\documentclass[
pdftex,
12pt,
a4paper,
twoside,
openright,
onecolumn,
bibliography=totoc,
fleqn,
headsepline,
titlepage
]{scrreprt}

\pdfoutput=1

\usepackage[a4paper,twoside]{geometry}

\usepackage[utf8]{inputenc}
\usepackage[T1]{fontenc}
\usepackage{lmodern}

\usepackage{amsmath}
\usepackage{amsfonts}
\usepackage{amssymb}

\usepackage{bm}
\usepackage{dlfltxbcodetips}
\usepackage{mathbbol}
\usepackage{stmaryrd}

\usepackage{appendix}
\usepackage{array}
\usepackage[english]{babel}
\usepackage{booktabs}
\usepackage{calc}
\usepackage{color}
\usepackage{enumerate}
\usepackage{float}
\usepackage{forloop}
\usepackage{framed}
\usepackage{graphicx}
\usepackage[unicode,bookmarks]{hyperref}
\usepackage{ifthen}
\usepackage{longtable}
\usepackage{makeidx}
\usepackage[numbers,square]{natbib}
\usepackage[neverdecrease]{paralist}
\usepackage{subfig}
\usepackage{tabularx}
\usepackage[absolute,overlay]{textpos}
\usepackage{tikz}
\usepackage{titletoc}
\usepackage[subfigure,titles]{tocloft}
\usepackage{url}
\usepackage{wasysym}
\usepackage{wrapfig}

\usepackage[automark,headsepline,plainheadsepline]{scrpage2}

\usetikzlibrary{3d}
\usetikzlibrary{arrows}
\usetikzlibrary{matrix}

\renewcaptionname{english}{\figurename}{Fig.}
\renewcaptionname{english}{\tablename}{Tab.}

\clearscrheadfoot

\ihead[\headmark]{\headmark}
\ohead[\pagemark]{\pagemark}

\setheadsepline{.4pt}

\setkomafont{caption}{\small\bfseries}
\setkomafont{captionlabel}{\small\bfseries}
\setcapindent{1em}

\definecolor{lightblue}{rgb}{.8 ,.85,1. }
\definecolor{midblue}  {rgb}{.15,.25,.65}
\definecolor{darkblue} {rgb}{.05,.1 ,.5 }
\definecolor{lightgray}{rgb}{.85,.85,.85}
\definecolor{darkgreen}{rgb}{0. ,.7 ,0. }
\definecolor{darkred}  {rgb}{.5 ,.1 ,.1 }
\definecolor{magenta}  {rgb}{1. ,0. ,.5 }
\definecolor{orange}   {rgb}{1. ,.5 ,0. }
\definecolor{yellow}   {rgb}{1. ,1. ,0. }

\numberwithin{figure}{chapter}
\numberwithin{equation}{chapter}

\setcapindent{0em}

\cftsetindents{section}{2em}{2.5em}
\cftsetindents{subsection}{3em}{2.5em}
\cftsetindents{subsubsection}{3em}{2.5em}

\titlecontents{chapter}[2em]
{\nopagebreak\vspace*{.2em}\hspace*{-2em}\rule{\textwidth}{.5mm}\vspace*{.3em}
\bfseries}
{\contentslabel{1.5em}}
{\hspace*{-1.5em}}
{\hfill\contentspage}
[\vspace*{-.1em}\rule{\textwidth}{.5mm}]

\newlength{\numwidth}
\setdefaultleftmargin{1.5em}{1.5em}{1em}{1em}{.5em}{.5em}

\newcolumntype{C}{>{\centering\arraybackslash}p}

\newenvironment{fshaded}{%
\MakeFramed {\FrameRestore}}{\endMakeFramed}

\newenvironment{impeq*}[1][\empty]{%
\definecolor{shadecolor}{rgb}{.8,.85,1}%
\definecolor{framecolor}{rgb}{1,1,1}%
\begin{fshaded}%
\ifthenelse{\equal{#1}{\empty}}{}{\textbf{#1}}\begin{align*}%
}{\end{align*}\end{fshaded}}

\newenvironment{example}[1][\empty]{%
\definecolor{framecolor}{rgb}{.85,.85,.85}%
\vspace{1em}
\begin{leftbar}%
\vspace{.25em}
\ifthenelse{\equal{#1}{\empty}}{}{\textbf{#1} \vspace{1em}\par}%
}{\vspace{.5em}
\end{leftbar}%
\vspace{.5em}
}

\newfloat{floaty}{thp}{lop}

\newenvironment{fexample}[1][\empty]{%
\definecolor{framecolor}{rgb}{.85,.85,.85}%
\begin{floaty}
\begin{leftbar}%
\vspace{.5em}
\ifthenelse{\equal{#1}{\empty}}{}{\textbf{#1} \vspace{1em}\par}%
}{\vspace{.5em}
\end{leftbar}%
\end{floaty}
}

\newcommand{\rimpeq}[2][\empty]{%
\definecolor{shadecolor}{rgb}{.8,.85,1}%
\definecolor{framecolor}{rgb}{1,1,1}%
\vspace{.75em}
\begin{tikzpicture}%
\hspace{-3pt}\node[fill=shadecolor, rounded corners=5pt, inner sep=3pt, align=justify, text width=\linewidth]{%
\ifthenelse{\equal{#1}{\empty}}%
{\vspace{-1em}\begin{align}#2\end{align}}%
{\textbf{#1}\begin{align}#2\end{align}}%
};%
\end{tikzpicture}%
\vspace{.5em}%
}

\newenvironment{rimpeq*}[2][\empty]{%
\definecolor{shadecolor}{rgb}{.8,.85,1}%
\definecolor{framecolor}{rgb}{1,1,1}%
\begin{tikzpicture}%
\node[fill=shadecolor,rounded corners=5pt, text width=\textwidth]{%
\ifthenelse{\equal{#1}{\empty}}{}{\textbf{#1}}\begin{equation*}}%
{\end{equation*}};%
\end{tikzpicture}%
}

\DeclareMathOperator{\sgn}{sgn}

\AtBeginDocument{

\let\obar\bar
\let\ohat\hat
\let\otilde\tilde

\let\bar\undefined
\let\hat\undefined
\let\tilde\undefined

\newcommand{\bar}[1]{\ensuremath{\overline{#1}}}						
\newcommand{\hat}[1]{\ensuremath{\widehat{#1}}}						
\newcommand{\tilde}[1]{\ensuremath{\widetilde{#1}}}					
}

\newcommand{\ten}[1]{\ensuremath{\boldsymbol{\mathsf{#1}}}}				
\newcommand{\form}[1]{\ensuremath{\boldsymbol{\mathsf{#1}}}}			

\newcommand{\mbb}[1]{\ensuremath{\mathbb{#1}}}							
\newcommand{\mcal}[1]{\ensuremath{\mathcal{#1}}}						
\newcommand{\msf}[1]{\ensuremath{\mathsf{#1}}}							
\newcommand{\mrm}[1]{\ensuremath{\mathrm{#1}}}							
\newcommand{\msp}[1]{\ensuremath{\mathbb{#1}}}							

\newcommand{\rsp}{\ensuremath{\mathbb{R}}}					

\newcommand{\mf}[1]{\ensuremath{\mathcal{#1}}}							
\newcommand{\tb}[2][\empty]{\ensuremath{\mathsf{T}_{#1}     #2}}		
\newcommand{\cb}[2][\empty]{\ensuremath{\mathsf{T}_{#1}^{*} #2}}		
\newcommand{\jb}[3][\empty]{\ensuremath{\mathsf{J}_{#1}^{#2} #3}}		

\newcommand{\lie}{\ensuremath{\pounds}}									
\newcommand{\ext}{\ensuremath{\boldsymbol{\msf{d}}}}					

\newcommand{\eps}{\ensuremath{\epsilon}}								
\newcommand{\veps}{\ensuremath{\varepsilon}}							
\newcommand{\phy}{\ensuremath{\varphi}}									

\newcommand{\abs}[1]{\ensuremath{\left \vert #1 \right \vert}}			
\newcommand{\norm}[1]{\ensuremath{\left \lVert #1 \right \rVert}}		

\newcommand{\contr}{\ensuremath{ \mathbin{\rule[0ex]{.4em}{.03em}\rule[0ex]{.03em}{1.5ex}} \, }}	
\newcommand{\iprod}{\ensuremath{\boldsymbol{\imath}}}					

\newcommand{\bracket}[1]{\ensuremath{\left < #1 \right >}}				

\newcommand{\degree}{\ensuremath{^\circ}}								

\let\div\undefined

\DeclareMathOperator{\div}{div}										

\DeclareMathOperator{\id}{id}										

\author{Michael Kraus \\ (michael.kraus@ipp.mpg.de)}
\title{Variational Integrators in Plasma Physics}

\begin{document}

\newgeometry{left=2cm,right=2cm,top=2.5cm,bottom=2.5cm}

\begin{titlepage}
\begin{flushleft}
\textsc{\LARGE{Michael Kraus}}\\
\rule{\textwidth}{3pt}\\
\vspace{0.05\textheight}
\textsf{\textbf{\Huge Variational Integrators in Plasma Physics}}
\end{flushleft}

\vspace{1cm}

\begin{textblock}{14}(1,3)
\begin{large}
\begin{tikzpicture}[scale=1]%

\draw [white] (0,18) node {.};

\draw [blue!20!white, draw opacity=0.5] (6.5,7.75) node {$\Omega_{L_{d}} = \dfrac{\partial^{2} L_{d}}{\partial q_{k} \, \partial q_{k+1}} (q_{k}, q_{k+1}) \, dq_{k} \wedge dq_{k+1}$};

\draw [blue!20!white, draw opacity=0.5] (10,3.5) node {$D_{1} L (q_{k}, q_{k+1}) + D_{2} L (q_{k-1}, q_{k}) = 0$};

\draw [blue!20!white, draw opacity=0.5] (10,14.75) node {$\sum \limits_{a} \dfrac{\partial \mcal{L}_{\square}}{\partial y^{a}} \big( \phy_{\square^{1}}, \phy_{\square^{2}}, \phy_{\square^{3}}, \phy_{\square^{4}} \big) = 0$};

\draw [blue!20!white, draw opacity=0.5] (-1,13.5) node {$\mcal{A}_{d} [\phy] = \sum \limits_{\square} \mcal{L}_{d} \circ j^{1} \phy (\square)$};

\draw [blue!20!white, draw opacity=0.5] (3,16.5) node {$\sum \limits_{\substack{\square\\ \square \cap \partial \mf{U} \neq \emptyset}} \bigg( \sum \limits_{\substack{a\\ \square^{a} \in \partial \mf{U}}} (j^{1} \phy)^{*} J^{a} (\square) \bigg) = 0$};

\draw [blue!20!white, draw opacity=0.5] (1, 5) node {$\ext \mcal{A}_{d} \cdot V = \sum \limits_{\substack{\square\\ \square \cap \mrm{int} \, \mf{U} \neq \emptyset}} \ext \mcal{L}_{\square} \cdot V + \theta_{L_{d}} (\phy) \cdot V$};

\draw [blue!20!white, draw opacity=0.5] (5,1) node {$\sum \limits_{\substack{\square\\ \square \cap \partial \mf{U} \neq \emptyset}} \bigg( \sum \limits_{\substack{a\\ \square^{a} \in \partial \mf{U}}} \Big[ (j^{1} \phy)^{*} ( V \contr W \contr \Omega_{L_{d}}^{a} \Big] (\square) \bigg) = 0$};

\end{tikzpicture}
\end{large}
\end{textblock}

\begin{center}
\includegraphics[width=.8\textwidth]{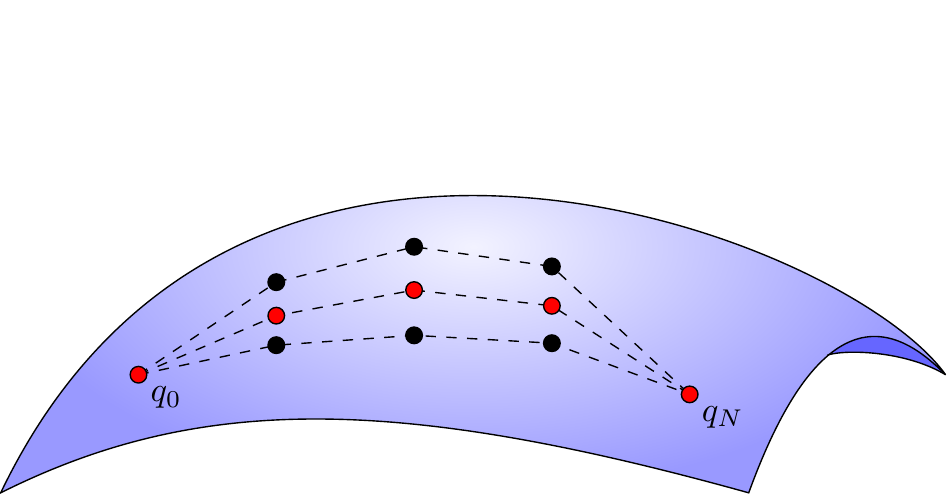}
\end{center}
\vfill
\rule{\textwidth}{3pt}

\vspace{0.2cm}
\begin{minipage}{0.14\textwidth}
\begin{flushleft}
\includegraphics[height=1.8cm]{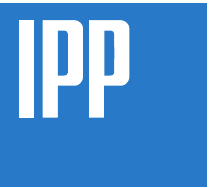}
\end{flushleft}
\end{minipage}
\hfill
\begin{minipage}{0.24\textwidth}
\begin{flushright}
\includegraphics[height=1.8cm]{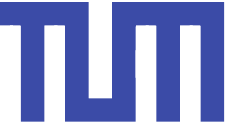}
\end{flushright}
\end{minipage}

\end{titlepage}

\restoregeometry

\begin{titlepage}

\begin{center}
\LARGE TECHNISCHE UNIVERSIT\"AT M\"UNCHEN

\vspace{3\lineskip}
\LARGE Fakult\"at f\"ur Mathematik\\

\vspace{3\lineskip}
\Large Lehrstuhl für Numerische Methoden der Plasmaphysik\\

\vspace{2.0cm}
\huge \textbf{Variational Integrators\\
in Plasma Physics}

\vspace{1.0cm}
\large \textbf{Michael Kraus}

\vspace{2.0cm}
Vollst\"andiger Abdruck der von der Fakult\"at f\"ur\\
Mathematik der Technischen Universit\"at M\"unchen\\
zur Erlangung des akademischen Grades eines\\
Doktors der Naturwissenschaften (Dr. rer. nat.)\\
genehmigten Dissertation.

\vspace{1.5cm}

\begin{tabular}{lrl}
Vorsitzende:			& & Univ.-Prof. Dr. Simone Warzel \\
\\
Pr\"ufer der Dissertation:
& 1. & Univ.-Prof. Dr. Eric Sonnendr\"ucker \\
& 2. & Univ.-Prof. Dr. Oliver Junge \\
& 3. & Prof. Philip J. Morrison, Ph.D., \\
&    & University of Texas at Austin/USA
\end{tabular}

\vspace{1.5cm}
Die Dissertation wurde am 03.04.2013 bei der\\
Technischen Universit\"at M\"unchen eingereicht\\
und durch die Fakult\"at f\"ur Mathematik\\
am 01.07.2013 angenommen.

\end{center}

\end{titlepage}

\newpage

\mbox{}
\thispagestyle{empty}
\newpage

\pagenumbering{roman}

\addcontentsline{toc}{section}{Abstract}
\section*{Abstract}

To a large extent, research in plasma physics is concerned with the description and analysis of energy and momentum transfer between different scales and different kinds of waves.
In the numerical modelling of such phenomena it appears to be crucial to describe the transfer processes preserving the underlying conservation laws in order to prevent physically spurious solutions.

In this work, special numerical methods, so called variational integrators, are developed for several models of plasma physics. Special attention is given to conservation properties like conservation of energy and momentum.

By design, variational integrators are applicable to all systems that have a Lagrangian formulation.
Usually, equations of motion are derived by Hamilton's action principle and then discretised. In the application of the variational integrator theory, the order of these steps is reversed.
At first, the Lagrangian and the accompanying variational principle are discretised, such that discrete equations of motion can be obtained directly by applying the discrete variational principle to the discrete Lagrangian.
The advantage of this approach is that the resulting discretisation automatically retains the conservation properties of the continuous system.

Following an overview of the geometric formulation of classical mechanics and field theory, which forms the basis of the variational integrator theory, variational integrators are introduced in a framework adapted to problems from plasma physics.
The applicability of variational integrators is explored for several important models of plasma physics: particle dynamics (guiding centre dynamics), kinetic theory (the Vlasov-Poisson system) and fluid theory (magnetohydrodynamics).

These systems, with the exception of guiding centre dynamics, do not possess a Lagrangian formulation to which the variational integrator methodology is directly applicable. Therefore the theory is extended by linking it to Ibragimov's theory of integrating factors and adjoint equations.
It allows us to find a Lagrangian for all ordinary and partial differential equations and systems thereof, for which the number of variables equals the number of equations.
Consequently, the applicability of the variational integrators is extended to a much larger family of systems as envisaged in the original theory.
The theory allows for the application of Noether's theorem to analyse the conservation properties of the system, both at the continuous and the discrete level.

In numerical examples, the conservation properties of the derived schemes are analysed.
In case of guiding centre dynamics, momentum in the toroidal direction of a tokamak is preserved exactly. The particle energy exhibits an error, but the absolute value of this error stays constant during the entire simulation. Therefore numerical dissipation is absent.
In case of the kinetic theory, the total number of particles, total linear momentum and total energy are preserved exactly, i.e., up to machine accuracy.
In case of magnetohydrodynamics, the total energy, cross helicity and the divergence of the magnetic field are preserved up to machine precision.

These conservation properties not only make the numerical schemes more stable than those obtained by traditional discretisation methods, but they also reduce unphysical behaviour like spurious loss of energy or momentum, thereby increasing the trustworthiness of numerical simulations.

\newpage

\mbox{}
\thispagestyle{empty}
\newpage

\addcontentsline{toc}{section}{Zusammenfassung}
\section*{Zusammenfassung}

Weite Teile der theoretischen Plasmaphysik beschäftigen sich mit der Beschreibung und Untersuchung des Transfers von Energie und Impuls zwischen verschiedenen Skalen und verschiedenen Arten von Wellen.
In der numerischen Modellierung dieser Phänomene erscheint es entscheidend, die den Transferprozessen zugrunde liegenden Erhaltungsgesetze zu bewahren, um unphysikalische Lösungen zu unterbinden.

In der vorliegenden Arbeit werden spezielle numerische Verfahren, sogenannte Variationsintegratoren, für verschiedene Modelle der Plasmaphysik entwickelt. Ein besonderer Augenmerk liegt dabei auf der Erhaltung physikalischer Größen wie Impuls und Energie.

Prinzipiell sind Variationsintegratoren auf alle Systeme anwendbar, die eine Lagrangesche Formulierung aufweisen. Die grundlegende Idee ist, dass man nicht, wie üblich, mit Hilfe des Hamiltonschen Variationsprinzips Bewegungsgleichungen ableitet und diese dann diskretisiert sondern die Abfolge dieser Schritte umkehrt. Man diskretisiert die Lagrange-Funktion und das Variationsprinzip und leitet damit direkt diskrete Bewegungsgleichungen ab.
Diese Vorgehensweise hat den Vorteil, dass man Diskretisierungen erhält, die automatisch die Erhaltungseigenschaften des kontinuierlichen Systems bewahren.

Im Folgenden wird zuerst die Theorie der Variationsintegratoren entwickelt. Da diese auf einer geometrischen Formulierung der klassischen Mechanik und Feldtheorie beruht, werden deren Grundlagen ebenfalls dargestellt. Anschließend wird die Anwendbarkeit der Variationsintegratoren auf verschiedene Systeme der Plasmaphysik untersucht: auf Teilchendynamik (Guiding Centre Dynamik), die kinetische Theorie (das Vlasov-Poisson-System) und die Fluiddynamik (Magnetohydrodynamik).

Da diese Systeme, mit Ausnahme der Guiding Centre Dynamik, keine passende Lagrangesche Formulierung aufweisen, wird die Theorie der Variationsintegratoren durch Kopplung mit Ibragimovs Theorie der integrierenden Faktoren und adjungierten Gleichungen erweitert. Diese erlaubt es, eine Lagrange-Funktion für alle gewöhnlichen und partiellen Differentialgleichungen und Differentialgleichungssysteme zu finden, bei denen die Anzahl der Variablen der Anzahl der Gleichungen entspricht. Dadurch ist es möglich, die Anwendbarkeit der Variationsintegratoren auf eine viel größere Familie von Systemen zu erweitern, als dies in der ursprünglichen Theorie vorgesehen ist.
Die Theorie ermöglicht die Anwendung des Noether-Theorems zur Untersuchung der Erhaltungseigenschaften des Systems, sowohl im Kontinuierlichen wie auch im Diskreten.

In numerischen Beispielen werden die Erhaltungseigenschaften der Variationsintegratoren untersucht. Im Falle der Guiding Centre Dynamik ist der Impuls in toroidaler Richtung im Tokamak exakt erhalten. Die Energie weist einen Fehler auf, dessen Absolutwert vom gewählten Zeitschritt abhängt, aber im Laufe einer Simulation konstant bleibt, d.h. es tritt keine numerische Dissipation auf.
Im Falle der kinetischen Theorie werden die Gesamtzahl der Teilchen, der Gesamtimpuls und die Gesamtenergie exakt erhalten (d.h. bis auf Maschinengenauigkeit). Im Falle der Magnetohydrodynamik wird ebenfalls die Gesamtenergie und zusätzlich die Kreuzhelizität und die Divergenz des Magnetfeldes exakt erhalten.

Diese Erhaltungseigenschaften führen zu Verfahren mit verbesserter numerischer Stabilität im Vergleich zu Diskretisierungen durch traditionelle Verfahren. Mindestens ebenso wichtig ist aber, dass sie unphysikalisches Verhalten wie die numerische Dissipation von Energie oder Impuls reduzieren und dadurch die Glaubwürdigkeit numerischer Ergebnisse erhöhen.

\tableofcontents
\clearpage{\pagestyle{empty}\cleardoublepage}

\chapter*{Notation}
\markboth{Notation}{}

\begin{longtable}{ c l }
$ \mcal{D} $		& Differential operator				\\
$ \mcal{A} $		& Action							\\
$ \mcal{A}_{d} $	& Discrete action					\\
$ H $				& Hamiltonian						\\
$ L $				& Lagrangian						\\
$ h $				& Particle Hamiltonian				\\
$ L_d $			& Discrete Lagrangian				\\
$ \mcal{H} $		& Hamiltonian density				\\
$ \mcal{L} $		& Lagrangian density				\\
$ \Theta $			& Canonical one-form				\\ 
$ \Omega $			& Symplectic two-form				\\
$ \Theta_{L} $		& Lagrangian one-form				\\
$ \Omega_{L} $		& Lagrangian two-form				\\
$ \Theta_{d} $		& Discrete one-form					\\
$ \Omega_{d} $		& Discrete two-form					\\
$ t $				& Time								\\
$ q $				& Generalised coordinates			\\
$ \dot{q} $		& Generalised velocities			\\
$ p $				& Generalised (conjugate) momenta	\\
$ E $				& Energy							\\
$ \mcal{E} $		& Energy functional					\\
$ f $				& Distribution function				\\
$ g $				& Ibragimov multiplier of $f$		\\
$ \phi $			& Electrostatic potential			\\
$ \psi $			& Ibragimov multiplier of $\phi$	\\
$ \Phi $			& Potential vorticity				\\
$ F $				& Functional of the distribution function $f$	\\
$ Z $				& $L^{2}$ norm of the distribution function $f$	\\
$ p $				& Kinetic gas pressure				\\
$ P $				& Generalised pressure				\\
$ V $				& Velocity field					\\
$ \omega $			& Vorticity							\\
$ \psi $			& Streaming function				\\
$ E $				& Electric field					\\
$ B $				& Magnetic field					\\
$ J $				& Current							\\
$ A $				& Magnetic vector potential			\\
$ A^{*} $			& Modified vector potential			\\
$ \msf{A} $		& Electromagnetic 4-potential		\\
$ \ten{F} $		& Electromagnetic field tensor		\\
$ \mcal{J} $						& Jacobian					\\
$ [ \cdot , \cdot ] $				& Poisson brackets			\\
$ \{ \cdot , \cdot \} $			& Lie-Poisson brackets		\\
$ [ \cdot , \cdot , \cdot ] $		& Nambu particle brackets	\\
$ \{ \cdot , \cdot , \cdot \} $	& Nambu field brackets		\\
$ \bracket{\cdot , \cdot} $		& Pairing					\\
$ r $				& Minor tokamak radius				\\
$ R $				& Major tokamak radius				\\
$ q $				& Safety factor						\\
$ \mu $			& Magnetic moment					\\
$ \phy $			& Toroidal angle					\\
$ \vartheta $		& Poloidal angle					\\
$ \mf{M} , \mf{N} $	& Manifold							\\
$ \Omega^{n} (\mf{M}) $ & $n$ forms on $\mf{M}$			\\
$ \mf{Q} $				& Configuration manifold with coordinates $(q)$					\\
$ \tb{Q} $				& Tangent bundle of $\mf{Q}$ with coordinates $(q, v)$, velocity phase manifold	\\
$ \cb{Q} $				& Cotangent bundle of $\mf{Q}$ with coordinates $(q, p)$			\\
$ \mf{C} ( \mf{Q} ) $							& Set of parametrised curves in $\mf{Q}$								\\
$ \mf{C} ( q_{1}, q_{2}, [t_{1}, t_{2}] ) $	& Set of parametrised curves in $\mf{Q}$ that start at $q_{1}$, end at $q_{2}$, and \\
& are parametrised by values in the interval $[t_{1}, t_{2}]$		\\
$ \mf{C}_{L} ( \mf{Q} ) $						& Set of parametrised curves in $\mf{Q}$ that are solutions of the \\
& Euler-Lagrange equations \\
$ \mf{C} ( \mf{Y} ) $							& Set of parametrised sections in $\mf{Y}$								\\
$ \mf{C}_{L} ( \mf{Y} ) $						& Set of parametrised sections in $\mf{Y}$ that are solutions of the \\
& Euler-Lagrange equations \\
$ \mf{X} $				& Base manifold (usually time, spacetime or phasespacetime)			\\
$ \mf{Y} $				& Fibre bundle										\\
$ \mf{F} $				& Typical fibre										\\
$ \mf{U}_{\mf{X}} $	& Closed submanifold of the base manifold $\mf{X}$	\\
$ \mf{U} $				& Parametrisation manifold of $\mf{U}_{\mf{X}}$		\\
$ \jb[]{k} $			& $k$'th jet bundle					\\
$ j^{k} $				& $k$'th jet prolongation			\\
$ \phy $				& Section							\\
$ V, X $				& Vector field						\\
$ \aleph $				& General geometric object (scalar function, vector field, differential form)	\\
$ \xi, \eta $			& Transformation map				\\
$ \eps, s, t $			& Group parameter					\\
$ \lie $				& Lie derivative					\\
$ \ext $				& Exterior derivative				\\
$ \wedge $				& Wedge product						\\
$ \iprod $				& Interior product					\\
$ \contr $				& Contraction						\\
$ \cong $				& is isomorphic to					\\
iff					& if and only if					\\
$ \mrm{int}(\mf{U}) $	& interior of $\mf{U}$				\\
$ \mrm{cl}(\mf{U}) $	& closure of $\mf{U}$				\\
$ \partial \mf{U} $	& boundary of $\mf{U}$				\\
\end{longtable}

Local deviations from this list cannot be excluded.
Some letters are defined twice, such as $\phy$ for the toroidal angle as well as for a general section. However, the current meaning should always be clear from the context in which the symbol is used.
Some symbols with varying meaning (like for general scalars, scalar fields, vector fields, differential forms, etc.) are not listed here or may overlap with notation defined here. \\

Coordinates indices $i,j,k$ run from $1$ to $n$, indices $\mu, \nu, \sigma$ run from $0$ to $n$,  where $n$ is the dimension of the space.
Indices $a,b$ correspond the components of fields or sections. \\

Throughout the whole thesis, all maps are assumed to be smooth, while all manifolds are assumed to be smooth as well as oriented. \\

\clearpage{\pagestyle{empty}\cleardoublepage}

\pagenumbering{arabic}

\chapter{Introduction}

Plasma physics is one of the most challenging fields in classical physics. Not only does it describe systems that consist of a vast number of particles, but these particles are charged and interact through the mean field they generate, leading to a collective behaviour and a tremendous complexity of the dynamics.
This complexity repeatedly provides us with new imponderabilities, not anticipated before.
To investigate the complicated behaviour inherent to any plasma system, pure theory is not sufficient.
Too great is the complexity of nowadays' problems to solve them by pen and paper alone.
Therefore computer simulations have become an essential part of plasma physics research, and have been for some time already.

With ever more powerful computers becoming available, ever larger simulations become feasible. Larger simulations in terms of the simulation domain but also in terms of the simulation time.
Consider for example simulations of the entire plasma of ITER, a new experimental device under construction in France, which has a volume of about $840 \, \mrm{m}^{3}$. And think of simulations of an entire plasma discharge, which in the case of ITER might last up to $400$ seconds, almost an eternity on the timescales of important plasma processes like small scale turbulence.

To be able to do such long simulations and still obtain accurate results, standard discretisation methods do not suffice. Most often they are based on the minimisation of local errors, but do not limit global error growth, thereby accumulating errors in each and every timestep, eventually leading to unphysical solutions.

Another example is turbulence, one of the large standing problems in classical physics and an ubiquitous topic in plasma physics. Its description involves the analysis of energy and momentum transfer processes between different scales and different kinds of waves. In the numerical modelling of such phenomena it appears to be crucial to describe these transfer processes while preserving the underlying conservation laws in order to prevent physically spurious solutions.
It cannot be expected to describe an energy cascade correctly, if energy is numerically created or dissipated.

Lastly, consider magnetic reconnection, a problem that will be dealt with in some more detail later on. It describes how magnetic field lines open up and reconnect in certain physical situations, resulting in a change of the topology of the magnetic field. Even under ideal conditions, that do not feature reconnection processes (i.e., the magnetic field line topology is fixed), most numerical schemes find reconnection events due to numerical dissipation.
If such methods are used to model real reconnection processes, one can never be absolutely certain to which extent the results are due to physical effects and to which extent they are just numerical artefacts.

\section{Geometric Discretisation}

To overcome these problems, the global structure of the equations, namely their geometry, has to be taken into account in the course of discretisation.
Following \citeauthor{Christiansen:2011} \cite{Christiansen:2011}, a geometric structure is a global property, that can be defined independently of particular coordinate representations of the differential equations at hand (see also \citeauthor{BuddPiggot:2000} \cite{BuddPiggot:2000}).
Examples for such structures encompass topology, like magnetic field line topology, conservation laws and symmetries, such as conservation of energy which arises through the invariance of a system under infinitesimal time translations, constraints like the divergence of the magnetic field which has to vanish, or identities like those from vector calculus and their generalisations from differential geometry.

The preservation of such geometric properties on the discrete level can have crucial influence on the quality of a simulation.
It affects stability and global error growth, reduces numerical artefacts, like spurious loss of energy or momentum, and thereby reduces the likelihood of inaccurate and unphysical behaviour. In difficult cases, simulations often only become possible by using geometric discretisation methods. This is especially true for long time simulations, where the unlimited growth of global errors, like in the energy of the system, can lead to numerical instabilities or at least physically wrong results.

\section{Symplecticity}

A geometric concept that plays an important role throughout this work is symplecticity. For a one-dimensional Hamiltonian system, the symplectic structure amounts to a skew-symmetric matrix that is a measure for phasespace area. Computing the product of this matrix with two phasespace vectors yields the area of the parallelogram spanned by the two vectors. For autonomous Hamiltonian systems, this area is always preserved. Consequently the Hamiltonian flow is called \emph{symplectic}. Maintaining this preservation of area on the discrete level has important consequences for the resulting integrators like very good energy behaviour (for more details see \citeauthor{SanzSernaCalvo:1994} \cite{SanzSernaCalvo:1994}, \citeauthor{LeimkuhlerReich:2004} \cite{LeimkuhlerReich:2004} and \citeauthor{HairerLubichWanner:2006} \cite{HairerLubichWanner:2006}).

In more than one dimension, the conserved quantity is the sum of the areas of the parallelograms that result by projecting the two phasespace vectors to the coordinate axes. As a consequence, phasespace volume is preserved under symplectic maps leading to other conservation laws like conservation of the total number of particles in a system.
In the framework of partial differential equations, the concept of symplecticity is generalised to multisymplecticity. Simply put, a multisymplectic map is symplectic with respect to both space and time.

It is noteworthy that a symplectic structure can also be defined on the Lagrangian side, indicating that the class of systems endowed with a symplectic structure is larger than the class of Hamiltonian systems. A fact that was already known to Lagrange and the details of which will be explained in chapter two.

\section{Variational Integrators}

One special geometric discretisation method is represented by variational integrators. They can be applied to any equation or system of equations that can be derived by means of a variational principle.
The general idea is simple. It can be described as discretising the theory instead of discretising the equations.
Part of the development of variational integrators was therefore the development of discrete counterparts of classical mechanics and classical field theory.
Although those are not complete counterparts, they are sufficient to derive geometric integration schemes and analyse their properties with respect to the observance of conservation laws.

In order to derive variational integrators one first has to discretise the basic constituents of the variational principle, the Lagrangian and the action integral. One has to approximate the particle positions or fields and their derivatives and select a quadrature rule.
Then a discrete variational principle is applied to the resulting discrete action, directly leading to discrete equations of motion (Euler-Lagrange equations).
There are several advantages of this method compared with a direct discretisation of the continuous Euler-Lagrange equations.
The obtained integrators preserve a discrete analogue of the symplectic or multisymplectic form (for finite or infinite dimensional systems, respectively).
This implies conservation of phasespace volume and a very good energy behaviour. In general, the energy is not preserved exactly, but it exhibits an oscillating behaviour about a fixed value. Consequently, the energy is not constant. But what is important is that the energy is not dissipated or growing unphysically, instead its error is bounded.
Furthermore, variational integrators conserve discrete momenta, that is conserved quantities corresponding to a symmetry of the system, practically exactly (up to machine precision).

Quite often, one can recover existing methods that are well known for their good conservative properties via a discrete variational principle. The Newmark scheme,
Störmer–Verlet, or symplectic Runge-Kutta methods are examples.

\section{Outline and Contributions}

In chapter two, an introduction to the geometric formulation of classical mechanics and classical field theory is given, together with an overview of the most important differential geometric tools.
The theory of variational principles is reviewed in a geometric setting. While the material presented in this section is not original, the presentation is detailed and self-contained. It should be accessible to non-specialists, applied mathematicians and theoretical physicists alike.

The theory of Ibragimov is presented. It allows us to find a variational formulation for certain systems that naturally do not have such a formulation. Furthermore, the Noether theorem can be applied in this framework to study symmetries and link them with conservation laws.
The combination of Ibragimov's theory with the discrete variational principle is a very important result of this work, as it allows to derive variational integrators for a much larger class of systems than had been foreseen in the original theory.

In chapter three, the basic theory and methodology of variational integrators is presented, both for finite dimensional systems (e.g., particle mechanics) and infinite dimensional systems (e.g., field theories). Proofs for the discrete conservation properties are given or at least sketched, including a discrete version of Noether's theorem.
The chapter is closed by an example, namely the advection equation. Ibragimov's theory is used to construct a Lagrangian and consecutively a variational integrator is derived. Both the continuous and the discrete Noether theorem are applied to obtain conservation laws for that equation.
Again, the general theory of variational integrators is not original and largely influenced by Marsden and coworkers \cite{MarsdenPatrick:1998, KouranbaevaShkoller:2000, MarsdenWest:2001, Lew:2003}, but the presentation is adapted to our framework.
As opposed to this, the application of variational integrators to Ibragimov's extended Lagrangians is proposed here for the first time. Its role for this work is crucial, since for most plasma physics problems a natural variational formulation in terms of Eulerian coordinates has not yet been found.

Chapters four to six explore the applicability of variational integrators to different systems which are important in plasma physics. Three classes of problems are distinguished: particle dynamics, kinetic theory, and plasma fluid theory.

In chapter four, several variational integrators are derived for the motion of particles in a non-uniform magnetic field. Specifically, the motion of centres of the helical trajectory of a particle (guiding centre) is considered, extending previous work in several aspects:
different discretisations of the Lagrangian are explored and dynamics in higher dimensions is considered.
We find that the variational integrators obtained here show excellent long-time behaviour, describing particle orbits correctly after millions of characteristic times and hundreds of millions of timesteps, while for standard methods like a fourth order Runge-Kutta scheme large deviations from the correct orbit are observed.
The chapter concludes with the sketch of a possible application of the derived integrators in particle-in-cell codes. Here, it is possible to employ a variational principle for the combined system of particles and fields, leading to schemes that consistently respect the conservation properties of the complete system.

In chapter five, variational integrators for the Vlasov-Poisson system in one-dimension are derived.
This is a typical test bed for kinetic problems, e.g., it has recently been employed to test the conservation properties of new schemes such as discontinuous Galerkin methods \cite{Ayuso:2009, Ayuso:2012a, Ayuso:2012b, HeathGamba:2012, ChengGamba:2012, ChengGamba:2013}.
Considering a one dimensional problem (one space plus one velocity coordinate) reduces the computational burden, while retaining the qualitative physical behaviour, including phase-mixing and collective effects.
One of the integrators for this system shows extraordinary conservation properties, preserving the total particle number, the total energy, total linear momentum and the $L^{2}$ norm exactly, i.e., up to machine accuracy. Problems only arise if the grid resolution is insufficient to resolve small scale structures in the distribution function.
To treat such cases, a velocity space collision operator is introduced. It dissipates the $L^{2}$ norm and removes subgrid modes, while retaining the conservation of the total particle number, energy and momentum.
Furthermore, a linear integrator is derived, which is computationally less demanding but keeps the conservation properties intact, albeit with less accuracy. Energy for example is not preserved to machine precision, but oscillating about some fixed value, as is typical for symplectic integrators. Still, no numerical dissipation is present.
The derived integrators are then applied to different standard benchmark cases like Landau damping, the twostream instability and the Jeans instability.

In chapter six, a variational discretisation of magnetohydrodynamics is obtained. The resulting integrator has similarly astonishing properties as the one for the Vlasov-Poisson system, namely exact conservation of the total energy and cross helicity. Here, a staggered grid approach has to be taken to avoid unphysical oscillations in the velocity and pressure fields, a typical problem in incompressible fluid dynamics.
The integrator is applied to a range of quite different examples like Alfvén waves, which appear to travel virtually forever through the computational domain, the passive advection of a magnetic loop by the velocity field, the emergence of current sheaths in the turbulent setting of a Orszag-Tang vortex, and several current sheath models as they are used in reconnection studies.

In the appendix, semi-discretisation strategies based on variational integrators or closely related methods are sketched.
In appendix A, a variational-spectral method for the vorticity equation and for the Vlasov-Poisson system is derived. Here, the spatial dimensions are transformed into Fourier space and only time or time and velocity are treated variationally.
In appendix B, discretisations of Poisson brackets and various generalisations thereof are considered. Here, only phasespace is discretised but not time. The derivations of these discretisations share many similarities with the derivation of variational integrators. Therefore it is not surprising that the resulting schemes are found to be similar.

\chapter{Geometric Mechanics and Field Theory}\label{ch:classical}

\textit{``Physicists have had a long-lasting love affair with the idea of generating physical laws by setting the derivative of some functional to zero. This is called an action principle. The most famous action principle is Hamilton’s principle, which produces Lagrange’s equations of mechanics upon variation.''}
- Philip Morrison \cite{Morrison:1998}

\vspace{2em}

In this chapter, a short overview of the geometric formulation of classical mechanics and field theory is given. It is important to understand some of the geometric underpinnings of the treated systems to appreciate the presented geometric discretisation methods that aim at preserving exactly these structures. Unfortunately, the geometric point of view is seldom treated in lectures on classical mechanics at university, nor are they common knowledge in the plasma physics community.
The elegance and beauty of the geometric formulation will certainly appeal to the reader yet unfamiliar with it.

We begin with the presentation of the geometric setting, i.e., some basic notions about manifolds, differential forms and fibre bundles. Hereafter, the formulation of Lagrangian mechanics and field theory is presented, at first in an analytic and then in a geometric way, applying the utilities introduced in the first section.
In this context, the theory of Ibragimov is reviewed. It allows us to find a Lagrangian for any ordinary or partial differential equation or any system of differential equations where the number of equations equals the number of dependent variables, which is usually the case in physical systems. Hence it allows us to find extended Lagrangian formulations for systems that do not posses a classical Lagrangian as it is often the case in plasma physics.

Some emphasis is put on Noether's theorem, which connects symmetries and conservation laws. Beginning with an analytic description of point transformations and one-parameter groups, the Noether theorem is presented for particles, fields and extended Lagrangian formulations according to Ibragimov. Hereafter, a geometric formulation of Noether's theorem is developed using the notion of momentum maps.

This treatment restrains itself mostly to the Lagrangian side, nevertheless connections with Hamiltonian mechanics and field theory are drawn to compare some results with discretisation methods developed on that side and to outline some alternative strategies. Central to these ideas are various kinds of brackets, namely the classical Poisson brackets and their generalisations in form of Nambu, Lie-Poisson and Dirac brackets.

If the reader is interested in more detailed treatments he can find some recommendations below. There are lots of classical as well as modern introductions to differential geometry and exterior calculus.
The more recent ones include \citeauthor{Lovett:2010} \cite{Lovett:2010}, \citeauthor{Castillo:2011} \cite{Castillo:2011} and \citeauthor{Epstein:2010} \cite{Epstein:2010} on the physics oriented side and \citeauthor{Lee:2012} \cite{Lee:2012} and \citeauthor{Tu:2011} \cite{Tu:2011} on the math oriented side.
Some classics are \citeauthor{Schutz:1980} \cite{Schutz:1980}, \citeauthor{Burke:1985} \cite{Burke:1985} and \citeauthor{AbrahamMarsdenRatiu:1988} \cite{AbrahamMarsdenRatiu:1988}.
A short and nevertheless comprehensive overview of differential forms are the lecture notes by \citeauthor{Sjamaar:2006} \cite{Sjamaar:2006} which are freely available on the internet.
Good introductions can quite often also be found in general relativity textbooks, e.g., \citeauthor{Ryder:2009} \cite{Ryder:2009}, \citeauthor{Hobson:2006} \cite{Hobson:2006} and \citeauthor{Carroll:2003} \cite{Carroll:2003}.

A basic introduction to the geometric formulation of classical mechanics is \citeauthor{JoseSaletan:1998} \cite{JoseSaletan:1998}.
More advanced treatments are \citeauthor{Holm:2009} \cite{Holm:2009} and \citeauthor{MarsdenRatiu:2002} \cite{MarsdenRatiu:2002}.
Some classics that are still very useful today are \citeauthor{Arnold:1989} \cite{Arnold:1989}, \citeauthor{AbrahamMarsden:1978} \cite{AbrahamMarsden:1978}, \citeauthor{SudarshanMukunda:1974} \cite{SudarshanMukunda:1974} and \citeauthor{SaletanCromer:1971} \cite{SaletanCromer:1971}.
Freely available lectures by \citeauthor{Holm:2011} \cite{Holm:2011}, \citeauthor{Marsden:2009} \cite{Marsden:2009} and \citeauthor{Ratiu:2005} \cite{Ratiu:2005} can be found on the internet.
Two review papers that focus on problems related to fluid dynamics but also have an introductory character are \citeauthor{Morrison:1998} \cite{Morrison:1998} and \citeauthor{Salmon:1988} \cite{Salmon:1988}.

The theory of jet bundles is introduced in the monographs of \citeauthor{KrasilshchikVinogradov:1999} \cite{KrasilshchikVinogradov:1999}, \citeauthor{Olver:1995} \cite{Olver:1995} and \citeauthor{Saunders:1989} \cite{Saunders:1989}.
Its application to classical mechanics and field theories is explained in \citeauthor{GotayMarsden:1998} \cite{GotayMarsden:1998}, \citeauthor{MarsdenPatrick:1998} \cite{MarsdenPatrick:1998}, \citeauthor{Aldaya:1980} \cite{Aldaya:1980}, \citeauthor{Echeverria-Enriquez:1996} \cite{Echeverria-Enriquez:1996, Echeverria-Enriquez:2000}, \citeauthor{Giachetta:2010} \cite{Giachetta:2010, Giachetta:1997, Giachetta:2009}, \citeauthor{Sardanashvily:2009} \cite{Sardanashvily:2009} and references therein.
The survey articles by \citeauthor{Saunders:2008} \cite{Saunders:2008} and \citeauthor{Krupka:1973} \cite{Krupka:1973} are also instructive.

This chapter is largely influenced by all of the aforementioned references. It makes no claim of originality, except for presenting the material in a mostly self-contained and coherent way.
However, the following presentation, especially the sections after the geometric introduction, should be easier accessible, as the level of detail and explanation often surpasses that of the original works which are quite challenging at times.

\section{Geometric Foundations}

This section tries to give a short overview of the geometric foundations underlying the theory presented below.
All of the geometric tools that are used later on should be covered, the only exception being basic Lie group theory.

\subsection{Smooth Manifolds}

Probably the simplest definition of a manifold $\mf{M}$ is a set of points that can be labelled by coordinates.
Locally, manifolds look like the Euclidean space. Globally, however, they might have a much more complicated structure.
Therefore it is often not possible to define a global coordinate system on a manifold (think of a circle or a sphere) and one has to find coordinate patches (charts) that together cover the whole manifold.
A \emph{chart} (local coordinate system) is a pair $(\mf{U}, \phi)$, where $\mf{U}$ is an open subset of $\mf{M}$ and $\phi$ is a one-to-one map from $\mf{U}$ onto some open subset of $\rsp^{n}$
\begin{align}
\phi (p) &= \big( x^{1} (p), x^{2} (p) , ..., x^{n} (p) \big) , &
p &\in \mf{U} . & &&
\end{align}

Hence, a chart labels each point $p$ in $\mf{U}$ by $n$ real numbers.
If more than one chart is necessary to cover the whole manifold, most likely some points will lie in the domain of more than one chart. In that case we demand that there exists a transition map as follows.
If $(\mf{U}_{1}, \phi_{1})$ and $(\mf{U}_{2}, \phi_{2})$ are two coordinate patches overlapping in $\mf{U} = \mf{U}_{1} \cup \mf{U}_{2}$, we request that $\psi = \phi_{2} \circ \phi_{1}^{-1}$ is smooth. In that case, the charts $\phi_{1}$ and $\phi_{2}$ are said to be compatible.
The set of compatible charts that covers all of a manifold is called an \emph{atlas}.

Consider as an example the configuration space $\mf{Q}$ of a mechanical system. The definition of a coordinate chart on $\mf{Q}$ amounts to a choice of generalised coordinates.
Fortunately, for the cases considered in this work, it is always possible to find a global coordinate patch, thereby avoiding the subtleties arising from having more than one coordinate patch.

Most manifolds in mathematical physics are smooth manifolds, continuous and infinitely often differentiable. We will always assume that this is the case.
Furthermore, we shall assume that all of the considered manifolds are orientable.

In the remainder of this section we will consider some intrinsic objects and operations that can be defined on manifolds and will be used in the subsequent treatment.

\subsection{Vector Fields}\label{sec:geometry_vector_fields}

On trivial manifolds, i.e., such manifolds that can be identified with a linear vector space like $\rsp^{n}$, the definition of vector fields is straight forward.
Indeed, it is customary to identify points $p$ of the space with the corresponding vector $x = (x^{1}, ..., x^{n})$, leading to the usual notion of vectors.
General manifolds, however, are not necessarily linear, so vectors cannot be defined by the usual means.
The simplest geometric way to describe a vector $V$ at a point $p$ on a nontrivial manifold is intuitively as the tangent to a parametrised curve $c(t)$ in $\mf{M}$, satisfying $c(0) = p$.

\subsubsection{Tangent Vectors}

A parametrised curve $c(t)$ in $\mf{M}$ is a smooth map from some interval $\mcal{I} \subseteq \rsp$ to the manifold $\mf{M}$
\begin{align}\label{eq:geometry_vectors_tangent_vectors_1}
c : \mcal{I} \rightarrow \mf{M} .
\end{align}

If coordinates on $\mf{M}$ are denoted $(x^{\mu})$, this can be explicitly written as
\begin{align}\label{eq:geometry_vectors_tangent_vectors_2}
c : t \mapsto x^{\mu} (t) .
\end{align}

Without loss of generality assume that $\mcal{I}$ contains the point $0 \in \rsp$ and that $c(0) = p$.
Consider the directional derivative of a function $f : \mf{M} \rightarrow \rsp$ along the curve $c$, that is
\begin{align}\label{eq:geometry_vectors_tangent_vectors_3}
V_{p} (f) = \dfrac{d}{dt} \Big[ f \circ c (t) \Big] \bigg\vert_{t=0} .
\end{align}

For the trivial case, $\mf{M} = \rsp^{n}$, this corresponds to
\begin{align}\label{eq:geometry_vectors_tangent_vectors_4}
V_{p} (f) &= \dfrac{dx^{\mu}}{dt} \, \dfrac{\partial f}{\partial x^{\mu}} &
& \text{with} &
\dfrac{dx^{\mu}}{dt} &\equiv V_{p}^{\mu} . & &&
\end{align}

$V_{p}^{\mu}$ are the components of the tangent vector of $c(t)$ at $t=0$.
In the general case, this is used as a definition.
As (\ref{eq:geometry_vectors_tangent_vectors_4}) is fully general and independent of $f$, the vector $V_{p}$ can be written as
\begin{align}\label{eq:geometry_vectors_tangent_vectors_5}
V_{p} &= V_{p}^{\mu} \, \partial_{\mu} &
& \text{where} &
\partial_{\mu} &\equiv \dfrac{\partial}{\partial x^{\mu}} . & &&
\end{align}

Vectors on a manifold therefore correspond to first order differential operators.
$\partial_{\mu}$ are the local basis in which the vector components are expressed.

\subsubsection{Tangent Bundle}

The tangent vector $V_{p}$ is an element of the local \emph{tangent space} $\tb[p]{\mf{M}}$, where $\tb[p]{\mf{M}}$ is the set of all tangent vectors (i.e., all possible directional derivatives) to $\mf{M}$ at $p$ and has the same dimension as $\mf{M}$.
That means, $\tb[p]{\mf{M}}$ can be obtained by considering the tangents to all possible curves passing through that point.
Coordinates $x^{\mu}$ on $\mf{M}$ induce a basis $\partial_{\mu}$ on $\tb[p]{\mf{M}}$.
Therefore $(\partial_{\mu})$ build a \emph{natural coordinate system} on $\tb[p]{\mf{M}}$.
One possible way of obtaining $\tb[p]{\mf{M}}$ is to consider all possible curves $c$ through $p$ and evaluate (\ref{eq:geometry_vectors_tangent_vectors_5}) for each single one.

An important consequence of vectors at different points $p$ of the manifold being elements of different vector spaces $\tb[p]{\mf{M}}$ is that they cannot be added or subtracted. This is only possible for vectors at the same point $p$, i.e., vectors which are elements of the same tangent space $\tb[p]{\mf{M}}$.

Collecting all the $\tb[p]{\mf{M}}$, for each point $p \in \mf{M}$, into one single object (which is a disjoint union) yields the \emph{tangent bundle}
\begin{align}\label{eq:geometry_vectors_tangent_bundle}
\tb{\mf{M}} = \bigcup \limits_{p \in \mf{M}} \tb[p]{\mf{M}} .
\end{align}

It is the set of all tangent vectors at all points of $\mf{M}$ and has the structure of a differentiable manifold.
More details on the tangent bundle will be presented in section \ref{sec:geometry_fibre_bundles} on fibre bundles.

A vector field $V$ on a manifold $\mf{M}$ is a function that assigns a vector $V_{p} \in \tb[p]{\mf{M}}$ to each point $p \in \mf{M}$.
Consequently, all vector fields $V$ on $\mf{M}$ lie in $\tb{\mf{M}}$.
A vector field $V$ is therefore a smooth, linear map
\begin{align}
V : \mf{M} \rightarrow \tb{\mf{M}} .
\end{align}

The tangent space $\tb[p]{\mf{M}}$ at each point $p \in \mf{M}$ is a real vector space.
Therefore two vector fields $V$ and $W$ may be added or multiplied by a scalar field $f : \mf{M} \rightarrow \rsp$ as follows
\begin{align}
(V + W) (p) &= V (p) + W (p) , &
(f V) (p) &= f(p) \, V (p) , &
& p \in \mf{M} . & &&
\end{align}

In the physics literature, the $V^{\mu}$ are usually referred to as contravariant components of the vector field $V$.

\subsection{Integral Curves and Flows}

We will now study a topic that is essential in the geometric formulation of the action principle and the study of symmetries.
It is based on the observation that vector fields induce, at least locally, a family of transformations of the manifold onto itself.

A one-parameter family of transformations is a differentiable map
\begin{align}
\phy : \mf{M} \times \rsp \rightarrow \mf{M}
\end{align}

that depends on a real parameter. Therefore they are called one-parameter-groups of transformations. They map points of $\mf{M}$ onto different points of $\mf{M}$
\begin{align}
\phy &: p \mapsto \phy (p, t) &
& \text{with} &
p &\in \mf{M} , \quad
t \in \rsp , & &&
\end{align}

such that
\begin{align}
\phy (p, 0) &= p &
& \text{and} &
\phy \big( \phy (p, s) , t \big) &= \phy ( p , s+t) &
& \text{for all} &
p &\in \mf{M} , \quad
s,t \in \rsp . & &&
\end{align}

Thus, upon defining $\phy_{t} (p) \equiv \phy (p, t)$ we can write
\begin{align}
\phy_{s+t} (p) &= \phy_{s} \circ \phy_{t} = \phy_{t} \circ \phy_{s} &
& \text{and} &
\phy_{0} &= \id . & &&
\end{align}

As
\begin{align}
\phy_{t} \circ \phy_{-t} = \phy_{-t} \circ \phy_{t} = \phy_{0} = \id
\end{align}

each map $\phy_{t}$ has an inverse $\phy_{t}^{-1} = \phy_{-t}$ that is also differentiable.
Therefore, each $\phy_{t}$ is a diffeomorphism of $\mf{M}$ onto itself, and the set of transformations $\{ \phy_{t} \; \big\vert \; t \in \rsp \}$ is a group of diffeomorphisms of $\mf{M}$ onto itself.

Each one-parameter-group of transformations $\phy$ on $\mf{M}$ determines a family of curves in $\mf{M}$ (referred to as the \emph{orbits} of the group).
The map
\begin{align}
\phy_{p} : \rsp \rightarrow \mf{M}
\end{align}

given by
\begin{align}
\phy_{p} (t) = \phy (p, t)
\end{align}

is a differentiable curve in $\mf{M}$ for each $p \in \mf{M}$.
The vector field tangent to these curves generated by the one-parameter-group of transformations
\begin{align}
V = \dfrac{d}{dt} \, \phy_{t} \, \bigg\vert_{t=0}
\end{align}

is called the \emph{infinitesimal generator} of $\phy$.
Since $\phy_{p} (0) = \phy (p, 0) = p$ the tangent vector to the curve $\phy_{p}$ belongs to $\tb[p]{\mf{M}}$.
The curves $\phy_{p}$ are integral curves of $V$.

\subsection{Fibre Bundles}\label{sec:geometry_fibre_bundles}

\begin{figure}[H]
\centering

\subfloat{
\begin{tikzpicture}[scale=1.4]%

\coordinate (A1) at (+0.5, +0.86602540378443871);
\coordinate (B1) at (-0.5, +0.86602540378443871);
\coordinate (C1) at (+0.5, -0.86602540378443871);
\coordinate (D1) at (-0.5, -0.86602540378443871);

\coordinate (A2) at (+0.86602540378443871, +0.5);
\coordinate (B2) at (+0.86602540378443871, -0.5);
\coordinate (C2) at (-0.86602540378443871, +0.5);
\coordinate (D2) at (-0.86602540378443871, -0.5);

\draw [line width=0.5mm, color=blue!60!white]	(0,+1) -- ++(  0:0.8cm);
\draw [line width=0.5mm, color=blue!60!white]	(0,+1) -- ++(180:0.8cm);

\draw [line width=0.5mm, color=blue!60!white]	(0,-1) -- ++(  0:0.8cm);
\draw [line width=0.5mm, color=blue!60!white]	(0,-1) -- ++(180:0.8cm);

\draw [line width=0.5mm, color=blue!60!white]	(+1,0) -- ++( 90:0.8cm);
\draw [line width=0.5mm, color=blue!60!white]	(+1,0) -- ++(270:0.8cm);

\draw [line width=0.5mm, color=blue!60!white]	(-1,0) -- ++( 90:0.8cm);
\draw [line width=0.5mm, color=blue!60!white]	(-1,0) -- ++(270:0.8cm);

\draw [line width=0.5mm, color=blue!60!white]	(A1) -- ++(150:0.8cm);
\draw [line width=0.5mm, color=blue!60!white]	(A1) -- ++(330:0.8cm);

\draw [line width=0.5mm, color=blue!60!white]	(B1) -- ++( 30:0.8cm);
\draw [line width=0.5mm, color=blue!60!white]	(B1) -- ++(210:0.8cm);

\draw [line width=0.5mm, color=blue!60!white]	(C1) -- ++(210:0.8cm);
\draw [line width=0.5mm, color=blue!60!white]	(C1) -- ++( 30:0.8cm);

\draw [line width=0.5mm, color=blue!60!white]	(D1) -- ++(330:0.8cm);
\draw [line width=0.5mm, color=blue!60!white]	(D1) -- ++(150:0.8cm);

\draw [line width=0.5mm, color=blue!60!white]	(A2) -- ++(120:0.8cm);
\draw [line width=0.5mm, color=blue!60!white]	(A2) -- ++(300:0.8cm);

\draw [line width=0.5mm, color=blue!60!white]	(B2) -- ++( 60:0.8cm);
\draw [line width=0.5mm, color=blue!60!white]	(B2) -- ++(240:0.8cm);

\draw [line width=0.5mm, color=blue!60!white]	(C2) -- ++( 60:0.8cm);
\draw [line width=0.5mm, color=blue!60!white]	(C2) -- ++(240:0.8cm);

\draw [line width=0.5mm, color=blue!60!white]	(D2) -- ++( 120:0.8cm);
\draw [line width=0.5mm, color=blue!60!white]	(D2) -- ++( 300:0.8cm);

\draw [line width=0.5mm]	(0,0)	circle	(+1.);

\filldraw [color=blue!60!white]		(0,+1) circle (.05);
\filldraw [color=blue!60!white]		(0,-1) circle (.05);
\filldraw [color=blue!60!white]		(+1,0) circle (.05);
\filldraw [color=blue!60!white]		(-1,0) circle (.05);

\filldraw [color=blue!60!white]		(A1) circle (.05);
\filldraw [color=blue!60!white]		(B1) circle (.05);
\filldraw [color=blue!60!white]		(C1) circle (.05);
\filldraw [color=blue!60!white]		(D1) circle (.05);

\filldraw [color=blue!60!white]		(A2) circle (.05);
\filldraw [color=blue!60!white]		(B2) circle (.05);
\filldraw [color=blue!60!white]		(C2) circle (.05);
\filldraw [color=blue!60!white]		(D2) circle (.05);

\end{tikzpicture}
}
\subfloat{\hspace{5em}}
\subfloat{
\begin{tikzpicture}[x={(-1cm, 0.5cm)}, y={(+1cm, 0.5cm)}, z={(0cm,1cm)}, scale=1.25]%

\begin{scope}[canvas is xy plane at z=0]
\draw [line width=0.5mm]	(0,0)	circle	(+1.);

\end{scope}

\begin{scope}[canvas is xz plane at y=0]
\draw [line width=0.5mm, color=blue!60!white]	(1,-1)	-- (1,+1);
\end{scope}

\begin{scope}[canvas is xz plane at y=1]
\draw [line width=0.5mm, color=blue!60!white]	(0,-1)	-- (0,+1);
\end{scope}

\begin{scope}[canvas is xz plane at y=+0.5]
\draw [line width=0.5mm, color=blue!60!white]	(+0.86602540378443871,-1)	-- (+0.86602540378443871,+1);
\end{scope}

\begin{scope}[canvas is xz plane at y=+0.86602540378443871]
\draw [line width=0.5mm, color=blue!60!white]	(+0.5,-1)	-- (+0.5,+1);
\end{scope}

\begin{scope}[canvas is xz plane at y=-0.70710678118654757]
\draw [line width=0.5mm, color=blue!60!white]	(+0.70710678118654757,-1)	-- (+0.70710678118654757,+1);
\end{scope}

\begin{scope}[canvas is xz plane at y=+0.70710678118654757]
\draw [line width=0.5mm, color=blue!60!white]	(-0.70710678118654757,-1)	-- (-0.70710678118654757,+1);
\end{scope}

\begin{scope}[canvas is xz plane at y=-0.70710678118654757]
\draw [line width=0.5mm, color=blue!60!white]	(-0.70710678118654757,-1)	-- (-0.70710678118654757,+1);
\end{scope}

\begin{scope}[canvas is xz plane at y=-0.96592582628906831]
\draw [line width=0.5mm, color=blue!60!white]	(-0.25881904510252074,-1)	-- (-0.25881904510252074,+1);
\draw [line width=0.5mm, color=blue!60!white]	(+0.25881904510252074,-1)	-- (+0.25881904510252074,+1);
\end{scope}

\begin{scope}[canvas is xz plane at y=+0.25881904510252074]
\draw [line width=0.5mm, color=blue!60!white]	(-0.96592582628906831,-1)	-- (-0.96592582628906831,+1);
\end{scope}

\begin{scope}[canvas is xz plane at y=-0.25881904510252074]
\draw [line width=0.5mm, color=blue!60!white]	(-0.96592582628906831,-1)	-- (-0.96592582628906831,+1);
\end{scope}

\end{tikzpicture}
}

\caption{Left: Tangent bundle of the circle $\msp{S}^{1}$. Each of the blue lines attached to a point of the circle depicts a fibre of the tangent bundle $\tb{\msp{S}^{1}}$. Right: To avoid spurious intersections, the same tangent bundle $\tb{\msp{S}^{1}}$ is drawn with the fibres parallel to each other. Now the circle is depicted in the horizontal plane and the fibres $\tb[p]{\msp{S}^{1}}$ are vertical lines.}
\label{fig:geometry_fibre_bundles1}
\end{figure}
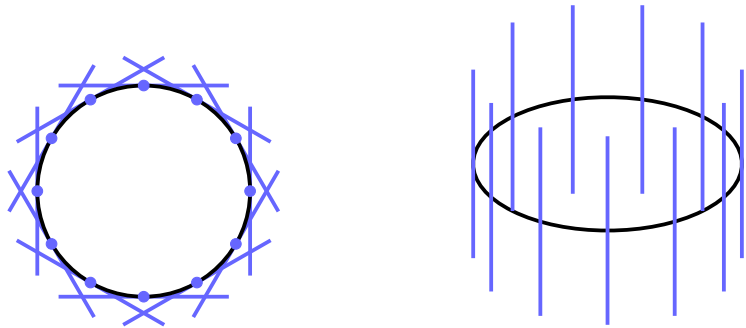

Reconsider the construction of the tangent bundle $\tb{\mf{M}}$ of a manifold $\mf{M}$ from section \ref{sec:geometry_vector_fields}.
$\tb{\mf{M}}$ was built by attaching to each point $p \in \mf{M}$ the tangent vector space $\tb[p]{\mf{M}}$ at that point (\ref{eq:geometry_vectors_tangent_bundle}).
The resulting object is generally referred to as a \emph{fibre bundle} with the vector spaces $\tb[p]{\mf{M}}$ being the \emph{fibres} that are attached to each point $p$ of the \emph{base space} $\mf{M}$.
See figure \ref{fig:geometry_fibre_bundles1} for a pictorial view of the tangent bundle of the circle $\msp{S}^{1}$.

In general, a fibre bundle is characterised by four quantities, the base space $\mf{X}$, the total space $\mf{Y}$, a typical fibre $\mf{F}$ and a projection $\pi$.
Fibre bundles are particular manifolds with the property of being decomposable into fibres. The points of a single fibre are related to one another while points of different fibres are not. This is formalised by defining a \emph{projection map}.

\subsubsection{Projections}

Consider a fibre bundle $\mf{Y}$ over $\mf{X}$, a point $x$ in the base space $\mf{X}$, and the fibre $\mf{Y}_{x} \cong \mf{F}$ at that point.
The \emph{natural} (or \emph{canonical}) \emph{projection} $\pi$ maps each element $y$ of $\mf{Y}_{x}$ to the point $x$, the fibre is attached to, i.e.,
\begin{align}
\pi_{\mf{X}} : \mf{Y} \rightarrow \mf{X}
\end{align}

or in coordinates
\begin{align}
\pi_{\mf{X}} : (x, y) \mapsto x .
\end{align}

To avoid confusion, e.g., if more than one projection appears in a treatment, we sometimes also write $\pi_{\mf{X} \mf{Y}}$ to denote both the source and the target space explicitly.

\subsubsection{Sections}

A section $\phy$ of a fibre bundle $\mf{Y}$ over $\mf{X}$ is a map that assigns to each point $x$ in the base manifold a point in the fibre bundle
\begin{align}
\phy &: \mf{X} \rightarrow \mf{Y} &
& \text{such that} &
\pi_{\mf{X} \mf{Y}} \circ \phy &= \id_{\mf{X}} , & &&
\end{align}

or in coordinates
\begin{align}
\phy : x \mapsto (x, y) .
\end{align}

Sections $\phy : \mf{X} \rightarrow \mf{Y}$ are also denoted as $\phy \in \Gamma (\pi_{\mf{X} \mf{Y}})$, where $\Gamma (\pi_{\mf{X} \mf{Y}})$ denotes the set of all sections in the fibre bundle $\mf{Y}$ over $\mf{X}$, characterised by the projection $\pi_{\mf{X} \mf{Y}}$.
Similarly, the fibre bundle itself is often just denoted as $\pi_{\mf{X} \mf{Y}}$.
If the model space $\mf{F}$ for fibres is a linear space, then $\Gamma (\pi_{\mf{X} \mf{Y}})$ has a natural linear space structure as well.

In this framework, vector fields are sections of $\tb{\mf{M}}$.
Later on, we will introduce jet bundles that are very practical when dealing with field theories. These are fibre bundles defined over spacetime that contain fields and their derivatives up to a given order.

\begin{figure}[H]
\centering

\begin{tikzpicture}[scale=1.25]%

\draw [line width=0.5mm, color=black]			(-3,0) -- (+3,0);

\foreach \y in {-2.50,-2.25,...,+2.50} {
\draw [line width=0.5mm, color=blue!60!white]	(\y,-0.5) -- (\y,+2);
}

\draw [line width=0.75mm, color=blue!80!white] plot [smooth, tension=.7] coordinates { (-2.75, 1.5) (-1.5, 1.7) (-0.2, 1.1) (+1.3, 0.8) (+2.75, 1.2)};

\end{tikzpicture}

\caption{A section in the tangent bundle $\tb{\rsp}$ of the real line.}
\end{figure}

\subsubsection{Vertical Lifts}

Some fibre bundles $\pi_{\mf{X} \mf{Y}}$, e.g., tangent bundles or jet bundles, have a natural way of \emph{lifting} a curve in $\mf{X}$ by a section $y$.
The tangent lift
\begin{align}
\tb{f} : \tb{\mf{M}} \rightarrow \tb{\mf{N}}
\end{align}

of a function
\begin{align}
f : \mf{M} \rightarrow \mf{N}
\end{align}

is defined by
\begin{align}
V_{p} \mapsto Df \, V_{p} .
\end{align}

The lift of a curve $c(t)$ in $\mf{X}$ to the tangent bundle $\tb{\mf{X}}$ is given as
\begin{align}
\tb{c} (t) : t \mapsto \bigg( x^{\mu} (t) , \dfrac{d x^{\mu}}{d t} (t) \bigg) .
\end{align}

Obviously, not every curve in the tangent bundle is the lift of a curve in the base space. Those curves that are lifted ones are called \emph{holonomic}.
If $\phy_{t}$ is a one-parameter group of diffeomorphisms on the basespace $\mf{X}$, then $\tb{\phy_{t}}$ will be a one-parameter group of diffeomorphisms on the tangent bundle $\tb{\mf{X}}$.
Thus, if $V$ is the infinitesimal generator of $\phy_{t}$, its tangent lift to $\tb{(\tb{\mf{X}})}$ is given as
\begin{align}
\tb{V} = \dfrac{d}{dt} \, \tb{\phy_{t}} \, \bigg\vert_{t=0}
\end{align}

in coordinates
\begin{align}
\tb{V} : (x^{\mu}, V^{\mu}) \mapsto \bigg( \bigg( x^{\mu} , \dfrac{d x^{\mu}}{d t} \bigg) , \bigg( \dfrac{d x^{\mu}}{d t} , \dfrac{d^{2} x^{\mu}}{d t^{2}} \bigg) \bigg) .
\end{align}

The space of such vector fields over $\tb{\mf{X}}$ plays an important role in the next section on Lagrangian dynamics.

\subsection{Differential Forms}

\textit{``Very loosely speaking, differential forms are whatever appears under an integral sign.''}

\hfill - Loring Tu \cite{Tu:2011}

\vspace{1em}

Somewhat less loosely speaking, they allow us to define integrands over manifolds, thereby providing a natural way to integrate over curved spaces.
Just as vector fields, they are intrinsic objects associated to any manifold, and in fact their simplest instance, differential one-forms,  are the dual concept to vector fields (which are first order differential operators).

\subsubsection{Differential One-Forms}

A differential one-form $\alpha$ (hereafter just referred to as one-form) on a manifold $\mf{M}$ assigns an element of the dual space $\cb[p]{\mf{M}}$ of the tangent space $\tb[p]{\mf{M}}$ to each point $p \in \mf{M}$. It is therefore a smooth, linear map
\begin{align}\label{eq:geometry_forms_one_1}
\alpha : \mf{M} \rightarrow \cb{\mf{M}} .
\end{align}

Each such one-form takes values in the \emph{cotangent space} $\cb[p]{\mf{M}}$ at that point, where $\cb[p]{\mf{M}}$ contains all dual vectors (covariant vectors) at that point $p$ and has the same dimension as $\mf{M}$. Collecting all the $\cb[p]{\mf{M}}$ for each point $p \in \mf{M}$ into one single object gives the \emph{cotangent bundle}
\begin{align}\label{eq:geometry_forms_one_2}
\cb{\mf{M}} = \bigcup \limits_{p \in \mf{M}} \cb[p]{\mf{M}} ,
\end{align}

such that the one-form $\alpha$ is a section of $\cb{\mf{M}}$. That way, $\cb{\mf{M}}$ is also the dual to $\tb{\mf{M}}$.
The space $\cb[p]{\mf{M}}$ is a vector space, such that two forms $\alpha$ and $\beta$ of the same order may be added or multiplied by a scalar field $f : \mf{M} \rightarrow \rsp$
\begin{align}\label{eq:geometry_forms_one_3}
(\alpha + \beta) (p) &= \alpha (p) + \beta (p) , &
(f \beta) (p) &= f(p) \, \beta (p) , &
& p \in \mf{M} . & &&
\end{align}

The simplest geometric way to describe a one-form $\alpha$ is as the differential of a function $f$ on $\mf{M}$.
Assign to each point $p$ of the curve local coordinates $x^{\mu} (t)$ and a scalar function $f(x^{\mu})$ on $\mf{M}$.
The differential of the function $f$ at $p$ is
\begin{align}\label{eq:geometry_forms_one_4}
df = \dfrac{\partial f}{\partial x^{\mu}} \, dx^{\mu} .
\end{align}

The first expression corresponds to the components of the gradient of $f$, and $(dx^{\mu})$ forms a local basis, dual to the basis $(\partial_{\mu})$ on $\tb{\mf{M}}$.
However, not all one-forms are differentials of a function.
In general, one-forms are written as
\begin{align}\label{eq:geometry_forms_one_5}
\alpha = \alpha_{\mu} \, dx^{\mu} ,
\end{align}

where $(dx^{\mu})$ is the aforementioned basis on $\cb{\mf{M}}$, defined by letting a basis one-form act on a basis vector, i.e.,
\begin{align}\label{eq:geometry_forms_one_6}
dx^{\mu} (\partial_{\nu}) = \dfrac{\partial x^{\mu}}{\partial x^{\nu}} = \delta^{\mu}_{\nu} ,
\end{align}

such that the result of a general one-form $\alpha$ acting on a general vector $v$ is given by
\begin{align}\label{eq:geometry_forms_one_7}
\alpha (v) = \alpha_{\mu} \, v^{\nu} \, \delta^{\mu}_{\nu} = \alpha_{\mu} \, v^{\mu} .
\end{align}

A one-form $\alpha$ is a linear functional, such that it acts on a linear combination of vectors $v, w \in \tb{\mf{M}}$ with scalars $a, b \in \rsp$ as
\begin{align}\label{eq:geometry_forms_one_8}
\alpha ( av + bw) = a \, \alpha(v) + b \, \alpha(w) .
\end{align}

In the physics literature, the $\alpha_{\mu}$ are usually referred to as covariant components of the covector field $\alpha$. Strictly speaking, one-forms can only be identified with covector fields if the underlying manifold is endowed with a metric, which defines a canonical isomorphism of $\cb[p]{\mf{M}}$ and $\tb[p]{\mf{M}}$, thus identifying vectors and there duals. In physical applications this is indeed most often the case.
To change between vectors and one-forms, the $\flat$ and $\sharp$ operators can be defined like
\begin{align}\label{eq:geometry_forms_one_9}
V^{\flat} &= V_{i} dx^{i} &
& \text{and} &
\alpha_{\sharp} &= \alpha^{i} \, \partial_{i} . & &&
\end{align}

The flat operator $\flat$ returns the one-form corresponding to a vector field, and the sharp operator $\sharp$ returns the vector field corresponding to a one-form. Therefore, the action is the same as in music, but with respect to indices instead of notes.

\subsubsection{Higher Order Differential Forms}

A differential two-form $\omega$ is a function
\begin{align}\label{eq:geometry_forms_higher_order_1}
\omega : \mf{M} \rightarrow \Omega^{2} (\mf{M})
\end{align}

where $\Omega^{2} (\mf{M})$ is the space of two-forms on $\mf{M}$.
It is generally written\footnote{
The factor of $1/2$ originates from the fact that in the sum over indices all contributions are taken into account twice.
Whether it is written or not depends on notational convention.
}
\begin{align}\label{eq:geometry_forms_higher_order_2}
\omega = \dfrac{1}{2} \, \omega_{\mu \nu} \, dx^{\mu} \wedge dx^{\nu} .
\end{align}

where $\wedge$ denotes the wedge product, which is defined in the next paragraph.
Two-forms are antisymmetric, such that
\begin{align}\label{eq:geometry_forms_higher_order_3}
\omega_{\mu \nu} = - \omega_{\nu \mu} .
\end{align}

A similar result as (\ref{eq:geometry_forms_one_8}) holds also for a two-form $\omega$, which is a bilinear functional, acting on vectors $v, w, z \in \tb{\mf{M}}$ with scalars $a,b \in \rsp$ as
\begin{align}\label{eq:geometry_forms_higher_order_4}
\omega ( av + bw, z) &= a \, \omega(v, z) + b \, \omega(w, z) &
& \text{and} &
\omega ( v, aw + bz) &= a \, \omega(v, w) + b \, \omega(v, z) .
\end{align}

Differential $n$-forms (differential forms of order $n$) are completely covariant, totally antisymmetric tensors.
A $n$-form $\theta$ is a function
\begin{align}\label{eq:geometry_forms_higher_order_5}
\theta : \mf{M} \rightarrow  \Omega^{n} (\mf{M}) ,
\end{align}

where $\Omega^{n} (\mf{M})$ is the space of $n$ forms on $\mf{M}$.
As a consequence of the antisymmetry property (\ref{eq:geometry_forms_higher_order_3}), the highest order forms that can exist on a manifold $\mf{M}$ of dimension $m$ are of order $m$.
The generalisation of (\ref{eq:geometry_forms_higher_order_4}) to higher order forms is straight forward.

\begin{example}[Examples: Differential Forms in a Three-Dimensional Manifold]

1-Form:
\begin{align}
A = A_{\mu} \, dx^{\mu} = A_{1} \, dx^{1} + A_{2} \, dx^{2} + A_{3} \, dx^{3}
\end{align}

2-Form:
\begin{align}
F = \dfrac{1}{2} \, F_{\mu \nu} \, dx^{\mu} \wedge dx^{\nu} = F_{12} \, dx^{1} \wedge dx^{2} + F_{23} \, dx^{2} \wedge dx^{3} + F_{31} \, dx^{3} \wedge dx^{1}
\end{align}

3-Form:
\begin{align}
\Omega = \dfrac{1}{3!} \, \Omega_{\mu \nu \sigma} \, dx^{\mu} \wedge dx^{\nu} \wedge dx^{\sigma} = \Omega_{123} \, dx^{1} \wedge dx^{2} \wedge dx^{3}
\end{align}

\vspace{-1em}

\end{example}

\vspace{-1em}

\subsubsection{Wedge Product}

The wedge product takes a $p$-form $\xi$ and a $q$-form $\eta$ and returns a $(p+q)$-form
\begin{align}
\xi \wedge \eta \, (v_{1}, v_{2}, ..., v_{p+q}) = \dfrac{1}{(p+q)!} \sum \limits_{\sigma \in S_{p+q}} \sgn (\sigma) \, \xi (v_{\sigma_{1}}, v_{\sigma_{2}}, ..., v_{\sigma_{p}}) \, \eta (v_{\sigma_{p+1}}, v_{\sigma_{p+2}}, ..., v_{\sigma_{p+q}}) ,
\end{align}

where $(\sigma_{1}, ..., \sigma_{p}, \sigma_{p+1}, ..., \sigma_{p+q})$ is an element of $S_{p+q}$, the group of all permutations of the numbers $\{ 1, 2, ..., p+q \}$, and $\sgn (\sigma)$ is the sign of the permutation, i.e.,
\begin{align}
\sgn (\sigma) =
\begin{cases}
1 & \text{odd permutation} , \\
0 & \text{even permutation} ,
\end{cases}
\end{align}

such that
\begin{align}\label{eq:geometry_forms_wedge_product1}
\xi \wedge \eta = (-1)^{pq} \, \eta \wedge \xi .
\end{align}

It is associative,
\begin{align}
( \alpha \wedge \beta ) \wedge \gamma = \alpha \wedge ( \beta \wedge \gamma ) ,
\end{align}

and bilinear,
\begin{align}
\begin{split}
(a \, \alpha + b \, \beta ) \wedge \gamma &= a \, ( \alpha \wedge \gamma) + b \, ( \beta \wedge \gamma) , \\
\alpha \wedge ( b \, \beta + c \, \gamma ) &= b \, ( \alpha \wedge \beta) + c \, ( \alpha \wedge \gamma) ,
\end{split}
\end{align}

but in general not commutative.
Due to the antisymmetry property (\ref{eq:geometry_forms_higher_order_3}), the wedge product of a basis form with itself vanishes,
\begin{align}
dx^{\mu} \wedge dx^{\mu} = 0 .
\end{align}

For that reason there can be no forms of higher order than the dimensionality of the space they are defined on.

\begin{example}[Examples: Wedge Products of Differential Forms]

Consider the examples from above, again defined on a three-dimensional manifold,
\begin{subequations}
\begin{align}
A &= A_{1} \, dx^{1} + A_{2} \, dx^{2} + A_{3} \, dx^{3} , \\
F &= F_{12} \, dx^{1} \wedge dx^{2} + F_{23} \, dx^{2} \wedge dx^{3} + F_{31} \, dx^{3} \wedge dx^{1} , \\
\Omega &= \Omega_{123} \, dx^{1} \wedge dx^{2} \wedge dx^{3} .
\end{align}
\end{subequations}

The wedge product of $A$ with itself is
\begin{align}
A \wedge A
\nonumber
&= A_{1} A_{2} \, dx^{1} \wedge dx^{2} + A_{2} A_{3} \, dx^{2} \wedge dx^{3} + A_{3} A_{1} \, dx^{3} \wedge dx^{1} \\
&+ A_{2} A_{1} \, dx^{2} \wedge dx^{1} + A_{3} A_{2} \, dx^{3} \wedge dx^{2} + A_{1} A_{3} \, dx^{1} \wedge dx^{3}
= 0 ,
\end{align}

which is obvious as for one-forms $A \wedge A = - A \wedge A$ (\ref{eq:geometry_forms_wedge_product1}).
The wedge product of $A$ and $F$ is
\begin{align}
A \wedge F
&= ( A_{1} F_{23} + A_{2} F_{31} + A_{3} F_{12} ) \, dx^{1} \wedge dx^{2} \wedge dx^{3} .
\end{align}

Let us try to compute the wedge product of $\Omega$ with the basis forms $dx^{\mu}$, that is
\begin{subequations}
\begin{align}
\Omega \wedge dx^{1} &= \Omega_{123} \, dx^{1} \wedge dx^{2} \wedge dx^{3} \wedge dx^{1} = \hphantom{-} \Omega_{123} \, \underbrace{dx^{1} \wedge dx^{1}}_{=0} \wedge dx^{2} \wedge dx^{3} = 0 , \\
\Omega \wedge dx^{2} &= \Omega_{123} \, dx^{1} \wedge dx^{2} \wedge dx^{3} \wedge dx^{2} = - \Omega_{123} \, dx^{1} \wedge \underbrace{dx^{2} \wedge dx^{2}}_{=0} \wedge dx^{3} = 0 , \\
\Omega \wedge dx^{3} &= \Omega_{123} \, dx^{1} \wedge dx^{2} \wedge \underbrace{dx^{3} \wedge dx^{3}}_{=0} = 0 .
\end{align}
\end{subequations}

We see that all of these vanish, which is no surprise as $\Omega$ is a form of maximum order.

\end{example}

\subsubsection{Interior Product}

One-forms $\alpha$ are linear functionals that map vector fields to functions
\begin{align}
\alpha : \tb{\mf{M}} \rightarrow \rsp .
\end{align}

In general, $n$-forms $\theta$ are $n$-linear functionals, mapping $n$ vector fields to functions
\begin{align}
\theta : \bigotimes \limits_{n} \tb{\mf{M}} \rightarrow \rsp .
\end{align}

The \emph{interior product} of a vector field $v$ and a one-form $\alpha$ is defined as their \emph{contraction}, denoted
\begin{align}
\iprod_{v} \alpha = v \contr \alpha = \bracket{ \alpha , v } .
\end{align}

The interior product $\iprod_{v}$ of some $n$-form $\theta$ and a vector field $v$ yields a $(n-1)$-form $\iprod_{v} \theta$,
\begin{align}
\iprod_{v} \theta ( \underbrace{ \, \dots \dots \, }_{\text{$n$ slots}} ) = \theta ( \, \overbrace{ v, \underbrace{ \, \dots \dots \, }_{\text{$n-1$ slots}} }^{\text{$n$ slots}} \, ) .
\end{align}

It can therefore be regarded as a map
\begin{align}
\iprod_{v} : \Omega^{n} (\mf{M}) \rightarrow \Omega^{n-1} (\mf{M})
\end{align}

or component-wise
\begin{align}
\iprod_{v} : \theta^{j_{1} \, j_{2} ... \, j_{p}} \mapsto \theta_{k \, j_{2} ... \, j_{p}} \, v^{k} .
\end{align}

The interior product of a vector field $v$ and a scalar function $f$ is zero by definition.
$\iprod_{v}$ is an anti-derivation: for a $p$-form $\xi$ and a $q$-form $\eta$
\begin{align}
\iprod_{v} ( \xi \wedge \eta ) = ( \iprod_{v} \xi ) \wedge \eta + (-1)^{p} \, \xi \wedge ( \iprod_{v} \eta ) .
\end{align}

For example, for $v = v^{\sigma} \partial_{\sigma}$ and the two-form $dx^{\mu} \wedge dx^{\nu}$, we have
\begin{align}
\iprod_{v} ( dx^{\mu} \wedge dx^{\nu} )
&= ( \iprod_{v} dx^{\mu} ) \, dx^{\nu} - dx^{\mu} \, ( \iprod_{v} dx^{\nu} )
= v^{\mu} \, dx^{\nu} - v^{\nu} \, dx^{\mu} .
\end{align}

\subsubsection{Exterior Derivative}

The exterior derivative $\ext$ maps $n$-forms into $(n+1)$-forms
\begin{align}
\ext : \Omega^{n} (\mf{M}) \rightarrow \Omega^{n+1} (\mf{M}) ,
\end{align}

thus taking functions (which are considered zero-forms) to one-forms, one-forms to two-forms, and so on.
It is axiomatically defined as follows.
If $f : \mf{M} \rightarrow \rsp$ is a function (zero-form), than $\ext f$ is the ordinary differential
\begin{align}
\ext f (v) = vf = v^{\mu} \partial_{\mu} f ,
\end{align}

equivalently
\begin{align}
\ext f = (\partial_{\mu} f) \, dx^{\mu} .
\end{align}

$\ext$ is an anti-derivation, i.e., if $\xi$ is a $p$-form and $\eta$ a $q$-form, than
\begin{align}
\ext ( \xi \wedge \eta ) = \ext \xi \wedge \eta + (-1)^{p} \, \xi \wedge \ext \eta .
\end{align}

When applied twice, the exterior derivative vanishes, i.e., $\ext^{2} = 0$ or $\ext (\ext \alpha) = 0$ for any $n$-form $\theta$.
The exterior derivative $\ext$ is linear, such that for every $a \in \rsp$
\begin{align}
\ext (a \alpha) &= a \, \ext \alpha &
& \text{and} &
\ext (\alpha + \beta) &= \ext \alpha + \ext \beta . & &&
\end{align}

The vanishing of the exterior derivative when applied twice, $\ext^{2} = 0$, leads to the notion of closed and exact forms.
A $n$-form $\theta$ is \emph{closed} if $\ext \theta = 0$.
A $n$-form $\theta$ is \emph{exact} if $\theta = \ext \eta$ for a $(n-1)$-form $\eta$.
An exact form is always closed, but a closed form is not necessarily exact.

\begin{example}[Example: Electromagnetic Field]

The natural description of the magnetic potential $A$ is as a one-form
\begin{align}
A = A_{1} \, dx^{1} + A_{2} \, dx^{2} + A_{3} \, dx^{3} .
\end{align}

The exterior derivative of $A$,
\begin{align}
\nonumber \ext A
&= \left( \dfrac{\partial A_{2}}{\partial x^{1}} - \dfrac{\partial A_{1}}{\partial x^{2}} \right) dx^{1} \wedge dx^{2}
+ \left( \dfrac{\partial A_{1}}{\partial x^{3}} - \dfrac{\partial A_{3}}{\partial x^{1}} \right) dx^{3} \wedge dx^{1}
+ \left( \dfrac{\partial A_{3}}{\partial x^{2}} - \dfrac{\partial A_{2}}{\partial x^{3}} \right) dx^{2} \wedge dx^{3} \\
&= \dfrac{1}{2} \, \big( \partial_{\mu} A_{\nu} - \partial_{\nu} A_{\mu} \big) \, dx^{\mu} \wedge dx^{\nu} \\
&\equiv \dfrac{1}{2} \, F_{\mu \nu} \, dx^{\mu} \wedge dx^{\nu} ,
\end{align}

yields the magnetic field tensor
\begin{align}
F = F_{12} \, dx^{1} \wedge dx^{2} + F_{23} \, dx^{2} \wedge dx^{3} + F_{31} \, dx^{3} \wedge dx^{1} .
\end{align}

Interestingly, $\ext A$ looks like a curl, and indeed, the components of $F$ correspond to the components of the magnetic field $B = \nabla \times A$,
\begin{align}
F =
\begin{pmatrix}
\hphantom{-} 0 & \hphantom{-} B_{3} & - B_{2} \\
-B_{3} & 0 & \hphantom{-} B_{1} \\
\hphantom{-} B_{2} & - B_{1} & \hphantom{-} 0
\end{pmatrix}
.
\end{align}

Therefore the natural representation of the magnetic field is a two-form. 
The exterior derivative of $F$,
\begin{align}
\ext F = \left( \dfrac{\partial F_{23}}{\partial x^{1}} + \dfrac{\partial F_{31}}{\partial x^{2}} + \dfrac{\partial F_{12}}{\partial x^{3}} \right) dx^{1} \wedge dx^{2} \wedge dx^{3} ,
\end{align}

does of course vanish as $\ext F = \ext^{2} A = 0$.
Interestingly, $\ext F$ looks like a divergence and indeed, it corresponds to $\nabla \cdot B = 0$.

\end{example}

These examples show that on a three-dimensional manifold, the exterior derivative corresponds to the operators from vector calculus. The exterior derivative of a zero-form corresponds to the gradient, the exterior derivative of a one-form corresponds to the curl, and the exterior derivative of a two-form corresponds to the divergence.
On manifolds of dimension other than three, the exterior derivative provides a generalisation of these operators.

\subsection{Pullback}

In all considerations of this subsection, $\phy$ is regarded as a diffeomorphism
\begin{align}
\phy : \mf{M} \rightarrow \mf{N} .
\end{align}

The pullback of $\phy$ allows us to \emph{pull back} geometric objects from the target manifold $\mf{N}$ to the source manifold $\mf{M}$. This includes functions, vector fields and differential forms.
The pullback of a scalar field $f : \mf{N} \rightarrow \rsp$ by $\phy$ is given by composition
\begin{align}
\phy^{*} f = f \circ \phy .
\end{align}

The result is a scalar field $\phy^{*} f : \mf{M} \rightarrow \rsp$.
The pullback of a $n$-form $\omega \in \Omega^{n} (\mf{N})$ by $\phy$ is a $n$-form $\phy^{*} \omega \in \Omega^{n} (\mf{M})$, defined point-wise by
\begin{align}
\big( \phy^{*} \omega \big)_{p} \big( v_{1}, ..., v_{n} \big) &= \omega_{\phy (p)} \big( D \phy (p) \cdot v_{1}, ..., D \phy (p) \cdot v_{n} \big) , &
p &\in \mf{N} . & &&
\end{align}

While $\phy^{*} \omega$ is acting on vectors $v_{i} \in \tb[p]{\mf{M}}$, $\omega$ is acting on vectors $D \phy (p) \cdot  v_{i} \in \tb[\phy(p)]{\mf{N}}$.
The following diagram should help to clarify this.

\begin{figure}[H]
\centering
\begin{large}
\begin{tikzpicture}
\node (m11) at ( 0,  0) {$\tb{\mf{M}}$};
\node (m12) at ( 3,  0) {$\tb{\mf{N}}$};

\node (m21) at ( 0, -2) {$\mf{M}$};
\node (m22) at ( 3, -2) {$\mf{N}$};

\node (m31) at ( 0, -2.5) {$\in$};
\node (m32) at ( 3, -2.5) {$\in$};

\node (m41) at ( 0, -3) {$\phy^{*} \omega$};
\node (m42) at ( 3, -3) {$         \omega$};

\path[-stealth, line width=.4mm]
(m11) edge node [above] {$\tb{\phy}$} (m12);
\path[-stealth, line width=.4mm]
(m21) edge node [above] {$\phy$}      (m22);
\path[-stealth, line width=.4mm]
(m11) edge node {} (m21);
\path[-stealth, line width=.4mm]
(m12) edge node {} (m22);
\end{tikzpicture}
\end{large}
\end{figure}

The pullback of a wedge product is the wedge product of the pullback
\begin{align}
\phy^{*} ( \alpha \wedge \beta ) = (\phy^{*} \alpha) \wedge (\phy^{*} \beta) .
\end{align}

The pullback of an exterior derivative is the exterior derivative of the pullback
\begin{align}
\phy^{*} ( \ext \omega ) = \ext ( \phy^{*} \omega )
\end{align}

where $\omega$ is any differential form.

\subsection{Lie Derivative}\label{sec:geometry_lie_derivative}

As already pointed out, on a manifold $\mf{M}$ it is generally not possible to add or subtract vectors at different points $p \in \mf{M}$, as those live in different vector spaces.
This causes a problem as for differentiation in the usual sense, one needs to do exactly that.
Nevertheless, there are ways to define derivatives on manifolds. Otherwise, they would not be fun to deal with.
The arguably most important derivative on manifolds is the Lie derivative.
It describes how a geometric object (a function, a vector, a form) changes when it is dragged along some vector field $X$.
It will be defined below in two different approaches that have been shown to be equivalent \cite{Holm:2009, MarsdenRatiu:2002}.

\subsubsection{Algebraic Definition}

The Lie derivative $\lie_{X}$ along a vector field $X \equiv X^{\mu} \partial_{\mu}$ is defined through its action on a scalar function $f$,
\begin{align}
\lie_{X} f \equiv X f = X^{\mu} \, f_{,\mu} ,
\end{align}

a vector field $Y \equiv Y^{\nu} \partial_{\nu}$,
\begin{align}
\lie_{X} Y \equiv XY - YX = \left( X^{\mu} ( \partial_{\mu} Y^{\nu} ) - Y^{\mu} ( \partial_{\mu} X^{\nu} ) \right) \partial_{\nu} ,
\end{align}

where in the second identity, the products $XY$ and $YX$ are viewed as composition of differential operators.
The Lie derivative of anything else is defined through the requirement that it is a derivative, i.e., that it fulfils Leibniz' rule.
The Lie derivative of differential forms can be obtained by following this rule. The result is a very beautiful relation called Cartan's magic formula or \emph{infinitesimal homotopy relation}

\rimpeq{\label{eq:geometry_lie_derivative_cartans_magic_formula}
\lie_{X} \omega &= \ext \left( \iprod_{X} \omega \right) + \iprod_{X} \left( \ext \omega \right) &
& \text{(\textbf{Cartan's Magic Formula})} .
}

\subsubsection{Dynamical Definition}

An alternative definition of the Lie derivative along a vector field $X$ with flow $\phy_{t}$ is given by
\begin{align}\label{eq:geometry_lie_derivative_dynamical_definition}
\lie_{X} \aleph = \dfrac{d}{dt} \phy_{t}^{*} \aleph \bigg\vert_{t=0} ,
\end{align}

where $\aleph$ can now be a scalar function, a vector field or a differential form.
So this definition, referred to as dynamical definition of the Lie derivative, is formally the same for all geometric entities (see for example \citeauthor{MarsdenRatiu:2002} \cite{MarsdenRatiu:2002} or \citeauthor{Holm:2009} \cite{Holm:2009}).

\subsubsection{Properties}

The Lie derivative does not change the tensorial character of the object it is acting on, e.g., a scalar stays a scalar, a vector stays a vector, a one-form stays a one-form, and so on.
It commutes with the exterior derivative,
\begin{align}
\lie_{X} \ext \omega = \ext \left( \lie_{X} \omega \right) ,
\end{align}

and if $\phy : \mf{M} \rightarrow \mf{N}$ is a diffeomorphism, the pullback of the Lie derivative by $\phy$ is
\begin{align}
\phy^{*} \lie_{X} \omega = \lie_{\phy^{*} X} ( \phy^{*} \omega ) .
\end{align}

Now we should be well equipped to approach the geometric formulation of Lagrangian dynamics.

\section{Lagrangian Dynamics}

At the age of 19, Lagrange found a solution to the long-standing isoperimetric problem\footnote{
Historical notes according to \citeauthor{Holm:2011} \cite{Holm:2011}.
}$^{,}$\footnote{
The isoperimetric problem asks, among all closed surfaces of a given fixed perimeter in the plane, which curve maximises the area that it encloses?
Lagrange sent his solution to this problem to Euler in 1755.
}. As it turned out, more important than the answer to this special problem was Lagrange's solution method, which lead to what we now call the \emph{Euler-Lagrange equations}
\begin{align}\label{eq:classical_lagrangian_euler_lagrange_equations}
\dfrac{\partial L}{\partial q} - \dfrac{d}{dt} \dfrac{\partial L}{\partial \dot{q}} = 0 ,
\end{align}

where $L (q, \dot{q})$ is the Lagrangian function, which often corresponds to the kinetic energy minus the potential energy, $q$ are generalised coordinates and $\dot{q}$ generalised velocities.
The great advantage of this formulation is that it is completely covariant. It does not depend on a specific coordinate representation.

Some years later, at the age of 18, Hamilton found that these equations can be derived by the \emph{principle of stationary action}\footnote{
Quite often Hamilton's principle is called ``principle of least action'', which is misleading. In fact, the action does not need to take a minimum but just a critical point. For the derivation of the equations of motion, it doesn't make a difference if the critical point is a minimum, a maximum or a saddle point. Admittedly, most often it is indeed a minimum, but there are counter-examples as well (e.g., under certain conditions the action of the harmonic oscillator takes neither a minimum nor a maximum).}.
It states that, considering all possible trajectories $q(t)$ a system could follow to get from state $a$ to state $b$, the following integral, called the action,
\begin{align}
\mcal{A} [q] = \int \limits_{a}^{b} L \big( q (t), \dot{q} (t) \big) \, dt ,
\end{align}

is stationary for the actual physical trajectory $q(t)$.
This means that the variation of $\mcal{A}$,
\begin{align}
\delta \mcal{A} [q] = \delta \int \limits_{a}^{b} L \big( q (t), \dot{q} (t) \big) \, dt ,
\end{align}

vanishes for the trajectory $q(t)$ that is actually taken by the system, i.e., the trajectory $q(t)$ that fulfils the Euler-Lagrange equations (\ref{eq:classical_lagrangian_euler_lagrange_equations}).

Starting from this description, the generalisation of classical mechanics to field theories is mostly straight forward.
The Lagrangian generalises from a function of position $q$ and velocity $\dot{q}$ to a function of the coordinates $x^{\mu}$ (independent variables; most often spacetime), the fields $\phy^{a} (x)$ (dependent variables) and their derivatives with respect to the coordinates.
The envisaged applications are all first order theories, i.e., their Lagrangians depend only on first order derivatives of the fields, $\phy^{a}_{\mu} = \partial \phy^{a} / \partial x^{\mu}$, and are thus of the form
\begin{align}
L = L \big( x^{\mu}, \phy^{a} (x), \phy^{a}_{\mu} (x) \big) .
\end{align}

In classical field theory, the Lagrangian density $\mcal{L}$ is often preferred over the Lagrangian function $L$, as it allows for more general notation.
The connection between the two,
\begin{align}
\mcal{L} = L \, \omega,
\end{align}

is drawn by the volume form $\omega$ of the base space (e.g., $\omega = dt \wedge dx \wedge dy \wedge dz$ for spacetime).
Finally, the action becomes an integral not only over time but over the whole base space $\mf{X}$,
\begin{align}
\mcal{A} = \int \limits_{\mf{X}} \mcal{L} (x^{\mu}, \phy^{a}, \phy^{a}_{\mu}) .
\end{align}

In the case of spacetime this is
\begin{align}
\mcal{A} = \int \limits_{\mf{X}} L  (x^{\mu}, \phy^{a}, \phy^{a}_{\mu}) \, dt \, dx \, dy \, dz .
\end{align}

In this section, a thorough derivation of the Euler-Lagrange equations based on Hamilton's action principle is presented. At first from an analytic point of view, utilising the notion of one-parameter families of transformations, thereby staying close to what is usually taught at university classes in classical mechanics.
This is followed by a presentation from a geometric point of view, namely on tangent and cotangent bundles, which is a much more natural description of the problem.
After a short comment on the phasespace Lagrangian, a popular object in the description of reduced kinetic theories in plasma physics, the theory on jet bundles is outlined. It allows to unify the theory for finite dimensional and infinite dimensional systems in one single framework and has many advantages over the tangent bundle theory.
Nevertheless, we also include the tangent bundle theory, which has the two-fold purpose of a gentle introduction to abstract formalism and a useful tool for problems where the full fledged framework of jet bundles is not needed (e.g., in large parts of our treatment of particle dynamics).
Finally, it is shown how to find what we call extended Lagrangians for systems that do not posses a classical Lagrangian, as it is often the case in plasma physics.
The geometric point of view as it is stressed in this section is essential in the derivation of the variational integrators and the analysis of symmetries and conservation properties.

\subsection{Hamilton's Action Principle}\label{ch:classical_lagrangian_action_principle}

\begin{figure}[tb]
\centering
\includegraphics[width=.6\textwidth]{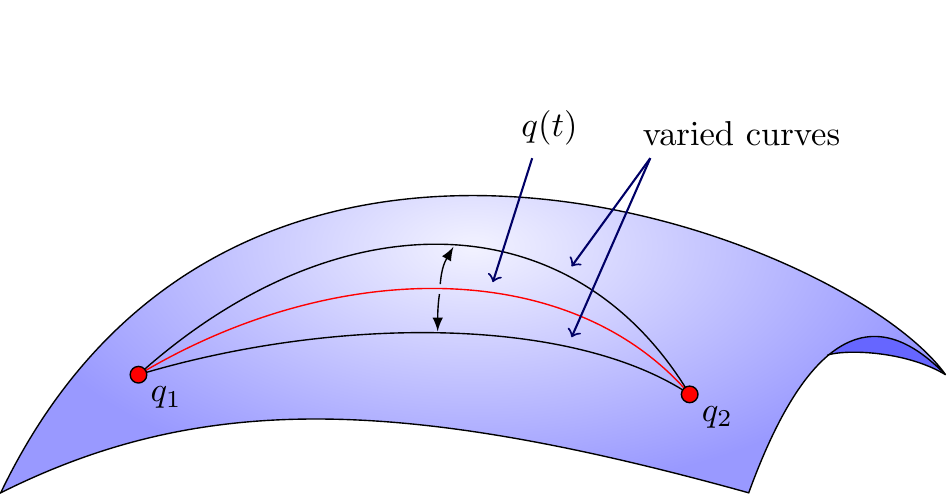}
\caption{Variations of the trajectory $q (t)$.}
\end{figure}

Hamilton's action principle answers the following question: of all possible paths $q(t)$ a system could choose to go from $q(t_{1})$ to $q(t_{2})$, which is the one it actually picks?\footnote{
The following derivation follows along the lines of \citeauthor{SaletanCromer:1971} \cite{SaletanCromer:1971} and \citeauthor{JoseSaletan:1998} \cite{JoseSaletan:1998}.
}

Inserting a given path $q(t)$ into the Lagrangian $L (q (t), \dot{q} (t))$, it becomes a function of time only.
This time dependency is different for all the possible $q(t)$, such that the integral
\begin{align}
\int \limits_{t_{1}}^{t_{2}} L \big( q(t), \dot{q}(t) \big) \, dt
\end{align}

takes different values for different paths $q(t)$. Hamilton's principle states that this integral takes a critical point for the physical path.

To make calculations tractable let us not consider all possible paths from $q(t_{1})$ to $q(t_{2})$ but a family in which each path is determined by a parameter $\eps$.
This family shall contain the actual, physical path for $\eps = 0$.
Each path is a function of time, labelled by $\eps$,
\begin{align}\label{eq:classical_lagrangian_action_principle_family_of_paths}
q_{\eps} (t) \equiv q (t, \eps)
\hspace{2em}
\text{with fixed $\eps$} .
\end{align}

The function $q_{\eps} (t)$ shall be differentiable for both, $t$ at fixed $\eps$ and $\eps$ at fixed $t$, such that mixed partial derivatives can be exchanged
\begin{align}\label{eq:classical_lagrangian_action_principle_partial_derivatives}
\dfrac{\partial^{2} q_{\eps}}{\partial t \, \partial \eps} = \dfrac{\partial^{2} q_{\eps}}{\partial \eps \, \partial t} .
\end{align}

All paths shall start at $q (t_{1})$ and end at $q (t_{2})$, such that
\begin{align}\label{eq:classical_lagrangian_action_principle_boundaries1}
\begin{split}
& q_{\eps} (t_{1}) = q_{0} (t_{1}) = q (t_{1}) \equiv q_{1} \\
& q_{\eps} (t_{2}) = q_{0} (t_{2}) = q (t_{2}) \equiv q_{2}
\end{split}
\end{align}

or
\begin{align}\label{eq:classical_lagrangian_action_principle_boundaries2}
\dfrac{\partial q_{\eps}}{\partial \eps} (t_{1}) =
\dfrac{\partial q_{\eps}}{\partial \eps} (t_{2}) = 0 .
\end{align}

One example of such a family, often considered exclusively in classical mechanics textbooks, is
\begin{align}
q_{\eps} (t) = q(t) + \eps \, \delta q .
\end{align}

This, however, requires that $q(t)$ takes values in a linear space, an assumption that generally cannot be made in the geometric framework on manifolds.
Therefore we consider more general transformations of the form (\ref{eq:classical_lagrangian_action_principle_family_of_paths}).
The action integral is given by
\begin{align}
\mcal{A} [q_{\eps}]
= \int \limits_{t_{1}}^{t_{2}} L \big( q_{\eps} (t), \dot{q}_{\eps} (t) \big) \, dt ,
\end{align}

and has different values for different $\eps$.
Hamilton's principle of stationary action states that for each one-parameter family $q_{\eps}$ that fulfils the above conditions (\ref{eq:classical_lagrangian_action_principle_partial_derivatives} - \ref{eq:classical_lagrangian_action_principle_boundaries2}), $q$ is a critical point of the action iff

\rimpeq{
\dfrac{d}{d\eps} \mcal{A} [q_{\eps}] \bigg\vert_{\eps=0}
= \bigg[ \dfrac{d}{d\eps} \int \limits_{t_{1}}^{t_{2}} L \big( q_{\eps} (t), \dot{q}_{\eps} (t) \big) \, dt \bigg] \bigg\vert_{\eps=0} &= 0 &
\text{(\textbf{Hamilton's Action Principle})} .
}

This means that the time integral of the Lagrangian becomes stationary for the actual, physical motion.
As $t_{1}$ and $t_{2}$ are arbitrary, this is a general statement of Newton's second law.
The integration limits are fixed, so that the derivative can be moved into the integral, such that
\begin{align}
\dfrac{d}{d\eps} \mcal{A} [q_{\eps}] \bigg\vert_{\eps=0}
= \int \limits_{t_{1}}^{t_{2}} \dfrac{dL}{d\eps} \big( q_{\eps} (t), \dot{q}_{\eps} (t) \big) \bigg\vert_{\eps=0} dt
= \int \limits_{t_{1}}^{t_{2}} \bigg[ \dfrac{\partial L}{\partial q} \dfrac{\partial q_{\eps}}{\partial \eps} + \dfrac{\partial L}{\partial \dot{q}} \dfrac{\partial \dot{q}_{\eps}}{\partial \eps} \bigg] \bigg\vert_{\eps=0} dt .
\end{align}

The derivative $\dot{q}_{\eps} \equiv dq_{\eps}/dt$ is the generalised velocity along a particular trajectory that is labelled by a specific value of $\eps$.
This time derivative is taken for fixed $\eps$, so it should better be denoted $\partial q_{\eps} / \partial t$.
However, the important point is that the derivatives with respect to time $t$ and with respect to $\eps$ can be exchanged as in (\ref{eq:classical_lagrangian_action_principle_partial_derivatives}) such that integration by parts can be performed
\begin{align}
\dfrac{d}{d\eps} \mcal{A} [q_{\eps}] \bigg\vert_{\eps=0}
&= \int \limits_{t_{1}}^{t_{2}} \bigg[ \dfrac{\partial L}{\partial q} \dfrac{\partial q_{\eps}}{\partial \eps} + \dfrac{\partial L}{\partial \dot{q}} \dfrac{\partial}{\partial t} \dfrac{\partial q_{\eps}}{\partial \eps} \bigg] \bigg\vert_{\eps=0} dt \\
&= \int \limits_{t_{1}}^{t_{2}} \bigg[ \dfrac{\partial L}{\partial q} - \dfrac{d}{dt} \dfrac{\partial L}{\partial \dot{q}} \bigg] \dfrac{\partial q_{\eps}}{\partial \eps} \bigg\vert_{\eps=0} dt + \int \limits_{t_{1}}^{t_{2}} \dfrac{d}{dt} \bigg[ \dfrac{\partial L}{\partial \dot{q}} \dfrac{\partial q_{\eps}}{\partial \eps} \bigg] \bigg\vert_{\eps=0} dt .
\end{align}

The second integral vanished according to (\ref{eq:classical_lagrangian_action_principle_boundaries2}),
\begin{align}
\bigg[ \dfrac{\partial L}{\partial \dot{q}} \dfrac{\partial q_{\eps}}{\partial \eps} \bigg] \bigg\vert_{t_{1}}^{t_{2}} = 0 ,
\end{align}

such that the first integral has to vanish as well.
At that, it has to vanish for each and every path. And as the $\partial q_{\eps} / \partial \eps$ are arbitrary functions of time (only restricted in that they have to vanish at the endpoints), the expression in square brackets has to vanish\footnote{
For a deeper discussion of this point see \citeauthor{JoseSaletan:1998} \cite{JoseSaletan:1998}, section 3.1, or \citeauthor{GelfandFomin:1963} \cite{GelfandFomin:1963}.
}. Of course, this expression corresponds to the Euler-Lagrange equations

\rimpeq{
\dfrac{\partial L}{\partial q} (q, \dot{q}) - \dfrac{d}{dt} \dfrac{\partial L}{\partial \dot{q}} (q, \dot{q}) &= 0 &
\text{(\textbf{Euler-Lagrange Equations})}.
}

The usual notation is retained by defining
\begin{align}
\delta \equiv \dfrac{d}{d\eps} \bigg\vert_{\eps=0}
\end{align}

and calling $\delta$ an infinitesimal variation. The previous derivation summarises
\begin{align}
\delta \mcal{A}
= \int \limits_{t_{1}}^{t_{2}} \delta L \big( q(t), \dot{q} (t) \big) \, dt
= \int \limits_{t_{1}}^{t_{2}} \bigg[ \dfrac{\partial L}{\partial q} \delta q + \dfrac{\partial L}{\partial \dot{q}} \delta \dot{q} \bigg] \, dt
= \int \limits_{t_{1}}^{t_{2}} \bigg[ \dfrac{\partial L}{\partial q} - \dfrac{d}{dt} \dfrac{\partial L}{\partial \dot{q}} \bigg] \, \delta q \, dt
= 0 .
\end{align}

By the above arguments it is clear that the variation of the time derivative of $q$ equals the time derivative of the variation of $q$, i.e.,
\begin{align}
\delta \dot{q} = \dfrac{d}{dt} \delta q ,
\end{align}

an important point that is often obfuscated by oversimplification.

\subsubsection{Hamilton's Action Principle for Fields}

In the infinite dimensional case (field theory), the Lagrangian can in principle depend on the fields, their derivatives, and also the coordinates.
The latter is however not the case in the envisaged applications, hence for simplicity it is not considered here. The theory on jet bundles includes this case without further complications.

With the restriction to Lagrangians that are only functions of the fields and their first derivatives, the derivation of Hamilton's action principle is not any more complicated than in the finite dimensional case (particle mechanics). Without loss of generality, consider as base space $\mf{X}$ only time plus one spatial dimension $(t,x)$ and a theory of just one scalar field $\phy (t,x)$. The Lagrangian density is thus a function
\begin{align}
\mcal{L} = \mcal{L} \big( \phy (t,x), \phy_{t} (t,x), \phy_{x} (t,x) \big)
\end{align}

and the action is given by
\begin{align}
\mcal{A} [\phy]
= \int \limits_{\mf{X}} \mcal{L} ( \phy, \phy_{t}, \phy_{x} )
= \int \limits_{\mf{X}} L ( \phy, \phy_{t}, \phy_{x} ) \, dt \, dx .
\end{align}

For instructive reasons, all considerations in this section will be taken with respect to the Lagrangian function $L$.
Consider a family of variations $\phy^{\eps}$ of $\phy$ that is defined by
\begin{align}
\phy^{\eps} (t,x) &= \phy (t, x, \eps) &
& \text{with} &
\phy^{0} &= \phy (t,x) . & &&
\end{align}

The variation of the action can be expressed as
\begin{align}\label{eq:classical_fields_action_principle_1}
\dfrac{d}{d \eps} \mcal{A} [\phy^{\eps}] \bigg\vert_{\eps=0} = \dfrac{d}{d\eps} \bigg[ \int \limits_{\mf{X}} L \big( \phy^{\eps} (t,x), \phy^{\eps}_{t} (t,x), \phy^{\eps}_{x} (t,x) \big) \, dt \, dx \bigg] \bigg\vert_{\eps=0}
\end{align}

and Hamilton's principle of stationary action states that $\phy$ is a critical point of the action iff (\ref{eq:classical_fields_action_principle_1}) vanishes for all variations $\phy^{\eps}$ of $\phy$.
The differentiation is carried out under the integral to give
\begin{align}
\dfrac{d}{d \eps} \mcal{A} [\phy^{\eps}] \bigg\vert_{\eps=0}
&= \int \limits_{\mf{X}} \bigg[ \dfrac{dL}{d\eps} \big( \phy^{\eps}, \phy^{\eps}_{t}, \phy^{\eps}_{x} \big) \bigg] \bigg\vert_{\eps=0} \, dt \, dx \\
&= \int \limits_{\mf{X}} \bigg[ \dfrac{\partial L}{\partial \phy} \dfrac{\partial \phy^{\eps}}{\partial \eps} + \dfrac{\partial L}{\partial \phy_{t}} \dfrac{\partial \phy_{t}^{\eps}}{\partial \eps} + \dfrac{\partial L}{\partial \phy_{x}} \dfrac{\partial \phy_{x}^{\eps}}{\partial \eps} \bigg] \bigg\vert_{\eps=0} \, dt \, dx
.
\end{align}

The second and third term are integrated by parts with respect to $t$ and $x$, respectively,
\begin{align}
\dfrac{d}{d \eps} \mcal{A} [\phy^{\eps}] \bigg\vert_{\eps=0}
&= \int \limits_{\mf{X}} \bigg[ \dfrac{\partial L}{\partial \phy} - \dfrac{\partial}{\partial t} \dfrac{\partial L}{\partial \phy_{t}} - \dfrac{\partial}{\partial x} \dfrac{\partial L}{\partial \phy_{x}} \bigg] \dfrac{\partial \phy^{\eps}}{\partial \eps} \bigg\vert_{\eps=0} \, dt \, dx ,
\end{align}

where it is assumed that the fields vanish at infinity such that the boundary terms vanish.
Except for this restriction, the functions $\partial_{\eps} \phy^{\eps}$ are arbitrary, such that the variation of the action vanishes, iff the expression in square brackets vanishes, which is what now leads to the Euler-Lagrange field equations in one spatial dimension,

\rimpeq{\label{eq:classical_fields_euler_lagrange_equations_1d}
\dfrac{\partial L}{\partial \phy} - \dfrac{\partial}{\partial t} \dfrac{\partial L}{\partial \phy_{t}} - \dfrac{\partial}{\partial x} \dfrac{\partial L}{\partial \phy_{x}} &= 0 &
\text{(\textbf{Euler-Lagrange Field Equations})}.
}

As expected, there were no surprises and the derivation was very similar to the one of the finite dimensional case, presented in section \ref{ch:classical_lagrangian_action_principle}.
The usual notation is retained by identifying
\begin{align}
\delta \phy = \dfrac{d}{d\eps} \phy^{\eps} \bigg\vert_{\eps=0}
\end{align}

such that the variation of the action reads
\begin{align}
\delta \mcal{A} [\phy]
&= \int \limits_{\mf{X}} \delta L \big( \phy^{\eps}, \phy^{\eps}_{t}, \phy^{\eps}_{x} \big) \, dt \, dx
\end{align}

and the usual manipulations amount to
\begin{align}
\delta \mcal{A} [\phy]
&= \int \limits_{\mf{X}} \bigg[\dfrac{\partial L}{\partial \phy} \, \delta \phy + \dfrac{\partial L}{\partial \phy_{t}} \, \delta \phy_{t} + \dfrac{\partial L}{\partial \phy_{x}} \delta \phy_{x} \bigg] \, dt \, dx
= \int \limits_{\mf{X}} \bigg[ \dfrac{\partial L}{\partial \phy} - \dfrac{\partial}{\partial t} \dfrac{\partial L}{\partial \phy_{t}} - \dfrac{\partial}{\partial x} \dfrac{\partial L}{\partial \phy_{x}} \bigg] \, \delta \phy \, dt \, dx
\end{align}

with the term in square brackets corresponding to the Euler-Lagrange field equations (\ref{eq:classical_fields_euler_lagrange_equations_1d}).

\subsection{Dynamics on the Tangent Bundle}\label{ch:classical_lagrangian_tangent_bundles}

In introductory textbooks on classical mechanics, the Lagrangian is often defined as a function on the configuration space $\mf{Q}$ with coordinates $q$, which for example might be the three-dimensional euclidean space $\msp{E}^{3}$. In that setting, the velocities $\dot{q}$ and accelerations $\ddot{q}$ correspond to the first and second time derivative of $q$, and if $q$ is a vector in $\msp{E}^{3}$, they are as well. The Euler-Lagrange equations are second order differential equations.

This setting, however, does not seem natural. The Lagrangian is defined with respect to $q$ and its first time derivative $\dot{q}$. So, strictly speaking, $L$ is not a function on $\mf{Q}$ but on a larger space.

Let us take a step back and ask what determines the state of a system. It is not just the position $q$ of all its constituents, but also their respective velocities $\dot{q}$. So the state of a system corresponds to a point in a state space labelled by $(q, \dot{q})$\footnote{
This state space is also called \emph{velocity phasespace} in analogy to the \emph{phasespace} in Hamiltonian dynamics.
}. When the system evolves in time, both $q$ and $\dot{q}$ change, consequently the evolution of both, $q$ and $\dot{q}$, has to be computed, not just the evolution of the coordinates $q$.
It is therefore natural to define the Lagrangian on exactly this space of states. This point of view has many advantages. The obvious one is that the Euler-Lagrange equations become first order differential equations for $q$ and $\dot{q}$.
The consequence of this first order nature of the equations is a separation of the trajectories in state space. There is only one trajectory passing through each point in state space, allowing for the construction of phase portraits. These are visual solutions of the dynamical equations and often useful in the analysis of a dynamical systems' behaviour.

In the next step, these ideas are translated into the geometric language of manifolds. The configuration space is regarded as a smooth manifold, still denoted $\mf{Q}$ and called the \emph{configuration manifold}, with points labelled by $q$. The velocity phasespace corresponds to the tangent bundle $\tb{\mf{Q}}$ of that configuration manifold $\mf{Q}$, called the \emph{velocity phase manifold}, with elements labelled by $(q, v)$\footnote{
At this point, a comment is in order. In the literature, points of $\tb{\mf{Q}}$ are often labelled $(q, \dot{q})$. This notation implies that all curves in $\tb{\mf{Q}}$ are the tangent lift of some curve $q(t)$ in $\mf{Q}$. That is of course not the case! The previous statement is only true for physical trajectories, i.e., solutions $(q, \dot{q})$ of the Euler-Lagrange equations (\ref{eq:classical_lagrangian_euler_lagrange_equations}). But there exist much more curves in $\tb{\mf{Q}}$ for which $v \neq \dot{q}$.
}.
The Lagrangian therefore comes naturally as a map
\begin{align}
L : \tb{\mf{Q}} \rightarrow \rsp .
\end{align}

The advantage of this point of view might not be apparent if one just thinks in term of Euclidean spaces. It will become clearer considering a particle whose motion is constrained to the two-dimensional surface of a sphere $\msp{S}^{2}$. The velocity vector of a particle moving in $\msp{E}^{3}$ is also a vector in $\msp{E}^{3}$.
In $\msp{S}^{2}$, however, the velocity vector of a particle is tangent to the sphere. It is not contained in the sphere, but reaches out into the $\msp{E}^{3}$ in which the sphere is embedded. So to describe the particle motion on $\msp{S}^{2}$, one has to consider an embedding space. It is not possible to do that only by means of $\msp{S}^{2}$ alone. In this example, there might not be much of an issue, but in other dynamical systems the embedding space might not be so easily found and if it can be found might not have any physical meaning.\footnote{
See \citeauthor{JoseSaletan:1998} \cite{JoseSaletan:1998}, section 2.4, for a more detailed discussion. This example is taken from there.
}

Of course, $\tb{\mf{Q}}$ is also a space embedding $\mf{Q}$, but in contrast to $\msp{S}^{2}$ and $\msp{E}^{3}$, there is an intrinsic relation between $\mf{Q}$ and $\tb{\mf{Q}}$, given by the tangent lift as discussed in section \ref{sec:geometry_fibre_bundles}, i.e., the tangent bundle $\tb{\mf{Q}}$ is obtained from $\mf{Q}$ by attaching to each point $q \in \mf{Q}$ the tangent space $\tb[q]{\mf{Q}}$ at that point. The linear space $\tb[q]{\mf{Q}}$ contains all possible velocities at $q$, which are of course tangent to $\mf{Q}$ at that point.

The Lagrangian maps points $(q, v)$ of $\tb{\mf{Q}}$ to the real numbers $\rsp$.
The resulting values are completely independent from the coordinates on $\tb{\mf{Q}}$. That way, the description of the dynamics is intrinsic, independent on any particular choice of coordinate systems.
Restricting the Lagrangian to solutions of the Euler-Lagrange equations (\ref{eq:classical_lagrangian_euler_lagrange_equations}), it becomes a function of $(q, \dot{q})$, as before. In the applications considered later on, it will always be possible to find a global coordinate system for the configuration manifold $\mf{Q}$ and its tangent bundle $\tb{\mf{Q}}$. This allows to circumvent the issues arising in the case, when the path $q(t)$ and its deformations $q_{\eps} (t)$ are not located in a single coordinate patch\footnote{
For a discussion of these issues see \citeauthor{Holm:2009} \cite{Holm:2009}, section 4.1.
}.

\subsubsection{Hamilton's Action Principle on the Tangent Bundle}

In this section, a derivation of Hamilton's action principle in a geometric framework is presented\footnote{
The following derivation follows along the lines of \citeauthor{MarsdenRatiu:2002} \cite{MarsdenRatiu:2002} and \citeauthor{MarsdenWest:2001} \cite{MarsdenWest:2001}. For proofs of some of the statements have a look at \cite{MarsdenRatiu:2002}, section 8.1.
}.
Consider the space of paths $\mf{C} ( \mf{Q} )$ that connect two points in $\mf{Q}$,
\begin{align}\label{eq:classical_lagrangian_tangent_bundle_path_space_1}
\mf{C} ( \mf{Q} ) = \big\{ c : \mf{I} \rightarrow \mf{Q} \; \big\vert \; \mf{I} \subset \rsp \; \text{smooth and bounded} \big\} .
\end{align}

Fixing two points $q_{1}$ and $q_{2}$ in $\mf{Q}$ as well as an interval $[t_{1}, t_{2}]$, the path space from $q_{1}$ to $q_{2}$ is defined as
\begin{align}\label{eq:classical_lagrangian_tangent_bundle_path_space_2}
\mf{C} ( q_{1}, q_{2}, [t_{1}, t_{2}] ) = \big\{ c : [t_{1}, t_{2}] \rightarrow \mf{Q} \; \big\vert \; c(t_{1}) = q_{1}, c(t_{2}) = q_{2} \big\} \subset \mf{C} ( \mf{Q} ) .
\end{align}

Elements $c$ of $\mf{C} ( q_{1}, q_{2}, [t_{1}, t_{2}] )$ are maps that relate points $q$ in configuration space $\mf{Q}$ to points $t$ in the time interval $[t_{1}, t_{2}]$, whereby the first and last points, $c(t_{1})$ and $c(t_{2})$, take fixed values, $q_{1}$ and $q_{2}$, respectively.
Consequently, the action can be written as a map $\mcal{A} : \mf{C} ( q_{1}, q_{2}, [t_{1}, t_{2}] ) \rightarrow \rsp$ assigning real values to each path $c$,
\begin{align}\label{eq:classical_lagrangian_tangent_bundle_action_1}
\mcal{A} [c] = \int \limits_{t_{1}}^{t_{2}} L \big( c (t), \dot{c} (t) \big) \, dt .
\end{align}

If $L$ is the Lagrangian on $\tb{\mf{Q}}$, and $c$ is a path $c : [t_{1}, t_{2}] \rightarrow \mf{Q}$ that connects $q_{1} = c(t_{1})$ with $q_{2} = c (t_{2})$, Hamilton's principle of stationary action states that $c$ obeys the Euler-Lagrange equations (\ref{eq:classical_lagrangian_euler_lagrange_equations}), iff  $c$ is a critical point of the function $\mcal{A} : \mf{C} ( q_{1}, q_{2}, [t_{1}, t_{2}] ) \rightarrow \rsp$, that is $\mcal{A}$ is stationary for $c$ or $\delta \mf{A} [c] = 0$.
Stationarity of $\mcal{A} [c]$ means that $\mcal{A} [c]$ does not change under infinitesimal variations of the path $c$.
Such infinitesimal variations of $c$ live in the tangent space $\tb[c]{\mf{C} ( q_{1}, q_{2}, [t_{1}, t_{2}] )}$ of $\mf{C} ( q_{1}, q_{2}, [t_{1}, t_{2}] )$ at $c$.
The tangent vector to some path $c_{\eps} \in \mf{C} ( q_{1}, q_{2}, [t_{1}, t_{2}] )$ is given by
\begin{align}\label{eq:classical_lagrangian_tangent_bundle_variational_vector_field}
V (t) = \dfrac{d}{d\eps} c_{\eps} (t) \bigg\vert_{\eps=0} .
\end{align}

For each fixed $t$, $c_{\eps}$ is a curve in $\mf{Q}$ through the point $c(t)$, such that $V (t)$ is a tangent vector to $\mf{Q}$ based at $c(t)$, i.e., $V (t) \in \tb[c(t)]{\mf{Q}}$ and thus $\pi_{\mf{Q}} \circ V = c$, where $\pi_{\mf{Q}}$ is the canonical projection $\pi_{\mf{Q}} : \tb{\mf{Q}} \rightarrow \mf{Q}$.
From the restrictions $c_{\eps} (t_{1}) = q_{1}$ and $c_{\eps} (t_{2}) = q_{2}$ follows that $V (t_{1}) = 0$ and $V (t_{2}) = 0$, but otherwise $V$ is an arbitrary function.
To summarise, the \emph{infinitesimal variation} of a path $c : [t_{1}, t_{2}] \rightarrow \mf{Q}$, is the set of maps
\begin{align}\label{eq:classical_lagrangian_tangent_bundle_action_principle_boundaries}
V : [t_{1}, t_{2}] &\rightarrow \tb{\mf{Q}} &
& \text{for which} &
\pi_{\mf{Q}} \circ V &= c &
& \text{and} &
V (t_{1}) = V (t_{2}) &= 0
.
\end{align}

$V$ is called an infinitesimal variation of the path $c$ with fixed endpoints and naturally denoted $V = \delta c$.
With the chain rule one obtains
\begin{align}\label{eq:classical_lagrangian_tangent_bundle_action_principle_1}
\delta \mcal{A} [c]
= \dfrac{d}{d\eps} \mcal{A} [c_{\eps}] \bigg\vert_{\eps=0}
= \dfrac{\partial \mcal{A}}{\partial c} [c] \cdot \dfrac{d c_{\eps}}{d \eps} \bigg\vert_{\eps=0}
= \ext \mcal{A} [c] \cdot V ,
\end{align}

where $\ext \mcal{A} [c]$ and $V = dc_{\eps} / d\eps \vert_{\eps=0}$ are regarded as elements of the cotangent and tangent spaces $\cb[c]{\mf{C}}$ and $\tb[c]{\mf{C}}$ on the manifold $\mf{C}$, respectively.
Therefore the variation of the action can be formulated as
\begin{align}\label{eq:classical_lagrangian_tangent_bundle_action_principle_2}
\ext \mcal{A} [c] \cdot V = \dfrac{d}{d\eps} \bigg[ \int \limits_{t_{1}}^{t_{2}} L \big( c_{\eps} (t), \dot{c}_{\eps} (t) \big) \, dt \bigg] \bigg\vert_{\eps=0} ,
\end{align}

where $L$ is a function of the tangent lift of $c_{\eps}$.
Computation of the derivative under the integral,
\begin{align}\label{eq:classical_lagrangian_tangent_bundle_action_principle_3}
\ext \mcal{A} [c] \cdot V
= \int \limits_{t_{1}}^{t_{2}} \dfrac{d}{d\eps} \bigg[ L \big( c_{\eps} (t), \dot{c}_{\eps} (t) \big) \bigg] \bigg\vert_{\eps=0} \, dt
= \int \limits_{t_{1}}^{t_{2}} \bigg[ \dfrac{\partial L}{\partial q} \cdot \dfrac{d c_{\eps}}{d\eps} \bigg\vert_{\eps=0} + \dfrac{\partial L}{\partial v} \cdot \dfrac{d \dot{c}_{\eps}}{d\eps} \bigg\vert_{\eps=0} \bigg] \, dt ,
\end{align}

leads to
\begin{align}\label{eq:classical_lagrangian_tangent_bundle_action_principle_4}
\ext \mcal{A} [c] \cdot V
= \int \limits_{t_{1}}^{t_{2}} \bigg[ \dfrac{\partial L}{\partial q} \cdot V + \dfrac{\partial L}{\partial v} \cdot \dot{V} \bigg] \, dt
= \int \limits_{t_{1}}^{t_{2}} \ext L \big( c (t), \dot{c} (t) \big) \cdot \tb{V(t)} \, dt ,
\end{align}

with $(V, \dot{V})$ the coordinates of the tangent lift of $V$,
\begin{align}
\tb{V} : (q, V) \mapsto \big( (q, \dot{q}) , ( V, \dot{V} ) \big) .
\end{align}

The second term can be integrated by parts, leading to
\begin{align}\label{eq:classical_lagrangian_tangent_bundle_action_principle_5}
\ext \mcal{A} [c] \cdot V
= \int \limits_{t_{1}}^{t_{2}} \bigg[ \dfrac{\partial L}{\partial q} - \dfrac{d}{dt} \dfrac{\partial L}{\partial v} \bigg] \cdot V \, dt + \bigg[ \dfrac{\partial L}{\partial v} \cdot V \bigg]_{t_{1}}^{t_{2}} ,
\end{align}

where the second term vanishes as $V$ vanishes on both ends of the trajectory (\ref{eq:classical_lagrangian_tangent_bundle_action_principle_boundaries}), such that
\begin{align}\label{eq:classical_lagrangian_tangent_bundle_action_principle_6}
\ext \mcal{A} [c] \cdot V = \int \limits_{t_{1}}^{t_{2}} \bigg[ \dfrac{\partial L}{\partial q} - \dfrac{d}{dt} \dfrac{\partial L}{\partial v} \bigg] \cdot V \, dt = \int \limits_{t_{1}}^{t_{2}} D_{\text{EL}} L \big( \ddot{c} (t) \big) \cdot V(t) \, dt ,
\end{align}

where
\begin{align}
D_{\text{EL}} L (c) : \ddot{\mf{Q}} \rightarrow \cb{\mf{Q}}
\end{align}

is an one-form valued function, defining the Euler-Lagrange operator, and $\ddot{\mf{Q}}$ is a submanifold of $\tb{(\tb{\mf{Q}})}$, such that
\begin{align}
\ddot{\mf{Q}} = \big\{ w \in \tb{(\tb{\mf{Q}})} \; \big\vert \; \tb{\pi_{\mf{Q}}} (w) = \pi_{\tb{\mf{Q}}} (w) \big\} \subset \tb{(\tb{\mf{Q}})} .
\end{align}

In other words, $\ddot{\mf{Q}}$ is the set of second derivatives $\ddot{c} (0)$ of curves $c : \mf{I} \rightarrow \mf{Q}$, which are of the form $((q, \dot{q}), (\dot{q}, \ddot{q})) \in \tb{(\tb{\mf{Q}})}$.
The requirement $\delta \mf{A} [c] = 0$ is equivalent to $\ext \mcal{A} [c] \cdot V = 0$ for all $V \in \tb[c]{\mf{C} ( q_{1}, q_{2}, [t_{1}, t_{2}] )}$ as well as to the Euler-Lagrange equations (\ref{eq:classical_lagrangian_euler_lagrange_equations}), that now are rewritten $D_{\text{EL}} L \big( c(t) \big) = 0$, as $V$ is arbitrary, except for it has to vanish at the end points of the trajectory.
The covariance of the Euler-Lagrange equation, first observed in the original work of Lagrange, is obtained here as a natural consequence of the geometric framework.

\subsubsection{Phasespace Lagrangian}

In plasma physics, there exists another notation that enjoys a certain prevalence, namely that of the \emph{phasespace Lagrangian} \cite{Littlejohn:1983}.
In this formulation, the Lagrangian is not defined as a function on the tangent bundle $\tb{\mf{Q}}$ of the configuration space $\mf{Q}$, but instead the tangent bundle takes the role of the configuration space, such that the Lagrangian is defined on $\tb{(\tb{\mf{Q}})}$.

In practice, the tangent bundle structure of the configuration space $\mf{Z} \cong \tb{\mf{Q}}$ is neglected, and the Lagrangian is defined as a function on $\tb{\mf{Z}} \cong \tb{(\tb{\mf{Q}})}$.
Everything else then follows in a straight forward way.

\subsection{Dynamics on the Jet Bundle}\label{ch:classical_lagrangian_jet_bundles}

Another view is offered by employing jet bundle theory.\footnote{
The derivations of this section follow along the lines of
\citeauthor{GotayMarsden:1998} \cite{GotayMarsden:1998},
\citeauthor{MarsdenPatrick:1998} \cite{MarsdenPatrick:1998, MarsdenPekarsky:2001},
\citeauthor{KouranbaevaShkoller:2000} \cite{KouranbaevaShkoller:2000},
\citeauthor{Kouranbaeva:1999} \cite{Kouranbaeva:1999}
and
\citeauthor{West:2004} \cite{West:2004}.
}
Rewriting the theory in this framework might at first sight seem to unnecessarily complicate things, but its great advantage is that it offers a concise notation that readily generalises to the case of field theories.
That way, jet bundle theory offers a general formulation of the variational problem that applies to finite as well as infinite dimensional systems.
Besides, this formulation is fully covariant, the analysis of symmetries with Noether's theorem is simpler, and it resembles the discrete setting quite nicely.

The idea of a jet is to combine the independent variables (coordinates), the dependent variables (trajectories, fields) and their partial derivatives up to a given order in one single geometric object.
Jets provide a coordinate-free description of differential equations, which is especially useful in the theory of partial differential equations, where they allow us to represent an infinite-dimensional space of maps by sections of a finite dimensional space of jets, thereby avoiding the intricacies of infinite dimensional manifolds. \\

Consider a function $\phy (x)$. It establishes a correspondence between each value $x \in \mf{X}$ and another value $\phy (x) \in \mf{F}$. This second value $\phy (x)$ can be considered as a point in the fibre $\mf{F}$ above $x$. So it seems natural to construct a fibre bundle $\mf{Y}$ over the base manifold $\mf{X}$ with fibres corresponding to $\mf{F}$.
In other words, $\mf{Y}$ is obtained by attaching a fibre $\mf{F}$ to each point $x \in \mf{X}$, such that the fibres of $\mf{Y}$ contain all possible values of functions $\phy (x)$, which can therefore be considered as sections $\phy$ in the bundle $\mf{Y}$,
\begin{align}
\phy : \mf{X} &\rightarrow \mf{Y} &
& \text{with} &
\pi_{\mf{X} \mf{Y}} \circ \phy &= \id_{\mf{X}} , &
&&
\end{align}

where $\pi_{\mf{X} \mf{Y}}$ is the canonical projection
\begin{align}
\pi_{\mf{X} \mf{Y}} : \mf{Y} \rightarrow \mf{X} .
\end{align}

If $(x^{\mu}, y^{a})$ are coordinates on $\mf{Y}$, a section $\phy$ is a map $x \mapsto (x^{\mu} , \phy^{a} (x))$, where we denote by $\phy^{a}$ the vertical components of $\phy$, i.e., the fibre coordinates of $\phy (x)$.
In this setting, the equivalent to the tangent bundle is the first jet bundle $\jb{1}{\mf{Y}}$, which contains the first order partial derivatives of each section $\phy \in \mf{Y}$.
In the same way, the $k$th jet bundle $\jb{k}{\mf{Y}}$ of $\mf{Y}$ is the space that contains the partial derivatives of each section $\phy \in \mf{Y}$ up to order $k$.
However, in what follows only the first order jet bundle $\jb{1}{\mf{Y}}$ is needed, so all considerations are restricted to that case.

Coordinates on $\jb{1}{\mf{Y}}$ are $(x^{\mu}, y^{a}, v^{a}_{\mu})$, where $x^{\mu}$ are the coordinates of the base manifold $\mf{X}$, $y^{a}$ are the values of fields at $x$, and $v^{a}_{\mu}$ are all possible values of the partial derivative of $y^{a}$ with respect to $x^{\mu}$.
$\jb{1}{\mf{Y}}$ has two natural projections.
It can be viewed as a fibre bundle over $\mf{X}$ with the source projection
\begin{align}
\pi_{\mf{X} , \jb{1}{\mf{Y}}} : \jb{1}{\mf{Y}} \rightarrow \mf{X}
\end{align}

as well as a fibre bundle over $\mf{Y}$ with the target projection
\begin{align}
\pi_{\mf{Y} , \jb{1}{\mf{Y}}} : \jb{1}{\mf{Y}} \rightarrow \mf{Y} .
\end{align}

\begin{figure}[H]
\centering
\begin{large}
\begin{tikzpicture}
\matrix (m) [matrix of math nodes,row sep=3em,column sep=4em,minimum width=2em] {
\mf{Y} & \\
& \jb{1}{\mf{Y}} \\
\mf{X} & \\
};
\path[-stealth, line width=.4mm]
(m-2-2) edge node [right] {$\;\pi_{\mf{Y} \jb{1}{\mf{Y}}}$} (m-1-1)
(m-2-2) edge node [right] {$\;\;\pi_{\mf{X} \jb{1}{\mf{Y}}}$} (m-3-1)
(m-1-1) edge node [left ] {$\pi_{\mf{X} \mf{Y}}$} (m-3-1);
\end{tikzpicture}
\end{large}
\end{figure}

The first point of view is especially important.
Consider a section $\phy : \mf{X} \rightarrow \mf{Y}$ of $\pi_{\mf{X} \mf{Y}}$. Its tangent map $\tb[x]{\phy}$ at $x \in \mf{X}$ is represented by the matrix $\partial \phy^{a} (x) / \partial x^{\mu} = \phy^{a}_{\mu} (x)$ and thus can be identified with an element of $\jb[\phy (x)]{1}{\mf{Y}}$.
The map $x \mapsto \tb[x]{\phy}$ is therefore a section of $\pi_{\mf{X} , \jb{1}{\mf{Y}}}$, i.e., a section of $\jb{1}{\mf{Y}}$ regarded as a bundle over $\mf{X}$.
This section, denoted $j^{1} \phy$, is called the \emph{first jet prolongation}\footnote{
The jet prolongation can be seen as producing a coordinate-free Taylor expansion to first (in general $k$th) order, as the jet bundle $\jb{1}{\mf{Y}}$ contains all functions that have the same Taylor series up to the first term.
} (also \emph{canonical prolongation}) of a section $\phy (x)$,
\begin{align}
j^{1} \phy : \mf{X} &\rightarrow \jb{1}{\mf{Y}} &
& \text{in coordinates} &
j^{1} \phy : x &\mapsto \big( x^{\mu}, \phy^{a} (x), \phy^{a}_{\mu} (x) \big) . &
\end{align}

Such sections $j^{1} \phy$ of $\jb{1}{\mf{Y}}$ that correspond to the canonical prolongation\footnote{
Not all sections of $\jb{1}{\mf{Y}}$ are prolongations of a section $\phy \in \mf{Y}$.
} of a section $\phy \in \mf{Y}$ are called holonomic. For them $v^{a}_{\mu}$ can be identified with $\phy^{a}_{\mu} = \partial \phy^{a} / \partial x^{\mu}$. \\

In this setting, a section $j^{1} \phy$ of $\pi_{\mf{X} , \jb{1}{\mf{Y}}}$ generalises the notion of a trajectory and a field.
The \emph{Lagrangian density} $\mcal{L}$ is a $n$-form on the jet bundle $\jb{1}{\mf{Y}}$,
\begin{align}
\mcal{L} : \jb{1}{\mf{Y}} \rightarrow \Omega^{n} (\mf{X}) ,
\end{align}

where $\Omega^{n} (\mf{X})$ denotes the $n$-forms on $\mf{X}$.
The Lagrangian $L$ is a function on the jet bundle $\jb{1}{\mf{Y}}$,
\begin{align}
L : \jb{1}{\mf{Y}} \rightarrow \rsp .
\end{align}

The connection between the two is drawn by the volume form $\omega$ of the base manifold $\mf{X}$,
\begin{align}
\mcal{L} = L \omega .
\end{align}

Here, $n$ is the dimension of the base space $\mf{X}$ and $\omega = dx^{1} \wedge dx^{2} \wedge ... \wedge dx^{n}$, e.g., for $\mf{X}$ corresponding to spacetime we have $\omega = dt \wedge dx \wedge dy \wedge dz$.

\subsubsection{Hamilton's Action Principle on the Jet Bundle}

In the framework of jet bundles, the action is given as the integral of the pullback of the Lagrangian density $\mcal{L}$ with the first jet prolongation $j^{1} \phy$ of a section $\phy : \mf{X} \rightarrow \mf{Y}$
\begin{align}\label{eq:classical_lagrangian_jet_bundle_action_1}
\mcal{A} [\phy] = \int \limits_{\mf{X}} (j^{1} \phy)^{*} \mcal{L} .
\end{align}

As $\mcal{L} = L \, \omega$ and $L$ is a smooth function, the following expressions are equivalent
\begin{align}
(j^{1} \phy)^{*} \mcal{L}
= L \big( j^{1} \phy \big) \, \omega
.
\end{align}

Writing the action (\ref{eq:classical_lagrangian_jet_bundle_action_1}) with respect to the last expression and in coordinates
\begin{align}\label{eq:classical_lagrangian_jet_bundle_action_2}
\mcal{A} [\phy]
= \int \limits_{\mf{X}} L \big( j^{1} \phy \big) \, \omega
= \int \limits_{\mf{X}} L ( x^{\mu}, \phy^{a}, \phy^{a}_{\mu} ) \, \omega
\end{align}

establishes a correspondence between (\ref{eq:classical_lagrangian_jet_bundle_action_1}) and previous formulation (\ref{eq:classical_lagrangian_tangent_bundle_action_1}).
Hamilton's principle states that $\phy$ is a critical point of the action iff
\begin{align}\label{eq:classical_lagrangian_jet_bundle_action_principle_1}
\dfrac{d}{d\eps} \mcal{A} [\phy_{\eps}] \bigg\vert_{\eps=0}
= \dfrac{d}{d\eps} \bigg[ \int \limits_{\mf{X}} (j^{1} \phy_{\eps})^{*} \mcal{L} \bigg] \bigg\vert_{\eps=0}
= 0
\end{align}

for all variations $\phy_{\eps}$ of $\phy$.
These variations are defined as a composition of the trajectory $\phy$ and the vertical transformation $\eta_{\eps}$ of the underlying fibre bundle, namely,
\begin{align}
\phy_{\eps} = \eta_{\eps} \circ \phy = \eta_{\eps} (\phy)
\end{align}

such that the variational vector field $V : \mf{X} \rightarrow \tb{\mf{Y}}$ is defined as
\begin{align}
V
&= \dfrac{d}{d\eps} ( \eta_{\eps} \circ \phy ) \bigg\vert_{\eps=0}
= \dfrac{d \eta_{\eps}}{d\eps} (\phy) \bigg\vert_{\eps = 0} ,
\end{align}

or explicitly
\begin{align}
V : x \mapsto \Big( \big( x, \phy^{a} (x) \big) , \big( 0 , V_{\eta}^{a} \big) \Big)
\end{align}

where $V_{\eta}$ is the generating vector field of the transformation $\eta_{\eps}$ with components
\begin{align}
V_{\eta}^{a} &= \dfrac{d}{d\eps} \eta^{a}_{\eps} \bigg\vert_{\eps=0} ,
\end{align}

$\eta_{\eps}^{a}$ being the $y^{a}$ component of $\eta_{\eps}$.
For the moment, we are considering only vertical transformations as that is sufficient for the derivation of the Euler-Lagrange equations, but the Euler-Lagrange equations are also obtained for general variations, not necessarily of the form $\eta_{\eps} \circ \phy$.

The flow map $\eta_{\eps}$ can be interpreted as dragging the path $\phy$ along $V_{\eta}$ through the configuration space.
From now on we drop the $\eta$ index on the field components of the generating vector field. As we do not consider transformations in the coordinates this is no origin of confusion.
The jet prolongation of $V$ to $\jb{1}{\mf{Y}}$ is given by
\begin{align}
j^{1} V &= \dfrac{d}{d\eps} j^{1} ( \eta_{\eps} \circ \phy ) \bigg\vert_{\eps=0}
\end{align}

or in coordinates
\begin{align}
j^{1} V &: x \mapsto \Big( \big( x^{\mu}, \phy^{a} (x) , \phy^{a}_{\nu} (x) \big) , \big( 0, V^{a} , V^{a}_{\nu} + V^{a}_{b} \, \phy^{b}_{\nu} \big) \Big) .
\end{align}

With this and
\begin{align}
j^{1} \phy_{\eps} = j^{1} ( \eta_{\eps} \circ \phy ) = j^{1} \eta_{\eps} \circ j^{1} \phy
\end{align}
such that
\begin{align}
(j^{1} \phy_{\eps})^{*} \mcal{L}
= (j^{1} \phy)^{*} (j^{1} \eta_{\eps})^{*} \mcal{L}
\end{align}

the action principle (\ref{eq:classical_lagrangian_jet_bundle_action_principle_1}) can be rewritten as
\begin{align}\label{eq:classical_lagrangian_jet_bundle_action_principle_2}
\dfrac{d}{d\eps} \mcal{A} [\phy_{\eps}] \bigg\vert_{\eps=0}
= \dfrac{d}{d\eps} \bigg[ \int \limits_{\mf{X}} (j^{1} \phy)^{*} (j^{1} \eta_{\eps})^{*} \mcal{L} \bigg] \bigg\vert_{\eps=0}
= 0 .
\end{align}

With the dynamical definition of the Lie derivative (\ref{eq:geometry_lie_derivative_dynamical_definition}),
\begin{align}\label{eq:classical_lagrangian_jet_bundle_lie_derivative}
\lie_{j^{1} V} \mcal{L} &= \dfrac{d}{d\eps} \Big[ (j^{1} \eta_{\eps})^{*} \mcal{L} \Big] \bigg\vert_{\eps=0} , &
& \text{and} &
\ext \mcal{A} \cdot V
&= \dfrac{d}{d\eps} \mcal{A} [\phy_{\eps}] \bigg\vert_{\eps=0} , &
\end{align}

(\ref{eq:classical_lagrangian_jet_bundle_action_principle_2}) becomes a beautiful, general, geometric formulation of Hamilton's action principle

\rimpeq{\label{eq:classical_lagrangian_jet_bundle_action_principle_3}
\ext \mcal{A} [\phy] \cdot V
= \int \limits_{\mf{X}} (j^{1} \phy)^{*} (\lie_{j^{1} V} \mcal{L})
&= 0 &
& \text{(\textbf{Hamilton's Action Principle})} .
}

This form of the action principle has several advantages.
First and most importantly, this equation is the very same for particles as well as for fields.
Second, it makes explicit the use of the jet prolongation of the trajectory $\phy$ to $j^{1} \phy$ and the variational vector field $V$ to $j^{1} V$, whereas in the tangent bundle formulation, the tangent lift of $\phy$ and $V$ is not explicit in the notation. \\
Last but not least, it is not too difficult to generalise from variations in the configuration space to variations in the full jet space. Thereby not only considering vertical variations, but general variations that might have both horizontal and vertical components.

Coming back to the derivation of the Euler-Lagrange equations, Cartan's magic formula
\begin{align}
\lie_{j^{1} V} \mcal{L}
= i_{j^{1} V} \ext \mcal{L} + \ext ( i_{j^{1} V} \mcal{L} )
\end{align}

needs to be employed to give
\begin{align}\label{eq:classical_lagrangian_jet_bundle_action_principle_4}
\ext \mcal{A} [\phy] \cdot V
= \int \limits_{\mf{X}} (j^{1} \phy)^{*} ( i_{j^{1} V} \ext \mcal{L} )
+ \int \limits_{\mf{X}} \ext \big( (j^{1} \phy)^{*} ( i_{j^{1} V} \mcal{L} ) \big)
= 0
\end{align}

where in the second integral we used that the pullback and the exterior derivative commute.
The second integral vanishes due to Stokes' theorem and the assumption that the variations of $\phy$ vanish at the boundary $\partial \mf{X}$.
Therefore, what is left is just
\begin{align}\label{eq:classical_lagrangian_jet_bundle_action_principle_5}
\ext \mcal{A} [\phy] \cdot V
= \int \limits_{\mf{X}} (j^{1} \phy)^{*} ( i_{j^{1} V} \ext \mcal{L} )
= 0 .
\end{align}

This expression will be the basis for deriving the actual Euler-Lagrange equations in the jet bundle framework for both, particle mechanics and field theory, below.

\subsubsection{Classical Mechanics on Jet Bundles}

In classical mechanics, the base manifold is just time, $\mf{X} = \rsp$, with coordinates $t$.
$\mf{Y}$ is a fibre bundle over time, with the fibres $\mf{Y}_{t}$ corresponding to the configuration space $\mf{Q}$, elements labelled by $q$ and coordinates $(t, q)$, i.e., time and the generalised coordinates.
The fibres $J^{1}_{(t,q)} \mf{Y}$ of the first jet bundle $\jb{1}{\mf{Y}}$ contain the time derivatives of all sections $c : \mf{X} \rightarrow \mf{Y}$. Its coordinates are $(t, q, v)$.
The jet prolongation is given in coordinates by
\begin{align}
j^{1} c : (t, c) \mapsto ( t, c, \dot{c} ) .
\end{align}

Observe that $\jb{1}{\mf{Y}}$ can be identified with (is isomorphic to) $\rsp \times \tb{\mf{Q}}$, sections of $\mf{Y}$ correspond to trajectories $q (t)$ in $\mf{Q}$, sections of $\jb{1}{\mf{Y}}$ to trajectories $(q (t), v (t))$ in $\tb{\mf{Q}}$, and that the jet prolongation is analogous to the tangent lift, such that holonomic sections $j^{1} c$ of $\jb{1}{\mf{Y}}$ are solutions $( q (t), \dot{q} (t) )$ of the Euler-Lagrange equations.
In classical mechanics the volume form is just $\omega = dt$.
This implies that the Lagrangian is a function
\begin{align}
L : \jb{1}{\mf{Y}} \rightarrow \rsp.
\end{align}

The coordinate expression of $\mcal{L}$ is
\begin{align}
\mcal{L} = L (q, v) \, dt ,
\end{align}

where we are considering a time-independent Lagrangian $L$.
Starting from (\ref{eq:classical_lagrangian_jet_bundle_action_principle_5}), compute the exterior derivative, contract with $j^{1} V$ and do the usual partial integration
\begin{align}\label{eq:classical_lagrangian_jet_bundle_action_principle_6}
\ext \mcal{A} [c] \cdot V
= \int \limits_{\mf{X}} \bigg[ \dfrac{\partial L}{\partial q} (j^{1} c) \cdot V + \dfrac{\partial L}{\partial v} (j^{1} c) \cdot \dot{V} \bigg] \, dt
= \int \limits_{\mf{X}} \bigg[ \dfrac{\partial L}{\partial q} (j^{1} c) - \dfrac{d}{dt} \dfrac{\partial L}{\partial v} (j^{1} c) \bigg] \cdot V \, dt
= 0 .
\end{align}

The usual arguments then again yield the Euler-Lagrange equations
\begin{align}
\dfrac{\partial L}{\partial q} (j^{1} c) - \dfrac{d}{dt} \dfrac{\partial L}{\partial v} (j^{1} c) = 0 .
\end{align}

Note that in the jet bundle framework, the case of an explicit time dependency of the Lagrangian is automatically included.

\subsubsection{Field Theory on Jet Bundles}

In field theory, the base manifold $\mf{X}$ is usually identified with spacetime. Its points are denoted $x$ and its coordinates are $(t, x, y, z)$ abbreviated as $x^{\nu}$.
$\mf{Y}$ is thus a fibre bundle over spacetime with coordinates $(x^{\nu}, y^{a})$, where $y^{a}$ are the different fields or field components of the theory, and the first jet bundle $\jb{1}{\mf{Y}}$ has coordinates $(x^{\nu}, y^{a}, v^{a}_{\nu})$.
Hence, the Lagrangian density is a function $\mcal{L} (x^{\nu}, y^{a}, v^{a}_{\nu})$.

Considering a field theory of a (possibly vector valued) field $\phy (x) : \mf{X} \rightarrow \mf{Y}$, defined over spacetime, one can directly start from the action principle as formulated in equation (\ref{eq:classical_lagrangian_jet_bundle_action_principle_3})
\begin{align}\label{eq:classical_fields_jet_bundles_action_principle_1}
\ext \mcal{A} [\phy] \cdot V
= \int \limits_{\mf{X}} (j^{1} \phy)^{*} (\lie_{j^{1} V} \mcal{L})
= 0 ,
\end{align}

as all considerations that lead to this equation were completely general.
All the hard work of section \ref{ch:classical_lagrangian_jet_bundles} is paying off now.
The coordinate expressions of $\phy$ and $V$ and their jet prolongations $j^{1} \phy$ and $j^{1} V$ are
\begin{subequations}
\begin{align}
\phy &: x \mapsto \big( x^{\mu}, \phy^{a} \big) , &
j^{1} \phy &: x \mapsto \big( x^{\mu}, \phy^{a}, \phy^{a}_{\mu} \big) , \\
V &: x \mapsto \big( ( x^{\mu}, \phy^{a} ) , ( 0, V^{a} ) \big) , &
j^{1} V &: x \mapsto \big( ( x^{\mu}, \phy^{a}, \phy^{a}_{\mu} ) , ( 0, V^{a}, V^{a}_{\mu} ) \big) .
\end{align}
\end{subequations}

Inserting this into (\ref{eq:classical_fields_jet_bundles_action_principle_1}) and making the exterior derivative, the contraction and the pullback explicit gives
\begin{align}\label{eq:classical_fields_jet_bundles_action_principle_2}
\ext \mcal{A} [\phy] \cdot V
&= \int \limits_{\mf{X}} \bigg[ \dfrac{\partial L}{\partial y^{a}} (j^{1} \phy) \, V^{a} + \dfrac{\partial L}{\partial v^{a}_{\nu}} (j^{1} \phy) \, V^{a}_{\nu} \bigg] \, \omega \\
&= \int \limits_{\mf{X}} \bigg[ \dfrac{\partial L}{\partial y^{a}} (j^{1} \phy) - \dfrac{\partial}{\partial x^{\nu}} \dfrac{\partial L}{\partial v^{a}_{\nu}} (j^{1} \phy) \bigg] \, V^{a} \, \omega
\end{align}

which leads to the Euler-Lagrange field equations for a theory of a field $\phy$ on spacetime

\rimpeq{\label{eq:classical_fields_jet_bundles_euler_lagrange_field_equations}
\dfrac{\partial L}{\partial y^{a}} (j^{1} \phy) - \dfrac{\partial}{\partial x^{\nu}} \dfrac{\partial L}{\partial v^{a}_{\nu}} (j^{1} \phy) &= 0 &
& \text{(\textbf{Euler-Lagrange Field Equations})} .
}

\subsection{Variational Route to the Cartan Form}

In this section we want to describe a variational derivation of the Cartan form, one of the two fundamental geometric structures of classical mechanics and classical field theories (the other one being the (multi)symplectic form, covered in the next section).

In most treatments, the Cartan form and the multisymplectic form are constructed by using the Legendre transformation to pull back the canonical forms from the Hamiltonian side (cotangent bundle) to the Lagrangian side (tangent bundle).
However, it has been shown by \citeauthor{MarsdenPatrick:1998} \cite{MarsdenPatrick:1998} that the Cartan form arises naturally in the boundary term of the variation of the action in Hamilton's action principle, thus allowing to obtain these structures while staying on the Lagrangian side, entirely.
The advantage of this approach is the possibility of a geometric treatment of theories for which a Hamiltonian cannot be defined. This is especially important in the light of extended Lagrangians as they will be introduced in section \ref{ch:classical_extended_lagrangians}.

After a short look at the Cartan one-form in the tangent bundle setting, which is restricted to autonomous systems of classical mechanics\footnote{
It is possible to derive the Cartan one-form for non-autonomous systems of classical mechanics in the tangent bundle framework (see \citeauthor{MarsdenWest:2001} \cite{MarsdenWest:2001}, section 4.2), but it is not much clearer than the more general derivation on jet bundles.
}, we will generalise the derivation to jet bundles and obtain an expression of the Cartan form that is valid for autonomous as well as non-autonomous systems of classical mechanics and also field theories.

\subsubsection{Lagrangian One- and Two-Form}

Besides leading to the equations of motion, the variational principle provides a direct and natural way to derive the fundamental geometric structures of classical mechanics\footnote{
The following derivation follows along the lines of \citeauthor{MarsdenPatrick:1998} \cite{MarsdenPatrick:1998}, section 2, and \citeauthor{MarsdenRatiu:2002} \cite{MarsdenRatiu:2002}, section 8.2.
}.
For this derivation, the boundary conditions $\delta q (t_{1}) = \delta q (t_{2}) = 0$ are removed, while the time interval is kept fixed.
Thus the variational principle reads
\begin{align}\label{eq:classical_lagrangian_tangent_bundle_action_principle_7}
\ext \mcal{A} [q(t)] \cdot \delta q(t) = \int \limits_{t_{1}}^{t_{2}} \bigg[ \dfrac{\partial L}{\partial q} - \dfrac{d}{dt} \dfrac{\partial L}{\partial v} \bigg] \cdot \delta q \, dt + \bigg[ \dfrac{\partial L}{\partial v} \cdot \delta q \bigg]_{t_{1}}^{t_{2}}
\end{align}

where the variations $\delta q$ do not vanish at the boundary point, so that the last term on the right hand side does not vanish.
This last term corresponds to a linear pairing of the function $\partial L / \partial v$, which is a function of $(q, \dot{q})$, with the tangent vector $\delta q$.
This term can be regarded as a one-form on $\tb{\mf{Q}}$ \footnote{
One could be tempted to regard $\partial L / \partial \dot{q}$ as a one-form on $\mf{Q}$ as it only has a component in $dq$. The same way $\delta q$ could be regarded as a tangent vector on $\mf{Q}$.
However, $\partial L / \partial \dot{q}$ is a function of $(q, \dot{q})$ and therefore clearly a function on $\tb{\mf{Q}}$.
$\delta q$ can also be replaced with a more general vector $\delta \ohat{q} \in \tb{(\tb{\mf{Q}})}$ that has non-vanishing components $(\delta q, \delta v)$.
}, referred to as the \emph{Lagrangian one-form} or \emph{Cartan one-form},
\begin{align}\label{eq:classical_particles_lagrangian_cartan_one_form}
\Theta_{L} = \dfrac{\partial L}{\partial v} \, dq .
\end{align}

This means that the Lagrangian one-form $\Theta_{L}$ is the boundary term of the functional derivative of the action, if the boundary is varied.
The negative of the exterior derivative of the Lagrangian one-form gives the \emph{Lagrangian two-form}, also referred to as the \emph{symplectic two-form}
\begin{align}\label{eq:classical_particles_lagrangian_symplectic_two_form}
\Omega_{L} \equiv - \ext \Theta_{L} ,
\end{align}

given in coordinates by
\begin{align}
\Omega_{L} 
= \dfrac{\partial^{2} L}{\partial q^{i} \, \partial v^{j}} \, d q^{i} \wedge d q^{j}
+ \dfrac{\partial^{2} L}{\partial v^{i} \, \partial v^{j}} \, d v^{i} \wedge d q^{j}
.
\end{align}

For details on the connection between the Lagrangian one-form $\Theta_{L}$ on $\tb{\mf{Q}}$ and the canonical one-form $\Theta$ on $\cb{\mf{Q}}$ as well as between the Lagrangian two-form $\Omega_{L}$ on $\tb{\mf{Q}}$ and the canonical symplectic two-form $\Omega$ on $\cb{\mf{Q}}$ the reader is referred to \citeauthor{MarsdenRatiu:2002} \cite{MarsdenRatiu:2002}.

\subsubsection{Cartan Form and Multisymplectic Form}

To derive the Cartan form in a general setting that applies to classical mechanics as well as to field theories, the action principle on the jet bundle has to be generalised a bit further\footnote{
The following derivation follows along the lines of
\citeauthor{MarsdenPatrick:1998} \cite{MarsdenPatrick:1998, MarsdenPekarsky:2001},
\citeauthor{KouranbaevaShkoller:2000} \cite{KouranbaevaShkoller:2000},
\citeauthor{Kouranbaeva:1999} \cite{Kouranbaeva:1999}, chapter 4,
and
\citeauthor{West:2004} \cite{West:2004}, chapter 5,.
}.
Before, only vertical variations of the action were considered, thereby implicitly restricting the treatment to Lagrangians, that are not explicitly time dependent in the case of particles, or do not explicitly depend on the coordinates in the case of fields, respectively.
But in order to obtain the correct Cartan form in the general case also horizontal variations need to be considered.
Whereas the Euler-Lagrange equations obtained by considering vertical or arbitrary variations are the same, the Cartan form is missing one term if only vertical variations are accounted for. \\

Allowing also for horizontal variations brings some complications.
A transformation $\eta$ acting on a section $\phy : \mf{U}_{\mf{X}} \rightarrow \mf{Y}$, defined over a bounded domain $\mf{U}_{\mf{X}} \subset \mf{X}$,
\begin{align}\label{eq:classical_fields_cartan_form_1}
\eta : \big( x, \phy(x) \big) \mapsto \big( \eta_{\mf{X}} (x) , \eta_{\mf{Y}} ( x, \phy (x) ) \big)
\end{align}

changes not only the section $\phy$ to $\eta \circ \phy$ but also the base space from $\mf{U}_{\mf{X}}$ to $\eta_{\mf{X}} (\mf{U}_{\mf{X}})$.
We explain now how to get around this issue.

Consider a smooth manifold $\mf{U}$ with smooth closed boundary $\partial \mf{U}$. $\mf{U}$ shall be a parametrisation of the space $\mf{U}_{\mf{X}} \subset \mf{X}$ on which the physical sections are defined. This is similar to the previous case, where a space of curves $\mf{C} (\mf{Q})$ was defined (\ref{eq:classical_lagrangian_tangent_bundle_path_space_1}), such that elements of $\mf{C} (\mf{Q})$ correspond to parametrisations of the physical trajectories.
Thus, in total analogy we define the set of smooth maps
\begin{align}\label{eq:classical_fields_cartan_form_2}
\mf{C} (\mf{Y}) = \big\{ \phi : \mf{U} \rightarrow \mf{Y} \; \big\vert \; \pi_{\mf{X} \mf{Y}} \circ \phi : \mf{U} \rightarrow \mf{X} \; \text{is an embedding} \big\}
\end{align}

in coordinates
\begin{align}\label{eq:classical_fields_cartan_form_3}
\phi : u &\mapsto \big( x^{\mu} (u), \phi^{a} (u) \big) &
& \text{with} &
& \text{$x (u)$ an embedding} , & && 
\end{align}

and $\phi^{a} (u)$ are the fibre coordinates of $\phi (u)$.
Points in $\mf{X}$ and $\mf{U}$ are denoted $x$ and $u$, respectively, and their coordinates are denoted $x^{\mu}$ and $u^{\mu}$, respectively.
For each $\phi \in \mf{C} (\mf{Y})$ define
\begin{align}\label{eq:classical_fields_cartan_form_4}
\phi_{\mf{X}} &\equiv \pi_{\mf{X} \mf{Y}} \circ \phi &
& \text{and} &
\mf{U}_{\mf{X}} &\equiv  \pi_{\mf{X} \mf{Y}} \circ \phi (\mf{U}) &
& \text{such that} &
\phi_{\mf{X}} : \mf{U} \rightarrow \mf{U}_{\mf{X}}
\end{align}

in coordinates
\begin{align}\label{eq:classical_fields_cartan_form_5}
\phi_{\mf{X}} : u \mapsto x^{\mu} (u) .
\end{align}

Since $\phi_{\mf{X}}$ is assumed to be an embedding, $\mf{U}_{\mf{X}}$ is a submanifold of $\mf{X}$ that has a smooth closed boundary, just like $\mf{U}$. It is the physical space on which the fields and trajectories, i.e., the physical sections, are defined. Closed boundaries are necessary as the term we are interested in of the variational principle, the one that yields the Cartan form, is the boundary term that arises from the partial integration.

The map $\phi_{\mf{X}}$ is a diffeomorphism between $\mf{U}$ and $\mf{U}_{\mf{X}}$. It maps between the physical space and its parametrisation, such that the composition
\begin{align}\label{eq:classical_fields_cartan_form_6}
\phy = \phi \circ \phi_{\mf{X}}^{-1}
\end{align}

corresponds to a parametrisation of the sections that are physical fields or trajectories.
These physical sections, defined on $\mf{U}_{\mf{X}}$, can be seen as sections of the fibre bundle represented by $\pi_{\mf{U_{\mf{X}}} , \mf{Y}}$.
These are maps
\begin{align}\label{eq:classical_fields_cartan_form_7}
\phy : \mf{U}_{\mf{X}} &\rightarrow \mf{Y} &
& \text{with} &
\pi_{\mf{U_{\mf{X}}} , \mf{Y}} \circ \phy &= \pi_{\mf{X} \mf{Y}} \circ \phy = \id  &
\end{align}

in coordinates
\begin{align}\label{eq:classical_fields_cartan_form_8}
\phy : x \mapsto ( x^{\mu}, \phy^{a} (x) ) .
\end{align}

The fields that are varied in the action principle are the $\phi$. Horizontal variations of the $\phy$ would change the base space $\mf{U}_{\mf{X}}$ on which the fields are defined.
However, a transformation $\mf{Y} \rightarrow \mf{Y}$ acts naturally on the space $\mf{C}$, mapping $\mf{C}$ into itself, even though $\mf{U}_{\mf{X}}$ is not mapped into itself.

\begin{figure}[H]
\centering
\begin{large}
\begin{tikzpicture}
\matrix (m) [matrix of math nodes,row sep=3em,column sep=4em,minimum width=2em] {
& \mf{Y} & \\
& & \\
\mf{X} & \mf{U}_{\mf{X}} & \mf{U} \\
};

\path[-stealth, line width=.4mm]
([yshift=+3pt]m-3-2.east) edge node [above] {$\phi_{\mf{X}}^{-1}$} ([yshift=+3pt]m-3-3.west)
([yshift=-2pt]m-3-3.west) edge node [below] {$\phi_{\mf{X}}$}      ([yshift=-2pt]m-3-2.east)
([xshift=-2pt]m-1-2.south) edge node [below left ] {$\pi_{\mf{U_{\mf{X}}} \mf{Y}}$} ([xshift=-2pt]m-3-2.north)
([xshift=+2pt]m-3-2.north) edge node [below right] {$\phy$}                         ([xshift=+2pt]m-1-2.south)
(m-1-2) edge node [above left ] {$\pi_{\mf{X} \mf{Y}}$} (m-3-1)
(m-3-3) edge node [above right] {$\phi$}                (m-1-2);
\path[left hook->, line width=.4mm]
(m-3-2) edge node [below] {} (m-3-1);
\end{tikzpicture}
\end{large}
\end{figure}

To clarify some of the previous statements, consider the analogous derivation on tangent bundles from section \ref{ch:classical_lagrangian_tangent_bundles}.
We defined $\mf{C} ( \mf{Q} )$ as the space of trajectories that connect two points in $\mf{Q}$.
These trajectories $c \in \mf{C} ( \mf{Q} )$ were considered as maps
\begin{align}\label{eq:classical_fields_cartan_form_9}
c : \mf{I} &\rightarrow \mf{Q} &
& \text{with} &
\mf{I} &\subset \rsp \; \text{smooth and bounded} . & &&
\end{align}

In the jet bundle framework they correspond to elements of $\mf{C} ( \mf{Y} )$
\begin{align}\label{eq:classical_fields_cartan_form_10}
c : \mf{U} \rightarrow \mf{Y} .
\end{align}

In the general case, the parameter space $\mf{U}$ does not have just one dimension but as many as the base space $\mf{X}$.
So the correspondence of the previous and the current notation is
\begin{align*}
\mf{I} &\leftrightarrow \mf{U} , &
\mf{Q} &\leftrightarrow \mf{Y} , &
\mf{C} ( \mf{Q} ) &\leftrightarrow \mf{C} ( \mf{Y} ) , &
c &\leftrightarrow \phi , &
q &\leftrightarrow \phy . &
\end{align*}

Going back to the general theory, the tangent space to $\mf{C}$ at a point $\phi$ is the set $\tb[\phi]{\mf{C}}$ defined as
\begin{align}\label{eq:classical_fields_cartan_form_11}
\tb[\phi]{\mf{C}} (\mf{Y}) = \big\{ V : \mf{U} \rightarrow \tb{\mf{Y}} \; \big\vert \; \pi_{ \mf{Y} , \tb{\mf{Y}} } \circ V = \phi  \big\} .
\end{align}

The elements $V$ of $\tb[\phi]{\mf{C}}$ are called variations of the sections $\phi$ and have coordinate expressions
\begin{align}\label{eq:classical_fields_cartan_form_12}
V : u \mapsto \Big( \big( x^{\mu} (u), \phy^{a} (u) \big) , \big( V^{\mu} (u), V^{a} (u) \big) \Big) ,
\end{align}

where the $V^{\mu}$ correspond to horizontal variations and the $V^{a}$ to vertical variations.
To each vector field $V$ on $\tb[\phi]{\mf{C}}$ belongs a vector field $V_{\mf{X}}$ on $\mf{X}$, given by projection
\begin{align}\label{eq:classical_fields_cartan_form_13}
V_{\mf{X}} \equiv \tb{ \pi_{\mf{X} \mf{Y}} } \circ V .
\end{align}

The projectors can be explicitly written as
\begin{align}\label{eq:classical_fields_cartan_form_14}
\pi_{ \mf{Y} , \tb{\mf{Y}} } &: \big( ( x^{\mu}, y^{a} ) , (V^{\mu}, V^{a}) \big) \mapsto ( x^{\mu}, y^{a} ) , \\
\tb{ \pi_{\mf{X} \mf{Y}} } & : \big( ( x^{\mu}, y^{a} ) , (V^{\mu}, V^{a}) \big) \mapsto ( x^{\mu}, V^{\mu} ) ,
\end{align}

such that $V_{\mf{X}}$ has the coordinate expression
\begin{align}\label{eq:classical_fields_cartan_form_15}
V_{\mf{X}} : u \mapsto \big( x^{\mu} (u), V^{\mu} (u) \big) .
\end{align}

The definition (\ref{eq:classical_fields_cartan_form_11}) of $\tb[\phi]{\mf{C}}$ can also be seen by considering the variation of a path in $\mf{C} (\mf{Y})$,
\begin{align}\label{eq:classical_fields_cartan_form_16}
\phi^{\eps} : u \mapsto \big( x^{\mu} (\eps, u) , \phi^{a} (\eps, u) \big) .
\end{align}

The derivatives of this expression and its projection to $\mf{X}$ are
\begin{align}\label{eq:classical_fields_cartan_form_17}
\dfrac{d\phi^{\eps}}{d\eps} \bigg\vert_{\eps=0} &= \big( V^{\mu} (u), V^{a} (u) \big) ,
\end{align}

which coincides with $V(u)$ in $\tb[\phi(u)]{\mf{Y}}$, and
\begin{align}
\dfrac{d\phi_{\mf{X}}^{\eps}}{d\eps} \bigg\vert_{\eps=0} &= \dfrac{d}{d\eps} \Big[ \pi_{\mf{X} \mf{Y}} \circ \phi^{\eps} \Big] \bigg\vert_{\eps=0} = \tb{\pi_{\mf{X} \mf{Y}}} V = V_{\mf{X}} .
\end{align}

To aid understanding some of these relations are depicted in the diagram below.

\begin{figure}[H]
\centering
\begin{large}
\begin{tikzpicture}
\matrix (m) [matrix of math nodes,row sep=3em,column sep=4em,minimum width=2em] {
& \tb{\mf{Y}} & \tb{\mf{X}} \\
\mf{U} & \mf{Y} & \mf{X} \\
};

\path[-stealth, line width=.4mm]
(m-1-2) edge node [above] {$\tb{\pi_{\mf{X} \mf{Y}}}$}  (m-1-3)
(m-2-1) edge node [above] {$V$}                         (m-1-2)
(m-2-1) edge node [below] {$\phi$}                      (m-2-2)
(m-2-2) edge node [below] {$\pi_{\mf{X} \mf{Y}}$}       (m-2-3)
(m-1-2) edge node [right] {$\pi_{\mf{Y}, \tb{\mf{Y}}}$} (m-2-2)
(m-1-3) edge node [right] {$\pi_{\mf{X}, \tb{\mf{X}}}$} (m-2-3);
\end{tikzpicture}
\end{large}
\end{figure}

Let us rephrase the action principle from the last section in this notation.
The action functional
\begin{align}\label{eq:classical_fields_cartan_form_18}
\mcal{A} : \mf{C} (\mf{Y}) \rightarrow \rsp
\end{align}

is expressed as
\begin{align}\label{eq:classical_fields_cartan_form_19}
\mcal{A} [\phi]
= \int \limits_{\mf{U}_{\mf{X}}} ( j^{1} \phy )^{*} \mcal{L}
= \int \limits_{\mf{U}_{\mf{X}}} \mcal{L} ( j^{1} \phy )
= \int \limits_{\mf{U}_{\mf{X}}} \mcal{L} \big(  j^{1} ( \phi \circ \phi_{\mf{X}}^{-1}) \big) .
\end{align}

As the action $\mcal{A} [\phi]$ depends on $\phi$ only through $\phy$, for any diffeomorphism $\gamma : \mf{U} \rightarrow \mf{U}$
\begin{align}\label{eq:classical_fields_cartan_form_20}
\mcal{A} [\phi \circ \gamma] = \mf{A} [\phi] .
\end{align}

As a consequence, the Euler-Lagrange equations only determine $\phy$ uniquely, not $\phi$. However, as $\phy$ corresponds to the physical fields or trajectories, that is all we need.
Hamilton's principle states that a section $\phi$ of $\mf{C} (\mf{Y})$ solves the Euler-Lagrange equations, iff the action is critical,
\begin{align}\label{eq:classical_fields_cartan_form_21}
\ext \mcal{A} [\phi] \cdot V = 0 ,
\end{align}

for all variations $V \in \tb[\phi]{\mf{C} (\mf{Y})}$ which are zero on the boundary $\partial \mf{U}$ of $\mf{U}$.
To derive the Cartan form, the last restriction has to be removed. The above expression corresponds to
\begin{align}\label{eq:classical_fields_cartan_form_22}
\ext \mcal{A} [\phi] \cdot V
= \dfrac{d}{d\eps} \mcal{A} [\phi^{\eps}] \bigg\vert_{\eps=0} = 0 .
\end{align}

A variation $\phi^{\eps}$ of a section $\phi$ is induced by a transformation $\eta_{\mf{Y}}^{\eps}$ on the configuration space
\begin{align}\label{eq:classical_fields_cartan_form_23}
\eta_{\mf{Y}}^{\eps} : \mf{Y} &\rightarrow \mf{Y} &
& \text{with} &
\eta_{\mf{Y}}^{0} &= \id & &&
\end{align}

through
\begin{align}\label{eq:classical_fields_cartan_form_24}
\phi^{\eps} &= \eta_{\mf{Y}}^{\eps} \circ \phi .
\end{align}

We impose the condition that $\eta_{\mf{Y}}^{\eps}$ covers a diffeomorphism
\begin{align}
\eta_{\mf{X}}^{\eps} : \mf{X} \rightarrow \mf{X} .
\end{align}

In coordinates
\begin{align}\label{eq:classical_fields_cartan_form_25}
\eta_{\mf{Y}}^{\eps} : \big( x, y \big) \mapsto \big( \eta_{\mf{X}}^{\mu} (x), \eta_{\mf{Y}}^{a} (x, y) \big) . &
\end{align}

The diffeomorphism on the base space $\mf{X}$ is obtained through the projection
\begin{align}\label{eq:classical_fields_cartan_form_27}
\eta_{\mf{X}}^{\eps} = \pi_{\mf{X} \mf{Y}} \circ \eta_{\mf{Y}}^{\eps} .
\end{align}

The following diagram should help clarify these relations.

\begin{figure}[H]
\centering
\begin{large}
\begin{tikzpicture}
\matrix (m) [matrix of math nodes,row sep=3em,column sep=4em,minimum width=2em] {
\mf{Y} & \mf{Y} \\
\mf{U}_{\mf{X}} & \eta_{\mf{X}}^{\eps} (\mf{U}_{\mf{X}}) \\
};

\path[-stealth, line width=.4mm]
(m-1-1) edge node [above] {$\eta_{\mf{Y}}^{\eps}$} (m-1-2)
(m-2-1) edge node [below] {$\eta_{\mf{X}}^{\eps}$} (m-2-2)
(m-2-1) edge node [left]  {$\phy$}                 (m-1-1)
(m-2-2) edge node [right] {$\phy^{\eps}$}          (m-1-2);
\end{tikzpicture}
\end{large}
\end{figure}

We see now why it is necessary to introduce a parameter space $\mf{U}$.
A physical section $\phy = \phi \circ \phi_{\mf{X}}^{-1}$ is a section of $\pi_{\mf{U}_{\mf{X}}, \mf{Y}}$.
But the transformation $\eta^{\eps} \circ \phi$ induces a section $\phy^{\eps} = \eta_{\mf{Y}}^{\eps} \circ ( \phi \circ \phi_{\mf{X}}^{-1} ) \circ (\eta_{\mf{X}}^{\eps})^{-1}$ of $\pi_{\eta_{\mf{X}} (\mf{U}_{\mf{X}}),\mf{Y}}$, i.e., the base space itself changes under the transformation.
This becomes more evident by looking at the coordinate expressions
\begin{subequations}\label{eq:classical_fields_cartan_form_28}
\begin{align}
\phi &: u \mapsto \big( x^{\mu} (u) , \phi^{a} (u) \big) , &
\phi^{\eps} &: u \mapsto \big( x^{\mu} (u) , \eta_{\mf{Y}}^{a} ( x (u), \phi (u) ) \big) , \\
\phy &: x \mapsto \big( x^{\mu} , \phi^{a} (x) \big) , &
\phy^{\eps} &: \otilde{x} \mapsto \Big( \otilde{x}^{\mu} , \eta_{\mf{Y}}^{a} ( x, \phy (x) ) \Big) , &
\end{align}
\end{subequations}

where $\otilde{x} \in \eta_{\mf{X}} (\mf{U}_{\mf{X}})$ and $x = (\eta_{\mf{X}}^{\eps})^{-1} ( \otilde{x} ) \in \mf{U}_{\mf{X}}$.
If we consider not variations of the physical sections $\phy : \mf{U}_{\mf{X}} \rightarrow \mf{Y}$ but variations of the sections $\phi : \mf{U} \rightarrow \mf{Y}$, the point $u$ in the base space $\mf{U}$ stays fixed.

Applying the transformation (\ref{eq:classical_fields_cartan_form_23}) to the action (\ref{eq:classical_fields_cartan_form_19}), the variation (\ref{eq:classical_fields_cartan_form_22}) becomes
\begin{align}\label{eq:classical_fields_cartan_form_29}
\ext \mcal{A} [\phi] \cdot V
&= \dfrac{d}{d\eps} \mcal{A} [\eta_{\mf{Y}}^{\eps} \circ  \phi] \bigg\vert_{\eps=0}
= \dfrac{d}{d\eps} \int \limits_{\eta_{\mf{X}}^{\eps} (\mf{U}_{\mf{X}})} \mcal{L} \big( j^{1} \phy^{\eps} \big) \bigg\vert_{\eps=0}
= \int \limits_{\mf{U}_{\mf{X}}} \dfrac{d}{d\eps} \big( \eta_{\mf{X}}^{\eps} \big)^{*} \mcal{L} \big( j^{1} \phy^{\eps} \big) \bigg\vert_{\eps=0} .
\end{align}

Application of the chain rule yields
\begin{align}\label{eq:classical_fields_cartan_form_30}
\ext \mcal{A} [\phi] \cdot V
&= \int \limits_{\mf{U}_{\mf{X}}} \dfrac{d}{d\eps} \bigg[ \big( \eta_{\mf{X}}^{\eps} \big)^{*} \mcal{L} \big( j^{1} \phy^{0} \big) \bigg] \bigg\vert_{\eps=0}
+ \int \limits_{\mf{U}_{\mf{X}}} \dfrac{d}{d\eps} \bigg[ \big( \eta_{\mf{X}}^{0} \big)^{*} \mcal{L} \big( j^{1} \phy^{\eps} \big) \bigg] \bigg\vert_{\eps=0} .
\end{align}

In the first integral, apply the dynamical definition of the Lie derivative as before in (\ref{eq:classical_lagrangian_jet_bundle_lie_derivative}), and in the second integral realise that $\eta_{\mf{X}}^{0}$ is just the identity
\begin{align}\label{eq:classical_fields_cartan_form_31}
\ext \mcal{A} [\phi] \cdot V
&= \int \limits_{\mf{U}_{\mf{X}}}  \lie_{V_{\mf{X}}} \mcal{L} \big( j^{1} \phy \big)
+ \int \limits_{\mf{U}_{\mf{X}}} \dfrac{d}{d\eps} \bigg[ \mcal{L} \big( j^{1} \phy^{\eps} \big) \bigg] \bigg\vert_{\eps=0} .
\end{align}

Use Cartan's magic formula (\ref{eq:geometry_lie_derivative_cartans_magic_formula}) in the first integral, and rewrite the second integral by making use of the identity
\begin{align}\label{eq:classical_fields_cartan_form_32}
\dfrac{d}{d\eps} \bigg[ \mcal{L} \big( j^{1} \phy^{\eps} \big) \bigg] \bigg\vert_{\eps=0}
= \dfrac{d}{d\eps} \bigg[ L \big( j^{1} \phy^{\eps} \big) \bigg] \bigg\vert_{\eps=0} \, \omega
= \Big( \iprod_{j^{1} V_{\phy}} \ext L \big( j^{1} \phy \big) \Big) \, \omega ,
\end{align}

such that
\begin{align}\label{eq:classical_fields_cartan_form_33}
\ext \mcal{A} [\phi] \cdot V
&= \int \limits_{\mf{U}_{\mf{X}}} \iprod_{V_{\mf{X}}} \ext \mcal{L} \big( j^{1} \phy \big)
+ \int \limits_{\mf{U}_{\mf{X}}} \ext \Big( \iprod_{V_{\mf{X}}} \mcal{L} \big( j^{1} \phy \big) \Big)
+ \int \limits_{\mf{U}_{\mf{X}}} \Big( \iprod_{j^{1} V_{\phy}} \ext L \big( j^{1} \phy \big) \Big) \, \omega .
\end{align}

The first integral vanishes as $\mcal{L} = L \omega$ and therefore
\begin{align}\label{eq:classical_fields_cartan_form_34}
\ext \mcal{L} (j^{1} \phy)
= \ext L(j^{1} \phy) \wedge \omega + L(j^{1} \phy) \, \ext \omega
= L_{\mu}(j^{1} \phy) \, dx^{\mu} \wedge \omega + L(j^{1} \phy) \, \ext \omega
= 0 ,
\end{align}

but $\omega = dx^{1} \wedge ... \wedge dx^{n}$ is a form of maximum order on the base space, such that $\ext \omega = 0$ and $dx^{\mu} \wedge \omega = 0$ for all $\mu$.
By Stokes' theorem, the second integral can be transformed into a surface integral, with the result that
\begin{align}\label{eq:classical_fields_cartan_form_35}
\ext \mcal{A} [\phi] \cdot V
&= \int \limits_{\partial \mf{U}_{\mf{X}}} \iprod_{V_{\mf{X}}} \mcal{L} \big( j^{1} \phy \big)
+ \int \limits_{\mf{U}_{\mf{X}}} \Big( \iprod_{j^{1} V_{\phy}} \ext L \big( j^{1} \phy \big) \Big) \, \omega .
\end{align}

Now we have to compute the vector field $V_{\phy}$ corresponding to the transformation of the physical section $\phy$
\begin{align}\label{eq:classical_fields_cartan_form_36}
V_{\phy}
= \dfrac{d}{d\eps} \phy^{\eps} \bigg\vert_{\eps=0}
= \dfrac{d}{d\eps} \Big[  \eta_{\mf{Y}}^{\eps} \circ ( \phi \circ \phi_{\mf{X}}^{-1} ) \circ (\eta_{\mf{X}}^{\eps})^{-1} \Big] \bigg\vert_{\eps=0} ,
\end{align}

and its jet prolongation $j^{1} V_{\phy}$.
Applying the chain rule, we get
\begin{align}\label{eq:classical_fields_cartan_form_37}
V_{\phy}
&= \dfrac{d}{d\eps} \Big[  \eta_{\mf{Y}}^{\eps} \circ ( \phi \circ \phi_{\mf{X}}^{-1} ) \circ (\eta_{\mf{X}}^{0})^{-1} \Big] \bigg\vert_{\eps=0}
+ \dfrac{d}{d\eps} \Big[  \eta_{\mf{Y}}^{0} \circ ( \phi \circ \phi_{\mf{X}}^{-1} ) \circ (\eta_{\mf{X}}^{\eps})^{-1} \Big] \bigg\vert_{\eps=0} .
\end{align}

In the first term we use (\ref{eq:classical_fields_cartan_form_17}) and in the second term we use the fact that $d/d\eps (\eta_{\mf{X}}^{\eps})^{-1} \vert_{\eps=0} = - V_{\mf{X}}$\footnotemark, such that
\begin{align}\label{eq:classical_fields_cartan_form_38}
V_{\phy}
&= V \circ \phi_{\mf{X}}^{-1} - \tb{(\phi \circ \phi_{\mf{X}}^{-1})} \circ V_{\mf{X}} .
\end{align}

The tangent lift of the vector field $V_{\mf{X}}$ is simply
\begin{align}\label{eq:classical_fields_cartan_form_39}
\tb{(\phi \circ \phi_{\mf{X}}^{-1})} \circ V_{\mf{X}} = \Big( \big( x^{\nu}, \phy^{a} (x) \big) , \big( V^{\nu} , \phy^{a}_{\mu} V^{\mu} \big) \Big)
\end{align}

such that
\begin{align}\label{eq:classical_fields_cartan_form_40}
V_{\phy} &= ( 0 , \delta \phy^{a} ) = \big( 0 , V^{a} - \phy^{a}_{\mu} V^{\mu} \big) .
\end{align}

\footnotetext{
This can be seen by the group property of the transformation
\begin{align*}
\eta_{\mf{X}}^{-\eps} &= ( \eta_{\mf{X}}^{\eps} )^{-1} &
& \rightarrow &
- V_{\mf{X}} = \dfrac{d}{d\eps} \eta_{\mf{X}}^{-\eps} \bigg\vert_{\eps=0}
&= \dfrac{d}{d\eps} ( \eta_{\mf{X}}^{\eps} )^{-1} \bigg\vert_{\eps=0} & &&
\end{align*}

or by a simple calculation as follows
\begin{align*}
\eta_{\mf{X}}^{\eps} \circ ( \eta_{\mf{X}}^{\eps} )^{-1} = \eta_{\mf{X}} \Big( \eps, \eta_{X}^{-1} ( \eps , x ) \Big) = \id
\end{align*}

such that
\begin{align*}
\dfrac{d}{d\eps} \Big[ \eta_{\mf{X}}^{\eps} \circ ( \eta_{\mf{X}}^{\eps} )^{-1} \Big] \bigg\vert_{\eps=0}
= \dfrac{d}{d\eps} \Big[ \eta_{\mf{X}}^{\eps} \circ ( \eta_{\mf{X}}^{0} )^{-1} \Big] \bigg\vert_{\eps=0}
+ \dfrac{d}{d\eps} \Big[ \eta_{\mf{X}}^{0} \circ ( \eta_{\mf{X}}^{\eps} )^{-1} \Big] \bigg\vert_{\eps=0}
= 0.
\end{align*}

$\eta_{\mf{X}}^{0} = \id$ and $( \eta_{\mf{X}}^{0} )^{-1} = \id^{-1} = \id$ as well, such that
\begin{align*}
\dfrac{d}{d\eps} \Big[ \eta_{\mf{X}}^{\eps} \circ ( \eta_{\mf{X}}^{\eps} )^{-1} \Big] \bigg\vert_{\eps=0}
= \dfrac{d}{d\eps} \Big[ \eta_{\mf{X}}^{\eps} \Big] \bigg\vert_{\eps=0}
+ \tb{\id} \circ \dfrac{d}{d\eps} \Big[ ( \eta_{\mf{X}}^{\eps} )^{-1} \Big] \bigg\vert_{\eps=0}
= V_{\mf{X}} + \dfrac{d}{d\eps} ( \eta_{\mf{X}}^{\eps} )^{-1} \bigg\vert_{\eps=0}
= 0 ,
\end{align*}

where the tangent lift of the identity is the identity on the tangent space.
}

This is just the vertical component of the vector field $V$ from (\ref{eq:classical_fields_cartan_form_12})\footnote{
Any vector $V \in \tb[\phi]{\mf{C}}$ can be decomposed into a horizontal and a vertical component $V = V^{h} + V^{v}$, where $V^{h} = \tb{\phy} \circ V_{\mf{X}}$ and $V^{v} = V - V^{h}$.
}.
The jet prolongation of $V_{\phy}$ along $j^{1} \phy$ is
\begin{align}\label{eq:classical_fields_cartan_form_41}
j^{1} V_{\phy}
&= ( 0 , \delta \phy^{a} , \delta \phy^{a}_{\nu} )
= \Big( 0 , V^{a} - \phy^{a}_{\mu} V^{\mu} , \partial_{\nu} (V^{a} - \phy^{a}_{\mu} V^{\mu}) \Big) .
\end{align}

With that we compute the action (\ref{eq:classical_fields_cartan_form_35}) as
\begin{align}\label{eq:classical_fields_cartan_form_42}
\ext \mcal{A} [\phi] \cdot V
&=  \int \limits_{\mf{U}_{\mf{X}}} \bigg[
\dfrac{\partial L}{\partial y^{a}} ( j^{1} \phy ) \, \delta \phy^{a} + \dfrac{\partial L}{\partial v^{a}_{\mu}} ( j^{1} \phy ) \, \delta \phy^{a}_{\mu}
\bigg] \, \omega
+ \int \limits_{\partial \mf{U}_{\mf{X}}} L ( j^{1} \phy ) \, V^{\mu} \, \omega_{\mu}
\end{align}

where $\omega_{\mu} = \partial_{\mu} \contr \omega$.
Integrate by parts the second term of the first integral
\begin{align}\label{eq:classical_fields_cartan_form_43}
\ext \mcal{A} [\phi] \cdot V
\nonumber
&= \int \limits_{\mf{U}_{\mf{X}}}
\bigg[
\dfrac{\partial L}{\partial y^{a}} ( j^{1} \phy )
- \dfrac{\partial}{\partial x^{\mu}} \bigg( \dfrac{\partial L}{\partial v^{a}_{\mu}} ( j^{1} \phy ) \bigg)
\bigg] \, \delta \phy^{a} \, \omega \\
&+ \int \limits_{\partial \mf{U}_{\mf{X}}}
\bigg[
L ( j^{1} \phy ) \, V^{\mu}
+ \dfrac{\partial L}{\partial v^{a}_{\mu}} ( j^{1} \phy ) \, V^{a}
- \dfrac{\partial L}{\partial v^{a}_{\mu}} ( j^{1} \phy ) \, \phy^{a}_{\nu} V^{\nu}
\bigg] \, \omega_{\mu} .
\end{align}

To bring this expression into a coordinate-free form, consider a general vector field $W = ( W^{\nu} , W^{a} , W^{a}_{\nu} )$.
Its contraction with $dy^{a} \wedge \omega$ is given by
\begin{align}\label{eq:classical_fields_cartan_form_44}
\iprod_{W} (dy^{a} \wedge \omega) = W^{a} \, \omega - (-1)^{p^{\nu}} \, W^{\nu} \, dy^{a} \wedge \omega_{\nu}
\end{align}

where $p^{\nu}$ is the number of permutations in the computation of $\omega_{\nu} = \partial_{\nu} \contr \omega$ (the additional minus results from the permutation with $dy^{a}$).
The pullback of this relation with $j^{1} \phy$ is
\begin{align}\label{eq:classical_fields_cartan_form_45}
(j^{1} \phy)^{*} \iprod_{W} (dy^{a} \wedge \omega)
= W^{a} (j^{1} \phy) \, \omega - (-1)^{p^{\nu}} \, W^{\nu} \, \phy^{a}_{\nu} \, dx^{\nu} \wedge \omega_{\nu}
= ( W^{a} (j^{1} \phy) - W^{\nu} \, \phy^{a}_{\nu} ) \, \omega .
\end{align}

Applying this result to the variation of the action (\ref{eq:classical_fields_cartan_form_43}), we find
\begin{align}\label{eq:classical_fields_cartan_form_46}
( j^{1} \phy )^{*} \iprod_{j^{1} V} (dy^{a} \wedge \omega) = \delta \phy^{a} \, \omega .
\end{align}

Further, consider the expression
\begin{align}\label{eq:classical_fields_cartan_form_47}
\iprod_{W} (dy^{a} \wedge \omega_{\nu})
= W^{a} \omega_{\nu} - (-1)^{p^{\mu}} \, W^{\mu} \, dy^{a} \wedge \omega_{\nu \mu} .
\end{align}

With
\begin{align}\label{eq:classical_fields_cartan_form_48}
dx^{\sigma} \wedge \omega_{\nu \mu} =
\begin{cases}
\quad 0 & \sigma \neq \mu, \nu \\
\hphantom{-} \omega_{\nu} & \mu = \sigma \\
-  \omega_{\mu} & \nu = \sigma
\end{cases}
\end{align}

the pullback of (\ref{eq:classical_fields_cartan_form_47}) with $j^{1} \phy$ is
\begin{align}\label{eq:classical_fields_cartan_form_49}
(j^{1} \phy)^{*} \iprod_{W} (dy^{a} \wedge \omega_{\mu})
\nonumber
&= W^{a} \omega_{\mu} - (-1)^{p^{\nu}} \, W^{\nu} \, \phy^{a}_{\sigma} \, dx^{\sigma} \wedge \omega_{\mu \nu} \\
&= W^{a} \omega_{\mu} - W^{\nu} \phy^{a}_{\nu} \omega_{\mu} + \phy^{a}_{\mu} \, W^{\nu} \omega_{\mu} .
\end{align}

Applying this result to the variation of the action (\ref{eq:classical_fields_cartan_form_43}), we find
\begin{align}\label{eq:classical_fields_cartan_form_50}
(j^{1} \phy)^{*} \iprod_{V} ( dy^{a} \wedge \omega_{\mu} )
\nonumber
&= V^{a} \omega_{\mu} - V^{\sigma} \phy^{a}_{\nu} \, dx^{\nu} \wedge \omega_{\mu \sigma}
= V^{a} \omega_{\mu} - V^{\sigma} \phy^{a}_{\nu} \, ( \omega_{\mu} \delta_{\sigma}^{\nu} - \omega_{\sigma} \delta_{\mu}^{\nu} ) \\
&= V^{a} \omega_{\mu}
- V^{\nu} \phy^{a}_{\nu} \, \omega_{\mu}
+ V^{\nu} \phy^{a}_{\mu} \, \omega_{\nu} .
\end{align}

A final but simple computation shows
\begin{align}\label{eq:classical_fields_cartan_form_51}
(j^{1} \phy)^{*} \iprod_{j^{1} V} \omega = V^{\mu} \, \omega_{\mu} .
\end{align}

Therefore, the variation of the action (\ref{eq:classical_fields_cartan_form_43}) can be written as
\begin{align}\label{eq:classical_fields_cartan_form_52}
\ext \mcal{A} [\phi] \cdot V
\nonumber
&= \int \limits_{\mf{U}_{\mf{X}}} ( j^{1} \phy )^{*} \iprod_{j^{1} V}
\bigg[ \dfrac{\partial L}{\partial y^{a}}
- \dfrac{\partial}{\partial x^{\mu}} \dfrac{\partial L}{\partial v^{a}_{\mu}}
\bigg] \, dy^{a} \wedge \omega \\
& \hspace{6em}
+ \int \limits_{\partial \mf{U}_{\mf{X}}} ( j^{1} \phy )^{*} \iprod_{j^{1} V}
\bigg[
\dfrac{\partial L}{\partial v^{a}_{\mu}} \, dy^{a} \wedge \omega_{\mu}
+ \bigg( L - \dfrac{\partial L}{\partial v^{a}_{\mu}} \, \phy^{a}_{\mu} \bigg) \omega
\bigg] .
\end{align}

In the first integral we find the Euler-Lagrange equations as we derived them before by considering only vertical variations (\ref{eq:classical_fields_jet_bundles_euler_lagrange_field_equations}).
The expression in square brackets in the second integral is the looked for Cartan form

\rimpeq{\label{eq:classical_fields_cartan_form_53}
\Theta_{L} &= \dfrac{\partial L}{\partial v^{a}_{\mu}} \, dy^{a} \wedge \omega_{\mu} + \bigg( L - \dfrac{\partial L}{\partial v^{a}_{\mu}} \, \phy^{a}_{\mu} \bigg) \, \omega &
& \hspace{5em} &
& \text{(\textbf{Cartan Form})} . &
}

In the case of particle mechanics, this is a one-form, otherwise it is a form of the order of the base manifold.
The \emph{symplectic form} is defined as the exterior derivative of the Cartan form

\rimpeq{\label{eq:classical_fields_cartan_form_54}
\Omega_{L} = - \ext \Theta_{L} = dy^{a} \wedge \ext \bigg ( \dfrac{\partial L}{\partial v^{a}_{\mu}} \bigg) \wedge \omega_{\mu} - \ext \bigg( L - \dfrac{\partial L}{\partial v^{a}_{\mu}} \, \phy^{a}_{\mu} \bigg) \wedge \omega  &
& \text{(\textbf{Multisymplectic Form})} . &
}

Both, the Cartan and the symplectic form, are defined on the first jet bundle, i.e., $\Theta_{L} \in \Omega^{1}$ and $\Omega_{L} \in \Omega^{2}$.
Using these expressions, the variation of the action (\ref{eq:classical_fields_cartan_form_52}) can be written as
\begin{align}\label{eq:classical_fields_cartan_form_55}
\ext \mcal{A} [\phi] \cdot V
&= \int \limits_{\mf{U}_{\mf{X}}} ( j^{1} \phy )^{*} ( \iprod_{j^{1} V} \Omega_{L} )
+ \int \limits_{\partial \mf{U}_{\mf{X}}} ( j^{1} \phy )^{*} ( \iprod_{j^{1} V} \Theta_{L} ) .
\end{align}

This will be the starting point to prove the preservation of the multisymplectic form along the Lagrangian flow in section (\ref{sec:classical_multisymplectic_form}).
The expressions (\ref{eq:classical_fields_cartan_form_53}) and (\ref{eq:classical_fields_cartan_form_54}) are more general than the ones we derived previously, (\ref{eq:classical_particles_lagrangian_cartan_one_form}) and (\ref{eq:classical_particles_lagrangian_symplectic_two_form}), in that they lift the restriction to a time-independent Lagrangian and describe field theories as well.

\subsection{Preservation of the Symplectic Form}

In this and the next section we want to prove the conservation of the symplectic and multisymplectic forms under Lagrangian flows.
At first we do so on the tangent bundle, thereby restricting ourselves to the case of particle dynamics.
The approach is then generalised to the framework of jet bundles, whereby we obtain a general proof that is valid for both particle and field systems.

\subsubsection{Euler-Lagrange Map and Lagrangian Vector Fields}

This section mostly aims at making the literature more easily accessible. It is not strictly necessary to understand the subsequent treatment. \\

Define the submanifold $\ddot{\mf{Q}}$ of $\tb{(\tb{\mf{Q}})}$ to be
\begin{align}
\ddot{\mf{Q}} \equiv \bigg\{ w \in \tb{(\tb{\mf{Q}})} \; \bigg\vert \; \tb{\pi_{\mf{Q}}} w = \pi_{\tb{\mf{Q}}} w \bigg\} \subset \tb{(\tb{\mf{Q}})} .
\end{align}

This states that $\ddot{\mf{Q}}$ contains those elements of $\tb{(\tb{\mf{Q}})}$ for which the two projections $\pi_{\tb{\mf{Q}}}$ and $\tb{\pi_{\mf{Q}}}$ coincide.

\begin{figure}[H]
\centering
\begin{large}
\begin{tikzpicture}
\matrix (m) [matrix of math nodes, row sep=3em, column sep=4em, minimum width=2em] {
\tb{\mf{Q}} & \tb{(\tb{\mf{Q}})} \\
\mf{Q} & \\
};
\path[-stealth, line width=.4mm]
(m-1-1) edge node [left ] {$\pi_{\mf{Q}}$} (m-2-1)
([yshift=+2pt]m-1-2.west) edge node [above] {$\pi_{\tb{\mf{Q}}}$} ([yshift=+2pt]m-1-1.east)
([yshift=-3pt]m-1-2.west) edge node [below] {$\tb{\pi_{\mf{Q}}}$} ([yshift=-3pt]m-1-1.east);
\end{tikzpicture}
\end{large}
\end{figure}

To see what that means, write both expressions in coordinates
\begin{subequations}
\begin{align}
\pi_{\tb{\mf{Q}}} : \big( (q, v), (\dot{q}, \dot{v}) \big) &\mapsto (q, v) , \\
\tb{\pi_{\mf{Q}}} : \big( (q, v), (\dot{q}, \dot{v}) \big) &\mapsto (q, \dot{q}) .
\end{align}
\end{subequations}

Requiring that both projections are equivalent therefore means singling out those elements of $\tb{(\tb{\mf{Q}})}$ for which $v = \dot{q}$ and therefore also $\dot{v} = \ddot{q}$.
These correspond to curves $q(t) \in \mf{Q}$ which are tangent lifted twice, first to $\tb{\mf{Q}}$, then to $\tb{(\tb{\mf{Q}})}$.
In other words, elements $w \in \ddot{\mf{Q}}$ are those elements of $\tb{(\tb{\mf{Q}})}$ that have the coordinate expression
\begin{align}
w = \big( (q, \dot{q}), (\dot{q}, \ddot{q}) \big) .
\end{align}

An alternative definition of the second order submanifold $\ddot{\mf{Q}}$ is therefore
\begin{align}
\ddot{\mf{Q}} \equiv \bigg\{ \dfrac{d^{2} q}{dt^{2}} (0) \in \tb{(\tb{\mf{Q}})} \; \bigg\vert \; \text{$q (t)$ a curve in $\mf{Q}$} \bigg\} \subset \tb{(\tb{\mf{Q}})} .
\end{align}

Given a Lagrangian $L$, there exists a map on $\ddot{\mf{Q}}$
\begin{align}
D_{\text{EL}} L : \ddot{\mf{Q}} \rightarrow \cb{\mf{Q}}
\end{align}

referred to as the the Euler-Lagrange map.
It defines a one-form in the dual space $\cb{\mf{Q}}$ of $\tb{\mf{Q}}$ with coordinate expression
\begin{align}
D_{\text{EL}} L = \dfrac{\partial L}{\partial q} - \dfrac{d}{dt} \dfrac{\partial L}{\partial v} .
\end{align}

It is a function on $\ddot{\mf{Q}} \subset \tb{(\tb{\mf{Q}})}$ as
\begin{align}
D_{\text{EL}} L
= \dfrac{\partial L}{\partial q} (q, \dot{q})
- \dfrac{\partial^{2} L}{\partial v \, \partial q} (q, \dot{q}) \cdot \dot{q}
- \dfrac{\partial^{2} L}{\partial v \, \partial v} (q, \dot{q}) \cdot \ddot{q} .
\end{align}

With that, the variational principle (\ref{eq:classical_lagrangian_tangent_bundle_action_principle_7}) can be written
\begin{align}\label{eq:classical_lagrangian_tangent_bundle_action_principle_8}
\ext \mcal{A} [q(t)] \cdot \delta q(t) = \int \limits_{t_{1}}^{t_{2}} D_{\text{EL}} L (\ddot{q}) \cdot \delta q \, dt + \bigg[ \Theta_{L} (\dot{q}) \cdot \delta \ohat{q} \bigg]_{t_{1}}^{t_{2}}
\end{align}

where $\ddot{q}$ refers to an element of $\ddot{\mf{Q}}$ and thus has coordinates $\big( (q, \dot{q}), (\dot{q}, \ddot{q}) \big)$. Similarly, $\dot{q}$ refers to an element of $\tb{\mf{Q}}$ with coordinates $(q, \dot{q})$.
The variation $\delta \ohat{q}$ is defined as
\begin{align}
\delta \ohat{q} \equiv \dfrac{d}{d\eps} \bigg\vert_{\eps = 0} \dfrac{d}{dt} \bigg\vert_{t=0} q_{\eps} (t)
\end{align}

or in coordinates
\begin{align}
\delta \ohat{q} (t) = \Big( \big( q (t) , \dot{q} (t) \big) , \big( \delta q (t) , \delta \dot{q} (t) \big) \Big) .
\end{align}

This Euler-Lagrange map can now be used to define the \emph{Lagrangian vector field}
\begin{align}
X_{L} : \tb{\mf{Q}} \rightarrow \ddot{\mf{Q}}
\end{align}

as a second order vector field on $\tb{\mf{Q}}$ satisfying
\begin{align}
D_{\text{EL}} L \circ X_{L} = 0 .
\end{align}

The flow of $X_{L}$ is called the \emph{Lagrangian flow}
\begin{align}
F_{L}^{t} : \tb{\mf{Q}} \rightarrow \tb{\mf{Q}} .
\end{align}

By construction $q \in \mf{C} ( \mf{Q} )$ is a solution of the Euler-Lagrange equations iff $(q, \dot{q})$ is an integral curve of $X_{L}$.
In the next section, the Lagrangian flow $F_{L}$ will be defined without referring to the Euler-Lagrange map but by using coordinate expressions instead.
The advantage is a somewhat easier treatment.

\subsubsection{Lagrangian Flows and Preservation of the Symplectic Form}

Denote the vector field on $\tb{\mf{Q}}$ that solves the Euler-Lagrange equations by $X_{L}$.
Its flow, referred to as the Lagrangian flow, is a map
\begin{align}\label{eq:classical_lagrangian_symplectic_form_flows_1}
F_{L}^{s} : \tb{\mf{Q}} \rightarrow \tb{\mf{Q}} ,
\end{align}

taking initial values $(q_{0}, \dot{q}_{0})$ to points of the corresponding phasespace trajectory at time $s$, that is
\begin{align}\label{eq:classical_lagrangian_symplectic_form_flows_2}
F_{L}^{s} : \big( q_{0}, \dot{q}_{0} \big) \mapsto \big( q (s), \dot{q} (s) \big) ,
\end{align}

such that
\begin{align}\label{eq:classical_lagrangian_symplectic_form_flows_3}
\dfrac{\partial L}{\partial q} \big( (q (s), \dot{q} (s) \big) - \dfrac{d}{ds} \dfrac{\partial L}{\partial \dot{q}} \big( (q (s), \dot{q} (s) \big) = 0 .
\end{align}

The Lagrangian vector field is accordingly defined as
\begin{align}\label{eq:classical_lagrangian_symplectic_form_flows_4}
X_{L} (q_{0}, v_{0}) = \dfrac{d}{ds} F_{L}^{s} (q_{0}, v_{0}) \bigg\vert_{s=0} .
\end{align}

Restrict the action $\mcal{A}$ to the subspace $\mf{C}_{L} \subset \mf{C} (\mf{Q})$ of solutions of the Euler-Lagrange equations.
Elements $q \in \mf{C}_{L}$ are integral curves of $X_{L}$, and therefore uniquely determined by the initial condition $ v_{q} = \big( q (0), \dot{q} (0) \big) \in \tb{\mf{Q}}$.
Consequently, $\mf{C}_{L}$ may be identified with the space of initial conditions, i.e., $\mf{C}_{L}$ is isomorphic to $\tb{\mf{Q}}$.

Associate to $v_{q}$ the integral curve $s \mapsto F_{L}^{s} (v_{q})$ with $s \in [0, t]$
\begin{align}\label{eq:classical_lagrangian_symplectic_form_flows_5}
q(s) &= \pi_{\mf{Q}} \big( F_{L}^{s} (v_{q}) \big) &
& \text{with} &
F_{L}^{s} (v_{q}) &= \big( q (s), \dot{q} (s) \big) . & &&
\end{align}

The restricted action $\mcal{A}_{t}$ corresponds to the value of $\mcal{A}$ on that curve.
It defines a map
\begin{align}\label{eq:classical_lagrangian_symplectic_form_flows_6}
\mcal{A}_{t} : \tb{\mf{Q}} \rightarrow \rsp
\end{align}

by
\begin{align}\label{eq:classical_lagrangian_symplectic_form_flows_7}
\mcal{A}_{t} [v_{q}] &= \mcal{A} [q] &
& \text{with} &
q &\in \mf{C}_{L} &
& \text{and} &
\big( q (0) , \dot{q} (0) \big) &= v_{q} . & &&
\end{align}

or explicitly
\begin{align}\label{eq:classical_lagrangian_symplectic_form_flows_8}
\mcal{A}_{t} [v_{q}]
= \int \limits_{0}^{t} L \big( q (s), \dot{q} (s) \big) \, ds
= \int \limits_{0}^{t} L \big( F_{L}^{s} (v_{q}) \big) \, ds
\end{align}

Calculating the variation of the restricted action, the first term in (\ref{eq:classical_lagrangian_tangent_bundle_action_principle_7}) vanishes, as $\mcal{A}$ is restricted to solutions of the Euler-Lagrange equations
\begin{align}\label{eq:classical_lagrangian_symplectic_form_flows_9}
\ext \mcal{A} [q] \cdot \delta q
&= \bigg[ \dfrac{\partial L}{\partial \dot{q}} (q, \dot{q}) \cdot \delta \ohat{q} \bigg]_{0}^{t} &
& \text{with} &
q &\in \mf{C}_{L} . & &&
\end{align}

As $\mcal{A}_{t}$ is considered a real-valued function on $\tb{\mf{Q}}$, this becomes
\begin{align}\label{eq:classical_lagrangian_symplectic_form_flows_10}
\ext \mcal{A}_{t} [v_{q}] \cdot w_{v_{q}}
&= \Theta_{L} \big( F_{L}^{t} (v_{q}) \big) \cdot \dfrac{d}{d\eps} F_{L}^{t} (v_{q}^{\eps}) \bigg\vert_{\eps=0} - \Theta_{L} (v_{q}) \cdot w_{v_{q}}
\end{align}

with $v_{q}^{\eps}$ an arbitrary curve in $\tb{\mf{Q}}$, namely
\begin{align}\label{eq:classical_lagrangian_symplectic_form_flows_11}
v_{q}^{\eps} : \rsp &\rightarrow \tb{\mf{Q}} &
& \text{such that} &
v_{q}^{0} &= v_{q} &
& \text{and} &
w_{v_{q}} &= \dfrac{d}{d\eps} v_{q}^{\eps} \bigg\vert_{\eps=0} . &
\end{align}

Since $w_{v_{q}}$ is arbitrary, (\ref{eq:classical_lagrangian_symplectic_form_flows_10}) is equivalent to
\begin{align}\label{eq:classical_lagrangian_symplectic_form_flows_12}
\ext \mcal{A}_{t} = (F_{L}^{t})^{*} \Theta_{L} - \Theta_{L} .
\end{align}

Taking the exterior derivative
\begin{align}\label{eq:classical_lagrangian_symplectic_form_flows_13}
0
= \ext^{2} \mcal{A}_{t}
= \ext (F_{L}^{t})^{*} \Theta_{L} - \Theta_{L}
= - (F_{L}^{t})^{*} \Omega_{L} + \Omega_{L}
.
\end{align}

leads to the conservation of the symplectic form $\Omega_{L}$ along the Lagrangian flow $X_{L}$
\begin{align}\label{eq:classical_lagrangian_symplectic_form_flows_14}
(F_{L}^{t})^{*} \Omega_{L} = \Omega_{L} .
\end{align}

These results can be considered the Lagrangian equivalent of the well-known conservation of the symplectic form by Hamiltonian flows \cite{MarsdenRatiu:2002, Arnold:1989, GotayMarsden:1998}.

\subsection{Preservation of the Multisymplectic Form}\label{sec:classical_multisymplectic_form}

In this section we want to show that the multisymplectic form $\Omega_{L}$ from (\ref{eq:classical_fields_cartan_form_54}) is preserved under the Lagrangian flow, a generalisation of the results from the previous section.
Therefore we recall equation (\ref{eq:classical_fields_cartan_form_55}) for the variation of the action, that is
\begin{align}\label{eq:classical_multisymplectic_form_1}
\ext \mcal{A} [\phi] \cdot V
&= \int \limits_{\mf{U}_{\mf{X}}} ( j^{1} \phy )^{*} ( \iprod_{j^{1} V} \Omega_{L} )
+ \int \limits_{\partial \mf{U}_{\mf{X}}} ( j^{1} \phy )^{*} ( \iprod_{j^{1} V} \Theta_{L} ) .
\end{align}

The action $\mcal{A}$ takes an extremum for $\phi \in \mf{C} (\mf{Y})$ if the first integral vanishes.
The corresponding integrand vanishes not only for vector fields $j^{1} V$, corresponding to vertical transformations, but for general vector fields $W$ on $\jb{1}{\mf{Y}}$, that can be tangent to any $j^{1} \phy$.
As a consequence, $\phi$ is an extremum of the action, if the variation of the action (\ref{eq:classical_multisymplectic_form_1}) vanishes for all vectors $W \in \tb{(\jb{1}{\mf{Y}})}$. Such $\phi$ are solutions of the Euler-Lagrange equations (\ref{eq:classical_fields_jet_bundles_euler_lagrange_field_equations}).

We define $\mf{C}_{L}$ to be the restriction of $\mf{C} (\mf{Y})$, defined in (\ref{eq:classical_fields_cartan_form_2}), to solutions of the Euler-Lagrange equations, i.e.,
\begin{align}\label{eq:classical_multisymplectic_form_2}
\mf{C}_{L} = \big\{ \phi \in \mf{C} (\mf{Y}) \; \big\vert \; (j^{1} \phy)^{*} [ W \contr \Omega_{L} ] = 0 \quad \text{for all} \;\; W \in \tb{(\jb{1}{\mf{Y}})} \big\} ,
\end{align}

such that $\phy^{a}$ is an element of $\mf{C}_{L}$ if
\begin{align}\label{eq:classical_multisymplectic_form_3}
\dfrac{\partial L}{\partial y^{a}} (j^{1} \phy) - \dfrac{\partial}{\partial x^{\mu}} \left( \dfrac{\partial L}{\partial v^{a}_{\mu}} (j^{1}\phy) \right) = 0
\quad \text{in} \;\; \mf{U}_{\mf{X}} .
\end{align}

A vector field $V \in \tb{\mf{C}_{L}}$ is called a first variation.
Its flow maps solutions $\phi$ of the Euler-Lagrange equations to other solutions of the Euler-Lagrange equations, such that sections $\phi \in \mf{C}_{L}$ are integral curves of $V$.

If we restrict the variation of the action to $\mf{C}_{L}$, the first integral in (\ref{eq:classical_multisymplectic_form_1}) becomes zero.
Computing the exterior derivative of (\ref{eq:classical_multisymplectic_form_1}) and restricting it to two first variations $V, W \in \tb{\mf{C}_{L}}$, gives
\begin{align}\label{eq:classical_multisymplectic_form_4}
0
= \ext^{2} \mcal{A} [\phi] \cdot V \cdot W
= \int \limits_{\partial \mf{U}_{\mf{X}}} ( j^{1} \phy )^{*} ( V \contr W \contr \ext \Theta_{L} ) .
\end{align}

This states that the multisymplectic form $\Omega_{L}$ is conserved

\rimpeq{\label{eq:classical_multisymplectic_form_formula}
\int \limits_{\partial \mf{U}_{\mf{X}}} ( j^{1} \phy )^{*} ( V \contr W \contr \Omega_{L} ) &= 0 &
& \text{(\textbf{Multisymplectic Form Formula})} . &
}

A detailed proof of this expression is omitted but can be found in \citeauthor{MarsdenPatrick:1998} \cite{MarsdenPatrick:1998}.

\subsection{Extended Lagrangians}\label{ch:classical_extended_lagrangians}

In the variational treatment of field theoretic problems from plasma physics one faces the problem that most systems do not have a natural Lagrangian formulation.
Similarly, even so most systems are Hamiltonian, they do not feature a canonical Hamiltonian formulation with respect to canonical conjugate variables, but only a so called noncanonical formulation. 
Therefore it is not possible to write a canonical Lagrangian for these systems.

However, to apply the variational integrator formalism, a Lagrangian is indispensable.
Salvation is brought by Ibragimov and his theory of integrating factors and adjoint equations \cite{Ibragimov:2006,Ibragimov:2007}. The basic idea is to extend the system by doubling the number of dependent variables. This enables to write down a Lagrangian that is the product of the original equations and the added auxiliary variables, such that the variation with respect to those new variables results in the original equations.

In the following, we summarise the important definitions and statements of \citeauthor{Ibragimov:2006} \cite{Ibragimov:2006}.

\subsubsection{Integrating Factors}

Integrating factors provide means to solve differential equations.
Any first order differential equation of the form
\begin{align}\label{eq:classical_ibragimov_integrating_factors_deq1}
a (x, y) \, \dfrac{dy}{dx} + b (x, y) = 0
\end{align}

can also be written in differential form, i.e.,
\begin{align}\label{eq:classical_ibragimov_integrating_factors_deq2}
a (x, y) \, dy + b (x, y) \, dx = 0 .
\end{align}

This equations is said to be \emph{exact} if its left hand side is the differential of some function $\alpha (x, y)$
\begin{align}\label{eq:classical_ibragimov_integrating_factors_deq3}
a (x, y) \, dy + b (x, y) \, dx = d \alpha (x, y) .
\end{align}

In general, (\ref{eq:classical_ibragimov_integrating_factors_deq2}) is not exact, but it can become exact upon multiplying by an appropriate function $\mu (x, y)$
\begin{align}\label{eq:classical_ibragimov_integrating_factors_deq4}
\mu (x, y) \, \big( a (x, y) \, dy + b (x, y) \, dx \big) = d \alpha (x, y) .
\end{align}

This function $\mu (x, y)$ is called an \emph{integrating factor} for (\ref{eq:classical_ibragimov_integrating_factors_deq2}).
As
\begin{align}\label{eq:classical_ibragimov_integrating_factors_differential}
\dfrac{\partial \alpha}{\partial x} &= \mu b &
& \text{and} &
\dfrac{\partial \alpha}{\partial y} &= \mu a & &&
\end{align}

the integrability condition for (\ref{eq:classical_ibragimov_integrating_factors_differential}), $\alpha_{xy} = \alpha_{yx}$, yields an equation for determining the integrating factor
\begin{align}
\dfrac{\partial (\mu a)}{\partial x} = \dfrac{\partial (\mu b)}{\partial y} .
\end{align}

In general, for an ordinary differential equation of order $s$,
\begin{align}\label{eq:classical_ibragimov_integrating_factors_deq_ho1}
a (x, y, y', y, ..., y^{(s-1)}) \, y^{(s)} + b (x, y, y' , ..., y^{(s-1)}) = 0 ,
\end{align}

a differential function $\mu (x, y, y', y, ..., y^{(s-1)})$ is an integrating factor, if the multiplication by $\mu$ converts the left hand side of (\ref{eq:classical_ibragimov_integrating_factors_deq_ho1}) into a total derivative of some function $\alpha (x, y, y', y, ..., y^{(s-1)})$,
\begin{align}\label{eq:classical_ibragimov_integrating_factors_deq_ho2}
\mu a y^{(s)} + \mu b = D_{x} ( \alpha ) .
\end{align}

Together, equations (\ref{eq:classical_ibragimov_integrating_factors_deq_ho1}) and (\ref{eq:classical_ibragimov_integrating_factors_deq_ho2}) imply $D_{x} (\alpha) = 0$, such that
\begin{align}
\alpha (x, y, y', y, ..., y^{(s-1)}) = const ,
\end{align}

thereby reducing the order of the differential equation to solve.
The integrating factor for (\ref{eq:classical_ibragimov_integrating_factors_deq_ho1}) is determined by
\begin{align}
\dfrac{\delta}{\delta y} \big( \mu a y^{(s)} + \mu b \big) = 0
\end{align}

where $\delta / \delta y$ is the variational derivative.

\subsubsection{Adjoint Equations}

Consider a first order linear partial differential equation for a scalar field $u(x)$
\begin{align}\label{eq:classical_ibragimov_adjoint_equations_deq}
\mcal{D} [u] = a^{\mu} (x) \, u_{\mu} + b(x) \, u = f(x) .
\end{align}

The first order linear differential operator $\mcal{D}$, corresponding to this equation is
\begin{align}\label{eq:classical_ibragimov_adjoint_equations_op}
\mcal{D} [u] = a^{\mu} (x) \, \partial_{\mu} + b(x) .
\end{align}

The adjoint operator to $\mcal{D}$ is a first-order linear differential operator $\mcal{D}^{*}$ such that
\begin{align}
v \, \mcal{D} [u] - u \, \mcal{D}^{*} [v] = \div \, C (x)
\end{align}

for all functions $u$ and $v$ and some vector field $C(x)$.
The adjoint operator is uniquely determined
\begin{align}\label{eq:classical_ibragimov_adjoint_equations_adj_op}
\mcal{D}^{*} [v] = - \partial_{\mu} ( a^{\mu} v) + b v .
\end{align}

It defines the adjoint equation to (\ref{eq:classical_ibragimov_adjoint_equations_deq}) by
\begin{align}\label{eq:classical_ibragimov_adjoint_equations_adj_deq}
\mcal{D}^{*} [v] = - \partial_{\mu} (a^{\mu} v) + b v = 0 .
\end{align}

If $\mcal{D} [u] = \mcal{D}^{*} [u]$ for any function $u$, the operator $\mcal{D}$ is called self-adjoint. \\

All of these statements and definitions generalise straight forwardly to systems of $m$ partial differential equations of arbitrary order $s$,
\begin{align}\label{eq:classical_ibragimov_adjoint_equations_pdes}
F_{a} (x, u, ..., u_{(s)}) &= 0, &
a &= 1, ..., m . & &&
\end{align}

$F_{a} (x, u, ..., u_{s})$ are differential functions with $n$ independent variables $x^{\mu}$ and $m$ dependent variables $u^{a}$.
The system of adjoint equations to (\ref{eq:classical_ibragimov_adjoint_equations_pdes}) is defined by
\begin{align}\label{eq:classical_ibragimov_adjoint_equations_adj_pdes}
F_{a}^{*} ( x, u, v, ..., u_{(s)}, v_{(s)} ) &= \dfrac{\delta (v^{b} F_{b})}{\delta u^{a}} = 0, &
a &= 1, ..., m , &
\end{align}

where $v = v^{a} (x)$ are $m$ new dependent variables, referred to as \emph{auxiliary variables}.
If the system obtained by substituting $v=u$ in the adjoint equations (\ref{eq:classical_ibragimov_adjoint_equations_adj_pdes}),
\begin{align}
F_{a}^{*} ( x, u, u, ..., u_{(s)}, u_{(s)} ) = 0 ,
\end{align}

is identical with the original system (\ref{eq:classical_ibragimov_adjoint_equations_pdes}), the system is called \emph{self-adjoint}.

\subsubsection{Extended Lagrangians}

The extended system of differential equations, composed of the system of $m$ partial differential equations of order $s$,
\begin{align}\label{eq:classical_ibragimov_lagrangians_pdes}
F_{a} (x, u, ..., u_{(s)}) &= 0 ,
\end{align}

together with its adjoint equations,
\begin{align}\label{eq:classical_ibragimov_lagrangians_adj_pdes}
F_{a}^{*} ( x, u, v, ..., u_{(s)}, v_{(s)} ) &\equiv \dfrac{\delta (v^{b} F_{b})}{\delta u^{a}} = 0 ,
\end{align}

has a Lagrangian given by
\begin{align}\label{eq:classical_ibragimov_lagrangians_lagrangian}
L = v^{a} \, F_{a} .
\end{align}

Obviously, the variation with respect to the auxiliary variables $v^{a}$, in this context also referred to as \emph{Ibragimov multipliers}, yields the original equations (\ref{eq:classical_ibragimov_lagrangians_pdes}),
\begin{align}
\dfrac{\delta L}{\delta v^{a}} = F_{a} (x, u, ..., u_{(s)}) .
\end{align}

Similarly, the variation with respect to the original variables $u^{a}$ yields the adjoint equations (\ref{eq:classical_ibragimov_lagrangians_adj_pdes}),
\begin{align}
\dfrac{\delta L}{\delta u^{a}} = F_{a}^{*} ( x, u, v, ..., u_{(s)}, v_{(s)} ) .
\end{align}

The definition of the adjoint equations (\ref{eq:classical_ibragimov_adjoint_equations_adj_pdes}) suggested this result already.

\subsubsection{Symmetries and Conserved Quantities}

Here, we do not want to anticipate results from the next section, but a short comment seems appropriate.
The original application Ibragimov had in mind for this method was the analysis of symmetries and conservation laws by applying Noether's theorem to systems without classical Lagrangian.
However, an analysis of (\ref{eq:classical_ibragimov_lagrangians_lagrangian}) will obviously lead to symmetries of the extended system of equations that are not necessarily symmetries of the original system. Therefore it is required to apply some appropriate restriction.

If the operator at hand is self-adjoint, the identification of the auxiliary variables $v$ with the original variables $u$ identifies the conserved flux of the extended system with the conserved flux of the original system at once. Most often, however, this is not the case. Ibragimov therefore defines the concept of quasi-self-adjointness \cite{Ibragimov:2007b}, but for us a simpler idea suffices.

The symmetries of the extended system can be reduced to symmetries of the original system, if the auxiliary variables can be expressed with respect to the original variables and their derivatives, i.e.,
\begin{align}\label{eq:classical_ibragimov_symmetries1}
v^{a} = v^{a} ( x, u , ..., u_{(s)} ) .
\end{align}

The adjoint equation (\ref{eq:classical_ibragimov_adjoint_equations_adj_pdes}) thereby becomes
\begin{align}\label{eq:classical_ibragimov_symmetries2}
F_{a}^{*} \big( x, u, v ( x, u , ..., u_{(s)} ), ..., u_{(s)}, v_{(s)} ( x, u , ..., u_{(s)} ) \big) = 0
\end{align}

and thus a function of $x$ and $u$ only.
If it is possible to select the relation (\ref{eq:classical_ibragimov_symmetries1}) such that (\ref{eq:classical_ibragimov_symmetries2}) is automatically respected when $u$ solves (\ref{eq:classical_ibragimov_lagrangians_pdes}), then a conservation law for the extended system amounts to a physical conservation law.

\subsubsection{Multisymplectic Form}

The proof of preservation of the multisymplectic form along the Lagrangian flow from section \ref{sec:classical_multisymplectic_form} can be applied directly to extended Lagrangians.
An open question is if the multisymplectic structure of the extended system can be restricted to the physical system, similar to the restriction of conservation laws of the extended system to conservation laws of the physical system.

It is however not clear if this is in general possible as the physical system is not necessarily Lagrangian and therefore might not even have a compatible multisymplectic structure.
But in some cases, such as Hamiltonian systems like the Vlasov-Poisson system or ideal magnetohydrodynamics, the original system certainly has a multisymplectic structure, such that the development of a restriction method, that establishes a relation between the multisymplectic forms of the extended and the original system, appears worthwhile.

\section{Noether Theorem}\label{sec:classical_noether_theorem}

The Noether theorem \cite{Noether:1918, KosmannSchwarzbach:2010, Neuenschwander:2010} is one of the deepest and most influential insights of mathematical physics.
It states that each continuous symmetry of a Lagrangian corresponds to a conservation law of the associated Euler-Lagrange equations and vice versa.

\subsection{Point Transformations and One Parameter Groups}

Consider an infinitesimal transformation on $\mf{Q}$ which is of the form
\begin{align}\label{eq:noether_point_transformations_1}
q (t) \rightarrow q^{\eps} (t) &= \xi \big( q(t), \eps \big) = \xi^{\eps} \big( q (t) \big) &
& \text{with} &
\xi^{0} &= \id &
& \text{such that} &
q^{0} (t) &= q (t) . &
\end{align}

If the functional dependency of $q^{\eps}$ on the two parameters $t$ and $\eps$ is of importance, we also write $q^{\eps} (t) = q (t, \eps)$.
The transformation $\xi$ maps each point $q \in \mf{Q}$ to a new point $q^{\eps} \in \mf{Q}$ and is therefore called a \emph{point transformation}.
It maps each trajectory $q (t) \in \mf{C} (\mf{Q})$ to a new trajectory $q^{\eps} (t) \in \mf{C} (\mf{Q})$. But what we are looking for are transformations of the Lagrangian and therefore transformations on $\tb{\mf{Q}}$.
Thus we need the tangent lift of $\xi^{\eps}$ which is defined as
\begin{align}\label{eq:noether_point_transformations_2}
\tb{\xi^{\eps}} : (q^{\eps}) &\mapsto ( q^{\eps} , \dot{q}^{\eps}) &
& \text{with} &
\dot{q}^{\eps} (t) &= \dfrac{\partial }{\partial t} q (t, \eps) . &
\end{align}

It maps each point $\dot{q} \in \tb[q]{\mf{Q}}$ to a new point $\dot{q}^{\eps} \in \tb[q^{\eps}]{\mf{Q}}$.
Hence a transformation $\xi^{\eps}$ on $\mf{Q}$ induces a transformation $\tb{\xi^{\eps}}$ on $\tb{\mf{Q}}$.

We shall always assume that $\xi^{\eps}$ corresponds to a continuous family of transformations, such that the trajectories $q^{\eps} (t)$ are continuous in both, $t$ and $\eps$, and the transformation $\xi^{\eps} (q)$ constitutes a one-parameter group of curves.

Instead of specifying the transformation (\ref{eq:noether_point_transformations_1}) directly, it can also be defined by its \emph{generating vector field}
\begin{align}\label{eq:noether_point_transformations_3}
X (q, t) &\equiv \dfrac{\partial q^{\eps} (t)}{\partial \eps} \bigg\vert_{\eps = 0} &
& \text{with} &
X &= X^{i} \, \dfrac{\partial}{\partial q^{i}} , &
\end{align}

sometimes referred to as the infinitesimal symmetry direction.

\subsection{Noether Theorem for Particle Systems}

A Lagrangian has a symmetry if it is invariant under a point transformation $q^{\eps} (t) = \xi (q, \eps)$, that is
\begin{align}\label{eq:noether_particles_1}
L \big( q^{\eps} (t), \dot{q}^{\eps} (t) \big) &= L \big( q (t), \dot{q} (t) \big) &
& \text{for all $\eps$} . &
\end{align}

This is equivalent to
\begin{align}\label{eq:noether_particles_2}
\dfrac{d}{d \eps} \bigg\vert_{\eps = 0} L \big( q^{\eps} (t), \dot{q}^{\eps} (t) \big) &= 0
\end{align}

or explicitly
\begin{align}\label{eq:noether_particles_3}
\dfrac{d}{d \eps} \bigg\vert_{\eps = 0} L \big( q^{\eps}, \dot{q}^{\eps} \big)
= \dfrac{\partial L}{\partial q} \big( q, \dot{q} \big) \cdot X
+ \dfrac{\partial L}{\partial v} \big( q, \dot{q} \big) \cdot \dot{X}
= 0 .
\end{align}

If $q(t)$ solves the Euler-Lagrange equations (\ref{eq:classical_lagrangian_euler_lagrange_equations}),
the first term on the right-hand side of (\ref{eq:noether_particles_3}) can be rewritten, such that the symmetry condition becomes
\begin{align}\label{eq:noether_particles_5}
\bigg[ \dfrac{d}{dt} \dfrac{\partial L}{\partial v} \big( q, \dot{q} \big) \bigg] \cdot X
+ \dfrac{\partial L}{\partial v} \big( q, \dot{q} \big) \cdot \dot{X}
= 0 .
\end{align}

This is a total time derivative and in fact a conservation law
\begin{align}\label{eq:noether_particles_6}
\dfrac{d}{dt} \bigg[ \dfrac{\partial L}{\partial v} \big( q, \dot{q} \big) \cdot X (q, t) \bigg] &= 0 .
\end{align}

It states that solutions $q$ of the Euler-Lagrange equations preserve $\partial L / \partial v$ in direction $X$.

\begin{example}[Example: Point Particle]

Consider a simple transformation that amounts to a time-independent translation
\begin{align}
q^{\eps} (t) &= q (t) + \eps X , &
\dot{q}^{\eps} (t) &= \dot{q} (t) . & &&
\end{align}

The corresponding transformed Lagrangian is
\begin{align}
L \big( q^{\eps} (t), \dot{q}^{\eps} (t) \big)
&= \dfrac{m}{2} \, \big( \dot{q} (t) \big)^{2}
= L ( q, \dot{q} )
\end{align}

which is obviously the same as the untransformed Lagrangian.
The symmetry condition is therefore trivially fulfilled
\begin{align}
\dfrac{\partial}{\partial \eps} \bigg\vert_{\eps = 0} L \big( q^{\eps} (t), \dot{q}^{\eps} (t) \big) = 0
\end{align}

and the corresponding conservation law
\begin{align}
& \dfrac{d}{dt} \bigg[ \dfrac{\partial L}{\partial \dot{q}} \big( q, \dot{q} \big) \cdot X \bigg]
= \dfrac{d}{dt} \bigg[ m \dot{q} \cdot X \bigg]
= 0
\end{align}

states that momentum is preserved in direction of $X$.

\end{example}

\subsection{Noether Theorem for Field Theories}

Now we want to generalise the Noether theorem for finite dimensional systems to infinite dimensional systems.
Still we restrict to vertical transformations, i.e., transformations on the configuration space alone, as that will be sufficient for the following treatment.

The transformation of a (possibly vector valued) field $y$ is described by a one-parameter group of transformations
\begin{align}\label{eq:noether_fields_1}
\eta^{\eps} (x, y) &= ( x, y^{\eps} ) &
& \text{such that} &
\eta^{0} (x, y) &= ( x, y ) . &
\end{align}

or in coordinates
\begin{align}\label{eq:noether_fields_3}
\eta : (x, y) \mapsto ( x^{\nu}, \eta^{a} ( x, y, \eps ) \big) .
\end{align}

The infinitesimal generator of the transformation $\eta$ is
\begin{align}\label{eq:noether_fields_4}
X (x, y) &= \dfrac{d \eta}{d \eps} \bigg\vert_{\eps = 0} = X^{a} \, \dfrac{\partial}{\partial y^{a}} &
& \text{with components} &
X^{a  } (x, y) &= \dfrac{d \eta^{a  }}{d\eps} \bigg\vert_{\eps = 0} ,
\end{align}

and its first jet prolongation is given by
\begin{align}\label{eq:noether_fields_5}
j^{1} X (x, y, v)
&= X^{a} \, \dfrac{\partial}{\partial y^{a}}
+ \dfrac{\partial X^{a}}{\partial y^{b}} \dfrac{\partial y^{b}}{\partial x^{\nu}} \dfrac{\partial}{\partial v^{a}_{\nu}} .
\end{align}

The Lagrangian has a symmetry if it is invariant under this transformation
\begin{align}\label{eq:noether_fields_6}
L \big( x^{\nu}, \, \eta^{a} (x, y), \, \eta^{a}_{\nu} (x, y, v) \big) &= L \big( x^{\nu}, y^{a}, v^{a}_{\nu} \big) &
& \text{for all $\eps$} . &
\end{align}

This is equivalent to
\begin{align}\label{eq:noether_fields_7}
j^{1} X (L)
&= \dfrac{d}{d \eps} \bigg\vert_{\eps = 0} L \big( x^{\nu}, \, \eta^{a} (x, y), \, \eta^{a}_{\nu} (x, y, v) \big) = 0
\end{align}

or explicitly
\begin{align}\label{eq:noether_fields_8}
j^{1} X (L)
&= \dfrac{\partial L}{\partial y^{a  }} \big( x, y, v \big) \cdot X^{a  }
+ \dfrac{\partial L}{\partial v^{a}_{\nu}} \big( x, y, v \big) \cdot X^{a}_{\nu}
= 0.
\end{align}

If $y(x)$ solves the Euler-Lagrange field equations (\ref{eq:classical_fields_jet_bundles_euler_lagrange_field_equations}),
we can replace the first term on right-hand-side and obtain 
\begin{align}\label{eq:noether_fields_10}
j^{1} X (L)
&= \bigg[ \dfrac{\partial}{\partial x^{\nu}} \left( \dfrac{\partial L}{\partial v^{a}_{\nu}} \right) \bigg] \cdot X^{a  }
+ \dfrac{\partial L}{\partial v^{a}_{\nu}} \cdot \bigg[ \dfrac{\partial}{\partial x^{\nu}} X^{a} \bigg]
= 0 .
\end{align}

We immediately see that this is a divergence
\begin{align}\label{eq:noether_fields_11}
j^{1} X (L)
&= \dfrac{\partial}{\partial x^{\nu}} \bigg[ \dfrac{\partial L}{\partial v^{a}_{\nu}} \cdot X^{a} \bigg]
= 0 .
\end{align}

The term in square brackets is called \emph{Noether field}.
Integration in the spatial dimensions yields a conservation law for solutions of the discrete Euler-Lagrange field equations,
\begin{align}\label{eq:noether_fields_12}
\dfrac{d}{dt} \int \dfrac{\partial L}{\partial v^{a}_{t}} \cdot X^{a} \, dx = 0 ,
\end{align}

assuming that the boundary terms vanish.

\subsection{Noether Theorem for Extended Lagrangians}

Consider the generator of a transformation of the configuration bundle for a field $u$ as in (\ref{eq:noether_fields_4}),
\begin{align}\label{eq:noether_extended_1}
X (x, u) &= X^{u} \, \dfrac{\partial}{\partial u} .
\end{align}

To be applied to the Lagrangian of an extended system of equations (\ref{eq:classical_ibragimov_lagrangians_lagrangian}), this generator has to be extended to the auxiliary variable $v$ by a certain function $X^{v}$ as
\begin{align}\label{eq:noether_extended_2}
Y (x, u, v) &= X^{u} \, \dfrac{\partial}{\partial u} + X^{v} \, \dfrac{\partial}{\partial v} .
\end{align}

Everything else follows by applying the theory of the last section to the extended Lagrangian, considered as describing a theory of two fields $(u,v)$.
Only in the last step of integrating the Noether field, one has to find a relation between the auxiliary field $v$ and the original field $u$, i.e., one has to find a functional expression for $v$ in terms of $u$ and its derivatives like in (\ref{eq:classical_ibragimov_symmetries1}).

\chapter{Variational Integrators}\label{ch:variational}

The seminal work in the development of a discrete equivalent of classical mechanics was presented by \citeauthor{Veselov:1988} \cite{Veselov:1988, Veselov:1991}.
His method, based on a discrete variational principle, leads to symplectic integration schemes that automatically preserve constants of motion.
A comprehensive review of discrete mechanics can be found in \citeauthor{MarsdenWest:2001} \cite{MarsdenWest:2001}. This includes also a more thorough account on the historical development.
The theory was extended to partial differential equations in form of first order Lagrangian theories by \citeauthor{MarsdenPatrick:1998} \cite{MarsdenPatrick:1998} and soon also to second order Lagrangian theories by \citeauthor{KouranbaevaShkoller:2000} \cite{KouranbaevaShkoller:2000, Kouranbaeva:1999}.
Another extension was that of asynchronous variational integrators by \citeauthor{Lew:2003} \cite{Lew:2003}, where each point in the spatial grid has its own timestep. This way, exact local energy conservation can be achieved in addition to an often observed speedup in runtime.
In more recent developments, the variational integrator method was applied to Maxwell's equations by \citeauthor{Stern:2009a} \cite{Stern:2009a, Stern:2009b} and fluid problems by \citeauthor{Pavlov:2011} \cite{Pavlov:2011, Pavlov:2009} and \citeauthor{Gawlik:2011} \cite{Gawlik:2011}.

This chapter gives an overview of the theory of variational integrators for finite-dimensional as well as for infinite-dimensional systems.
It follows mostly along the lines of \citeauthor{MarsdenWest:2001} \cite{MarsdenWest:2001}, \citeauthor{MarsdenPatrick:1998} \cite{MarsdenPatrick:1998}, \citeauthor{KouranbaevaShkoller:2000} \cite{KouranbaevaShkoller:2000}, and \citeauthor{Kouranbaeva:1999} \cite{Kouranbaeva:1999}.

\section{Discrete Particle Dynamics}\label{sec:vi_finite}

The derivation of the discrete theory follows along the lines of the derivation of the continuous theory.
The starting point is the discretisation of the action integral and the Lagrangian.
There is some degree of freedom in the choice of the discrete quadrature rule as well as in the approximation of the generalised coordinates $q$ and the generalised velocities $\dot{q}$.
Everything else follows in a straight forward way, so that these choices are determining the respective form of the discrete Lagrangian as well as the resulting discrete equations of motion.

Time will be discretised uniformly, i.e., the timestep $h$ is constant, $q_{k}$ denotes the generalised coordinates at timepoint $k$, and $\dot{q}_{k}$ the generalised velocities at timepoint $k$.
The discrete Lagrangian approximates the time integral of the continuous Lagrangian between two consecutive points in time, $k$ and $k+1$
\begin{align}\label{eq:vi_finite_quadrature}
L_{d} (q_{k}, q_{k+1}) \approx \int \limits_{t_{k}}^{t_{k+1}} L (q, \dot{q}) \, dt .
\end{align}

Its exact expression is determined by the quadrature rule used to approximate the integral.
Here, we assume that the quadrature rule depends only on $(q_{k}, q_{k+1})$.
The discrete action thus becomes merely a sum over the time index of discrete Lagrangians
\begin{align}\label{eq:vi_finite_action}
\mcal{A}_{d} = \sum \limits_{k=0}^{N-1} L_{d} (q_{k}, q_{k+1}) .
\end{align}

The generalised velocities will usually be discretised by simple finite-difference expressions\footnote{
In the first term of the trapezoidal rule (\ref{eq:vi_finite_trapezoidal}), this corresponds to a forward finite-difference, in the second term to a backward finite-difference, and in the midpoint rule (\ref{eq:vi_finite_midpoint}) to a centred finite-difference.
}, i.e.
\begin{align}
\dot{q} &\approx \dfrac{q_{k+1} - q_{k}}{h} &
& \text{for} &
t &\in [ t_{k} , t_{k+1} ] . &
\end{align}

The quadrature (\ref{eq:vi_finite_quadrature}) is most often realised by either the trapezoidal rule
\begin{align}\label{eq:vi_finite_trapezoidal}
L_{d}^{\text{tr}} (q_{k}, q_{k+1}) = \dfrac{h}{2} \, L \bigg( q_{k}, \dfrac{q_{k+1} - q_{k}}{h} \bigg) + \dfrac{h}{2} \, L \bigg( q_{k+1}, \dfrac{q_{k+1} - q_{k}}{h} \bigg)
\end{align}

or the midpoint rule
\begin{align}\label{eq:vi_finite_midpoint}
L_{d}^{\text{mp}} (q_{k}, q_{k+1}) = h \, L \bigg( \dfrac{q_{k} + q_{k+1}}{2}, \dfrac{q_{k+1} - q_{k}}{h} \bigg) .
\end{align}

The configuration manifold of the discrete theory is still $\mf{Q}$, but the discrete state space is $\mf{Q} \times \mf{Q}$ instead of $\tb{\mf{Q}}$,
such that the discrete Lagrangian $L_{d}$ is a function
\begin{align}\label{eq:vi_finite_discrete_lagrangian}
L_{d} : \mf{Q} \times \mf{Q} \rightarrow \rsp .
\end{align}

\subsection{Discrete Action Principle}

\begin{figure}[tb]
\centering
\includegraphics[width=.6\textwidth]{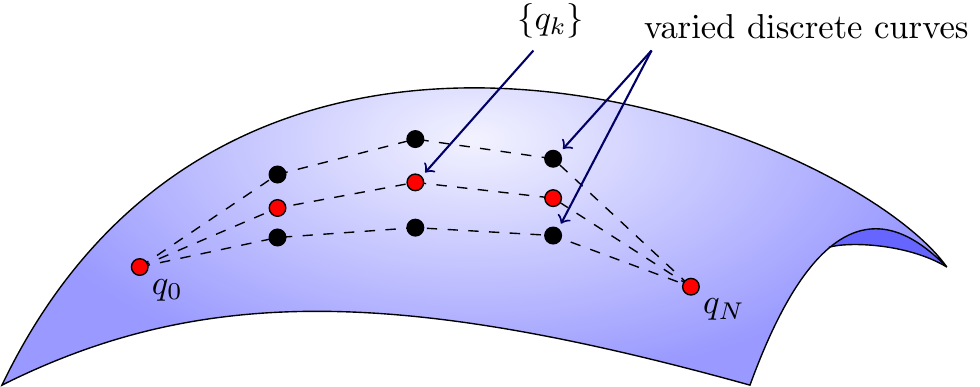}
\caption{Variations of the discrete trajectory $\{ q_{k} \}_{k=0}^{N}$.}
\end{figure}

The discrete trajectories $q_{d} = \{ q_{k} \}_{k=0}^{N}$ are required to satisfy a discrete version of Hamilton's principle of least action
\begin{align}
\delta \mcal{A}_{d} [q_{d}] = \delta \sum \limits_{k=0}^{N-1} L_{d} (q_{k}, q_{k+1}) = 0 .
\end{align}

The variation of the action is
\begin{align}\label{eq:vi_finite_variation_01}
\delta \mcal{A}_{d} [q_{d}]
&= \sum \limits_{k=0}^{N-1} \big[ D_{1} \, L_{d} (q_{k}, q_{k+1}) \cdot \delta q_{k} + D_{2} \, L_{d} (q_{k}, q_{k+1}) \cdot \delta q_{k+1} \big]
\end{align}
where $D_{i}$ denotes the derivative with respect to to the $i$th argument.
What follows corresponds to a discrete integration by parts, i.e., a reordering of the summation.
The $k=0$ term is removed from the first part of the sum and the $k=N-1$ term is removed from the second part
\begin{align}\label{eq:vi_finite_variation_02}
\delta \mcal{A}_{d} [q_{d}]
\nonumber
&= D_{1} \, L_{d} (q_{0}, q_{1}) \cdot \delta q_{0} + \sum \limits_{k=1}^{N-1} D_{1} \, L_{d} (q_{k}, q_{k+1}) \cdot \delta q_{k} \\
&+ \sum \limits_{k=0}^{N-2} D_{2} \, L_{d} (q_{k}, q_{k+1}) \cdot \delta q_{k+1} + D_{2} \, L_{d} (q_{N-1}, q_{N}) \cdot \delta q_{N}
.
\end{align}

As the variations at the endpoints, $\delta q_{0}$ and $\delta q_{N}$, are kept fixed, the corresponding terms vanish.
At last, the summation range of the second sum is shifted upwards by one with the arguments of the discrete Lagrangian adapted correspondingly
\begin{align}\label{eq:vi_finite_variation_03}
\delta \mcal{A}_{d} [q_{d}]
&=  \sum \limits_{k=1}^{N-1} \big[ D_{1} \, L_{d} (q_{k}, q_{k+1}) + D_{2} \, L_{d} (q_{k-1}, q_{k}) \big] \cdot \delta q_{k}
.
\end{align}

Hamilton's principle of least action requires the variation of the discrete action $\delta \mcal{A}_{d}$ to vanish for any choice of $\delta q_{k}$. Consequently, the expression in the square brackets of (\ref{eq:vi_finite_variation_03}) has to vanish. This defines the

\rimpeq[Discrete Euler-Lagrange Equations]{\label{eq:vi_finite_deleqs}
D_{1} \, L_{d} (q_{k}, q_{k+1}) + D_{2} \, L_{d} (q_{k-1}, q_{k}) = 0
}

The discrete Euler-Lagrange equations (\ref{eq:vi_finite_deleqs}) define an evolution map
\begin{align}\label{eq:vi_finite_evolution_map}
F_{L_{d}}
\; : \; \mf{Q} \times \mf{Q} \rightarrow \mf{Q} \times \mf{Q}
\; : \; ( q_{k-1}, q_{k} ) \mapsto ( q_{k}, q_{k+1} ) .
\end{align}

Starting from two configurations, $q_{0} \approx q (t_{0})$ and $q_{1} \approx q (t_{1} = t_{0} + h)$, the successive solution of the discrete Euler-Lagrange equations (\ref{eq:vi_finite_deleqs}) for $q_{2}$, $q_{3}$, etc., up to $q_{N}$, determines the discrete trajectory $\{ q_{k} \}_{k=0}^{N}$.
Quite often, however, it is more practical to prescribe an initial position and momentum instead of the configuration of the first two timesteps.
We therefore define the discrete momentum at timestep $k$ as\footnote{
The two expressions are equal by the discrete Euler-Lagrange equations (\ref{eq:vi_finite_deleqs}).
}
\begin{align}\label{eq:vi_finite_momentum}
p_{k} = - D_{1} L_{d} (q_{k}, q_{k+1}) = D_{2} L_{d} (q_{k-1}, q_{k}) ,
\end{align}

such that the variational integrator (\ref{eq:vi_finite_deleqs}) can be rewritten in

\begin{subequations}\label{eq:vi_finite_position_momentum}
\rimpeq[Position-Momentum Form]{
\label{eq:vi_finite_position_momentum_a}
p_{k  } &=           -  D_{1} L_{d} (q_{k}, q_{k+1}) \\
\label{eq:vi_finite_position_momentum_b}
p_{k+1} &= \hphantom{-} D_{2} L_{d} (q_{k}, q_{k+1})
}
\end{subequations}

Given $(q_{k}, p_{k})$, the first equation can be solved for $q_{k+1}$. This is generally a nonlinearly implicit equation that has to be solved by some iterative technique like Newton's method.
The second equation is an explicit function, so to obtain $p_{k+1}$ we merely have to plug in $q_{k}$ and $q_{k+1}$.
The corresponding Hamiltonian evolution map is
\begin{align}\label{eq:vi_finite_evolution_map_position_momentum}
\otilde{F}_{L_{d}}
\; : \; \cb{\mf{Q}} \rightarrow \cb{\mf{Q}}
\; : \; ( q_{k}, p_{k} ) \mapsto ( q_{k+1}, p_{k+1} ) .
\end{align}

Thus, starting with an initial position $q_{0}$ and an initial momentum $p_{0}$, the repeated solution of (\ref{eq:vi_finite_position_momentum}) gives the same discrete trajectory $\{ q_{k} \}_{k=0}^{N}$ as (\ref{eq:vi_finite_deleqs}).
The position-momentum form, as a one-step method, is usually easier to implement than the discrete Euler-Lagrange equations (\ref{eq:vi_finite_deleqs}). And for most problems, initial conditions are more naturally prescribed via the position and momentum of the particle at a given point in time, $(q_{0}, p_{0})$.
If, however, only the position of the particle at two points in time, $(q_{0}, q_{1})$, is known, the Euler-Lagrange equations (\ref{eq:vi_finite_deleqs}) are the more natural way of describing the dynamics.

This of course is just reflecting the difference in the Lagrangian and Hamiltonian point of view. For $n$ degrees of freedom, the variational principle leads to $n$ differential equations of second order. Hamilton's equations, on the other hand, are $2n$ differential equations of first order.
Which form eventually is used largely depends on the problem at hand.

\begin{fexample}[Example: Point Particle]
Consider a particle with mass $m$, moving in some potential $V$.
Its continuous Lagrangian is
\begin{align}\label{eq:vi_finite_example1_lagrangian}
L (q, \dot{q}) = \dfrac{1}{2} \, m \dot{q}^{2} - V(q)
\end{align}

Approximated by the trapezoidal rule, the discrete Lagrangian reads
\begin{align}\label{eq:vi_finite_example1_lagrangian_trapezoidal}
L_{d}^{\text{tr}} (q_{k}, q_{k+1}) = h \, \bigg[ \dfrac{m}{2} \bigg( \dfrac{q_{k+1} - q_{k}}{h} \bigg)^{2} - \dfrac{V (q_{k}) + V (q_{k+1})}{2} \bigg] .
\end{align}

Applying the discrete Euler-Lagrange equations (\ref{eq:vi_finite_deleqs}) to this expression results in discrete equations of motion
\begin{align}\label{eq:vi_finite_example1_deleqs_trapezoidal}
m \, \dfrac{q_{k+1} - 2 \, q_{k} + q_{k-1}}{h^{2}} = - \nabla V (q_{k})
\end{align}

which clearly are a discrete version of Newton's second law
\begin{align}\label{eq:vi_finite_example1_newton}
m \ddot{q} = - \nabla V = F .
\end{align}

For comparison, consider also the midpoint approximation
\begin{align}\label{eq:vi_finite_example1_lagrangian_midpoint}
L_{d}^{\text{mp}} (q_{k}, q_{k+1}) = h \, \bigg[ \dfrac{m}{2} \, \bigg( \dfrac{q_{k+1} - q_{k}}{h} \bigg)^{2} - V \bigg( \dfrac{q_{k} + q_{k+1}}{2} \bigg) \bigg]
\end{align}

which leads to
\begin{align}\label{eq:vi_finite_example1_deleqs_midpoint}
m \, \dfrac{q_{k+1} - 2 \, q_{k} + q_{k-1}}{h^{2}} = - \dfrac{1}{2} \, \bigg[ \nabla V \bigg( \dfrac{q_{k-1} + q_{k}}{2} \bigg) + \nabla V \bigg( \dfrac{q_{k} + q_{k+1}}{2} \bigg) \bigg]
\end{align}

and thus a different discretisation of (\ref{eq:vi_finite_example1_newton}).
The position-momentum form (\ref{eq:vi_finite_position_momentum}) of the trapezoidal Lagrangian (\ref{eq:vi_finite_example1_lagrangian_trapezoidal}) can be written as
\begin{subequations}
\begin{align}
\dfrac{q_{k+1} - q_{k}}{h} &= \; \dfrac{1}{m} \, \bigg[ p_{k} - \dfrac{h}{2} \, \nabla V (q_{k}) \bigg] \\
\dfrac{p_{k+1} - p_{k}}{h} &= -  \dfrac{1}{2} \, \bigg[ \nabla V (q_{k}) + \nabla V (q_{k+1}) \bigg]
\end{align}
\end{subequations}

and the one of the midpoint Lagrangian (\ref{eq:vi_finite_example1_lagrangian_midpoint}) reads
\begin{subequations}
\begin{align}
\dfrac{q_{k+1} - q_{k}}{h} &= \dfrac{1}{m} \, \bigg[ p_{k} - \dfrac{h}{2} \, \nabla V \bigg( \dfrac{q_{k} + q_{k+1}}{2} \bigg) \bigg] \\
\dfrac{p_{k+1} - p_{k}}{h} &= - \nabla V \bigg( \dfrac{q_{k} + q_{k+1}}{2} \bigg) .
\end{align}
\end{subequations}

This bears a close resemblance of Hamilton's equations of motion, where the additional term in the first equations can be interpreted as extrapolating the momentum $p_{k}$ to $t_{k+1/2}$.
As already noted, it is not always so easy to solve (\ref{eq:vi_finite_position_momentum_a}) for $q_{k+1}$, but in general this is an implicit equation.

\end{fexample}

\subsection{Discrete Tangent Space}

The discrete path space is defined as
\begin{align}
\mf{C}_{d} ( \mf{Q} )
= \big\{ c_{d} : \{ t_{k} \}_{k=0}^{N} \rightarrow \mf{Q} \big\}
\end{align}

where
\begin{align}
\{ t_{k} \}_{k=0}^{N} = \{ t_{k} = kh \; \vert \; k = 0, ..., N \} \subset \rsp
\end{align}

is an increasing sequence of times and $h$ is the discrete timestep.
$\mf{C}_{d} ( \mf{Q} )$ contains all possible discrete trajectories $c_{d}$ in $\mf{Q}$ and is isomorphic to $\mf{Q} \times ... \times \mf{Q}$ ($N+1$ copies)
\begin{align}
\mf{C}_{d} ( \mf{Q} ) \cong \bigtimes \limits_{N+1} \mf{Q} .
\end{align}

The subspace of $\mf{C}_{d} ( \mf{Q} )$ that contains all discrete trajectories with fixed endpoints $q_{0}$ and $q_{N}$ is defined as
\begin{align}
\mf{C}_{d} ( q_{0}, q_{N}, \{ t_{k} \}_{k=0}^{N} ) = \big\{ c_{d} : \{ t_{k} \}_{k=0}^{N} \rightarrow \mf{Q} \; \big\vert \; c(t_{0}) = q_{0}, c(t_{N}) = q_{N} \big\} .
\end{align}

The discrete action map $\mcal{A}_{d} : \mf{C}_{d} ( q_{0}, q_{N}, \{ t_{k} \}_{k=0}^{N} ) \rightarrow \rsp$ is defined as
\begin{align}
\mcal{A}_{d} [c_{d}]
&= \sum \limits_{k=0}^{N-1} L_{d} \big( c(t_{k}), c(t_{k+1}) \big)
= \sum \limits_{k=0}^{N-1} L_{d} ( q_{k}, q_{k+1} ) &
& \text{where} &
q_{k} &= c (t_{k}) . &
\end{align}

The tangent space $\tb[c_{d}]{\mf{C}_{d} ( q_{0}, q_{N}, \{ t_{k} \}_{k=0}^{N} )}$ to $\mf{C}_{d} ( q_{0}, q_{N}, \{ t_{k} \}_{k=0}^{N} )$ at $c_{d}$ contains the variations of $c_{d}$. It is defined as the set of maps
\begin{align}
v_{c_{d}} : \{ t_{k} \}_{k=0}^{N} &\rightarrow \tb{\mf{Q}} &
& \text{such that} &
\pi_{\mf{Q}} \circ v_{c_{d}} &= c_{d} &
& \text{and} &
v (t_{0}) = v (t_{N}) &= 0 , &
\end{align}

where $\pi_{\mf{Q}}$ is the canonical projection $\pi_{\mf{Q}} : \tb{\mf{Q}} \rightarrow \mf{Q}$ and local coordinates are given by
\begin{align}
v_{c_{d}} = \{ (q_{k}, v_{k}) \}_{k=0}^{N} .
\end{align}

$v_{c_{d}}$ is called a discrete variation of the discrete path $c_{d}$ and sometimes denoted $v_{c_{d}} = \delta c_{d}$.
The variation of the discrete action can therefore be formulated as
\begin{align}
\ext \mcal{A}_{d} [c_{d}] \cdot v_{c_{d}} = \sum \limits_{k=0}^{N-1} \big[ D_{1} L (q_{k}, q_{k+1}) \cdot v_{k} + D_{2} L (q_{k}, q_{k+1}) \cdot v_{k+1} \big] .
\end{align}

A reordering of the sum (discrete partial integration) leads to
\begin{align}
\ext \mcal{A}_{d} [c_{d}] \cdot v_{c_{d}}
\nonumber
&= \sum \limits_{k=1}^{N-1} \big[ D_{1} L (q_{k}, q_{k+1}) + D_{2} L (q_{k-1}, q_{k}) \big] \cdot v_{k} \\
&+ D_{1} L (q_{0}, q_{1}) \cdot v_{0} + D_{2} L (q_{N-1}, q_{N}) \cdot v_{N} .
\end{align}

where the terms in the second line vanish as $v_{0} = v (t_{0}) = 0$ as well as $v_{N} = v (t_{N}) = 0$ and thus
\begin{align}
\ext \mcal{A}_{d} [c_{d}] \cdot v_{c_{d}}
&= \sum \limits_{k=1}^{N-1} \big[ D_{1} L (q_{k}, q_{k+1}) + D_{2} L (q_{k-1}, q_{k}) \big] \cdot v_{k} .
\end{align}

The arbitrariness of the $v_{k}$ once more yields the discrete Euler-Lagrange equations.

\subsection{Discrete One- and Two-Form}\label{sec:vi_finite_cartan_form}

As in the continuous case, the discrete one-form is obtained by computing the variation of the action for varying endpoints
\begin{align}
\ext \mcal{A}_{d} [q_{d}] \cdot \delta q_{d}
\nonumber
&= \sum \limits_{k=0}^{N-1} \big[ D_{1} \, L_{d} (q_{k}, q_{k+1}) \cdot \delta q_{k} + D_{2} \, L_{d} (q_{k}, q_{k+1}) \cdot \delta q_{k+1} \big] \\
\nonumber
&= \sum \limits_{k=1}^{N-1} \big[ D_{1} \, L_{d} (q_{k}, q_{k+1}) + D_{2} \, L_{d} (q_{k-1}, q_{k}) \big] \cdot \delta q_{k} \\
&\hspace{3em}
+ D_{1} \, L_{d} (q_{0}, q_{1}) \cdot \delta q_{0} + D_{2} \, L_{d} (q_{N-1}, q_{N}) \cdot \delta q_{N} .
\end{align}

The two latter terms originate from the variation at the boundaries. They form the discrete counterpart of the Lagrangian one-form.
However, there are two boundary terms that define two distinct one-forms on $\mf{Q} \times \mf{Q}$
{%
\setlength{\arraycolsep}{2pt}%
\begin{align}
\begin{split}
\begin{array}{ll}
\Theta_{L_{d}}^{-} ( q_{0}  , q_{1} ) \cdot ( \delta q_{0}   , \delta q_{1} ) &\equiv           -  D_{1} L_{d} (q_{0},   q_{1}) \cdot \delta q_{0} , \\
\Theta_{L_{d}}^{+} ( q_{N-1}, q_{N} ) \cdot ( \delta q_{N-1} , \delta q_{N} ) &\equiv \hphantom{-} D_{2} L_{d} (q_{N-1}, q_{N}) \cdot \delta q_{N} .
\end{array}
\end{split}
\end{align}
}%

In general, these one-forms are defined as
\begin{align}\label{eq:vi_finite_discrete_one_form}
\begin{split}
\Theta_{L_{d}}^{-} ( q_{k} , q_{k+1} ) &\equiv           -  D_{1} L (q_{k}, q_{k+1}) , \\
\Theta_{L_{d}}^{+} ( q_{k} , q_{k+1} ) &\equiv \hphantom{-} D_{2} L (q_{k}, q_{k+1}) .
\end{split}
\end{align}

As $\ext L_{d} = \Theta_{L_{d}}^{+} - \Theta_{L_{d}}^{-}$ and $\ext^{2} L_{d} = 0$ one observes that
\begin{align}
\ext \Theta_{L_{d}}^{+} = \ext \Theta_{L_{d}}^{-}
\end{align}

such that the exterior derivative of both discrete one-forms defines the same \emph{discrete Lagrangian two-form} or \emph{discrete symplectic form}
\begin{align}\label{eq:vi_finite_discrete_two_form}
\Omega_{L_{d}}
&= \ext \Theta_{L_{d}}^{+}
= \ext \Theta_{L_{d}}^{-}
= \dfrac{\partial^{2} L_{d}}{\partial q_{k} \, \partial q_{k+1}} (q_{k}, q_{k+1}) \, dq_{k} \wedge dq_{k+1} &
& \text{(no summation)} .
\end{align}

\subsection{Preservation of the Discrete Symplectic Form}\label{sec:vi_finite_symplectic_form}

Consider the exterior derivative of the discrete action (\ref{eq:vi_finite_action}).
Upon insertion of the discrete Euler-Lagrange equations (\ref{eq:vi_finite_deleqs}) it becomes
\begin{align}\label{eq:vi_finite_symplectic_ext_action}
\ext \mcal{A}_{d}
&= D_{1} L_{d} (q_{0}, q_{1}) \cdot \form{d} q_{0} + D_{2} L_{d} (q_{N-1}, q_{N}) \cdot \form{d} q_{N}
= \Theta_{L_{d}}^{+} (q_{N-1}, q_{N}) - \Theta_{L_{d}}^{-} (q_{0}, q_{1}) .
\end{align}

On the right hand side we find the just defined Lagrangian one-forms (\ref{eq:vi_finite_discrete_one_form}).
Taking the exterior derivative of (\ref{eq:vi_finite_symplectic_ext_action}) gives
\begin{align}\label{eq:vi_finite_symplectic_preservation}
\Omega_{d} (q_{0}, q_{1}) = \Omega_{d} (q_{N-1}, q_{N}) ,
\end{align}

where $q_{N-1}$ and $q_{N}$ are connected with $q_{0}$ and $q_{1}$ through the discrete Euler-Lagrange equations (\ref{eq:vi_finite_deleqs}).
Therefore, (\ref{eq:vi_finite_symplectic_preservation}) implies that the discrete symplectic structure $\Omega_{d}$ is preserved while the system advances from $t=0$ to $t=Nh$ according to the discrete equations of motion (\ref{eq:vi_finite_deleqs}).
As the number of timesteps $N$ is arbitrary, the discrete symplectic form $\Omega_{d}$ is preserved at all times of the simulation.
Note that this does not automatically imply that the continuous symplectic structure $\Omega$ is preserved under the discrete map $F_{L_{d}}$.

\subsection{Composition Methods}\label{vi:finite_composition_methods}

The composition of a one-step variational integrator $\otilde{F}_{L_{d}}$ with different step sizes is a simple method of obtaining higher order schemes.
We assume that the initial scheme $\otilde{F}_{L_{d}}$ is symmetric, that is
\begin{align}
\otilde{F}_{L_{d}}^{h} \circ \otilde{F}_{L_{d}}^{-h} = \id ,
\end{align}

as this simplifies the construction.
Alternatively, efficient methods can also be built by combining a non-symmetric method with its adjoint.
The interested reader can find more information on these issues in \citeauthor{MarsdenWest:2001} \cite{MarsdenWest:2001} and \citeauthor{HairerLubichWanner:2006} \cite{HairerLubichWanner:2006}.
If a numerical method
\begin{align}
\otilde{F}_{L_{d}} : \cb{\mf{Q}} \times \rsp \rightarrow \cb{\mf{Q}}
\end{align}

is symmetric, it can be used to compose higher order methods by splitting up each timestep into $s$ substeps \cite{HairerLubichWanner:2006,McLachlan:1995,MarsdenWest:2001}
\begin{align}
\ohat{F}_{L_{d}}^{h} = \otilde{F}_{L_{d}}^{\gamma_{s} h} \circ ... \circ \otilde{F}_{L_{d}}^{\gamma_{i} h} \circ ... \circ \otilde{F}_{L_{d}}^{\gamma_{1} h}
\end{align}

where the careful selection of the $\gamma_{i}$ is crucial for the performance of the resulting scheme.
In this section, we show that a variational integrator is self-adjoint and thereby symmetric if its discrete Lagrangian is self-adjoint, a condition that is easily checked. We present some fourth and sixth order composition methods that can be applied in most situations.

It is worth mentioning that the composition can already be implemented at the level of the Lagrangian.
We will outline this at the end of the section, and later on, when we come to the discrete Noether theorem, the consequences for the discrete conservation laws are described.

\subsubsection{Adjoint of a Method and Adjoint Lagrangians}

In this subsection, we show that a variational integrator is symmetric if its Lagrangian is self-adjoint.
The adjoint $\otilde{F}_{L_{d}}^{*}$ of a method $\otilde{F}_{L_{d}}$ is defined as
\begin{align}
\big( \otilde{F}_{L_{d}}^{*} (h) \big) \circ \big( \otilde{F}_{L_{d}} (-h) \big) = \id.
\end{align}

A method $\otilde{F}_{L_{d}}$ is self-adjoint if $\otilde{F}_{L_{d}}^{*} = \otilde{F}_{L_{d}}$, therefore a self-adjoint method is also symmetric.
We would like to establish a condition of the discrete Lagrangian that tells us if the resulting method is self-adjoint or not.
We therefore define the adjoint Lagrangian $L_{d}^{*}$ of a discrete Lagrangian $L_{d}$ as
\begin{align}\label{eq:vi_finite_composition_adjoint_lagrangian}
L_{d}^{*} (q_{k}, q_{k+1}, h) \equiv - L_{d} (q_{k+1}, q_{k}, -h) .
\end{align}

Hence the Lagrangian $L_{d}$ is self-adjoint if
\begin{align}\label{eq:vi_finite_composition_selfadjoint_lagrangian}
L_{d} (q_{k}, q_{k+1}, h) = - L_{d} (q_{k+1}, q_{k}, -h) .
\end{align}

We want to show that if a discrete Lagrangian is self-adjoint so is the resulting method.
We start by establishing that adjoint Lagrangians admit adjoint methods, i.e., if $\otilde{F}_{L_{d}}$ is the Hamiltonian map resulting from $L_{d}$ and $\otilde{F}_{L_{d}^{*}}$ is the map resulting from $L_{d}^{*}$ then $\otilde{F}_{L_{d}}^{*} = \otilde{F}_{L_{d}^{*}}$.
In position momentum-form (\ref{eq:vi_finite_position_momentum}), the map $\otilde{F}_{L_{d}}$ is defined as
\begin{align}\label{eq:vi_finite_composition_original_method}
\otilde{F}_{L_{d}} :
\begin{cases}
\, p_{k  } &=           -  D_{1} L_{d} (q_{k}, q_{k+1}, h) \\
\, p_{k+1} &= \hphantom{-} D_{2} L_{d} (q_{k}, q_{k+1}, h)
\end{cases} .
\end{align}

Its adjoint method $\otilde{F}_{L_{d}}^{*} (h) = \big( \otilde{F}_{L_{d}} (-h) \big)^{-1}$ is the map
\begin{align}\label{eq:vi_finite_composition_adjoint_method1}
\otilde{F}_{L_{d}}^{*} :
\begin{cases}
\, p_{k  } &= \hphantom{-} D_{2} L_{d} (q_{k+1}, q_{k}, -h) \\
\, p_{k+1} &=           -  D_{1} L_{d} (q_{k+1}, q_{k}, -h)
\end{cases} .
\end{align}

And the map $\otilde{F}_{L_{d}^{*}}$ corresponding to the adjoint Lagrangian (\ref{eq:vi_finite_composition_adjoint_lagrangian}) is
\begin{align}\label{eq:vi_finite_composition_adjoint_method2}
\otilde{F}_{L_{d}^{*}} :
\begin{cases}
\, p_{k  } &=           -  D_{1} L_{d}^{*} (q_{k}, q_{k+1}, h) \\
\, p_{k+1} &= \hphantom{-} D_{2} L_{d}^{*} (q_{k}, q_{k+1}, h)
\end{cases} .
\end{align}

Computing the derivatives of the definition of the adjoint Lagrangian (\ref{eq:vi_finite_composition_adjoint_lagrangian})
\begin{align}
\begin{split}
\hphantom{-} D_{1} L_{d}^{*} (q_{k}, q_{k+1}, h) &=           -  D_{2} L_{d} (q_{k+1}, q_{k}, -h) \\
-  D_{2} L_{d}^{*} (q_{k}, q_{k+1}, h) &= \hphantom{-} D_{1} L_{d} (q_{k+1}, q_{k}, -h)
\end{split}
\end{align}

establishes the equality of (\ref{eq:vi_finite_composition_adjoint_method1}) and (\ref{eq:vi_finite_composition_adjoint_method2}), i.e., if two Lagrangians are adjoint so are the resulting methods.
Computing the derivatives of the definition of the self-adjoint Lagrangian (\ref{eq:vi_finite_composition_selfadjoint_lagrangian})
\begin{align}
\hphantom{-} D_{1} L_{d} (q_{k}, q_{k+1}, h) &=           -  D_{2} L_{d} (q_{k+1}, q_{k}, -h) \\
-  D_{2} L_{d} (q_{k}, q_{k+1}, h) &= \hphantom{-} D_{1} L_{d} (q_{k+1}, q_{k}, -h)
\end{align}

establishes the equality of (\ref{eq:vi_finite_composition_original_method}) and (\ref{eq:vi_finite_composition_adjoint_method1}), i.e., if the Lagrangian is self-adjoint so is the resulting method.
We have thereby obtained a condition for symmetry of a variational integrator that can easily be checked.

\subsubsection{Fourth Order Composition Methods}

If $\otilde{F}_{L_{d}}$ is a method of order $r$, a method $\ohat{F}_{L_{d}}$ of order $r+2$ is obtained by the composition \cite{HairerLubichWanner:2006}
\begin{align}\label{eq:vi_composition_4o3s}
\ohat{F}_{L_{d}}^{h} &= \otilde{F}_{L_{d}}^{\gamma h} \circ \otilde{F}_{L_{d}}^{(1-2\gamma) h} \circ \otilde{F}_{L_{d}}^{\gamma h} &
& \text{with} &
\gamma &= (2 - 2^{1/(r+1)})^{-1} . &
\end{align}

Hence, if $\otilde{F}_{L_{d}}$ is of second order, the resulting method $\ohat{F}_{L_{d}}$ will be of fourth order.
Note that symmetric methods are always of even order (for details see \citeauthor{MarsdenWest:2001} \cite{MarsdenWest:2001}).
A method of the same order but with generally smaller errors is obtained by considering five steps
\begin{align}\label{eq:vi_composition_4o5s}
\ohat{F}_{L_{d}}^{h} &= \otilde{F}_{L_{d}}^{\gamma h} \circ \otilde{F}_{L_{d}}^{\gamma h} \circ \otilde{F}_{L_{d}}^{(1-4\gamma) h} \circ \otilde{F}_{L_{d}}^{\gamma h} \circ \otilde{F}_{L_{d}}^{\gamma h} &
& \text{with} &
\gamma &= (4 - 4^{1/(r+1)})^{-1} . &
\end{align}

Multiple application of these compositions yields methods of orders higher than four.

\subsubsection{Sixth Order Composition Methods}

Higher order compositions can also be constructed directly (see \citeauthor{HairerLubichWanner:2006} \cite{HairerLubichWanner:2006}, section 3.2).
A sixth order method with seven substeps is given by
\begin{align}\label{eq:vi_composition_6o7s}
\begin{split}
\gamma_{1} = \gamma_{7} &= + 0.78451361047755726381949763 , \\
\gamma_{2} = \gamma_{6} &= + 0.23557321335935813368479318 , \\
\gamma_{3} = \gamma_{5} &= - 1.17767998417887100694641568 , \\
\gamma_{4} &= + 1.31518632068391121888424973 ,
\end{split}
\end{align}

but again smaller errors can be achieved by using nine steps
\begin{align}\label{eq:vi_composition_6o9s}
\begin{split}
\gamma_{1} = \gamma_{9} &= + 0.39216144400731413927925056 , \\
\gamma_{2} = \gamma_{8} &= + 0.33259913678935943859974864 , \\
\gamma_{3} = \gamma_{7} &= - 0.70624617255763935980996482 , \\
\gamma_{4} = \gamma_{6} &= + 0.08221359629355080023149045 , \\
\gamma_{5} &= + 0.79854399093482996339895035 .
\end{split}
\end{align}

The computational effort of these high order methods is quite large. Each step requires the solution of a nonlinear system of equations.
Given the outstanding performance already second order variational integrators are able to deliver, the necessity for such high order methods is probably rarely found.
Nevertheless, if extremely high accuracy is indispensable, these methods can be applied.

\subsubsection{Composite Discrete Lagrangians}

The composition schemes presented can all be derived as Euler-Lagrange equations from a composite discrete Lagrangian.
There are several equivalent possibilities of constructing such a Lagrangian and the corresponding discrete Euler-Lagrange equations. We will present only one, for details on the alternatives see \citeauthor{MarsdenWest:2001} \cite{MarsdenWest:2001}, section 2.5.

The composite discrete Lagrangian of a method
\begin{align}\label{eq:vi_composite_lagrangian_method}
\ohat{F}_{L_{d}}^{h} = \otilde{F}_{L_{d}^{s}}^{\gamma_{s} h} \circ ... \circ \otilde{F}_{L_{d}^{i}}^{\gamma_{i} h} \circ ... \circ \otilde{F}_{L_{d}^{1}}^{\gamma_{1} h} .
\end{align}

with $s$ substeps can be written as
\begin{align}\label{eq:vi_composite_lagrangian_lagrangian}
\ohat{L}_{d} (q_{k}^{0}, q_{k}^{1}, ..., q_{k}^{s} ) = \sum \limits_{i=1}^{s} L_{d}^{i} (q_{k}^{i-1}, q_{k}^{i}, \gamma_{i} h)
\end{align}

where we identify $q_{k} = q_{k}^{0}$ and $q_{k+1} = q_{k}^{s}$ such that $q_{k}^{s} = q_{k+1}^{0}$.
The discrete action becomes
\begin{align}
\mcal{A}_{d} \big( \{ q_{k}^{0}, q_{k}^{1}, ..., q_{k}^{s} \}_{k=0}^{N-1} ) = \sum \limits_{k=0}^{N-1} \ohat{L}_{d} (q_{k}^{0}, q_{k}^{1}, ..., q_{k}^{s} )
\end{align}

and we obtain discrete Euler-Lagrange equations
\begin{align}
D_{2} L_{d}^{s} (q_{k-1}^{s-1}, q_{k-1}^{s}, \gamma_{s} h) + D_{1} L_{d}^{i} (q_{k}^{0} ,q_{k}^{1}, \gamma_{1} h) &= 0 \\
& \; \; \vdots \\
D_{2} L_{d}^{i} (q_{k}^{i-1}, q_{k}^{i}, \gamma_{i} h) + D_{1} L_{d}^{i+1} (q_{k}^{i} ,q_{k}^{i+1}, \gamma_{i+1} h) &= 0 \\
& \; \; \vdots \\
D_{2} L_{d}^{s} (q_{k}^{s-1}, q_{k}^{s}, \gamma_{s} h) + D_{1} L_{d}^{1} (q_{k+1}^{0} ,q_{k+1}^{1}, \gamma_{1} h) &= 0
.
\end{align}

The maps $\otilde{F}_{L_{d}^{i}}^{\gamma_{i} h}$ in the composition method (\ref{eq:vi_composite_lagrangian_method}) can therefore be written as
\begin{align}
\otilde{F}_{L_{d}^{i}}^{\gamma_{i} h} : (q_{k}^{i-1} , p_{k}^{i-1}) \mapsto (q_{k}^{i} , p_{k}^{i})
\end{align}

with
\begin{subequations}
\begin{align}
p_{k}^{i-1}            &=           -  D_{1} L_{d}^{i} ( q_{k}^{i-1} , q_{k}^{i} , \gamma_{i} h ) \\
p_{k}^{i\hphantom{-1}} &= \hphantom{-} D_{2} L_{d}^{i} ( q_{k}^{i-1} , q_{k}^{i} , \gamma_{i} h ) .
\end{align}
\end{subequations}

The existence of a composite Lagrangian $\ohat{L}_{d}$ corresponding to a composite method $\ohat{F}_{L_{d}}$ is important in the analysis of conserved quantities. The discrete Noether theorem (next section) has to be applied to the composite Lagrangian to determine the quantities that are discretely conserved to the order of the composition method.
It cannot be expected that the errors of the conserved quantities of the discrete Lagrangians $L_{d}^{i}$, which are used to build the composition scheme, scale with the order of the composition scheme.

\subsection{Discrete Noether Theorem}\label{sec:vi_finite_noether_theorem}

The discrete Noether theorem, just as the continous Noether theorem, draws the connection between symmetries of a discrete Lagrangian and quantities that are conserved by the discrete Euler-Lagrange equations or, equivalently, the discrete Lagrangian flow.
The continuous theory translates straight forwardly to the discrete case. Therefore, we repeat just the important steps, translated to the discrete setting.

\subsubsection{Discrete Noether Theorem for Particle Systems}

Consider a one parameter group of discrete curves $\{ q_{k}^{\eps} \}_{k=0}^{N}$ such that $q_{k}^{0} (q_{k}) = q_{k}$.
The discrete Lagrangian $L_{d}$ has a symmetry if it is invariant under this transformation
\begin{align}\label{eq:vi_noether_finite_1}
L_{d} \big( q_{k}^{\eps}, q_{k+1}^{\eps} \big) &= L_{d} \big( q_{k}, q_{k+1} \big) &
& \text{for all $\eps$ and $k$} . &
\end{align}

The direction of such a symmetry is
\begin{align}\label{eq:vi_noether_finite_2}
X_{k} = \dfrac{\partial q_{k}^{\eps}}{\partial \eps} \bigg\vert_{\eps = 0}
\end{align}

such that
\begin{align}\label{eq:vi_noether_finite_3}
\dfrac{d}{d \eps} L_{d} \big( q_{k}^{\eps}, q_{k+1}^{\eps} \big) \bigg\vert_{\eps = 0}
&= D_{1} L_{d} \big( q_{k}, q_{k+1} \big) \cdot X_{k  }
+ D_{2} L_{d} \big( q_{k}, q_{k+1} \big) \cdot X_{k+1}
.
\end{align}

If $\{ q_{k} \}$ solves the discrete Euler-Lagrange equations
\begin{align}\label{eq:vi_noether_finite_4}
D_{1} \, L_{d} (q_{k}, q_{k+1}) + D_{2} \, L_{d} (q_{k-1}, q_{k}) = 0
\end{align}

we can replace the first term on the right hand side of (\ref{eq:vi_noether_finite_3}) to get
\begin{align}\label{eq:vi_noether_finite_5}
0
&= - D_{2} L_{d} \big( q_{k-1}, q_{k} \big) \cdot X_{k}
+ D_{2} L_{d} \big( q_{k}, q_{k+1} \big) \cdot X_{k+1}
.
\end{align}

This amounts to a discrete conservation law of the form

\rimpeq{\label{eq:vi_noether_finite_6}
D_{2} L_{d} \big( q_{k-1}, q_{k} \big) \cdot X_{k} &= D_{2} L_{d} \big( q_{k}, q_{k+1} \big) \cdot X_{k+1} &
& \text{(\textbf{Discrete Noether Theorem})} .
}

It states that solutions $\{ q_{k} \}$ of the discrete Euler-Lagrange equations preserve the components of the momentum map $p_{k} = D_{2} L_{d} \big( q_{k-1}, q_{k} \big)$ in direction $X_{k}$.

\begin{fexample}[Example: Free Point Particle]

Consider a transformation that amounts to an infinitesimal spatial translation
\begin{align}
q^{\eps}_{k} = q_{k} + \eps \, X .
\end{align}

The discrete Lagrangian is invariant under this transformation
\begin{align}
L_{d} ( q^{1,\eps}, q^{2\eps} )
= \dfrac{h}{2} \, \bigg( \dfrac{q^{2} + \eps \, X - q^{1} - \eps \, X}{h} \bigg)^{2}
= \dfrac{h}{2} \, \bigg( \dfrac{q^{2} - q^{1}}{h} \bigg)^{2}
= L_{d} ( q^{1}, q^{2} ) ,
\end{align}

such that the symmetry condition is trivially fulfilled
\begin{align}
\dfrac{\partial}{\partial \eps} L_{d} \big( q_{k}^{\eps}, q_{k+1}^{\eps} \big) \bigg\vert_{\eps = 0} = 0 .
\end{align}

The discrete conservation law following from the symmetry of the Lagrangian under spatial translation
\begin{align}
\bigg( \dfrac{q_{k  } - q_{k-1}}{h} \bigg) \cdot X
= \bigg( \dfrac{q_{k+1} - q_{k  }}{h} \bigg) \cdot X
\end{align}

amounts to the preservation of the discrete momentum in direction of $X$.

\end{fexample}

\subsubsection{Energy}

In continuous particle dynamics, the conservation of energy follows from translational symmetry of the Lagrangian with respect to time.
In discrete particle dynamics, with a fixed timestep $h$, it is not possible to consider infinitesimal translations with respect to time. In the setting we described, it is therefore not possible to prove conservation of the discrete energy by applying Noether's theorem.
Indeed, most often we find that energy is not conserved exactly, but only approximately, in that the energy error is bounded by some threshold value. This behaviour is typical for symplectic methods (see e.g. \citeauthor{HairerLubichWanner:2006} \cite{HairerLubichWanner:2006} and references therein).

Nevertheless, it is possible to achieve and prove exact energy conservation by making the timestep $h$ a dynamical variable. It is thereby determined by the variational principle, such that energy is conserved exactly \cite{Kane:1999}. And in the Noether theorem, infinitesimal transformations of time can be considered as well.

However, we do not follow this path. Still we are interested in the energy conserving properties of our variational integrators.
We therefore ``read'' the expression for the discrete energy from the Lagrangian. In the particle case, the Hamiltonian is an explicit part of the Lagrangian, such that its discrete counterpart follows directly from the discretisation of the Lagrangian.

\subsubsection{Discrete Noether Theorem for Composite Lagrangians}

As already pointed out in section \ref{vi:finite_composition_methods}, special care has to be taken in the case of composite discrete Lagrangians (\ref{eq:vi_composite_lagrangian_lagrangian})
\begin{align}\label{eq:vi_noether_composition_1}
\ohat{L}_{d} (q_{k}^{0}, q_{k}^{1}, ..., q_{k}^{s} ) = \sum \limits_{i=1}^{s} L_{d}^{i} (q_{k}^{i-1}, q_{k}^{i}, \gamma_{i} h) .
\end{align}

$\ohat{L}_{d}$ might have different discrete expressions of the conserved momenta than the $L_{d}^{i}$.
And it might even have different conservational properties, i.e., not all conserved momenta of $\ohat{L}_{d}$ might be conserved by the $L_{d}^{i}$ or vice versa.
Therefore the discrete Noether theorem (\ref{eq:vi_noether_finite_6}) has to be applied to $\ohat{L}_{d}$.

Similarly, the discrete expression for the energy is different for $\ohat{L}_{d}$ and the $L_{d}^{i}$.
To clarify this, let us look at an example. Consider the midpoint Lagrangian (\ref{eq:vi_finite_midpoint})
\begin{align}\label{eq:vi_noether_composition_2}
L_{d}^{\text{mp}} (q_{k}, q_{k+1}) = h \, L \bigg( \dfrac{q_{k} + q_{k+1}}{2}, \dfrac{q_{k+1} - q_{k}}{h} \bigg) ,
\end{align}

and the fourth order, three step composition method (\ref{eq:vi_composition_4o3s}),
\begin{align}\label{eq:vi_noether_composition_3}
\ohat{F}_{L_{d}}^{h} &= \otilde{F}_{L_{d}}^{\gamma h} \circ \otilde{F}_{L_{d}}^{(1-2\gamma) h} \circ \otilde{F}_{L_{d}}^{\gamma h} &
& \text{with} &
\gamma &= (2 - 2^{1/(r+1)})^{-1} . &
\end{align}

The discrete Hamiltonian corresponds to
\begin{align}\label{eq:vi_noether_composition_4}
H_{d} = h \, H \bigg( \dfrac{q_{k} + q_{k+1}}{2}, \dfrac{q_{k+1} - q_{k}}{h} \bigg) ,
\end{align}

such that the discrete Hamiltonian of the composite Lagrangian is
\begin{align}\label{eq:vi_noether_composition_5}
\ohat{H}_{d}
= h \, \bigg[
\gamma \, H \bigg( \dfrac{q_{k}^{0} + q_{k}^{1}}{2}, \dfrac{q_{k}^{1} - q_{k}^{0}}{h} \bigg)
+ (1-2\gamma) \, H \bigg( \dfrac{q_{k}^{1} + q_{k}^{2}}{2}, \dfrac{q_{k}^{2} - q_{k}^{1}}{h} \bigg)
+ \gamma \, H \bigg( \dfrac{q_{k}^{2} + q_{k}^{3}}{2}, \dfrac{q_{k}^{3} - q_{k}^{2}}{h} \bigg)
\bigg] ,
\end{align}

where $q_{k}^{0} = q_{k}$ and $q_{k}^{3} = q_{k+1}$. Only the error of this composite Hamiltonian will scale with the order $r$ of the scheme.
That might at first seem surprising, as symplectic methods are supposed to conserve a direct discretisation of the continuous Hamiltonian to at least order $r$, i.e.,
\begin{align}\label{eq:vi_noether_composition_6}
H (q_{n}, p_{n}) = H (q_{0}, p_{0}) + \mcal{O} (h^{r})
\end{align}

for exponentially long time intervals $kh \leq e^{h_{0} / 2h}$ with some constant $h_{0}$.
However, in the proof of this relation (see \citeauthor{HairerLubichWanner:2006} \cite{HairerLubichWanner:2006} and references therein) it is assumed that the discrete flow map $\ohat{F}_{L_{d}}$ preserves the continuous symplectic form $\Omega$. But we have only proved that the discrete symplectic form $\Omega_{d}$ is preserved. We can therefore not assume that this result translates directly. What we can always assume (for conservative systems) is that there exists a discrete energy that is preserved to the the order of the scheme, and in the case of the composition schemes this is an expression analogous to (\ref{eq:vi_noether_composition_5}).

\section{Discrete Field Theory}\label{sec:vi_infinite}

The derivation of the discrete field theory is a straight forward generalisation of the derivation for particle dynamics.
The only difference is that the basic physical quantity is not the Lagrangian, defined on a one-dimensional ``grid'' of time, but the Lagrangian density, defined over a multidimensional grid of spacetime or phasespacetime.

Again, the starting point is the discretisation of the action integral and the Lagrangian density.
The discrete Lagrangian density approximates the (phase)spacetime integral of the continuous Lagrangian density over one cell of the (phase)spacetime grid, e.g. with one spatial dimension this is
\begin{align}\label{eq:vi_infinite_lagrangian1}
\mcal{L}_{d} (\phy_{i,k}, \phy_{i+1,k}, \phy_{i,k+1}, \phy_{i+1,k+1}) \approx \int \limits_{t_{k}}^{t_{k+1}} \int \limits_{x_{i}}^{x_{i+1}} \mcal{L} \big( \phy, \phy_{t}, \phy_{x} \big)
\end{align}

with the corresponding action being a sum over the whole grid
\begin{align}\label{eq:vi_infinite_action1}
\mcal{A}_{d} = \sum \limits_{k=0}^{N_{t}-1} \sum \limits_{i=0}^{N_{x}-1} \mcal{L}_{d} (\phy_{i,k}, \phy_{i+1,k}, \phy_{i,k+1}, \phy_{i+1,k+1}) .
\end{align}

\begin{figure}[b]
\centering
\includegraphics[width=.3\textwidth]{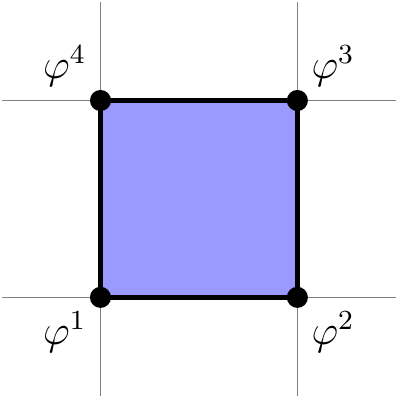}
\caption{Basic element of a two-dimensional spacetime grid.}
\label{fig:vi_infinite_gridbox}
\end{figure}

To make manipulations more tractable, the discrete Lagrangian density is rewritten in a slightly more abstract way, namely in terms of cells rather than grid points. Let us consider a cell determined by its vertices $(\phy^{1}, \phy^{2}, \phy^{3}, \phy^{4})$, like it is depicted in Fig. \ref{fig:vi_infinite_gridbox}.
For now, the horizontal axis shall be space, denoted by $x$, and the vertical axis shall be time, denoted by $t$.

Here, only a midpoint rule is considered, but the application of other quadrature rules is straight forward.
The fields $\phy$ are thus approximated by
\begin{align}\label{eq:vi_infinite_fields}
\phy (t,x) \approx \dfrac{1}{4} \big( \phy^{1} + \phy^{2} + \phy^{3} + \phy^{4} \big) .
\end{align}

For the approximation of the derivatives, there are in principle two possibilities for each coordinate, e.g. an $x$ derivative can be defined as $(\phy^{2} - \phy^{1}) / h_{x}$ as well as $(\phy^{3} - \phy^{4}) / h_{x}$.
Again, a midpoint-like averaging over the two possibilities is employed, such that
\begin{align}\label{eq:vi_infinite_derivatives}
\phy_{t} (t,x) &\approx \dfrac{1}{2} \bigg( \dfrac{\phy^{4} - \phy^{1}}{h_{t}} + \dfrac{\phy^{3} - \phy^{2}}{h_{t}} \bigg) , &
\phy_{x} (t,x) &\approx \dfrac{1}{2} \bigg( \dfrac{\phy^{2} - \phy^{1}}{h_{x}} + \dfrac{\phy^{3} - \phy^{4}}{h_{x}} \bigg) . &
\end{align}

Applying this to (\ref{eq:vi_infinite_lagrangian1}), the resulting discrete Lagrangian density reads
\begin{align}\label{eq:vi_infinite_lagrangian2}
\mcal{L}_{d} (\phy^{1}, \phy^{2}, \phy^{3}, \phy^{4}) \approx h_{t} \, h_{x} \, \mcal{L} \bigg( \dfrac{\phy^{1} + \phy^{2} + \phy^{3} + \phy^{4}}{4}, \dfrac{\phy^{4} - \phy^{1}}{2 h_{t}} + \dfrac{\phy^{3} - \phy^{2}}{2 h_{t}}, \dfrac{\phy^{2} - \phy^{1}}{2 h_{x}} + \dfrac{\phy^{3} - \phy^{4}}{2 h_{x}} \bigg)
\end{align}

and the discrete action becomes
\begin{align}\label{eq:vi_infinite_action2}
\mathcal{A}_{d} [ \phy_{d} ] = \sum \limits_{\substack{\text{grid}\\\text{boxes}}} \mcal{L}_{d} \big( \varphi^{1}, \varphi^{2}, \varphi^{3}, \varphi^{4} \big)
\end{align}

where $\phy_{d} = \{ \{ \phy_{i,k} \}_{i=0}^{N_{x}-1} \}_{k=0}^{N_{t}-1} $ is the discrete field.
The application of Hamilton's principle
\begin{align}\label{eq:vi_infinite_action3}
\dfrac{d}{d\eps} \mathcal{A}_{d} [ \phy_{d}^{\eps} ] \bigg\vert_{\eps = 0} = \dfrac{d}{d\eps} \sum \limits_{\substack{\text{grid}\\\text{boxes}}} \mcal{L}_{d} \big( \varphi^{1,\eps}, \varphi^{2,\eps}, \varphi^{3,\eps}, \varphi^{4,\eps} \big) \bigg\vert_{\eps = 0}
\end{align}
leads to discrete Euler-Lagrange field equations (DELFEQs) just as it lead to Euler-Lagrange equations in the continuous case.
With
\begin{align}\label{eq:vi_infinite_action4}
\delta \mcal{A}_{d} &\equiv \dfrac{d}{d\eps} \mathcal{A}_{d} [ \phy_{d}^{\eps} ] \bigg\vert_{\eps = 0} &
& \text{and} &
\delta \phy_{i,k} &\equiv \dfrac{d}{d\eps} \phy_{i,k}^{\eps} \bigg\vert_{\eps = 0} &
\end{align}

the variation of the action can be written as
\begin{align}\label{eq:vi_infinite_action5}
\delta \mathcal{A}_{d} = \sum \limits_{\substack{\text{grid}\\\text{boxes}}} \dfrac{\partial \mathcal{L}_{d}}{\partial \varphi^{a}} (\varphi^{1}, \varphi^{2}, \varphi^{3}, \varphi^{4}) \cdot \delta \varphi^{a} &&
(1 \leq a \leq 4) .
\end{align}

\begin{figure}[tb]
\centering
\includegraphics[width=.5\textwidth]{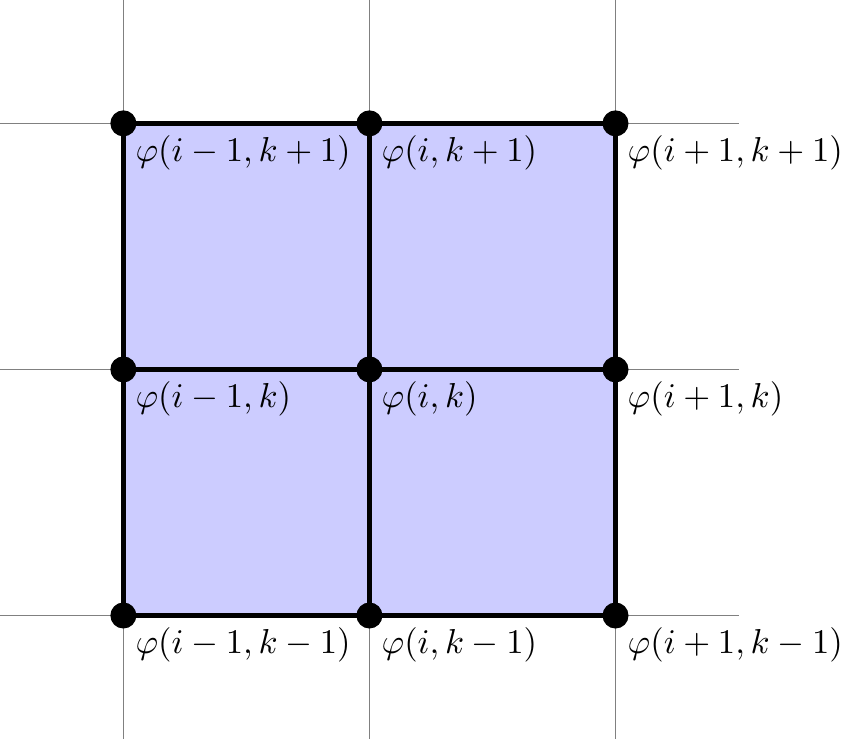}
\caption{Contributions of a specific grid point $(i,k)$ to the variation of the discrete action.}
\label{fig:vi_infinite_variation}
\end{figure}

As the variation of the action has to vanish for each $\delta \phy_{i,k}$ on the spacetime grid, it is sufficient to consider the contributions of $\delta \mathcal{A}_{d}$ that are multiplied by the variation of $\phy$ at a given grid point $(i,k)$
\begin{align}\label{eq:vi_infinite_action6}
\delta \mathcal{A}_{d} =
\nonumber
& \hspace{1em}
... + \dfrac{\partial \mcal{L}_d}{\partial \phy^1} \Big( \phy_{i,  k  }, \phy_{i+1,k  }, \phy_{i+1,k+1}, \phy_{i,  k+1} \Big) \cdot \delta \phy_{i,k} + ... \\
\nonumber
& \hspace{4em}
... + \dfrac{\partial \mcal{L}_d}{\partial \phy^2} \Big( \phy_{i-1,k  }, \phy_{i,  k  }, \phy_{i,  k+1}, \phy_{i-1,k+1} \Big) \cdot \delta \phy_{i,k} + ... \\
\nonumber
& \hspace{7em}
... + \dfrac{\partial \mcal{L}_d}{\partial \phy^3} \Big( \phy_{i-1,k-1}, \phy_{i,  k-1}, \phy_{i,  k  }, \phy_{i-1,k  } \Big) \cdot \delta \phy_{i,k} + ... \\
& \hspace{10em}
... + \dfrac{\partial \mcal{L}_d}{\partial \phy^4} \Big( \phy_{i,  k-1}, \phy_{i+1,k-1}, \phy_{i+1,k  }, \phy_{i,  k  } \Big) \cdot \delta \phy_{i,k} + ...
\hspace{1em} = 0 .
\end{align}

In total there are four such contributions, originating from the Lagrangian densities $\mcal{L}_{d}$ of the four squares that touch the point  $(i,k)$ as is depicted in Fig. \ref{fig:vi_infinite_variation}. The sum of the factors of $\delta \phy_{ik}$ corresponds to the

\rimpeq[Discrete Euler-Lagrange Field Equations]{\label{eq:vi_infinite_delfeqs}
0
\nonumber
&= \dfrac{\partial \mcal{L}_d}{\partial \phy^1} \Big( \phy_{i,  k  }, \phy_{i+1,k  }, \phy_{i+1,k+1}, \phy_{i,  k+1} \Big)
+ \dfrac{\partial \mcal{L}_d}{\partial \phy^2} \Big( \phy_{i-1,k  }, \phy_{i,  k  }, \phy_{i,  k+1}, \phy_{i-1,k+1} \Big) \\
&+ \dfrac{\partial \mcal{L}_d}{\partial \phy^3} \Big( \phy_{i-1,k-1}, \phy_{i,  k-1}, \phy_{i,  k  }, \phy_{i-1,k  } \Big)
+ \dfrac{\partial \mcal{L}_d}{\partial \phy^4} \Big( \phy_{i,  k-1}, \phy_{i+1,k-1}, \phy_{i+1,k  }, \phy_{i,  k  } \Big) .
}

\begin{example}[Example: Wave Equation]

The continuous Lagrangian density for the wave equations is
\begin{align}
\mathcal{L} \big( u_t(x,t), u_x(x,t) \big) = \dfrac{1}{2} \, \bigg( \dfrac{\partial u (x,t)}{\partial t} \bigg)^{2} - \dfrac{1}{2} \, \bigg( \dfrac{\partial u (x,t)}{\partial x} \bigg)^{2} .
\end{align}

A discretisation with the midpoint rule (\ref{eq:vi_infinite_derivatives}) leads to the discrete Lagrangian density
\begin{align}
\mathcal{L}_{d} \big( u^{1}, u^{2}, u^{3}, u^{4} \big)
&= \dfrac{1}{2} \, \bigg( \dfrac{u^{4} - u^{1}}{2 h_{t}} + \dfrac{u^{3} - u^{2}}{2 h_{t}} \bigg)^{2}
- \dfrac{1}{2} \, \bigg( \dfrac{u^{2} - u^{1}}{2 h_{x}} + \dfrac{u^{3} - u^{4}}{2 h_{x}} \bigg)^{2} .
\end{align}

The four contributions to the discrete Euler-Lagrange field equations (\ref{eq:vi_infinite_delfeqs}) are
\begin{subequations}
\begin{align}
\dfrac{\partial \mathcal{L}_d}{\partial u^1} \Big( u_{i,  k  }, u_{i+1,k  }, u_{i+1,k+1}, u_{i,  k+1} \Big) =
\nonumber
&- \dfrac{u_{i  , k+1} - u_{i  , k  }}{4 h_{t}^{2}} - \dfrac{u_{i+1, k+1} - u_{i+1, k  }}{4 h_{t}^{2}} \\
&+ \dfrac{u_{i+1, k  } - u_{i  , k  }}{4 h_{x}^{2}} + \dfrac{u_{i+1, k+1} - u_{i  , k+1}}{4 h_{x}^{2}} ,
\end{align}
\begin{align}
\dfrac{\partial \mathcal{L}_d}{\partial u^2} \Big( u_{i-1,k  }, u_{i,  k  }, u_{i,  k+1}, u_{i-1,k+1} \Big) =
\nonumber
&- \dfrac{u_{i-1, k+1} - u_{i-1, k  }}{4 h_{t}^{2}} - \dfrac{u_{i  , k+1} - u_{i  , k  }}{4 h_{t}^{2}} \\
&- \dfrac{u_{i  , k  } - u_{i-1, k  }}{4 h_{x}^{2}} - \dfrac{u_{i  , k+1} - u_{i-1, k+1}}{4 h_{x}^{2}} ,
\end{align}
\begin{align}
\dfrac{\partial \mathcal{L}_d}{\partial u^3} \Big( u_{i-1,k-1}, u_{i,  k-1}, u_{i,  k  }, u_{i-1,k  } \Big) =
\nonumber
&+ \dfrac{u_{i-1, k  } - u_{i-1, k-1}}{4 h_{t}^{2}} + \dfrac{u_{i  , k  } - u_{i  , k-1}}{4 h_{t}^{2}} \\
&- \dfrac{u_{i  , k-1} - u_{i-1, k-1}}{4 h_{x}^{2}} - \dfrac{u_{i  , k  } - u_{i-1, k  }}{4 h_{x}^{2}} ,
\end{align}
\begin{align}
\dfrac{\partial \mathcal{L}_d}{\partial u^4} \Big( u_{i,  k-1}, u_{i+1,k-1}, u_{i+1,k  }, u_{i,  k  } \Big) =
\nonumber
&+ \dfrac{u_{i  , k  } - u_{i  , k-1}}{4 h_{t}^{2}} + \dfrac{u_{i+1, k  } - u_{i+1, k-1}}{4 h_{t}^{2}} \\
&+ \dfrac{u_{i+1, k-1} - u_{i  , k-1}}{4 h_{x}^{2}} + \dfrac{u_{i+1, k  } - u_{i  , k  }}{4 h_{x}^{2}} .
\end{align}
\end{subequations}

Summing up all these terms, the discrete wave equation is obtained
\begin{multline}\label{eq:vi_infinite_wave_equation_discrete}
\dfrac{u_{i-1, k+1} - 2 \, u_{i-1, k  } + u_{i-1, k-1}}{4 h_{t}^{2}}
+ 2 \, \dfrac{u_{i  , k+1} - 2 \, u_{i  , k  } + u_{i  , k-1}}{4 h_{t}^{2}}
+      \dfrac{u_{i+1, k+1} - 2 \, u_{i+1, k  } + u_{i+1, k-1}}{4 h_{t}^{2}} = \\
=      \dfrac{u_{i+1, k+1} - 2 \, u_{i  , k+1} + u_{i-1, k+1}}{4 h_{x}^{2}}
+ 2 \, \dfrac{u_{i+1, k  } - 2 \, u_{i  , k  } + u_{i-1, k  }}{4 h_{x}^{2}}
+      \dfrac{u_{i+1, k-1} - 2 \, u_{i  , k-1} + u_{i-1, k-1}}{4 h_{x}^{2}} .
\end{multline}

This clearly is a discrete version of the continuous wave equation
\begin{align*}
\dfrac{\partial^{2} u (x,t)}{\partial t^{2}} = \dfrac{\partial^{2} u (x,t)}{\partial x^{2}} ,
\end{align*}

with the following stencil
\begin{align}
\dfrac{1}{4 h_{t}^{2}}
\begin{bmatrix}
\hphantom{-} 1 & \hphantom{-} 2 & \hphantom{-} 1 \\
-2 &             -4 &             -2 \\
\hphantom{-} 1 & \hphantom{-} 2 & \hphantom{-} 1
\end{bmatrix}
u
=
\dfrac{1}{4 h_{x}^{2}}
\begin{bmatrix}
\hphantom{-} 1 & -2 & \hphantom{-} 1 \\
\hphantom{-} 2 & -4 & \hphantom{-} 2 \\
\hphantom{-} 1 & -2 & \hphantom{-} 1
\end{bmatrix}
u .
\end{align}

We observe that the derivative with respect to one direction is averaged in the other direction, i.e., the time derivative is averaged over three neighbouring points in space, and the spatial derivative is averaged over three neighbouring points in time.
This averaging of derivatives is a common feature often found in variational integrators of field theories. It appears to be one of the decisive features that account for the superior performance of variational integrators.

\end{example}

\newpage 

\subsection{Discrete Jet Space}

The discrete phasespace, on which the Lagrangian density is defined, is the discrete first jet bundle $\jb{1}{\mf{Y}}$
\begin{align}\label{vi_infinite_discrete_jet_space_1}
\mcal{L}_{d} : \jb{1}{\mf{Y}} \rightarrow \rsp .
\end{align}

To understand its structure, some considerations are in order.
The discrete spacetime of dimension two is the grid of points
\begin{align}\label{vi_infinite_discrete_jet_space_2}
\mf{X} = \mbb{Z} \times \mbb{Z} = \big\{ (i,k) \big\} \cong \big\{ ( i h_{x} , k h_{t} ) \; \big\vert \; i, k \in \mbb{Z} \big\} .
\end{align}

It corresponds to a grid with elements $x_{i,j}$ in continuous spacetime.
The discrete fibre bundle over $\mf{X}$ is
\begin{align}\label{vi_infinite_discrete_jet_space_3}
\mf{Y} = \mf{X} \times \mf{F} .
\end{align}

where $\mf{F}$ is a smooth manifold. Elements of $\mf{Y}$ over the point $(i,k)$ are denoted by $y_{i,k}$ and the projection $\pi_{\mf{X} \mf{Y}}$ is given by
\begin{align}\label{vi_infinite_discrete_jet_space_4}
\pi_{\mf{X} \mf{Y}} (y_{i,k}) = (i,k) .
\end{align}

A square $\square$ on $\mf{X}$ is an ordered quadruplet
\begin{align}\label{vi_infinite_discrete_jet_space_5}
\square = \big( (i, k), (i+1, k), (i+1, k+1), (i, k+1) \big) ,
\end{align}

defining a grid cell, c.f. figure \ref{fig:vi_infinite_gridbox}.
The first component of $\square$, denoted $\square^{1}$, is the first vertex of the square, with equivalent definitions for the other three vertices
\begin{align}\label{vi_infinite_discrete_jet_space_6}
\square^{1} &= (i, k) , &
\square^{2} &= (i+1, k) , &
\square^{3} &= (i+1, k+1) , &
\square^{4} &= (i, k+1) . &
\end{align}

A section $\phy$ of $\mf{Y}$ is a map
\begin{align}\label{vi_infinite_discrete_jet_space_7}
\phy : \mf{U} \subseteq \mf{X} \rightarrow \mf{Y}
\end{align}

such that
\begin{align}\label{vi_infinite_discrete_jet_space_8}
\pi_{\mf{X} \mf{Y}} \circ \phy = \id_{\mf{U}}
\end{align}

and the $\phy^{a}$ with $a \in \{ 1, 2, 3, 4 \}$ of figure \ref{fig:vi_infinite_gridbox} correspond to
\begin{align}\label{vi_infinite_discrete_jet_space_9}
\phy^{1} &= \phy (\square^{1}) = \phy_{\square^{1}} , &
\phy^{2} &= \phy (\square^{2}) = \phy_{\square^{2}} , &
\phy^{3} &= \phy (\square^{3}) = \phy_{\square^{3}} , &
\phy^{4} &= \phy (\square^{4}) = \phy_{\square^{4}} .
\end{align}

The set of squares on $\mf{X}$ is denoted $\mf{X}^{\square}$.
The first jet bundle of $\mf{Y}$ is given by
\begin{align}\label{vi_infinite_discrete_jet_space_10}
\jb{1}{\mf{Y}}
\nonumber
&\equiv \big\{ \big( \square , (y_{i, k}, y_{i+1, k}, y_{i+1, k+1}, y_{i, k+1}) \big) \; \big\vert \; \square \in \mf{X}^{\square}, \; (i,k) = \square^{1} , \; y_{i, k}, y_{i+1, k}, y_{i+1, k+1}, y_{i, k+1} \in \mf{F} \big\} \\
&\equiv \mf{X}^{\square} \times \mf{F}^{4} .
\end{align}

The first jet prolongation of a section $\phy$ on $\mf{Y}$ is the map
\begin{align}\label{vi_infinite_discrete_jet_space_11}
j^{1} \phy : \mf{X}^{\square} \mapsto \jb{1}{\mf{Y}}
\end{align}

defined by
\begin{align}\label{vi_infinite_discrete_jet_space_12}
j^{1} \phy (\square) \equiv \big( \square, \phy (\square^{1}), \phy (\square^{2}), \phy (\square^{3}), \phy (\square^{4}) \big) .
\end{align}

The first jet is defined to include first order derivatives, that at the discrete level are functions of $\square$ and $\phy (\square^{a})$.
The restriction of a vector field $V$ on $\mf{Y}$ to the fibre $\mf{Y}_{i,k}$ is denoted $V_{i,k}$, and similarly for vector fields on $\jb{1}{\mf{Y}}$.
The first jet prolongation of a vector field $V$ on $\mf{Y}$ is the vector field $j^{1} V$ on $\jb{1}{\mf{Y}}$, defined by
\begin{align}\label{vi_infinite_discrete_jet_space_13}
j^{1} V ( \square ) \equiv \big( \square, V_{\square^{1}} (y_{\square^{1}}), V_{\square^{2}} (y_{\square^{2}}), V_{\square^{3}} (y_{\square^{3}}), V_{\square^{4}} (y_{\square^{4}}) \big) ,
\end{align}

for any square $\square$.

A point $(i,k) \in \mf{X}$ is \emph{touched} by a square, if it is a vertex of that square. A point $(i,k) \in \mf{U} \subseteq \mf{X}$ is an \emph{interior point} of $\mf{U}$, if $\mf{U}$ contains all four squares of $\mf{X}$ that touch $(i,k)$. The \emph{interior} $\mrm{int} \, \mf{U}$ of $\mf{U}$ is the collection of all interior points of $\mf{U}$. The closure $\mrm{cl} \, \mf{U}$ of $\mf{U}$ is the union of all squares touching interior points of $\mf{U}$. A point $(i,k)$ is a \emph{boundary point} of $\mf{U}$ if it is a point in both, $\mf{U}$ and $\mrm{cl} \, \mf{U}$, which is not an interior point. The boundary $\partial \mf{U}$ of $\mf{U}$ is the set of boundary points of $\mf{U}$, such that
\begin{align}
\partial \mf{U} \equiv ( \mf{U} \cap \mrm{cl} \, \mf{U} ) \, \backslash \, \mrm{int} \, \mf{U} .
\end{align}

\subsubsection{Discrete Action Principle}

The Lagrangian density on a given square is a function
\begin{align}\label{eq:vi_infinite_discrete_jet_space_action_principle_1}
\mcal{L}_{\square} : \mf{F}^{4} \rightarrow \rsp
\end{align}

defined as
\begin{align}\label{eq:vi_infinite_discrete_jet_space_action_principle_2}
\mcal{L}_{\square} \big( y^{1}, y^{2}, y^{3}, y^{4} \big) \equiv \mcal{L}_{d} \big( \square, y^{1}, y^{2}, y^{3}, y^{4} \big) .
\end{align}

Thus, the discrete Lagrangian density $\mcal{L}_{d}$ can be regarded as the choice of a function $\mcal{L}_{\square}$ on each square $\square$ of $\mf{X}$.
The variables on the domain of $\mcal{L}_{\square}$ are denoted $y^{1}, y^{2}, y^{3}, y^{4}$, independently of the actual $\square$.

If $\mf{C} (\mf{Y})$ is the set of sections of $\mf{Y}$ on a subset $\mf{U} \subseteq \mf{X}$, the discrete action $\mcal{A}_{d}$ is a real-valued function on $\mf{C} (\mf{Y})$, defined by
\begin{align}\label{eq:vi_infinite_discrete_jet_space_action_principle_3}
\mcal{A}_{d} [\phy] = \sum \limits_{\substack{\square\\ \square \subseteq \mf{U}}} \mcal{L}_{d} \circ j^{1} \phy (\square) .
\end{align}

The variations $\phy^{\eps}$ of a section $\phy \in \mf{C} (\mf{Y})$ are described by a vertical map $\eta^{\eps}$ and its generating vector field $V \in \tb[\phy]{\mf{C} (\mf{Y})}$.
The map $\phy^{\eps}$ corresponds to a one-parameter family of sections
\begin{align}\label{eq:vi_infinite_discrete_jet_space_action_principle_4}
\phy^{\eps} (i,k) \equiv ( \eta^{\eps} \circ \phy ) (i,k) \equiv \eta^{\eps}_{ik} \big( \phy (i,k) \big) ,
\end{align}

where $\eta_{ik}^{\eps}$ is the flow of $V_{ik}$ on $\mf{Y}_{ik}$.
The action principle is to seek those sections $\phy$ for which
\begin{align}\label{eq:vi_infinite_discrete_jet_space_action_principle_5}
\ext \mcal{A}_{d} = \dfrac{d}{d\eps} \mcal{A}_{d} [\phy^{\eps}] \bigg\vert_{\eps = 0} = 0
\end{align}

for all vector fields $V$ on $\mf{C} (\mf{Y})$.
By focusing upon a fixed $(i,k) \in \mrm{int}(\mf{U})$ and the same arguments as in the previous section, c.f. equations (\ref{eq:vi_infinite_action5}) and (\ref{eq:vi_infinite_action6}), we find the discrete Euler-Lagrange field equations for all $(i,k) \in \mrm{int}(\mf{U})$

\rimpeq{\label{eq:vi_infinite_discrete_jet_space_action_principle_6}
\hspace{-1.5em}
\sum \limits_{\substack{a\\ (i,k) = \square^{a}}} \dfrac{\partial \mcal{L}_{\square}}{\partial y^{a}} \big( \phy_{\square^{1}}, \phy_{\square^{2}}, \phy_{\square^{3}}, \phy_{\square^{4}} \big) &= 0 &
& \text{(\textbf{Discrete Euler-Lagrange Equations})} .
}

\subsection{Discrete Cartan Form}

Allowing for nonzero variations on the boundary $\partial \mf{U}$ will lead us to the discrete Cartan form.
In that case, the vector field $V$ does not necessarily vanish on $\partial \mf{U}$.

For each point $(i,k)$ of the boundary $\partial \mf{U}$, find the squares in $\mf{U}$ that touch $(i,k)$.
There is at least one such square since $(i,k) \in \mrm{cl} (\mf{U})$, but not four such squares since $(i,k) \notin \mrm{int} (\mf{U})$.
For each of the touching squares, $(i,k)$ occurs at the $a$'th vertex for at most three values of $a = \{ 1, 2, 3, 4 \}$, such that the expressions
\begin{subequations}\label{eq:vi_infinite_discrete_jet_space_cartan_form_1}
\begin{align}
\dfrac{\partial \mcal{L}_d}{\partial y^1} \Big( \phy_{i,  k  }, \phy_{i+1,k  }, \phy_{i+1,k+1}, \phy_{i,  k+1} \Big) \, V_{ik} (\phy_{ik}) , \\
\dfrac{\partial \mcal{L}_d}{\partial y^2} \Big( \phy_{i-1,k  }, \phy_{i,  k  }, \phy_{i,  k+1}, \phy_{i-1,k+1} \Big) \, V_{ik} (\phy_{ik}) , \\
\dfrac{\partial \mcal{L}_d}{\partial y^3} \Big( \phy_{i-1,k-1}, \phy_{i,  k-1}, \phy_{i,  k  }, \phy_{i-1,k  } \Big) \, V_{ik} (\phy_{ik}) , \\
\dfrac{\partial \mcal{L}_d}{\partial y^4} \Big( \phy_{i,  k-1}, \phy_{i+1,k-1}, \phy_{i+1,k  }, \phy_{i,  k  } \Big) \, V_{ik} (\phy_{ik}) ,
\end{align}
\end{subequations}

give at most three contributions.
The total contribution to $\ext \mcal{A}_{d}$ from the boundary is the sum of all such terms.
In discrete particle mechanics, we found two one-forms. The above list suggests that there are four Cartan forms, which we define to be
\begin{subequations}\label{eq:vi_infinite_discrete_jet_space_cartan_form_2}
\begin{multline}
\Theta_{L_{d}}^{1} \big( \phy_{i,  k  }, \phy_{i+1,k  }, \phy_{i+1,k+1}, \phy_{i,  k+1} \big) \cdot \big( v_{y_{i,k}}, v_{y_{i+1,k}}, v_{y_{i+1,k+1}}, v_{y_{i,k+1}} \big) = \\
= \dfrac{\partial \mcal{L}_d}{\partial y^1} \big( \phy_{i,  k  }, \phy_{i+1,k  }, \phy_{i+1,k+1}, \phy_{i,  k+1} \big) \cdot \big( v_{y_{ij}}, 0, 0, 0 \big) ,
\end{multline}
\vspace{-2em}
\begin{multline}
\Theta_{L_{d}}^{2} \big( \phy_{i-1,k  }, \phy_{i,  k  }, \phy_{i,  k+1}, \phy_{i-1,k+1} \big) \cdot \big( v_{y_{i-1,k}}, v_{y_{i,k}}, v_{y_{i,k+1}}, v_{y_{i-1,k+1}} \big) = \\
= \dfrac{\partial \mcal{L}_d}{\partial y^2} \big( \phy_{i-1,k  }, \phy_{i,  k  }, \phy_{i,  k+1}, \phy_{i-1,k+1} \big) \cdot \big( 0, v_{y_{ij}}, 0, 0 \big) ,
\end{multline}
\vspace{-2em}
\begin{multline}
\Theta_{L_{d}}^{3} \big( \phy_{i-1,k-1}, \phy_{i,  k-1}, \phy_{i,  k  }, \phy_{i-1,k  } \big) \cdot \big( v_{y_{i-1,k-1}}, v_{y_{i,k-1}}, v_{y_{i,k}}, v_{y_{i-1,k}} \big) = \\
= \dfrac{\partial \mcal{L}_d}{\partial y^3} \big( \phy_{i-1,k-1}, \phy_{i,  k-1}, \phy_{i,  k  }, \phy_{i-1,k  } \big) \cdot \big( 0, 0, v_{y_{ij}}, 0 \big) ,
\end{multline}
\vspace{-2em}
\begin{multline}
\Theta_{L_{d}}^{4} \big( \phy_{i,  k-1}, \phy_{i+1,k-1}, \phy_{i+1,k  }, \phy_{i,  k  } \big) \cdot \big( v_{y_{i,k-1}}, v_{y_{i+1,k-1}}, v_{y_{i+1,k}}, v_{y_{i,k}} \big) = \\
= \dfrac{\partial \mcal{L}_d}{\partial y^4} \big( \phy_{i,  k-1}, \phy_{i+1,k-1}, \phy_{i+1,k  }, \phy_{i,  k  } \big) \cdot \big( 0, 0, 0, v_{y_{ij}} \big) .
\end{multline}
\end{subequations}

This quadruple $\big( \Theta_{L_{d}}^{1}, \Theta_{L_{d}}^{2}, \Theta_{L_{d}}^{3}, \Theta_{L_{d}}^{4} \big)$ is regarded as the discrete counterpart of the Cartan form $\Theta_{L}$ from (\ref{eq:classical_fields_cartan_form_53}).
For a vector field $V$ from $\tb{\mf{C} (\mf{Y})}$, the expressions from the list become
\begin{multline}\label{eq:vi_infinite_discrete_jet_space_cartan_form_3}
\big[ (j^{1} \phy)^{*} ( \iprod_{j^{1} V} \Theta_{L_{d}}^{a} ) \big] (\square) = \\
= \Theta_{L_{d}}^{a} \big( \phy_{\square^{1}}, \phy_{\square^{2}}, \phy_{\square^{3}}, \phy_{\square^{4}} \big) \cdot \big( V_{\square^{1}} (y_{\square^{1}}), V_{\square^{2}} (y_{\square^{2}}), V_{\square^{3}} (y_{\square^{3}}), V_{\square^{4}} (y_{\square^{4}}) \big) .
\end{multline}

We collect all the contributions from the boundary into one single object $\theta_{L_{d}} (\phy) \cdot V$, where $\theta_{L_{d}}$ is the one-form on the space of sections $\mf{C} (\mf{Y})$, defined by
\begin{align}\label{eq:vi_infinite_discrete_jet_space_cartan_form_4}
\theta_{L_{d}} (\phy) \cdot V \equiv
\sum \limits_{\substack{\square\\ \square \cap \partial \mf{U} \neq \emptyset}} \bigg(
\sum \limits_{\substack{a\\ \square^{a} \in \partial \mf{U}}}
\Big[ (j^{1} \phy)^{*} ( \iprod_{j^{1} V} \Theta_{L_{d}}^{a} ) \Big] (\square)
\bigg) .
\end{align}

Here we sum over all squares $\square$ that touch the boundary, and all vertices $\square^{a}$ of those squares that are elements of the boundary.
With this, the variation of the discrete action (\ref{eq:vi_infinite_discrete_jet_space_action_principle_5}) can be written as
{%
\setlength{\arraycolsep}{2pt}%
\begin{align}
\label{eq:vi_infinite_discrete_jet_space_cartan_form_5}
\ext \mcal{A}_{d} \cdot V
\nonumber
&= \sum \limits_{\substack{\square\\ \square \cap \mrm{int} \, \mf{U} \neq \emptyset}} \bigg(
\sum \limits_{\substack{a\\ \square^{a} \in \mrm{int} \, \mf{U}}}
\Big[ (j^{1} \phy)^{*} ( \iprod_{j^{1} V} \Theta_{L_{d}}^{a} ) \Big] (\square)
\bigg) \\
& \hspace{12em}
+ \sum \limits_{\substack{\square\\ \square \cap \partial \mf{U} \neq \emptyset}} \bigg(
\sum \limits_{\substack{a\\ \square^{a} \in \partial \mf{U}}}
\Big[ (j^{1} \phy)^{*} ( \iprod_{j^{1} V} \Theta_{L_{d}}^{a} ) \Big] (\square)
\bigg)
\\
\label{eq:vi_infinite_discrete_jet_space_cartan_form_6}
&= \sum \limits_{\substack{\square\\ \square \cap \mrm{int} \, \mf{U} \neq \emptyset}} \ext \mcal{L}_{\square} \cdot V + \theta_{L_{d}} (\phy) \cdot V = 0 ,
\end{align}
}%

which is similar to the continuous result (\ref{eq:classical_fields_cartan_form_55}).

\subsection{Discrete Multisymplectic Form}

The four Cartan forms (\ref{eq:vi_infinite_discrete_jet_space_cartan_form_2}) correspond to the exterior derivative of the discrete Lagrangian
\begin{align}\label{eq:vi_infinite_discrete_jet_space_multisymplectic_1}
\ext \mcal{L}_{\square}
= \sum \limits_{\substack{a\\ \square \subseteq \mf{U}}}
\dfrac{\partial \mcal{L}_{\square}}{\partial y^{a}} \big( \phy_{\square^{1}}, \phy_{\square^{2}}, \phy_{\square^{3}}, \phy_{\square^{4}} \big) \, dy_{\square^{a}}
= \sum \limits_{\substack{a\\ \square \subseteq \mf{U}}}
\Theta_{L_{d}}^{a} \big( \phy_{\square^{1}}, \phy_{\square^{2}}, \phy_{\square^{3}}, \phy_{\square^{4}} \big) ,
\end{align}

such that upon defining the discrete multisymplectic form as the exterior derivative of the discrete Cartan form,
\begin{align}\label{eq:vi_infinite_discrete_jet_space_multisymplectic_2}
\Omega_{L_{d}}^{a} = - \ext \Theta_{L_{d}}^{a} ,
\end{align}

due to $\ext^{2} \mcal{L} = 0$ we get
\begin{align}\label{eq:vi_infinite_discrete_jet_space_multisymplectic_3}
\sum \limits_{\substack{a\\ \square \subseteq \mf{U}}}
\Omega_{L_{d}}^{a} \big( \phy_{\square^{1}}, \phy_{\square^{2}}, \phy_{\square^{3}}, \phy_{\square^{4}} \big)
= \Omega_{L_{d}}^{1} (\square) + \Omega_{L_{d}}^{2} (\square) + \Omega_{L_{d}}^{3} (\square) + \Omega_{L_{d}}^{4} (\square)
= 0 .
\end{align}

For a square $\square$ in $\mf{X}$, define the projection
\begin{align}\label{eq:vi_infinite_discrete_jet_space_multisymplectic_4}
\pi_{\square} : \mf{C} (\mf{Y}) \rightarrow \jb{1}{\mf{Y}}
\end{align}

in analogy to the jet prolongation (\ref{vi_infinite_discrete_jet_space_11},\ref{vi_infinite_discrete_jet_space_12}) by
\begin{align}\label{eq:vi_infinite_discrete_jet_space_multisymplectic_5}
\pi_{\square} (\phy) \equiv \big( \square , \phy (\square^{1}) , \phy (\square^{2}) , \phy (\square^{3}) , \phy (\square^{4}) \big) ,
\end{align}

such that the forms $\pi_{\square}^{*} \Theta_{L}^{a}$ are computed as
\begin{align}\label{eq:vi_infinite_discrete_jet_space_multisymplectic_6}
(\pi_{\square}^{*} \Theta_{L_{d}}^{a}) (\phy) \cdot V = \dfrac{\partial \mcal{L}_d}{\partial y^a} \big( \phy_{\square^{1}}, \phy_{\square^{2}}, \phy_{\square^{3}}, \phy_{\square^{4}} \big) \cdot V ( \square^{a} ) ,
\end{align}

and the one-form (\ref{eq:vi_infinite_discrete_jet_space_cartan_form_4}) becomes
\begin{align}\label{eq:vi_infinite_discrete_jet_space_multisymplectic_7}
\theta_{L_{d}} =
\sum \limits_{\substack{\square\\ \square \cap \partial \mf{U} \neq \emptyset}} \bigg(
\sum \limits_{\substack{a\\ \square^{a} \in \partial \mf{U}}}
\pi_{\square}^{*} \Theta_{L_{d}}^{a} \bigg) .
\end{align}

Consider the subspace $\mf{C}_{L} \subset \mf{C} (\mf{Y})$ of the space of sections $\phy$, that solve the discrete Euler-Lagrange equations (\ref{eq:vi_infinite_discrete_jet_space_action_principle_6}).
A first variation at a solution $\phy$ of the discrete Euler-Lagrange field equations (\ref{eq:vi_infinite_discrete_jet_space_action_principle_6}) corresponds to a vector field $V \in \tb[\phy]{\mf{C}_{L}}$ such that the associated flow maps $\phy$ to other solutions of the discrete Euler-Lagrange field equations, i.e., sections $\phy \in \mf{C}_{L}$ are integral curves of $V$.
Restricting the action (\ref{eq:vi_infinite_discrete_jet_space_cartan_form_5}) to the subspace $\mf{C}_{L}$, the first sum in (\ref{eq:vi_infinite_discrete_jet_space_cartan_form_5}) becomes zero and only the one-form (\ref{eq:vi_infinite_discrete_jet_space_cartan_form_4}) is retained.
Computing the exterior derivative of the variation of the action and restricting to two first variations $V, W \in \tb[\phy]{\mf{C}_{L}}$, we obtain
\begin{align}
0= \ext^{2} \mcal{A}_{d} \cdot V \cdot W
&= \ext \theta_{L_{d}} (\phi) (V,W)
= \sum \limits_{\substack{\square\\ \square \cap \partial \mf{U} \neq \emptyset}} \bigg(
\sum \limits_{\substack{a\\ \square^{a} \in \partial \mf{U}}}
V \contr W \contr \pi_{\square}^{*} \Omega_{L_{d}}^{a} \bigg) ,
\end{align}

which is equivalent to
\begin{align}
\sum \limits_{\substack{\square\\ \square \cap \partial \mf{U} \neq \emptyset}} \bigg(
\sum \limits_{\substack{a\\ \square^{a} \in \partial \mf{U}}}
\Big[ (j^{1} \phy)^{*} ( V \contr W \contr \Omega_{L_{d}}^{a} \Big] (\square)
\bigg) = 0 .
\end{align}

This is the discrete analogue to the multisymplectic form formula (\ref{eq:classical_multisymplectic_form_formula}).

\subsection{Discrete Noether Theorem}\label{sec:vi_infinite_noether_theorem}

We restrict our treatment to the case of a scalar field theory, one spatial dimension, and vertical transformations, but at least the first two restrictions are easily lifted \cite{MarsdenPatrick:1998}.
Consider a one-parameter group of
\begin{align}\label{eq:vi_infinite_noether_1}
\phy^{\eps}_{i,k} &= \eta^{\eps} \circ \phy_{i,k} &
& \text{such that} &
\phy^{0}_{i,k} &= \phy_{i,k} . & &&
\end{align}

The infinitesimal generator of the transformation $\phy^{\eps}_{i,k}$ is
\begin{align}\label{eq:vi_infinite_noether_2}
X_{i,k} &= \dfrac{d}{d \eps} \phy^{\eps}_{i,k} \bigg\vert_{\eps = 0} &
& \text{or in abstract notation} &
X_{\square^{a}} (\phy_{\square^{a}}) &= \dfrac{d}{d \eps} \phy^{\eps}_{\square^{a}} \bigg\vert_{\eps = 0} . &
\end{align}

The discrete Lagrangian has a symmetry if it is invariant under this transformation
\begin{align}\label{eq:vi_infinite_noether_3}
\mathcal{L}_{d} \big( \phy_{\square^{1}}^{\eps}, \phy_{\square^{2}}^{\eps}, \phy_{\square^{3}}^{\eps}, \phy_{\square^{4}}^{\eps} \big)
&= \mathcal{L}_{d} \big( \phy_{\square^{1}}, \phy_{\square^{2}}, \phy_{\square^{3}}, \phy_{\square^{4}} \big) &
& \text{for all $\eps$} . & &&
\end{align}

This is equivalent to
\begin{align}\label{eq:vi_infinite_noether_4}
\ext \mathcal{L}_{\square} \cdot X
&= \dfrac{d}{d \eps} \bigg\vert_{\eps = 0} \mathcal{L}_{d} \big( \phy_{\square^{1}}^{\eps}, \phy_{\square^{2}}^{\eps}, \phy_{\square^{3}}^{\eps}, \phy_{\square^{4}}^{\eps} \big) \\
&= \sum \limits_{\substack{a\\ \square \subseteq \mf{U}}} \bigg[
\dfrac{\partial \mcal{L}_d}{\partial y^a} \big( \phy_{\square^{1}}, \phy_{\square^{2}}, \phy_{\square^{3}}, \phy_{\square^{4}} \big) \cdot X_{\square^{a}} (\phy_{\square^{a}}) \bigg] = 0 ,
\end{align}

or explicitly in grid coordinates,
\begin{align}\label{eq:vi_infinite_noether_5}
\ext \mathcal{L}_{\square} \cdot X
\nonumber
&= \dfrac{\partial \mcal{L}_d}{\partial y^1} \Big( \phy_{i,  k  }, \phy_{i+1,k  }, \phy_{i+1,k+1}, \phy_{i,  k+1} \Big) \cdot X_{i,k}
\\
\nonumber
& \hspace{4em}
+ \dfrac{\partial \mcal{L}_d}{\partial y^2} \Big( \phy_{i,  k  }, \phy_{i+1,k  }, \phy_{i+1,k+1}, \phy_{i,  k+1} \Big) \cdot X_{i+1,k}
\\
\nonumber
& \hspace{8em}
+ \dfrac{\partial \mcal{L}_d}{\partial y^3} \Big( \phy_{i,  k  }, \phy_{i+1,k  }, \phy_{i+1,k+1}, \phy_{i,  k+1} \Big) \cdot X_{i+1,k+1}
\\
& \hspace{12em}
+ \dfrac{\partial \mcal{L}_d}{\partial y^4} \Big( \phy_{i,  k  }, \phy_{i+1,k  }, \phy_{i+1,k+1}, \phy_{i,  k+1} \Big) \cdot X_{i,k+1}
= 0 .
\end{align}

Since the $\phy$ are solutions of the discrete Euler-Lagrange field equations (\ref{eq:vi_infinite_discrete_jet_space_action_principle_6}), the generating vector field $X$ is a first variation, i.e., $X \in \tb[\phy]{\mf{C}_{L}}$.
This means that the sum in (\ref{eq:vi_infinite_discrete_jet_space_cartan_form_6}) vanishes, and we obtain
\begin{align}\label{eq:vi_infinite_noether_6}
\theta_{L_{d}} (\phy) \cdot X \equiv
\sum \limits_{\substack{\square\\ \square \cap \partial \mf{U} \neq \emptyset}} \bigg(
\sum \limits_{\substack{a\\ \square^{a} \in \partial \mf{U}}}
\Big[ (j^{1} \phy)^{*} ( \iprod_{j^{1} X} \Theta_{L_{d}}^{a} ) \Big] (\square)
\bigg) = 0,
\end{align}

or explicitly in grid coordinates,
\begin{align}\label{eq:vi_infinite_noether_7}
0
\nonumber
= \sum \limits_{\substack{i\\ (i,k) \in \partial \mf{U}}} \bigg[ &
\dfrac{\partial \mcal{L}_d}{\partial y^1} \Big( \phy_{i,   1}, \phy_{i+1, 1}, \phy_{i+1, 2}, \phy_{i,   2} \Big) \cdot X_{i,1} \\
\nonumber
& \hspace{2em}
+ \dfrac{\partial \mcal{L}_d}{\partial y^2} \Big( \phy_{i-1, 1}, \phy_{i,   1}, \phy_{i,   2}, \phy_{i-1, 2} \Big) \cdot X_{i,1} \\
\nonumber
& \hspace{6em}
+ \dfrac{\partial \mcal{L}_d}{\partial y^3} \Big( \phy_{i-1, n_{t}-1}, \phy_{i,   n_{t}-1}, \phy_{i,   n_{t}}, \phy_{i-1, n_{t}} \Big) \cdot X_{i,n_{t}} \\
& \hspace{10em}
+ \dfrac{\partial \mcal{L}_d}{\partial y^4} \Big( \phy_{i,   n_{t}-1}, \phy_{i+1, n_{t}-1}, \phy_{i+1, n_{t}}, \phy_{i,   n_{t}} \Big) \cdot X_{i,n_{t}}
\bigg] .
\end{align}

\begin{figure}[tb]
\centering
\begin{tikzpicture}[scale=0.7]%

\tikzstyle{every node}=[font=\small]

\draw [style=help lines, step=2]	(-9,-3)	grid		(+9,+5);

\draw [dashed, line width=0.4mm]	(-8,-2)	--	(-6,-2);
\draw [solid,  line width=0.4mm]	(-6,-2)	--	(-4,-2);
\draw [dashed, line width=0.4mm]	(-4,-2)	--	(-2,-2);
\draw [solid,  line width=0.4mm]	(-2,-2)	--	(+2,-2);
\draw [dashed, line width=0.4mm]	(+2,-2)	--	(+4,-2);
\draw [solid,  line width=0.4mm]	(+4,-2)	--	(+6,-2);
\draw [dashed, line width=0.4mm]	(+6,-2)	--	(+8,-2);

\draw [dashed, line width=0.4mm]	(-8,+4)	--	(-6,+4);
\draw [solid,  line width=0.4mm]	(-6,+4)	--	(-4,+4);
\draw [dashed, line width=0.4mm]	(-4,+4)	--	(-2,+4);
\draw [solid,  line width=0.4mm]	(-2,+4)	--	(+2,+4);
\draw [dashed, line width=0.4mm]	(+2,+4)	--	(+4,+4);
\draw [solid,  line width=0.4mm]	(+4,+4)	--	(+6,+4);
\draw [dashed, line width=0.4mm]	(+6,+4)	--	(+8,+4);

\draw [solid,  line width=0.4mm]	(-6,+4)	--	(-6,+2);
\draw [dashed, line width=0.4mm]	(-6,+2)	--	(-6, 0);
\draw [solid,  line width=0.4mm]	(-6, 0)	--	(-6,-2);

\draw [solid,  line width=0.4mm]	(+6,+4)	--	(+6,+2);
\draw [dashed, line width=0.4mm]	(+6,+2)	--	(+6, 0);
\draw [solid,  line width=0.4mm]	(+6, 0)	--	(+6,-2);

\fill [fill=blue!40!white, fill opacity=0.5]	(-2,-2)	rectangle	(+0,+0);
\fill [fill=blue!40!white, fill opacity=0.5]	(-2,+2)	rectangle	(+0,+4);
\fill [fill=blue!40!white, fill opacity=0.5]	(+2,+2)	rectangle	(+0,+4);
\fill [fill=blue!40!white, fill opacity=0.5]	(+2,-2)	rectangle	(+0,+0);

\filldraw [color=black]	(0, -2) circle (.1);
\filldraw [color=black]	(0, +4) circle (.1);

\draw [black]	(+9, -2)	node[anchor=west]	{$k=1$};
\draw [black]	(+9, +4)	node[anchor=west]	{$k=n_t$};

\draw [black]	(-8, -3)	node[anchor=north]	{$n_{x}\vphantom{1}$};
\draw [black]	(-6, -3)	node[anchor=north]	{$1$};
\draw [black]	(-2, -3)	node[anchor=north]	{$i-1$};
\draw [black]	( 0, -3)	node[anchor=north]	{$i\vphantom{1}$};
\draw [black]	(+2, -3)	node[anchor=north]	{$i+1$};
\draw [black]	(+6, -3)	node[anchor=north]	{$n_{x}\vphantom{1}$};
\draw [black]	(+8, -3)	node[anchor=north]	{$1$};

\draw [black]	(-1, -1)	node	{(b)};
\draw [black]	(+1, -1)	node	{(a)};
\draw [black]	(-1, +3)	node	{(c)};
\draw [black]	(+1, +3)	node	{(d)};

\end{tikzpicture}

\caption{Boundary terms in a spatially periodic integration domain.}
\label{fig:vi_infinite_noether_grid}
\end{figure}
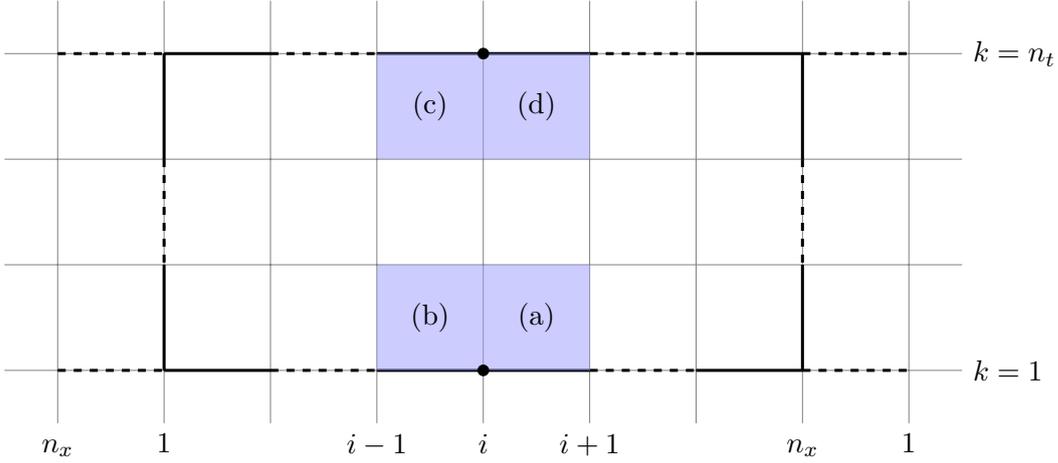

The various contributions arise as depicted in figure \ref{fig:vi_infinite_noether_grid}.
If we fix the spatial index $i$, there are four cells that touch the boundary, two at $(i,1)$ and two at $(i,n_{t})$, respectively.
The first contribution arises from the cell marked (a). The point $(i,1)$ corresponds to $\square^{1}$, such that the derivative of the Lagrangian is computed with respect to $y^{1} = y(\square^{1})$. The other contributions follow in the same way.
In principle, there are also contributions to (\ref{eq:vi_infinite_noether_6}) that arise from the spatial boundary at $i=1$ and $i=n_{x}$. One either has to account for those or select boundary conditions, that automatically take care of these contributions. We shall always use periodic boundary conditions (as depicted), such that practically there is no boundary in the spatial dimension.

With the help of the discrete symmetry condition (\ref{eq:vi_infinite_noether_5}), we can replace the first two lines in (\ref{eq:vi_infinite_noether_7}), such that
\begin{multline}\label{eq:vi_infinite_noether_8}
\sum \limits_{\substack{i\\ (i,k) \in \partial \mf{U}}} \bigg[
\dfrac{\partial \mcal{L}_d}{\partial y^3} \Big( \phy_{i-1, n_{t}-1}, \phy_{i,   n_{t}-1}, \phy_{i,   n_{t}}, \phy_{i-1, n_{t}} \Big)
+ \dfrac{\partial \mcal{L}_d}{\partial y^4} \Big( \phy_{i,   n_{t}-1}, \phy_{i+1, n_{t}-1}, \phy_{i+1, n_{t}}, \phy_{i,   n_{t}} \Big)
\bigg] \cdot X_{i,n_{t}} =
\\
= \sum \limits_{\substack{i\\ (i,k) \in \partial \mf{U}}} \bigg[
\dfrac{\partial \mcal{L}_d}{\partial y^3} \Big( \phy_{i-1, 1}, \phy_{i,   1}, \phy_{i,   2}, \phy_{i-1, 2} \Big)
+ \dfrac{\partial \mcal{L}_d}{\partial y^4} \Big( \phy_{i,   1}, \phy_{i+1, 1}, \phy_{i+1, 2}, \phy_{i,   2} \Big)
\bigg] \cdot X_{i,2} .
\end{multline}

This is a conservation law, and as the number of timesteps $n_{t}$ is arbitrary, can be written as
\begin{align}\label{eq:vi_infinite_noether_9}
\sum \limits_{i=1}^{n_{x}-1} \bigg[
\dfrac{\partial \mcal{L}_d}{\partial y^3} \Big( \phy_{i-1, k-1}, \phy_{i,   k-1}, \phy_{i,   k}, \phy_{i-1, k} \Big)
+ \dfrac{\partial \mcal{L}_d}{\partial y^4} \Big( \phy_{i,   k-1}, \phy_{i+1, k-1}, \phy_{i+1, k}, \phy_{i,   k} \Big)
\bigg] \cdot X_{i,k} = const .
\end{align}

The structure of this conservation law implies that if the continuous Lagrangian has a vertical symmetry, any consistent discretisation of the Lagrangian will lead to a discrete conservation law corresponding to that symmetry.

\subsection{Discrete Momentum Maps}

Geometrically, a conserved quantity is described by a momentum map 
\begin{align}\label{eq:vi_infinite_momentum_maps_1}
J^{a} &\equiv \iprod_{j^{1} X} \Theta_{L_{d}}^{a} &
& \text{with} &
X_{\square^{a}} (\phy_{\square^{a}}) &= \dfrac{d}{d \eps} \phy^{\eps}_{\square^{a}} \bigg\vert_{\eps = 0} , &
\end{align}

such that (\ref{eq:vi_infinite_noether_6}) takes the form
\begin{align}\label{eq:vi_infinite_momentum_maps_2}
\sum \limits_{\substack{\square\\ \square \cap \partial \mf{U} \neq \emptyset}} \bigg(
\sum \limits_{\substack{a\\ \square^{a} \in \partial \mf{U}}}
(j^{1} \phy)^{*} J^{a} (\square)
\bigg) = 0,
\end{align}

or explicitly in grid coordinates,
\begin{multline}\label{eq:vi_infinite_momentum_maps_3}
\sum \limits_{i=1}^{n_{x}-1} \bigg[
J^{1} \big( \phy_{i,   1}, \phy_{i+1, 1}, \phy_{i+1, 2}, \phy_{i,   2} \big)
+ J^{2} \big( \phy_{i,   1}, \phy_{i+1, 1}, \phy_{i+1, 2}, \phy_{i,   2} \big)
\\
+ J^{3} \big( \phy_{i,   n_{t}-1}, \phy_{i+1, n_{t}-1}, \phy_{i+1, n_{t}}, \phy_{i,   n_{t}} \big)
+ J^{4} \big( \phy_{i,   n_{t}-1}, \phy_{i+1, n_{t}-1}, \phy_{i+1, n_{t}}, \phy_{i,   n_{t}} \big)
\bigg] = 0.
\end{multline}

The invariance of the Lagrangian (\ref{eq:vi_infinite_noether_3}) and (\ref{eq:vi_infinite_discrete_jet_space_multisymplectic_1}) implies
\begin{multline}\label{eq:vi_infinite_momentum_maps_4}
J^{1} \big( \phy_{i,   k}, \phy_{i+1, k}, \phy_{i+1, k+1}, \phy_{i,   k+1} \big)
+ J^{2} \big( \phy_{i,   k}, \phy_{i+1, k}, \phy_{i+1, k+1}, \phy_{i,   k+1} \big) \\
+ J^{3} \big( \phy_{i,   k}, \phy_{i+1, k}, \phy_{i+1, k+1}, \phy_{i,   k+1} \big)
+ J^{4} \big( \phy_{i,   k}, \phy_{i+1, k}, \phy_{i+1, k+1}, \phy_{i,   k+1} \big) = 0 ,
\end{multline}

so that we can replace the last line of (\ref{eq:vi_infinite_momentum_maps_3}) to get
\begin{multline}\label{eq:vi_infinite_momentum_maps_5}
\sum \limits_{i=1}^{n_{x}-1} \bigg[
J^{3} \big( \phy_{i,   n_{t}-1}, \phy_{i+1, n_{t}-1}, \phy_{i+1, n_{t}}, \phy_{i,   n_{t}} \big)
+ J^{4} \big( \phy_{i,   n_{t}-1}, \phy_{i+1, n_{t}-1}, \phy_{i+1, n_{t}}, \phy_{i,   n_{t}} \big)
\bigg]
=
\\
=
\sum \limits_{i=1}^{n_{x}-1} \bigg[
J^{3} \big( \phy_{i,   1}, \phy_{i+1, 1}, \phy_{i+1, 2}, \phy_{i,   2} \big)
+ J^{4} \big( \phy_{i,   1}, \phy_{i+1, 1}, \phy_{i+1, 2}, \phy_{i,   2} \big)
\bigg]
.
\end{multline}

As this is true for any $n_{t}$, we can infer the general expression
\begin{align}\label{eq:vi_infinite_momentum_maps_6}
\sum \limits_{i=1}^{n_{x}-1} \bigg[
J^{3} \big( \phy_{i,   k}, \phy_{i+1, k}, \phy_{i+1, k+1}, \phy_{i,   k+1} \big)
+ J^{4} \big( \phy_{i,   k}, \phy_{i+1, k}, \phy_{i+1, k+1}, \phy_{i,   k+1} \big)
\bigg] = const
.
\end{align}

This is equivalent to the conservation law (\ref{eq:vi_infinite_noether_9}).

\section{Example: The Advection Equation}

In this section, we apply the theory of the previous chapter to an interesting and instructive example, namely the advection equation.
Its structure is similar to that of the Vlasov equation which we will study in chapter \ref{ch:kinetic_theory}.
The linear advection equation for a field $u(t,x)$ (in the following referred to as \emph{density}) reads
\begin{align}\label{eq:vi_advection_equation}
\partial_{t} u + c \, \partial_{x} u = 0 ,
\end{align}

where the velocity, $c$, is a constant parameter.
For initial conditions $u (t = 0, x) = u_{0} (x)$, the analytic solution is
\begin{align}
u (t,x) = u_{0} (x - ct) .
\end{align}

The advection equation does not have a natural variational formulation, but we can apply the theory from section \ref{ch:classical_extended_lagrangians} to obtain an extended Lagrangian.

\subsection{Extended Lagrangian}

The extended Lagrangian for the advection equation is obtained by multiplying (\ref{eq:vi_advection_equation}) with the auxiliary variable $v(t,x)$.
The solution vector of the extended system is denoted $w = (u,v)$, such that the Lagrangian can be written as
\begin{align}\label{eq:vi_advection_lagrangian}
\mcal{L} (j^{1} w) = v \big( u_{t} + c u_{x} \big) \, \omega ,
\end{align}

where $\omega = dt \wedge dx$ is the volume form, so that the corresponding action is $\mcal{A} = \int \mcal{L}$.
The variational derivative of the action with respect to the Ibragimov multiplier $v$ retains the advection equation
\begin{align}
\dfrac{\delta \mcal{A}}{\delta v} = + u_{t} + c u_{x} = 0 .
\end{align}

The variation with respect to the original variable $u$ yields the adjoint equation
\begin{align}
\dfrac{\delta \mcal{A}}{\delta u} = - v_{t} - c v_{x} = 0 .
\end{align}

It is immediately observed that the adjoint equation has the same solution as the original equation, such that if $u$ is a solution of the advection equation, then $w = (u,u)$ solves the Euler-Lagrange equations of the extended Lagrangian (\ref{eq:vi_advection_lagrangian}).
Here and in the remaining part of the section, $w$ denotes the combined fields $(u,v)$.

\subsection{Variational Integrator}

We discretise the Lagrangian (\ref{eq:vi_advection_lagrangian}) on a grid cell as depicted in figure \ref{fig:vi_infinite_gridbox} by
\begin{multline}\label{eq:vi_advection_integrator_1}
\mcal{L}_{d} ( w^{1}, w^{2}, w^{3}, w^{4} ) = h_{t} \, h_{x} \, \dfrac{1}{4} \Big( v^{1} + v^{2} + v^{3} + v^{4} \Big) \times \\
\times \bigg[ \dfrac{1}{2} \bigg( \dfrac{u^{4} - u^{1}}{h_{t}} + \dfrac{u^{3} - u^{2}}{h_{t}} \bigg) + \dfrac{c}{2} \bigg( \dfrac{u^{2} - u^{1}}{h_{x}} + \dfrac{u^{3} - u^{4}}{h_{x}} \bigg) \bigg]
\end{multline}

The discrete Euler-Lagrange field equations (\ref{eq:vi_infinite_delfeqs})
\begin{align}\label{eq:vi_advection_integrator_2}
0
\nonumber
&= \dfrac{\partial \mathcal{L}_d}{\partial g^1} \Big( w_{i,  k  }, w_{i+1,k  }, w_{i+1,k+1}, w_{i,  k+1} \Big)
+ \dfrac{\partial \mathcal{L}_d}{\partial g^2} \Big( w_{i-1,k  }, w_{i,  k  }, w_{i,  k+1}, w_{i-1,k+1} \Big) \\
&+ \dfrac{\partial \mathcal{L}_d}{\partial g^3} \Big( w_{i-1,k-1}, w_{i,  k-1}, w_{i,  k  }, w_{i-1,k  } \Big)
+ \dfrac{\partial \mathcal{L}_d}{\partial g^4} \Big( w_{i,  k-1}, w_{i+1,k-1}, w_{i+1,k  }, w_{i,  k  } \Big)
\end{align}

are computed as
\begin{align}\label{eq:vi_advection_integrator_3}
0
\nonumber
=& \dfrac{1}{4} \bigg[ \dfrac{u_{i+1, k+1} - u_{i+1, k-1}}{2 h_{t}} + 2 \, \dfrac{u_{i, k+1} - u_{i, k-1}}{2 h_{t}} + \dfrac{u_{i-1, k+1} - u_{i-1, k-1}}{2 h_{t}} \bigg] \\
+& \dfrac{c}{4} \bigg[ \dfrac{u_{i+1, k+1} - u_{i-1, k+1}}{2 h_{x}} + 2 \, \dfrac{u_{i+1, k} - u_{i-1, k}}{2 h_{x}} + \dfrac{u_{i+1, k-1} - u_{i-1, k-1}}{2 h_{x}} \bigg] .
\end{align}

As in the continuous case, the discrete adjoint equation has the exact same form as the discrete advection equation.
Again, we find the same kind of averaging we have already observed in the example of the wave equation (\ref{eq:vi_infinite_wave_equation_discrete}).

\subsection{Continuous Conservation Laws}

We will shortly prove some conservation laws of the advection equation, namely the conservation of the average density, the $L^{1}$ and $L^{2}$ norms, the total linear momentum and kinetic energy.
We consider vertical transformations with generators of the form
\begin{align}\label{eq:vi_advection_continuous_generator}
X &= X^{u} \, \dfrac{\partial}{\partial u} + X^{v} \, \dfrac{\partial}{\partial v} &
& \text{with} &
X^{a} &= X^{a} (u, v) , &
\end{align}

as discussed in section \ref{sec:classical_noether_theorem}.
The jet prolongation (\ref{eq:noether_fields_5}) of such vector fields are
\begin{align}\label{eq:vi_advection_continuous_generator_prolongation}
j^{1} X
= X^{u} \, \dfrac{\partial}{\partial u} + X^{v} \, \dfrac{\partial}{\partial v}
+ \bigg[ \dfrac{\partial X^{u}}{\partial u} \dfrac{\partial u}{\partial x^{\mu}} + \dfrac{\partial X^{v}}{\partial u} \dfrac{\partial u}{\partial x^{\mu}} \bigg] \, \dfrac{\partial}{\partial u_{\mu}}
+ \bigg[ \dfrac{\partial X^{u}}{\partial v} \dfrac{\partial v}{\partial x^{\mu}} + \dfrac{\partial X^{v}}{\partial v} \dfrac{\partial v}{\partial x^{\mu}} \bigg] \, \dfrac{\partial}{\partial v_{\mu}} ,
\end{align}

such that the invariance condition (\ref{eq:noether_fields_8}) becomes
\begin{align}\label{eq:vi_advection_continuous_invariance}
j^{1} X (L)
&= X^{v} \, \big( u_{t} + c \, u_{x} \big)
+ v  \, \big( u_{t} + c \, u_{x} \big) \bigg( \dfrac{\partial X^{u}}{\partial u} + \dfrac{\partial X^{v}}{\partial u} \bigg) .
\end{align}

Conservation laws (\ref{eq:noether_fields_12}) take the form
\begin{align}\label{eq:vi_advection_continuous_conservation_law}
\dfrac{d}{d t} \int \bigg[ \dfrac{\partial L}{\partial u_{t}} \big( j^{1} w \big) \, X^{u} + \dfrac{\partial L}{\partial v_{t}} \big( j^{1} w \big) \, X^{v} \bigg] \, dx
= \dfrac{d}{d t} \int v \, X^{u} \, dx .
\end{align}

We can use these expressions directly to test for several possible symmetries.

\subsubsection{Conservation of the Average Density}

Consider a vertical transformation generated by $X = (X^{u}, X^{v})$ with
\begin{align}\label{eq:vi_advection_noether_continuous_1}
X^{u} = 1,
\hspace{2em}
X^{v} = 0 .
\end{align}

The Lagrangian is trivially invariant under this transformation
\begin{align}\label{eq:vi_advection_noether_continuous_2}
j^{1} X (L) = 0 .
\end{align}

The corresponding conservation law (\ref{eq:vi_advection_continuous_conservation_law}) is
\begin{align}\label{eq:vi_advection_noether_continuous_3}
\dfrac{d}{d t} \int v \, dx = 0 .
\end{align}

If $u$ is a solution of the advection equation, the pair $(u,v) = (u,u)$ solves the extended system.
We can therefore reduce the conserved quantity in (\ref{eq:vi_advection_noether_continuous_3}) by specialising it to $v = u$.
This gives the conservation of the average density, namely
\begin{align}\label{eq:vi_advection_noether_continuous_4}
\dfrac{d}{d t} \int u \, dx = 0 .
\end{align}

If, in addition, the advection equation (\ref{eq:vi_advection_equation}) is provided with positive initial conditions, $u_{0} (x) \geq 0$, the positivity is preserved at later times, and thus (\ref{eq:vi_advection_noether_continuous_4}) turns into a conservation law for the $L^{1}$ norm of $u$.

\subsubsection{Conservation of the $L^{2}$ Norm}

Consider a different vertical transformation generated by $X = (X^{u}, X^{v})$ with
\begin{align}\label{eq:vi_advection_noether_continuous_5}
X^{u} = u,
\hspace{2em}
X^{v} = -v .
\end{align}

The Lagrangian density is invariant under this transformation as well
\begin{align}\label{eq:vi_advection_noether_continuous_6}
j^{1} X (L) = - v \big( u_{t} + c \, u_{x} \big) + v \big( u_{t} + c \, u_{x} \big) = 0 .
\end{align}

The corresponding conservation law (\ref{eq:vi_advection_continuous_conservation_law}) is
\begin{align}\label{eq:vi_advection_noether_continuous_7}
\dfrac{d}{d t} \int v u \, dx = 0 .
\end{align}

Upon identifying $v$ with $u$, this gives the conservation of the $L^{2}$ norm of $u$
\begin{align}\label{eq:vi_advection_noether_continuous_8}
\dfrac{d}{d t} \int u^{2} \, dx = 0 .
\end{align}

\subsubsection{Conservation of Linear Momentum and Kinetic Energy}

Lastly, consider the following vertical transformation generated by $X = (X^{u}, X^{v})$ with
\begin{align}\label{eq:vi_advection_noether_continuous_9}
X^{u} = c,
\hspace{2em}
X^{v} = 0 .
\end{align}

The Lagrangian density is trivially invariant also under this transformation
\begin{align}\label{eq:vi_advection_noether_continuous_10}
j^{1} X (L) = 0 .
\end{align}

The corresponding conservation law is
\begin{align}\label{eq:vi_advection_noether_continuous_11}
\dfrac{d}{d t} \int c \, v \, dx = 0 .
\end{align}

Upon identifying $v$ with $u$, this gives the conservation of linear momentum
\begin{align}\label{eq:vi_advection_noether_continuous_12}
\dfrac{d}{d t} \int c \, u \, dx = 0 .
\end{align}

Conservation of kinetic energy follows exactly the same way by choosing $X^{u} = \tfrac{1}{2} c^{2}$, i.e.,
\begin{align}\label{eq:vi_advection_noether_continuous_13}
\dfrac{d}{d t} \int \tfrac{1}{2} \, c^{2} \, v \, dx
= \dfrac{d}{d t} \int \tfrac{1}{2} \, c^{2} \, u \, dx
= 0 .
\end{align}

We have therefore proved conservation of the most important quantities related to the advection equation.

\subsection{Discrete Conservation Laws}

The discrete generator $X_{i,k} = (X_{i,k}^{u}, X_{i,k}^{v})$ of the transformation given by (\ref{eq:vi_advection_noether_continuous_1}) is
\begin{align}\label{eq:vi_advection_noether_discrete_1}
X_{i,k}^{u} = 1,
\hspace{2em}
X_{i,k}^{v} = 0 .
\end{align}

The discrete Lagrangian (\ref{eq:vi_advection_integrator_1}) is invariant under this transformation
\begin{align}\label{eq:vi_advection_noether_discrete_2}
\ext \mathcal{L}_{\square} \cdot X
\nonumber
&= \dfrac{\partial \mcal{L}_d}{\partial u^1} \Big( w_{i,  k  }, w_{i+1,k  }, w_{i+1,k+1}, w_{i,  k+1} \Big) 
\\
\nonumber
& \hspace{3em}
+ \dfrac{\partial \mcal{L}_d}{\partial u^2} \Big( w_{i,  k  }, w_{i+1,k  }, w_{i+1,k+1}, w_{i,  k+1} \Big) 
\\
\nonumber
& \hspace{6em}
+ \dfrac{\partial \mcal{L}_d}{\partial u^3} \Big( w_{i,  k  }, w_{i+1,k  }, w_{i+1,k+1}, w_{i,  k+1} \Big) 
\\
\nonumber
& \hspace{9em}
+ \dfrac{\partial \mcal{L}_d}{\partial u^4} \Big( w_{i,  k  }, w_{i+1,k  }, w_{i+1,k+1}, w_{i,  k+1} \Big) 
\\
\nonumber
&= \dfrac{h_{t} \, h_{x}}{8} \, \Big[ v_{i,k} + v_{i+1,k} + v_{i+1,k+1} + v_{i,k+1} \Big] \Big[ - \dfrac{1}{h_{t}} - \dfrac{c}{h_{x}} \Big] \, X^{u}_{i,k}
\\
\nonumber
& \hspace{3em}
+ \dfrac{h_{t} \, h_{x}}{8} \, \Big[ v_{i,k} + v_{i+1,k} + v_{i+1,k+1} + v_{i,k+1} \Big] \Big[ - \dfrac{1}{h_{t}} + \dfrac{c}{h_{x}} \Big] \, X^{u}_{i+1,k}
\\
\nonumber
& \hspace{6em}
+ \dfrac{h_{t} \, h_{x}}{8} \, \Big[ v_{i,k} + v_{i+1,k} + v_{i+1,k+1} + v_{i,k+1} \Big] \Big[ + \dfrac{1}{h_{t}} + \dfrac{c}{h_{x}} \Big] \, X^{u}_{i+1,k+1}
\\
\nonumber
& \hspace{9em}
+ \dfrac{h_{t} \, h_{x}}{8} \, \Big[ v_{i,k} + v_{i+1,k} + v_{i+1,k+1} + v_{i,k+1} \Big] \Big[ + \dfrac{1}{h_{t}} - \dfrac{c}{h_{x}} \Big] \, X^{u}_{i,k+1}
\\
&= 0 .
\end{align}

The corresponding conservation law (\ref{eq:vi_infinite_noether_9}) is
\begin{align}\label{eq:vi_advection_noether_discrete_3}
\nonumber
& \sum \limits_{i=1}^{n_{x}-1} \bigg[
\dfrac{\partial \mcal{L}_d}{\partial u^3} \Big( w_{i-1, k-1}, w_{i,   k-1}, w_{i,   k}, w_{i-1, k} \Big)
+ \dfrac{\partial \mcal{L}_d}{\partial u^4} \Big( w_{i,   k-1}, w_{i+1, k-1}, w_{i+1, k}, w_{i,   k} \Big)
\bigg] \cdot X^{u}_{i,k} \\
\nonumber
& \hspace{4em}
= h_{x} \, \sum \limits_{i=1}^{n_{x}-1} \dfrac{1}{4} \Big[ v_{i,k} + v_{i+1,k} + v_{i+1,k+1} + v_{i,k+1} \Big] \\
& \hspace{4em}
= h_{x} \, \sum \limits_{i=1}^{n_{x}-1} \dfrac{1}{4} \Big[ u_{i,k} + u_{i+1,k} + u_{i+1,k+1} + u_{i,k+1} \Big]
,
\end{align}

where the last equality arises from identifying $v$ with $u$.
The conservation of momentum and energy follows along the same lines with $X_{i,k}^{u} = c$ and $X_{i,k}^{u} = \tfrac{1}{2} c^{2}$, respectively. 

\chapter{Charged Particle Motion}

This section addresses a reduced description of the motion of charged particles in a plasma, the so called guiding centre dynamics.
It can be seen as a limit of gyrokinetic theory, the predominant model used in plasma physical particle-in-cell codes, that is valid when the magnetic field is very strong, such that the gyration orbit of the particle is very small, or when the electromagnetic field is almost spatially uniform, such that it varies only very little along the gyration orbit.
It is also closely related to drift kinetic theory which is a reduced model of kinetic theory (see next chapter).

After a short summary of the Lagrangian formulation for this set of dynamical equations, a set of variational integrators based on different quadrature rules is derived. This integrators is adapted to the reduced dynamics in the poloidal plane similar to \cite{QinGuanTang:2009} as well as to the dynamics in full tokamak geometry. Several higher order methods are derived by composition methods. The advantages of the variational discretisations compared to widely used Runge-Kutta schemes are demonstrated.

\section{Guiding Centre Dynamics}

In a magnetic field, charged particles move along a helix. This motion can be decomposed into the gyration about a magnetic field line, and the motion of the centre of the gyration (guiding centre) along the field line. Mathematically, this amounts to a coordinate transformation from spatial coordinates $x$ and the corresponding velocities $\dot{x}$ to guiding centre coordinates $(X, \Theta, u, \mu)$, where $X$ is the position of the guiding centre, $\Theta$ is the angle of the gyration (gyrophase), $\mu$ the magnetic moment, and $u = \dot{x} \cdot b$ is the velocity along the magnetic field lines (parallel velocity).

\citeauthor{Littlejohn:1983} \cite{Littlejohn:1983} was the first to devise a variational principle for the guiding centre motion and thereby find simple proofs for conservation of energy and angular momentum. His Lagrangian reads
\begin{align}\label{eq:particles_guiding_centre_lagrangian}
L = A^{*} \cdot \dot{X} + \mu \dot{\Theta} - H
\end{align}
with the Hamiltonian
\begin{align}\label{eq:particles_guiding_centre_hamiltonian}
H = \dfrac{1}{2} u^{2} + \mu B + \phi
\end{align}
and the so called ``modified vector potential'' (first discovered by \citeauthor{MorozovSolovev:1966} \cite{MorozovSolovev:1966})
\begin{align}\label{eq:particles_modified_vector_potential}
A^{*} = A + ub .
\end{align}

Here, $B$ is the magnetic field strength, $b$ its unit vector, and $\phi$ is the electrostatic potential. Units are chosen such that $e = m = c = 1$, with charge $e$, particle mass $m$, and speed of light $c$.

The Lagrangian (\ref{eq:particles_guiding_centre_lagrangian}) is regarded as a function of the guiding centre variables $(X, \Theta, u, \mu)$ and their time derivatives $(\dot{X}, \dot{\Theta}, \dot{u}, \dot{\mu})$.
For the spatial components, the Euler-Lagrange equations are
\begin{align}
\dfrac{d}{dt} \bigg( \dfrac{\partial L}{\partial \dot{X}^{i}} \bigg) - \dfrac{\partial L}{\partial X^{i}} = 0
\end{align}

which explicitly amounts to
\begin{align}\label{eq:particles_guiding_centre_equations_of_motion_spatial_1}
\dfrac{d}{dt} A^{*}_{i} - A^{*}_{j,i} \dot{X}^{j} + \mu B_{,i} + \phi_{,i} = 0 .
\end{align}

Computing the time derivative, this becomes
\begin{align}\label{eq:particles_guiding_centre_equations_of_motion_spatial_2}
\big( A^{*}_{i,j} - A^{*}_{j,i} \big) \, \dot{X}^{j} + \dot{u} b_{i} + \mu B_{,i} + \phi_{,i} = 0 .
\end{align}

From the equation of the gyrophase
\begin{align}\label{eq:particles_guiding_centre_equations_of_motion_gyrophase}
\dfrac{d}{dt} \bigg( \dfrac{\partial L}{\partial \dot{\Theta}} \bigg) = 0
\hspace{3em}
\text{or}
\hspace{3em}
\dot{\mu} = 0
\end{align}

the conservation of the magnetic moment $\mu$ is obtained.
The Euler-Lagrange equation of the parallel velocity is just
\begin{align}\label{eq:particles_guiding_centre_equations_of_motion_parallel_velocity}
\dfrac{\partial L}{\partial u} = 0
\hspace{3em}
\text{or}
\hspace{3em}
u = b \cdot \dot{X} ,
\end{align}

i.e. the definition of the parallel velocity.
And from the equation of $\mu$
\begin{align}\label{eq:particles_guiding_centre_equations_of_motionetic_moment}
\dfrac{\partial L}{\partial \mu} = 0
\hspace{3em}
\text{or}
\hspace{3em}
\dot{\Theta} = B
\end{align}

we obtain, upon restoration of physical units, the definition of the gyro frequency
\begin{align}
\dot{\Theta} = \omega = \dfrac{e B}{m} .
\end{align}

If the variation of the background electromagnetic fields is small along the radius of the gyration, the particle's motion can be approximated by the motion of just the guiding centre, averaging over the gyrophase. The corresponding reduced Lagrangian is
\begin{align}\label{eq:particles_guiding_centre_reduced_lagrangian}
L &= A^{*} \cdot \dot{X} - H , &
H &= \dfrac{1}{2} u^{2} + \mu B + \phi . &
\end{align}

This is the starting point for the derivation of a set of variational integrators for the guiding centre motion of charged particles in a tokamak.

\section{Variational Discretisation}

At first, the general derivation of a variational integrator for guiding centre motion is reproduced similar to \citeauthor{QinGuanTang:2009} \cite{QinGuanTang:2009}, where the trapezoidal rule is used to discretise the Lagrangian. In addition, we provide the derivation of a second integrator based on the midpoint rule that appears to be more stable at small timesteps and yields more accurate results.
In the last section the construction of higher order schemes by composition of low order schemes is sketched.

To allow for a compact notation, we introduce the generalised coordinates $q^{i} = X^{i}$ with $i = \{ 1,2,3 \}$ and $q^{u} = u$. Together, they are denoted $q^{\nu}$ with $\nu = \{ 1,2,3,u \}$ or just $q$. Correspondingly, the conjugate momenta are denoted $p^{\nu}$ or just $p$.

\subsection{Trapezoidal Discretisation}

Applying a trapezoidal discretisation to Littlejohn's guiding centre Lagrangian (\ref{eq:particles_guiding_centre_lagrangian}) gives
\begin{align}\label{eq:particles_guiding_centre_trapezoidal_lagrangian}
L_{d} (q_{k}, q_{k+1})
&= \dfrac{h}{2} \, L \bigg( q_{k  }, \dfrac{q_{k+1} - q_{k}}{h} \bigg)
+ \dfrac{h}{2} \, L \bigg( q_{k+1}, \dfrac{q_{k+1} - q_{k}}{h} \bigg) \\
\nonumber
&= h \, \bigg[ \dfrac{A_{i}^{*} (q_{k}) + A_{i}^{*} (q_{k+1})}{2} \dfrac{q^{i}_{k+1} - q^{i}_{k}}{h} - \dfrac{(q^{u}_{k})^{2} + (q^{u}_{k+1})^{2}}{4} \\
& \hspace{8em}
- \mu \, \dfrac{B (q_{k}) + B (q_{k+1})}{2} - \dfrac{\phi (q_{k}) + \phi (q_{k+1})}{2} \bigg]
.
\end{align}

This Lagrangian, however, results in a scheme with small stability region, as the expression resulting from the $(q^{u})^{2}$ term is explicit in $q^{u}$.
\citeauthor{QinGuanTang:2009} \cite{QinGuanTang:2009} replace this term with $q^{k} q^{k+1}$.
We explore a different modification in this term to make the resulting expression implicit in $q^{u}$, that is
\begin{align}\label{eq:particles_guiding_centre_trapezoidal_lagrangian_modified}
L_{d} (q_{k}, q_{k+1})
\nonumber
&= h \, \bigg[ \dfrac{A_{i}^{*} (q_{k}) + A_{i}^{*} (q_{k+1})}{2} \dfrac{q^{i}_{k+1} - q^{i}_{k}}{h} - \dfrac{1}{2} \, \bigg( \dfrac{q^{u}_{k} + q^{u}_{k+1}}{2} \bigg)^{2} \\
& \hspace{8em}
- \mu \, \dfrac{B (q_{k}) + B (q_{k+1})}{2} - \dfrac{\phi (q_{k}) + \phi (q_{k+1})}{2} \bigg]
.
\end{align}

Questions about which discretisations of the Lagrangian can be regarded as ``good'', i.e., produce well working, stable schemes, and which discretisations work less well remain largely unanswered.
Unfortunately, there exist no clear guidelines for the discretisation of the Lagrangian, but the preservation of symmetries in the course of discretisation certainly plays an important role.

Continuing with the derivation of the discrete Euler-Lagrange equations, these are defined as
\begin{subequations}\label{eq:particles_guiding_centre_trapezoidal_deleqs}
\begin{align}
\nonumber
\dfrac{\partial}{\partial q^{j}_{k}} & \left[ L_{d} (q_{k-1}, q_{k}) + L_{d} (q_{k}, q_{k+1}) \right] = \\
\label{eq:particles_guiding_centre_trapezoidal_deleqs_a}
&= h \, \bigg[ A^{*}_{i,j} (q_{k}) \, \dfrac{q^{i}_{k+1} - q^{i}_{k-1}}{2h} - \dfrac{A^{*}_{j} (q_{k+1}) - A^{*}_{j} (q_{k-1})}{2h} - \mu B_{,j} (q_{k}) - \phi_{,j} (q_{k}) \bigg] = 0 ,
\\
\nonumber
\dfrac{\partial}{\partial q^{u}_{k}} & \left[ L_{d} (q_{k-1}, q_{k}) + L_{d} (q_{k}, q_{k+1}) \right] = \\
\label{eq:particles_guiding_centre_trapezoidal_deleqs_b}
&= h \, \bigg[ b_{i} (q_{k}) \, \dfrac{q^{i}_{k+1} - q^{i}_{k-1}}{2h} - \dfrac{q^{u}_{k-1} + 2 \, q^{u}_{k} + q^{u}_{k+1}}{4} \bigg] = 0 .
\end{align}
\end{subequations}

This set of equations forms an implicit system for the solution of the guiding centre dynamics depending on data at three points in time, $q_{k+1}$, $q_{k}$ and $q_{k-1}$. It constitutes a discrete map
\begin{align}
( q_{k-1}, q_{k} ) \mapsto ( q_{k}, q_{k+1} ) .
\end{align}

Solving the system for $q_{k+1}$ yields nonlinearly implicit iteration rules for integrating the discrete phasespace trajectory $\{ q_{k} \}_{k=0}^{N}$ of the particle.
A possible solution strategy is to use a Newton solver for the nonlinear iteration. If the initial guess, e.g., by the linearised scheme derived below, is sufficiently close to the solution, no more than two or three iterations are needed. Using an analytic solution in the Newton iteration and a fixed number of iterations per timestep, such that an evaluation of the residual becomes unnecessary, the computational effort is about the same as for a standard fourth order Runge-Kutta method. The variational integrator is therefore computationally competitive to an explicit standard method while yielding superior results.

\subsubsection{Linearisation}

The linearisation of (\ref{eq:particles_guiding_centre_trapezoidal_deleqs}) allows for an easier comparison with the continuous equations of motion (\ref{eq:particles_guiding_centre_equations_of_motion_spatial_2}-\ref{eq:particles_guiding_centre_equations_of_motionetic_moment}) and highlights the differences between the variational integrator and a direct discretisation.
It is also possible to use this linearised scheme to compute an initial guess for a nonlinear solver applied to the above scheme.

Expand the $[ A^{*}_{j} (q_{k+1}) - A^{*}_{j} (q_{k-1}) ] / 2h$ term in (\ref{eq:particles_guiding_centre_trapezoidal_deleqs_a}) into a Taylor series about $q_{k}$
\begin{align}\label{eq:particles_guiding_centre_trapezoidal_linearisation_taylor}
\dfrac{A^{*}_{j} (q_{k+1}) - A^{*}_{j} (q_{k-1})}{2h} \approx A^{*}_{j,i} (q_{k}) \, \dfrac{q^{i}_{k+1} - q^{i}_{k-1}}{2h} + b_{j} (q_{k}) \, \dfrac{q^{u}_{k+1} - q^{u}_{k-1}}{2h}
.
\end{align}

Upon insertion into (\ref{eq:particles_guiding_centre_trapezoidal_deleqs_a}) we obtain a linearised set of equations
\begin{subequations}\label{eq:particles_guiding_centre_trapezoidal_linearisation_deleqs}
\begin{align}
\label{eq:particles_guiding_centre_trapezoidal_linearisation_deleqs_a}
& \Big[ A^{*}_{i,j} (q_{k}) - A^{*}_{j,i} (q_{k}) \Big] \, \dfrac{q^{i}_{k+1} - q^{i}_{k-1}}{2h} - b_{j} (q_{k}) \, \dfrac{q^{u}_{k+1} - q^{u}_{k-1}}{2h} - \mu B_{,j} (q_{k}) - \phi_{,j} (q_{k}) = 0 , \\
\label{eq:particles_guiding_centre_trapezoidal_linearisation_deleqs_b}
& b_{i} (q_{k}) \, \dfrac{q^{i}_{k+1} - q^{i}_{k-1}}{2h} - \dfrac{q^{u}_{k-1} + 2 \, q^{u}_{k} + q^{u}_{k+1}}{4} = 0
.
\end{align}
\end{subequations}

These equations resemble (\ref{eq:particles_guiding_centre_equations_of_motion_spatial_2}) and (\ref{eq:particles_guiding_centre_equations_of_motion_parallel_velocity}), while the original equation (\ref{eq:particles_guiding_centre_trapezoidal_deleqs_a}) resembles (\ref{eq:particles_guiding_centre_equations_of_motion_spatial_1}).

\subsubsection{Position Momentum Form}

The position-momentum form (\ref{eq:vi_finite_position_momentum}) of the trapezoidal integrator (\ref{eq:particles_guiding_centre_trapezoidal_deleqs}) is given by
\begin{subequations}
\begin{align}
p^{j}_{k  } &= - \dfrac{\partial}{\partial q^{j}_{k  }} L_{d} (q_{k}, q_{k+1}) , &
p^{j}_{k+1} &=   \dfrac{\partial}{\partial q^{j}_{k+1}} L_{d} (q_{k}, q_{k+1}) , & \\
p^{u}_{k  } &= - \dfrac{\partial}{\partial q^{u}_{k  }} L_{d} (q_{k}, q_{k+1}) , &
p^{u}_{k+1} &=   \dfrac{\partial}{\partial q^{u}_{k+1}} L_{d} (q_{k}, q_{k+1}) . &
\end{align}
\end{subequations}

Explicitly computing these expressions gives
\begin{subequations}
\begin{align}
p^{j}_{k  }
&=
- \dfrac{1}{2} \, A_{i,j}^{*} (q_{k}) \, \big[ q^{i}_{k+1} - q^{i}_{k} \big]
+ \dfrac{1}{2} \big[ A_{j}^{*} (q_{k}) + A_{j}^{*} (q_{k+1}) \big]
+ \dfrac{h}{2} \, \mu B_{,j} (q_{k})
+ \dfrac{h}{2} \, \phi_{,j} (q_{k}) ,
\\
p^{u}_{k  }
&=
- \dfrac{1}{2} \, b_{i} (q_{k}) \, \big[ q^{i}_{k+1} - q^{i}_{k} \big]
+ \dfrac{h}{4} \, \big[ q^{u}_{k} + q^{u}_{k+1} \big] ,
\\
p^{j}_{k+1}
&=
+ \dfrac{1}{2} \, A_{i,j}^{*} (q_{k+1}) \, \big[ q^{j}_{k+1} - q^{j}_{k} \big]
+ \dfrac{1}{2} \big[ A_{j}^{*} (q_{k}) + A_{j}^{*} (q_{k+1}) \big]
- \dfrac{h}{2} \, \mu B_{,j} (q_{k+1})
- \dfrac{h}{2} \, \phi_{,j} (q_{k+1}) ,
\\
p^{u}_{k+1}
&=
+ \dfrac{1}{2} \, b_{i} (q_{k+1}) \, \big[ q^{i}_{k+1} - q^{i}_{k} \big]
- \dfrac{h}{4} \, \big[ q^{u}_{k} - q^{u}_{k+1} \big] .
\end{align}
\end{subequations}

This set of equations forms a nonlinearly implicit system for the solution of the guiding centre dynamics depending only on data at two points in time, $q_{k+1}$ and $q_{k}$.
The first two equations have to be solved for the $q_{k+1}$. Afterwards, the $q_{k+1}$ are straight forwardly computed, as the other two equations are merely explicit functions.
The discrete map corresponding to this formulation is
\begin{align}
( q_{k}, p_{k} ) \mapsto ( q_{k+1}, p_{k+1} ) .
\end{align}

Solving this system, results in the same trajectory $\{ q_{k} \}_{k=0}^{N}$ as solving the system (\ref{eq:particles_guiding_centre_trapezoidal_deleqs}).

\subsubsection{Discrete Cartan Form}

The discrete Cartan one-forms corresponding to the Lagrangian (\ref{eq:particles_guiding_centre_trapezoidal_lagrangian_modified}) are given by
\begin{align}\label{eq:particles_guiding_centre_trapezoidal_cartan_oneform}
\Theta_{L_{d}}^{-} (q_{k}, q_{k+1}) &= - \dfrac{\partial}{\partial q^{\nu}_{k  }} L_{d} (q_{k}, q_{k+1}) \, d q^{\nu}_{k  } , &
\Theta_{L_{d}}^{+} (q_{k}, q_{k+1}) &= + \dfrac{\partial}{\partial q^{\nu}_{k+1}} L_{d} (q_{k}, q_{k+1}) \, d q^{\nu}_{k+1} , &
\end{align}

explicitly computed to be
\begin{align}
\Theta_{L_{d}}^{-} (q_{k}, q_{k+1})
\nonumber
&=
- \bigg\lgroup \dfrac{1}{h} \, A_{i,j}^{*} (q_{k}) \, \big[ q^{i}_{q_{k+1}} - q^{i}_{k} \big] - \dfrac{1}{h} \big[ A_{j}^{*} (q_{k}) + A_{j}^{*} (q_{k+1}) \big]  \\
\nonumber
& \hspace{4em}
- \dfrac{h}{2} \mu B_{,j} (q_{k}) - \dfrac{h}{2} \, \phi_{,j} (q_{k}) \bigg\rgroup \, dq^{j}_{k} \\
& \hspace{8em}
- \bigg\lgroup \dfrac{1}{2} \, b_{j} (q_{k}) \, \big[ q^{j}_{k+1} - q^{j}_{k} \big] - \dfrac{h}{4} \, \big[ q^{u}_{k} + q^{u}_{k+1} \big] \bigg\rgroup \, dq^{u}_{k} ,
\end{align}
\begin{align}
\Theta_{L_{d}}^{+} (q_{k}, q_{k+1})
\nonumber
&=
+ \bigg\lgroup \dfrac{1}{2} \, A_{i,j}^{*} (q_{k+1}) \, \big[ q^{i}_{k+1} - q^{i}_{k} \big] + \dfrac{1}{2} \big[ A_{j}^{*} (q_{k}) + A_{j}^{*} (q_{k+1}) \big] \\
\nonumber
& \hspace{4em}
- \dfrac{h}{2} \mu B_{,j} (q_{k+1}) - \dfrac{h}{2} \, \phi_{,j} (q_{k+1}) \bigg\rgroup \, dq^{j}_{k+1} \\
& \hspace{8em}
+ \bigg\lgroup \dfrac{1}{2} \, b_{j} (q_{k+1}) \, \big[ q^{j}_{k+1} - q^{j}_{k} \big] - \dfrac{h}{4} \, \big[ q^{u}_{k} + q^{u}_{k+1} \big] \bigg\rgroup \, dq^{u}_{k+1}
.
\end{align}

Together, these two forms determine the exterior derivative of the Lagrangian $L_{d} (q_{k}, q_{k+1})$
\begin{align}\label{eq:particles_guiding_centre_trapezoidal_cartan_ext_Lagrangian}
\ext L_{d} (q_{k}, q_{k+1}) = \Theta_{L_{d}}^{+} (q_{k}, q_{k+1}) - \Theta_{L_{d}}^{-} (q_{k}, q_{k+1}) .
\end{align}

As $\ext \ext L_{d} = 0$, the exterior derivative of both one-forms defines the same two-form
\begin{align}\label{eq:particles_guiding_centre_trapezoidal_cartan_symplectic_twoform}
\Omega_{d} (q_{k}, q_{k+1}) = \ext \Theta_{L_{d}}^{+} = \ext \Theta_{L_{d}}^{-} .
\end{align}

This is the discrete symplectic two-form $\Omega_{d}$.
Its preservation along the Lagrangian flow is given by construction and has been shown in section \ref{sec:vi_finite_symplectic_form}.

\subsection{Midpoint Discretisation}

We now derive an alternative method based on a midpoint discretisation. The main advantage over the trapezoidal scheme from the last section is a higher accuracy of the method.
We restrict ourselves to deriving the Euler-Lagrange equations and the position-momentum form of the equations of motion and do not repeat the derivation of the Cartan one-form.

Applying a midpoint discretisation to Littlejohn's guiding centre Lagrangian (\ref{eq:particles_guiding_centre_lagrangian}) gives
\begin{align}\label{eq:particles_guiding_centre_midpoint_lagrangian}
L_{d} (q_{k}, q_{k+1})
&= h \, L \bigg( \dfrac{q_{k} + q_{k+1}}{2}, \dfrac{q_{k+1} - q_{k}}{h} \bigg)
= h \, L \bigg( q_{k+1/2}, \dfrac{q_{k+1} - q_{k}}{h} \bigg) \\
&= h \, \bigg[ A_{i}^{*} (q_{k+1/2}) \, \dfrac{q^{i}_{k+1} - q^{i}_{k}}{h} - \left( \dfrac{q^{u}_{k} + q^{u}_{k+1}}{2} \right)^{2} - \mu B (q_{k+1/2}) - \phi (q_{k+1/2}) \bigg]
,
\end{align}

where $q_{k+1/2} = \tfrac{1}{2} ( q_{k} + q_{k+1} )$.
The discrete Euler-Lagrange equations of this midpoint discretisation are computed as
\begin{subequations}\label{eq:particles_guiding_centre_midpoint_deleqs}
\begin{align}
\nonumber
\label{eq:particles_guiding_centre_midpoint_deleqs_a}
\dfrac{\partial}{\partial q_{k}^{j}} & \Big[ L_{d} (q_{k-1}, q_{k}) + L_{d} (q_{k}, q_{k+1}) \Big] = \\
\nonumber
&= \dfrac{1}{2} \, \Big[ A^{*}_{i,j} (q_{k+1/2}) \Big] \Big[ q^{i}_{k+1} - q^{i}_{k} \Big]
+ \dfrac{1}{2} \, \Big[ A^{*}_{i,j} (q_{k-1/2}) \Big] \Big[ q^{i}_{k} - q^{i}_{k-1} \Big]
- \Big[ A^{*}_{j} (q_{k+1/2}) - A^{*}_{j} (q_{k-1/2}) \Big] \\
& \hspace{4em}
- \dfrac{h}{2} \, \mu \, \Big[ B_{,j} (q_{k-1/2}) + B_{,j} (q_{k+1/2}) \Big]
- \dfrac{h}{2} \, \Big[ \phi_{,j} (q_{k-1/2}) + \phi_{,j} (q_{k+1/2}) \Big] = 0 ,
\\
\nonumber
\label{eq:particles_guiding_centre_midpoint_deleqs_b}
\dfrac{\partial}{\partial q^{u}_{k}} & \left[ L_{d} (q_{k-1}, q_{k}) + L_{d} (q_{k}, q_{k+1}) \right] = \\
&= \dfrac{1}{2} \, \Big[ b_{i} (q_{k-1/2}) \left[ q^{i}_{k} - q^{i}_{k-1} \right] + b_{i} (q_{k+1/2}) \left[ q^{i}_{k+1} - q^{i}_{k} \right] \Big]
- \dfrac{h}{4} \Big[ q^{u}_{k-1} + 2 \, q^{u}_{k} + q^{u}_{k+1} \Big] = 0
.
\end{align}
\end{subequations}

The position-momentum form of the midpoint integrator is computed as
\begin{subequations}
\begin{align}
p^{i}_{k  }
&=
- \dfrac{1}{2} \, \big[ A^{*}_{i,j} (q_{k+1/2}) \big] \big[ q^{i}_{k+1} - q^{i}_{k} \big]
+ A^{*}_{j} (q_{k+1/2})
+ \dfrac{h}{2} \, \mu B_{,j} (q_{k+1/2})
+ \dfrac{h}{2} \, \phi_{,j} (q_{k+1/2})
\\
p^{u}_{k  }
&=
- \dfrac{1}{2} \, b_{i} (q_{k+1/2}) \Big[ q^{i}_{k+1} - q^{i}_{k} \Big]
+ \dfrac{h}{4} \, \Big[ q^{u}_{k} + q^{u}_{k+1} \Big]
\\
p^{i}_{k+1}
&=
+ \dfrac{1}{2} \, \big[ A^{*}_{i,j} (q_{k+1/2}) \big] \big[ q^{i}_{k+1} - q^{i}_{k} \big]
+ A^{*}_{j} (q_{k+1/2})
- \dfrac{h}{2} \, \mu B_{,j} (q_{k+1/2})
- \dfrac{h}{2} \, \phi_{,j} (q_{k+1/2})
\\
p^{u}_{k+1}
&=
+ \dfrac{1}{2} \, b_{i} (q_{k+1/2}) \Big[ q^{i}_{k+1} - q^{i}_{k} \Big]
- \dfrac{h}{4} \, \Big[ q^{u}_{k} + q^{u}_{k+1} \Big]
.
\end{align}
\end{subequations}

The first two equations have to be solved for $q_{k+1}$. The solution of the second two equations for $p_{k+1}$ is then straight forward.

\subsection{Higher Order Schemes}

If the discrete Lagrangian is self-adjoint, c.f. section \ref{vi:finite_composition_methods},
\begin{align}
L_{d}^{*} ( q_{k}, q_{k+1}, h ) \equiv - L_{d} ( q_{k+1}, q_{k}, -h ) = L_{d} ( q_{k+1}, q_{k}, h ) ,
\end{align}

the resulting variational integrator can be composed to yield higher order methods.
Both, the trapezoidal Lagrangian (\ref{eq:particles_guiding_centre_trapezoidal_lagrangian})
\begin{align*}
(L^{\text{tr}}_{d})^{*} ( q_{k}, q_{k+1}, h )
&= - (-h) \, \dfrac{1}{2} \, \bigg[ L \bigg( q_{k+1}, \dfrac{q_{k} - q_{k+1}}{-h} \bigg) - L \bigg( q_{k}, \dfrac{q_{k} - q_{k+1}}{-h} \bigg) \bigg] \\
&=   \dfrac{ h}{2} \, \bigg[ L \bigg( q_{k  }, \dfrac{q_{k+1} - q_{k}}{h} \bigg) - L \bigg( q_{k+1}, \dfrac{q_{k+1} - q_{k}}{h} \bigg) \bigg] \\
&= L_{d} ( q_{k}, q_{k+1}, h )
\end{align*}

and the midpoint Lagrangian (\ref{eq:particles_guiding_centre_midpoint_lagrangian})
\begin{align*}
(L^{\text{mp}}_{d})^{*} ( q_{k}, q_{k+1}, h )
&= - (-h) \, L \bigg( \dfrac{q_{k+1} + q_{k}}{2}, \dfrac{q_{k} - q_{k+1}}{-h} \bigg) \\
&= h \, L \bigg( \dfrac{q_{k} + q_{k+1}}{2}, \dfrac{q_{k+1} - q_{k}}{h} \bigg) \\
&= L^{\text{mp}}_{d} ( q_{k}, q_{k+1}, h )
\end{align*}

fulfil this property. The same is true for the modified trapezoidal Lagrangian (\ref{eq:particles_guiding_centre_trapezoidal_lagrangian_modified}), so that we can apply the composition rules from section \ref{vi:finite_composition_methods} to the schemes just derived.

\section{Particle Motion in the Poloidal Plane}

In this section we want to apply the derivations from the previous section to the motion of a charged particle in axisymmetric tokamaks. The toroidal symmetry allows us to reduce the dynamics to the poloidal plane of a tokamak.
At first, the derivation of \citeauthor{QinGuanTang:2009} \cite{QinGuanTang:2009} is reproduced, with the difference that cylindrical coordinates $(R,Z)$ are used instead of toroidal coordinates $(r, \theta)$.
The toroidal symmetry implies that the toroidal momentum $p_{\phy}$ is conserved and can be used to express the parallel velocity $u$ as a function of $(R,Z)$.
Hence, upon prescribing the values of the toroidal momentum $p_{\phy}$ and the magnetic moment $\mu$, only the coordinates of the poloidal plane $(R,Z)$ are treated as dynamical variables.

In \cite{QinGuanTang:2009}, the Lagrangian is discretised with the trapezoidal rule exclusively. We will also derive an integrator based on the midpoint rule. Finally we will compare all three of these schemes with a standard explicit Runge-Kutta method.
We restrict our treatment to the position-momentum form as that appears more natural with respect to the specification of initial conditions.

\subsubsection{Magnetic Field and Vector Potential}

For the magnetic field $B$ and the vector potential $A$ we will use analytic expressions following \citeauthor{QinGuanTang:2009} \cite{QinGuanTang:2009}.
The vector potential is given as
\begin{align}
A_{R} &= \dfrac{B_{0} R_{0} Z}{2 R} , &
A_{Z} &= - \ln \bigg( \dfrac{R}{R_{0}} \bigg) \, \dfrac{B_{0} R_{0}}{2} , &
A_{\phy} &= - \dfrac{B_{0} r^{2}}{2 q R} , &
r &= \sqrt{ (R-R_{0})^{2} + Z^{2} } ,
\end{align}

where subscripts $R$, $Z$ and $\phy$ denote the radial, vertical and toroidal components, respectively.
The magnetic field $B = \nabla \times A$ is
\begin{align*}
B_{R   } &= - \dfrac{B_{0} Z}{q R} , &
B_{Z   } &=   \dfrac{B_{0} \, (R - R_{0})}{q R} , &
B_{\phy} &= - \dfrac{B_{0} R_{0}}{R} , &
B &= \dfrac{B_{0} S}{q R} , &
\end{align*}

with the normalised magnetic field being
\begin{align}
b_{R   } &= - \dfrac{Z}{S} , &
b_{Z   } &= \dfrac{R-R_{0}}{S} , &
b_{\phy} &= - \dfrac{q R_{0}}{S} , &
S &= \sqrt{ r^{2} + q^{2} R_{0}^{2} } . &
\end{align}

Here, $R_{0}$ is the radial position of the magnetic axis, $B_{0}$ is the magnetic field at $R_{0}$, and $q$ is the safety factor, regarded as constant.
All the derivatives of the above expressions, which will be needed in the derived schemes, are listed in section \ref{sec:guiding_centre_derivatives}.

\subsubsection{Toroidal Momentum and Parallel Velocity}

It will be practical to express the Lagrangian with respect to the momenta $(p_{R}, p_{Z}, p_{\phy})$, such that in cylinder coordinates $(R, Z, \phy)$ we have
\begin{align}
L = p_{R} \dot{R} + p_{Z} \dot{Z} + p_{\phy} \dot{\phy} - H ,
\hspace{5em}
H = \dfrac{1}{2} u^{2} + \mu B - \phi ,
\end{align}

with
\begin{align}
p_{R}
&= A_{R}^{*} = A_{R} + u b_{R}, &
p_{Z}
&= A_{Z}^{*} = A_{Z} + u b_{Z}, &
p_{\phy}
&= R \, A_{\phy}^{*} = R \, \big( A_{\phy} + u b_{\phy} \big) .
\end{align}

As $p_{\phy}$ is conserved ($\partial L / \partial \phy = 0$),
\begin{align}
p_{\phy} = - R \, \bigg[ \dfrac{B_{0} r^{2}}{2 q R} + u \, \dfrac{q R_{0}}{S} \bigg] = const ,
\end{align}

we can compute a functional expression for the parallel velocity $u$ in which $p_{\phy}$ is a parameter
\begin{align}\label{eq:particles_parallel_velocity}
u
= - \bigg[ p_{\phy} + \dfrac{B_{0} r^{2}}{2 q} \bigg] \, \dfrac{S}{q R R_{0}}
= - \bigg[ p_{\phy} + \dfrac{B_{0} r^{2}}{2 q} \bigg] \, \dfrac{B}{B_{0} R_{0}}
.
\end{align}

\subsubsection{Reduced Lagrangian and Generalised Coordinates}

Projecting the motion to the poloidal plane and assuming the absence of any electrostatic field, the Lagrangian reduces to
\begin{align}
L = p_{R} \dot{R} + p_{Z} \dot{Z} - H
= A^{*}_{R} \dot{R} + A^{*}_{Z} \dot{Z} - \dfrac{1}{2} u^{2} - \mu B ,
\end{align}

where the components of the generalised magnetic potential read
\begin{align}
A^{*}_{R} &= \dfrac{B_{0} R_{0} Z}{2 R} + \bigg[ p_{\phy} + \dfrac{B_{0} r^{2}}{2 q} \bigg] \, \dfrac{Z}{q R R_{0}} , &
A^{*}_{Z} &= - \ln \bigg( \dfrac{R}{R_{0}} \bigg) \, \dfrac{B_{0} R_{0}}{2} - \bigg[ p_{\phy} + \dfrac{B_{0} r^{2}}{2 q} \bigg] \, \dfrac{R - R_{0}}{q R R_{0}} ,
\end{align}

and the parallel velocity $u$ is given by (\ref{eq:particles_parallel_velocity}).
We introduce generalised coordinates $q^{R}_{k}$ and $q^{Z}_{k}$ with discrete conjugate momenta $p^{R}_{k}$ and $p^{Z}_{k}$
\begin{align}\label{eq:particles_poloidal_plane_lagrangian}
L = A^{*}_{R} \dot{q}^{R} + A^{*}_{Z} \dot{q}^{Z} - \dfrac{1}{2} u^{2} - \mu B .
\end{align}

This is the basis for the following discretisations.

\subsection{Trapezoidal Discretisation}

Applying a trapezoidal discretisation to the reduced guiding centre Lagrangian (\ref{eq:particles_poloidal_plane_lagrangian}) gives
\begin{align}\label{eq:particles_poloidal_plane_lagrangian_trapezoidal}
L_{d}^{\text{tr}} (q_{k}, q_{k+1})
\nonumber
&= \dfrac{h}{2} \, L \bigg( q_{k  }, \dfrac{q_{k+1} - q_{k}}{h} \bigg)
+ \dfrac{h}{2} \, L \bigg( q_{k+1}, \dfrac{q_{k+1} - q_{k}}{h} \bigg) \\
\nonumber
&= h \, \bigg[
\dfrac{A^{*}_{R} (q_{k}) + A^{*}_{R} (q_{k+1})}{2} \dfrac{q^{R}_{k+1} - q^{R}_{k}}{h}
+ \dfrac{A^{*}_{Z} (q_{k}) + A^{*}_{Z} (q_{k+1})}{2} \dfrac{q^{Z}_{k+1} - q^{Z}_{k}}{h} \\
& \hspace{4em}
- \dfrac{1}{2} \, \bigg( \dfrac{u (q_{k}) + u (q_{k+1})}{2} \bigg)^{2}
- \mu \, \dfrac{B (q_{k}) + B (q_{k+1})}{2}
\bigg]
.
\end{align}

The position-momentum form (\ref{eq:vi_finite_position_momentum}) of the trapezoidal integrator is computed as
\begin{subequations}\label{eq:particles_poloidal_plane_position_momentum_trapezoidal}
\begin{align}
\label{eq:particles_poloidal_plane_position_momentum_trapezoidal_a}
p^{R}_{k}
\nonumber
&= \dfrac{1}{2} \bigg[
- A^{*}_{R,R} (q_{k}) \, \Big[ q^{R}_{k+1} - q^{R}_{k} \Big]
- A^{*}_{Z,R} (q_{k}) \, \Big[ q^{Z}_{k+1} - q^{Z}_{k} \Big]
+ \Big[ A^{*}_{R} (q_{k}) + A^{*}_{R} (q_{k+1}) \Big] \\
& \hspace{8em}
+ \dfrac{h}{2} \, u_{,R} (q_{k}) \, \Big[ u (q_{k}) + u (q_{k+1}) \Big]
+ h \, \mu B_{,R} (q_{k})
\bigg] ,
\\
\label{eq:particles_poloidal_plane_position_momentum_trapezoidal_b}
p^{Z}_{k}
\nonumber
&= \dfrac{1}{2} \bigg[
- A^{*}_{R,Z} (q_{k}) \, \Big[ q^{R}_{k+1} - q^{R}_{k} \Big]
- A^{*}_{Z,Z} (q_{k}) \, \Big[ q^{Z}_{k+1} - q^{Z}_{k} \Big]
+ \Big[ A^{*}_{Z} (q_{k}) + A^{*}_{Z} (q_{k+1}) \Big] \\
& \hspace{8em}
+ \dfrac{h}{2} \, u_{,Z} (q_{k}) \, \Big[ u (q_{k}) + u (q_{k+1}) \Big]
+ h \, \mu B_{,Z} (q_{k})
\bigg] ,
\\
\label{eq:particles_poloidal_plane_position_momentum_trapezoidal_c}
p^{R}_{k+1}
\nonumber
&= \dfrac{1}{2} \bigg[
+ A^{*}_{R,R} (q_{k+1}) \, \Big[ q^{R}_{k+1} - q^{R}_{k} \Big]
+ A^{*}_{Z,R} (q_{k+1}) \, \Big[ q^{Z}_{k+1} - q^{Z}_{k} \Big]
+ \Big[ A^{*}_{R} (q_{k}) + A^{*}_{R} (q_{k+1}) \Big] \\
& \hspace{8em}
- \dfrac{h}{2} \, u_{,R} (q_{k+1}) \, \Big[ u (q_{k}) + u (q_{k+1}) \Big]
- h \, \mu B_{,R} (q_{k+1})
\bigg] ,
\\
\label{eq:particles_poloidal_plane_position_momentum_trapezoidal_d}
p^{Z}_{k+1}
\nonumber
&= \dfrac{1}{2} \bigg[
+ A^{*}_{R,Z} (q_{k+1}) \, \Big[ q^{R}_{k+1} - q^{R}_{k} \Big]
+ A^{*}_{Z,Z} (q_{k+1}) \, \Big[ q^{Z}_{k+1} - q^{Z}_{k} \Big]
+ \Big[ A^{*}_{Z} (q_{k}) + A^{*}_{Z} (q_{k+1}) \Big] \\
& \hspace{8em}
- \dfrac{h}{2} \, u_{,Z} (q_{k+1}) \, \Big[ u (q_{k}) + u (q_{k+1}) \Big]
- h \, \mu B_{,Z} (q_{k+1})
\bigg]
.
\end{align}
\end{subequations}

The first two equations constitute a nonlinear system determining $q^{R}_{k+1}$ and $q^{Z}_{k+1}$ and can be written as a function $F (q_{k}, p_{k}, q_{k+1}) = 0$.
We solve it by Newton iteration with analytic Jacobian $\mcal{J}$, which is determined by computing the variation of the two equations with respect to $q^{R}_{k+1}$ and $q^{Z}_{k+1}$, such that
\begin{align}
\mcal{J} \, \delta q_{k+1}^{n+1}
=
\dfrac{1}{2}
\begin{pmatrix}
J_{11} & J_{12} \\
J_{21} & J_{22}
\end{pmatrix}
\begin{pmatrix}
\delta q^{R}_{k+1} \\
\delta q^{Z}_{k+1}
\end{pmatrix}
=
-
\begin{pmatrix}
F_{R} (q_{k}, p_{k}, q_{k+1}^{n}) \\
F_{Z} (q_{k}, p_{k}, q_{k+1}^{n})
\end{pmatrix} ,
\end{align}

where $n$ denotes the Newton step, such that
\begin{align}
q_{k+1}^{n+1} = q_{k+1}^{n} + \delta q_{k+1}^{n+1} .
\end{align}

The components of the Jacobian can be found in section \ref{sec:vi_finite_jacobians}.
After we solved for $q^{R}_{k+1}$ and $q^{Z}_{k+1}$, it is straight forward to evaluate (\ref{eq:particles_poloidal_plane_position_momentum_trapezoidal_c}) and (\ref{eq:particles_poloidal_plane_position_momentum_trapezoidal_d}) as these are explicit functions for $p^{R}_{k+1} (q_{k}, p_{k}, q_{k+1})$ and $p^{Z}_{k+1} (q_{k}, p_{k}, q_{k+1})$.

\subsection{Midpoint Discretisation}

Applying a midpoint discretisation to the reduced guiding centre Lagrangian (\ref{eq:particles_poloidal_plane_lagrangian}) gives
\begin{align}\label{eq:particles_poloidal_plane_lagrangian_midpoint}
L_{d}^{\text{mp}} & (q_{k}, q_{k+1})
\nonumber
= h \, L \bigg( \dfrac{q_{k} + q_{k+1}}{2}, \dfrac{q_{k+1} - q_{k}}{h} \bigg)
= h \, L \bigg( q_{k+1/2}, \dfrac{q_{k+1} - q_{k}}{h} \bigg) \\
& \hspace{2em}
= h \, \bigg[
A^{*}_{R} (q_{k+1/2}) \, \dfrac{q^{R}_{k+1} - q^{R}_{k}}{h}
+ A^{*}_{Z} (q_{k+1/2}) \, \dfrac{q^{Z}_{k+1} - q^{Z}_{k}}{h}
- \dfrac{u^{2} (q_{k+1/2})}{2}
- \mu B (q_{k+1/2})
\bigg]
.
\end{align}

The position-momentum form (\ref{eq:vi_finite_position_momentum}) of the midpoint integrator is computed as
\begin{subequations}\label{eq:particles_poloidal_plane_position_momentum_midpoint}
\begin{align}
\label{eq:particles_poloidal_plane_position_momentum_midpoint_a}
p^{R}_{k}
\nonumber
&= \dfrac{1}{2} \bigg[
- A^{*}_{R,R} (q_{k+1/2}) \, \Big[ q^{R}_{k+1} - q^{R}_{k} \Big]
- A^{*}_{Z,R} (q_{k+1/2}) \, \Big[ q^{Z}_{k+1} - q^{Z}_{k} \Big]
+ 2 \, A^{*}_{R} (q_{k+1/2}) \\
& \hspace{8em}
+ h \, u (q_{k+1/2}) \, u_{,R} (q_{k+1/2})
+ h \, \mu B_{,R} (q_{k+1/2})
\bigg] ,
\\
\label{eq:particles_poloidal_plane_position_momentum_midpoint_b}
p^{Z}_{k}
\nonumber
&= \dfrac{1}{2} \bigg[
- A^{*}_{R,Z} (q_{k+1/2}) \, \Big[ q^{R}_{k+1} - q^{R}_{k} \Big]
- A^{*}_{Z,Z} (q_{k+1/2}) \, \Big[ q^{Z}_{k+1} - q^{Z}_{k} \Big]
+ 2 \, A^{*}_{Z} (q_{k+1/2}) \\
& \hspace{8em}
+ h \, u (q_{k+1/2}) \, u_{,Z} (q_{k+1/2})
+ h \, \mu B_{,Z} (q_{k+1/2})
\bigg] ,
\\
\label{eq:particles_poloidal_plane_position_momentum_midpoint_c}
p^{R}_{k+1}
\nonumber
&= \dfrac{1}{2} \bigg[
+ A^{*}_{R,R} (q_{k+1/2}) \, \Big[ q^{R}_{k+1} - q^{R}_{k} \Big]
+ A^{*}_{Z,R} (q_{k+1/2}) \, \Big[ q^{Z}_{k+1} - q^{Z}_{k} \Big]
+ 2 \, A^{*}_{R} (q_{k+1/2}) \\
& \hspace{8em}
- h \, u (q_{k+1/2}) \, u_{,R} (q_{k+1/2})
- h \, \mu B_{,R} (q_{k+1/2})
\bigg] ,
\\
\label{eq:particles_poloidal_plane_position_momentum_midpoint_d}
p^{Z}_{k+1}
\nonumber
&= \dfrac{1}{2} \bigg[
+ A^{*}_{R,Z} (q_{k+1/2}) \, \Big[ q^{R}_{k+1} - q^{R}_{k} \Big]
+ A^{*}_{Z,Z} (q_{k+1/2}) \, \Big[ q^{Z}_{k+1} - q^{Z}_{k} \Big]
+ 2 \, A^{*}_{Z} (q_{k+1/2}) \\
& \hspace{8em}
- h \, u (q_{k+1/2}) \, u_{,Z} (q_{k+1/2})
- h \, \mu B_{,Z} (q_{k+1/2})
\bigg] .
\end{align}
\end{subequations}

The solution strategy is the same as for the trapezoidal method.

\subsection{Numerical Results}

We want to compare the different variational integrators with each other and a standard fourth order Runge-Kutta discretisation. Our main focus lies on the energy error and the geometry of the particle orbit for long time integration.
In the following, we consider the trapped particle example from \citeauthor{QinGuanTang:2009} \cite{QinGuanTang:2009}, which is initialised by
\begin{align*}
R &= R_{0} + 0.05 , &
Z &= 0 , &
\mu &= 2.25 \times 10^{-6} , &
p_{\phy} &= -1.077 \times 10^{-3} , &
\end{align*}

with
\begin{align*}
R_{0} &= 1 , &
B_{0} &= 1 , &
q &= 2 , &
\tau_{b} &= 43107 , &&
\end{align*}

where $\tau_{b}$ is the bounce time, determining the timestep. We use the normalisation proposed by \citeauthor{QinGuanTang:2009} \cite{QinGuanTang:2009}, where the parameters of the tokamak geometry are normalised by $R_{0}$ and $B_{0}$, to be able to compare with their results.

\begin{figure}[p]
\centering
\subfloat[Runge-Kutta]{
\includegraphics[width=.32\textwidth]{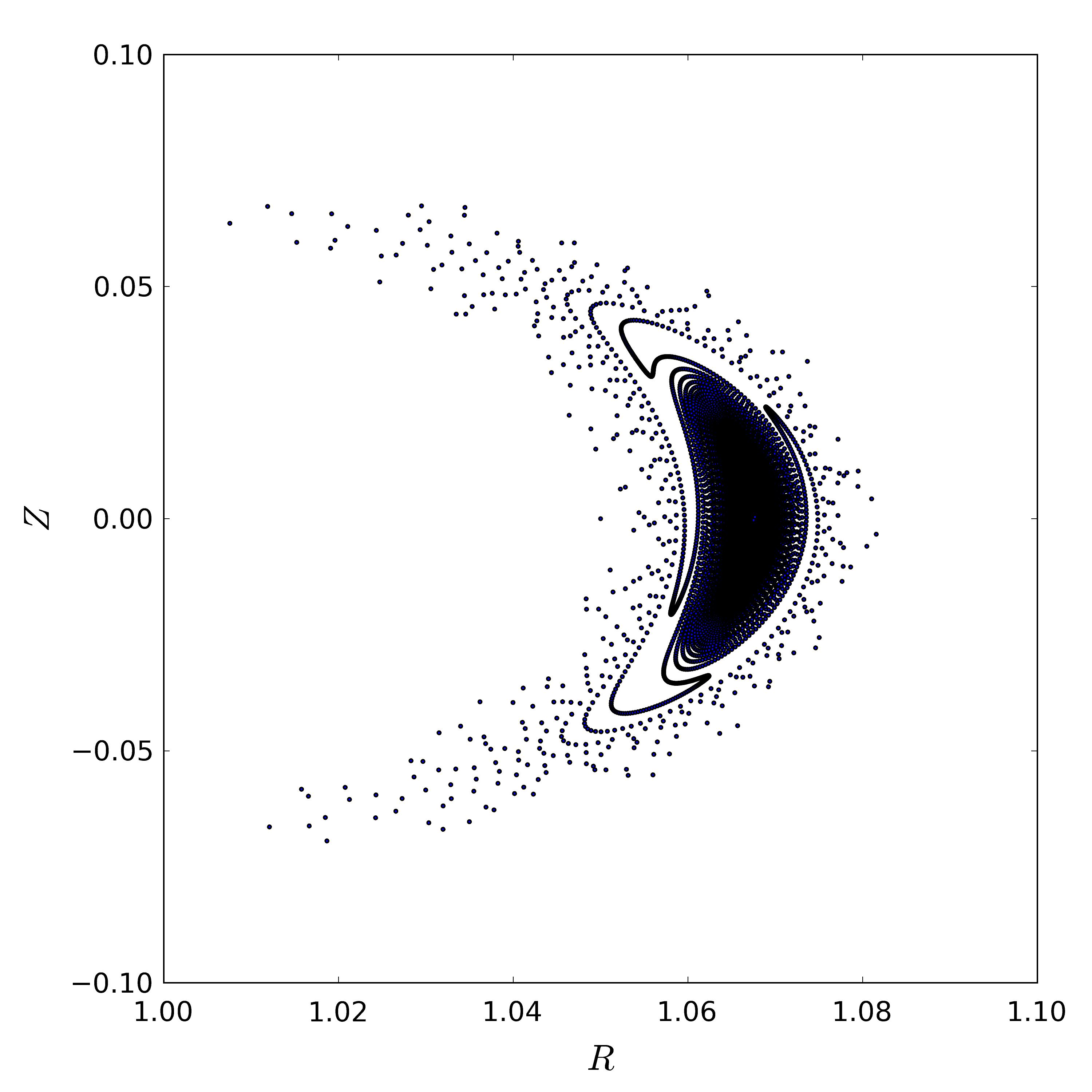}
}
\subfloat[Variational Midpoint]{
\includegraphics[width=.32\textwidth]{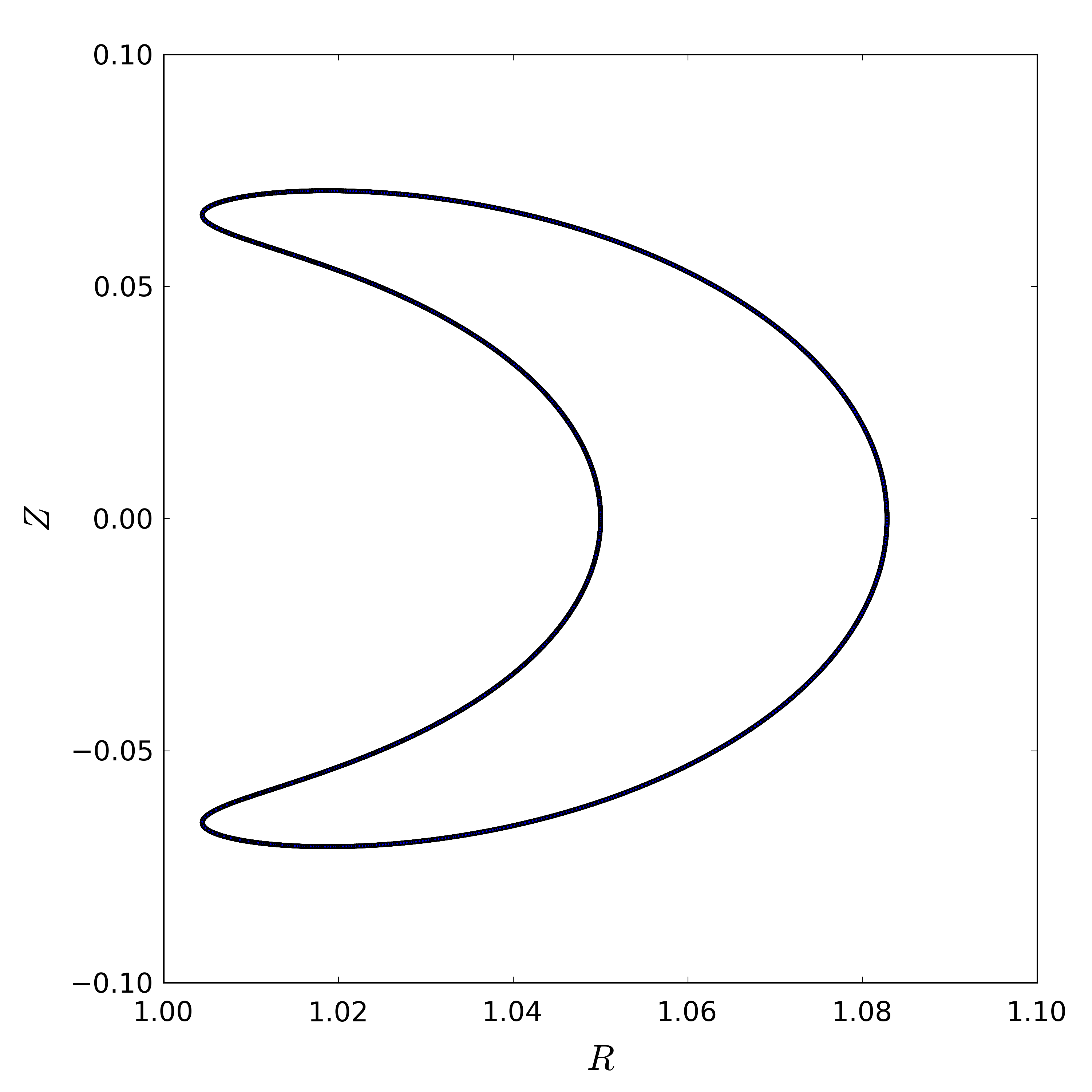}
}
\caption{Trapped particle after 8.000 bounce periods with 25 steps per bounce period. The trapezoidal integrator is not stable for the timestep considered in this example.}
\label{fig:particles_trapped_2d_nb25_orbits}
\end{figure}

\begin{figure}[p]
\centering
\subfloat[Runge-Kutta]{
\includegraphics[width=.32\textwidth]{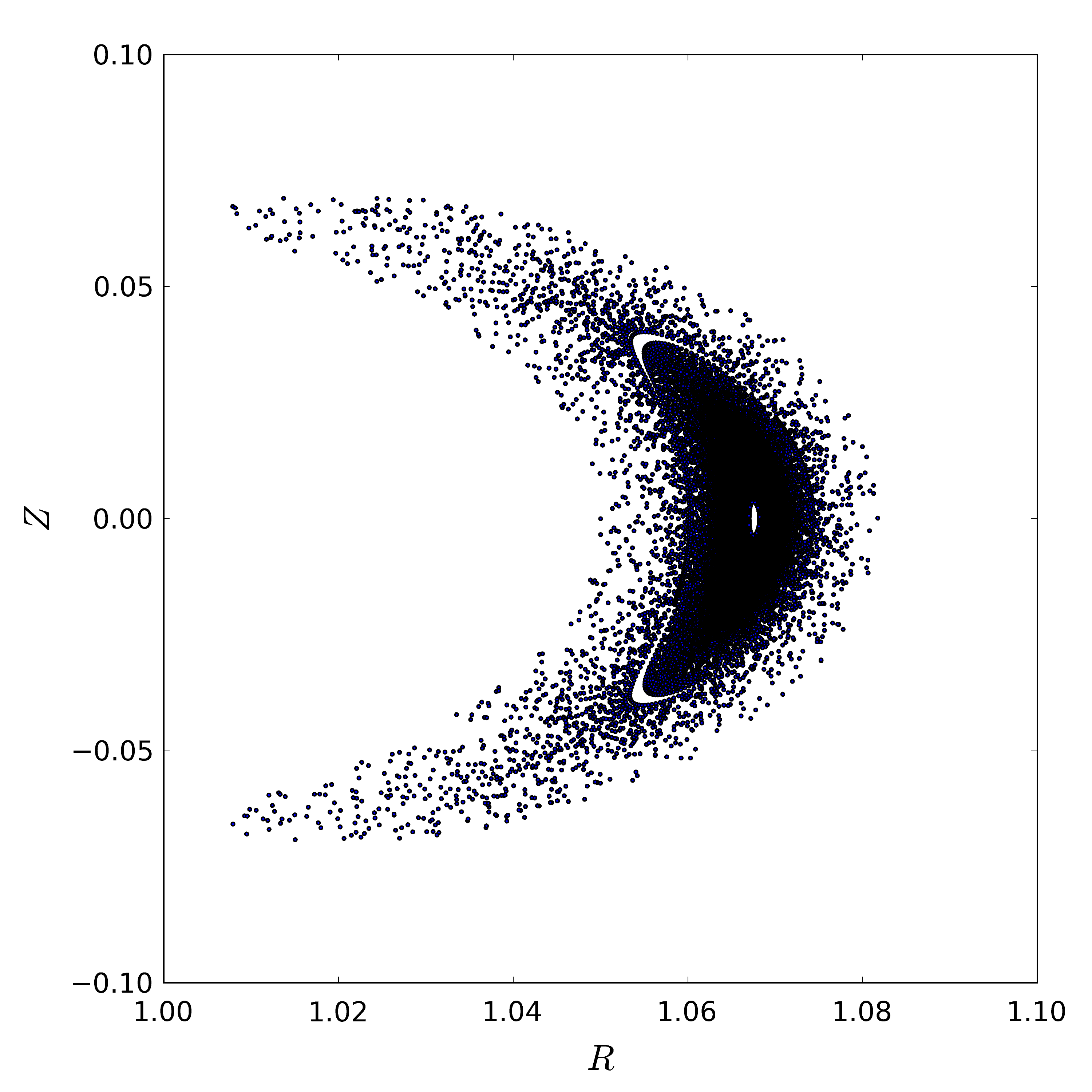}
}
\subfloat[Variational Trapezoidal]{
\includegraphics[width=.32\textwidth]{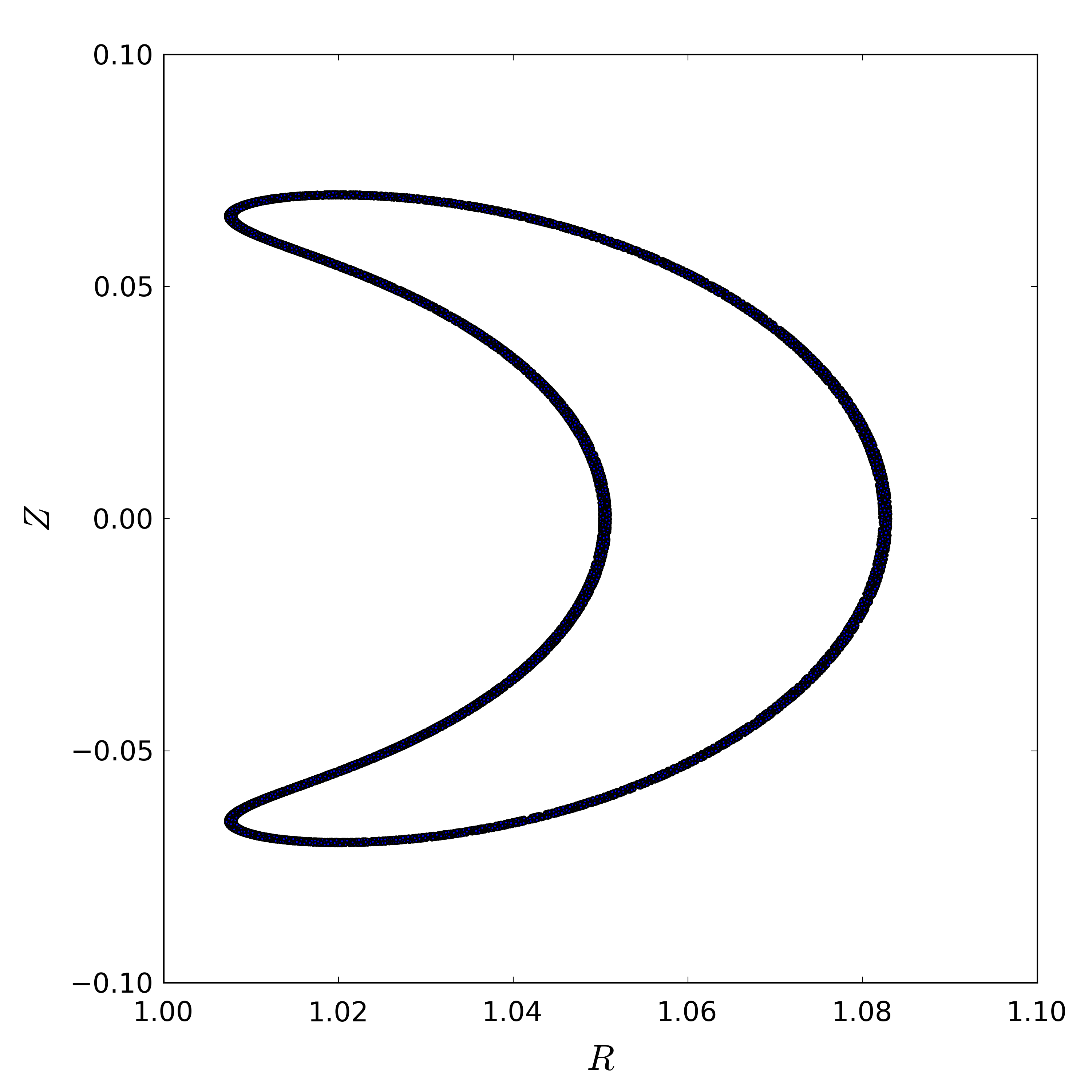}
}
\subfloat[Variational Midpoint]{
\includegraphics[width=.32\textwidth]{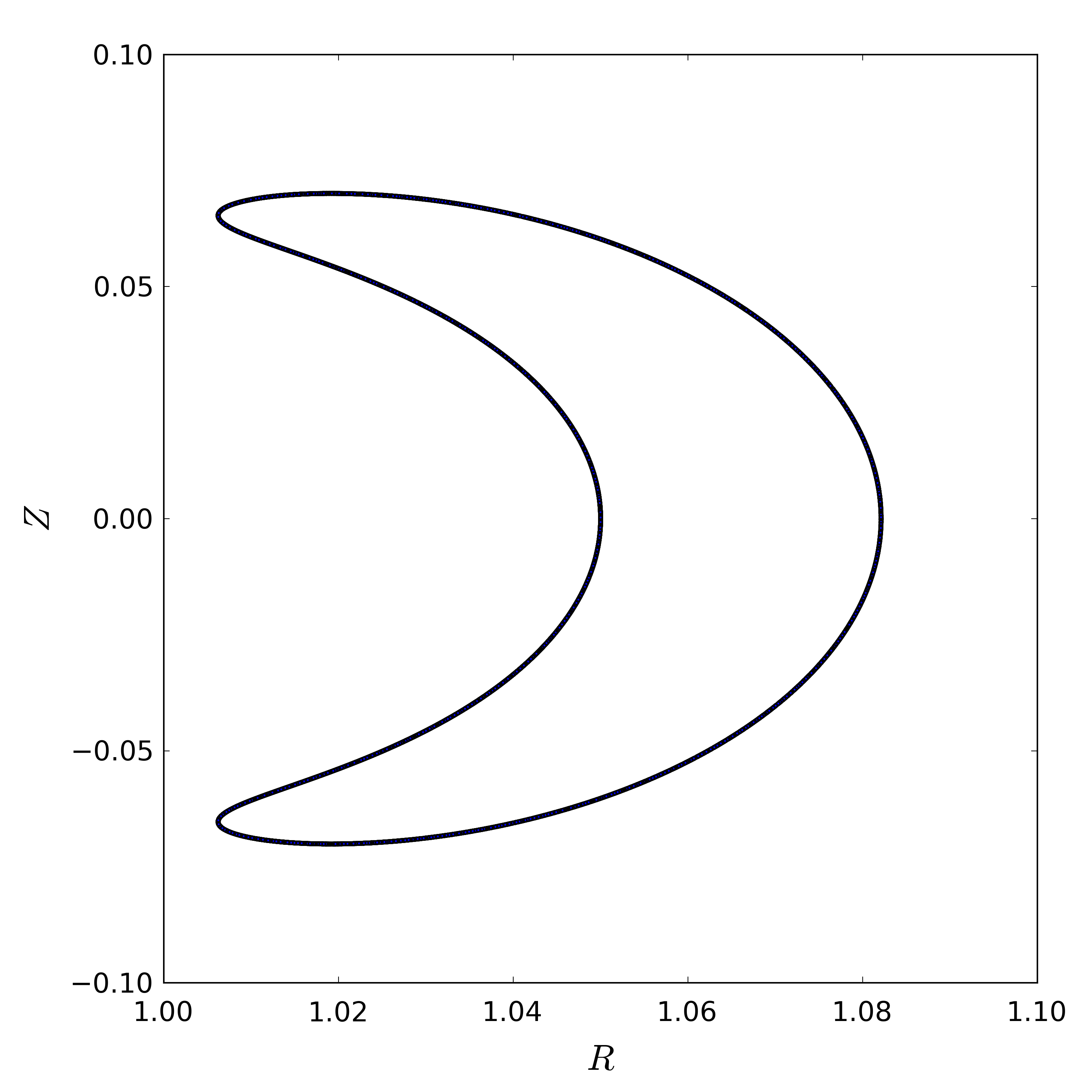}
}
\caption{Trapped particle after 200.000 bounce periods with 50 steps per bounce period.}
\label{fig:particles_trapped_2d_nb50_orbits}
\end{figure}

\begin{figure}[p]
\centering
\subfloat[Runge-Kutta]{
\includegraphics[width=.32\textwidth]{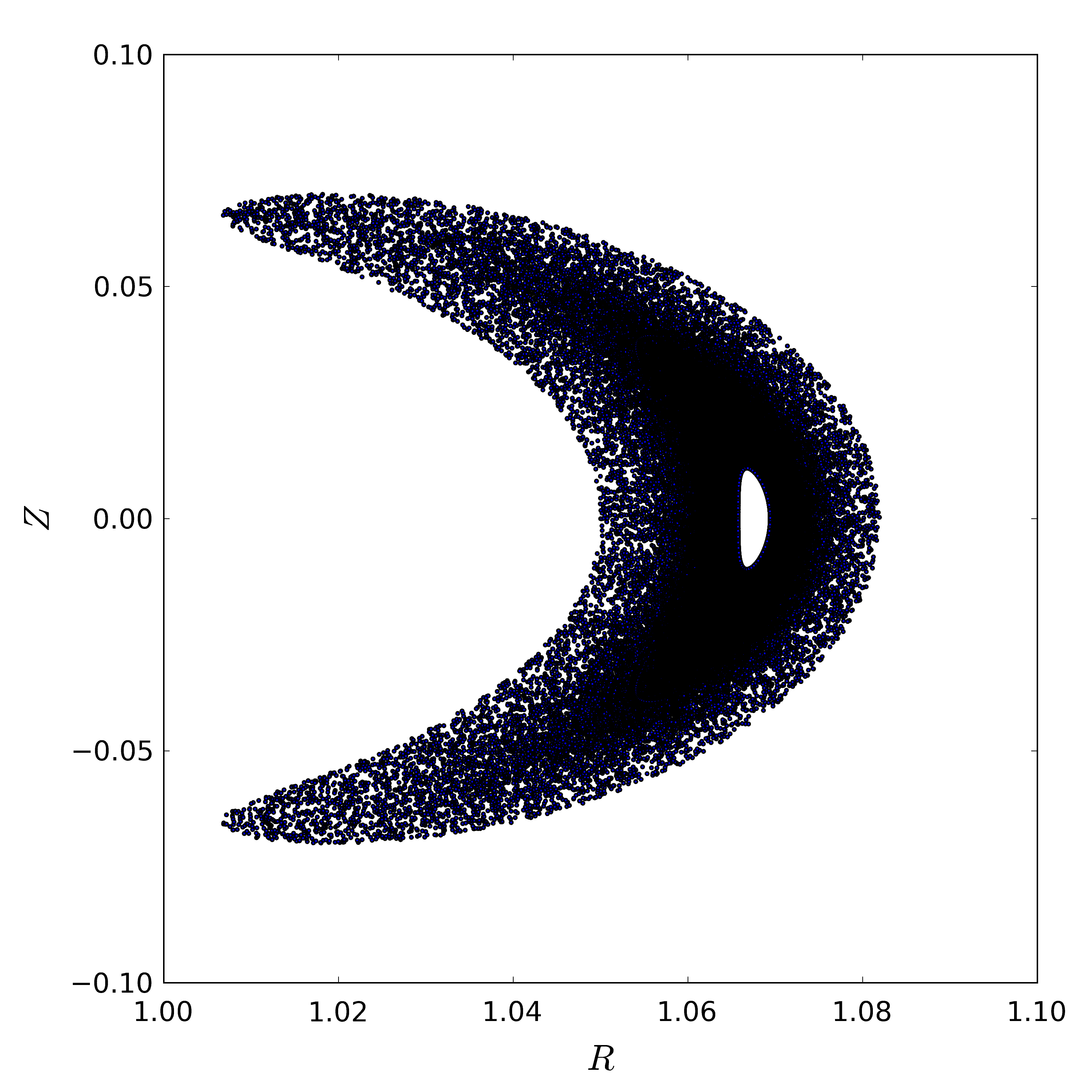}
}
\subfloat[Variational Trapezoidal]{
\includegraphics[width=.32\textwidth]{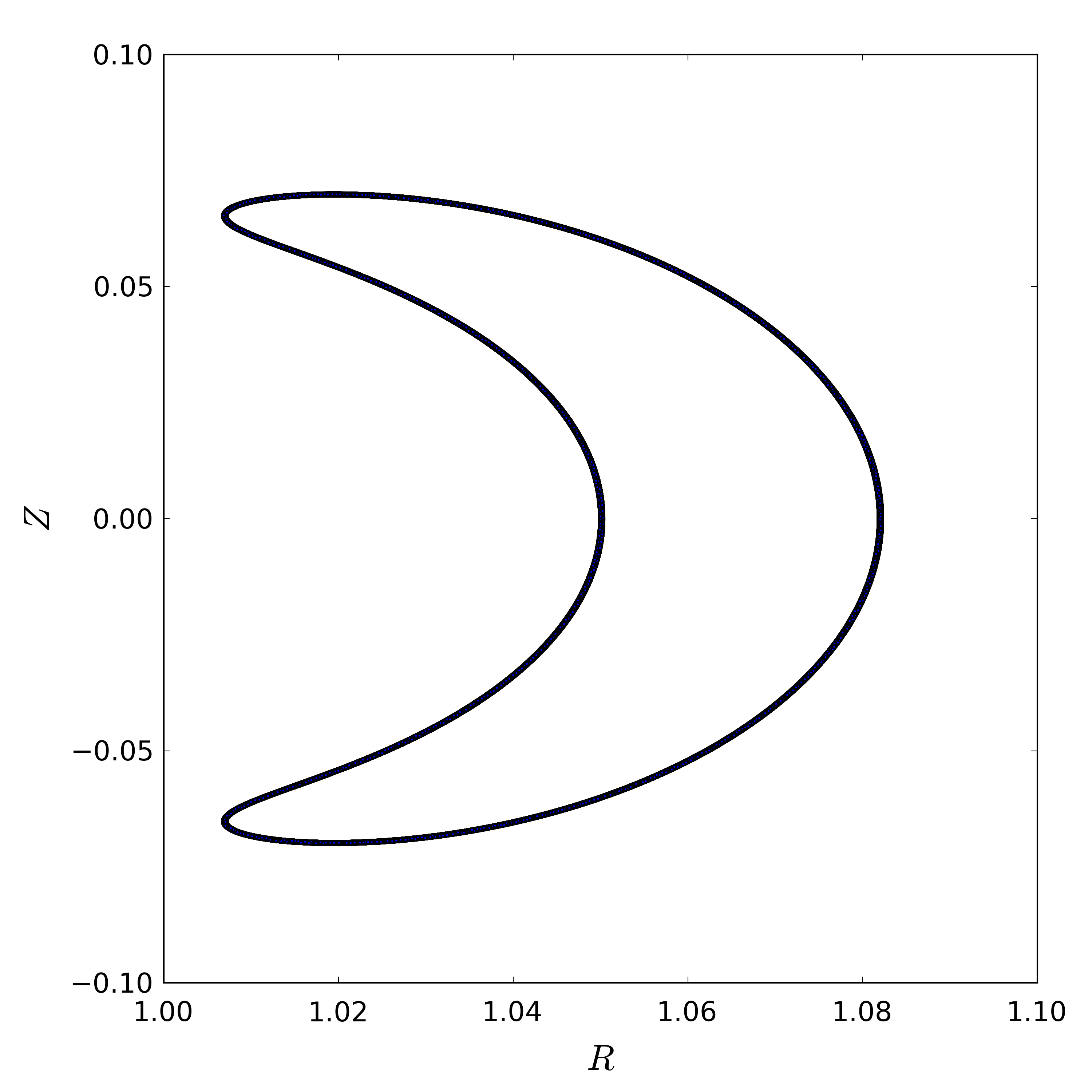}
}
\subfloat[Variational Midpoint]{
\includegraphics[width=.32\textwidth]{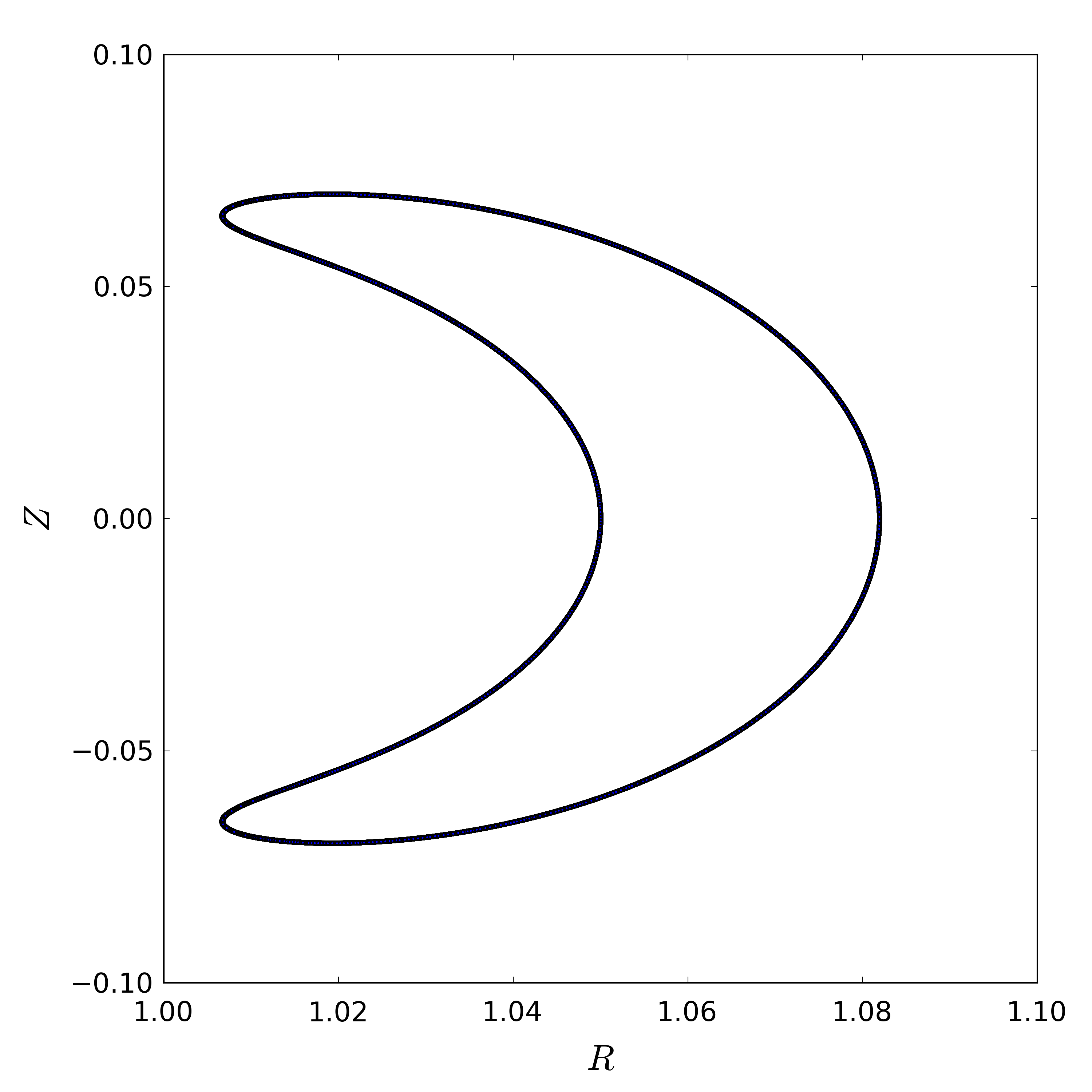}
}
\caption{Trapped particle after two million bounce periods with 100 steps per bounce period.}
\label{fig:particles_trapped_2d_nb100_orbits}
\end{figure}

\begin{figure}[p]
\centering
\includegraphics[width=.75\textwidth]{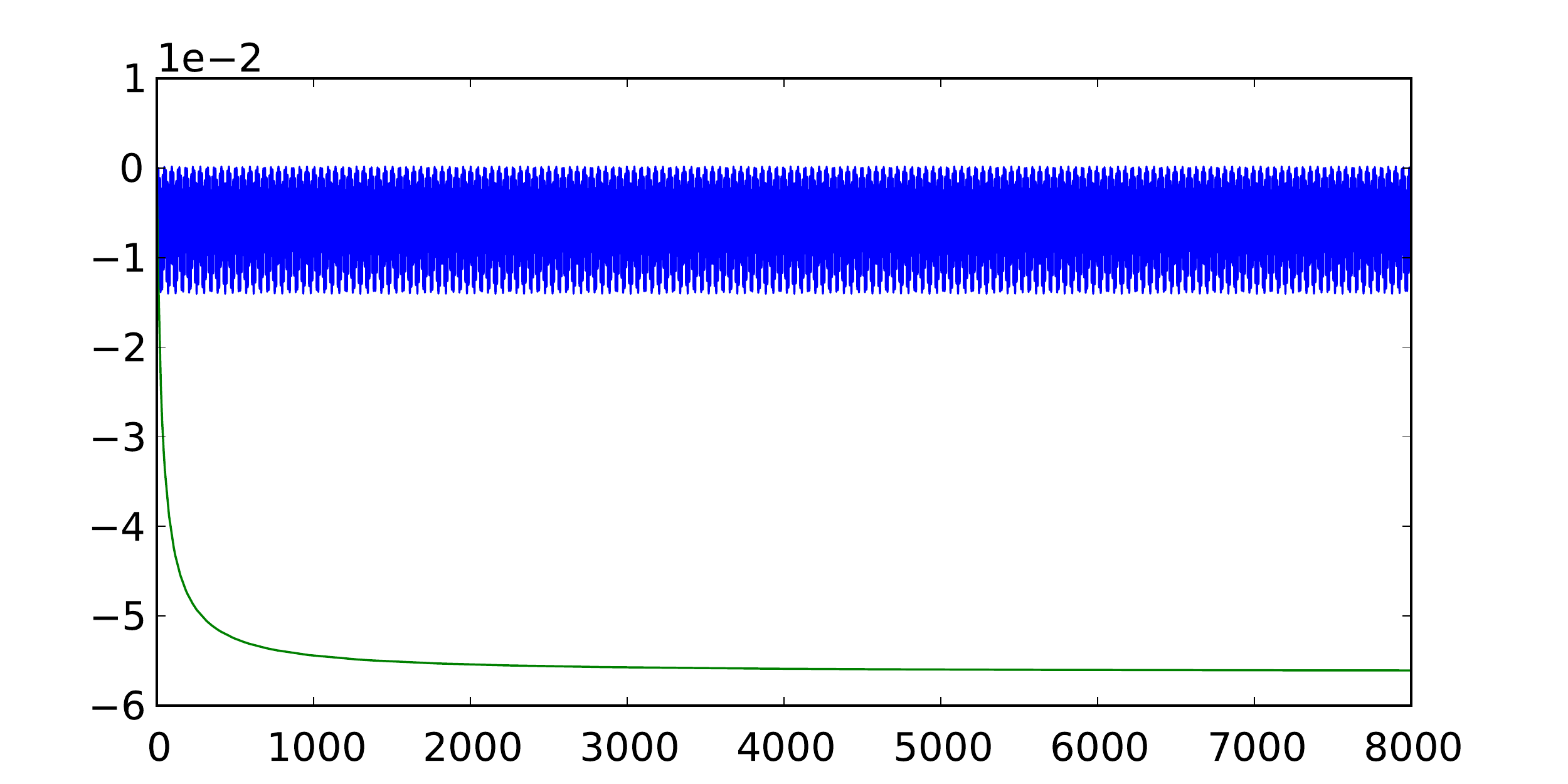}
\caption{Energy error for trapped particle with 25 timesteps per bounce period\\ (green: Runge-Kutta, blue: variational midpoint).}
\label{fig:particles_trapped_2d_nb25_energy}
\end{figure}

\begin{figure}[p]
\centering
\includegraphics[width=.75\textwidth]{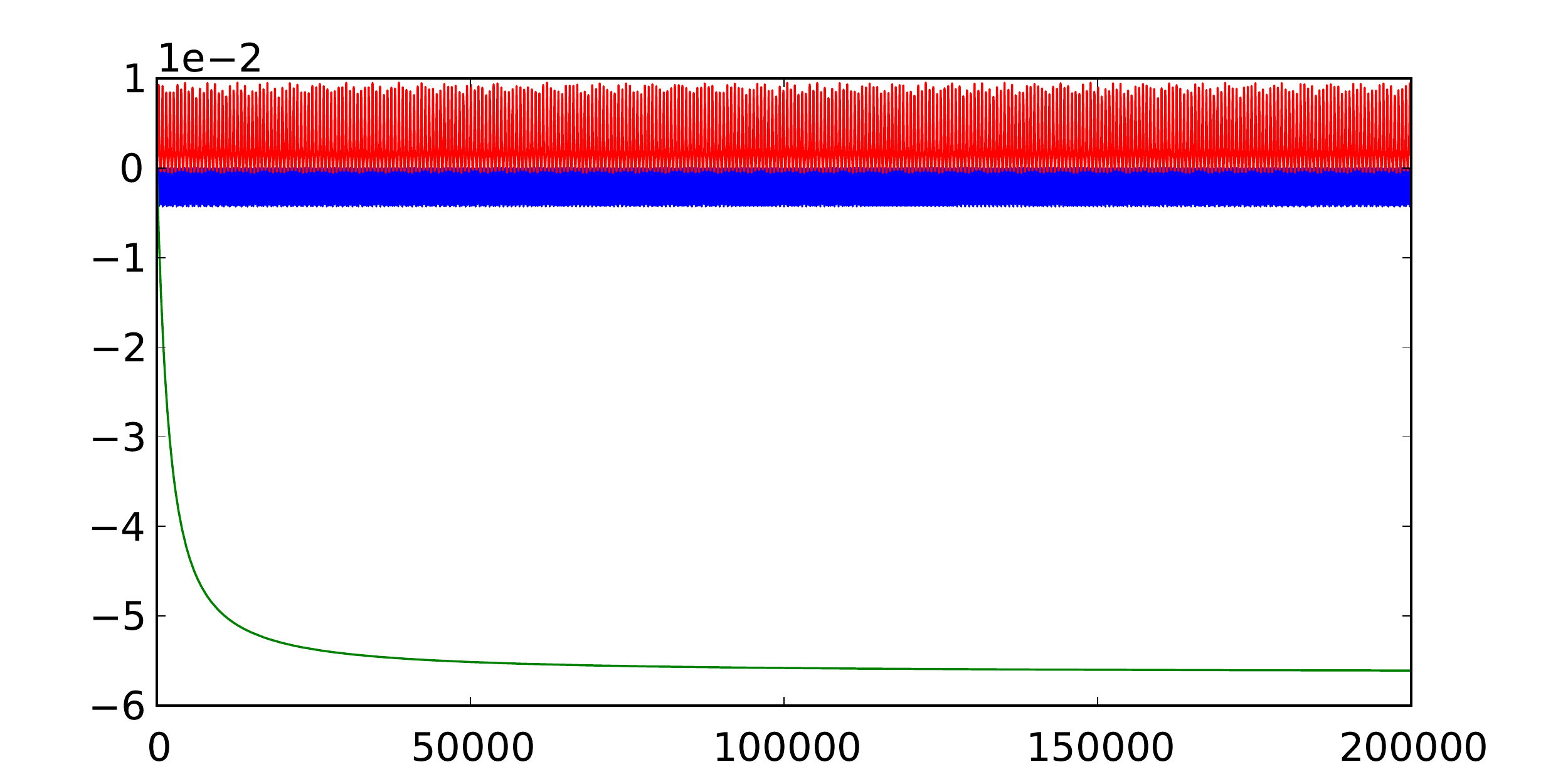}
\caption{Energy error for trapped particle with 50 timesteps per bounce period\\ (green: Runge-Kutta, blue: variational trapezoidal, red: variational midpoint).}
\label{fig:particles_trapped_2d_nb50_energy}
\end{figure}

\begin{figure}[p]
\centering
\includegraphics[width=.75\textwidth]{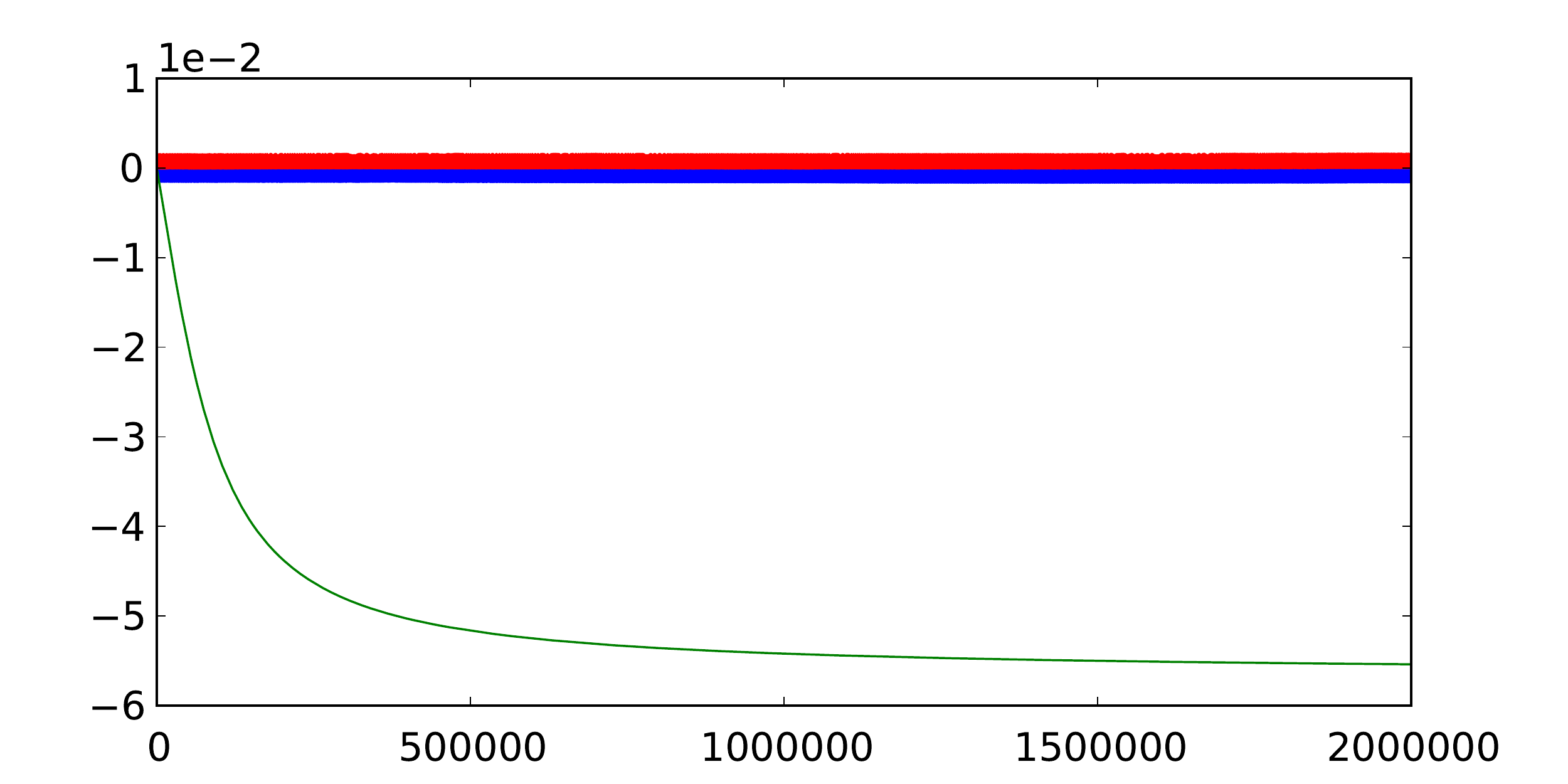}
\caption{Energy error for trapped particle with 100 timesteps per bounce period\\ (green: Runge-Kutta, blue: variational trapezoidal, red: variational midpoint).}
\label{fig:particles_trapped_2d_nb100_energy}
\end{figure}

We compare the evolution of the particle orbit for the two variational integrators and the Runge-Kutta method for three different timestep lengths that correspond to $25$, $50$ and $100$ timesteps per bounce period, respectively. In all three cases, we observe that the variational integrator follows the expected orbit accurately for long times, while the Runge-Kutta method exhibits substantial deviations (see figures \ref{fig:particles_trapped_2d_nb25_orbits} - \ref{fig:particles_trapped_2d_nb100_orbits}). Due to the loss of energy in the Runge-Kutta simulations (see figures \ref{fig:particles_trapped_2d_nb25_energy} - \ref{fig:particles_trapped_2d_nb100_energy}), the particle orbits shrink until eventually they almost contract to a point. On the contrary, the variational integrators exhibit an oscillating energy error, with a constant amplitude of the oscillation.
The amplitude of the energy error scales according to the order of the scheme, which is second-order accurate, i.e., halving the timestep results in a reduction of the error by a factor of four.

The stability region of the trapezoidal integrator appears to be smaller than that of the midpoint integrator as simulations with $25$ steps per bounce period are not possible.
At $50$ steps per bounce period, right from the beginning of the simulation, the trapezoidal integrator shows slight deviations of the orbit, but the global topology of the orbit is respected.
These inaccuracies are recognisable as a slight smearing of the orbit and are not observed for the midpoint integrator (see figure \ref{fig:particles_trapped_2d_nb25_orbits}).

Note that figures \ref{fig:particles_trapped_2d_nb100_orbits} and \ref{fig:particles_trapped_2d_nb100_energy} result from a simulation of two million characteristic times, corresponding to 200 million timesteps.
These results agree well with those found by \citeauthor{QinGuanTang:2009} \cite{QinGuanTang:2009}, but extend them in considering an alternative discretisation that appears to lead to more accurate results.

\section{Particle Motion in the Tokamak}

In this section, we want to discretise the particle motion in full tokamak geometry $(R, Z, \phy, u)$, such that also the toroidal coordinate $\phy$ and the parallel velocity $u$ are treated dynamically, while $\mu$ is still regarded as a parameter.
The resulting scheme is expected to be somewhat simpler than the previous one as we avoid the complicated functional expression for $u$.
Although this requires the solution of a larger system of equations describing the same dynamics, in this formulation both energy and toroidal momentum are independently conserved quantities. This allows us to test the variational integrators when the dynamics is constrained by more than one constant of motion.

As before, we introduce generalised coordinates $q^{R}_{k}$, $q^{Z}_{k}$, $q^{\phy}_{k}$, $q^{u}_{k}$ with discrete conjugate momenta $p^{R}_{k}$, $p^{Z}_{k}$, $p^{\phy}_{k}$, $p^{u}_{k}$.
For simplicity, we neglect the electrostatic field, could, however, be added with only minor complications.

\subsection{Trapezoidal Discretisation}

The discrete trapezoidal Lagrangian is
\begin{align}\label{eq:particles_tokamak_lagrangian_trapezoidal}
L_{d}^{\text{tr}} (q_{k}, q_{k+1})
\nonumber
&= h \, \bigg[
\dfrac{A_{R}^{*} (q_{k}) + A_{R}^{*} (q_{k+1})}{2} \, \dfrac{q^{R}_{k+1} - q^{R}_{k}}{h}
+ \dfrac{A_{Z}^{*} (q_{k}) + A_{Z}^{*} (q_{k+1})}{2} \, \dfrac{q^{Z}_{k+1} - q^{Z}_{k}}{h} \\
\nonumber
& \hspace{6em}
+ \dfrac{q^{R}_{k} \, A_{\phy}^{*} (q_{k}) + q^{R}_{k+1} \, A_{\phy}^{*} (q_{k+1})}{2} \, \dfrac{q^{\phy}_{k+1} - q^{\phy}_{k}}{h} \\
& \hspace{12em}
- \dfrac{q^{u}_{k} \, q^{u}_{k+1}}{2}
- \mu \, \dfrac{B (q_{k}) + B (q_{k+1})}{2}
\bigg] .
\end{align}

The position-momentum form (\ref{eq:vi_finite_position_momentum}) of the trapezoidal integrator is computed as
\begin{subequations}\label{eq:particles_tokamak_position_momentum_trapezoidal}
\begin{align}
\label{eq:particles_tokamak_position_momentum_trapezoidal_a}
p^{R}_{k}
\nonumber
&= \dfrac{1}{2} \bigg[
- A^{*}_{R,R} (q_{k}) \, \Big[ q^{R}_{k+1} - q^{R}_{k} \Big]
- A^{*}_{Z,R} (q_{k}) \, \Big[ q^{Z}_{k+1} - q^{Z}_{k} \Big]
- q^{R}_{k} \, A^{*}_{\phy,R} (q_{k}) \, \Big[ q^{\phy}_{k+1} - q^{\phy}_{k} \Big] \\
& \hspace{4em}
+ \Big[ A^{*}_{R} (q_{k}) + A^{*}_{R} (q_{k+1}) \Big]
- A_{\phy}^{*} (q_{k}) \, \Big[ q^{\phy}_{k+1} - q^{\phy}_{k} \Big]
+ h \, \mu B_{,R} (q_{k})
\bigg] ,
\\
\label{eq:particles_tokamak_position_momentum_trapezoidal_b}
p^{Z}_{k}
\nonumber
&= \dfrac{1}{2} \bigg[
- A^{*}_{R,Z} (q_{k}) \, \Big[ q^{R}_{k+1} - q^{R}_{k} \Big]
- A^{*}_{Z,Z} (q_{k}) \, \Big[ q^{Z}_{k+1} - q^{Z}_{k} \Big]
- q^{R}_{k} \, A^{*}_{\phy,Z} (q_{k}) \, \Big[ q^{\phy}_{k+1} - q^{\phy}_{k} \Big] \\
& \hspace{4em}
+ \Big[ A^{*}_{Z} (q_{k}) + A^{*}_{Z} (q_{k+1}) \Big]
+ h \, \mu B_{,Z} (q_{k})
\bigg] ,
\end{align}
\begin{align}
\label{eq:particles_tokamak_position_momentum_trapezoidal_c}
p^{\phy}_{k}
&= \dfrac{1}{2} \bigg[
+ \Big[ q^{R}_{k} \, A_{\phy}^{*} (q_{k}) + q^{R}_{k+1} \, A_{\phy}^{*} (q_{k+1}) \Big]
\bigg] ,
\\
\label{eq:particles_tokamak_position_momentum_trapezoidal_d}
p^{u}_{k}
&= \dfrac{1}{2} \bigg[
+ h \, q^{u}_{k+1}
- b_{R} (q_{k}) \, \Big[ q^{R}_{k+1} - q^{R}_{k} \Big]
- b_{Z} (q_{k}) \, \Big[ q^{Z}_{k+1} - q^{Z}_{k} \Big]
- q^{R}_{k} \, b_{\phy} (q_{k}) \, \Big[ q^{\phy}_{k+1} - q^{\phy}_{k} \Big]
\bigg] ,
\\
\label{eq:particles_tokamak_position_momentum_trapezoidal_e}
p^{R}_{k+1}
\nonumber
&= \dfrac{1}{2} \bigg[
+ A^{*}_{R,R} (q_{k+1}) \, \Big[ q^{R}_{k+1} - q^{R}_{k} \Big]
+ A^{*}_{Z,R} (q_{k+1}) \, \Big[ q^{Z}_{k+1} - q^{Z}_{k} \Big]
+ q^{R}_{k+1} \, A^{*}_{\phy,R} (q_{k+1}) \, \Big[ q^{\phy}_{k+1} - q^{\phy}_{k} \Big] \\
& \hspace{4em}
+ \Big[ A^{*}_{R} (q_{k}) + A^{*}_{R} (q_{k+1}) \Big]
+ A_{\phy}^{*} (q_{k+1}) \, \Big[ q^{\phy}_{k+1} - q^{\phy}_{k} \Big]
- h \, \mu B_{,R} (q_{k+1})
\bigg] ,
\\
\label{eq:particles_tokamak_position_momentum_trapezoidal_f}
p^{Z}_{k+1}
\nonumber
&= \dfrac{1}{2} \bigg[
+ A^{*}_{R,Z} (q_{k+1}) \, \Big[ q^{R}_{k+1} - q^{R}_{k} \Big]
+ A^{*}_{Z,Z} (q_{k+1}) \, \Big[ q^{Z}_{k+1} - q^{Z}_{k} \Big]
+ q^{R}_{k+1} \, A^{*}_{\phy,Z} (q_{k+1}) \, \Big[ q^{\phy}_{k+1} - q^{\phy}_{k} \Big] \\
& \hspace{4em}
+ \Big[ A^{*}_{Z} (q_{k}) + A^{*}_{Z} (q_{k+1}) \Big]
- h \, \mu B_{,Z} (q_{k+1})
\bigg] ,
\\
\label{eq:particles_tokamak_position_momentum_trapezoidal_g}
p^{\phy}_{k+1}
&= \dfrac{1}{2} \bigg[
+ \Big[ q^{R}_{k} \, A_{\phy}^{*} (q_{k}) + q^{R}_{k+1} \, A_{\phy}^{*} (q_{k+1}) \Big]
\bigg] ,
\\
\label{eq:particles_tokamak_position_momentum_trapezoidal_h}
p^{u}_{k+1}
&= \dfrac{1}{2} \bigg[
- h \, q^{u}_{k}
+ b_{R} (q_{k+1}) \, \Big[ q^{R}_{k+1} - q^{R}_{k} \Big]
+ b_{Z} (q_{k+1}) \, \Big[ q^{Z}_{k+1} - q^{Z}_{k} \Big]
+ q^{R}_{k+1} \, b_{\phy} (q_{k+1}) \, \Big[ q^{\phy}_{k+1} - q^{\phy}_{k} \Big]
\bigg]
.
\end{align}
\end{subequations}

As before, we employ Newton's method to solve the system. The Jacobian is now a $4 \times 4$ matrix, given by
\begin{align}
\mcal{J} \, \delta q_{k+1}^{n+1}
=
\dfrac{1}{2}
\begin{pmatrix}
J_{11} & J_{12} & J_{13} & J_{14} \\
J_{21} & J_{22} & J_{23} & J_{24} \\
J_{31} & J_{32} & J_{33} & J_{34} \\
J_{41} & J_{42} & J_{43} & J_{44}
\end{pmatrix}
\begin{pmatrix}
\delta q^{R}_{k+1} \\
\delta q^{Z}_{k+1} \\
\delta q^{\phy}_{k+1} \\
\delta q^{u}_{k+1}
\end{pmatrix}
=
-
\begin{pmatrix}
F_{R} (q_{k}, q_{k+1}^{n}) \\
F_{Z} (q_{k}, q_{k+1}^{n}) \\
F_{\phy} (q_{k}, q_{k+1}^{n}) \\
F_{u} (q_{k}, q_{k+1}^{n})
\end{pmatrix} ,
\end{align}

with components listed in section \ref{sec:vi_finite_jacobians}.

\subsection{Midpoint Discretisation}

The discrete midpoint Lagrangian is
\begin{multline}\label{eq:particles_tokamak_lagrangian_midpoint}
L_{d}^{\text{mp}} (q_{k}, q_{k+1})
= h \, \bigg[
A_{R}^{*} (q_{k+1/2}) \, \dfrac{q^{R}_{k+1} - q^{R}_{k}}{h}
+ A_{Z}^{*} (q_{k+1/2}) \, \dfrac{q^{Z}_{k+1} - q^{Z}_{k}}{h}
\\
+ q^{R}_{k+1/2} \, A_{\phy}^{*} (q_{k+1/2}) \, \dfrac{q^{\phy}_{k+1} - q^{\phy}_{k}}{h}
- \dfrac{1}{2} \, \bigg( \dfrac{q^{u}_{k} + q^{u}_{k+1}}{2} \bigg)^{2}
- \mu B (q_{k+1/2})
\bigg] .
\end{multline}

The position-momentum form (\ref{eq:vi_finite_position_momentum}) of the midpoint integrator is computed as
\begin{subequations}\label{eq:particles_tokamak_position_momentum_midpoint}
\begin{align}
\label{eq:particles_tokamak_position_momentum_midpoint_a}
p^{R}_{k}
\nonumber
&= \dfrac{1}{2} \bigg[
- A^{*}_{R,R} (q_{k+1/2}) \, \Big[ q^{R}_{k+1} - q^{R}_{k} \Big]
- A^{*}_{Z,R} (q_{k+1/2}) \, \Big[ q^{Z}_{k+1} - q^{Z}_{k} \Big]
- q^{R}_{k+1/2} \, A^{*}_{\phy,R} (q_{k+1/2}) \, \Big[ q^{\phy}_{k+1} - q^{\phy}_{k} \Big] \\
& \hspace{4em}
+ 2 \, A^{*}_{R} (q_{k+1/2})
- A^{*}_{\phy} (q_{k+1/2}) \, \Big[ q^{\phy}_{k+1} - q^{\phy}_{k} \Big]
+ h \, \mu B_{,R} (q_{k+1/2})
\bigg] ,
\\
\label{eq:particles_tokamak_position_momentum_midpoint_b}
p^{Z}_{k}
\nonumber
&= \dfrac{1}{2} \bigg[
- A^{*}_{R,Z} (q_{k+1/2}) \, \Big[ q^{R}_{k+1} - q^{R}_{k} \Big]
- A^{*}_{Z,Z} (q_{k+1/2}) \, \Big[ q^{Z}_{k+1} - q^{Z}_{k} \Big]
- q^{R}_{k+1/2} \, A^{*}_{\phy,Z} (q_{k+1/2}) \, \Big[ q^{\phy}_{k+1} - q^{\phy}_{k} \Big] \\
& \hspace{4em}
+ 2 \, A^{*}_{Z} (q_{k+1/2})
+ h \, \mu B_{,Z} (q_{k+1/2})
\bigg] ,
\end{align}
\begin{align}
\label{eq:particles_tokamak_position_momentum_midpoint_c}
p^{\phy}_{k}
&= \dfrac{1}{2} \bigg[
+ 2 \, q^{R}_{k+1/2} \, A_{\phy}^{*} (q_{k+1/2})
\bigg] ,
\\
\label{eq:particles_tokamak_position_momentum_midpoint_d}
p^{u}_{k}
\nonumber
&= \dfrac{1}{2} \bigg[
- b_{R} (q_{k+1/2}) \, \Big[ q^{R}_{k+1} - q^{R}_{k} \Big]
- b_{Z} (q_{k+1/2}) \, \Big[ q^{Z}_{k+1} - q^{Z}_{k} \Big]
- q^{R}_{k+1/2} \, b_{\phy} (q_{k+1/2}) \, \Big[ q^{\phy}_{k+1} - q^{\phy}_{k} \Big] \\
& \hspace{4em}
+ h \, q^{u}_{k+1/2}
\bigg] ,
\\
\label{eq:particles_tokamak_position_momentum_midpoint_e}
p^{R}_{k+1}
\nonumber
&= \dfrac{1}{2} \bigg[
+ A^{*}_{R,R} (q_{k+1/2}) \, \Big[ q^{R}_{k+1} - q^{R}_{k} \Big]
+ A^{*}_{Z,R} (q_{k+1/2}) \, \Big[ q^{Z}_{k+1} - q^{Z}_{k} \Big]
+ q^{R}_{k+1/2} \, A^{*}_{\phy,R} (q_{k+1/2}) \, \Big[ q^{\phy}_{k+1} - q^{\phy}_{k} \Big] \\
& \hspace{4em}
+ 2 \, A^{*}_{R} (q_{k+1/2})
+ A_{\phy}^{*} (q_{k+1/2}) \, \Big[ q^{\phy}_{k+1} - q^{\phy}_{k} \Big]
- h \, \mu B_{,R} (q_{k+1/2})
\bigg] ,
\\
\label{eq:particles_tokamak_position_momentum_midpoint_f}
p^{Z}_{k+1}
\nonumber
&= \dfrac{1}{2} \bigg[
+ A^{*}_{R,Z} (q_{k+1/2}) \, \Big[ q^{R}_{k+1} - q^{R}_{k} \Big]
+ A^{*}_{Z,Z} (q_{k+1/2}) \, \Big[ q^{Z}_{k+1} - q^{Z}_{k} \Big]
+ q^{R}_{k+1/2} \, A^{*}_{\phy,Z} (q_{k+1/2}) \, \Big[ q^{\phy}_{k+1} - q^{\phy}_{k} \Big] \\
& \hspace{4em}
+ 2 \, A^{*}_{Z} (q_{k+1/2})
- h \, \mu B_{,Z} (q_{k+1/2})
\bigg] ,
\\
\label{eq:particles_tokamak_position_momentum_midpoint_g}
p^{\phy}_{k+1}
&= \dfrac{1}{2} \bigg[
+ 2 \, q^{R}_{k+1/2} \, A_{\phy}^{*} (q_{k+1/2})
\bigg] ,
\\
\label{eq:particles_tokamak_position_momentum_midpoint_h}
p^{u}_{k+1}
\nonumber
&= \dfrac{1}{2} \bigg[
+ b_{R} (q_{k+1/2}) \, \Big[ q^{R}_{k+1} - q^{R}_{k} \Big]
+ b_{Z} (q_{k+1/2}) \, \Big[ q^{Z}_{k+1} - q^{Z}_{k} \Big]
+ q^{R}_{k+1/2} \, b_{\phy} (q_{k+1/2}) \, \Big[ q^{\phy}_{k+1} - q^{\phy}_{k} \Big] \\
& \hspace{4em}
- h \, q^{u}_{k+1/2}
\bigg]
.
\end{align}
\end{subequations}

\subsection{Discrete Noether Theorem}

In the four-dimensional treatment, the interesting question is that of conservation of the toroidal momentum $p_{\phy}$. As we have discussed in the last section, $p_{\phy}$ is a conserved quantity of the continuous system and should therefore be exactly conserved by the variational integrator.

The corresponding transformation is
\begin{align}
\phy_{k}^{\eps} = \phy_{k} + \eps X_{k}^{\phy} .
\end{align}

Both, the trapezoidal (\ref{eq:particles_tokamak_lagrangian_trapezoidal}) and the midpoint (\ref{eq:particles_tokamak_lagrangian_midpoint}) Lagrangian are invariant under this transformation, as can easily be seen.
The discrete conserved momenta are
\begin{align}
\dfrac{\partial L_{d}^{\text{tr}}}{\partial q^{\phy}_{k+1}} (q_{k}, q_{k+1}) \cdot X_{k}^{\phy}
&= \dfrac{1}{2} \Big( q^{R}_{k} \, A_{\phy}^{*} (q_{k}) + q^{R}_{k+1} \, A_{\phy}^{*} (q_{k+1}) \Big) \cdot X_{k}^{\phy}
= p^{\phy}_{k} \cdot X_{k}^{\phy}
\end{align}

for the trapezoidal discretisation and
\begin{align}
\dfrac{\partial L_{d}^{\text{mp}}}{\partial q^{\phy}_{k+1}} (q_{k}, q_{k+1}) \cdot X_{k}^{\phy}
&= q^{R}_{k+1/2} \, A_{\phy}^{*} (q_{k+1/2}) \cdot X_{k}^{\phy}
= p^{\phy}_{k} \cdot X_{k}^{\phy}
\end{align}

for the midpoint discretisation.

\subsection{Numerical Results}

We use the same initial conditions as in the previous section, that is
\begin{align*}
R &= R_{0} + 0.05 , &
Z &= 0 , &
\phy &= 0, &
\mu &= 2.25 \times 10^{-6} , &
\end{align*}

with
\begin{align*}
R_{0} &= 1 , &
B_{0} &= 1 , &
q &= 2 , &
\tau_{b} &= 43107 . & &&
\end{align*}

The initial parallel velocity $u$ is computed by evaluating equation (\ref{eq:particles_parallel_velocity}) for the initial values of $R$, $Z$ and $p^{\phy}$.
The initial momenta are given below, where $p^{Z} = A^{*}_{Z} (q_{0})$, i.e.,
\begin{align}
p^{R} &= 0 , &
p^{Z} &= - 2.438 \times 10^{-2} , &
p^{\phy} &= - 1.077 \times 10^{-3} , &
p^{u} &= 0 . &
\end{align}

\begin{figure}[p]
\centering
\subfloat[Runge-Kutta]{
\includegraphics[width=.32\textwidth]{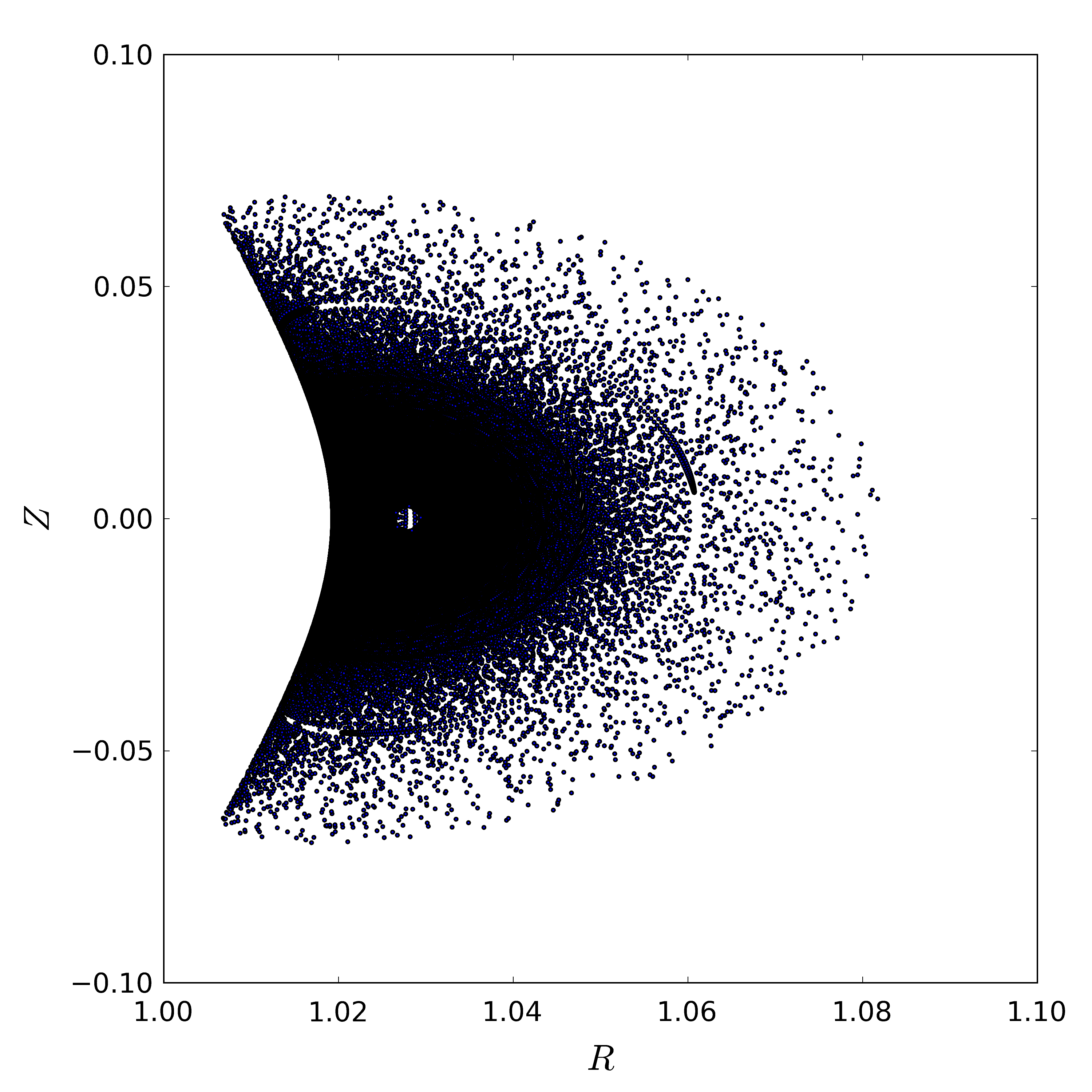}
}
\subfloat[Variational Midpoint]{
\includegraphics[width=.32\textwidth]{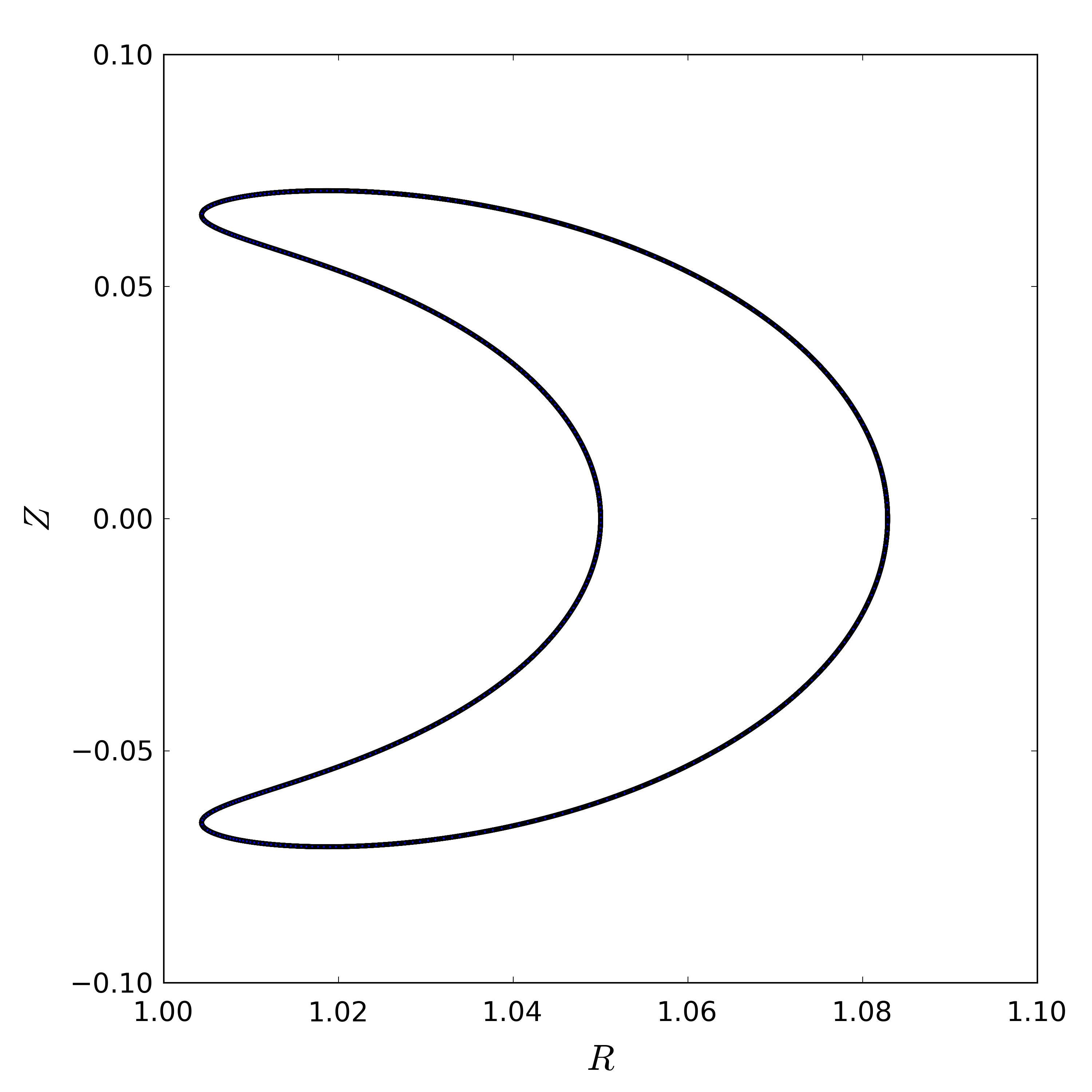}
}
\caption{Trapped particle after 50.000 bounce periods with 25 steps per bounce period. The trapezoidal integrator is not stable for the timestep considered in this example.}
\label{fig:particles_4d_trapped_nb25_orbits}
\end{figure}

\begin{figure}[p]
\centering
\subfloat[Runge-Kutta]{
\includegraphics[width=.32\textwidth]{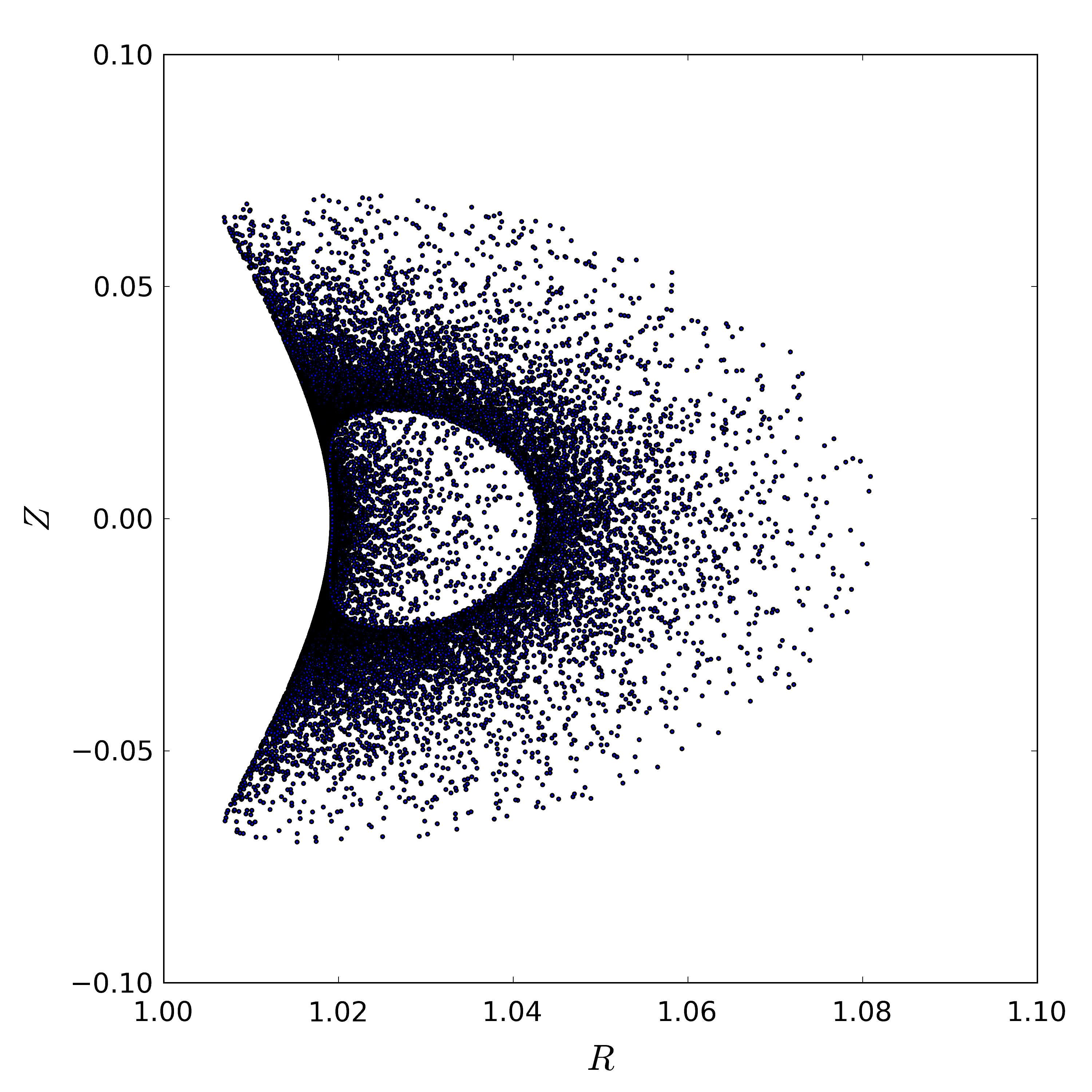}
}
\subfloat[Variational Trapezoidal]{
\includegraphics[width=.32\textwidth]{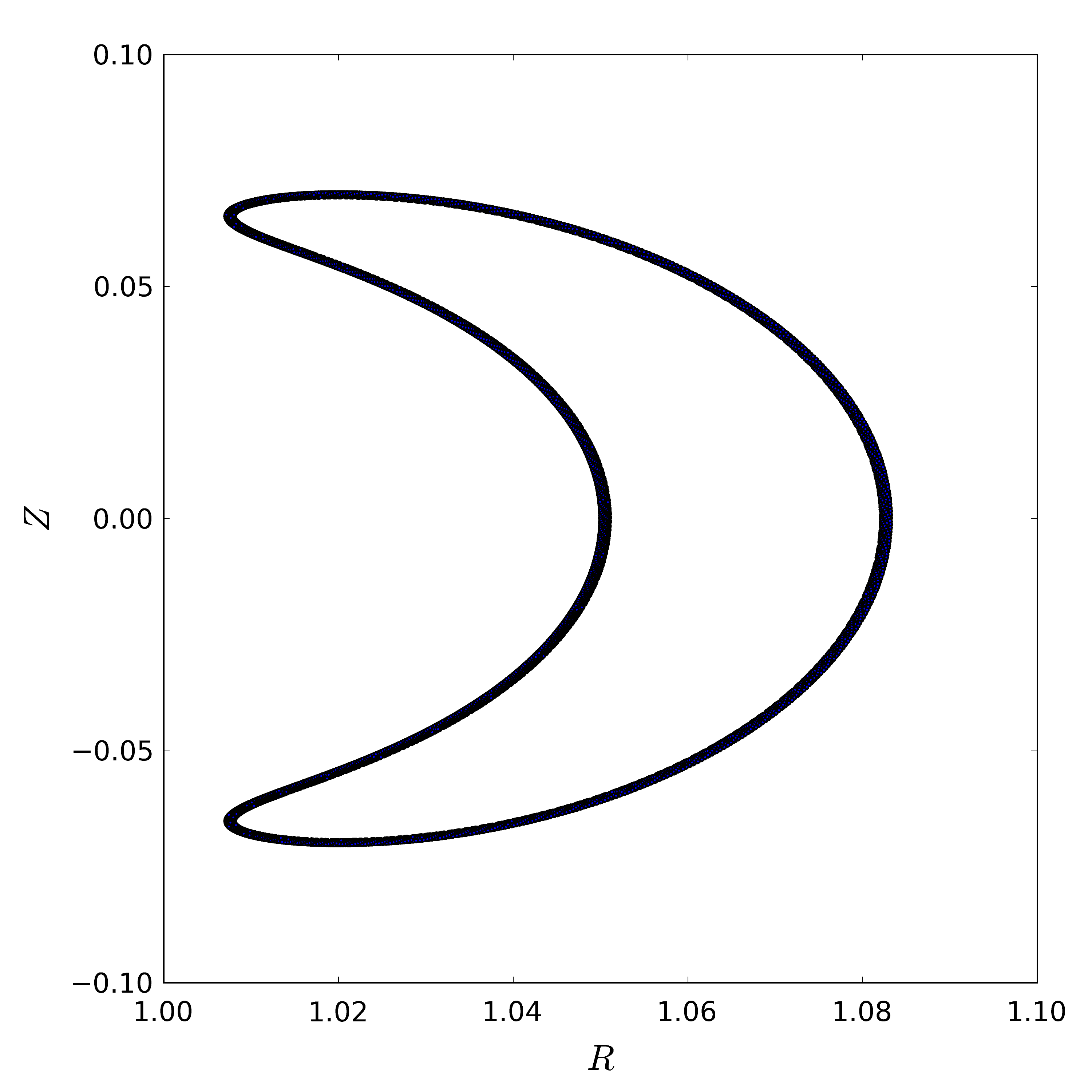}
}
\subfloat[Variational Midpoint]{
\includegraphics[width=.32\textwidth]{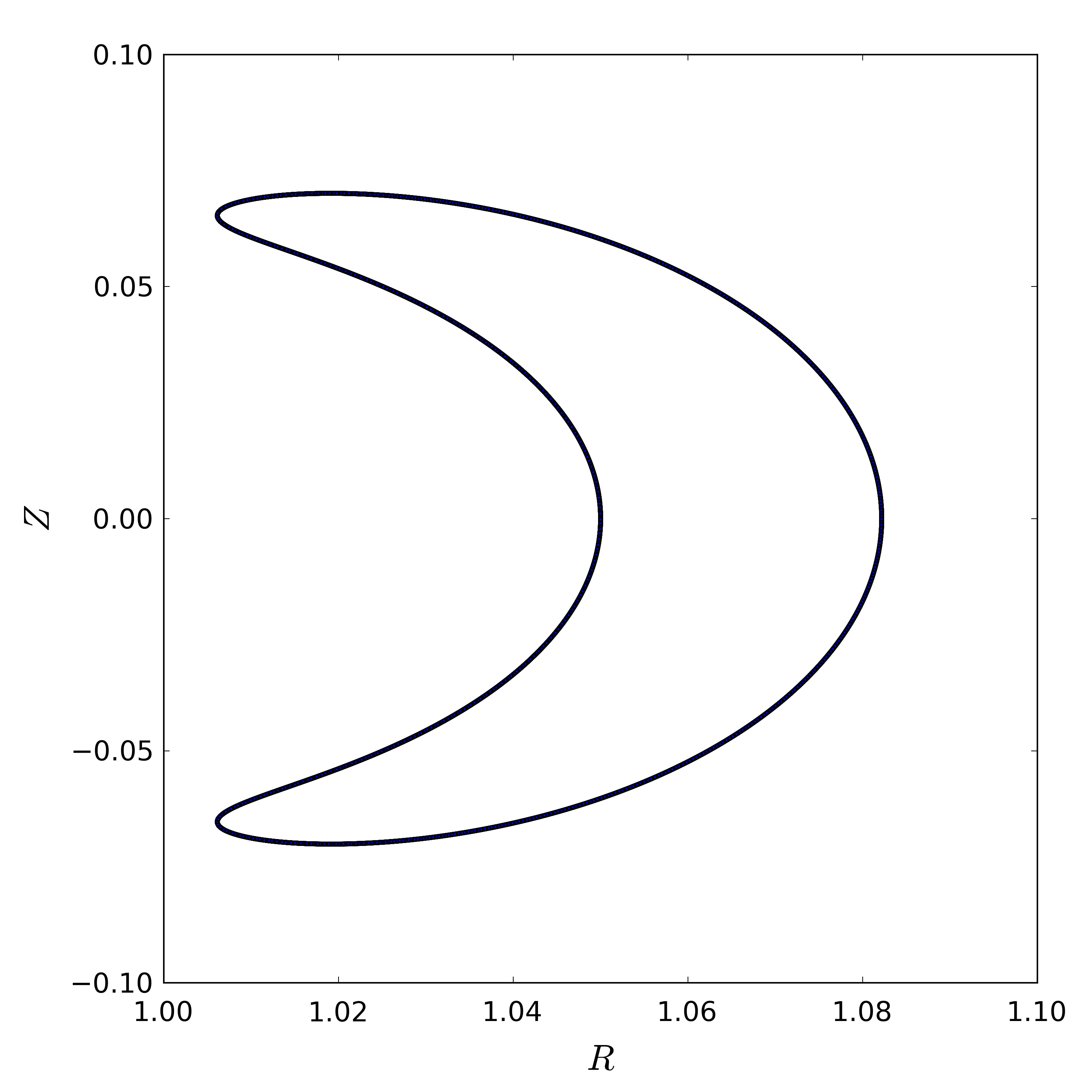}
}
\caption{Trapped particle after 200.000 bounce periods with 50 steps per bounce period.}
\label{fig:particles_4d_trapped_nb50_orbits}
\end{figure}

\begin{figure}[p]
\centering
\subfloat[Runge-Kutta]{
\includegraphics[width=.32\textwidth]{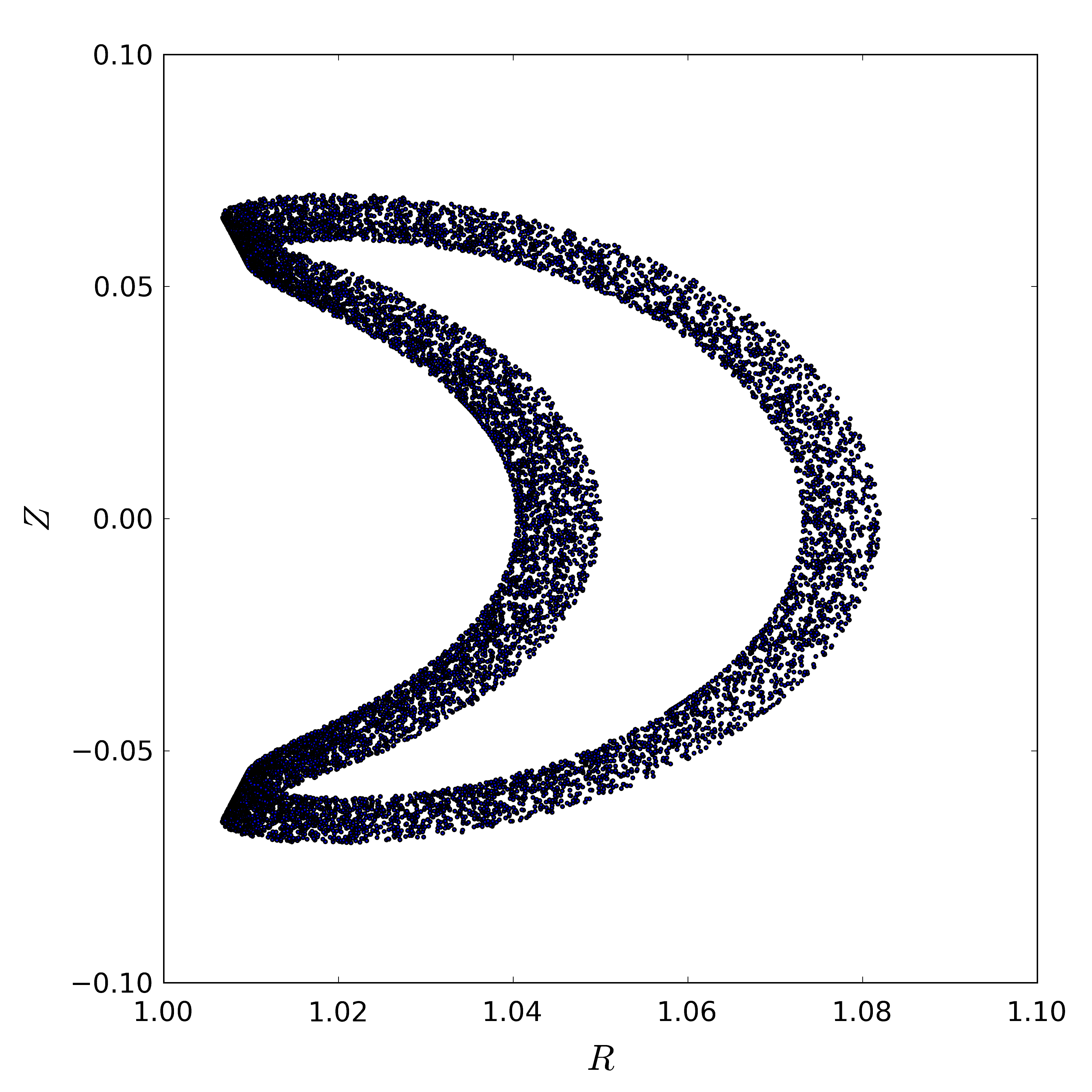}
}
\subfloat[Variational Trapezoidal]{
\includegraphics[width=.32\textwidth]{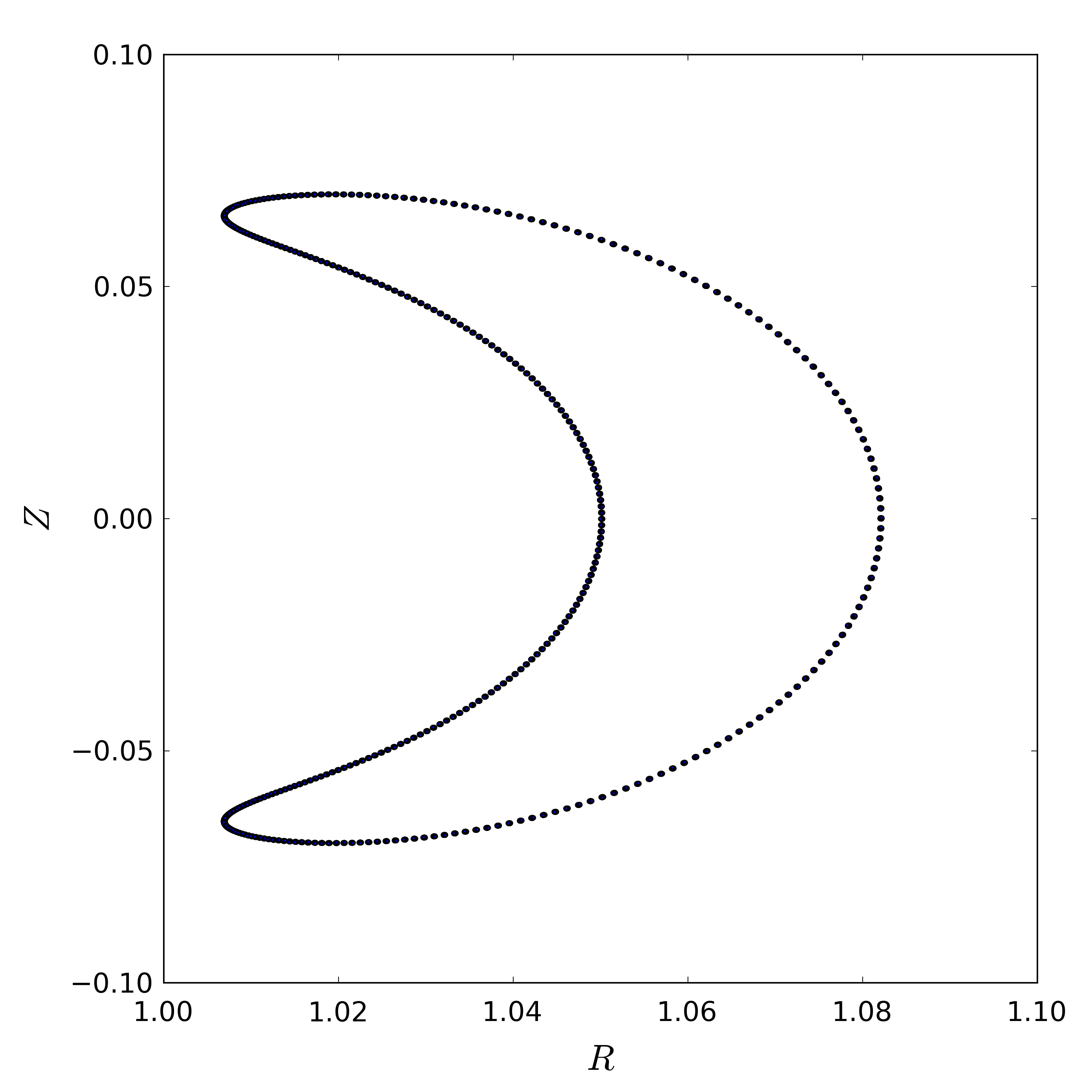}
}
\subfloat[Variational Midpoint]{
\includegraphics[width=.32\textwidth]{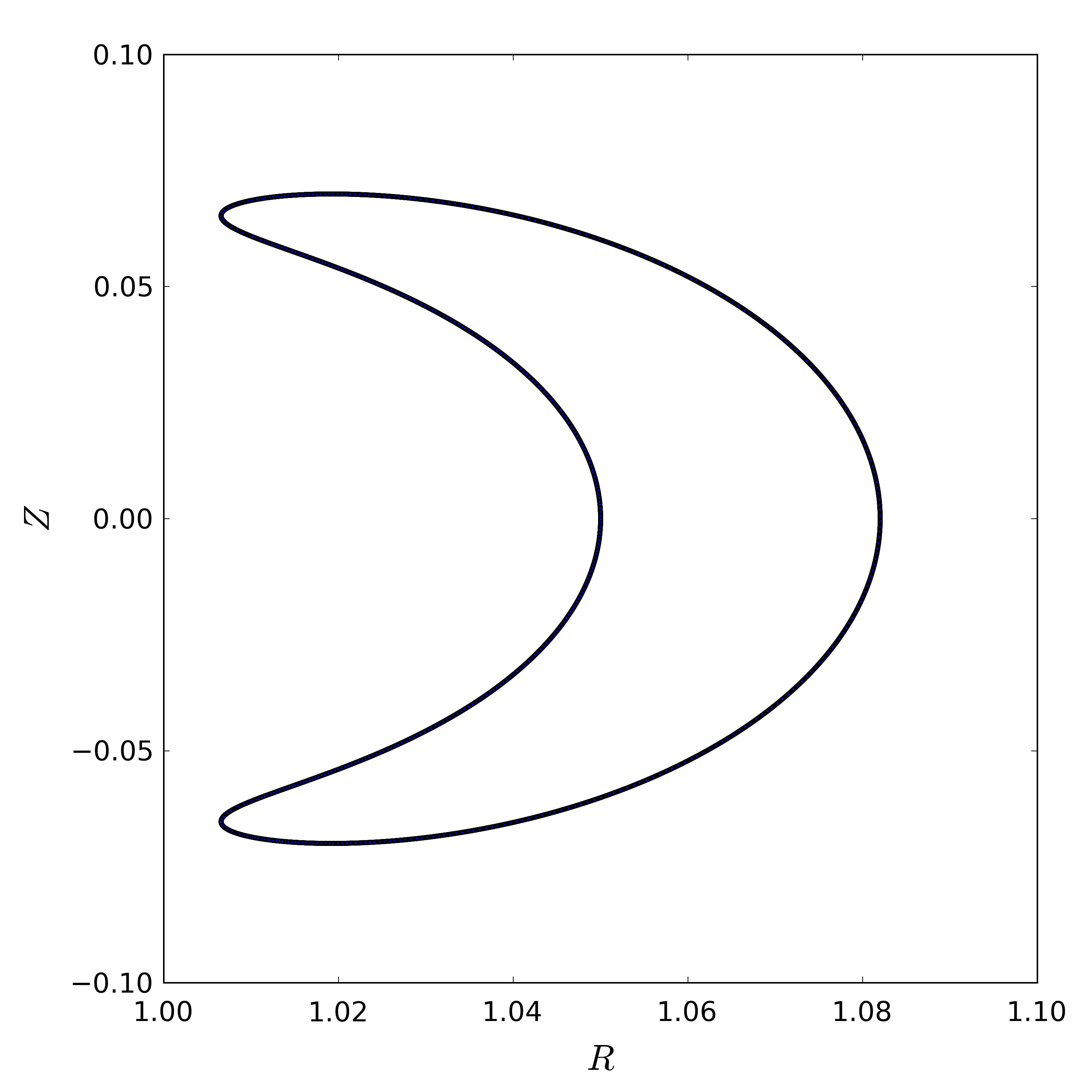}
}
\caption{Trapped particle after 250.000 bounce periods with 100 steps per bounce period.}
\label{fig:particles_4d_trapped_nb100_orbits}
\end{figure}

As in the two-dimensional case, we compare the two variational integrators with the explicit fourth order Runge-Kutta method with different timesteps corresponding to $25$, $50$ and $100$ steps per bounce period. The qualitative behaviour is the same as before. With the Runge-Kutta method, the particle orbits deviate severely from their expected shape while the variational integrators find the correct result (see figures \ref{fig:particles_4d_trapped_nb25_orbits} - \ref{fig:particles_4d_trapped_nb100_orbits}). The trapezoidal integrator is not stable for $25$ steps per bounce period and shows slight inaccuracies for $50$ steps per bounce period. The midpoint integrator is stable also for $25$ steps per bounce period and appears accurate already at this large timestep.

\begin{figure}[p]
\centering
\includegraphics[width=\textwidth]{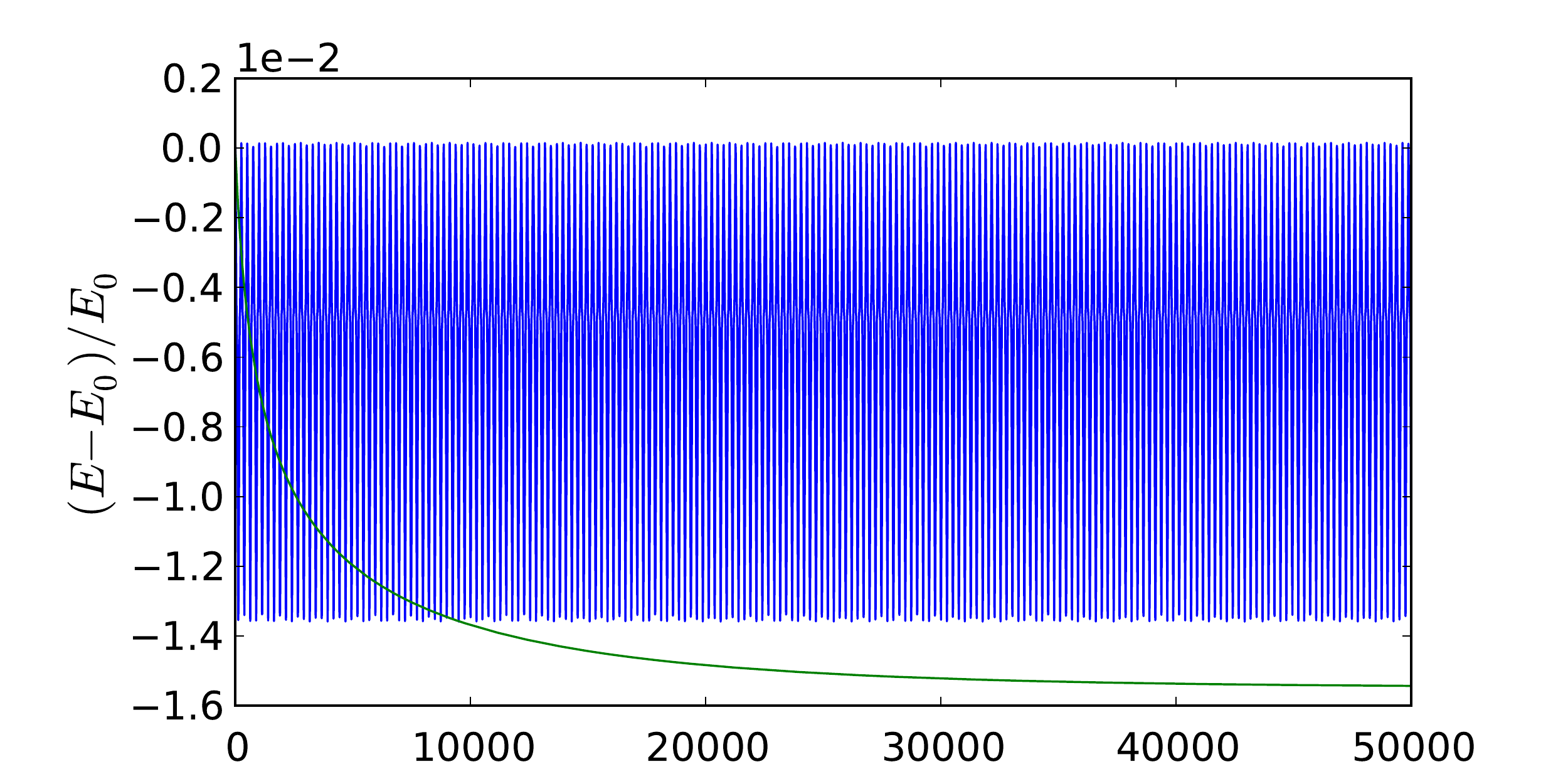}
\caption{Energy error for trapped particle with 25 steps per bounce period\\ (green: Runge-Kutta, blue: variational midpoint).}
\label{fig:particles_4d_trapped_nb25_energy}
\end{figure}

\begin{figure}[p]
\centering
\includegraphics[width=.9\textwidth]{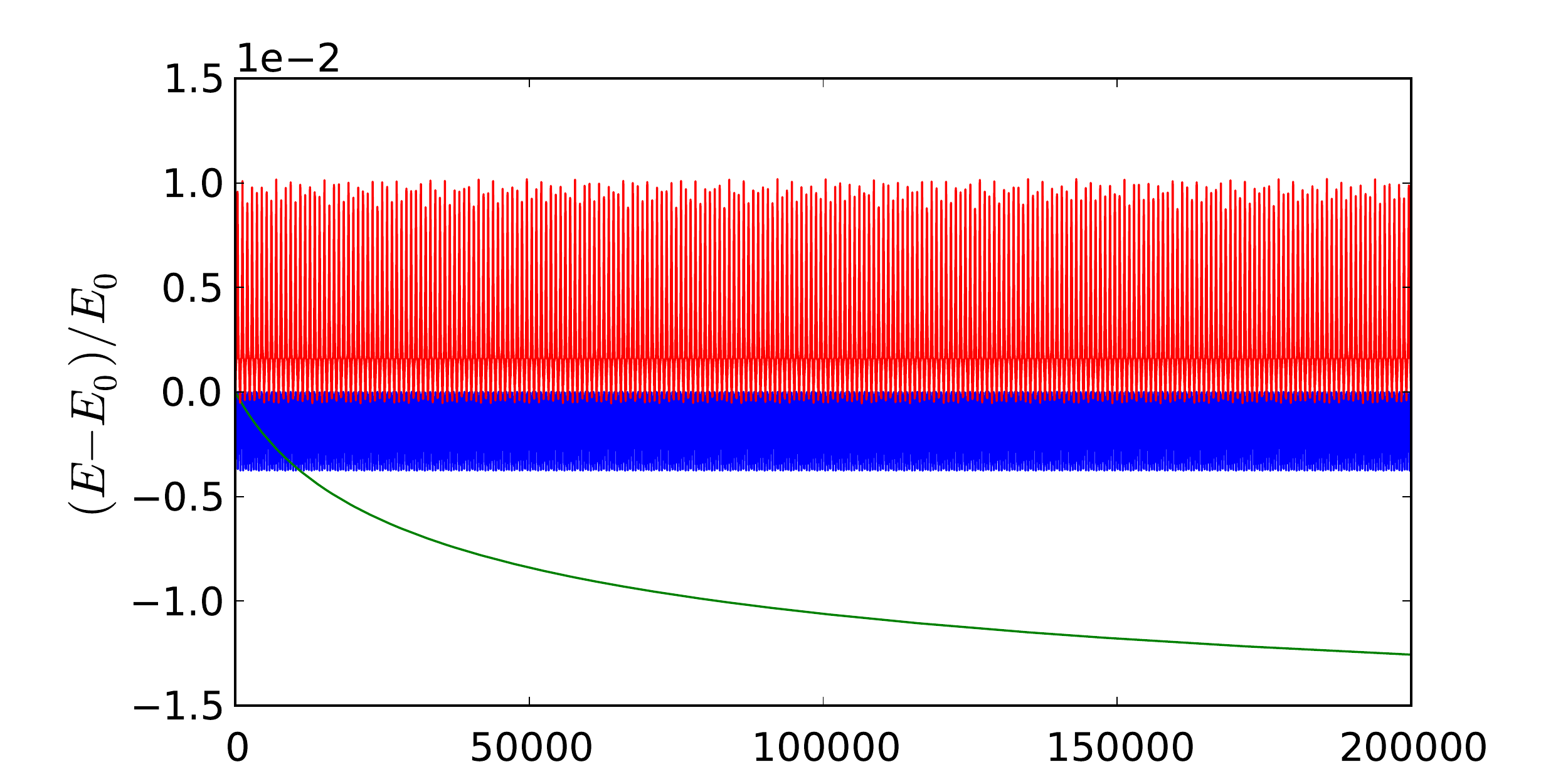}
\caption{Energy error for trapped particle with 50 timesteps per bounce period\\ (green: Runge-Kutta, blue: variational trapezoidal, red: variational midpoint).}
\label{fig:particles_4d_trapped_nb50_energy}
\end{figure}

While the variational integrators exhibit an oscillating energy error with a bounded amplitude of the oscillation, the Runge-Kutta method dissipates energy monotonically (see figures \ref{fig:particles_4d_trapped_nb25_energy}, \ref{fig:particles_4d_trapped_nb50_energy} and \ref{fig:particles_4d_trapped_nb25_energy}).
For the variational integrators, the amplitude of the error oscillation scales with the order of the scheme which is second order accurate.

\begin{figure}[p]
\centering
\includegraphics[width=.8\textwidth]{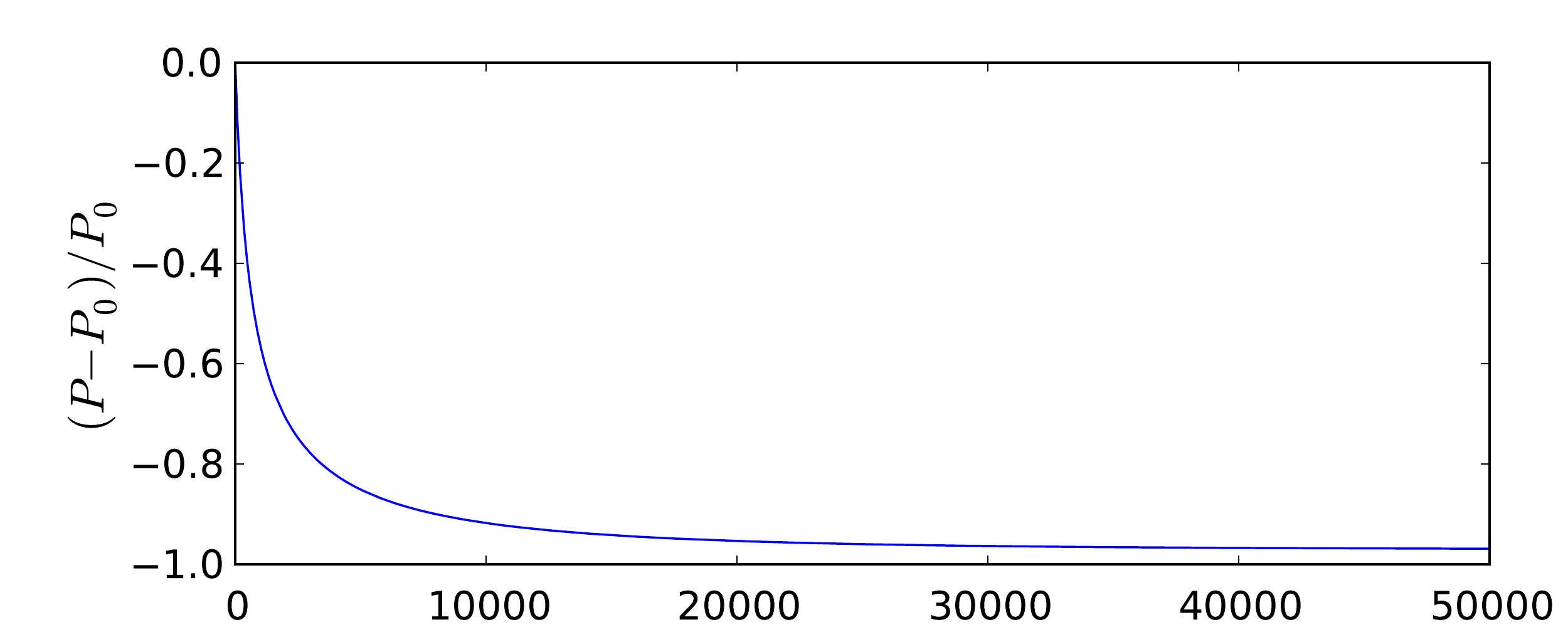}
\includegraphics[width=.8\textwidth]{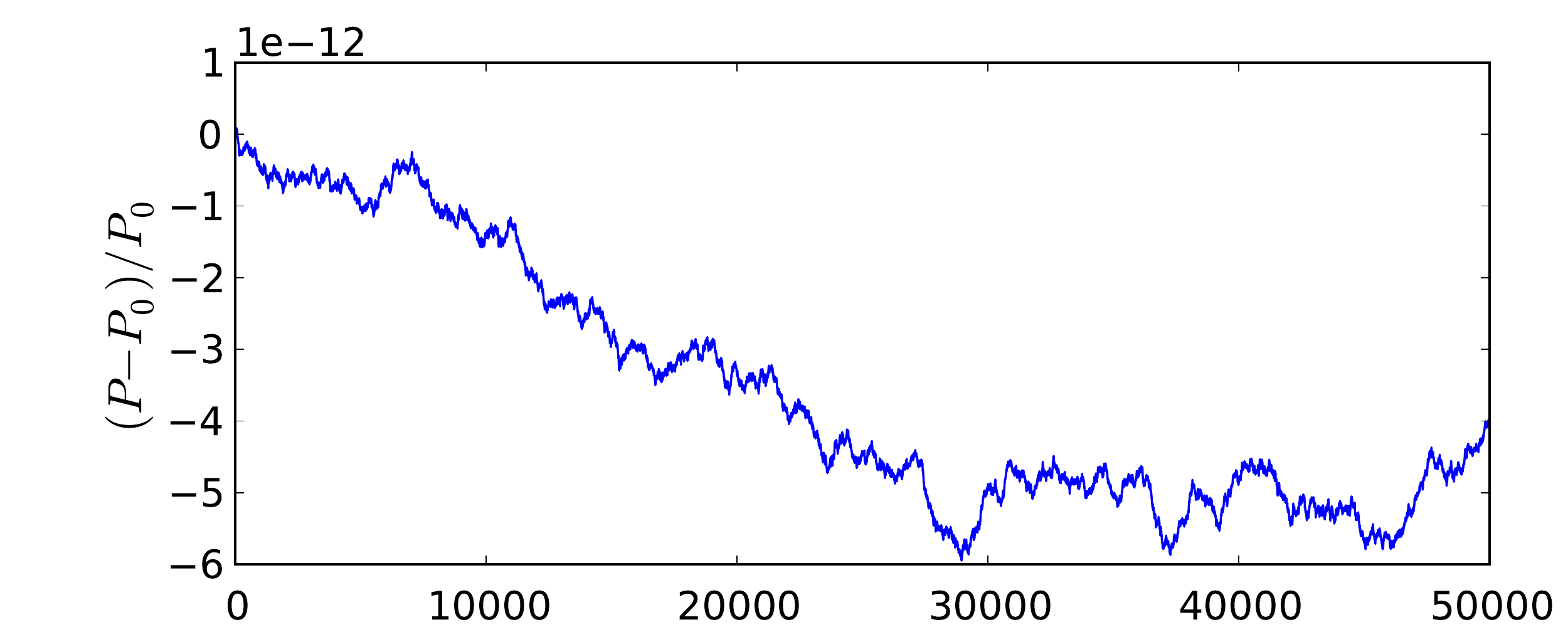}
\caption{Total linear momentum for a trapped particle with 25 steps per bounce period (top: Runge-Kutta, bottom: variational midpoint).}
\label{fig:particles_4d_trapped_nb25_momentum}
\end{figure}

\begin{figure}[p]
\centering
\includegraphics[width=.6\textwidth]{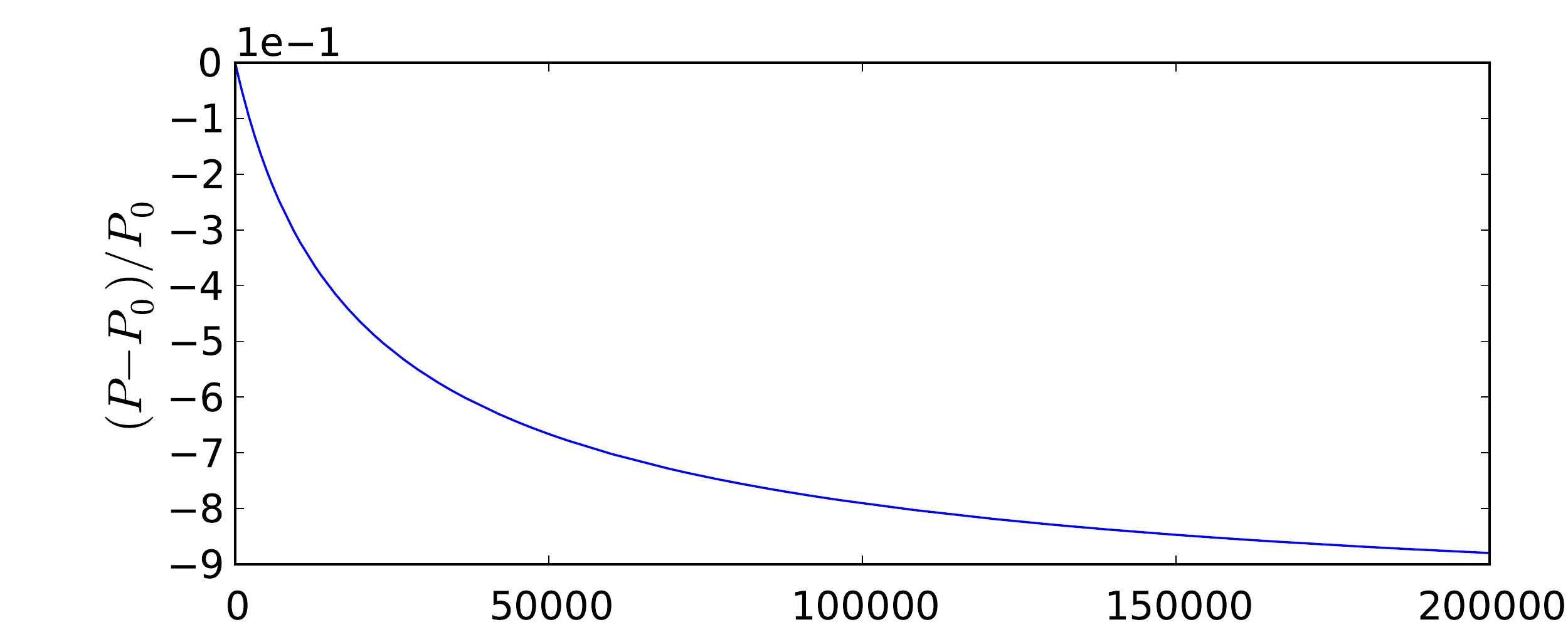}
\includegraphics[width=.6\textwidth]{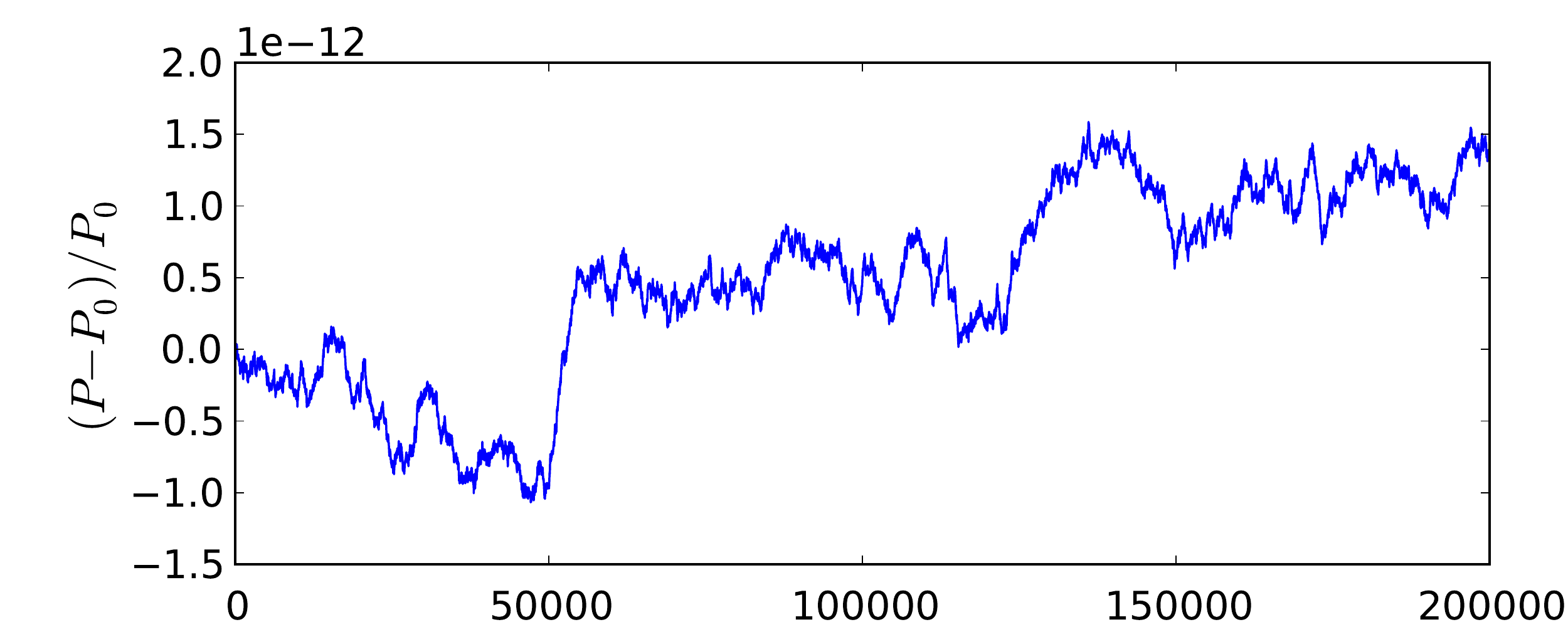}
\includegraphics[width=.6\textwidth]{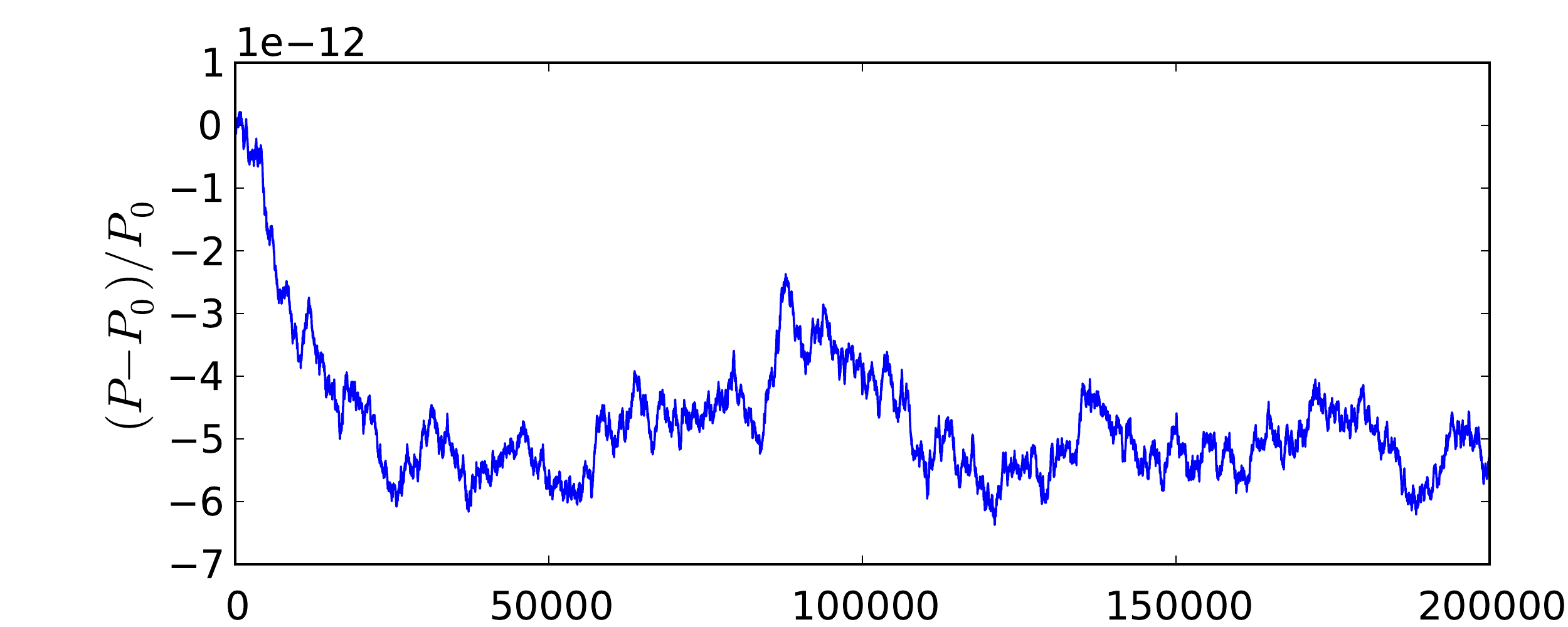}
\caption{Total linear momentum for a trapped particle with 50 timesteps per bounce period (top: Runge-Kutta, middle: variational trapezoidal, bottom: variational midpoint).}
\label{fig:particles_4d_trapped_nb50_momentum}
\end{figure}

\begin{figure}[p]
\centering
\includegraphics[width=.6\textwidth]{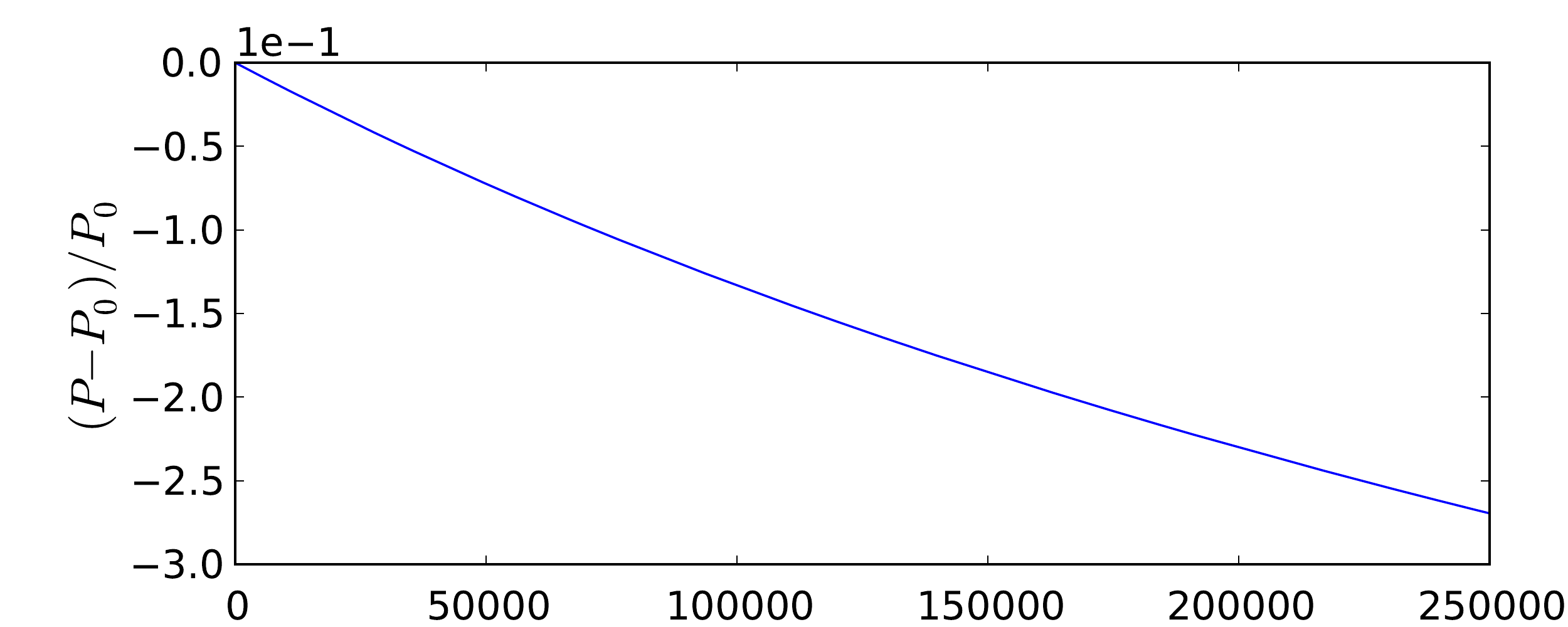}
\includegraphics[width=.6\textwidth]{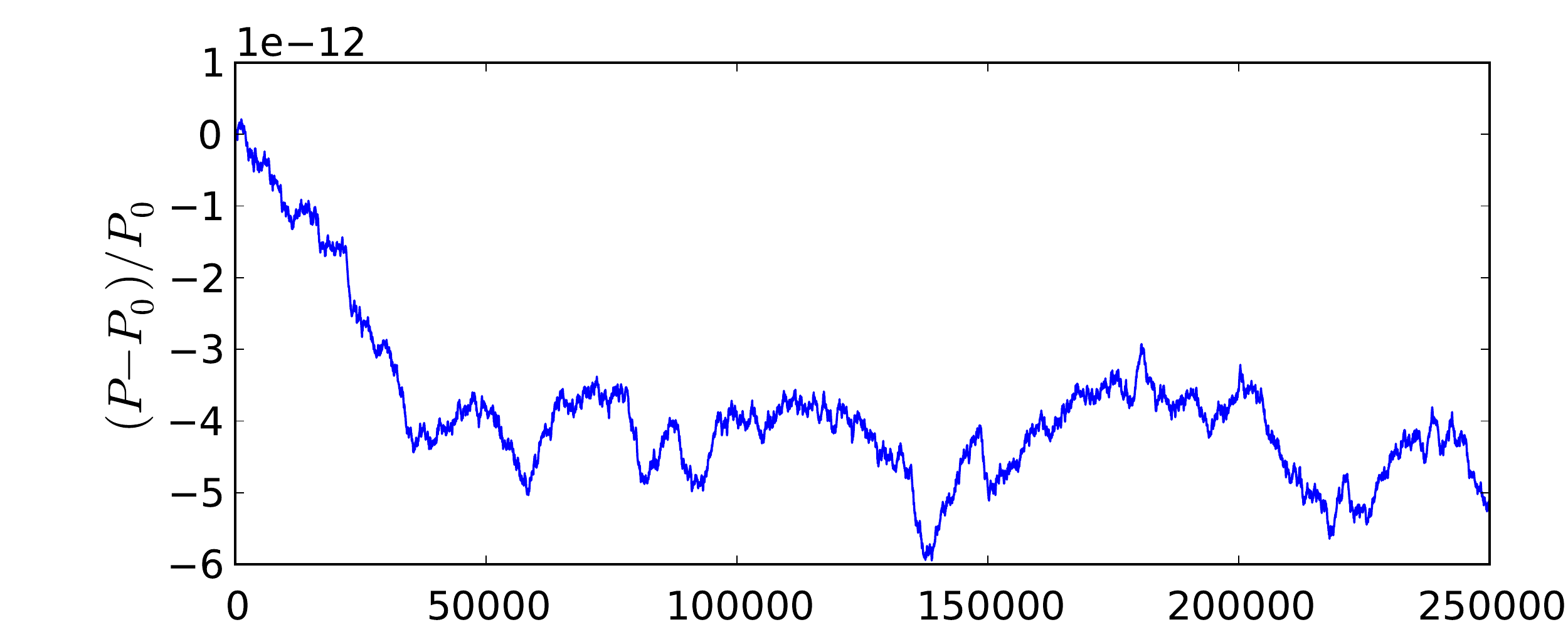}
\includegraphics[width=.6\textwidth]{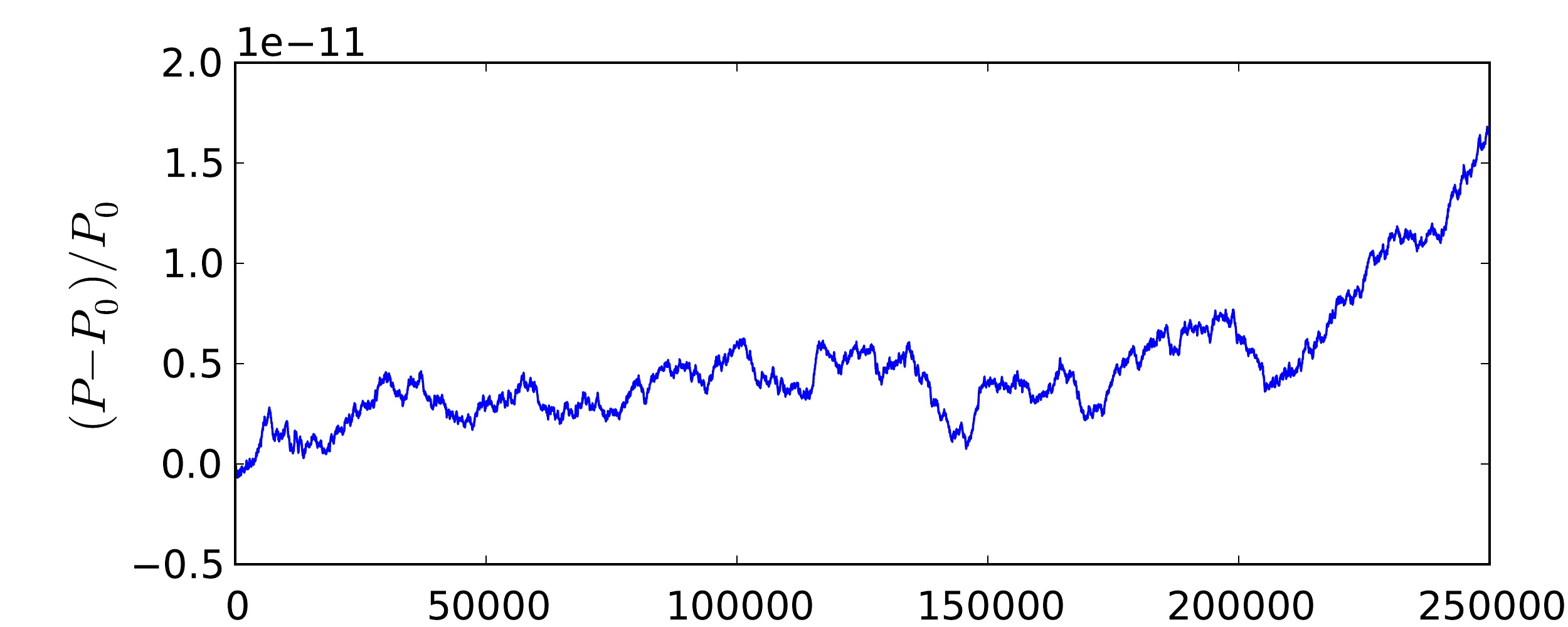}
\caption{Total linear momentum for a trapped particle with 100 timesteps per bounce period (top: Runge-Kutta, middle: variational trapezoidal, bottom: variational midpoint).}
\label{fig:particles_4d_trapped_nb100_momentum}
\end{figure}

The most interesting question about the four-dimensional integrators is that of toroidal momentum conservation.
As expected, the variational integrators exhibit only very small errors in the toroidal momentum, close to the machine accuracy (figures \ref{fig:particles_4d_trapped_nb25_momentum}, \ref{fig:particles_4d_trapped_nb50_momentum} and \ref{fig:particles_4d_trapped_nb100_momentum}) while the Runge-Kutta method dissipates the toroidal momentum as it appears monotonically. While the energy error of the Runge-Kutta method seems too approach a stationary value of order one percent or smaller, depending on the timestep, the toroidal momentum is dissipated almost completely during the course of the simulations, independently from the timestep.

The absolute value of the momentum error of the variational integrators seems to be mostly determined by the residual of the Newton iteration.
If the residual is too large, errors tend to accumulate, leading to an almost monotonic growths of the momentum error during the simulation. But even in that case, the toroidal momentum error is of orders $\mcal{O} (10^{-6}) ... \mcal{O} (10^{-3})$ for the simulation times considered here, and thus much smaller than with the Runge-Kutta method.

\begin{figure}[p]
\centering
\includegraphics[width=.9\textwidth]{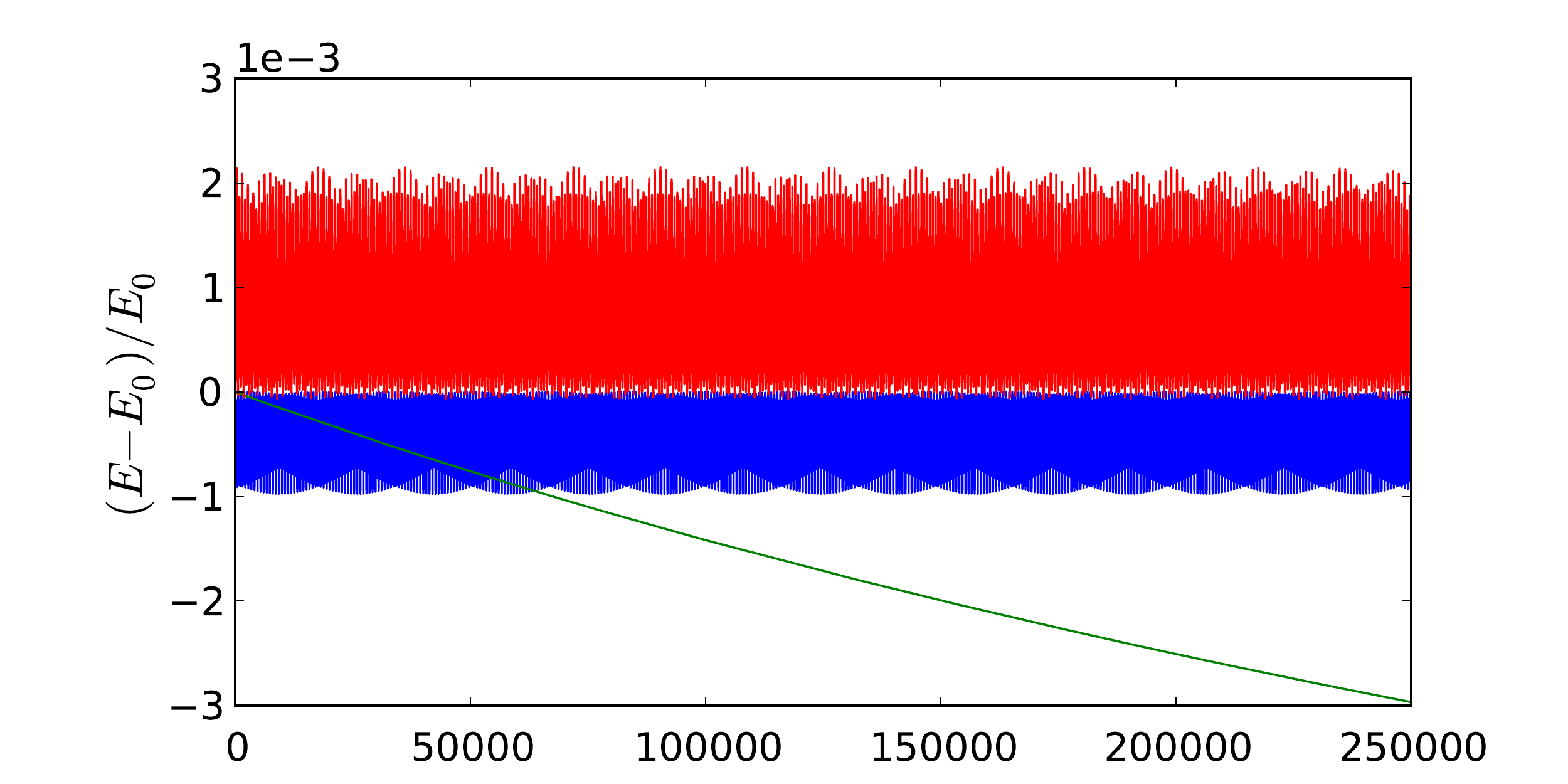}
\caption{Energy error for trapped particle with 100 timesteps per bounce period\\ (green: Runge-Kutta, blue: variational trapezoidal, red: variational midpoint).}
\label{fig:particles_4d_trapped_nb100_energy}
\end{figure}

\section{Variational PIC Scheme}

The reason for deriving variational integrators for guiding centre dynamics is of course not to compute particle orbits in the poloidal plane but the aim of finding better integration techniques for large scale particle-in-cell codes as they are used in the simulation of plasma turbulence.
In such codes, the electromagnetic fields are computed on a fixed grid, while the particles move in a mesh-free space.
To be self-consistent, the electromagnetic fields have to be computed given the particle positions.

As the field dynamics can also be derived from an action principle, it is possible to combine the particle Lagrangian and the field Lagrangian, together with an interaction term, to get a Lagrangian description of the full system. This can then be used to derive a variational integrator for the complete system of particles and fields, thereby not only conserving the symplectic form of each particle, but the multisymplectic form of the full system as well.

A similar idea has recently been presented by \citeauthor{Squire:2012} \cite{Squire:2012}. In that work, however, the electromagnetic fields are represented by discrete differential forms, which is the geometrically correct approach, but not necessary if the electrostatic potential alone is considered.
In the following we sketch the derivation of a variational PIC scheme based on the particle integrators from this chapter and a simple discretisation of the electrostatic potential. We do not compute the actual Euler-Lagrange equations as the aim of this section is merely to show the potential of this formulation and outline possible future directions of research.

\subsection{Total Lagrangian and Euler-Lagrange Equations}

If we restrict ourselves to the electrostatic case, the action reads
\begin{align}\label{eq:vi_finitec_action1}
\mcal{A} = \sum \limits_{p} \int \dot{x}_{p}^{2} \, dt + \dfrac{1}{2} \int \bigg( \dfrac{\partial \phi}{\partial x} \bigg)^{2} \, dt \, dx - \int \rho (t, x) \, \phi (t, x) \, dt \, dx
\end{align}

where $\rho$ is the charge density
\begin{align}\label{eq:vi_finitec_rho}
\rho (t, x) = \sum \limits_{s} q_{s} \int f_{s} (x,v,t) \, dv .
\end{align}

The distribution function of the species $s$ is computed as the sum of the distribution functions of each particle of that species
\begin{align}
f_{s} (t,x,v) = \sum \limits_{p} f_{p} (t,x,v) ,
\end{align}

and the particle distribution function $f_{p}$ is determined by
\begin{align}
f_{p} (t,x,v) = N_{p} \, S_{x} \big( x - x_{p} (t) \big) \, S_{v} \big( v - v_{p} (t) \big) ,
\end{align}

where $S_{x}$ and $S_{v}$ are the shape functions of the particle in space and velocity.
The simplest choice for $S_{v}$ are just $\delta$ functions
\begin{align}
S_{v} (v - v_{p}) = \delta (v^{x} - v^{x}_{p}) \, \delta (v^{y} - v^{y}_{p}) \, \delta (v^{z} - v^{z}_{p}) .
\end{align}

Smoother shape functions are often used for $S_{x}$, e.g., given by B-splines of order $l$,
\begin{align}
S_{x} (x - x_{p}) = \dfrac{1}{\Delta x_{p} \, \Delta y_{p} \, \Delta z_{p}} \, b_{l} \bigg( \dfrac{x - x_{p}}{\Delta x_{p}} \bigg) \, b_{l} \bigg( \dfrac{y - y_{p}}{\Delta y_{p}} \bigg) \, b_{l} \bigg( \dfrac{z - z_{p}}{\Delta z_{p}} \bigg) ,
\end{align}

where the $b_{l}$ are recursively defined as
\begin{subequations}\label{eq:vi_finitec_splines}
\begin{align}
b_{l} (\xi) &= \int \limits_{-\infty}^{+\infty} d\xi' \, b_{0} (\xi - \xi') \, b_{l-1} (\xi') \\
b_{0} (\xi) &=
\begin{cases}
1 & \text{if} \; \abs{\xi} < \tfrac{1}{2} , \\
0 & \text{otherwise} .
\end{cases}
\end{align}
\end{subequations}

Using these expressions, the action (\ref{eq:vi_finitec_action1}) becomes
\begin{align}\label{eq:vi_finitec_action2}
\mcal{A}
= \sum \limits_{p} \int \dot{x}_{p}^{2} \, dt
+ \dfrac{1}{2} \int \bigg( \dfrac{\partial \phi}{\partial x} \bigg)^{2} \, dt \, dx
- \sum \limits_{p} q_{p} \int S_{x} \big( x - x_{p} (t) \big) \, \phi (t, x) \, dt \, dx .
\end{align}

Computing the variation of the action with respect to $x_{p}$ gives the equations of motion for the particle
\begin{align}
\dfrac{\delta \mcal{A}}{\delta x_{p}} = \ddot{x}_{p} + q_{p} \int \dfrac{\partial S_{x}}{\partial x} (x - x_{p}) \, \phi(t, x) \, dx = 0 ,
\end{align}

and the variation with respect to $\phi$ gives the Poisson equation
\begin{align}
\dfrac{\partial \mcal{A}}{\partial \phi} = - \Delta \phi (t, x) - \sum \limits_{p} q_{p} \, S_{x} (x - x_{p}) .
\end{align}

This approach can now be used to obtain a fully variational discretisation of the system consisting of particles and fields.

\subsection{Variational Integrator}

The action (\ref{eq:vi_finitec_action2}) is discretised by a midpoint rule according to section (\ref{sec:vi_finite}) and (\ref{sec:vi_infinite}) as
\begin{align}
\mcal{A}_{d}
= \sum \limits_{k} \bigg[ & \dfrac{h_{t}}{2} \bigg( \dfrac{x_{p}^{k+1} - x_{p}^{k}}{h_{t}} \bigg)^{2} + \dfrac{h_{t}}{2} \sum \limits_{i} \bigg( \dfrac{\phi^{i+1, k} - \phi^{i,k}}{2 h_{x}} + \dfrac{\phi^{i+1, k+1} - \phi^{i, k+1}}{2 h_{x}} \bigg) \\
&- \dfrac{h_{t}}{16} \sum \limits_{i} \Big( \rho^{i,k} + \rho^{i+1,k} + \rho^{i+1,k+1} + \rho^{i,k+1} \Big) \Big( \phi^{i,k} + \phi^{i+1,k} + \phi^{i+1,k+1} + \phi^{i,k+1} \Big) \bigg]
\end{align}

where the discrete expression of $\rho$, i.e., the quadrature rule to compute (\ref{eq:vi_finitec_rho}), depends on the order of the splines (\ref{eq:vi_finitec_splines}) that are used.
For linear B-splines, a trapezoidal rule suffices, for quadratic B-splines, the Simpson rule should be used, and for cubic B-splines, Gauss' quadrature rule should be employed. The splines are integrated exactly by the corresponding quadrature rules.
The integration domain $V^{i,j}$ for $\rho^{i,k}$ in a two-dimensional setting is selected as depicted below.

\begin{figure}[H]
\centering
\begin{tikzpicture}[scale=0.7]%

\tikzstyle{every node}=[font=\small]

\draw [style=help lines, step=3]				(-4,-4)	grid		(+4,+4);

\draw [dashed, line width=0.5mm, fill=blue!20!white, fill opacity=0.5]	(-1.5,-1.5)	rectangle	(+1.5,+1.5);

\filldraw [color=black]	( 0,  0) circle (.1);

\filldraw [color=black]	(-1.5,    0) circle (.1);
\filldraw [color=black]	(+1.5,    0) circle (.1);
\filldraw [color=black]	(   0, -1.5) circle (.1);
\filldraw [color=black]	(   0, +1.5) circle (.1);

\draw [black]	(   0,    0)	node[anchor=south west]	{$(i, j)$};

\draw [black]	(-1.5,    0)	node[anchor=south east]	{$(i-1/2, j    )$};
\draw [black]	(+1.5,    0)	node[anchor=south west]	{$(i+1/2, j    )$};
\draw [black]	(   0, +1.5)	node[anchor=south west]	{$(i    , j-1/2)$};
\draw [black]	(   0, -1.5)	node[anchor=north west]	{$(i    , j+1/2)$};

\end{tikzpicture}
\end{figure}

So in one spatial dimension, $V_{i} = [ x^{i - 1/2} , x^{i+1/2} ]$, and the charge density becomes
\begin{align}
\rho^{i,k}
= \sum \limits_{p} \dfrac{q_{p}}{V^{i}} \int \limits_{V^{i}} S_{x} (x - x_{p}^{k}) \, dx .
\end{align}

Setting $\ohat{x} = x - x_{i}$ gives
\begin{align}
\rho^{i,k}
= \sum \limits_{p} \dfrac{q_{p}}{h_{x}} \int \limits_{-h_{x}/2}^{+h_{x}/2} S_{x} (x_{i} - x_{p}^{k} + \ohat{x}) \, d\ohat{x} .
\end{align}

For some of the quadrature rules the integration domain has to be $(-1,+1)$, therefore we introduce another transformation $\otilde{x} = 2 \ohat{x} / h_{x}$ with $d\ohat{x} = (h_{x} / 2) \, d\otilde{x}$, so that
\begin{align}
\rho^{i,k}
= \sum \limits_{p} \dfrac{q_{p}}{h_{x}} \int \limits_{-1}^{+1} S_{x} (x_{i} - x_{p}^{k} + h_{x} \, \otilde{x} / 2) \, d\otilde{x} .
\end{align}

We can now write the fully discrete expressions for the charge density, according to the different quadrature rules
\begin{subequations}
\begin{alignat}{2}
&\rho^{i,k} \big\vert_{\text{Trapezoidal}}
&&= \dfrac{1}{2} \sum \limits_{p} q_{p} \, \bigg[ S_{x} \bigg( x_{i} - x_{p}^{k} - \dfrac{h_{x}}{2} \bigg) + S_{x} \bigg( x_{i} - x_{p}^{k} + \dfrac{h_{x}}{2} \bigg) \bigg] ,
\\
&\rho^{i,k} \big\vert_{\text{Simpson}}
&&= \dfrac{1}{6} \sum \limits_{p} q_{p} \, \bigg[ S_{x} \bigg( x_{i} - x_{p}^{k} - \dfrac{h_{x}}{2} \bigg) + 4 \, S_{x} \bigg( x_{i} - x_{p}^{k} \bigg) + S_{x} \bigg( x_{i} - x_{p}^{k} + \dfrac{h_{x}}{2} \bigg) \bigg] ,
\\
&\rho^{i,k} \big\vert_{\text{Gauss}}
&&= \dfrac{1}{2} \sum \limits_{p} q_{p} \, \bigg[ S_{x} \bigg( x_{i} - x_{p} - \dfrac{h_{x}}{2 \sqrt{3}} \bigg) + S_{x} \bigg( x_{i} - x_{p} + \dfrac{h_{x}}{2 \sqrt{3}} \bigg) \bigg] .
\end{alignat}
\end{subequations}

Now we have all the necessary ingredients for the derivation of a fully variational PIC scheme.
The last step is to decide on the order of the B-splines and to actually compute the discrete Euler-Lagrange equations.
This, however, is left for future research.

\begin{subappendices}

\section{Calculation of Transit and Bounce Times}

In the calculation of the bounce time, we follow \citeauthor{Brizard:2011} \cite{Brizard:2011}.
We merely collect the necessary formulae. For details on the actual derivation we refer to \citeauthor{Brizard:2011}'s original work.
We assume a circular tokamak with large aspect ratio, such that
\begin{align}
\eps &\equiv \dfrac{r}{R} \ll 1 &
& \text{(inverse aspect ratio)} . &
\end{align}

For a passing particle, the transit frequency is given by
\begin{align}
\omega_{t} = \dfrac{\pi \, \omega_{\parallel} \, \sqrt{\kappa}}{\msf{K} (\kappa^{-1})}
\end{align}

while for a trapped particle, the bounce frequency is
\begin{align}
\omega_{b} = \dfrac{\pi \, \omega_{\parallel}}{2 \, \msf{K} (\kappa)} .
\end{align}

Here, $\omega_{\parallel}$ is the characteristic parallel frequency, defined as
\begin{align}
\omega_{\parallel} = \dfrac{1}{qR} \sqrt{\eps \, \dfrac{\mu B_{0}}{m	}} ,
\end{align}

$\kappa$ is the bounce-transit parameter, defined as
\begin{align}
\kappa \equiv \dfrac{E - \mu B_{0} (1-\eps)}{2 \eps \, \mu B_{0}} ,
\end{align}

and $\msf{K}$ is the complete elliptic integral of first kind.
Furthermore, $m$ is the mass of the particle, $E = \tfrac{1}{2} u^{2} + \mu B$ its energy, $\mu$ the magnetic moment, and $B$ the magnetic field.
The transit and bounce times are accordingly computed as
\begin{align}
\tau_{t} &= \dfrac{2\pi}{\omega_{t}} = \dfrac{2 \, \msf{K} (\kappa^{-1})}{\omega_{\parallel} \, \sqrt{\kappa}} , &
\tau_{b} &= \dfrac{2\pi}{\omega_{b}} = \dfrac{4 \, \msf{K} (\kappa)     }{\omega_{\parallel}} . &
\end{align}

For the trapped particle from the examples in this chapter, the bounce time is computed to be $\tau_{b} = 43107$.

\section{Jacobians}\label{sec:vi_finite_jacobians}

\subsubsection{2D Trapezoidal Method}

The Jacobian is defined as
\begin{align}
\mcal{J}^{\text{2DTR}}
=
\dfrac{1}{2}
\begin{pmatrix}
J_{11} & J_{12} \\
J_{21} & J_{22}
\end{pmatrix} ,
\end{align}

with components
\begin{subequations}
\begin{align}
J_{11} &= A^{*}_{R,R} (q_{k}) - A^{*}_{R,R} (q_{k+1}^{n}) - \tfrac{1}{2} \, h \, u_{,R} (q_{k}) \, u_{,R} (q_{k+1}^{n}) , \\
J_{12} &= A^{*}_{Z,R} (q_{k}) - A^{*}_{R,Z} (q_{k+1}^{n}) - \tfrac{1}{2} \, h \, u_{,R} (q_{k}) \, u_{,Z} (q_{k+1}^{n}) , \\
J_{21} &= A^{*}_{R,Z} (q_{k}) - A^{*}_{Z,R} (q_{k+1}^{n}) - \tfrac{1}{2} \, h \, u_{,Z} (q_{k}) \, u_{,R} (q_{k+1}^{n}) , \\
J_{22} &= A^{*}_{Z,Z} (q_{k}) - A^{*}_{Z,Z} (q_{k+1}^{n}) - \tfrac{1}{2} \, h \, u_{,Z} (q_{k}) \, u_{,Z} (q_{k+1}^{n}) .
\end{align}
\end{subequations}

\subsubsection{2D Midpoint Method}

The Jacobian is defined as
\begin{align}
\mcal{J}^{\text{2DMP}}
=
\dfrac{1}{4}
\begin{pmatrix}
J_{11} & J_{12} \\
J_{21} & J_{22}
\end{pmatrix} ,
\end{align}

with components
\begin{subequations}
\begin{align}
J_{11}
\nonumber
&=
2 \, \Big[ A^{*}_{R,R} (q_{k+1/2}^{n}) - A^{*}_{R,R} (q_{k+1/2}^{n}) \Big]
+ A^{*}_{R,RR} (q_{k+1/2}) \, \Big[ q^{R}_{k+1} - q^{R}_{k} \Big]
+ A^{*}_{Z,RR} (q_{k+1/2}) \, \Big[ q^{Z}_{k+1} - q^{Z}_{k} \Big] \\
& \hspace{2em}
- h \, u_{,R} (q_{k+1/2}) \, u_{,R} (q_{k+1/2})
- h \, u (q_{k+1/2}) \, u_{,RR} (q_{k+1/2})
- h \, \mu B_{,RR} (q_{k+1/2})
, \\
J_{12}
\nonumber
&=
2 \, \Big[ A^{*}_{Z,R} (q_{k+1/2}^{n}) - A^{*}_{R,Z} (q_{k+1/2}^{n}) \Big] \\
\nonumber
& \hspace{2em}
+ A^{*}_{R,RZ} (q_{k+1/2}) \, \Big[ q^{R}_{k+1} - q^{R}_{k} \Big]
+ A^{*}_{Z,RZ} (q_{k+1/2}) \, \Big[ q^{Z}_{k+1} - q^{Z}_{k} \Big] \\
& \hspace{4em}
- h \, u_{,R} (q_{k+1/2}) \, u_{,Z} (q_{k+1/2})
- h \, u (q_{k+1/2}) \, u_{,RZ} (q_{k+1/2})
- h \, \mu B_{,RZ} (q_{k+1/2})
, \\
J_{21}
\nonumber
&=
2 \, \Big[ A^{*}_{R,Z} (q_{k+1/2}) - A^{*}_{Z,R} (q_{k+1/2}) \Big] \\
\nonumber
& \hspace{2em}
+ A^{*}_{R,RZ} (q_{k+1/2}) \, \Big[ q^{R}_{k+1} - q^{R}_{k} \Big]
+ A^{*}_{Z,RZ} (q_{k+1/2}) \, \Big[ q^{Z}_{k+1} - q^{Z}_{k} \Big] \\
& \hspace{4em}
- h \, u_{,R} (q_{k+1/2}) \, u_{,Z} (q_{k+1/2})
- h \, u (q_{k+1/2}) \, u_{,RZ} (q_{k+1/2})
- h \, \mu B_{,RZ} (q_{k+1/2})
, \\
J_{22}
\nonumber
&=
2 \, \Big[ A^{*}_{Z,Z} (q_{k+1/2}) - A^{*}_{Z,Z} (q_{k+1/2}) \Big]
+ A^{*}_{R,ZZ} (q_{k+1/2}) \, \Big[ q^{R}_{k+1} - q^{R}_{k} \Big]
+ A^{*}_{Z,ZZ} (q_{k+1/2}) \, \Big[ q^{Z}_{k+1} - q^{Z}_{k} \Big] \\
& \hspace{2em}
- h \, u_{,Z} (q_{k+1/2}) \, u_{,Z} (q_{k+1/2})
- h \, u (q_{k+1/2}) \, u_{,ZZ} (q_{k+1/2})
- h \, \mu B_{,ZZ} (q_{k+1/2})
.
\end{align}
\end{subequations}

\subsubsection{4D Trapezoidal Method}

The Jacobian is defined as
\begin{align}
\mcal{J}^{\text{4DTR}}
=
\dfrac{1}{2}
\begin{pmatrix}
J_{11} & J_{12} & J_{13} & J_{14} \\
J_{21} & J_{22} & J_{23} & J_{24} \\
J_{31} & J_{32} & J_{33} & J_{34} \\
J_{41} & J_{42} & J_{43} & J_{44}
\end{pmatrix} ,
\end{align}

with components
\begin{subequations}
\begin{align}
J_{11} &= A^{*}_{R,R} (q_{k}) - A^{*}_{R,R} (q_{k+1}^{n}) , \\
J_{12} &= A^{*}_{Z,R} (q_{k}) - A^{*}_{R,Z} (q_{k+1}^{n}) , \\
J_{13} &= q^{R}_{k} \, A^{*}_{\phy,R} (q_{k}) + A_{\phy}^{*} (q_{k}) , \\
J_{14} &= - b_{R} (q_{k+1}^{n}) ,
\end{align}
\begin{align}
J_{21} &= A^{*}_{R,Z} (q_{k}) - A^{*}_{Z,R} (q_{k+1}^{n}) , \\
J_{22} &= A^{*}_{Z,Z} (q_{k}) - A^{*}_{Z,Z} (q_{k+1}^{n}) , \\
J_{23} &= q^{R}_{k} \, A^{*}_{\phy,Z} (q_{k}) , \\
J_{24} &= - b_{Z} (q_{k+1}^{n}) ,
\end{align}
\begin{align}
J_{31} &= - q^{R}_{k+1} \, A^{*}_{\phy,R} (q_{k+1}^{n})
- A_{\phy}^{*} (q_{k+1}^{n}) , \\
J_{32} &= - q^{R}_{k+1} \, A^{*}_{\phy,Z} (q_{k+1}^{n}) , \\
J_{33} &= 0 , \\
J_{34} &= - q^{R}_{k+1} \, b_{\phy} (q_{k+1}^{n}) ,
\end{align}
\begin{align}
J_{41} &= b_{R} (q_{k}) , \\
J_{42} &= b_{Z} (q_{k}) , \\
J_{43} &= q^{R}_{k} \, b_{\phy} (q_{k}) , \\
J_{44} &= - \tfrac{1}{2} \, h
.
\end{align}
\end{subequations}

\subsubsection{4D Midpoint Method}

The Jacobian is defined as
\begin{align}
\mcal{J}^{\text{4DMP}}
=
\dfrac{1}{4}
\begin{pmatrix}
J_{11} & J_{12} & J_{13} & J_{14} \\
J_{21} & J_{22} & J_{23} & J_{24} \\
J_{31} & J_{32} & J_{33} & J_{34} \\
J_{41} & J_{42} & J_{43} & J_{44}
\end{pmatrix} ,
\end{align}

with components
\begin{subequations}
\begin{align}
\nonumber
J_{11} &= 2 \, \big[ A^{*}_{R,R} (q_{k+1/2}) - A^{*}_{R,R} (q_{k+1/2}) \big]
+ 2 \, A^{*}_{\phy,R} (q_{k+1/2}) \, \Big[ q^{\phy}_{k+1} - q^{\phy}_{k} \Big] \\
\nonumber
& \hspace{4em}
+ A^{*}_{R,RR} (q_{k+1/2}) \, \Big[ q^{R}_{k+1} - q^{R}_{k} \Big]
+ A^{*}_{Z,RR} (q_{k+1/2}) \, \Big[ q^{Z}_{k+1} - q^{Z}_{k} \Big] \\
& \hspace{4em}
+ q^{R}_{k+1/2} \, A^{*}_{\phy,RR} (q_{k+1/2}) \, \Big[ q^{\phy}_{k+1} - q^{\phy}_{k} \Big]
- h \, \mu B_{,RR} (q_{k+1/2}) ,
\\
\nonumber
J_{12} &= 2 \, \big[ A^{*}_{Z,R} (q_{k+1/2}) - A^{*}_{R,Z} (q_{k+1/2}) \big]
+ A^{*}_{\phy,Z} (q_{k+1/2}) \, \Big[ q^{\phy}_{k+1} - q^{\phy}_{k} \Big] \\
\nonumber
& \hspace{4em}
+ A^{*}_{R,RZ} (q_{k+1/2}) \, \Big[ q^{R}_{k+1} - q^{R}_{k} \Big]
+ A^{*}_{Z,RZ} (q_{k+1/2}) \, \Big[ q^{Z}_{k+1} - q^{Z}_{k} \Big] \\
& \hspace{4em}
+ q^{R}_{k+1/2} \, A^{*}_{\phy,RZ} (q_{k+1/2}) \, \Big[ q^{\phy}_{k+1} - q^{\phy}_{k} \Big]
- h \, \mu B_{,RZ} (q_{k+1/2}) ,
\\
J_{13} &= 2 \, q^{R}_{k+1/2} \, A^{*}_{\phy,R} (q_{k+1/2}) + 2 \, A_{\phy}^{*} (q_{k}) ,
\\
\nonumber
J_{14} &= b_{R,R} (q_{k+1/2}) \, \Big[ q^{R}_{k+1} - q^{R}_{k} \Big]
+ b_{Z,R} (q_{k+1/2}) \, \Big[ q^{Z}_{k+1} - q^{Z}_{k} \Big]
+ q^{R}_{k+1/2} \, b_{\phy,R} (q_{k+1/2}) \, \Big[ q^{\phy}_{k+1} - q^{\phy}_{k} \Big] \\
& \hspace{4em}
+ b_{\phy} (q_{k+1/2}) \, \Big[ q^{\phy}_{k+1} - q^{\phy}_{k} \Big]
- 2 \, b_{R} (q_{k+1/2}) ,
\end{align}
\begin{align}
\nonumber
J_{21} &= 2 \, \Big[ A^{*}_{R,Z} (q_{k+1/2}) - A^{*}_{Z,R} (q_{k+1/2}) \Big]
+ A^{*}_{\phy,Z} (q_{k+1/2}) \, \Big[ q^{\phy}_{k+1} - q^{\phy}_{k} \Big] \\
\nonumber
& \hspace{4em}
+ A^{*}_{R,RZ} (q_{k+1/2}) \, \Big[ q^{R}_{k+1} - q^{R}_{k} \Big]
+ A^{*}_{Z,RZ} (q_{k+1/2}) \, \Big[ q^{Z}_{k+1} - q^{Z}_{k} \Big] \\
& \hspace{4em}
+ q^{R}_{k+1/2} \, A^{*}_{\phy,RZ} (q_{k+1/2}) \, \Big[ q^{\phy}_{k+1} - q^{\phy}_{k} \Big]
- h \, \mu B_{,RZ} (q_{k+1/2}) ,
\\
\nonumber
J_{22} &= 2 \, \Big[ A^{*}_{Z,Z} (q_{k+1/2}) - A^{*}_{Z,Z} (q_{k+1/2}) \Big] \\
\nonumber
& \hspace{4em}
+ A^{*}_{R,ZZ} (q_{k+1/2}) \, \Big[ q^{R}_{k+1} - q^{R}_{k} \Big]
+ A^{*}_{Z,ZZ} (q_{k+1/2}) \, \Big[ q^{Z}_{k+1} - q^{Z}_{k} \Big] \\
& \hspace{4em}
+ q^{R}_{k+1/2} \, A^{*}_{\phy,ZZ} (q_{k+1/2}) \, \Big[ q^{\phy}_{k+1} - q^{\phy}_{k} \Big]
- h \, \mu B_{,ZZ} (q_{k+1/2}) ,
\\
J_{23} &= 2 \, q^{R}_{k+1/2} \, A^{*}_{\phy,Z} (q_{k+1/2}) ,
\\
\nonumber
J_{24} &= b_{R,Z} (q_{k+1/2}) \, \Big[ q^{R}_{k+1} - q^{R}_{k} \Big]
+ b_{Z,Z} (q_{k+1/2}) \, \Big[ q^{Z}_{k+1} - q^{Z}_{k} \Big]
+ q^{R}_{k+1/2} \, b_{\phy,Z} (q_{k+1/2}) \, \Big[ q^{\phy}_{k+1} - q^{\phy}_{k} \Big] \\
& \hspace{4em}
- 2 \, b_{Z} (q_{k+1/2}) ,
\end{align}
\begin{align}
J_{31} &= - 2 \, q^{R}_{k+1/2} \, A^{*}_{\phy,R} (q_{k+1/2}) - 2 \, A_{\phy}^{*} (q_{k+1/2}) , \\
J_{32} &= - 2 \, q^{R}_{k+1/2} \, A^{*}_{\phy,Z} (q_{k+1/2}) , \\
J_{33} &= 0 , \\
J_{34} &= - 2 \, q^{R}_{k+1/2} \, b_{\phy} (q_{k+1/2}) ,
\end{align}
\begin{align}
\nonumber
J_{41} &= 2 \, b_{R} (q_{k+1/2})
+ b_{\phy} (q_{k+1/2}) \, \Big[ q^{\phy}_{k+1} - q^{\phy}_{k} \Big]
+ b_{R,R} (q_{k+1/2}) \, \Big[ q^{R}_{k+1} - q^{R}_{k} \Big] \\
& \hspace{4em}
+ b_{Z,R} (q_{k+1/2}) \, \Big[ q^{Z}_{k+1} - q^{Z}_{k} \Big]
+ q^{R}_{k+1/2} \, b_{\phy,R} (q_{k+1/2}) \, \Big[ q^{\phy}_{k+1} - q^{\phy}_{k} \Big] ,
\\
\nonumber
J_{42} &= 2 \, b_{Z} (q_{k+1/2})
+ b_{R,Z} (q_{k+1/2}) \, \Big[ q^{R}_{k+1} - q^{R}_{k} \Big] \\
& \hspace{4em}
+ b_{Z,Z} (q_{k+1/2}) \, \Big[ q^{Z}_{k+1} - q^{Z}_{k} \Big]
+ q^{R}_{k+1/2} \, b_{\phy,Z} (q_{k+1/2}) \, \Big[ q^{\phy}_{k+1} - q^{\phy}_{k} \Big] ,
\\
J_{43} &= 2 \, q^{R}_{k+1/2} \, b_{\phy} (q_{k+1/2}) ,
\\
J_{44} &= - h .
\end{align}
\end{subequations}

\section{Derivatives}\label{sec:guiding_centre_derivatives}

In this section, the reader can find an overview of the derivatives of the magnetic potential and the magnetic field that appear in the various integrators for guiding centre dynamics.
Of course, all derivatives with respect to $\phy$ vanish as we assume axisymmetry.

\subsubsection{Generalised Magnetic Potential in the Poloidal Plane}

\begin{subequations}
\begin{align}
A^{*}_{R,R} &= A_{R,R} + u_{,R} \, b_{R} + u \, b_{R,R} &
A^{*}_{Z,R} &= A_{Z,R} + u_{,R} \, b_{Z} + u \, b_{Z,R}
\\
A^{*}_{R,Z} &= A_{R,Z} + u_{,Z} \, b_{R} + u \, b_{R,Z} &
A^{*}_{Z,Z} &= A_{Z,Z} + u_{,Z} \, b_{Z} + u \, b_{Z,Z}
\end{align}
\end{subequations}
\vspace{-2em}
\begin{subequations}
\begin{align}
A^{*}_{R,RR} &= A_{R,RR} + u_{,RR} \, b_{R} + u_{,R} \, b_{R,R} + u_{,R} \, b_{R,R} + u \, b_{R,RR} \\
A^{*}_{R,RZ} &= A_{R,RZ} + u_{,RZ} \, b_{R} + u_{,R} \, b_{R,Z} + u_{,Z} \, b_{R,R} + u \, b_{R,RZ} \\
A^{*}_{R,ZZ} &= A_{R,ZZ} + u_{,ZZ} \, b_{R} + u_{,Z} \, b_{R,Z} + u_{,Z} \, b_{R,Z} + u \, b_{R,ZZ} \\
A^{*}_{Z,RR} &= A_{Z,RR} + u_{,RR} \, b_{Z} + u_{,R} \, b_{Z,R} + u_{,R} \, b_{Z,R} + u \, b_{Z,RR} \\
A^{*}_{Z,RZ} &= A_{Z,RZ} + u_{,RZ} \, b_{Z} + u_{,R} \, b_{Z,Z} + u_{,Z} \, b_{Z,R} + u \, b_{Z,RZ} \\
A^{*}_{Z,ZZ} &= A_{Z,ZZ} + u_{,ZZ} \, b_{Z} + u_{,Z} \, b_{Z,Z} + u_{,Z} \, b_{Z,Z} + u \, b_{Z,ZZ}
\end{align}
\end{subequations}

\subsubsection{Generalised Magnetic Potential in Tokamak Geometry}

\begin{subequations}
\begin{align}
A^{*}_{R,R} &= A_{R,R} + u \, b_{R,R} &
A^{*}_{Z,R} &= A_{Z,R} + u \, b_{Z,R} &
A^{*}_{\phy,R} &= A_{\phy,R} + u \, b_{\phy,R}
\\
A^{*}_{R,Z} &= A_{R,Z} + u \, b_{R,Z} &
A^{*}_{Z,Z} &= A_{Z,Z} + u \, b_{Z,Z} &
A^{*}_{\phy,Z} &= A_{\phy,Z} + u \, b_{\phy,Z}
\\
A^{*}_{R,u} &= b_{R} &
A^{*}_{Z,u} &= b_{Z} &
A^{*}_{\phy,u} &= b_{\phy}
\end{align}
\end{subequations}
\vspace{-2em}
\begin{subequations}
\begin{align}
A^{*}_{R,RR} &= A_{R,RR} + u \, b_{R,RR} &
A^{*}_{Z,RR} &= A_{Z,RR} + u \, b_{Z,RR} &
A^{*}_{\phy,RR} &= A_{\phy,RR} + u \, b_{\phy,RR}
\\
A^{*}_{R,RZ} &= A_{R,RZ} + u \, b_{R,RZ} &
A^{*}_{Z,RZ} &= A_{Z,RZ} + u \, b_{Z,RZ} &
A^{*}_{\phy,RZ} &= A_{\phy,RZ} + u \, b_{\phy,RZ}
\\
A^{*}_{R,ZZ} &= A_{R,ZZ} + u \, b_{R,ZZ} &
A^{*}_{Z,ZZ} &= A_{Z,ZZ} + u \, b_{Z,ZZ} &
A^{*}_{\phy,ZZ} &= A_{\phy,ZZ} + u \, b_{\phy,ZZ}
\\
A^{*}_{R,uR} &= b_{R,R} &
A^{*}_{Z,uR} &= b_{Z,R} &
A^{*}_{\phy,uR} &= b_{\phy,R}
\\
A^{*}_{R,uZ} &= b_{R,Z} &
A^{*}_{Z,uZ} &= b_{Z,Z} &
A^{*}_{\phy,uZ} &= b_{\phy,Z}
\\
A^{*}_{R,uu} &= 0 &
A^{*}_{Z,uu} &= 0 &
A^{*}_{\phy,uu} &= 0
\end{align}
\end{subequations}

\subsubsection{Parallel Velocity in the Poloidal Plane}

\begin{subequations}
\begin{align}
u_{,R} &= - \bigg[ p_{\phi} + \dfrac{r^{2} B_{0}}{2q} \bigg] \dfrac{B_{,R}}{B_{0} R_{0}}
- \dfrac{B \, (R-R_{0})}{q R_{0}}  \\
u_{,Z} &= - \bigg[ p_{\phi} + \dfrac{r^{2} B_{0}}{2q} \bigg] \dfrac{B_{,Z}}{B_{0} R_{0}}
- \dfrac{B Z}{q R_{0}}
\end{align}
\end{subequations}
\vspace{-1em}
\begin{subequations}
\begin{align}
u_{,RR} &= - \bigg[ p_{\phi} + \dfrac{r^{2} B_{0}}{2q} \bigg] \dfrac{B_{,RR}}{B_{0} R_{0}}
- \dfrac{2 \, (R-R_{0}) \, B_{,R} + B}{q R_{0}}  \\
u_{,RZ} &= - \bigg[ p_{\phi} + \dfrac{r^{2} B_{0}}{2q} \bigg] \dfrac{B_{,RZ}}{B_{0} R_{0}}
- \dfrac{(R-R_{0}) \, B_{,Z} + Z B_{,R}}{q R_{0}}  \\
u_{,ZZ} &= - \bigg[ p_{\phi} + \dfrac{r^{2} B_{0}}{2q} \bigg] \dfrac{B_{,ZZ}}{B_{0} R_{0}}
- \dfrac{2 Z B_{,Z} + B}{q R_{0}}
\end{align}
\end{subequations}

\subsubsection{Magnetic Potential}

\begin{subequations}
\begin{align}
A_{R,R} &= - \dfrac{B_{0} R_{0} Z}{2 R^{2}} &
A_{Z,R} &= - \dfrac{B_{0} R_{0}}{2 R} &
A_{\phy,R} &= - \dfrac{A_{\phy}}{R} - \dfrac{B_{0} \, (R - R_{0})}{qR}
\\
A_{R,Z} &= \dfrac{B_{0} R_{0}}{2 R} &
A_{Z,Z} &= 0 &
A_{\phy,Z} &= - \dfrac{B_{0} Z}{qR}
\end{align}
\end{subequations}
\vspace{-1em}
\begin{subequations}
\begin{align}
A_{R,RR} &= \dfrac{B_{0} R_{0} Z}{R^{3}} &
A_{Z,RR} &= \dfrac{B_{0} R_{0}}{2 R^{2}} &
A_{\phy,RR} &= - 2 \, \dfrac{A_{\phy,R}}{R} - \dfrac{B_{0}}{qR}
\\
A_{R,RZ} &= - \dfrac{B_{0} R_{0}}{2 R^{2}} &
A_{Z,RZ} &= 0 &
A_{\phy,RZ} &= \dfrac{B_{0} Z}{qR^{2}}
\\
A_{R,ZZ} &= 0 &
A_{Z,ZZ} &= 0 &
A_{\phy,ZZ} &= - \dfrac{B_{0}}{qR}
\end{align}
\end{subequations}

\subsubsection{Magnetic Field}

\begin{subequations}
\begin{align}
B_{,R} &= B \, \bigg[ \dfrac{R-R_{0}}{S^{2}} - \dfrac{1}{R} \bigg] &
S_{,R} &= \dfrac{R - R_{0}}{S}
\\
B_{,Z} &= B \, \dfrac{Z}{S^{2}} &
S_{,Z} &= \dfrac{Z}{S}
\end{align}
\end{subequations}
\vspace{-1em}
\begin{subequations}
\begin{align}
B_{,RR} &= B_{,R} \, \bigg[ \dfrac{R-R_{0}}{S^{2}} - \dfrac{1}{R} \bigg] + B \bigg[ \dfrac{1}{R^{2}} + \dfrac{1}{S^{2}} - 2 \, \dfrac{(R - R_{0})^{2}}{S^{4}} \bigg] \\
B_{,RZ} &= \dfrac{1}{S^{2}} \bigg[ Z \, B_{,R} - 2 \, (R-R_{0}) \, B_{,Z} \bigg] \\
B_{,ZZ} &= \dfrac{B}{S^{2}} \, \bigg[ 1 - \dfrac{Z^{2}}{S^{2}} \bigg]
\end{align}
\end{subequations}

\subsubsection{Normalised Magnetic Field}

\begin{subequations}
\begin{align}
b_{R,R} &= - \dfrac{1}{S} \, b_{R} b_{Z} &
b_{Z,R} &= - \dfrac{1}{S} \, \bigg[ b_{Z}^{2} - 1 \bigg] &
b_{\phy,R} &= - \dfrac{1}{S} \, b_{Z} b_{\phy}
\\
b_{R,Z} &= \hphantom{-} \dfrac{1}{S} \, \bigg[ b_{R}^{2} - 1 \bigg] &
b_{Z,Z} &= \hphantom{-} \dfrac{1}{S} \, b_{R} b_{Z} &
b_{\phy,Z} &= \hphantom{-} \dfrac{1}{S} \, b_{R} b_{\phy}
\end{align}
\end{subequations}
\vspace{-1em}
\begin{subequations}
\begin{align}
b_{R,RR} &= \hphantom{-} \dfrac{b_{R}}{S^{2}} \, \bigg[ 3 \, b_{Z}^{2} - 1 \bigg] &
b_{Z,RR} &= \hphantom{-} \dfrac{b_{Z}}{S^{2}} \, \bigg[ 3 \, b_{Z}^{2} - 3 \bigg] &
b_{\phy,RR} &= \hphantom{-} \dfrac{b_{\phy}}{S^{2}} \, \bigg[ 3 \, b_{Z}^{2} - 1 \bigg]
\\
b_{R,RZ} &= - \dfrac{b_{Z}}{S^{2}} \, \bigg[ 3 \, b_{R}^{2} - 1 \bigg] &
b_{Z,RZ} &= - \dfrac{b_{R}}{S^{2}} \, \bigg[ 3 \, b_{Z}^{2} - 1 \bigg] &
b_{\phy,RZ} &= - \dfrac{3}{S^{2}} \, b_{R} b_{Z} b_{\phy}
\\
b_{R,ZZ} &= \hphantom{-} \dfrac{b_{R}}{S^{2}} \, \bigg[ 3 \, b_{R}^{2} - 3 \bigg] &
b_{Z,ZZ} &= \hphantom{-} \dfrac{b_{Z}}{S^{2}} \, \bigg[ 3 \, b_{R}^{2} - 1 \bigg] &
b_{\phy,ZZ} &= \hphantom{-} \dfrac{b_{\phy}}{S^{2}} \, \bigg[ 3 \, b_{R}^{2} - 1 \bigg]
\end{align}
\end{subequations}

\end{subappendices}

\chapter{Kinetic Theory}\label{ch:kinetic_theory}

In the kinetic theory of plasma dynamics \cite{LandauLiftshitz10, Swanson:2011}, the system of charged particles constituting the plasma is described by a distribution function $f (t, x, v)$ that can be seen as a phasespace density. The integral of $f$ over some phase space region $\Omega$
\begin{align}\label{eq:vlasov_kinetic_1}
\int \limits_{\Omega} f (t, x, v) \, dx \, dv .
\end{align}

gives the number of particles in that region, such that its velocity integral yields the particle density
\begin{align}\label{eq:vlasov_kinetic_density}
n(x) = \int \limits_{-\infty}^{+\infty} f (t, x, v) \, dv ,
\end{align}

at a given point in space, and the integral over full phase space gives the total number of particles $N$ in the system
\begin{align}\label{eq:vlasov_kinetic_particle_number}
N = \int \limits_{-\infty}^{+\infty} \int \limits_{-\infty}^{+\infty} f (t, x, v) \, dx \, dv .
\end{align}

Conservation of phasespace volume along the particle trajectories implies that the distribution function evolves according to the Liouville equation
\begin{align}\label{eq:vlasov_kinetic_liouville}
\dfrac{\partial f}{\partial t} + \dot{x} \cdot \dfrac{\partial f}{\partial x} + \dot{v} \cdot \dfrac{\partial f}{\partial v} = 0
\end{align}

which is a linear advection equation in phasespace.
When the particle motion can be described by a canonical Hamiltonian system, this equation can be expressed as
\begin{align}\label{eq:vlasov_vlasov_1}
\dfrac{\partial f}{\partial t} + [ f, h ]_{xp} = 0
\end{align}

where
\begin{align}\label{eq:vlasov_poisson_brackets}
[ f, h ]_{xp} = \dfrac{\partial f}{\partial x} \dfrac{\partial h}{\partial p} - \dfrac{\partial f}{\partial p} \dfrac{\partial h}{\partial x} ,
\end{align}

are the canonical Poisson brackets, $h$ is the particle Hamiltonian for the system under consideration, and $p$ is the canonical momentum conjugate to $x$.

\section{The Vlasov-Poisson and Vlasov-Maxwell Systems}

Replacing the acceleration with the Lorentz force
\begin{align}
\dot{v}
= \dfrac{1}{m} F_{\text{Lorentz}}
= \dfrac{q}{m} \Big\lgroup E +  \dfrac{1}{c} \, v \times B \Big\rgroup
\end{align}

we obtain the Vlasov equation (also referred to as the collisionless Boltzmann equation)
\begin{align}\label{eq:vlasov_equation}
\dfrac{\partial f}{\partial t} + v \cdot \dfrac{\partial f}{\partial x} + \dfrac{q}{m} \Big\lgroup E + \dfrac{1}{c} \, v \times B \Big\rgroup \cdot \dfrac{\partial f}{\partial v} = 0 .
\end{align}

The fields $E$ and $B$ are computed self-consistently with respect to the particle distribution given by $f$.

\subsection{The Vlasov-Maxwell System}

The Vlasov-Maxwell system consists of the Vlasov equation (\ref{eq:vlasov_equation})
\begin{align}\label{eq:vlasov_maxwell}
\dfrac{\partial f}{\partial t} + v \cdot \dfrac{\partial f}{\partial x} + \dfrac{q}{m} \, \bigg( E + \dfrac{1}{c} \, v \times B \bigg) \cdot \dfrac{\partial f}{\partial v} = 0
\end{align}

and the dynamical Maxwell's equations
\begin{align}
\nabla \times E &= - \dfrac{1}{c} \, \dfrac{\partial B}{\partial t} . & &&
\nabla \times B &= \dfrac{1}{c} \, \bigg( \dfrac{\partial E}{\partial t} + j \bigg) . & &&
\end{align}

The charge density $\rho$ and the current density $j$ are given as moments of the distribution function
\begin{align}
\rho &= q \int f (t, x, v) \, dv , &
j    &= q \int v \, f (t, x, v) \, dv , & &&
\end{align}

and $E$ and $B$ satisfy the constraints
\begin{align}
\nabla \cdot E &= \rho , &
\nabla \cdot  B &= 0 . &
\end{align}

The particle Hamiltonian $h$ has the form
\begin{align}\label{eq:vlasov_vm_hamiltonian}
h = \dfrac{1}{2m} \, \bigg( p - \dfrac{q}{c} \, A (x) \bigg)^{2} + q \phi (x) .
\end{align}

where the canonical momentum $p$ is given by
\begin{align}
p = mv + \dfrac{q}{c} \, A(x) ,
\end{align}

$A$ is the magnetic vector potential and $\phi$ is the electrostatic potential, such that the electromagnetic fields are obtained by
\begin{align}
E &= - \nabla \phi - \dfrac{\partial A}{\partial t} , &
B &= \nabla \times A . & &&
\end{align}

\subsection{The Vlasov-Poisson System}

In the non-relativistic case, when $v \ll c$, the $v \times B$ in the Lorentz force is very small an can be neglected.
If, moreover, $B$ exhibits only little change over time, the electric and magnetic field decouple, and we can describe the system by the electrostatic Vlasov equation
\begin{align}\label{eq:vlasov_vp_vlasov_1}
\dfrac{\partial f}{\partial t} + v \cdot \dfrac{\partial f}{\partial x} + \dfrac{q}{m} \, E \cdot \dfrac{\partial f}{\partial v} = 0
\end{align}

The electric field $E$ can be replaced by the electrostatic potential $E = - \nabla \phi$, such that the Vlasov-equation becomes
\begin{align}\label{eq:vlasov_vp_vlasov_2}
\dfrac{\partial f}{\partial t} + v \cdot \dfrac{\partial f}{\partial x} - \dfrac{q}{m} \, \dfrac{\partial \phi}{\partial x} \cdot \dfrac{\partial f}{\partial v} = 0 .
\end{align}

The electrostatic potential $\phi$ is determined through the Poisson equation
\begin{align}\label{eq:vlasov_vp_poisson_1}
- \Delta \phi = \rho
\end{align}

where $\rho$ is the charge density
\begin{align}
\rho = qn = q \int f \, dv .
\end{align}

With the particle Hamiltonian $h$, consisting of the kinetic energy of the particles and their potential energy in the electrostatic field,
\begin{align}\label{eq:vlasov_vp_hamiltonian}
h =\dfrac{1}{2m} \, p^{2} + q \phi ,
\end{align}

the Vlasov equation (\ref{eq:vlasov_vp_vlasov_2}) can be expressed with Poisson brackets as in (\ref{eq:vlasov_vlasov_1}).
In the following treatment, it will however be more practical to express the Vlasov equation with respect to velocity phasespace variables $(x,v)$ instead of $(x,p)$ as in the canonical Poisson brackets (\ref{eq:vlasov_poisson_brackets}). In that form, the Vlasov equation reads
\begin{align}
\dfrac{\partial f}{\partial t} + [f, h]_{xv} = 0
\end{align}

with particle Hamiltonian
\begin{align}
h = \dfrac{m}{2} \, v^{2} + q \phi
\end{align}

and noncanonical Poisson brackets
\begin{align}
[ \cdot , \cdot ]_{x,p} = \dfrac{1}{m} [ \cdot , \cdot ]_{x,v} .
\end{align}

The additional mass factor will disappear in the normalisation procedure.

\subsubsection{Dimensionless Equations}

For the numerical treatment we normalise the Vlasov-Poisson equation to a dimensionless equation. 
Charges are normalised to $e$, where the electron charge is $q_{e} = - e$. Masses are normalised to the electron mass $m_{e}$.
The speed of light $c$ is set to one, and the average densities are also normalised to one,
\begin{align}
\dfrac{1}{L_{x}} \int f (x,v) \, dx \, dv = 1,
\end{align}

with $L_{x}$ the size of the system.
If both, electrons and ions (assumed to have opposite charge $\pm e$), are treated dynamically, the mass ratio must be accounted for in the Hamiltonian of the ions
\begin{subequations}
\begin{align}
\dfrac{\partial f_{e}}{\partial t} + [f_{e} , h_{e}] &= 0 , &
h_{e} &= \dfrac{1}{2} \, v^{2} - \phi , \\
\dfrac{\partial f_{i}}{\partial t} + [f_{i} , h_{i}] &= 0 , &
h_{i} &= \dfrac{1}{2} \dfrac{m_{i}}{m_{e}} \, v^{2} + \phi . &
\end{align}
\end{subequations}

In that case, the Poisson equation, of course, has contributions from both species
\begin{align}
- \Delta \phi &= \int ( f_{i} -  f_{e} ) \, dv .
\end{align}

If the scale of the dynamics we are interested in is very fast the ions can be considered as fixed in space and only the electrons have to be treated dynamically.
The ions therefore constitute a neutralising background that has to be taken into account in the Poisson equation in terms of the ions' average charge density (which is normalised to $1$). Thereby we obtain the following set of equations
\begin{align}\label{eq:vlasov_dimensionsless_system}
\dfrac{\partial f}{\partial t} + [f , h] &= 0 , &
h &= \dfrac{1}{2} \, v^{2} - \phi , &
- \Delta \phi &= 1 - \int f \, dv . &
\end{align}

We will restrict the following treatment to this case.

\subsection{Conservation Properties}

The Vlasov-Poisson (and Vlasov-Maxwell) system conserves a number of quantities that should in principle also be preserved in any numerical simulation.
If this is not possible, their error can give a hint of the validity of a numerical simulation and should therefore always be monitored.
These conserved quantities are\footnote{
Here, we merely list the conserved quantities, for proofs see e.g. the lecture notes of \citeauthor{Sonnendruecker:2013} \cite{Sonnendruecker:2013}.
}
\begin{itemize}
\item positivity and maximum principle (follows from characteristics theory)
\begin{align}\label{eq:vlasov_conservation_1}
0 \leq f (x,v,t) \; \leq \; \max \limits_{(x,v)} f_{0} (x,v) ,
\end{align}

\item total linear momentum $P$
\begin{align}\label{eq:vlasov_conservation_3}
\dfrac{d}{dt} \bigg( \int v \, f (x,v,t) \, dx \, dv \bigg) &= 0 ,
\end{align}

\item total energy $E$
\begin{align}\label{eq:vlasov_conservation_4}
\dfrac{d}{dt} \bigg( \int v^{2} \, f (x,v,t) \, dx \, dv + \dfrac{1}{2} \int \left( \dfrac{\partial \phi}{\partial x} \right)^{2} dx \bigg) &= 0 ,
\end{align}

\item any functional of the form
\begin{align}
\dfrac{d}{dt} \bigg( \int F(f) \, dx \, dv \bigg) &= 0 ,
\end{align}

with the especially important cases of integral norms $L_{p}$ and entropy

\item integral norms $L_{p}$
\begin{align}\label{eq:vlasov_conservation_5}
\dfrac{d}{dt} \bigg( \int \big( f (x,v,t) \big)^{p} \, dx \, dv \bigg) &= 0 &
& \text{for} &
& 1 \leq p \leq \infty , & &&
\end{align}

where $p=1$ corresponds to the total particle number $N$ and $p = \infty$ is the maximum norm

\item entropy $S$
\begin{align}
\dfrac{d}{dt} \bigg( \int f \, \ln f \, dx \, dv \bigg) &= 0 .
\end{align}

\end{itemize}

In the derivation of the variational integrators, we concentrate on preserving the total particle number, the total linear momentum, and the total energy.

\section{Review of Action Principles}

As pointed out in the chapter on classical mechanics and field theory, action principles are very powerful tools for the description of physical theories.
Not only do they allow us to derive equations of motions in a general and covariant way, but they also provide a machinery for finding conserved quantities for free, namely, the Noether theorem.

The very first action principle for the Vlasov-Maxwell system is due to \citeauthor{Low:1958} and was published in 1958 \cite{Low:1958}:
\begin{align}\label{eq:vlasov_action_low}
L = \int \bigg[ \dfrac{m}{2} \left( \dfrac{\partial x}{\partial t} \right)^2 - q \phi (x, t) + q \, \dfrac{\partial x}{\partial t} \cdot A (x, t) \bigg] \, f (x, v) \, dx \, dv + \dfrac{1}{2} \int \Big( E^2 - B^2 \Big) \, dx .
\end{align}
Since then, a plethora of different action principles for the Vlasov-Poisson and Vlasov-Maxwell systems have been proposed. These action principles can be classified by the variables they use for the particles part; there are Lagrangian descriptions \cite{Low:1958, Sugama:2000}, Eulerian descriptions \cite{YeMorrison:1992, Larsson:1992, Larsson:1993, Fla:1994, Brizard:2000a, Brizard:2000b}, and mixed Lagrangian-Eulerian (i.e., Hamilton-Jacobi) descriptions \cite{Pfirsch:1984, PfirschMorrison:1985, MorrisonPfirsch:1989, PfirschMorrison:1991}. For the electromagnetic fields, of course, in all action principles Eulerian variables are used.

While Lagrangian action principles are suitable for the derivation of numerical schemes for particle-in-cell (PIC) codes, the natural basis for a Vlasov code is a purely Eulerian action principle. These, however. have a severe difficulty, namely that the Vlasov system is inherently noncanonical - the distribution function $f$ does not have a canonical conjugate field variable. As we will see below, there are several possibilities to circumvent this shortcoming. Unfortunately, none of these possibilities, while being quite elegant analytically, do lend themselves to a straight forward discretisation.

What all action principles have in common is the electrostatic Lagrangian $\mcal{L}_{\phi}$ or the electromagnetic Lagrangian $\mcal{L}_{\text{EM}}$, depending on the Lagrangian describing the Vlasov-Poisson or the Vlasov-Maxwell system. These Lagrangian densities are given by
\begin{align}
\mcal{L}_{\phi} &= \dfrac{1}{2} \bigg( \dfrac{\partial \phi}{\partial x} \bigg)^{2}, &
\mcal{L}_{\text{EM}} &= \dfrac{1}{2} \Big( E^{2} + B^{2} \Big) . & && 
\end{align}

\subsection{Parametrisation of the Distribution Function}

An approach well known in fluid dynamics for a long time is the use of Clebsch variables, first applied to the Vlasov system by \citeauthor{YeMorrison:1992} \cite{YeMorrison:1992}. Here the field variable, in our case the distribution function $f$, is parametrised as the Poisson bracket of two Clebsch potentials $\alpha$ and $\beta$ that constitute canonically conjugate field variables:
\begin{align}
f = [ \alpha , \beta ] .
\end{align}

The corresponding Lagrangian density reads
\begin{align}
\mcal{L}_{f} = \alpha \dfrac{\partial \beta}{\partial t} - [ \alpha , \beta ] \, h ,
\end{align}

with the Hamiltonian density
\begin{align}
\mcal{H}_{f} = f \, h = [ \alpha , \beta ] \, h .
\end{align}

and $h$ the particle Hamiltonian.
The variations with respect to $\beta$ and $\alpha$ lead to Vlasov equations for $\alpha$ and $\beta$, respectively
\begin{align}
\dfrac{\partial \alpha}{\partial t} + [ \alpha, h ] &= 0 , &
- \dfrac{\partial \beta}{\partial t} - [ \beta,  h ] &= 0 . & &&
\end{align}

It is easy to show, that if $\alpha$ and $\beta$ obey the Vlasov equation, so does $f$. Just insert the parametrisation of $f$ into the Vlasov equation,
\begin{align}
\dfrac{\partial [ \alpha, \beta ]}{\partial t} + \Big[ [ \alpha, \beta ] , h \Big] = 0,
\end{align}

apply the Leibniz rule to the time derivative,
\begin{align}
\dfrac{\partial [ \alpha, \beta ]}{\partial t} =
\left[ \dfrac{\partial \alpha}{\partial t}, \beta \right] + \left[ \alpha , \dfrac{\partial \beta}{\partial t} \right],
\end{align}

and the Jacobi identity to the Poisson bracket,
\begin{align}
\Big[ [ \alpha, \beta ] , h \Big] + \Big[ [ \beta , h ] , \alpha \Big]  + \Big[ [ h , \alpha ] , \beta \Big] = 0,
\end{align}

to get
\begin{align}
\bigg[ \dfrac{\partial \beta}{\partial t} + [ \alpha, h ] , \beta \bigg] - \bigg[ \dfrac{\partial \beta}{\partial t} + [ \beta,  h ] , \alpha \bigg] = 0 .
\end{align}
The outer Poisson brackets on the left hand side vanish, if $\alpha$ and $\beta$ obey the Vlasov equation.

The problems with this parametrisation are the behaviour of the Clebsch potentials at the boundaries, especially for spatially periodic boundary conditions of $f$, their continuity over the domain, as well as the question of how to initialise $\alpha$ and $\beta$ for a given distribution function $f$. \\

A related approach was taken by \citeauthor{Fla:1994} \cite{Fla:1994, FlaKraus:2011}. He parametrises $f$ with respect to some reference distribution function $f_0$ from the same symplectic leaf, i.e., a distribution function with the same number of particles and the same energy.
$f_{0}$ does not refer to the initial conditions, but is in general evolving along with $f$.
The Poisson bracket of $f_0$ with some generator $S$ then gives the deformation of $f_0$ towards the actual distribution function $f$:
\begin{align}
f = f_0 + [ S , f_0 ] .
\end{align}

Here, $f_0$ might be an equilibrium solution (e.g. a local Maxwellian distribution) while $f$ describes the turbulent state of the system. In this description, however, there is no limitation on the difference $\Delta f = f - f_0$, besides that both, $f$ and $f_0$, have to lie on the same symplectic leaf. So this is not comparable to a $\delta f$ method.

The corresponding action is somewhat more complicated, as $f$ and $f_0$ generate different electromagnetic fields, which has to be taken into account. 
The Lagrangian for the distribution function reads
\begin{align}
\mcal{L}_{f} = f_0 \, \dfrac{\partial S}{\partial t} - [ S , f_0 ] \, h - e \, f_0 \, (\phi - \phi_0) ,
\end{align}

and the Lagrangian for the electrostatic potentials is
\begin{align}
\mcal{L}_{\phi} = \dfrac{1}{2} \Bigg\lgroup \bigg( \dfrac{\partial \phi}{\partial x} \bigg)^2 - \bigg( \dfrac{\partial \phi_0}{\partial x} \bigg)^2 \Bigg\rgroup .
\end{align}

So that the action reads
\begin{align}
\mcal{A} [f_{0}, f, \phi_{0}, \phi]
\nonumber
&= \int \Bigg\lgroup f_0 \, \dfrac{\partial S}{\partial t} - [ S , f_0 ] \, h - e \, f_0 \, (\phi - \phi_0) \Bigg\rgroup \, dt \, dx \, dp \\
& \hspace{8em}
+ \dfrac{1}{2} \int \Bigg\lgroup \bigg( \dfrac{\partial \phi}{\partial x} \bigg)^2 - \bigg( \dfrac{\partial \phi_0}{\partial x} \bigg)^2 \Bigg\rgroup \, dt \, dx .
\end{align}

The variation with respect to $S$ yields the Vlasov equation for $f_0$,
\begin{align}
\dfrac{\partial f_0}{\partial t} + [ f_0, h ] = 0 .
\end{align}

The variation with respect to $f_0$ yields a Vlasov-like equation for $S$,
\begin{align}
\dfrac{\partial S}{\partial t} + [ S, h ] = e \, ( \phi - \phi_0 ) .
\end{align}

The variations with respect to $\phi$ and $\phi_0$ yield the corresponding Poisson equations
\begin{align}
\Delta \phi\hphantom{_0} &= - e \int \Big( f_0 + [ S , f_0 ] \Big) \, dp = - e \int f \, dp \\
\Delta \phi_0 &= - e \int f_0 \, dp .
\end{align}

\citeauthor{Fla:1994}'s parametrisation is subject to similar problems as the Clebsch parametrisation. It is not obvious how to find the generating function $S$ for given distribution functions $f$ and $f_0$. Besides, it is much more suggestive to prescribe $f_0$ and $S$ and obtain the initial $f$ through the parametrisation as this procedure gives a physical meaning to $S$ that would be lost when proceeding the other way around. Unfortunately, this strategy can not be followed, unless a set of generating functions $S$ for some standard scenarios in plasma physical simulation is found.

\subsection{Constrained Variations}

\citeauthor{Brizard:2000a} \cite{Brizard:2000a, Brizard:2000b} suggested an action for the Vlasov-Maxwell system that is defined on an eight-dimensional extended phasespace, adding time and energy to position and momentum.
He uses constraint variations $\delta f = [ \delta S, f ]_{\mrm{cov}}$, where $\delta S$ is the infinitesimal generator of the variation and $[ \cdot , \cdot ]_{\mrm{cov}}$ denote Poisson brackets in the extended phasespace
\begin{align}\label{eq:vlasov_action_brizard_1}
[ f , h ]_{\mrm{cov}} = [f,h] + \dfrac{\partial f}{\partial \eps} \dfrac{\partial h}{\partial t} - \dfrac{\partial f}{\partial t} \dfrac{\partial h}{\partial \eps} .
\end{align}

The action is written as
\begin{align}\label{eq:vlasov_action_brizard_2}
\mcal{A}_{f} &= \int f \, H_{\mrm{cov}} \, dt \, d\eps \, dx \, dp ,
\end{align}

with the covariant Hamiltonian $H_{\mrm{cov}} (t, \eps; x, p) = h (t, x, p) - \eps$, thus
\begin{align}\label{eq:vlasov_action_brizard_3}
\mcal{A}_{f} &= \int f \, ( h - \eps) \, dt \, d\eps \, dx \, dp .
\end{align}

The variation of $\mcal{A}$ is computed as follows
\begin{align}\label{eq:vlasov_action_brizard_4}
\delta \mcal{A}_{f}
&= \int \delta f \, ( h - \eps) \, dt \, d\eps \, dx \, dp \\
&= \int [ \delta S, f ]_{\mrm{cov}} \, ( h - \eps) \, dt \, d\eps \, dx \, dp \\
&= \int \delta S \, [ f, ( h - \eps) ]_{\mrm{cov}} \, dt \, d\eps \, dx \, dp .
\end{align}

Introducing a space-time-split in the Poisson bracket, we get
\begin{align}\label{eq:vlasov_action_brizard_5}
\delta \mcal{A}_{f}
&= \int \delta S \, \Bigg\lgroup \dfrac{\partial f}{\partial \eps} \dfrac{\partial (h - \eps)}{\partial t} - \dfrac{\partial f}{\partial t} \dfrac{\partial (h - \eps)}{\partial \eps} + [ f, h ] - [ f, \eps ] \Bigg\rgroup \, dt \, d\eps \, dx \, dp .
\end{align}

As $\partial h / \partial \eps = \partial h / \partial t = 0$ and of course $\partial \eps / \partial t = 0$ as $(t, \eps; x, p)$ are independent variables, we get
\begin{align}\label{eq:vlasov_action_brizard_6}
\delta \mcal{A}_{f}
&= \int \delta S \, \Bigg\lgroup \dfrac{\partial f}{\partial t} + [ f, h ] \Bigg\rgroup \, dt \, dx \, dp ,
\end{align}

which holds for any variation $\delta S$ and thus yields the Vlasov equation, i.e.,
\begin{align}\label{eq:vlasov_action_brizard_7}
\dfrac{\delta \mcal{A}_{f}}{\delta S} = 0
\hspace{1em}
\Rightarrow
\hspace{1em}
\dfrac{\partial f}{\partial t} + [ f, h ] = 0 .
\end{align}

This action principle does not fit the formalism of chapter \ref{ch:variational} very well.
While it certainly is possible to use the variational integrator framework to compute a discrete variational derivative of the action (\ref{eq:vlasov_action_brizard_2}), constrained to the form $\delta f = [ \delta S, f ]$, it is not easy to see what properties the resulting discrete equations will have. The application of the discrete Noether theorem from section \ref{sec:vi_infinite_noether_theorem} does not seem straight forward. And it is not obvious how to incorporate a symmetrisation of the Poisson brackets in a natural way. We will see in section \ref{sec:vlasov_variational} that this is crucial to retain some of the symmetries of the continuous system on the discrete level and obtain a robust numerical scheme (see also appendix \ref{ch:brackets}).
Besides, the use of extended phasespace adds further complications, as the additional dimensions have to be removed by restriction of the dynamics to a hyperplane of constant energy after the application of the discrete action principle.

\subsection{Euler-Poincaré Reduction}\label{sec:kinetic_theory_euler_poincare}

Another action principle based on constrained variations is the one by \citeauthor{Cendra:1998} \cite{Cendra:1998}, which is obtained from Low's action principle (\ref{eq:vlasov_action_low}) by Euler-Poincaré reduction \cite{Holm:2011, Holm:2009, MarsdenRatiu:2002, Holm:1998}.
Albeit this formulation is not directly applicable to the variational integrator approach as well, it deserves some attention as it constitutes the most natural geometric description of the Vlasov-Maxwell system that has been found so far. The formulation becomes even more interesting as it has recently been extended to gyrokinetics by \citeauthor{Squire:2013} \cite{Squire:2013}.
Furthermore, work by \citeauthor{Pavlov:2011} \cite{Pavlov:2011} suggests that a descendent of the variational integrator method as it is presented in this work can be applied to this formulation (for more details see section \ref{sec:outlook_euler_poincare}).

The basic idea is to reduce the system by using the invariance of the Lagrangian under particle relabelling
\begin{align}\label{eq:vlasov_action_euler_poincare_1}
\psi (x_{0}, v_{0}) = \big( x(x_{0}, v_{0}), v(x_{0}, v_{0}) \big) .
\end{align}

$\psi$ is the particle evolution map. It maps particles with initial phasespace position $(x_{0}, v_{0})$ to their current phasespace position $(x, v)$.
The distribution function is therefore given as
\begin{align}\label{eq:vlasov_action_euler_poincare_2}
f = f_{0} \circ \psi^{-1} ,
\end{align}

or explicitly
\begin{align}\label{eq:vlasov_action_euler_poincare_3}
f \big( x( x_{0}, v_{0}, t) , v (x_{0}, v_{0}, t) \big) = f_{0} (x_{0}, v_{0}) ,
\end{align}

i.e., $f$ is just carried along the particle flow.
The Lagrangian can be written
\begin{align}\label{eq:vlasov_action_euler_poincare_4}
L_{f_{0}} (\psi, \dot{\psi}, \phi, \dot{\phi}, A, \dot{A} )
\nonumber
&= \int f_{0} (x_{0}, v_{0}) \, \bigg[ \bigg( \dfrac{e}{c} \, A(x) + mv \bigg) \cdot \dot{x} - \dfrac{m}{2} \, v^{2} - e \phi (x) \bigg] \, dx_{0} \, dv_{0} \\
& \hspace{8em}
+ \dfrac{1}{8\pi} \int \bigg[ \bigg( - \nabla \phi - \dfrac{\partial A}{\partial t} \bigg)^{2} - \bigg( \nabla \times A \bigg)^{2} \bigg] \, dx .
\end{align}

Invariance of the Lagrangian under the particle relabelling transformation $\psi$ means
\begin{align}\label{eq:vlasov_action_euler_poincare_5}
L_{f_{0}} (\psi, \dot{\psi}, \phi, \dot{\phi}, A, \dot{A} )
=
L_{f_{0} \psi^{-1}} (\psi\psi^{-1}, \dot{\psi}\psi^{-1}, \phi, \dot{\phi}, A, \dot{A} )
\equiv l (u, \dot{\psi}, \phi, \dot{\phi}, A, \dot{A} ) ,
\end{align}

where $u = ( \dot{x}, \dot{v} )$ is the phasespace velocity field
\begin{align}\label{eq:vlasov_action_euler_poincare_6}
u (x,v) \equiv \dot{\psi} \circ \psi^{-1} (x,v) ,
\end{align}

such that
\begin{multline}\label{eq:vlasov_action_euler_poincare_7}
l (u, \dot{\psi}, \phi, \dot{\phi}, A, \dot{A} )
= \int f (x, v) \, \bigg[ \bigg( \dfrac{e}{c} \, A(x) + mv \bigg) \cdot u_{x} - \dfrac{m}{2} \, v^{2} - e \phi (x) \bigg] \, dx \, dv \\
+ \dfrac{1}{8\pi} \int \bigg[ \bigg( - \nabla \phi - \dfrac{\partial A}{\partial t} \bigg)^{2} - \bigg( \nabla \times A \bigg)^{2} \bigg] \, dx ,
\end{multline}

where $x$ and $v$ are now considered as coordinates rather than fields and $u_{x}$ is the spatial component of the phasespace velocity $u$.

We now have to compute the variations of $l$ with respect to $u$ and $f$. Variations with respect to $\phi$ and $A$ yield Maxwell's equation in the usual way.
Variations of the particle evolution map $\psi$ lead to variations in the phasespace velocity,
\begin{align}\label{eq:vlasov_action_euler_poincare_8}
\delta u = \dfrac{\partial \eta}{\partial t} + [u, \eta] ,
\end{align}

which have the form of Lin constraints, well know in fluid dynamics \cite{SeligerWhitham:1968}. Here, $[ \cdot , \cdot ]$ denotes not the Poisson but the Lie bracket, i.e.,
\begin{align}\label{eq:vlasov_action_euler_poincare_9}
[u, \eta] = (u \cdot \nabla_{z}) \eta - (\eta \cdot \nabla_{z}) u ,
\end{align}

where $\nabla_{z}$ denotes the nabla operator in phasespace.
Variations of $\psi$ also induce variations of the distribution function,
\begin{align}\label{eq:vlasov_action_euler_poincare_10}
\delta f = - \nabla_{z} \cdot (f\eta) .
\end{align}

According to (\ref{eq:vlasov_action_euler_poincare_2}), the evolution of $f$ is is determined by the phasespace advection equation,
\begin{align}\label{eq:vlasov_action_euler_poincare_11}
\dfrac{\partial f}{\partial t} + (u \cdot \nabla_{z} ) f = 0 ,
\end{align}

which obviously is the Vlasov equation in conservation form.
Computing the variational derivative of $\int l \, dt$ with respect to $\delta u$ and $\delta f$ leads to the Euler-Poincaré equations,
\begin{align}\label{eq:vlasov_action_euler_poincare_12}
\dfrac{\partial}{\partial t} \dfrac{\delta l}{\delta u} + (u \cdot \nabla_{z} ) \dfrac{\delta l}{\delta u} = f \, \nabla_{z} \, \dfrac{\delta l}{\delta f} .
\end{align}

With the reduced Lagrangian (\ref{eq:vlasov_action_euler_poincare_7}) this leads to
\begin{align}
u_{x} &= v , &
u_{v} &= E + \dfrac{1}{c} \, v \times B , & &&
\end{align}

such that (\ref{eq:vlasov_action_euler_poincare_11}) takes the expected form (\ref{eq:vlasov_maxwell}).

\subsection{Lie Action Principles}

A number of action principles have been derived based on Lie group methods \cite{YeMorrison:1992, Larsson:1992, Larsson:1993, Fla:1994}.
These are not only more complicated than the actions introduced above, they also suffer from similar problems, i.e., they employ constrained variations, auxiliary variables, generating functions, etc. We do not want to go into detail here, but just mention that all of these do not seem to be applicable for our purposes.

\section{Variational Discretisation}\label{sec:vlasov_variational}

We have seen in the previous section that, even though a variety of action principles for the Vlasov-Poisson and Vlasov-Maxwell systems exist, none of them appears directly applicable within the variational integrator framework.
We therefore have to build an extended Lagrangian as described in section \ref{ch:classical_extended_lagrangians}.

\subsection{Extended Lagrangian}

To write the action for the dimensionless Vlasov-Poisson system (\ref{eq:vlasov_dimensionsless_system}), we need two Ibragimov multipliers, $g (t, x, v)$ for the Vlasov equation and $\psi (t, x)$ for the Poisson equation,
\begin{align}\label{eq:vlasov_vi_action}
\mathcal{A} [f, g, \phi, \psi]
&= \int g \, \bigg( \dfrac{\partial f}{\partial t} + [f, h] \bigg) \, dt \, dx \, dv
+ \int \psi \, \bigg( \Delta \phi + 1 - \int f \, dv \bigg) \, dt \, dx
\end{align}

Computing the variations results in the following equations of motion
\begin{subequations}
\begin{align}
\dfrac{\delta \mathcal{A}}{\delta g} &= + \dfrac{\partial f}{\partial t} + [f, h] = 0, &
\dfrac{\delta \mathcal{A}}{\delta \psi} &= \Delta \phi + \int f \, dv - 1 = 0 \\
\dfrac{\delta \mathcal{A}}{\delta f} &= - \dfrac{\partial g}{\partial t} - [g, h] - \psi = 0, &
\dfrac{\delta \mathcal{A}}{\delta \phi} &= \Delta \psi + \int [g,f] \, dv = 0
\end{align}
\end{subequations}

A compatible solution of the auxiliary variables is given by $g = f$ and $\psi = 0$, such that the adjoint equation for the Vlasov equation becomes the Vlasov equation itself, and the Poisson bracket in the adjoint equation of the Poisson equation is identical zero as $[f,f] = 0$, thereby admitting a constant solution for $\psi$, where we choose $\psi = 0$ to obtain the first equality.
The solution vector of the extended system of equations is thus $(f, f, \phi, 0)$.

\subsection{Variational Integrator}

\begin{figure}[H]
\centering
\begin{minipage}[c]{.06\textwidth}
$\dfrac{\partial \varphi}{\partial t}$:
\end{minipage}
\begin{minipage}[c]{.22\textwidth}
\includegraphics[width=\textwidth]{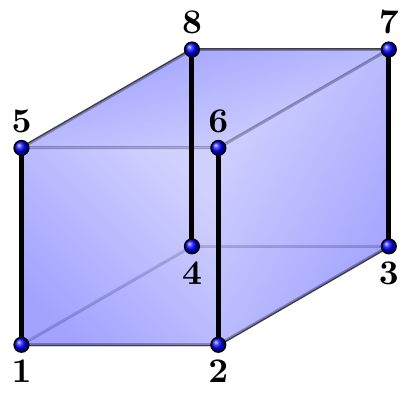}
\end{minipage}
\begin{minipage}[c]{.06\textwidth}
$\dfrac{\partial \varphi}{\partial x}$:
\end{minipage}
\begin{minipage}[c]{.22\textwidth}
\includegraphics[width=\textwidth]{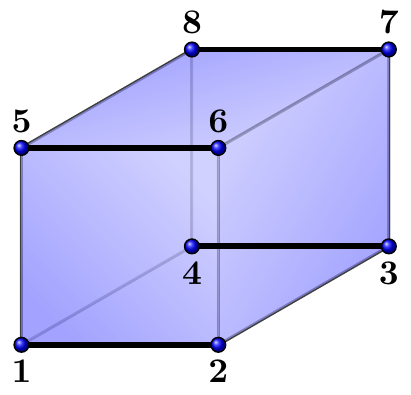}
\end{minipage}
\begin{minipage}[c]{.06\textwidth}
$\dfrac{\partial \varphi}{\partial v}$:
\end{minipage}
\begin{minipage}[c]{.22\textwidth}
\includegraphics[width=\textwidth]{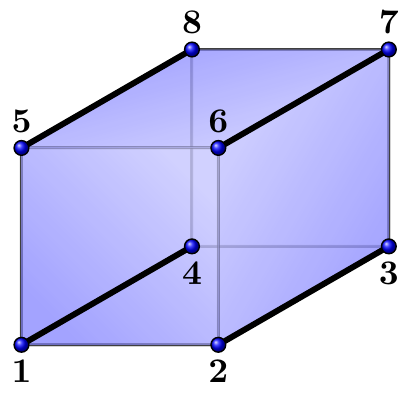}
\end{minipage}
\begin{minipage}[c]{.10\textwidth}
\includegraphics[width=\textwidth]{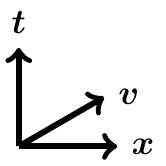}
\end{minipage}

\caption{For a given phasespacetime grid cell, there are four possible ways of defining derivatives in the different coordinate directions $(t,x,v)$, namely along the black lines. The actual discrete derivatives correspond to the averages over all four corresponding possibilities.}
\label{fig:vlasov_vi_derivatives}
\end{figure}

The discretisation of the action (\ref{eq:vlasov_vi_action}) follows exactly along the lines of section \ref{sec:vi_infinite}, the only difference being that we have three dimensions now, time $t$, space $x$, and velocity $v$ (see figure \ref{fig:vlasov_vi_derivatives}).

The approximations of the fields and derivatives are therefore
\begin{subequations}\label{eq:vlasov_vi_derivatives}
\begin{align}
\varphi && & \rightarrow &
& \phy_{\obar{i},\obar{j},\obar{k}} &
& \equiv
\dfrac{1}{8} \Big\lgroup \varphi^1 + \varphi^2 + \varphi^3 + \varphi^4 + \varphi^5 + \varphi^6 + \varphi^7 + \varphi^8 \Big\rgroup
\\
\dfrac{\partial \varphi}{\partial t} && & \rightarrow &
& \bigg( \dfrac{\partial \phy}{\partial t} \bigg)_{\obar{i},\obar{j},k} &
& \equiv
\dfrac{1}{4} \Bigg\lgroup \dfrac{\varphi^5 - \varphi^1}{h_t} + \dfrac{\varphi^6 - \varphi^2}{h_t} + \dfrac{\varphi^7 - \varphi^3}{h_t} + \dfrac{\varphi^8 - \varphi^4}{h_t} \Bigg\rgroup
\\
\dfrac{\partial \varphi}{\partial x} && & \rightarrow &
& \bigg( \dfrac{\partial \phy}{\partial x} \bigg)_{i,\obar{j},\obar{k}} &
& \equiv
\dfrac{1}{4} \Bigg\lgroup \dfrac{\varphi^2 - \varphi^1}{h_x} + \dfrac{\varphi^3 - \varphi^4}{h_x} + \dfrac{\varphi^6 - \varphi^5}{h_x} + \dfrac{\varphi^7 - \varphi^8}{h_x} \Bigg\rgroup
\\
\dfrac{\partial \varphi}{\partial v} && & \rightarrow &
& \bigg( \dfrac{\partial \phy}{\partial v} \bigg)_{\obar{i},j,\obar{k}} &
& \equiv
\dfrac{1}{4} \Bigg\lgroup \dfrac{\varphi^4 - \varphi^1}{h_v} + \dfrac{\varphi^3 - \varphi^2}{h_v} +  \dfrac{\varphi^8 - \varphi^5}{h_v} + \dfrac{\varphi^7 - \varphi^6}{h_v} \Bigg\rgroup . &
& \hspace{5em} &
\end{align}
\end{subequations}

The bar over the indices indicates in which dimension averaging is applied as depicted in figure \ref{fig:vlasov_vi_derivatives}. Considering a given phasespacetime grid cell, the time derivative, for example, can be defined along each of the four highlighted edges of that cell. Our discrete time derivative is the average over all four of that possibilities, denoted by overbars $\obar{i},\obar{j}$, but not $\obar{k}$ as that is the coordinate direction of the derivative.
The fields themselves are averaged in all three dimensions, such that their indices have overbars $\obar{i},\obar{j}$,$\obar{k}$.

In the Poisson equation, we do not have a $v$ dimension, so we have to define a reduced field average and a reduced $x$ derivative as follows
\begin{align}\label{eq:vlasov_vi_derivatives_2}
\phy_{\obar{i},\obar{k}} & \equiv \dfrac{1}{4} \Big\lgroup \varphi^1 + \varphi^2 + \varphi^5 + \varphi^6 \Big\rgroup, &
\bigg( \dfrac{\partial \phy}{\partial x} \bigg)_{i,\obar{k}} & \equiv \dfrac{1}{4} \Bigg\lgroup \dfrac{\varphi^2 - \varphi^1}{h_x} + \dfrac{\varphi^6 - \varphi^5}{h_x} \Bigg\rgroup . &
\end{align}

Some care has to be taken when discretising the Poisson bracket (see \citeauthor{SalmonTalley:1989} \cite{SalmonTalley:1989}).
To retain the properties of the continuous bracket (antisymmetry and the Jacobi identity) at the discrete level, a symmetrisation has to be introduced in the Lagrangian.
One has to realise that by partial integration the even permutations in the integrand are all identical (assuming boundary conditions such that the boundary terms of the partial integration vanish)
\begin{align*}
\int g \, [ f, h ] \, dx \, dv
= \int f \, [ h, g ] \, dx \, dv
= \int h \, [ g, f ] \, dx \, dv .
\end{align*}

Hence, instead of one of the permutations a convex combination can be used just as well
\begin{align}
\int g \, [ f, h ] \, dx \, dv &= \int \Big\lgroup \alpha \, g \, [ f, h ] + \beta \, f \, [ h, g ] + \gamma \, h \, [ g, f ] \Big\rgroup \, dx \, dv &
& \text{with} &
\alpha + \beta + \gamma &= 1 .
\end{align}

The symmetric case, i.e., the one that retains the properties of the bracket at the discrete level, corresponds to $\alpha = \beta = \gamma = 1/3$.
We therefore write the action
\begin{align}
\mathcal{A} [f, g, \phi, \psi]
\nonumber
&= \int \bigg[ g \, \dfrac{\partial f}{\partial t} + \dfrac{1}{3} \Big( g \, [f, h] + f \, [h, g] + h \, [g, f] \Big) \bigg] \, dt \, dx \, dv \\
& \hspace{12em}
- \int \bigg[ \dfrac{\partial \psi}{\partial x} \dfrac{\partial \phi}{\partial x} + \psi \, \bigg( \int f \, dv - 1 \bigg) \bigg] \, dt \, dx ,
\end{align}

where we also did a partial integration in the second integral to avoid second order derivatives.
As the two integrals in the action have different integration domains, $(t,x,v)$ for the Vlasov equations and $(t,x)$ for the Poisson equation, we split the discrete Lagrangian into two parts
\begin{align}
\mcal{L}_{d}^{V} &= g_{\obar{i},\obar{j},\obar{k}} \, \bigg( \dfrac{\partial f}{\partial t} \bigg)_{\obar{i},\obar{j},k} + \dfrac{1}{3} \Big( g_{\obar{i},\obar{j},\obar{k}} \, [f, h]_{i,j,k} + f_{\obar{i},\obar{j},\obar{k}} \, [h, g]_{i,j,k} + h_{\obar{i},\obar{j},\obar{k}} \, [g, f]_{i,j,k} \Big) \\
\mcal{L}_{d}^{P} &= - \bigg( \dfrac{\partial \psi}{\partial x} \bigg)_{i,\obar{k}} \bigg( \dfrac{\partial \phi}{\partial x} \bigg)_{i,\obar{k}} - \psi_{\obar{i},\obar{k}} \, \bigg( \sum \limits_{j} f_{\obar{i},\obar{j},\obar{k}} - 1 \bigg)
\end{align}

with the discrete Poisson bracket
\begin{align}
[f, h]_{i,j,k}
=
\bigg( \dfrac{\partial f}{\partial x} \bigg)_{i,\obar{j},\obar{k}}
\bigg( \dfrac{\partial h}{\partial v} \bigg)_{\obar{i},j,\obar{k}}
-
\bigg( \dfrac{\partial f}{\partial v} \bigg)_{\obar{i},j,\obar{k}}
\bigg( \dfrac{\partial h}{\partial x} \bigg)_{i,\obar{j},\obar{k}}
.
\end{align}

With these definitions the discrete action becomes
\begin{align}
\mathcal{A}_{d}
&= h_{t} \, h_{x} \, h_{v} \, \sum \limits_{i,j,k} \mcal{L}_{d}^{V} + h_{t} \, h_{x} \, \sum \limits_{i,k} \mcal{L}_{d}^{P} \\
\nonumber
&= h_{t} \, h_{x} \, h_{v} \, \sum \limits_{i,j,k} \bigg[ g_{\obar{i},\obar{j},\obar{k}} \, \bigg( \dfrac{\partial f}{\partial t} \bigg)_{\obar{i},\obar{j},k} + \dfrac{1}{3} \Big( g_{\obar{i},\obar{j},\obar{k}} \, [f, h]_{i,j,\obar{k}} + f_{\obar{i},\obar{j},\obar{k}} \, [h, g]_{i,j,\obar{k}} + h_{\obar{i},\obar{j},\obar{k}} \, [g, f]_{i,j,\obar{k}} \Big) \bigg] \\
& \hspace{8em}
- h_{t} \, h_{x} \, \sum \limits_{i,k} \bigg[ \bigg( \dfrac{\partial \psi}{\partial x} \bigg)_{i,\obar{k}} \bigg( \dfrac{\partial \phi}{\partial x} \bigg)_{i,\obar{k}} + \psi_{\obar{i},\obar{k}} \, \bigg( \sum \limits_{j} f_{\obar{i},\obar{j},\obar{k}} - 1 \bigg) \bigg] .
\end{align}

The discrete Euler-Lagrange field equations (\ref{eq:vi_infinite_delfeqs}) are computed as
\begin{align}
0
\nonumber
&= \dfrac{\partial \mcal{L}^{V}_d}{\partial g^1} \Big( y_{i,  j  ,k  }, y_{i+1,j  ,k  }, y_{i+1,j+1,k  }, y_{i,  j+1,k  }, y_{i,  j  ,k+1}, y_{i+1,j  ,k+1}, y_{i+1,j+1,k+1}, y_{i,  j+1,k+1} \Big) \\
\nonumber
&+ \dfrac{\partial \mcal{L}^{V}_d}{\partial g^2} \Big( y_{i-1,j  ,k  }, y_{i,  j  ,k  }, y_{i,  j+1,k  }, y_{i-1,j+1,k  }, y_{i-1,j  ,k+1}, y_{i,  j  ,k+1}, y_{i,  j+1,k+1}, y_{i-1,j+1,k+1} \Big) \\
\nonumber
&+ \dfrac{\partial \mcal{L}^{V}_d}{\partial g^3} \Big( y_{i-1,j-1,k  }, y_{i,  j-1,k  }, y_{i,  j  ,k  }, y_{i-1,j  ,k  }, y_{i-1,j-1,k+1}, y_{i,  j-1,k+1}, y_{i,  j  ,k+1}, y_{i-1,j  ,k+1} \Big) \\
\nonumber
&+ \dfrac{\partial \mcal{L}^{V}_d}{\partial g^4} \Big( y_{i,  j-1,k  }, y_{i+1,j-1,k  }, y_{i+1,j  ,k  }, y_{i,  j  ,k  }, y_{i,  j-1,k+1}, y_{i+1,j-1,k+1}, y_{i+1,j  ,k+1}, y_{i,  j  ,k+1} \Big) \\
\nonumber
&+ \dfrac{\partial \mcal{L}^{V}_d}{\partial g^5} \Big( y_{i,  j  ,k-1}, y_{i+1,j  ,k-1}, y_{i+1,j+1,k-1}, y_{i,  j+1,k  }, y_{i,  j  ,k  }, y_{i+1,j  ,k  }, y_{i+1,j+1,k  }, y_{i,  j+1,k  } \Big) \\
\nonumber
&+ \dfrac{\partial \mcal{L}^{V}_d}{\partial g^6} \Big( y_{i-1,j  ,k-1}, y_{i,  j  ,k-1}, y_{i,  j+1,k-1}, y_{i-1,j+1,k  }, y_{i-1,j  ,k  }, y_{i,  j  ,k  }, y_{i,  j+1,k  }, y_{i-1,j+1,k  } \Big) \\
\nonumber
&+ \dfrac{\partial \mcal{L}^{V}_d}{\partial g^7} \Big( y_{i-1,j-1,k-1}, y_{i,  j-1,k-1}, y_{i,  j  ,k-1}, y_{i-1,j  ,k-1}, y_{i-1,j-1,k  }, y_{i,  j-1,k  }, y_{i,  j  ,k  }, y_{i-1,j  ,k  } \Big) \\
&+ \dfrac{\partial \mcal{L}^{V}_d}{\partial g^8} \Big( y_{i,  j-1,k-1}, y_{i+1,j-1,k-1}, y_{i+1,j  ,k-1}, y_{i,  j  ,k-1}, y_{i,  j-1,k  }, y_{i+1,j-1,k  }, y_{i+1,j  ,k  }, y_{i,  j  ,k  } \Big)
\end{align}

for the discrete Vlasov equation and
\begin{align}
0
\nonumber
&= \dfrac{\partial \mcal{L}^{P}_d}{\partial \psi^1} \Big( y_{i,  k  }, y_{i+1,k  }, y_{i+1,k+1}, y_{i,  k+1} \Big)
+ \dfrac{\partial \mcal{L}^{P}_d}{\partial \psi^2} \Big( y_{i-1,k  }, y_{i,  k  }, y_{i,  k+1}, y_{i-1,k+1} \Big) \\
&+ \dfrac{\partial \mcal{L}^{P}_d}{\partial \psi^3} \Big( y_{i-1,k-1}, y_{i,  k-1}, y_{i,  k  }, y_{i-1,k  } \Big)
+ \dfrac{\partial \mcal{L}^{P}_d}{\partial \psi^4} \Big( y_{i,  k-1}, y_{i+1,k-1}, y_{i+1,k  }, y_{i,  k  } \Big)
\end{align}

for the discrete Poisson equation.

These discrete variations yield the following discrete Vlasov-Poisson system
\begin{align}
\label{eq:vlasov_vi_discrete_vlasov_equation}
0
\nonumber
&= \dfrac{\obar{f}_{k+1} - \obar{f}_{k-1}}{2 h_{t}} + \dfrac{1}{8} \Big( A (f_{k+1}, h_{k+1}) + A (f_{k+1}, h_{k}) + A (f_{k}, h_{k+1}) \\
& \hspace{8em}
+ 2 A(f_{k}, h_{k}) + A(f_{k}, h_{k-1}) + A (f_{k-1}, h_{k}) + A (f_{k-1}, h_{k-1}) \Big)
\\
\label{eq:vlasov_vi_discrete_poisson_equation}
0
\nonumber
&= \Delta_{fd} \phi_{ik+1} - \dfrac{1}{4} \bigg( \sum \limits_{j=0}^{n_{v}} f_{i-1jk+1} + 2 \sum \limits_{j=0}^{n_{v}} f_{ijk+1} + \sum \limits_{j=0}^{n_{v}} f_{i+1jk+1} \bigg) + 1 \\
\nonumber
& \hspace{6em}
+ 2 \Delta_{fd} \phi_{ik} - \dfrac{1}{2} \bigg( \sum \limits_{j=0}^{n_{v}} f_{i-1jk  } + 2 \sum \limits_{j=0}^{n_{v}} f_{ijk  } + \sum \limits_{j=0}^{n_{v}} f_{i+1jk  } \bigg) + 2 \\
& \hspace{8em}
+ \Delta_{fd} \phi_{ik-1} - \dfrac{1}{4} \bigg( \sum \limits_{j=0}^{n_{v}} f_{i-1jk-1} + 2 \sum \limits_{j=0}^{n_{v}} f_{ijk-1} + \sum \limits_{j=0}^{n_{v}} f_{i+1jk-1} \bigg) + 1 .
\end{align}

The time derivative of the distribution function $f$ is an average of centred-finite-differences over 9 grid points in phasespace, weighted as depicted below.

\begin{figure}[H]
\centering
\includegraphics[width=.2\textwidth]{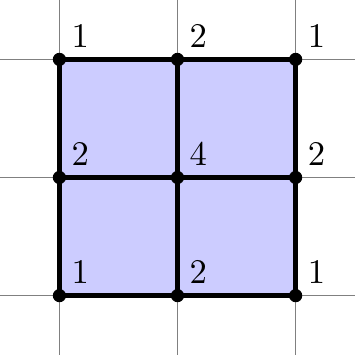}
\end{figure}

The Poisson bracket $A(\cdot,\cdot)$ is discretised by the well known Arakawa scheme \cite{Arakawa:1966} (see also appendix \ref{ch:brackets}). It is noteworthy that the Arakawa discretisation arises naturally from the variational principle. This is of course a consequence of the applied symmetrisation, but that symmetrisation is imperative to retain the symmetries of the continuous Lagrangian on the discrete level.

The discrete Laplace operator $\Delta_{fd}$ is just the standard finite difference stencil $[+1 \; -2 \; +1] / h^{2}$. The Poisson equation is, however, averaged over three points in time, and the charge density is averaged over three points in space.

\subsubsection{Simplifications}\label{sec:kinetic_theory_vi_simplifications}

Overall, we obtain a well working, nonlinearly implicit scheme. It allows, however, for some simplifications.
The first issue is the representation of the time derivative by a second order discretisation. This requires the prescription of initial conditions at two consecutive points in time, which is unnatural as the Vlasov equation requires only one initial condition.
Rewriting the discrete Vlasov equation (\ref{eq:vlasov_vi_discrete_vlasov_equation}) as
\begin{align}\label{eq:vlasov_vi_discrete_vlasov_equation_rewritten}
0
\nonumber
&= \dfrac{\obar{f}_{k+1} - \obar{f}_{k}}{2 h_{t}} + \dfrac{1}{8} \Big( A (f_{k+1}, h_{k+1}) + A (f_{k+1}, h_{k}) + A (f_{k}, h_{k+1}) + A (f_{k}, h_{k}) \Big) \\
&+ \dfrac{\obar{f}_{k} - \obar{f}_{k-1}}{2 h_{t}} + \dfrac{1}{8} \Big( A (f_{k}, h_{k}) + A (f_{k}, h_{k-1}) + A (f_{k-1}, h_{k}) + A (f_{k-1}, h_{k-1}) \Big)
\end{align}

the solution to this issue becomes apparent. Consider the reduced equation
\begin{align}\label{eq:vlasov_vi_discrete_vlasov_equation_reduced}
0
&= \dfrac{\obar{f}_{k+1} - \obar{f}_{k}}{h_{t}} + \dfrac{1}{4} \Big( A (f_{k+1}, h_{k+1}) + A (f_{k+1}, h_{k}) + A (f_{k}, h_{k+1}) + A (f_{k}, h_{k}) \Big) .
\end{align}

If we use this equation to determine $f_{1}$ for given initial conditions $f_{0}$, and use this as initial conditions for (\ref{eq:vlasov_vi_discrete_vlasov_equation_rewritten}), the solution of (\ref{eq:vlasov_vi_discrete_vlasov_equation_rewritten}) will always also be a solution of (\ref{eq:vlasov_vi_discrete_vlasov_equation_reduced}) and vice versa.
We can therefore solve the reduced system (\ref{eq:vlasov_vi_discrete_vlasov_equation_reduced}) instead of (\ref{eq:vlasov_vi_discrete_vlasov_equation_rewritten}), retaining all conservation properties, but replacing the time derivative of the distribution function with a discrete derivative of first order.

By a similar argument we can remove the time average in the discrete Poisson equation (\ref{eq:vlasov_vi_discrete_poisson_equation}).
If we prescribe an initial distribution function at $k=0$ and only use the third line of (\ref{eq:vlasov_vi_discrete_poisson_equation}) to compute the corresponding potential, all three parts of the discrete Poisson equation have to hold separately at all times.
We can thus replace the discrete Poisson equation by
\begin{align}\label{eq:vlasov_vi_discrete_poisson_equation_reduced}
0
&= \Delta_{fd} \phi_{ik+1} - \dfrac{1}{4} \bigg( \sum \limits_{j=0}^{n_{v}} f_{i-1jk+1} + 2 \sum \limits_{j=0}^{n_{v}} f_{ijk+1} + \sum \limits_{j=0}^{n_{v}} f_{i+1jk+1} \bigg) + 1 .
\end{align}

Interestingly, this reduction can already be implemented on the level of the discrete Lagrangian by modifying some of the averaging, i.e., by removing the time average from $g$ in the Vlasov equation and all time averages in the Poisson equation
\begin{align}
\nonumber
\mcal{L}_{d}^{V} &= g_{\obar{i},\obar{j},k} \, \bigg( \dfrac{\partial f}{\partial t} \bigg)_{\obar{i},\obar{j},k}
+
\dfrac{1}{3} \, g_{\obar{i},\obar{j},k} \, \bigg[
\bigg( \dfrac{\partial f}{\partial x} \bigg)_{i,\obar{j},\obar{k}}
\bigg( \dfrac{\partial h}{\partial v} \bigg)_{\obar{i},j,\obar{k}}
-
\bigg( \dfrac{\partial f}{\partial v} \bigg)_{\obar{i},j,\obar{k}}
\bigg( \dfrac{\partial h}{\partial x} \bigg)_{i,\obar{j},\obar{k}}
\bigg] \\
\nonumber
& \hspace{7em}
+
\dfrac{1}{3} \, f_{\obar{i},\obar{j},\obar{k}} \, \bigg[
\bigg( \dfrac{\partial h}{\partial x} \bigg)_{i,\obar{j},\obar{k}}
\bigg( \dfrac{\partial g}{\partial v} \bigg)_{\obar{i},j,k}
-
\bigg( \dfrac{\partial h}{\partial v} \bigg)_{\obar{i},j,\obar{k}}
\bigg( \dfrac{\partial g}{\partial x} \bigg)_{i,\obar{j},k}
\bigg] \\
& \hspace{7em}
+
\dfrac{1}{3} \, h_{\obar{i},\obar{j},\obar{k}} \, \bigg[
\bigg( \dfrac{\partial g}{\partial x} \bigg)_{i,\obar{j},k}
\bigg( \dfrac{\partial f}{\partial v} \bigg)_{\obar{i},j,\obar{k}}
-
\bigg( \dfrac{\partial g}{\partial v} \bigg)_{\obar{i},j,k}
\bigg( \dfrac{\partial f}{\partial x} \bigg)_{i,\obar{j},\obar{k}}
\bigg]
\\
\mcal{L}_{d}^{P} &= - \bigg( \dfrac{\partial \psi}{\partial x} \bigg)_{i,k} \bigg( \dfrac{\partial \phi}{\partial x} \bigg)_{i,k} - \psi_{\obar{i},k} \, \bigg( \sum \limits_{j} f_{\obar{i},\obar{j},k} - 1 \bigg)
.
\end{align}

This is an important point as it allows us to study and compare the discrete symmetries and conservation laws of both the original and the reduced scheme.
Without having done the analysis (which should be addressed in future work), we observe numerically that both schemes preserve the total particle number, the total linear momentum, and the total energy of the system exactly, i.e., up to machine accuracy.
Furthermore, integral norms of the distribution function $f$ are preserved, and as the scheme is symplectic by construction, the phasespace volume is preserved as well.
Solely, positivity of the distribution function and the maximum principle are not preserved automatically.

\subsection{Linearised Lagrangian}\label{sec:vlasov_linearised}

Another simplification that can be introduced on the level of the Lagrangian is a linearisation with respect to time.
The above scheme has extraordinary conservation properties, but in some situations the solution of a nonlinear system of equations might be too demanding in terms of computational time and might thus not be an option.
Moreover, to get quick convergence of the nonlinear iteration (e.g. in a Newton solver) a good predictor or initial guess is necessary.
We therefore derive a linearised scheme that can be used in both of these cases and examine which of the properties of the fully nonlinear scheme are retained.
Again, we apply a different averaging strategy, defining separate averages, space and velocity derivatives for timepoints $k$ and $k+1$
\begin{subequations}
\begin{align}
\varphi_{\obar{i}, \obar{j}, k  } &\equiv \dfrac{1}{4} \Big\lgroup \varphi^1 + \varphi^2 + \varphi^3 + \varphi^4 \Big\rgroup, &
\varphi_{\obar{i}, \obar{j}, k+1} &\equiv \dfrac{1}{4} \Big\lgroup \varphi^5 + \varphi^6 + \varphi^7 + \varphi^8 \Big\rgroup \\
\left( \dfrac{\partial \varphi}{\partial x} \right)_{i, \obar{j}, k  } &\equiv \dfrac{1}{2} \Bigg\lgroup \dfrac{\varphi^2 - \varphi^1}{h_x} + \dfrac{\varphi^3 - \varphi^4}{h_x} \Bigg\rgroup, &
\left( \dfrac{\partial \varphi}{\partial x} \right)_{i, \obar{j}, k+1} &\equiv \dfrac{1}{2} \Bigg\lgroup \dfrac{\varphi^6 - \varphi^5}{h_x} + \dfrac{\varphi^7 - \varphi^8}{h_x} \Bigg\rgroup \\
\left( \dfrac{\partial \varphi}{\partial v} \right)_{\obar{i}, j, k  } &\equiv \dfrac{1}{2} \Bigg\lgroup \dfrac{\varphi^4 - \varphi^1}{h_v} + \dfrac{\varphi^3 - \varphi^2}{h_v} \Bigg\rgroup, &
\left( \dfrac{\partial \varphi}{\partial v} \right)_{\obar{i}, j, k+1} &\equiv \dfrac{1}{2} \Bigg\lgroup \dfrac{\varphi^8 - \varphi^5}{h_v} + \dfrac{\varphi^7 - \varphi^6}{h_v} \Bigg\rgroup
\end{align}
\end{subequations}

and replace the symmetrised Poisson bracket by
\begin{align}
\nonumber
& \hspace{-5em} \overline{ \big( g \, [f, h] + f \, [h, g] + h \, [g, f] \big) }_{i,j,k} = \\
\nonumber
\hspace{3em}
= \dfrac{1}{12} \bigg(
&  g_{\obar{i},\obar{j},k  } \, \big( [f_{k}, h_{k+1}]_{ij} + [f_{k+1}, h_{k}]_{ij} \big)
+ g_{\obar{i},\obar{j},k+1} \, \big( [f_{k}, h_{k+1}]_{ij} + [f_{k+1}, h_{k}]_{ij} \big) \\
&
\nonumber
+ f_{\obar{i},\obar{j},k  } \, \big( [h_{k+1}, g_{k}]_{ij} + [h_{k+1}, g_{k+1}]_{ij} \big)
+ f_{\obar{i},\obar{j},k+1} \, \big( [h_{k  }, g_{k}]_{ij} + [h_{k  }, g_{k+1}]_{ij} \big) \\
&
+ h_{\obar{i},\obar{j},k  } \, \big( [g_{k}, f_{k+1}]_{ij} + [g_{k+1}, f_{k+1}]_{ij} \big)
+ h_{\obar{i},\obar{j},k+1} \, \big( [g_{k}, f_{k  }]_{ij} + [g_{k+1}, f_{k  }]_{ij} \big)
\bigg)
\end{align}

with
\begin{align}
[f_{k}, h_{k+1}]_{i,j}
=
\bigg( \dfrac{\partial f}{\partial x} \bigg)_{i,\obar{j},k  }
\bigg( \dfrac{\partial h}{\partial v} \bigg)_{\obar{i},j,k+1}
-
\bigg( \dfrac{\partial f}{\partial v} \bigg)_{\obar{i},j,k  }
\bigg( \dfrac{\partial h}{\partial x} \bigg)_{i,\obar{j},k+1}
\end{align}

such that in combinations of $f$ and $h$, both fields are always taken at different times.
In the time derivative and in the Poisson equation we apply the simplifications from the previous section, thereby obtaining the linearised discrete Lagrangians
\begin{align}
\mcal{L}_{d}^{V} &= g_{\obar{i},\obar{j},k} \, \bigg( \dfrac{\partial f}{\partial t} \bigg)_{\obar{i},\obar{j},k} + \overline{ \big( g \, [f, h] + f \, [h, g] + h \, [g, f] \big) }_{i,j,k}
\\
\mcal{L}_{d}^{P} &=
- \bigg( \dfrac{\partial \psi}{\partial x} \bigg)_{ik} \bigg( \dfrac{\partial \phi}{\partial x} \bigg)_{ik}
- \psi_{\obar{i}k} \, \bigg( \sum \limits_{j} f_{\obar{i}\obar{j}k} - 1 \bigg)
.
\end{align}

The resulting scheme is
\begin{align}
\label{eq:vlasov_vi_discrete_vlasov_equation_linearised_reduced}
0
&= \dfrac{\obar{f}_{k+1} - \obar{f}_{k}}{h_{t}} + \dfrac{1}{2} \Big( A (f_{k+1}, h_{k}) + A (f_{k}, h_{k+1}) \Big)
\\
\label{eq:vlasov_vi_discrete_poisson_equation_linearised_reduced}
0
\nonumber
&= \Delta_{fd} \phi_{ik+1} - \dfrac{1}{4} \bigg( \sum \limits_{j=0}^{n_{v}} f_{i-1jk+1} + 2 \sum \limits_{j=0}^{n_{v}} f_{ijk+1} + \sum \limits_{j=0}^{n_{v}} f_{i+1jk+1} \bigg) + 1 .
\end{align}

As we will see in the numerical examples, this scheme still preserves the total particle number, linear momentum, and integral norms of the distribution function, but it does not preserve the energy exactly. Instead the usual energy behaviour of symplectic methods is observed, i.e., the energy error oscillates about zero with a bounded amplitude of the oscillation.

The loss of exact energy conservation is almost certainly explained by destruction of some symmetry (namely the one responsible for energy conservation) in the discrete Lagrangian in the course of the linearisation procedure. Again, a detailed analysis of the discrete symmetries and discrete conservation laws should clarify this point.

\section{Velocity Space Collision Operator}\label{sec:vlasov_collision_operator}

A well known problem with low order finite difference schemes like ours is the development of oscillations when phasespace filaments of the order of the grid size develop. In other discretisation techniques, e.g. finite elements or semi-Lagrangian methods, interpolation procedures are employed which, as a side effect, damp these oscillations.
An alternative is to add a velocity space collision operator.

\subsection{Continuous Collision Operator}

We start by considering the Lenard-Bernstein operator \cite{LenardBernstein:1958} which conserves the total particle number but not momentum and energy
\begin{align}
C_{\text{LB}}[f] = \nu \, \dfrac{\partial}{\partial v} \bigg[ \dfrac{\partial f}{\partial v} + v f \bigg] .
\end{align}

We try to fix this by adding correction terms that restore these conservation properties and obtain the same collision operator as \citeauthor{FilbetSonnendruecker:2003} \cite{FilbetSonnendruecker:2003}.
The general expression of the collision operator with collision frequency $\nu$ is
\begin{align}\label{eq:vlasov_collision_operator_general}
C[f] = \nu \, \dfrac{\partial}{\partial v} \bigg[ \dfrac{\partial f}{\partial v} + A(v, f) \, f \bigg]
\end{align}

where the correction term is of the form
\begin{align}\label{eq:vlasov_collision_correction}
A (v, f) = \sum \limits_{n=1}^{K} A_{n} (f) \, v^{n-1} .
\end{align}

If the operator shall preserve the total particle number, linear momentum, and energy (i.e., the zeroth, first and second moment of the distribution function), the velocity integral of the collision operator, multiplied with $\{1, v, v^{2} \}$ has to vanish.
In general, to preserve the first $K$ moments, $C[f]$ has to fulfil
\begin{align}
\int \limits_{-v_{\text{max}}}^{+v_{\text{max}}} v^{k} \, C[f] \, dv
= \nu \int \limits_{-v_{\text{max}}}^{+v_{\text{max}}} v^{k} \, \dfrac{\partial}{\partial v} \bigg[ \dfrac{\partial f}{\partial v} + A(v, f) \, f \bigg] \, dv &= 0 &
& \text{for} &
k &= 0, ..., K . &
\end{align}

The integration domain should be the whole real line, $(-\infty, +\infty)$, but in the discrete case, it suffices if the velocity domain is large enough to ensure that $f$ and its derivatives vanish or are at least very small at the boundaries.
A partial integration with respect to $v$ (neglecting the collision frequency $\nu$) gives
\begin{align}
\bigg[ v^{k} \, \bigg( \dfrac{\partial f}{\partial v} + A(v, f) \, f \bigg) \bigg] \bigg\vert_{-v_{\text{max}}}^{+v_{\text{max}}} - k \int \limits_{-v_{\text{max}}}^{+v_{\text{max}}} v^{k-1} \, \bigg[ \dfrac{\partial f}{\partial v} + A (v, f) \, f \bigg] \, dv
\end{align}

which means that $f$ and $\partial_{v} f$ have to be (close to) zero at $v \pm v_{\text{max}}$ such that
\begin{align}
\bigg[ v^{k} \, \bigg( \dfrac{\partial f}{\partial v} + A(v, f) \, f \bigg) \bigg] \bigg\vert_{-v_{\text{max}}}^{+v_{\text{max}}} \approx 0 .
\end{align}

Assuming that this is fulfilled, conservation of the moments of $f$ requires the following expression to vanish
\begin{align}
\int \limits_{-v_{\text{max}}}^{+v_{\text{max}}} v^{k-1} \, \bigg[ \dfrac{\partial f}{\partial v} + A (v, f) \, f \bigg] \, dv .
\end{align}

Partial integration of the first term gives
\begin{align}
\big[ v^{k-1} f \big] \big\vert_{-v_{\text{max}}}^{+v_{\text{max}}}
- \int \limits_{-v_{\text{max}}}^{+v_{\text{max}}} \bigg[ (k-1) \, v^{k-2} - A (v, f) \, v^{k-1} \bigg] \, f \, dv
\end{align}

where the surface term vanishes (approximately) for $f \approx 0$ at $v \pm v_{\text{max}}$.
Plugging (\ref{eq:vlasov_collision_correction}) into the integral and writing $M_{i}$ for the $i$th moment of $f$, we get the set of conditions
\begin{align}
(k-1) \, M_{k-2} (f) &= \sum \limits_{n=1}^{K} A_{n} (f) \, M_{k+n-2} (f) &
& \text{for} &
k &= 1, ..., K . &
\end{align}

We want to preserve all moments up to $K=2$, so we compute
\begin{align*}
k=1 : && 0 &= A_{1} M_{0} + A_{2} M_{1} , & \\
k=2 : && M_{0} &= A_{1} M_{1} + A_{2} M_{2} . &
& \hspace{24em} &
\end{align*}

With the definition of the moments
\begin{align}
M_{0} &= n = \int f \, dv , &
M_{1} &= nu = \int vf \, dv , &
M_{2} &= n \veps = \int v^{2} f \, dv &
\end{align}

we get the system
\begin{align*}
0 &= A_{1} n + A_{2} nu , \\
n &= A_{1} nu + A_{2} n \veps . &
\end{align*}

the solution of which determines the correction factors to be
\begin{align}\label{eq:vlasov_collision_correction_factors}
A_{1} &= \dfrac{u}{u^{2} - \veps}, &
A_{2} &= - \dfrac{1}{u^{2} - \veps} . &
\hspace{15em}
\end{align}

where the factors $A_{i}$ depend on the distribution function $f$ through the momenta $n$, $u$ and $\eps$.
The full expression of our operator is thus
\begin{align}\label{eq:vlasov_collision_operator}
C[f] = \nu \, \dfrac{\partial}{\partial v} \bigg[ \dfrac{\partial}{\partial v} + \dfrac{v - u}{\veps - u^{2}} \bigg] \, f .
\end{align}

The denominator represents the temperature $T = \veps - u^{2}$ of the plasma, such that the correction factor $(v-u)/T$ corresponds to the thermal spread of the particles about the average velocity.

\subsubsection{Comparison With Other Operators}

The operator (\ref{eq:vlasov_collision_operator}) can be shown to be related to the one presented by \citeauthor{Oppenheim:1965} \cite{Oppenheim:1965}, \citeauthor{OngYu:1969} \cite{OngYu:1969}, as well as \citeauthor{Clemmow:1969} \cite{Clemmow:1969}, by multiplying the right hand side of (\ref{eq:vlasov_collision_operator}) with the denominator of the second term
\begin{align}
C'[f] &= \nu \, \dfrac{\partial}{\partial v} \bigg[ (\veps - u^{2}) \, \dfrac{\partial}{\partial v} + (v - u) \bigg] \, f .
\end{align}

Upon insertion of the definition of $\veps$ we get
\begin{align}\label{eq:vlasov_collision_operator_comparison}
C'[f]
\nonumber
&= \nu \, \dfrac{\partial}{\partial v} \bigg[ \bigg( \dfrac{1}{n} \int v^{2} \, f \, dv - u^{2} \bigg) \dfrac{\partial}{\partial v} + (v - u) \bigg] \, f \\
&= \nu \, \dfrac{\partial}{\partial v} \bigg[ \bigg( \dfrac{1}{n} \int ( v - u )^{2} \, f \, dv \bigg) \dfrac{\partial}{\partial v} + (v - u) \bigg] \, f
\end{align}

where the last equality holds as
\begin{align*}
\dfrac{1}{n} \int ( v - u )^{2} \, f \, dv
= \dfrac{1}{n} \int ( v^{2} - 2 u v + u^{2} ) \, f \, dv
= \dfrac{1}{n} \bigg( \int  v^{2} \, f \, dv - 2 n u^{2} + n u^{2} \bigg)
= \dfrac{1}{n} \int v^{2} \, f \, dv - u^{2} .
\end{align*}

Equation (\ref{eq:vlasov_collision_operator_comparison}) is the expression presented in the above references \cite{Oppenheim:1965, OngYu:1969, Clemmow:1969}. It has similar properties as our operator, i.e., it preserves the total particle number, the total linear momentum, and the total energy. It also relaxes towards a Maxwellian, but the coefficient of the diffusion term modifies its behaviour, such that the strength of the diffusion scales with the thermal energy $\eps - u^{2}$.

\subsection{Discrete Collision Operator}

To obtain the discrete collision operator we repeat the derivation at the discrete level in the same spirit we derived the discrete action principle, i.e., by mimicking the continuous derivation.
We discretise the derivatives in the collision operator (\ref{eq:vlasov_collision_operator_general}) by
\begin{align}\label{eq:vlasov_discrete_collision_operator_general}
C_{d} [j]
\nonumber
&= \nu \, \bigg[ \dfrac{f(j-1) - 2 f(j) + f(j+1)}{h_{v}^{2}} + A_{1,d} \, \dfrac{f(j+1) - f (j-1)}{2 h_{v}} \\
& \hspace{10em}
+ A_{2,d} \, \dfrac{v(j+1) \, f( j+1) - v(j-1) \, f(j-1)}{2 h_{v}} \bigg] .
\end{align}

where we dropped the spatial and time indices as the collision operator is always computed at a single point in spacetime $(i,k)$.
The coefficients $A_{1,d}$ and $A_{2,d}$ are choosen to enforce the discrete conservation properties
\begin{subequations}
\begin{align}
0 &= \sum \limits_{j} \Big( v (j) + v (j+1) \Big) \Big( C_{d} (j) + C_{d} (j+1) \Big) , \\
0 &= \sum \limits_{j} \Big( v^{2} (j) + v^{2} (j+1) \Big) \Big( C_{d} (j) + C_{d} (j+1) \Big) .
\end{align}
\end{subequations}

Discrete partial integration, i.e., reordering of the sums, then leads to the following expressions for the correction factors $A_{1,d}$ and $A_{2,d}$
\begin{align}
A_{1,d} &= - \dfrac{u_{d}}{\veps_{d} - u_{d}^{2}}, &
A_{2,d} &= \dfrac{1}{\veps_{d} - u_{d}^{2}} &
\hspace{12em}
\end{align}

which, no surprises, are exactly the same as their continuous counter parts (\ref{eq:vlasov_collision_correction_factors}).
The important result is that we automatically obtain the correct energy and momentum preserving discretisation of the moments
\begin{align}
n_{d} &= h_{v} \sum \limits_{j} f(j), &
u_{d} &= \dfrac{h_{v}}{n_{d}} \sum \limits_{j} v(j) f(j), &
\veps_{d} &= \dfrac{h_{v}}{n_{d}} \sum \limits_{j} v^{2}(j) f(j). &
\end{align}

The complete discretised collision operator is (replacing the subscript $d$ with grid coordinates)
\begin{align}\label{eq:vlasov_discrete_collision_operator}
C_{i,j,k}
\nonumber
&= \nu \, \bigg[ \dfrac{f(i, j-1, k) - 2 f(i, j, k) + f(i, j+1, k)}{h_{v}^{2}} \\
& \hspace{4em}
+ \dfrac{ [ v(j+1) - u (i,k) ] \, f(i, j+1, k) - [ v(j-1) - u (i,k) ] \, f(i, j-1, k)}{2 h_{v} \, [ \veps (i,k) - u^{2} (i,k) ]} \bigg] .
\end{align}

We add the discrete operator to the simplified nonlinear Vlasov equation (\ref{eq:vlasov_vi_discrete_vlasov_equation_reduced}) by employing a spacetime averaging approach mimicking the result of the discrete variational principle
\begin{align}
\dfrac{\obar{f}_{k+1} - \obar{f}_{k}}{h_{t}} + \dfrac{1}{4} \Big( A (f_{k+1}, h_{k+1}) + A (f_{k+1}, h_{k}) + A (f_{k}, h_{k+1}) + A(f_{k}, h_{k}) \Big) = C_{\obar{i}, j, \obar{k}} [f]
\end{align}

where
\begin{align}
\nonumber
C_{\obar{i}, j, \obar{k}} &= \dfrac{1}{8} \Big[ C_{i-1,j,k} + 2 \, C_{i,j,k} + C_{i+1,j,k} + C_{i-1,j,k+1} + 2 \, C_{i,j,k+1} + C_{i+1,j,k+1} \Big] .
\end{align}

This averaging does not have an effect on the conservation properties as the discrete collision operator is designed to conserve the total particle number, the linear momentum and energy locally, i.e., for each spacetime grid point $(i,k)$ separately.

\section{Numerical Examples}

In this section we consider several numerical examples that can be considered standard benchmark cases \cite{HeathGamba:2012, WatanabeSugama:2005, BesseSonnendruecker:2003, FilbetSonnendruecker:2003, ArberVann:2002, NakamuraYabe:1999, ChengKnorr:1976}.

If not noted otherwise, the simulation domain is $[0, 2\pi / k] \times [ - v_{\text{max}} , + v_{\text{max}} ]$.
The resolution is always $n_{x} = 201$, $n_{v} = 401$, and except for the simulations with the linear integrator, the timestep is $h_{t} = 0.1$ in units of the inverse plasma frequency. For the linear integrator the timestep is $h_{t} = 0.01$.

Most of the examples are initialised as a perturbation of a Maxwellian distribution, given by
\begin{align}\label{eq:vlasov_numerical_exmaples_maxwellian}
f_{M} (x,v) = \dfrac{1}{\sqrt{2\pi}} \, \exp \left\{ - \tfrac{1}{2} \, v^{2} \right\} .
\end{align}

The temperature is set to one such that the thermal velocity is also one. The density is normalised to one.
The initial potential is determined by the initial distribution function via the Poisson equation.

\subsection{Simulation Code}

The variational integrator for the Vlasov-Poisson system, equations (\ref{eq:vlasov_vi_discrete_vlasov_equation_reduced}) and (\ref{eq:vlasov_vi_discrete_poisson_equation_reduced}), constitutes a nonlinearly implicit system of equations.
The nonlinearity is solved by Newton's method where in each Newton step a direct linear solver based on LU decomposition and a GMRES correction is employed.
The initial guess for the Newton solver is either computed by the linear variational integrator from section \ref{sec:vlasov_linearised}, or, in linear or weakly nonlinear examples, the last timestep is used.
Depending on the problem, the Newton solver usually needs 1-3 iterations to converge with a residual smaller $10^{-11}$. In most cases, the LU decomposition of the Jacobian needs only be carried out once per timestep.

The implementation of efficient solvers is a topic left for of future research, but preliminary results suggest that for a sufficiently good initial guess, the LU decomposition can be replaced by an iterative method (GMRES), where only a few iterations are needed to solve the linear system.

The code is implemented in Python and Cython using PETSc \cite{petsc-web-page, petsc-user-ref} to solve the nonlinear system and take care of the parallel communication and MUMPS \cite{mumps-web-page} for the LU decomposition.

\subsection{Diagnostics}

We have not yet carried out a detailed analysis of the discrete conservation laws of the Vlasov-Poisson system. We therefore assume a discrete representation of the conservation properties according to the discretisation of the Lagrangian, i.e., a midpoint representation.

The total particle number is computed as
\begin{align}\label{eq:vlasov_diagnostics_1}
N_{k} = \dfrac{1}{4} \sum \limits_{i=1}^{n_{x}-1} \sum \limits_{j=1}^{n_{v}-1} \big( f_{i,j,k} + f_{i+1,j,k} + f_{i+1,j+1,k} + f_{i,j+1,k} \big) \, h_{x} h_{v} ,
\end{align}

where $i=n_{x}$ corresponds to $i=1$ as we use periodic boundary conditions in space. Furthermore, the velocity domain should always be chosen large enough, such that $f_{i,1} = f_{i, n_{v}} = 0$, the above expression is really just a sum of $f$ over the whole phasespace grid $(i,j)$,
\begin{align}\label{eq:vlasov_diagnostics_2}
N_{k} = \dfrac{1}{4} \sum \limits_{i=1}^{n_{x}} \sum \limits_{j=1}^{n_{v}} f_{i,j,k} \, h_{x} h_{v} .
\end{align}

Similar to (\ref{eq:vlasov_diagnostics_1}), the $L^{2}$ norm is computed as
\begin{align}\label{eq:vlasov_diagnostics_3}
L^{2}_{k} = \sum \limits_{i=1}^{n_{x}-1} \sum \limits_{j=1}^{n_{v}-1} \bigg[ \dfrac{1}{4} \Big( f_{i,j,k} + f_{i+1,j,k} + f_{i+1,j+1,k} + f_{i,j+1,k} \Big) \Bigg]^{2} \, h_{x} h_{v} ,
\end{align}

momentum is computed as
\begin{align}
P_{k} = \dfrac{1}{8} \sum \limits_{i=1}^{n_{x}-1} \sum \limits_{j=1}^{n_{v}-1} \big( f_{i,j,k} + f_{i+1,j,k} + f_{i+1,j+1,k} + f_{i,j+1,k} \big) \big( v_{j} + v_{j+1} \big) \, h_{x} h_{v} ,
\end{align}

energy is computed as
\begin{align}
E_{k} = \dfrac{1}{16} \sum \limits_{i=1}^{n_{x}-1} \sum \limits_{j=1}^{n_{v}-1}
\nonumber
& \big( f_{i,j,k} + f_{i+1,j,k} + f_{i+1,j+1,k} + f_{i,j+1,k} \big) \times \\
& \hspace{4em}
\times \big( h_{i,j,k} + h_{i+1,j,k} + h_{i+1,j+1,k} + h_{i,j+1,k} \big) \, h_{x} h_{v} ,
\end{align}

and entropy is computed as
\begin{align}
S_{k} = \dfrac{1}{4} \sum \limits_{i=1}^{n_{x}-1} \sum \limits_{j=1}^{n_{v}-1}
\nonumber
& \big( f_{i,j,k} + f_{i+1,j,k} + f_{i+1,j+1,k} + f_{i,j+1,k} \big) \times \\
& \hspace{4em}
\times \log \bigg( \dfrac{1}{4} \big( f_{i,j,k} + f_{i+1,j,k} + f_{i+1,j+1,k} + f_{i,j+1,k} \big) \bigg) \, h_{x} h_{v} .
\end{align}

It is expected that a rigorous calculation of the discrete conservation laws from the discrete Noether theorem (section \ref{sec:vi_infinite_noether_theorem}) may improve on the quality of the discrete conservation laws.
Even with the foregoing heuristic diagnostics, we obtain very satisfying results so that a more precise analysis is left for future work.

\subsection{Landau Damping}

Landau damping is probably the most popular benchmark for the Vlasov equation, first because it is a purely kinetic effect involving phase mixing, and second because there are analytical results available to compare with (at least in the linear case).
The initial distribution function is given by
\begin{align}
f (x,v) = f_{M} \, \big( 1 + A \, \cos (kx) \big) ,
\end{align}

where $f_{M}$ is a Maxwellian distribution (\ref{eq:vlasov_numerical_exmaples_maxwellian}).
With $k=0.5$, the spatial simulation domain is $[0, 4\pi]$, and $v_{\text{max}} = 10$.
The resolution is $n_{x} = 201$, $n_{v} = 401$, and the timestep is $h_{t} = 0.1$ in units of the inverse plasma frequency.
The spatial step width $h_{x} = 2\pi / n_{x} k$ depends on the chosen wave number $k$.

\subsubsection{Linear Landau Damping}

\begin{figure}
\centering
\includegraphics[width=.95\textwidth]{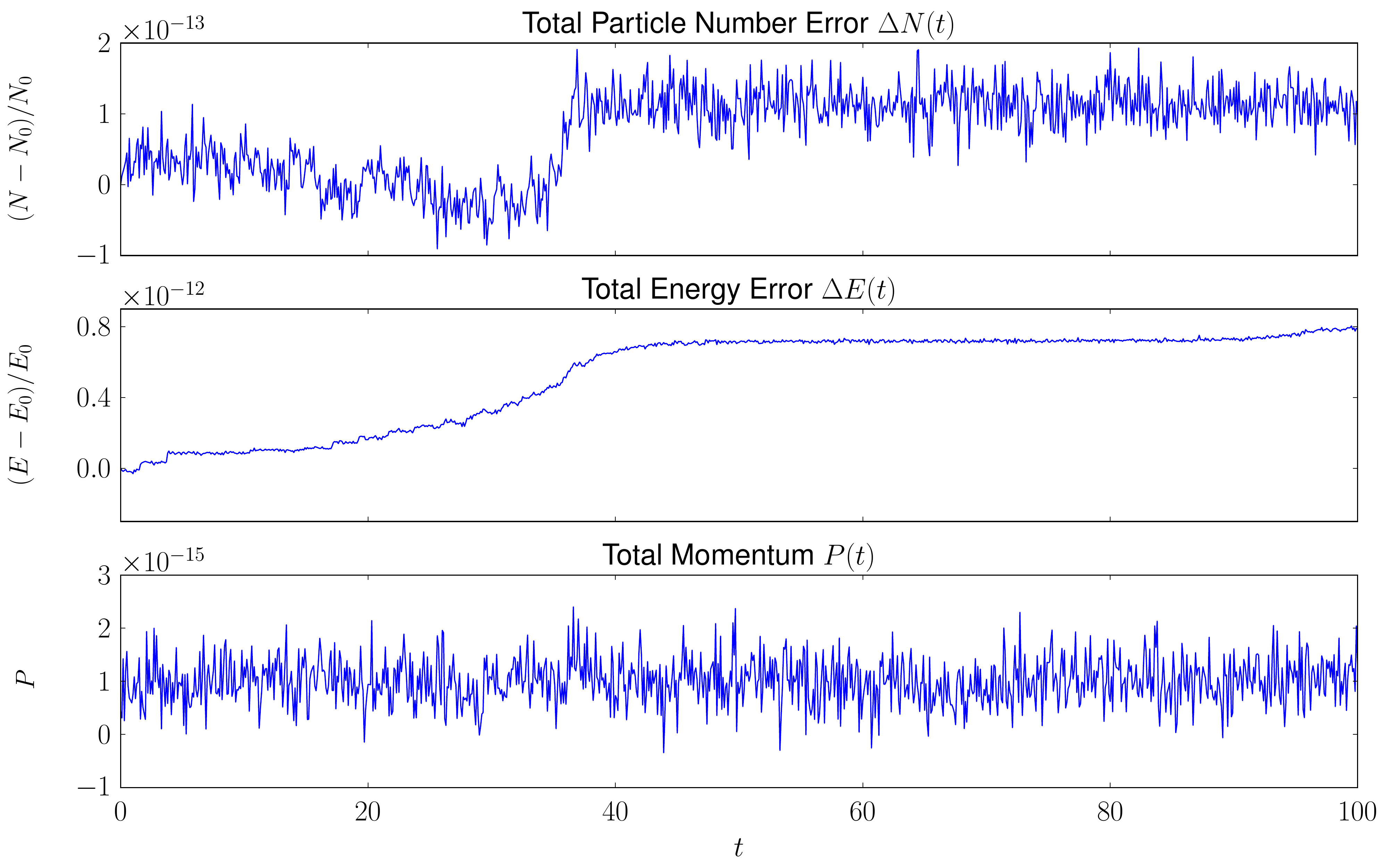}
\caption{Linear Landau damping without collisions. The total particle number and the total linear momentum are well preserved. Conservation of the total energy is slightly violated due to subgrid mode effects.}
\label{fig:vlasov_landau_linear_nu0_NEP}
\end{figure}

At first, we consider the linear case, which can be compared with theoretical results, that is an initial perturbation with $A = 0.01$, $k = 0.5$, and without collisions ($\nu = 0$).
Figure \ref{fig:vlasov_landau_linear_nu0_NEP}) shows the time traces of the errors of the total particle number, the total energy, and linear momentum, while the evolution of the electrostatic energy is displayed in figure \ref{fig:vlasov_landau_linear_nu0_potential}.
Using only the marked maxima, the damping rate is computed to be $\gamma = - 0.152$, which is very close to the theoretical value of $\gamma = - 0.153$.
Using only the first ten maxima, we obtain the predicted value $\gamma = - 0.153$.
The total particle number and the total linear momentum are preserved optimally (see figure \ref{fig:vlasov_landau_linear_nu0_NEP}) and exhibit the expected oscillatory behaviour about a constant value. The error in the total energy is very small but seems to grow during the simulation.
Indeed, the error of the particle number oscillates about zero for $t < 40$, and then jumps to $10^{-13}$.
The energy error appears to grow monotonically, although it remains very small.
This is attributed to the formation of structures on a scale length shorter than the grid step size, namely, subgrid modes which are discussed below.

We can therefore conclude that, without additional ad hoc devices (such as hyperdiffusion), the integrator shows remarkable conservation properties and accuracy, as far as the grid is sufficient to resolve the phase-mixing structure of the distribution function.
For long-time integration, special care of subgrid modes should be taken.

\subsubsection{Subgrid Modes and Collision Operator}

\begin{figure}
\centering
\includegraphics[width=\textwidth]{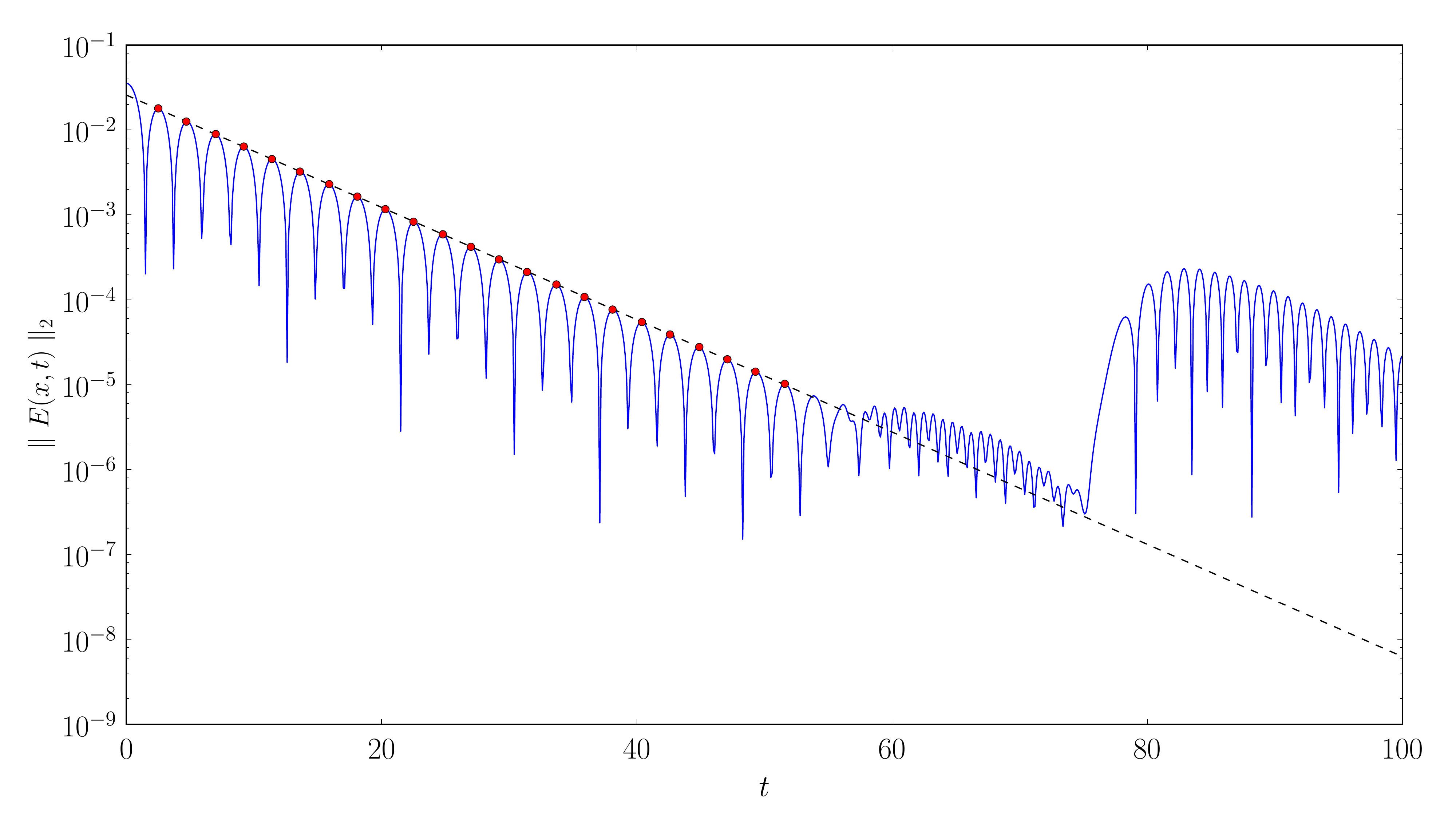}
\caption{Linear Landau damping without collisions. Evolution of the electrostatic energy. At about $t=40$, subgrid modes start to develop, eventually spoiling the damping of the initial perturbation.}
\label{fig:vlasov_landau_linear_nu0_potential}
\end{figure}

At about $t=40$, subgrid modes start to develop. Consequently, large gradients in the distribution function appear, which in turn lead to an unphysically large electrostatic potential. Therefore, the total energy error increases, as can be seen in figure \ref{fig:vlasov_landau_linear_nu0_NEP}, and the damping rate becomes spurious (figure \ref{fig:vlasov_landau_linear_nu0_potential}).
To remove these subgrid modes, we employ the collision operator described in section \ref{sec:vlasov_collision_operator}.
It dissipates the $L^{2}$ norm but retains the conservation of total particle number, total energy, and total linear momentum.

\begin{figure}
\centering
\includegraphics[width=.95\textwidth]{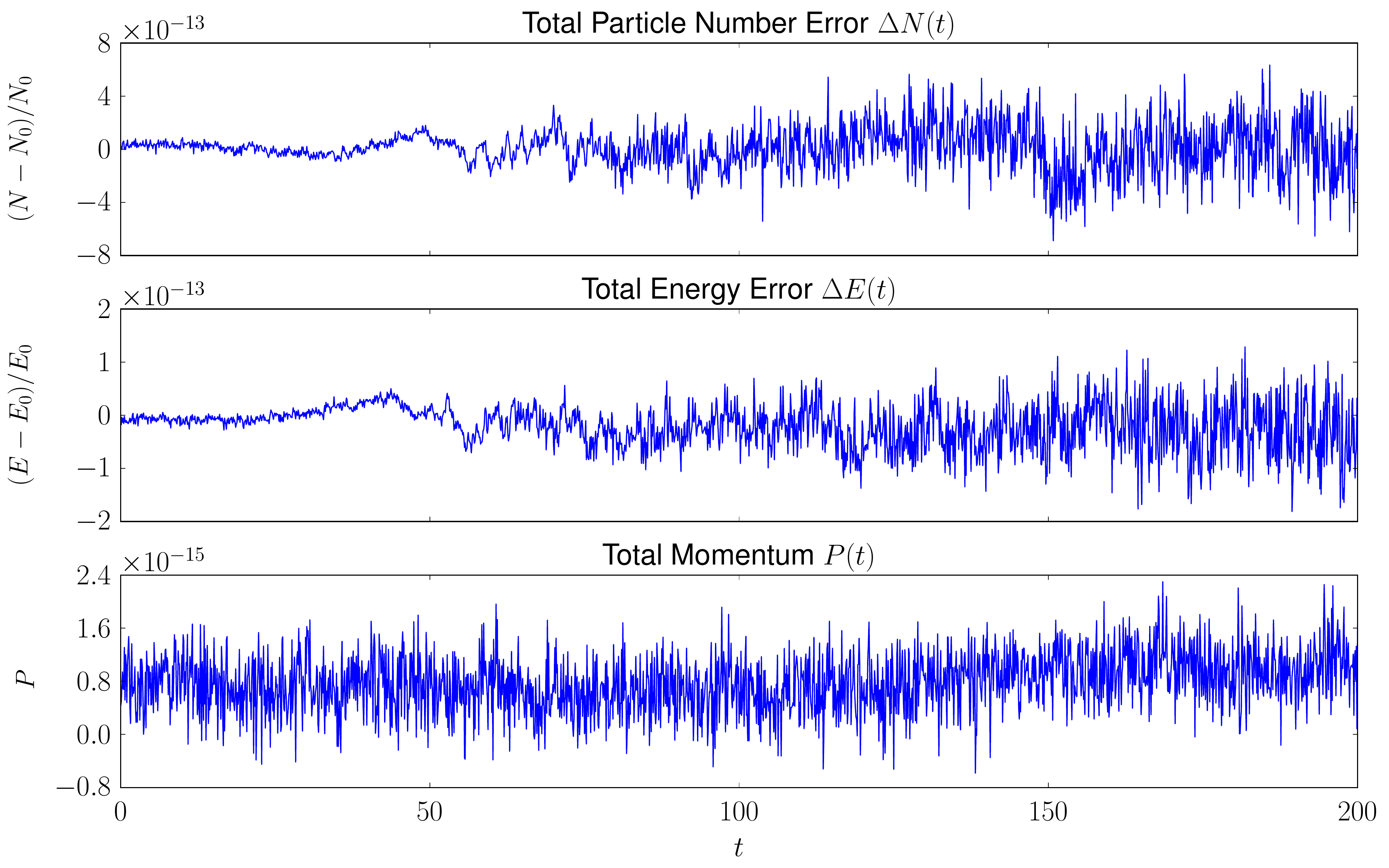}
\caption{Linear Landau damping with collision frequency $\nu = 10^{-4}$. With collisions, the total particle number, energy, and linear momentum are all preserved up to machine precision.}
\label{fig:vlasov_landau_linear_nu1E-4_NEP}
\end{figure}

\begin{figure}
\centering
\includegraphics[width=\textwidth]{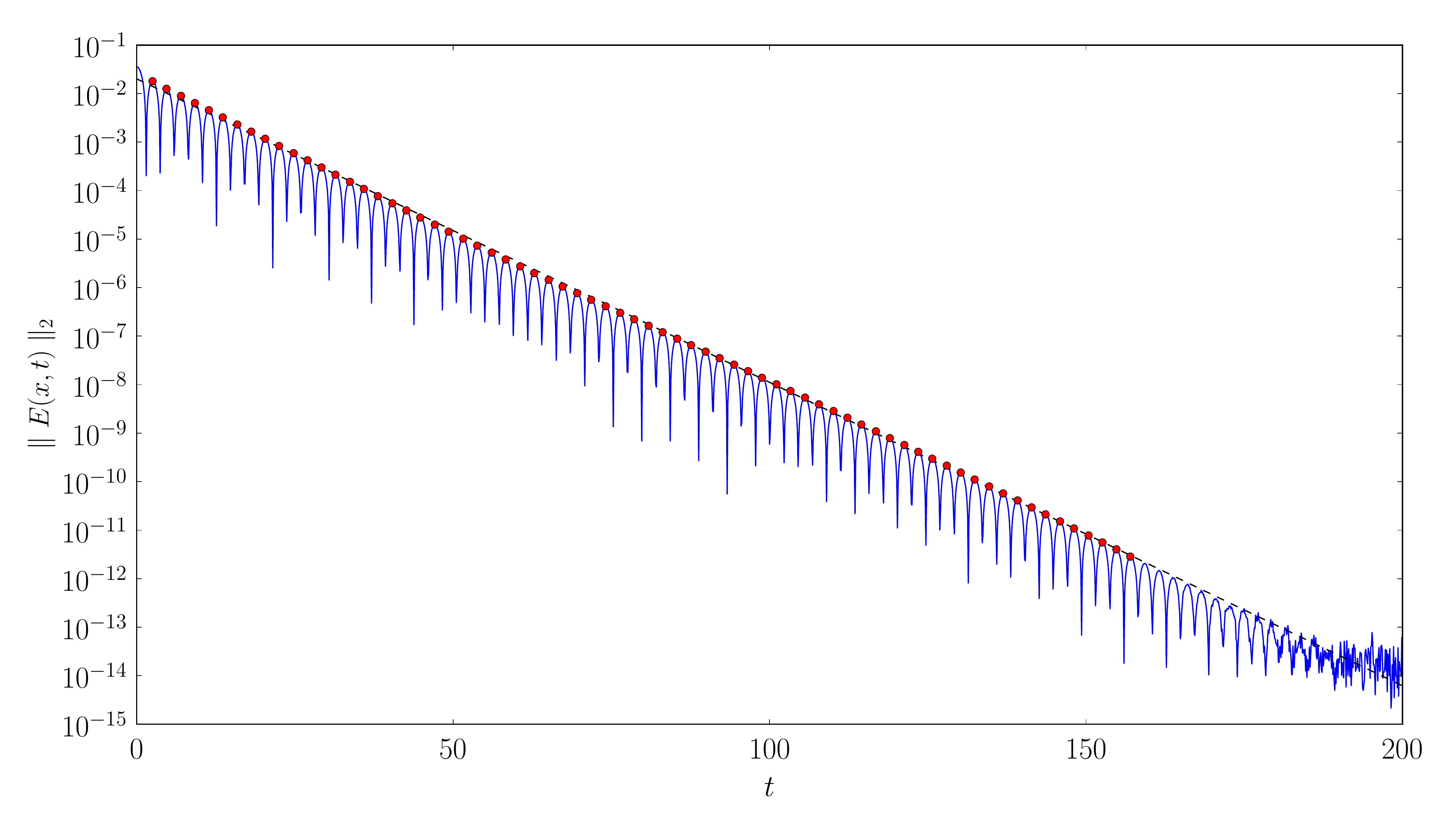}
\caption{Linear Landau damping with collision frequency $\nu = 10^{-4}$. Evolution of the electrostatic energy. The qualitative behaviour of the linear Landau damping is restored by the collisions, but the damping rate, $\gamma = - 0.144$, is to low.}
\label{fig:vlasov_landau_linear_nu1E-4_potential}
\end{figure}

At a collision frequency of $\nu = 10^{-4}$, the error of the conserved quantities is almost optimal (figure \ref{fig:vlasov_landau_linear_nu1E-4_NEP}), and the electrostatic field is damped up to the machine accuracy (figure \ref{fig:vlasov_landau_linear_nu1E-4_potential}).
When measuring the absorption coefficient $\gamma$ from the first timesteps, a good agreement with the theoretical value is observed.
Long-time measurements of $\gamma$ result, however, in values too small compared with the theoretical value. For the marked maxima in figure \ref{fig:vlasov_landau_linear_nu1E-4_potential}, we obtain $\gamma = - 0.144$.
This behaviour is explained by subgrid modes which are not damped completely by the collision operator and are therefore still active.

\begin{figure}
\centering
\includegraphics[width=.95\textwidth]{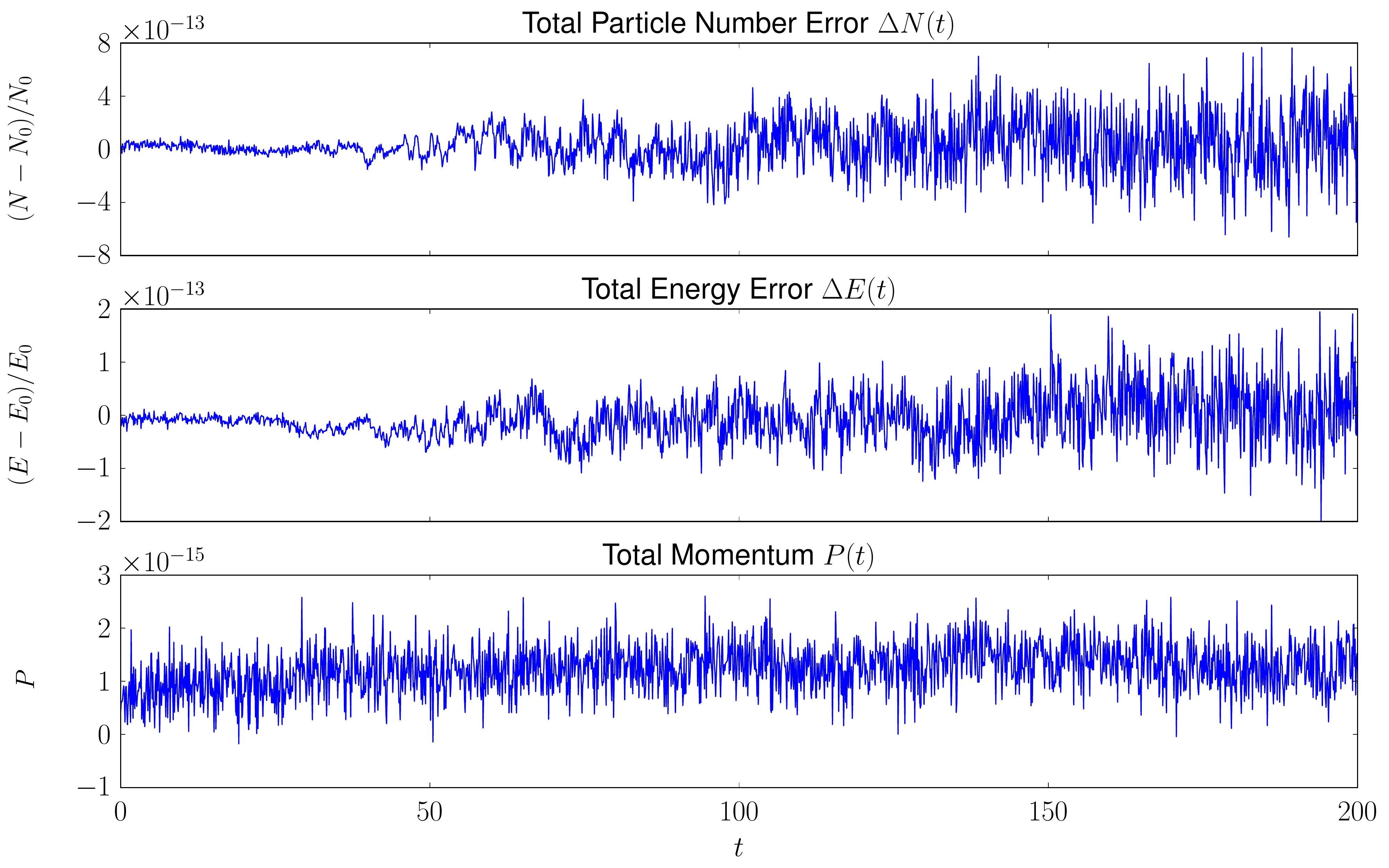}
\caption{Linear Landau damping with collision frequency $\nu = 4 \times 10^{-4}$. With collisions, the total particle number, energy, and linear momentum are all preserved up to machine precision. Hardly any differences to the case with $\nu = 10^{-4}$ are visible.}
\label{fig:vlasov_landau_linear_nu4E-4_NEP}
\end{figure}

To obtain the approximately correct value of $\gamma = - 0.152$, we have to increase the collision frequency to $\nu = 4 \times 10^{-4}$ (figure  \ref{fig:vlasov_landau_linear_nu4E-4_NEP}). Even so, there is almost no visible difference in the time traces of the energy error for $\nu = 1 \times 10^{-4}$ and  $\nu = 4 \times 10^{-4}$, figures \ref{fig:vlasov_landau_linear_nu1E-4_NEP} and \ref{fig:vlasov_landau_linear_nu4E-4_NEP}, respectively, the difference is obvious in the damping of the electrostatic potential.

These results suggests that the linear case, for which an analytical solution is known, can be used to tune the collision frequency for a given step width $h_{v}$ in velocity space. We will see in the following nonlinear examples, that in all cases a collision frequency of $\nu = 4 \times 10^{-4}$ is necessary to obtain accurate conservation of particle number, energy and momentum on long timescales.

\subsubsection{Nonlinear Landau Damping}

In the case of nonlinear Landau damping, $A = 0.5$ and $k = 0.5$, the previous observations manifest more clearly.
The effects of the subgrid modes are much more pronounced as the nonlinear character of the dynamics tends to develop smaller scale structures in phasespace. In particular, the phase mixing that comes along with Landau damping quickly develops very small phasespace structures that cannot be resolved.

\begin{figure}
\centering
\includegraphics[width=.95\textwidth]{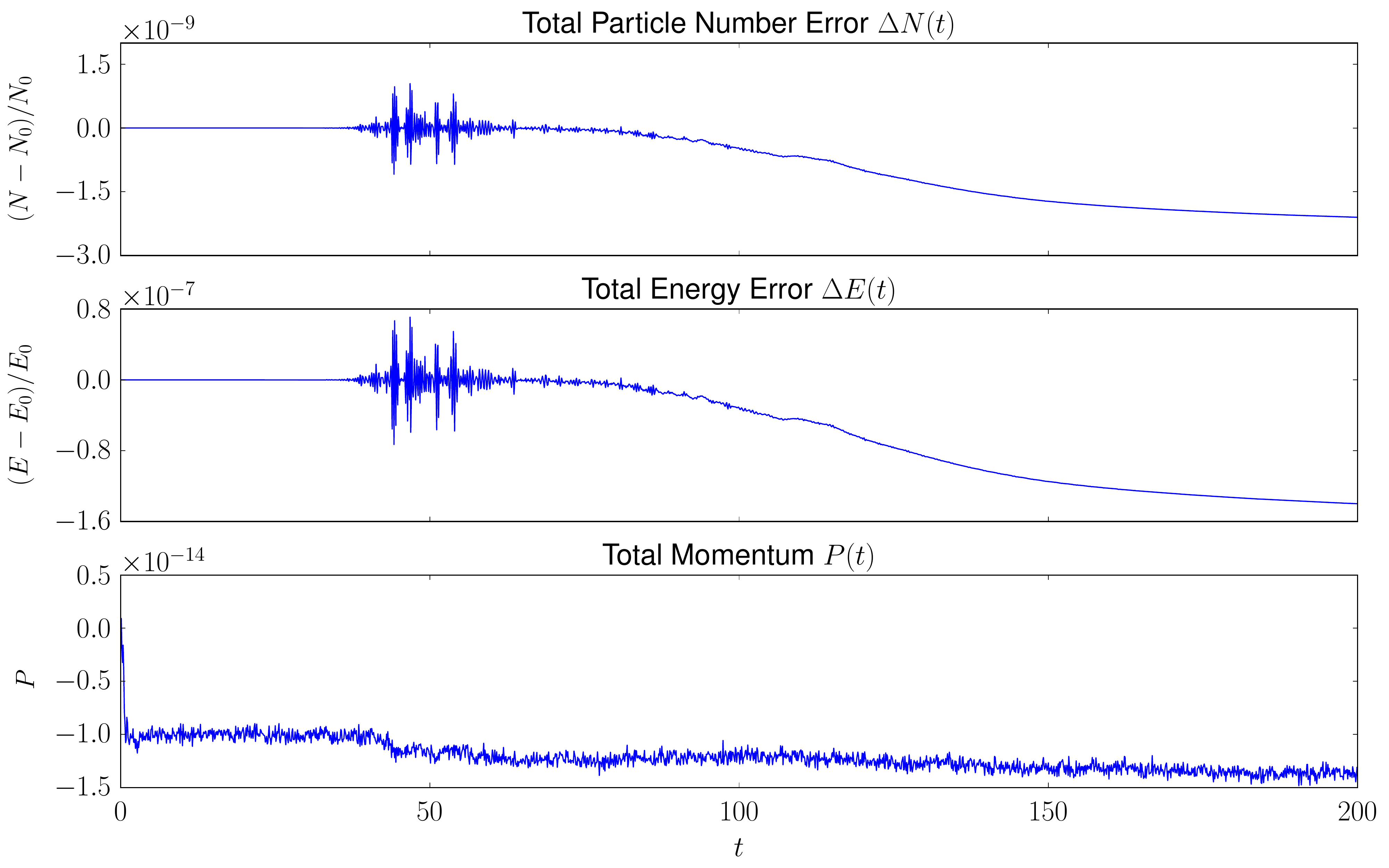}
\caption{Nonlinear Landau damping with collision frequency $\nu = 10^{-4}$. Conservation of the total particle number and energy is violated due to subgrid mode effects. Total linear momentum is preserved exactly.}
\label{fig:vlasov_landau_nonlinear_nu1E-4_NEP}
\end{figure}

\begin{figure}
\centering
\includegraphics[width=.95\textwidth]{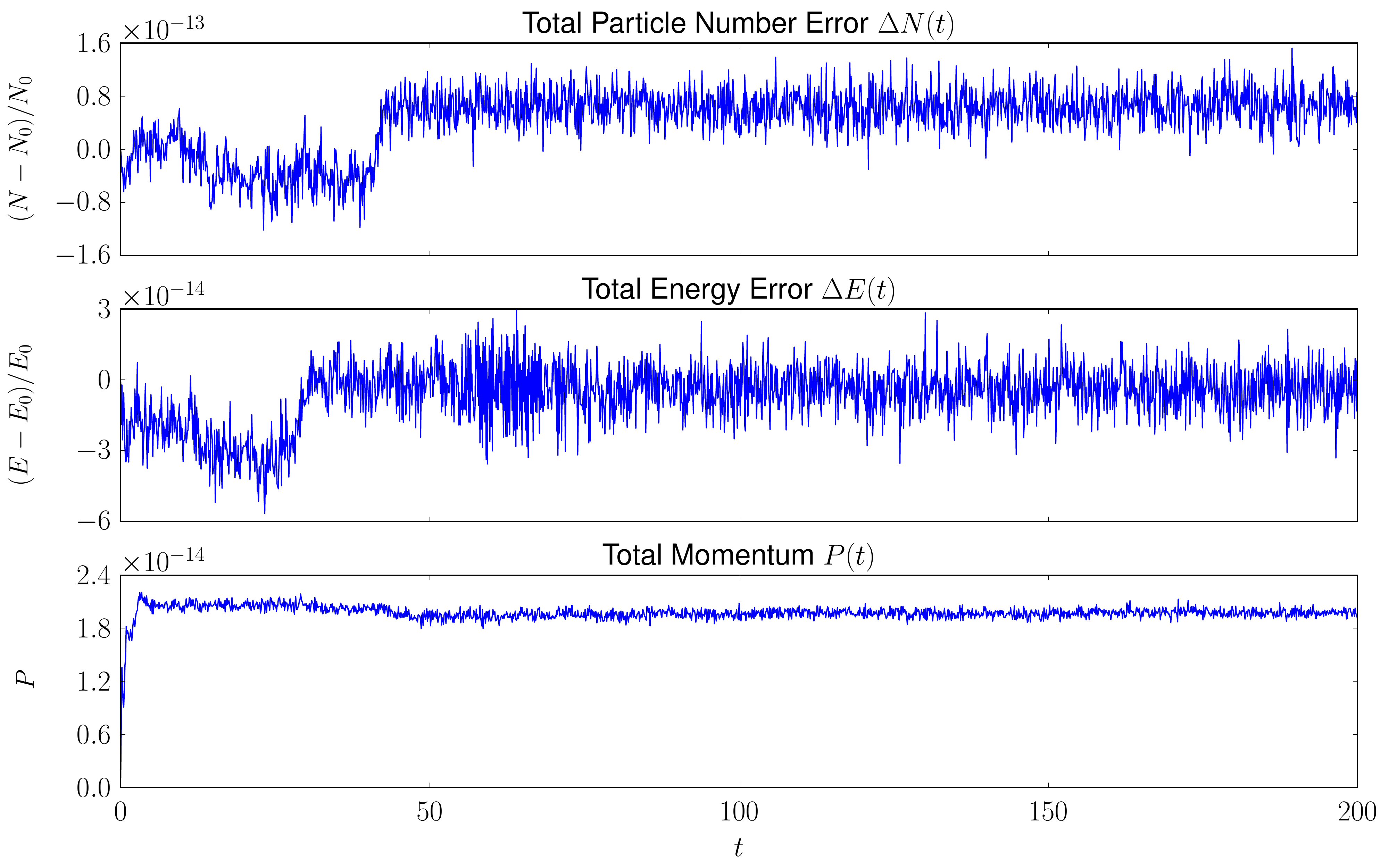}
\caption{Nonlinear Landau damping with collision frequency $\nu = 4 \times 10^{-4}$. Collisions retain the conservation of the total particle number and energy in addition to  exact preservation of the linear momentum.}
\label{fig:vlasov_landau_nonlinear_nu4E-4_NEP}
\end{figure}

With a collision frequency of $\nu = 10^{-4}$, the conservation of energy and the total particle number is severely violated. Only the error in the linear momentum is very small (see figure \ref{fig:vlasov_landau_nonlinear_nu1E-4_NEP}).
In contrast to the previous example, 
here a larger contribution to the error seems to come from the kinetic part.
The relatively large error in the total particle number is directly reflected in the error of the kinetic energy.
With a collision frequency of $\nu = 4 \times 10^{-4}$, instead, conservation of the total particle number, energy as well as the total linear momentum is optimal, see figure \ref{fig:vlasov_landau_nonlinear_nu4E-4_NEP}.

\begin{figure}
\centering
\includegraphics[width=\textwidth]{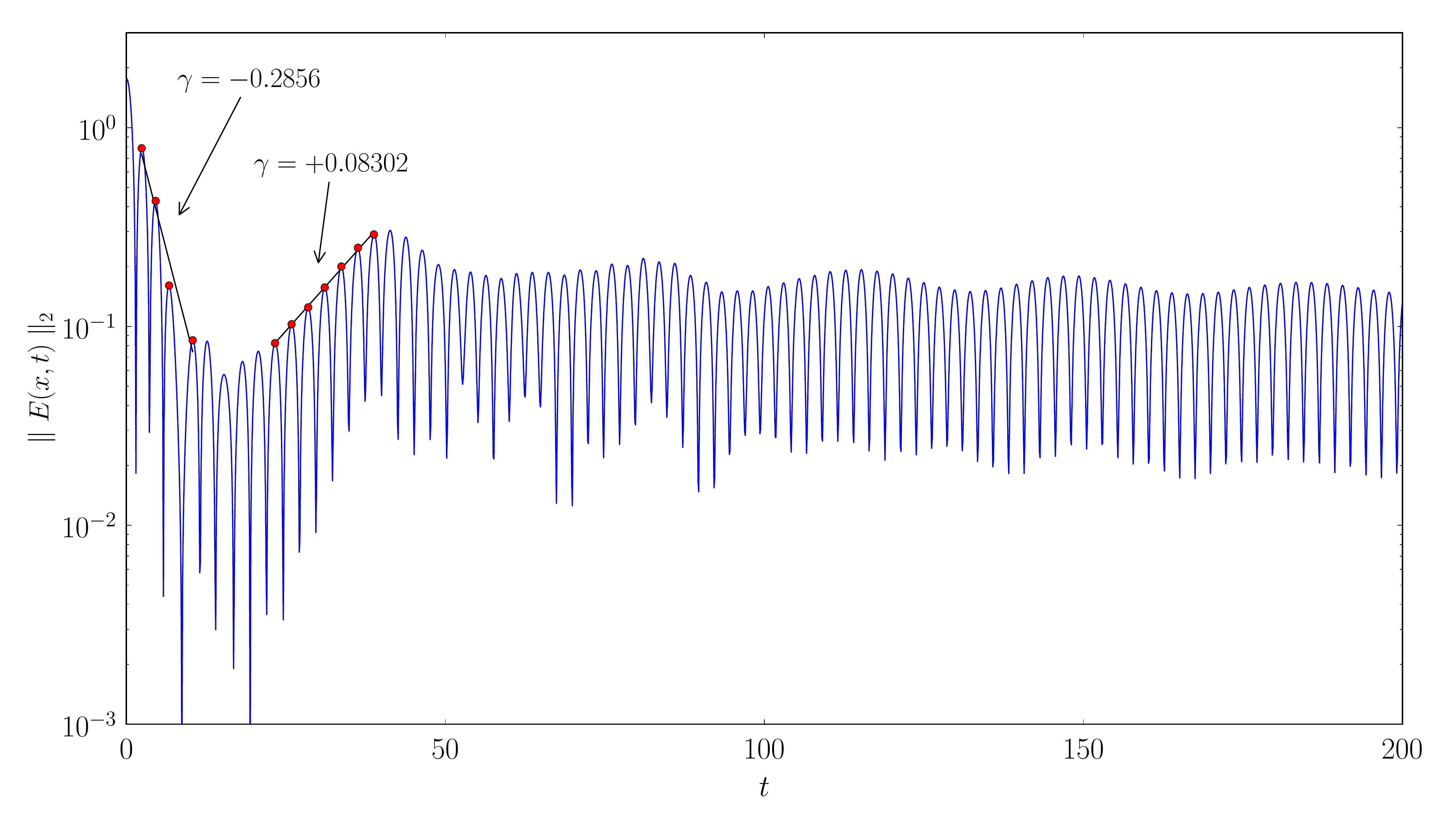}
\caption{Nonlinear Landau damping with collision frequency $\nu = 10^{-4}$. Evolution of the electrostatic energy.}
\label{fig:vlasov_landau_nonlinear_nu1E-4_potential}
\end{figure}

\begin{figure}
\centering
\includegraphics[width=\textwidth]{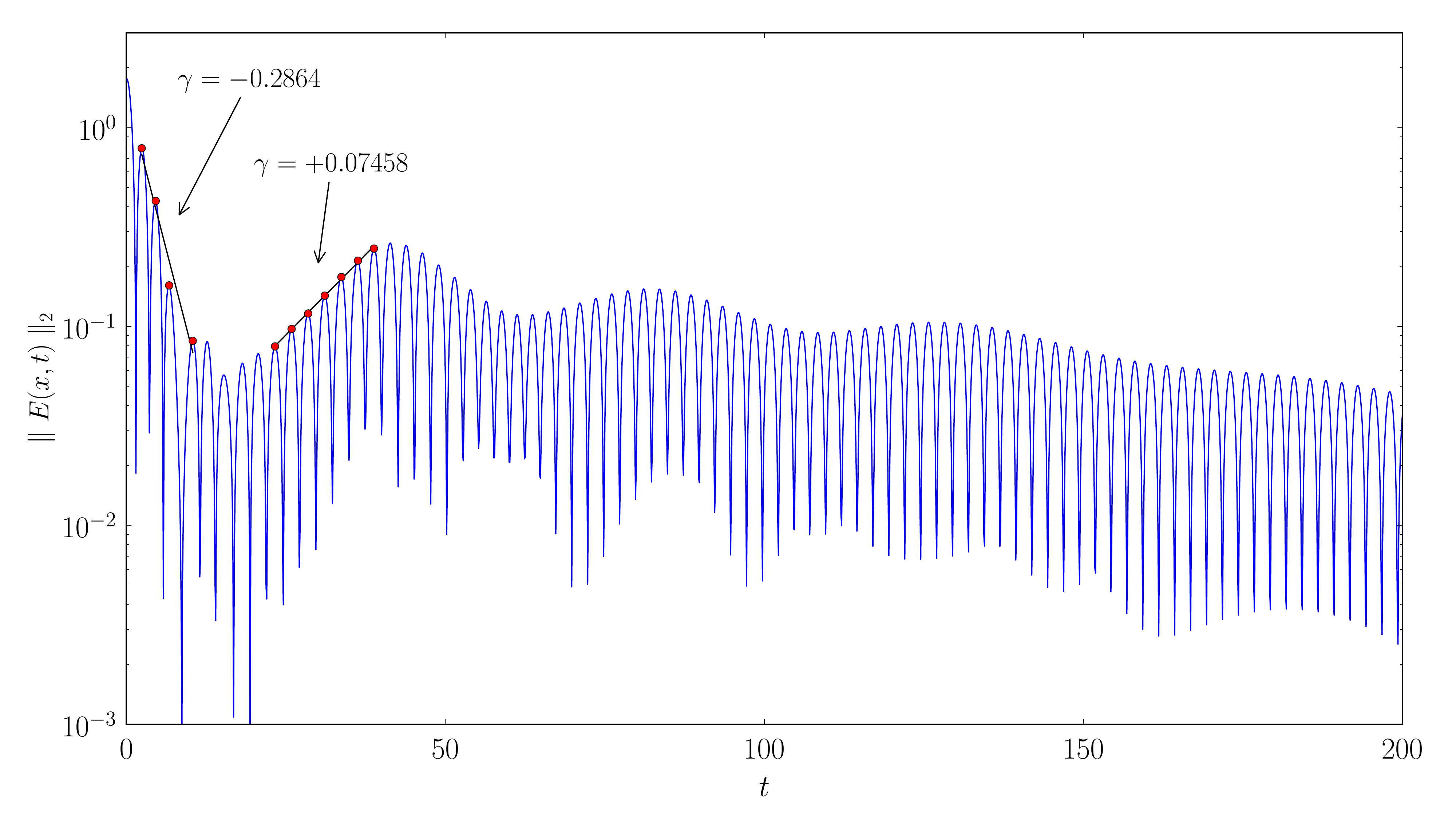}
\caption{Nonlinear Landau damping with collision frequency $\nu = 4 \times 10^{-4}$. Evolution of the electrostatic energy.}
\label{fig:vlasov_landau_nonlinear_nu4E-4_potential}
\end{figure}

The initial damping rate is hardly changed by the collisions. For $\nu=0$ we find $\gamma_{1} = - 0.2854$ (no figure), for $\nu = 10^{-4}$ we find $\gamma_{1} = - 0.2856$ (figure \ref{fig:vlasov_landau_nonlinear_nu1E-4_potential}) and for $\nu = 4 \times 10^{-4}$ we find $\gamma_{1} = - 0.2864$ (figure \ref{fig:vlasov_landau_nonlinear_nu4E-4_potential}). All numbers are equal to two digits and compare well with the existing literature, e.g., \citeauthor{ChengKnorr:1976} computed $\gamma_{1} = - 0.281$, \citeauthor{NakamuraYabe:1999} computed $\gamma_{1} = - 0.280$, and \citeauthor{HeathGamba:2012} computed $\gamma_{1} = - 0.287$.

The effect of the collisions on the second phase, where the electrostatic potential is growing again, is more pronounced. Without collisions, we obtain $\gamma_{2} = 0.0860$, with $\nu = 10^{-4}$ we obtain $\gamma_{2} = 0.0830$ and for $\nu = 4 \times 10^{-4}$ the we find the growth rate to be $\gamma_{2} = 0.0746$ and hence considerably reduced.
The results of \citeauthor{ChengKnorr:1976}, who computed $\gamma_{2} = 0.084$, and \citeauthor{NakamuraYabe:1999}, who computed $\gamma_{2} = 0.0845$, are closer to our results with less or no collisions. The result of \citeauthor{HeathGamba:2012}, $\gamma_{2} = 0.0746$, on the other side, matches ours exactly (for a summary see table \ref{tab:vlasov_nonlinear_landau_damping}).

Not surprisingly, the collisions damp the electrostatic field, and more so for larger collision frequencies $\nu$. The important question is whether they are just removing unphysical contributions to the field energy that originate from subgrid modes, or whether they damp the field too much.
As already described, the relatively large error in the energy and particle number for $\nu = 10^{-4}$ suggests that subgrid modes are not sufficiently damped. Therefore the electrostatic potential is likely to be affected as well, such that part of the electrostatic energy is due to subgrid modes and therefore unphysical. Consequently, the electrostatic energy is likely to be overestimated in that case. On the other hand, it cannot be anticipated that the smaller growth and consecutively stronger damping for $\nu = 4 \times 10^{-4}$ is closer to the real situation as we are counteracting a numerical effect with an effective collision operator, which neither allows us to draw conclusions for the collisionless case, nor represents the physical collision process. Furthermore, the electrostatic energy is much larger in the nonlinear case than in the linear case. Therefore the error in the total energy can not be attributed to the kinetic or the potential part without ambiguity. The correlation between the
errors in the total particle number and the total energy suggests, however, that the energy error arises mainly from the kinetic energy.
A definite conclusion is not possible, but simulations with higher resolution and higher order integrators should indicate which effects are physical and which are numerical.

\begin{figure}
\centering
\subfloat{
\includegraphics[width=.48\textwidth]{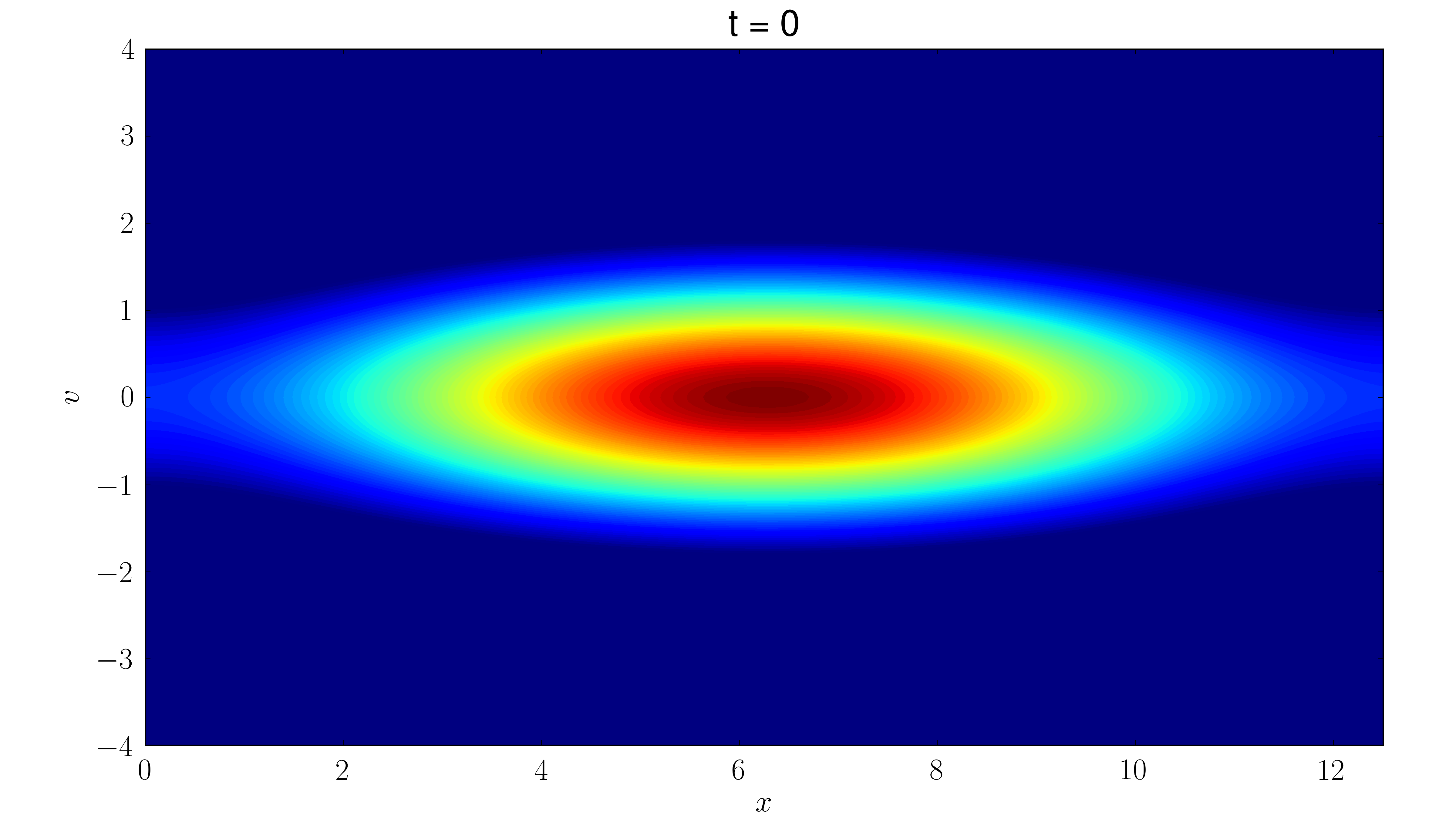}
}
\subfloat{
\includegraphics[width=.48\textwidth]{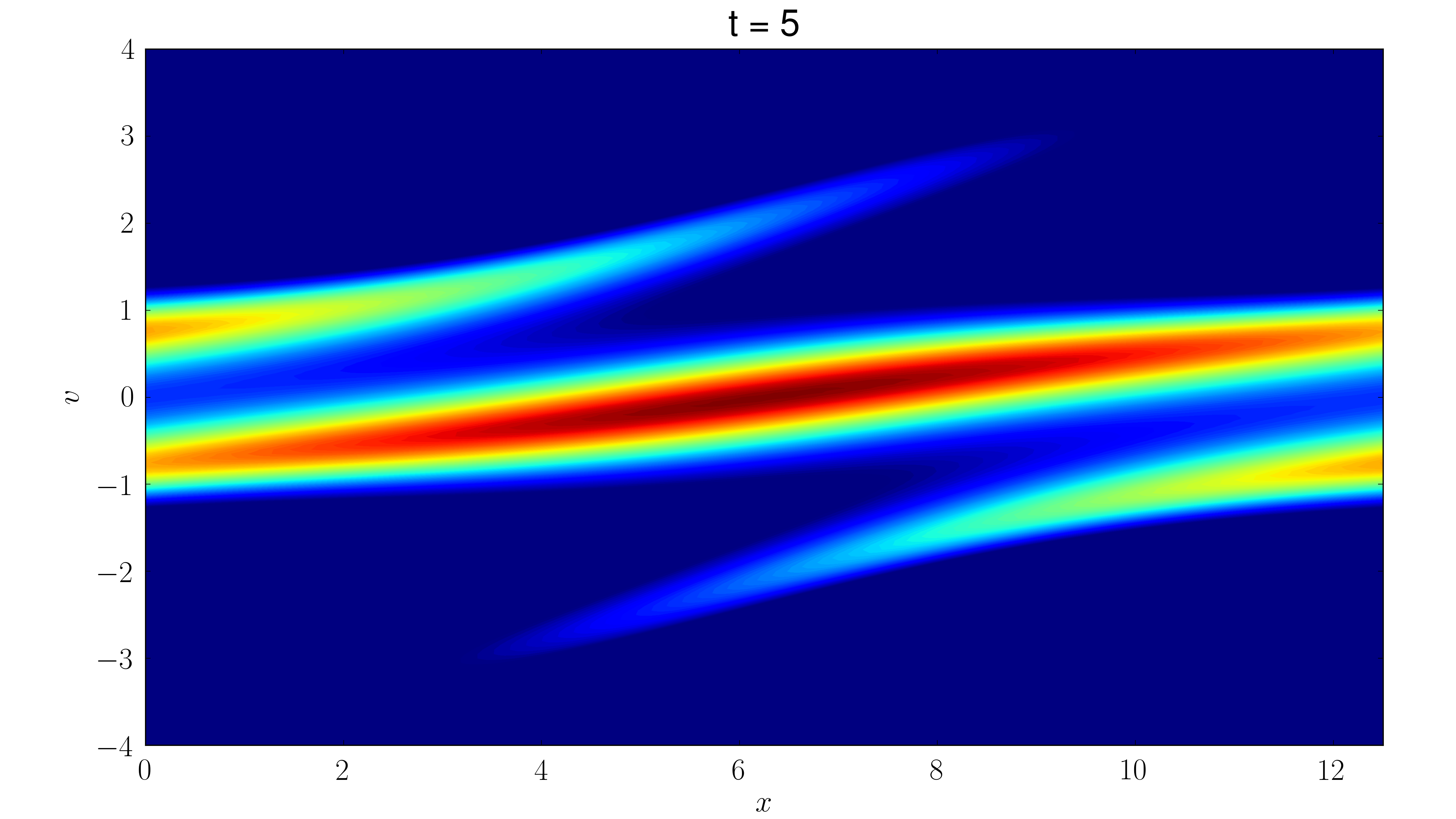}
}

\subfloat{
\includegraphics[width=.48\textwidth]{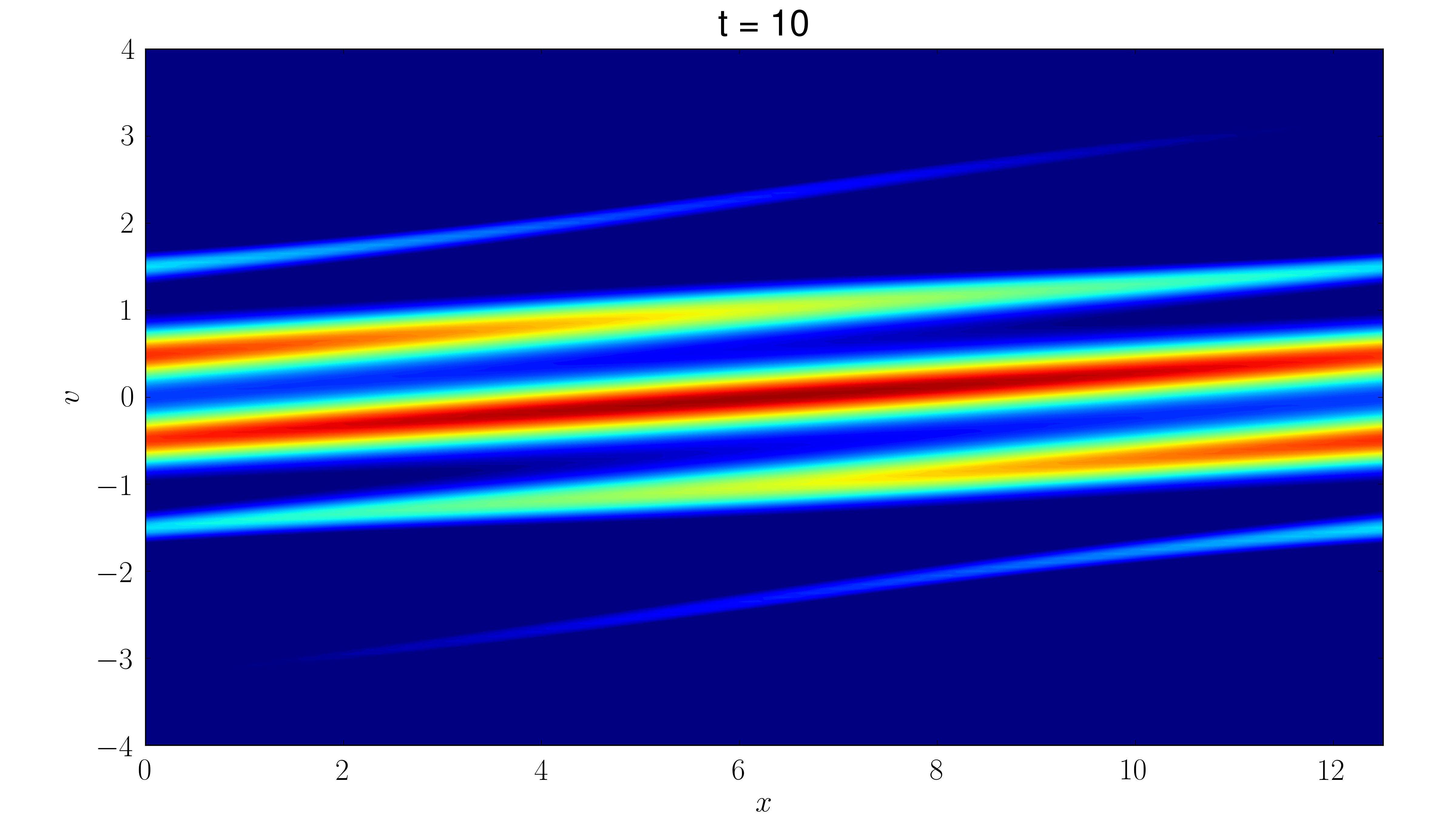}
}
\subfloat{
\includegraphics[width=.48\textwidth]{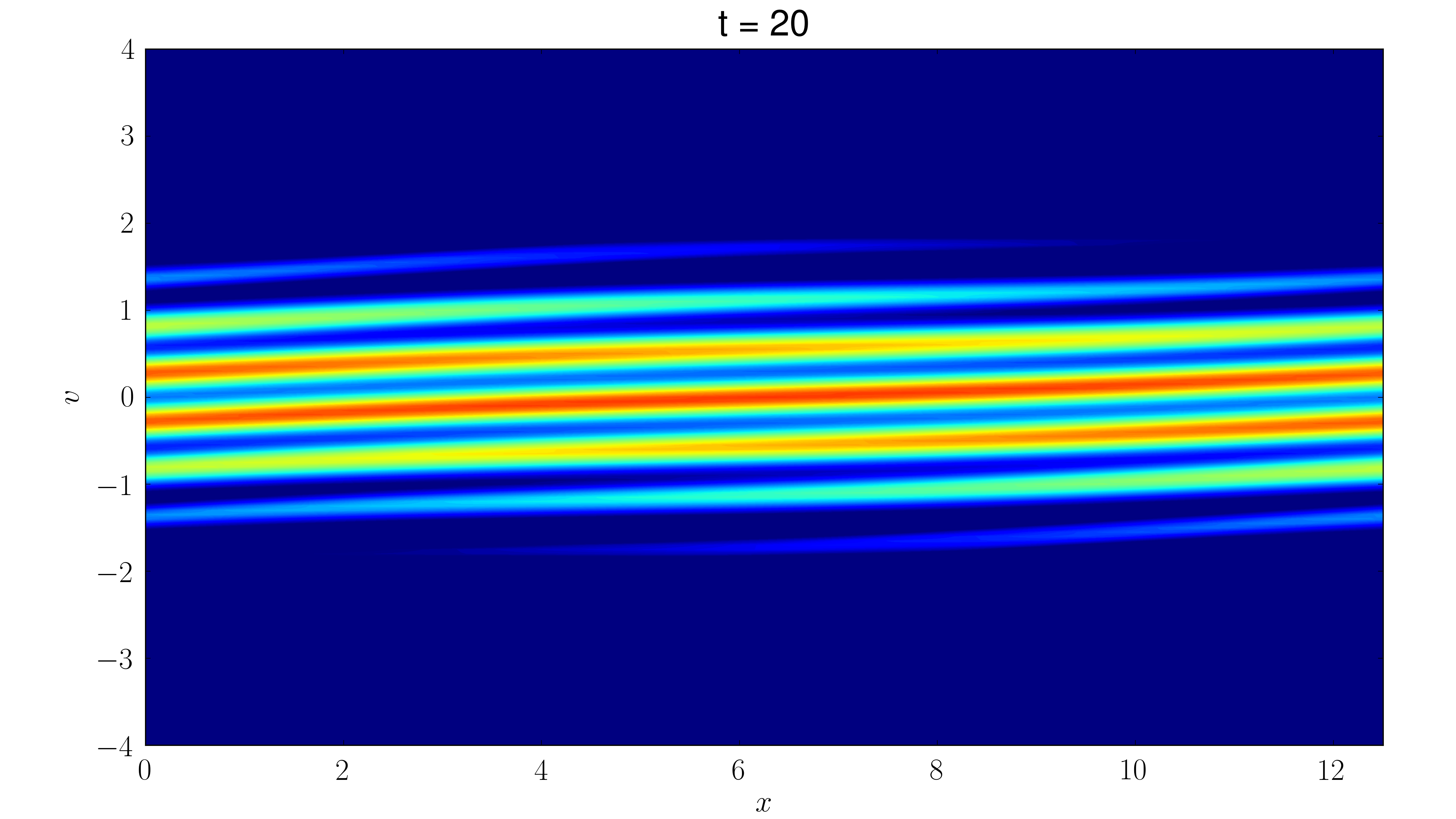}
}

\subfloat{
\includegraphics[width=.48\textwidth]{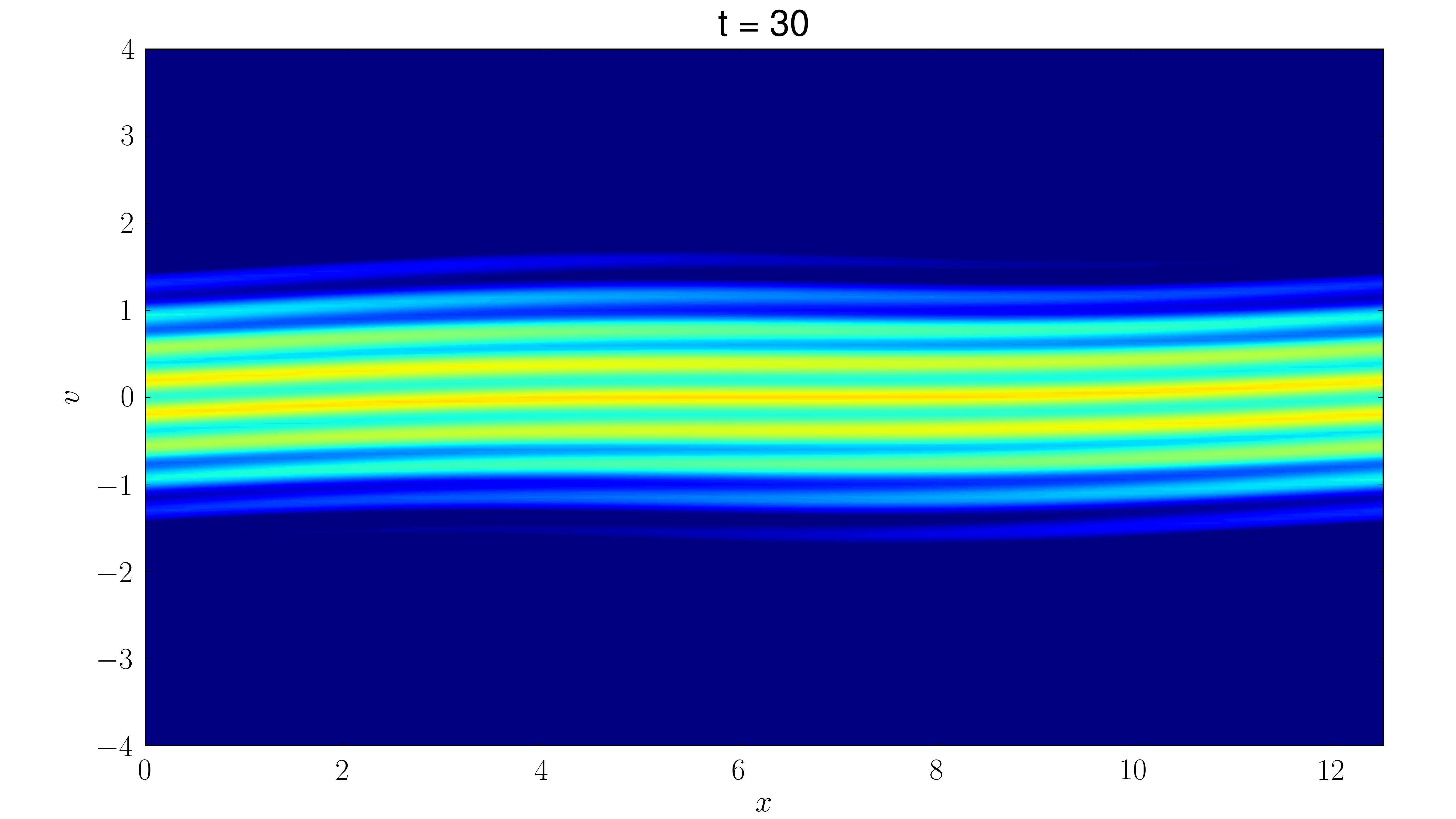}
}
\subfloat{
\includegraphics[width=.48\textwidth]{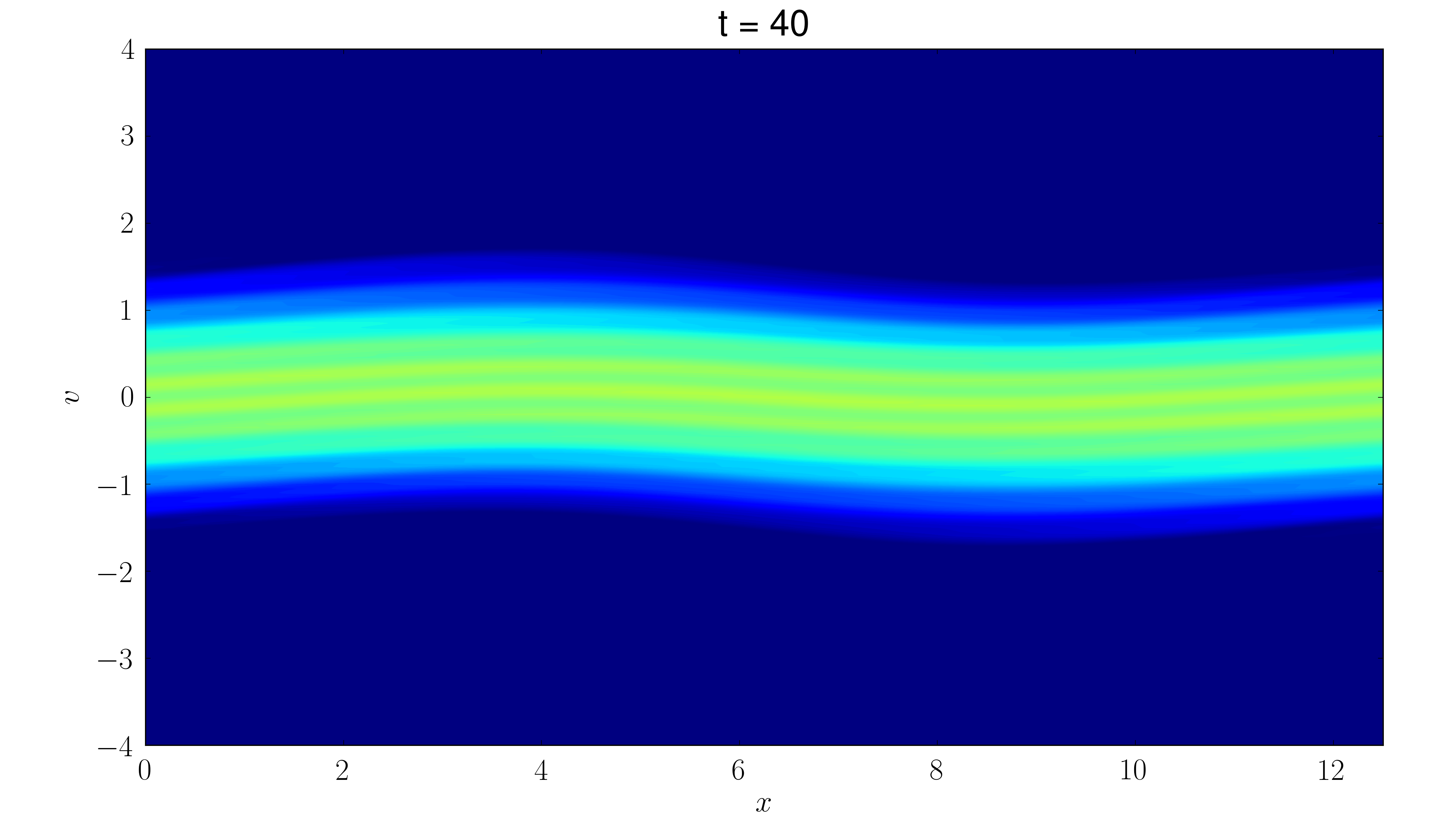}
}

\subfloat{
\includegraphics[width=.48\textwidth]{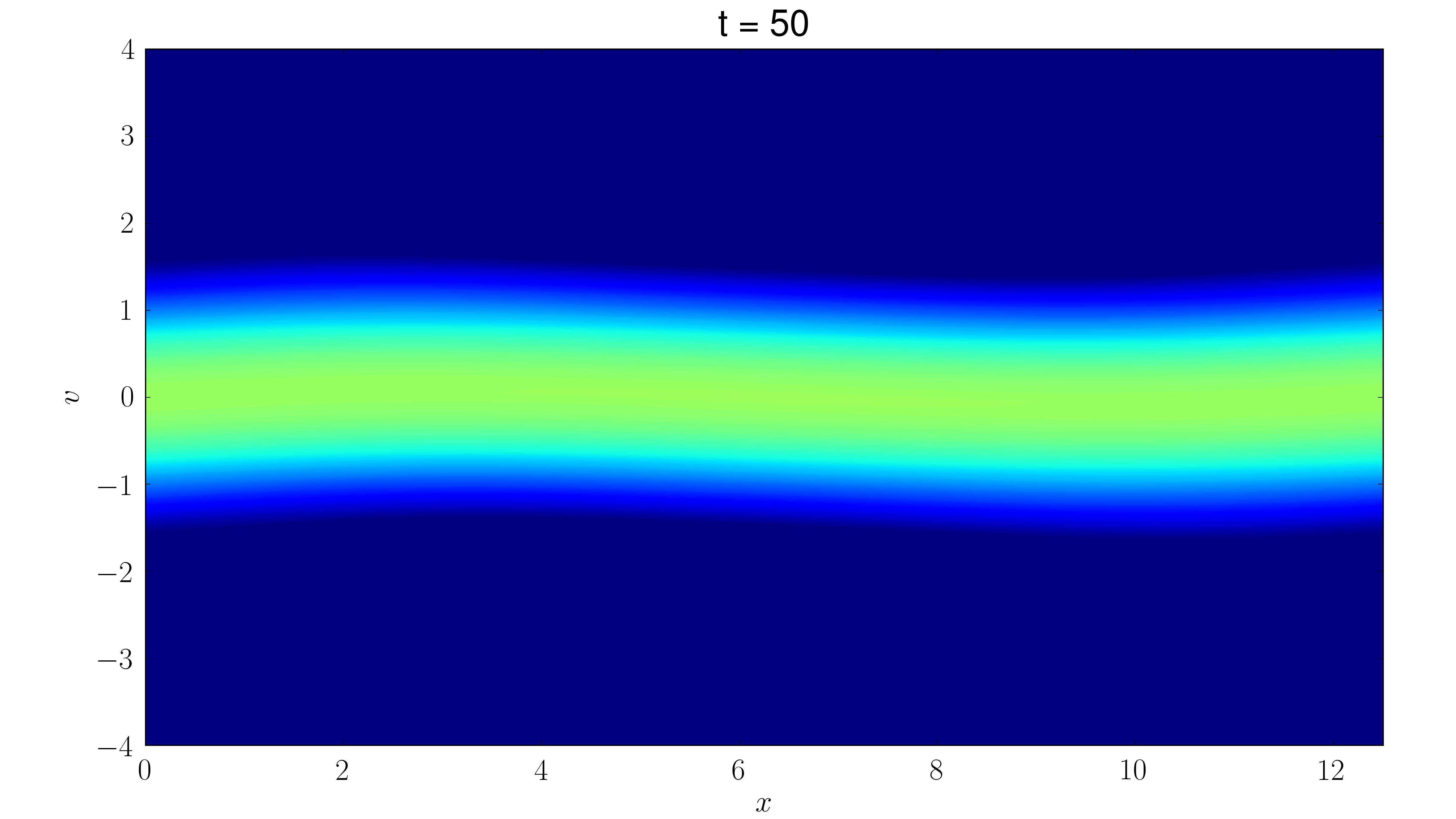}
}
\subfloat{
\includegraphics[width=.48\textwidth]{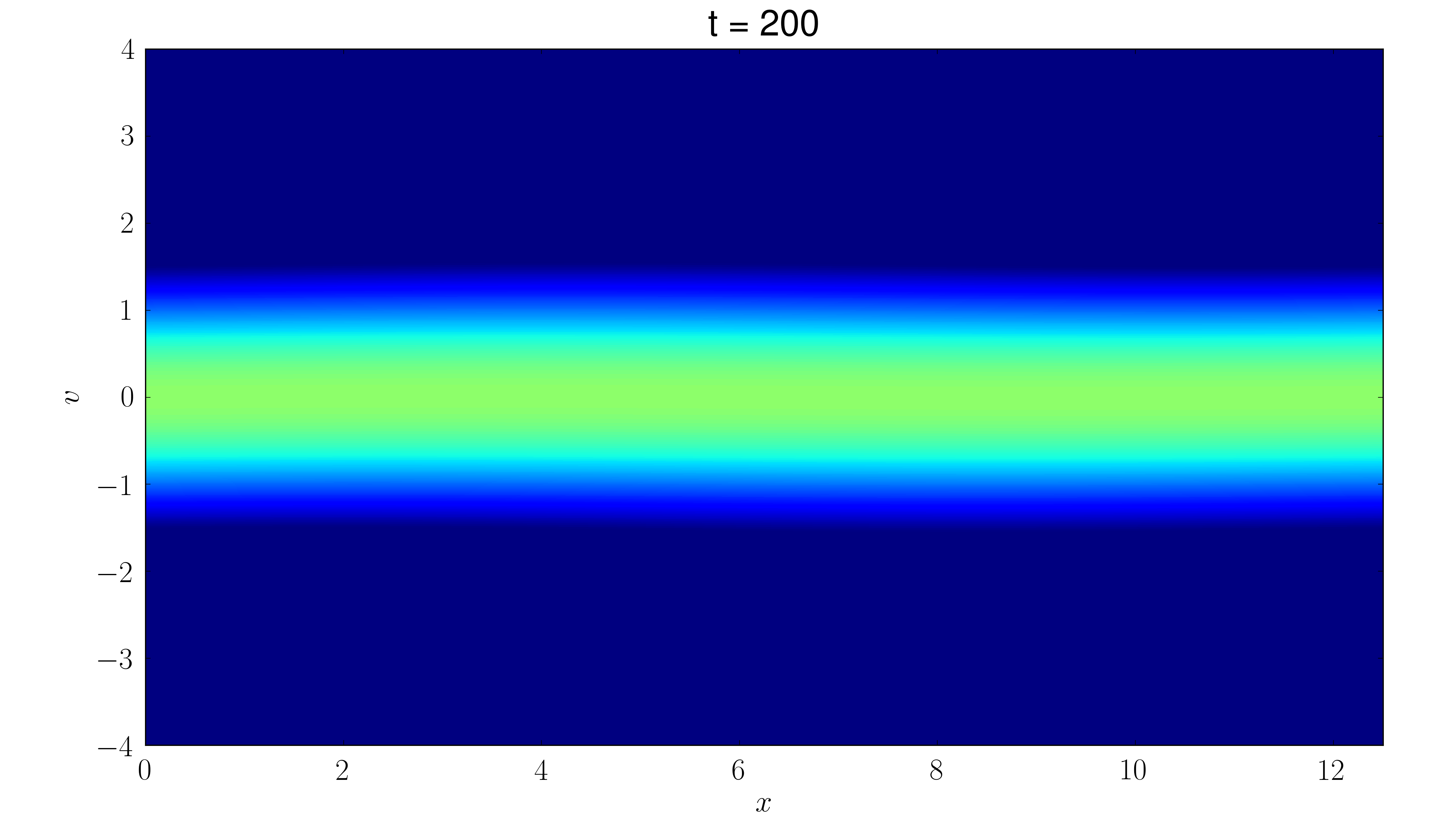}
}
\caption{Nonlinear Landau damping with collision frequency $\nu = 4 \times 10^{-4}$. Contours of the distribution function in phasespace. Contours are linear and constant.}
\label{fig:vlasov_landau_nonlinear_nu4E-4_F}
\end{figure}

In figure \ref{fig:vlasov_landau_nonlinear_nu4E-4_F}, the time evolution of the distribution function is plotted. The phase mixing is nicely visible as is the action of the collision operator. At about $t = 30$ the phasespace structures start to become too small to be resolvable and get therefore damped by the collisions. At about $t = 50$, the fine scale structures have disappeared almost completely, but a large scale oscillation is still visible. At $t = 200$ this oscillation has been further damped, such that it does not appear in the plot anymore.

\subsubsection{Linear Integrator}

For both, linear and nonlinear Landau damping, we also did simulations with the linear integrator from section \ref{sec:vlasov_linearised}. The simulation parameters are the same, except for the timestep, which was chosen as $h_{t} = 0.01$. We did only simulations without collisions, as in the linear scheme, the collision operator can only be treated explicitly (as it is inherently nonlinear) and therefore even smaller timesteps would be necessary.
For both, linear and nonlinear Landau damping, we find the same behaviour with the linear method as we do with the fully nonlinear method.

\begin{figure}
\centering
\includegraphics[width=.95\textwidth]{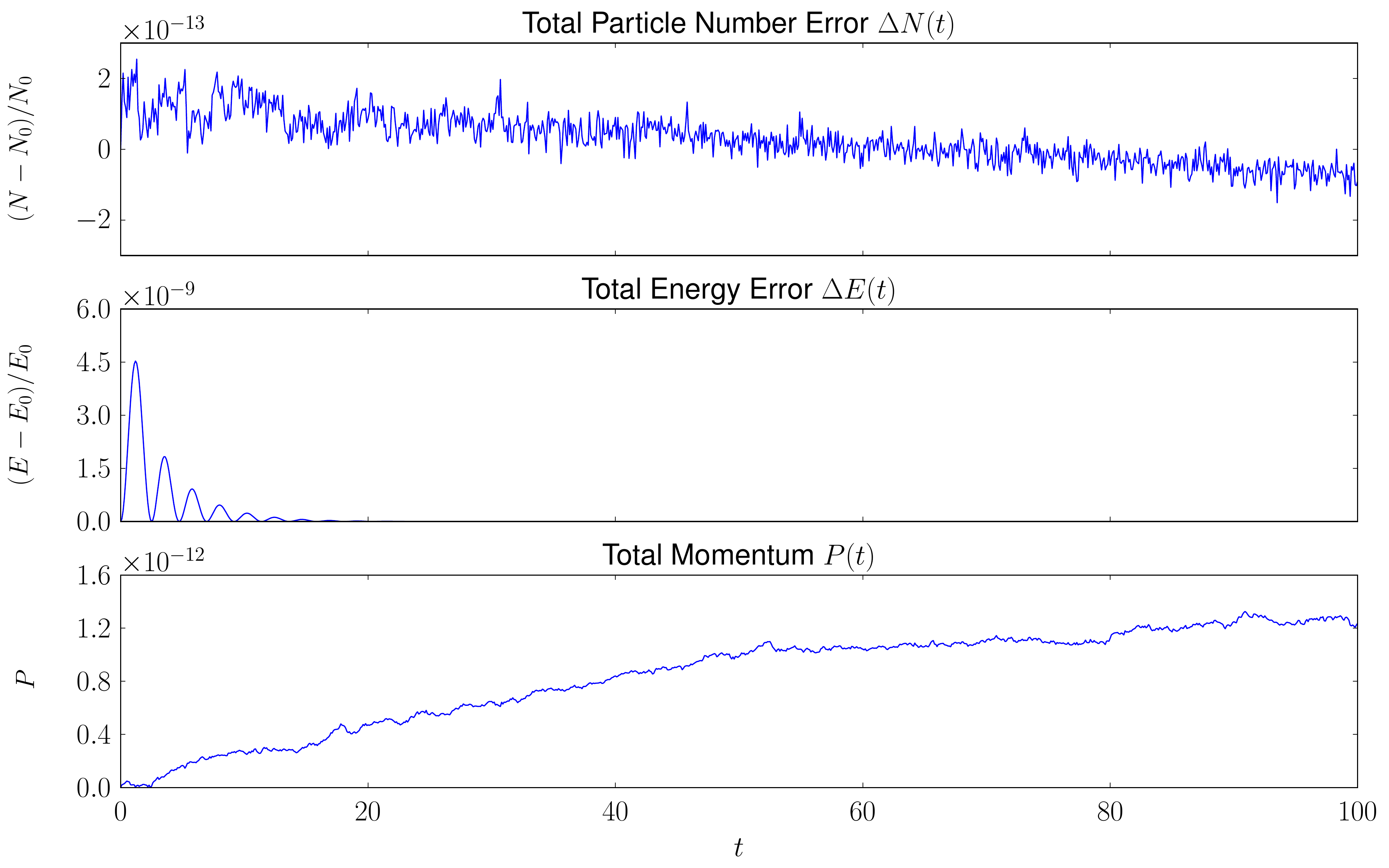}
\caption{Linear Landau damping with linear integrator and without collisions. The total particle number and the total linear momentum are well preserved. Conservation of the total energy is good, but the error is larger than with the nonlinear integrator.}
\label{fig:vlasov_landau_linear_nu0_linear_NEP}
\end{figure}

In the case of linear Landau damping, the total particle number and the total linear momentum are well preserved (see figure \ref{fig:vlasov_landau_linear_nu0_linear_NEP}). The conservation of the total energy is good, but the error is larger than with the nonlinear integrator. In fact, we observe a behaviour of the energy error that is typical for multisymplectic integrators, namely, the energy is not preserved exactly, but its error is bounded, often oscillating, where the amplitude of the oscillation depends on the timestep. This is the reason why we choose a smaller timestep for the simulations with the linear integrator, i.e., to still get good energy conservation.
In the case of linear Landau damping, the amplitude of the oscillation is $\mcal{O} (10^{-9})$. However, after the initial perturbation is damped, it becomes much smaller again.

\begin{figure}
\centering
\includegraphics[width=.95\textwidth]{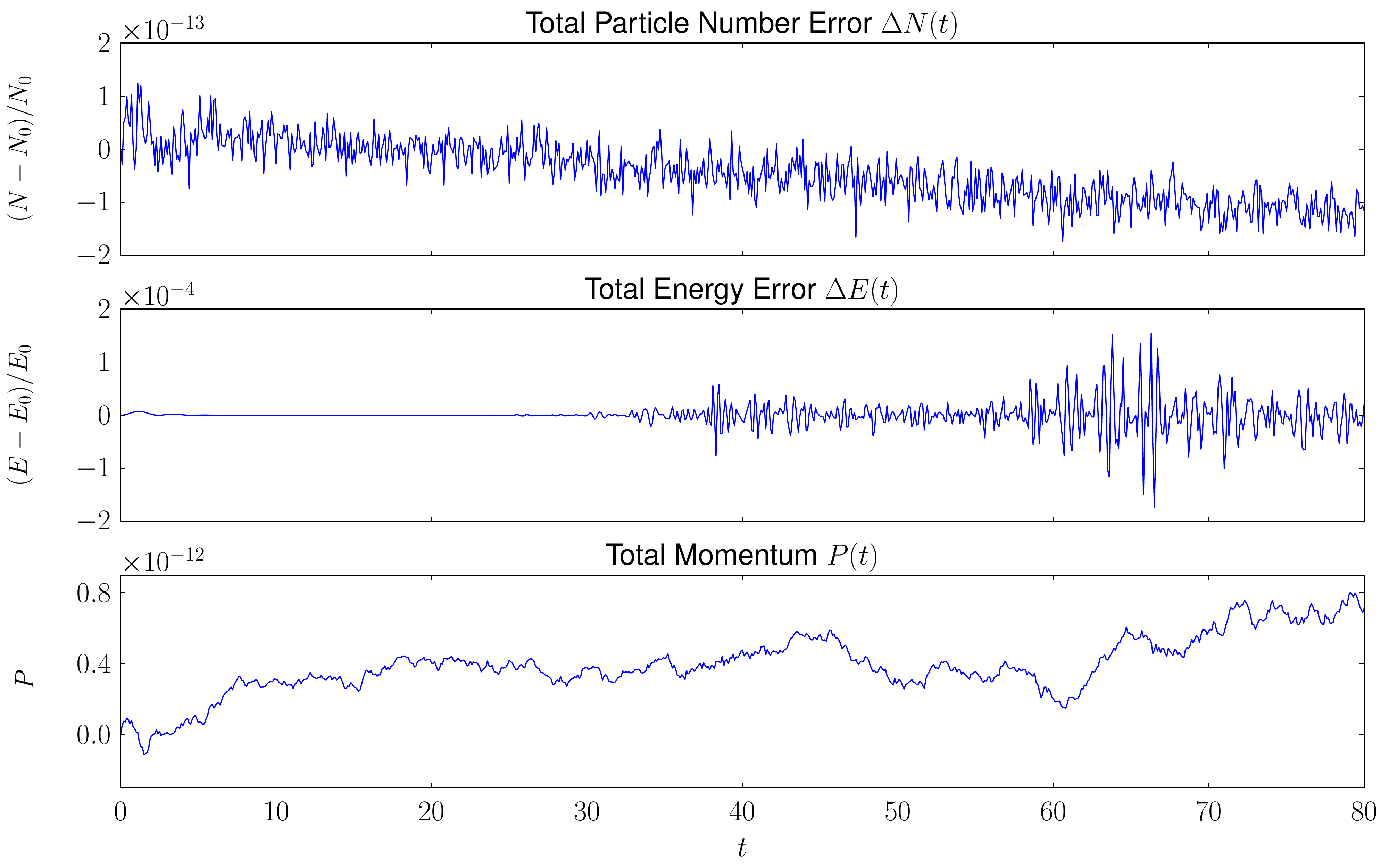}
\caption{Nonlinear Landau damping with linear integrator and without collisions. Total particle number and total linear momentum are well preserved, but the energy error is soon dominated by subgrid mode effects.}
\label{fig:vlasov_landau_nonlinear_nu0_linear_NEP}
\end{figure}

In the case of nonlinear Landau damping, initially, the total particle number and total linear momentum are well preserved and the energy error shows a similar behaviour as in the linear case, albeit with a larger amplitude of the error which is $\mcal{O} (10^{-5})$ (see figure \ref{fig:vlasov_landau_nonlinear_nu0_linear_NEP}). Eventually the energy error grows larger due to subgrid modes, which will also spoil the momentum and particle number conservation when running for longer times.

We see that the linear integrator poses a viable alternative if solving a nonlinear system is not an option.
That the linear integrator is working rather well is probably attributed to the fact that the Vlasov-Poisson system consists of two linear equations. The nonlinearity arises only through the coupling of the two equations.

\begin{table}
\centering

\begin{tabular}{|l|c|c|c|}
\hline
Integrator		& $\nu$				 	& $\gamma_{1}$	& $\gamma_{2}$ \\
\hline
Linear VI		& $0$					& $-0.285$		& $+0.087$ \\
Nonlinear VI	& $0$					& $-0.285$		& $+0.086$ \\
Nonlinear VI	& $1 \times 10^{-4}$	& $-0.286$		& $+0.083$ \\
Nonlinear VI	& $4 \times 10^{-4}$	& $-0.286$		& $+0.075$ \\
\hline
\citeauthor{ChengKnorr:1976}	\cite{ChengKnorr:1976}		& - & $-0.281$ & $+0.084$ \\
\citeauthor{NakamuraYabe:1999}	\cite{NakamuraYabe:1999}	& - & $-0.280$ & $+0.085$ \\
\citeauthor{HeathGamba:2012}	\cite{HeathGamba:2012}		& - & $-0.287$ & $+0.075$ \\
\hline
\end{tabular}
\caption{Damping and growth rates in nonlinear Landau damping simulations with variational integrators and comparison with previous works.}
\label{tab:vlasov_nonlinear_landau_damping}
\end{table}

\subsubsection{$L^{2}$ Norm and Entropy}

\begin{figure}
\centering
\includegraphics[width=.8\textwidth]{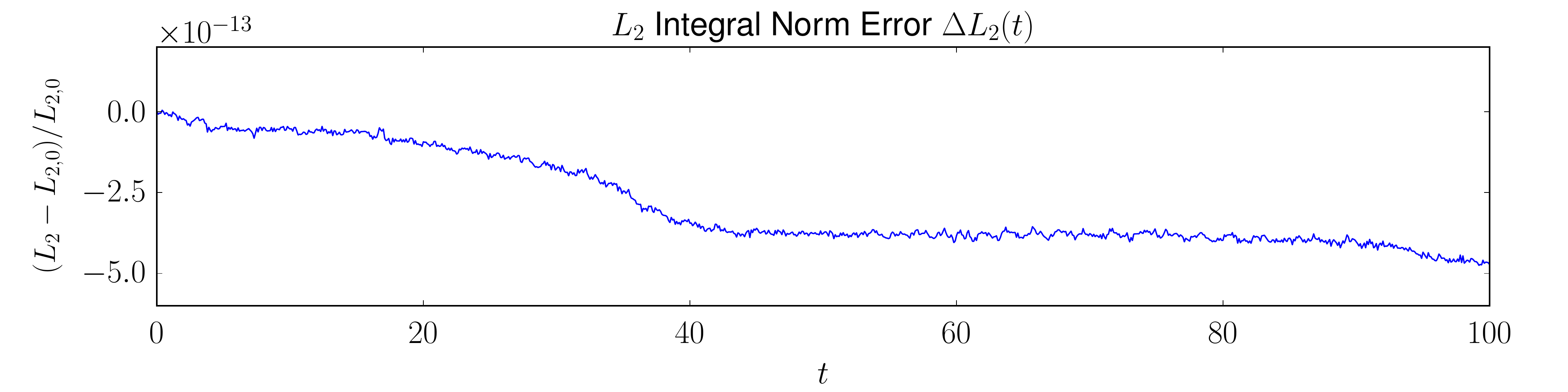}
\includegraphics[width=.8\textwidth]{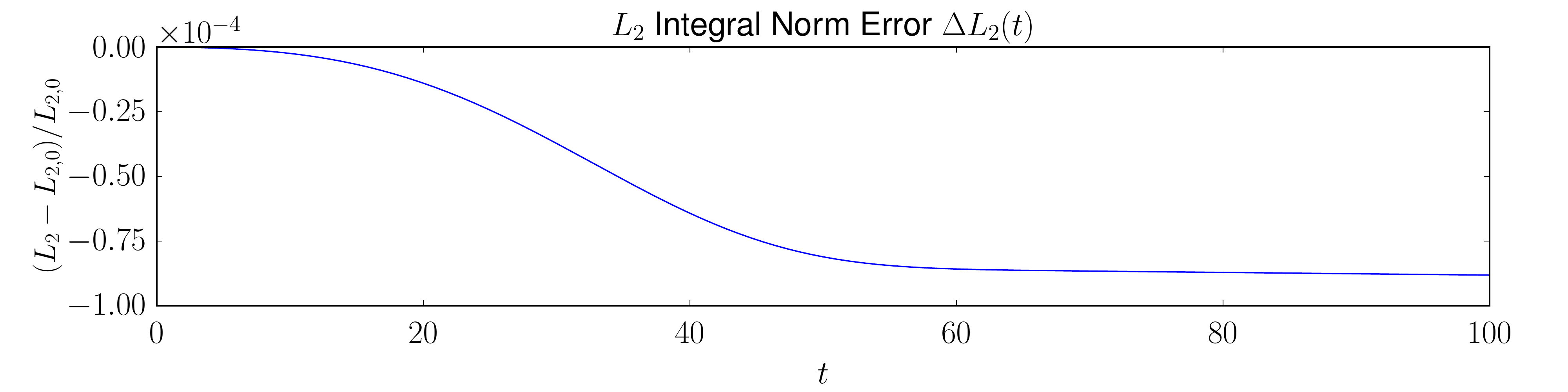}
\includegraphics[width=.8\textwidth]{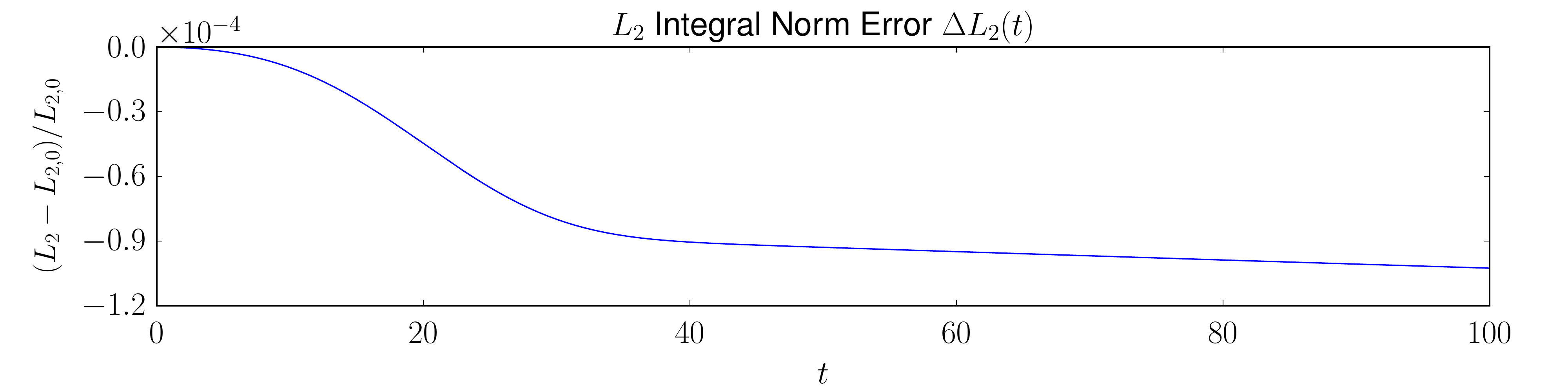}
\caption{Linear Landau damping. Evolution of the $L^{2}$ norm.\\ Top: $\nu = 0$, Middle: $\nu = 10^{-4}$, Bottom: $\nu = 4 \times 10^{-4}$.}
\label{fig:vlasov_landau_linear_L2}
\end{figure}

\begin{figure}
\centering
\includegraphics[width=.8\textwidth]{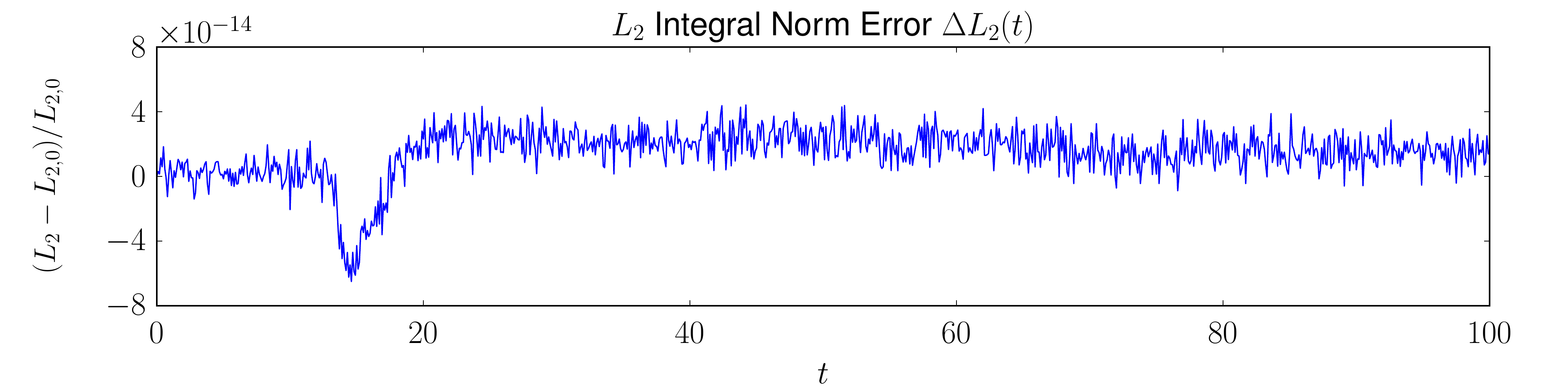}
\includegraphics[width=.8\textwidth]{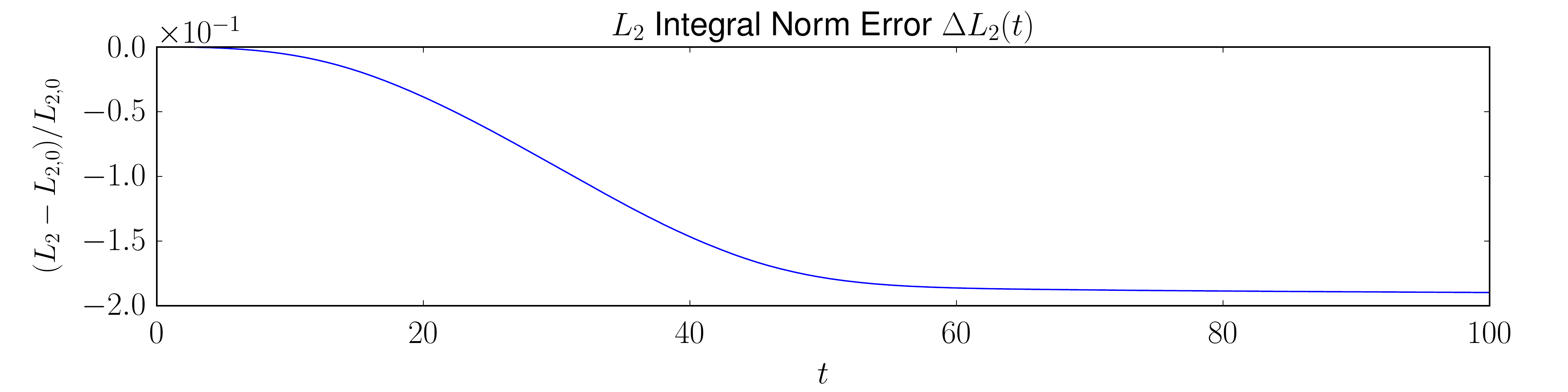}
\includegraphics[width=.8\textwidth]{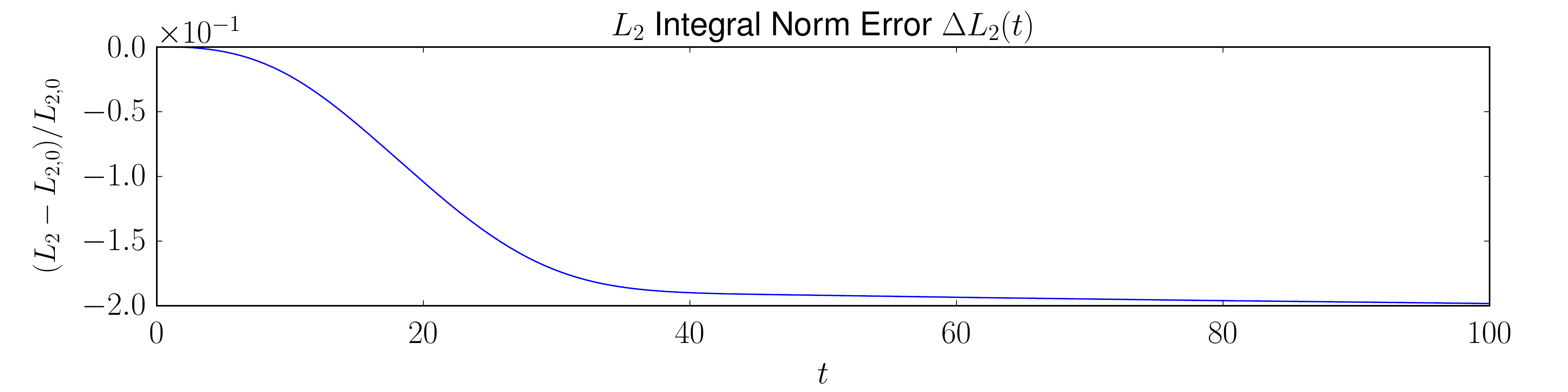}
\caption{Nonlinear Landau damping. Evolution of the $L^{2}$ norm.\\ Top: $\nu = 0$, Middle: $\nu = 10^{-4}$, Bottom: $\nu = 4 \times 10^{-4}$.}
\label{fig:vlasov_landau_nonlinear_L2}
\end{figure}

Before we move to the next example, a comment on the evolution of entropy and the $L^{2}$ norm is in order.
The variational integrator for the Vlasov-Poisson system (\ref{eq:vlasov_vi_discrete_vlasov_equation_reduced}, \ref{eq:vlasov_vi_discrete_poisson_equation_reduced}) preserves the $L^{2}$ norm of the distribution function exactly (see figures \ref{fig:vlasov_landau_linear_L2} and \ref{fig:vlasov_landau_nonlinear_L2}, top).
Through the application of the collision operator, the $L^{2}$ norm is dissipated (figures \ref{fig:vlasov_landau_linear_L2} and \ref{fig:vlasov_landau_nonlinear_L2}, middle and bottom).

\begin{figure}
\centering
\includegraphics[width=.8\textwidth]{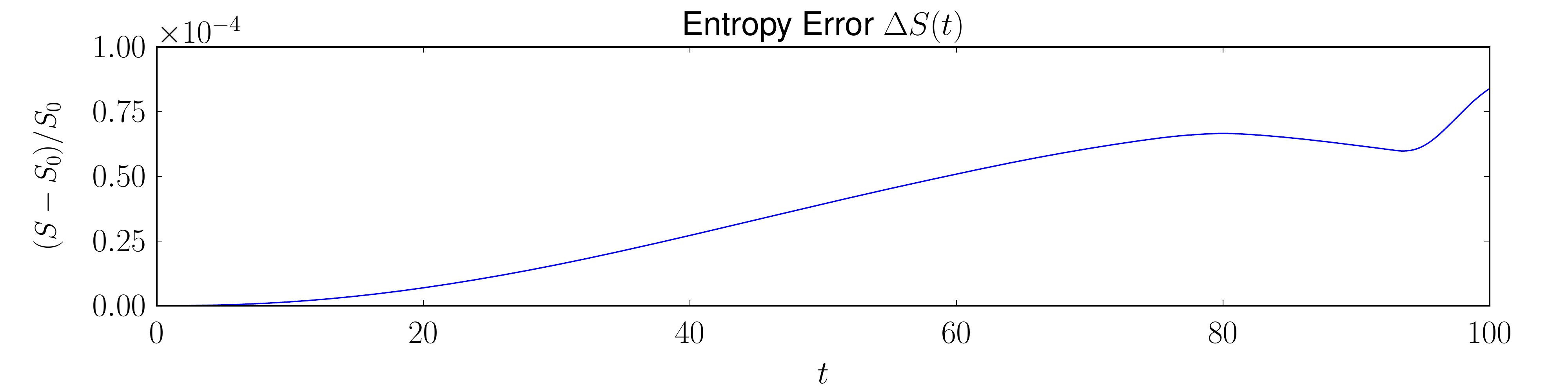}
\includegraphics[width=.8\textwidth]{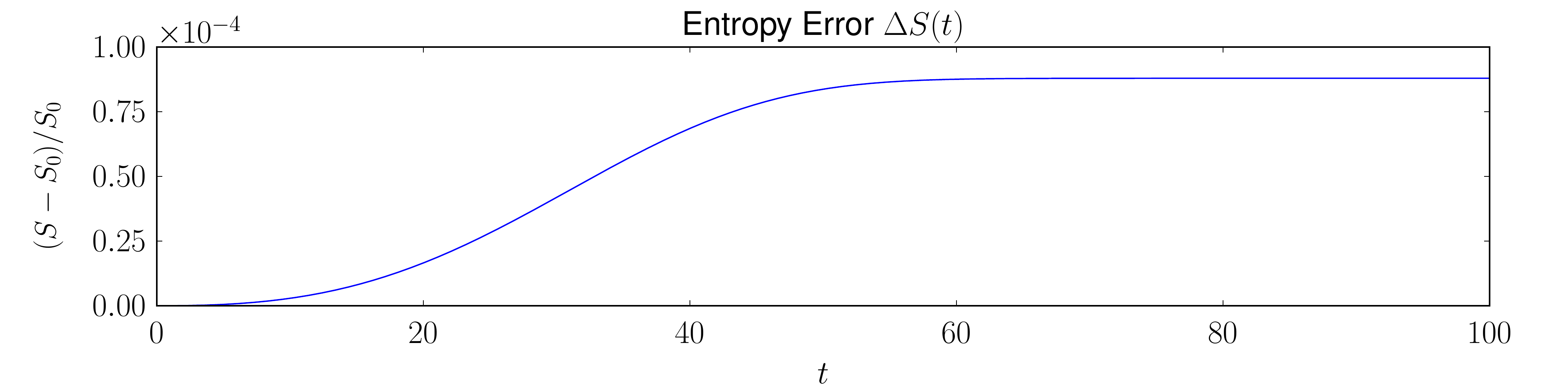}
\includegraphics[width=.8\textwidth]{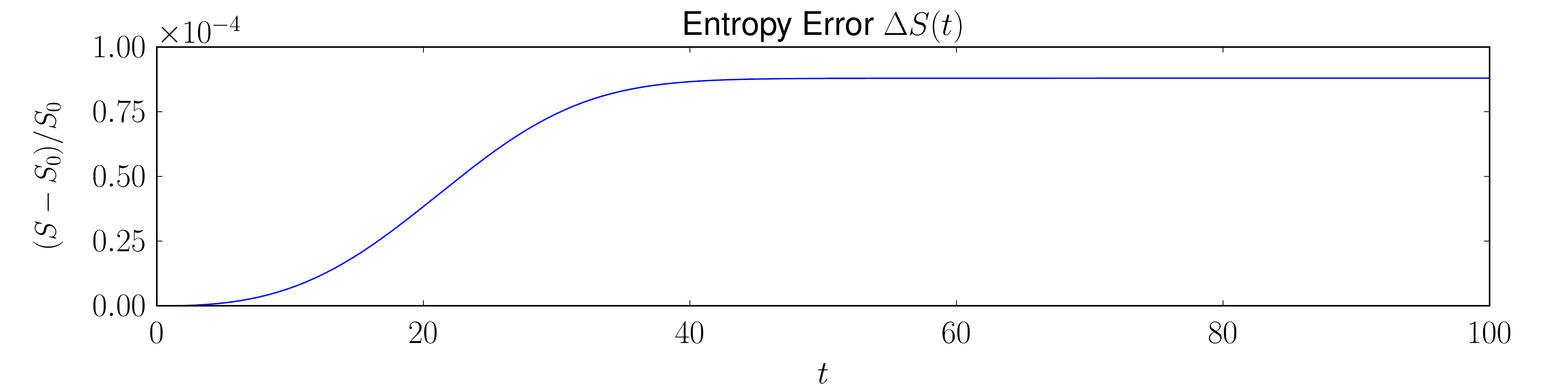}
\caption{Linear Landau damping. Evolution of the entropy $S$.\\ Top: $\nu = 0$, Middle: $\nu = 10^{-4}$, Bottom: $\nu = 4 \times 10^{-4}$.}
\label{fig:vlasov_landau_linear_S}
\end{figure}

\begin{figure}
\centering
\includegraphics[width=.8\textwidth]{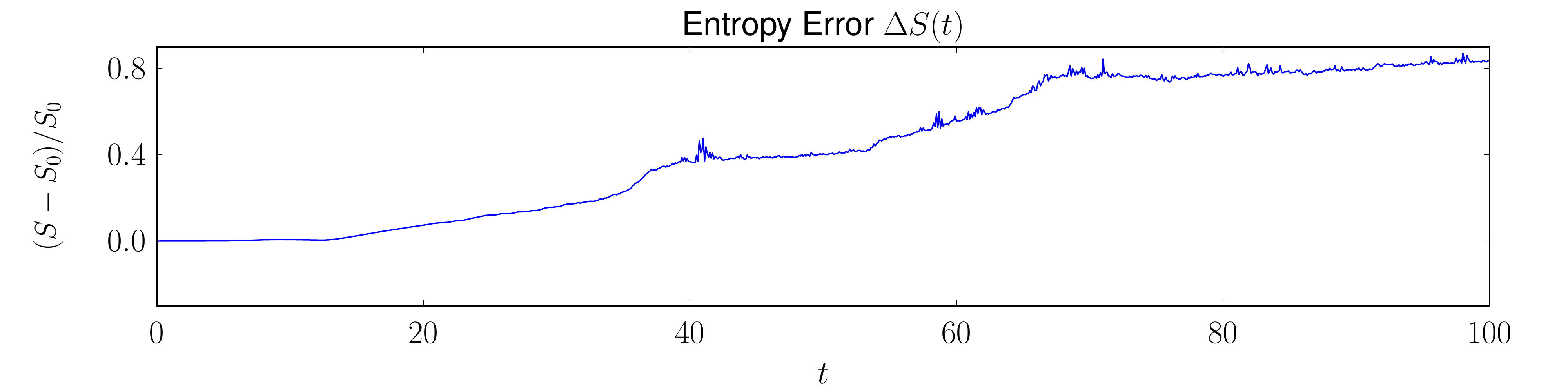}
\includegraphics[width=.8\textwidth]{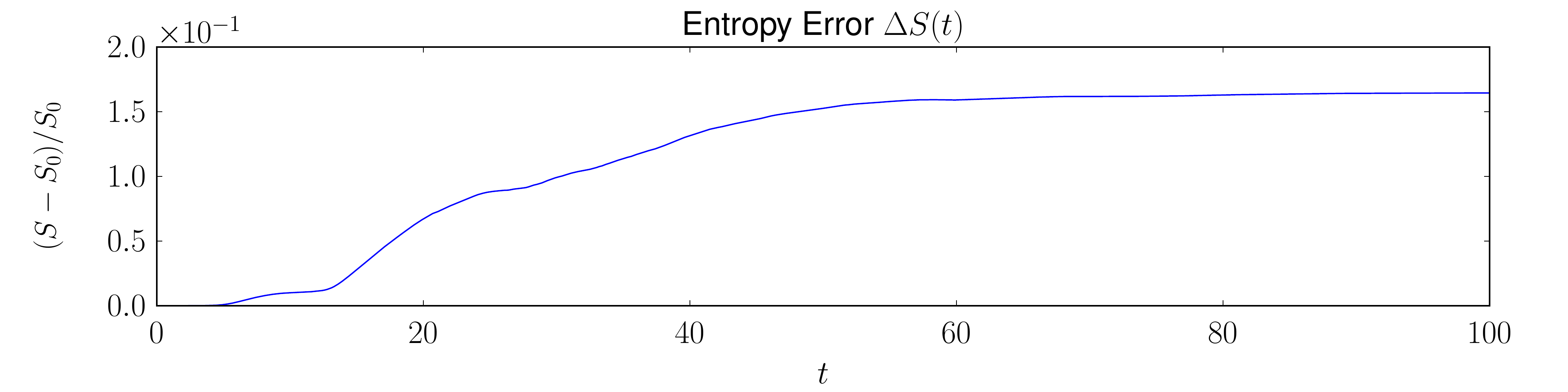}
\includegraphics[width=.8\textwidth]{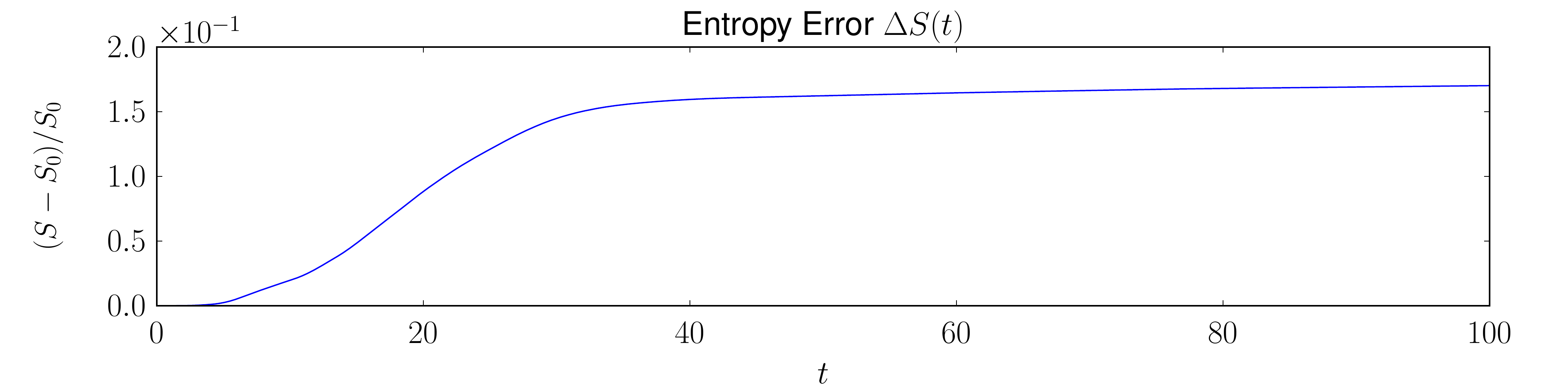}
\caption{Nonlinear Landau damping. Evolution of the entropy $S$.\\ Top: $\nu = 0$, Middle: $\nu = 10^{-4}$, Bottom: $\nu = 4 \times 10^{-4}$.}
\label{fig:vlasov_landau_nonlinear_S}
\end{figure}

Entropy is not conserved by our variational integrator, but without collisions it can be used as a diagnostic for the appearance of subgrid modes. In the case of linear Landau damping, where the effect of the subgrid modes onto the distribution function is rather small, the entropy grows slowly but steadily as can be seen in the top of figure \ref{fig:vlasov_landau_linear_S}.
In the case of nonlinear Landau damping, where the effect of the subgrid modes onto the distribution is more severe, the entropy is growing slowly at first, just as in the linear case, but starting from about $t = 15$ it is growing much more rapidly (top of figure \ref{fig:vlasov_landau_nonlinear_S}).
This sudden growth indicates the appearance of subgrid modes, long-time before they are visible in the energy diagnostics.
In the simulations with collisions (bottom of figures \ref{fig:vlasov_landau_linear_S} and \ref{fig:vlasov_landau_nonlinear_S}) the entropy is increasing more smoothly, not showing such sudden jumps.

\subsection{Twostream Instability}

The distribution function is initialised as
\begin{align}
f (x,v) = v^{2} \, f_{M} (x,v) \, \big( 1 + A \, \cos(kx) \big) ,
\end{align}

with amplitude $A = 0.05$ and wave number $k = 0.5$.
The simulation parameters are the same as before. The spatial domain is $[0, 4\pi]$, $v_{\text{max}} = 10$, $n_{x} = 201$, $n_{v} = 401$, $h_{t} = 0.1$.
The collision frequency is taken to be either $\nu = 0$ or $\nu = 4 \times 10^{-4}$, following the tuning by linear Landau damping, described in the previous section.
The distribution function describes two particle beams, propagating in opposite direction and having a small perturbation imposed on them.

\begin{figure}
\centering
\includegraphics[width=.95\textwidth]{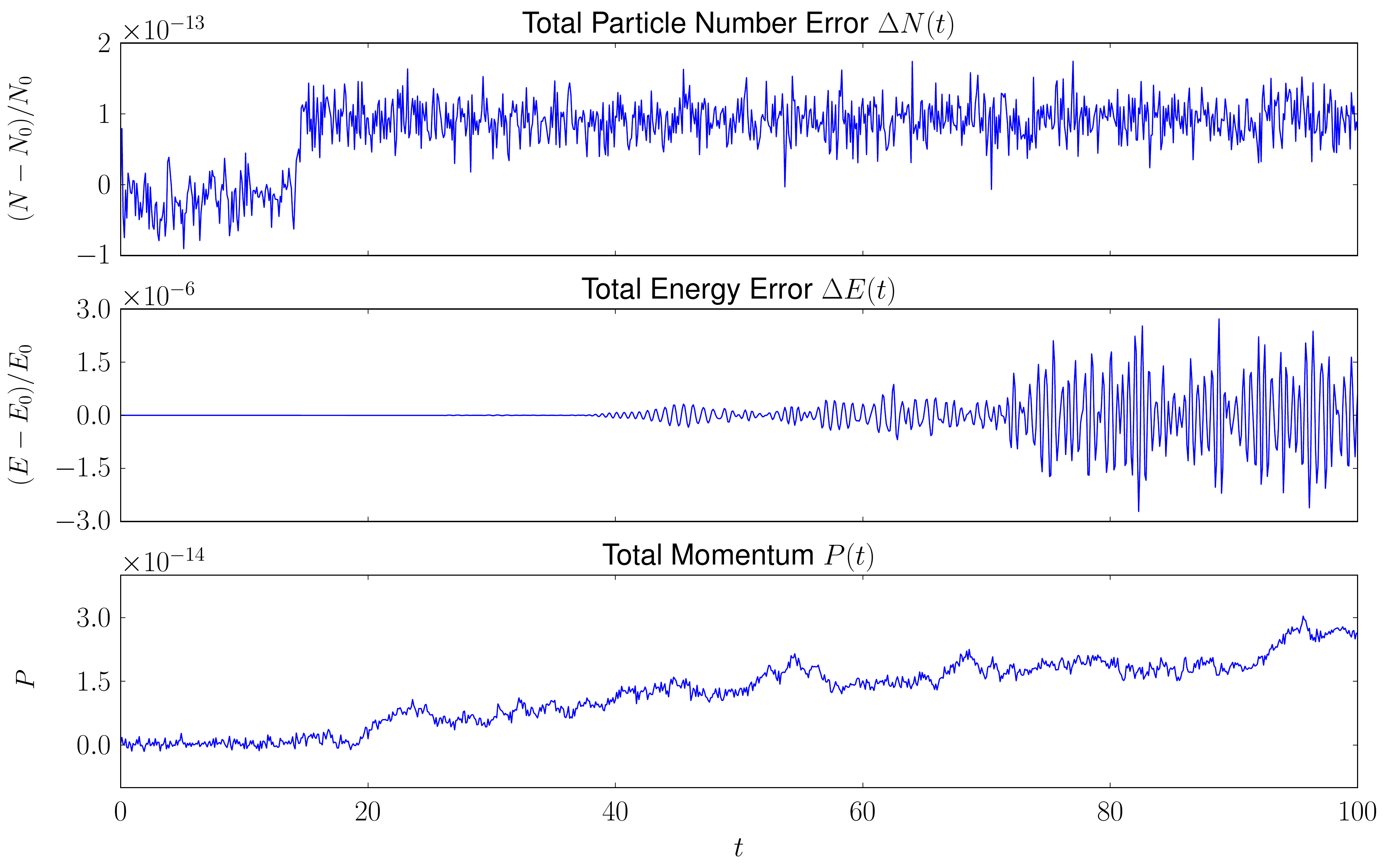}
\caption{Twostream instability without collisions. Total particle number and linear momentum are well preserved but energy conservation is violated due to subgrid modes.}
\label{fig:vlasov_twostream_nu0_NEP}
\end{figure}

\begin{figure}
\centering
\includegraphics[width=.95\textwidth]{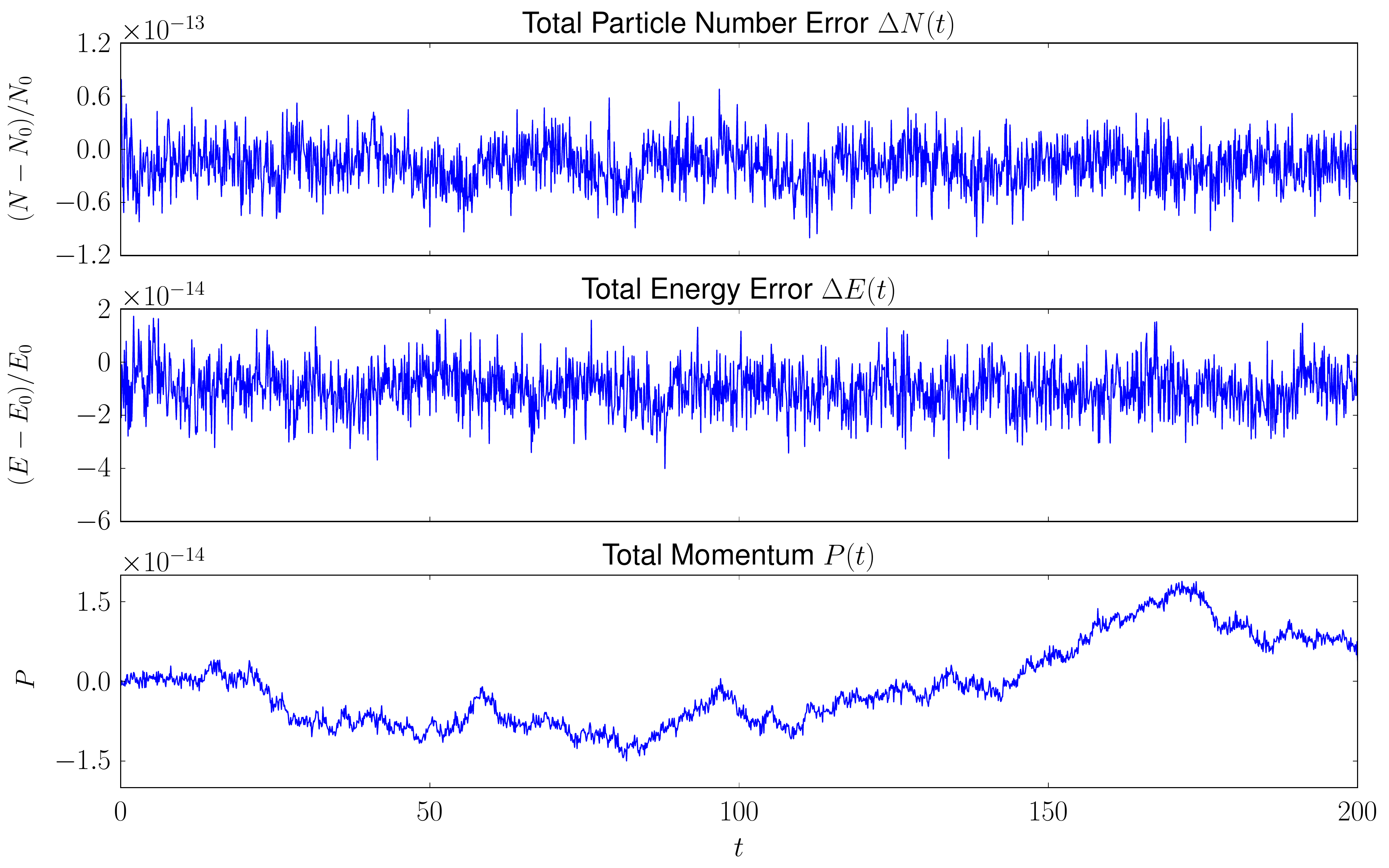}
\caption{Twostream instability with collision frequency $\nu = 4 \times 10^{-4}$. Collisions retain exact energy conservation in addition to exact preservation of the total particle number and linear momentum.}
\label{fig:vlasov_twostream_nu4E-4_NEP}
\end{figure}

In simulations without collisions ($\nu = 0$), we find a good conservation of the total particle number and the total linear momentum, but conservation of the total energy is violated (figure \ref{fig:vlasov_twostream_nu0_NEP}).
Employing the collision operator, optimal energy conservation can be restored (figure \ref{fig:vlasov_twostream_nu4E-4_NEP}). The collision frequency necessary to retain the correct energy throughout the whole simulation is $\nu = 4 \times 10^{-4}$, as it was estimated in the linear Landau damping simulations.

\begin{figure}
\centering
\subfloat{
\includegraphics[width=.48\textwidth]{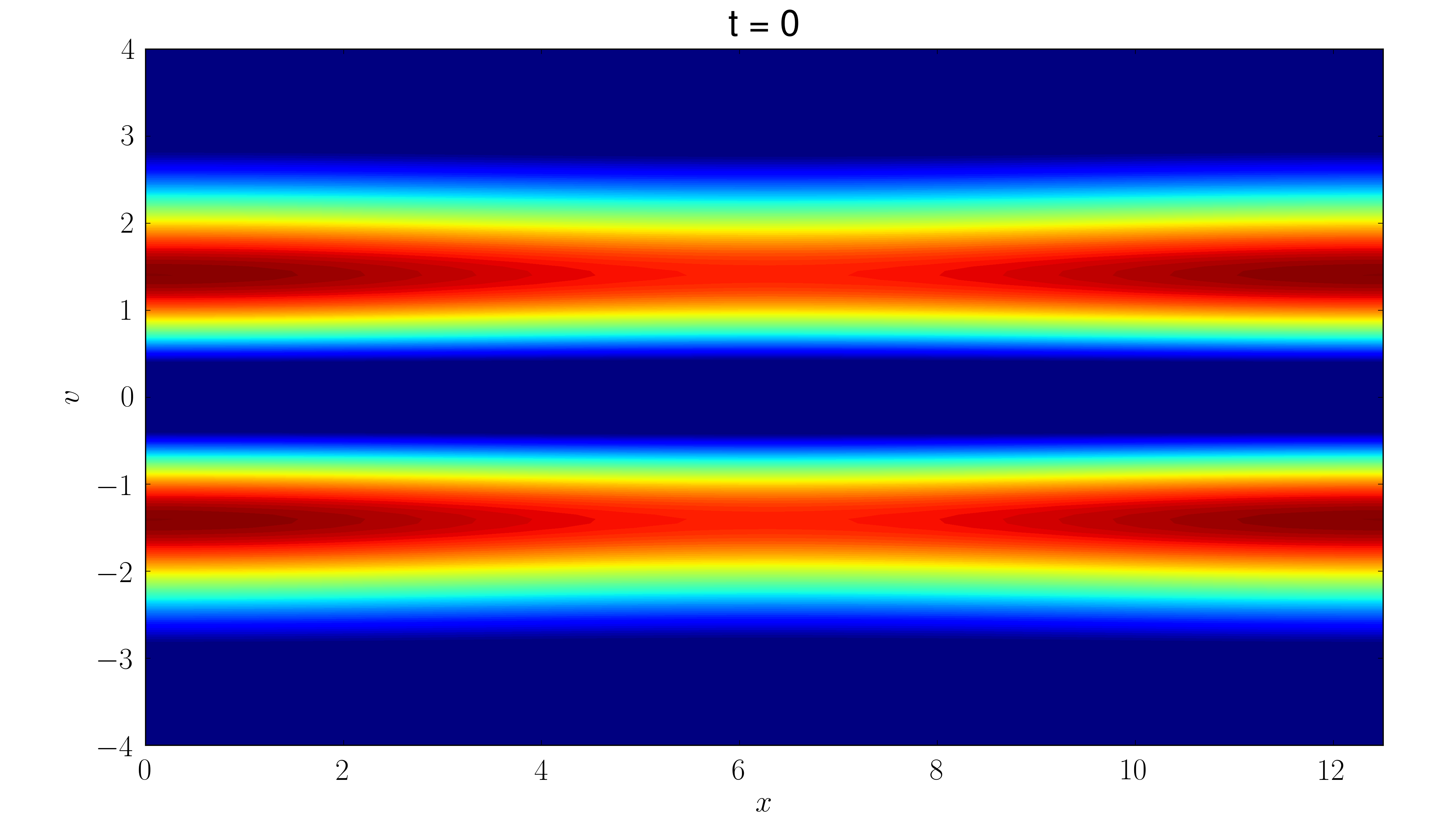}
}
\subfloat{
\includegraphics[width=.48\textwidth]{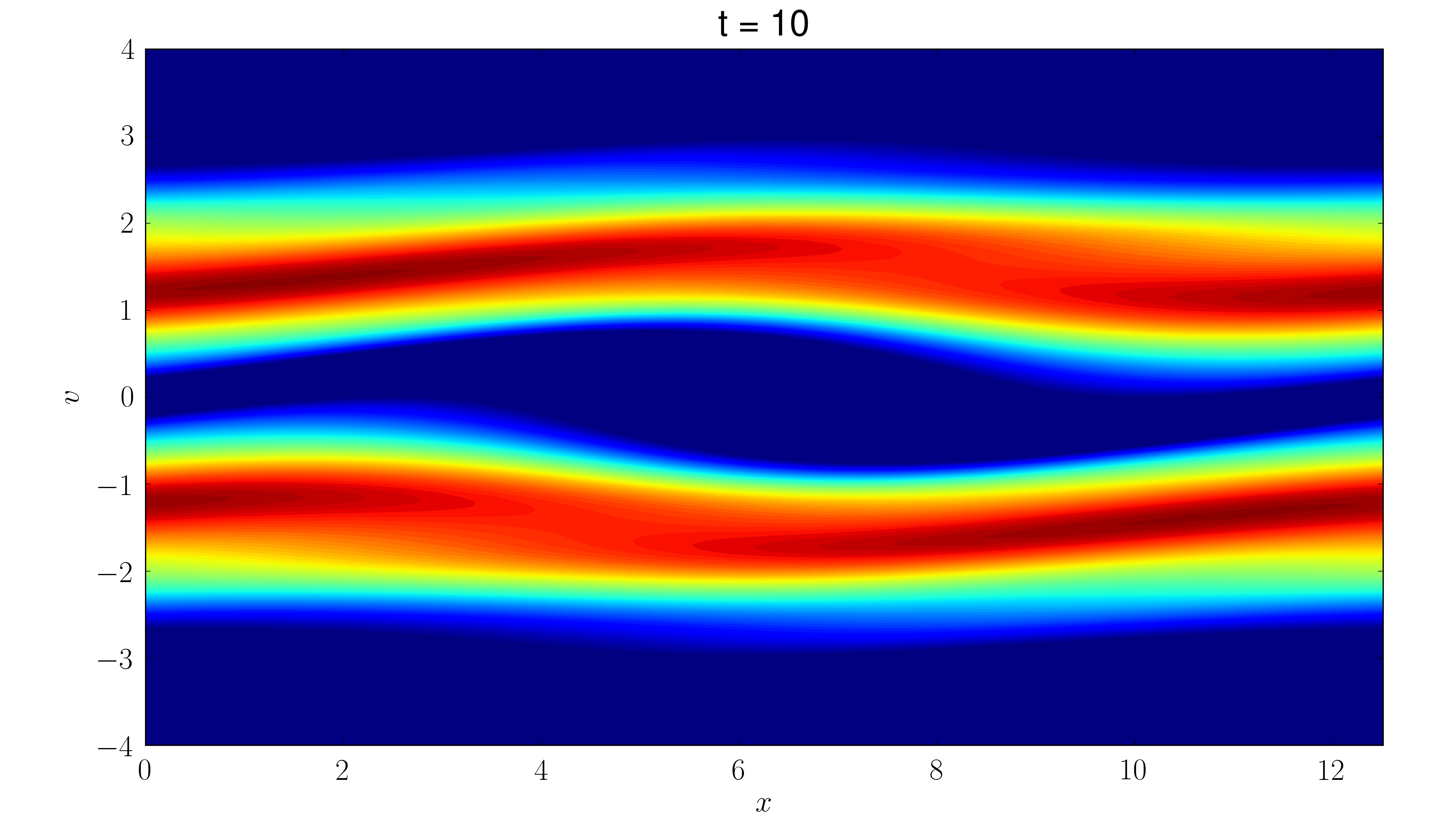}
}

\subfloat{
\includegraphics[width=.48\textwidth]{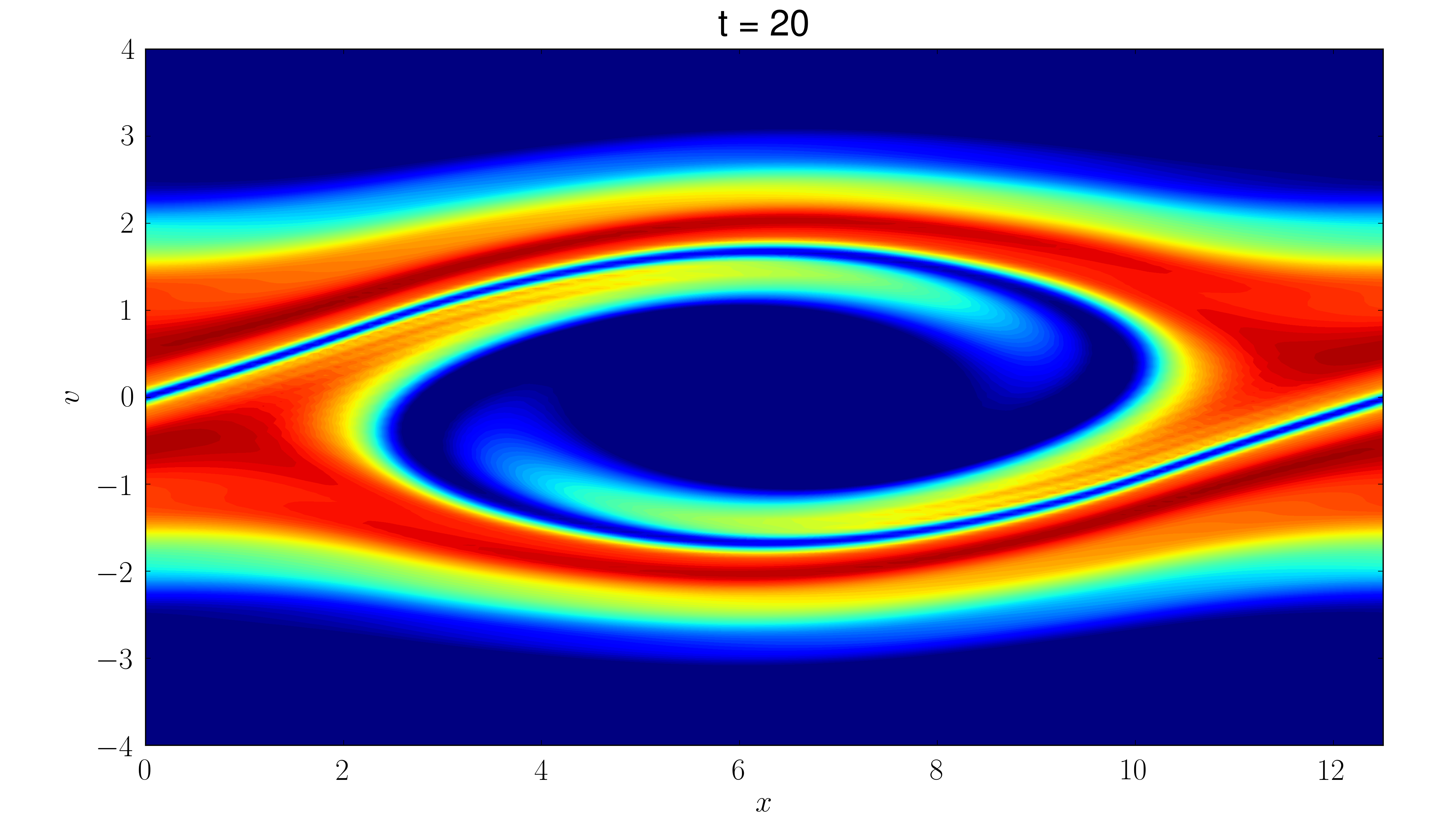}
}
\subfloat{
\includegraphics[width=.48\textwidth]{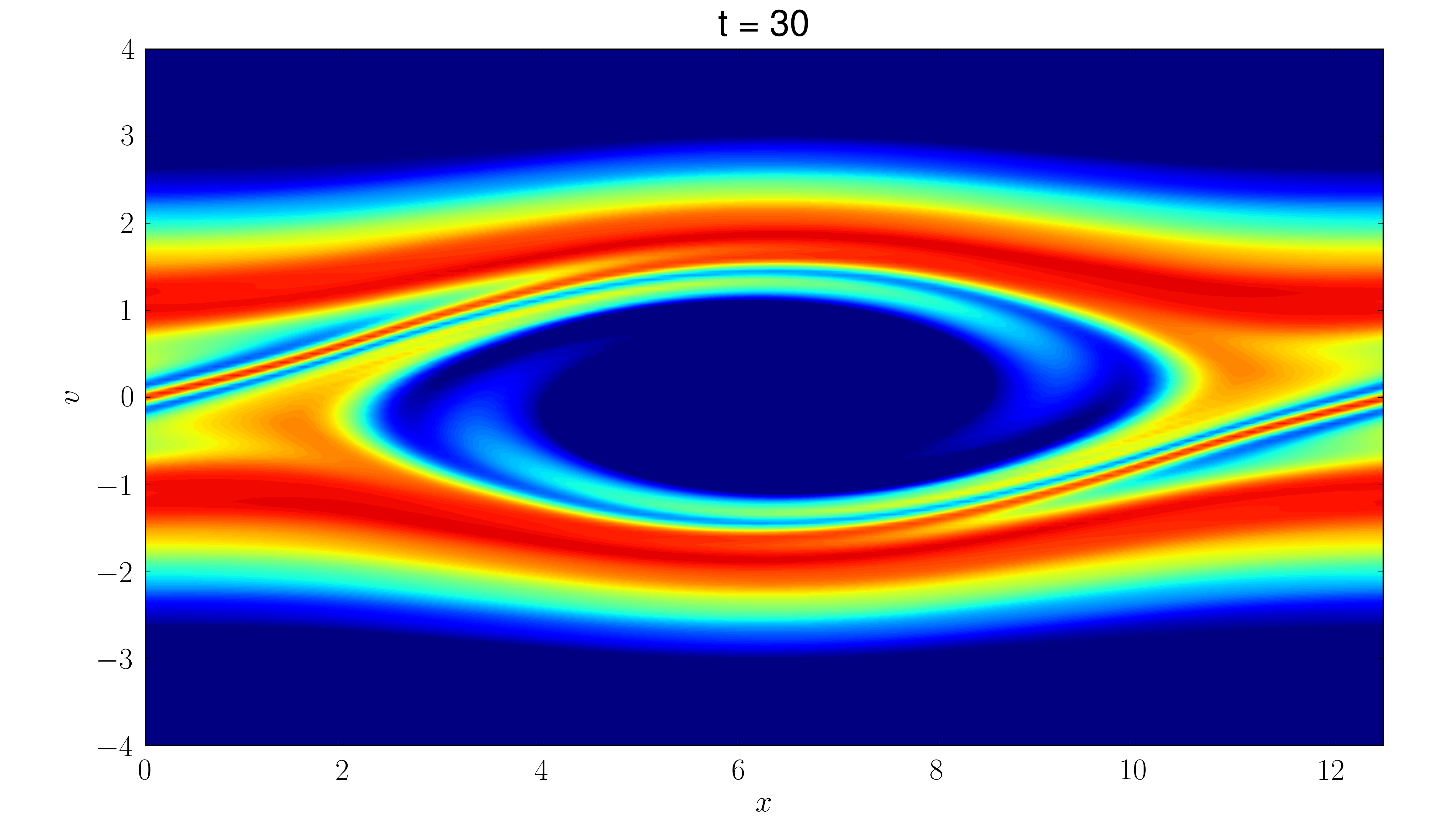}
}

\subfloat{
\includegraphics[width=.48\textwidth]{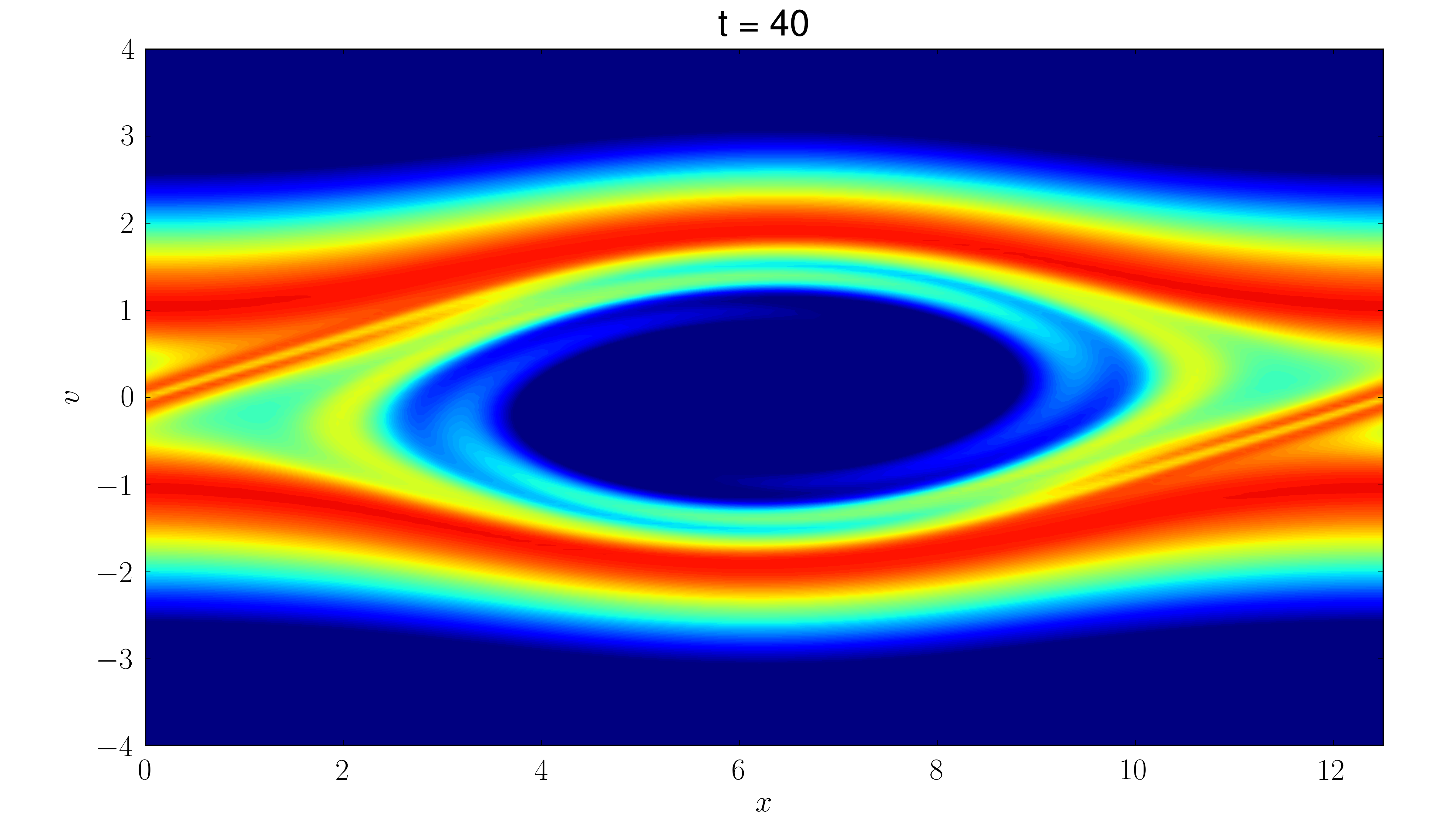}
}
\subfloat{
\includegraphics[width=.48\textwidth]{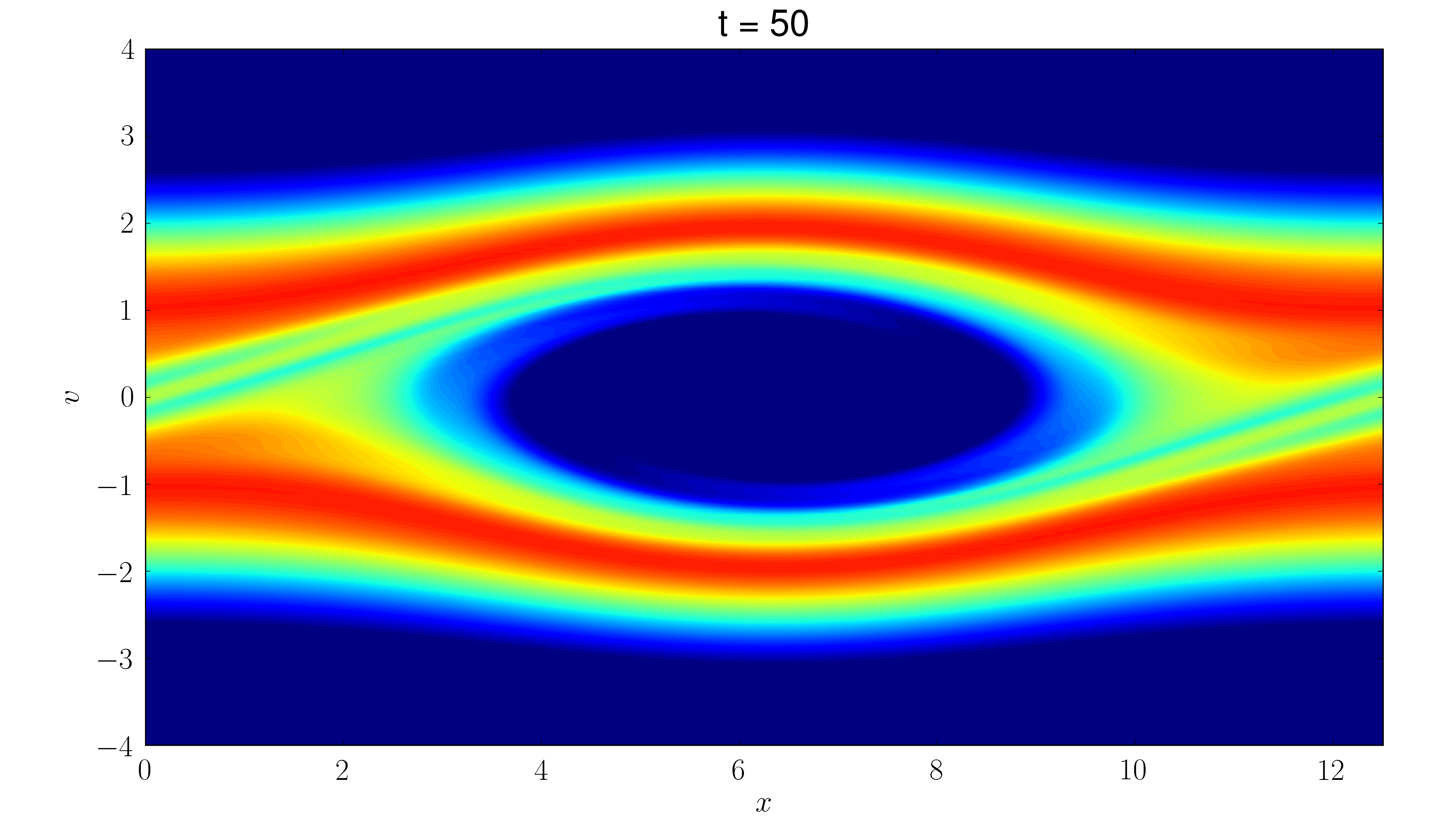}
}

\subfloat{
\includegraphics[width=.48\textwidth]{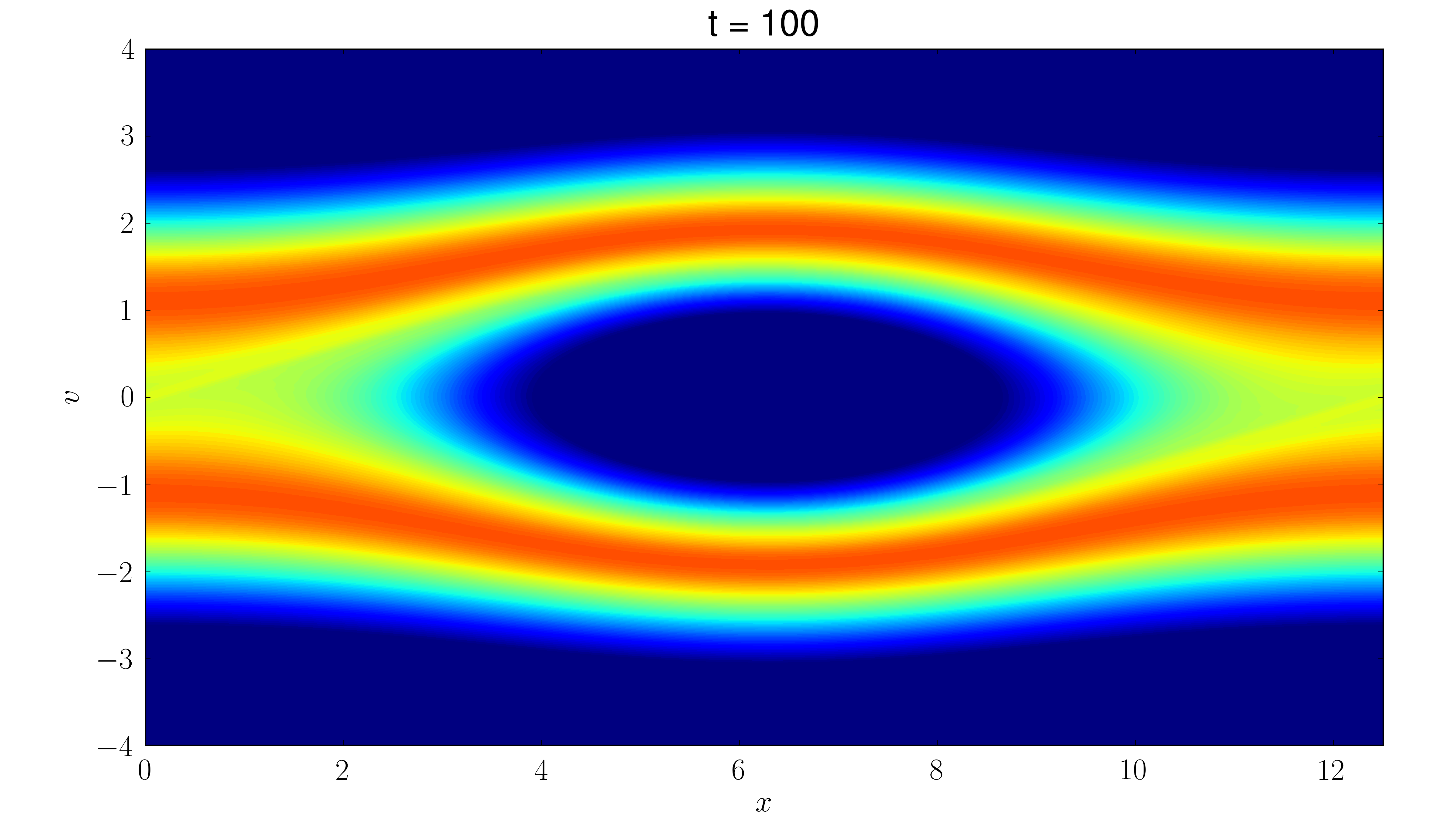}
}
\subfloat{
\includegraphics[width=.48\textwidth]{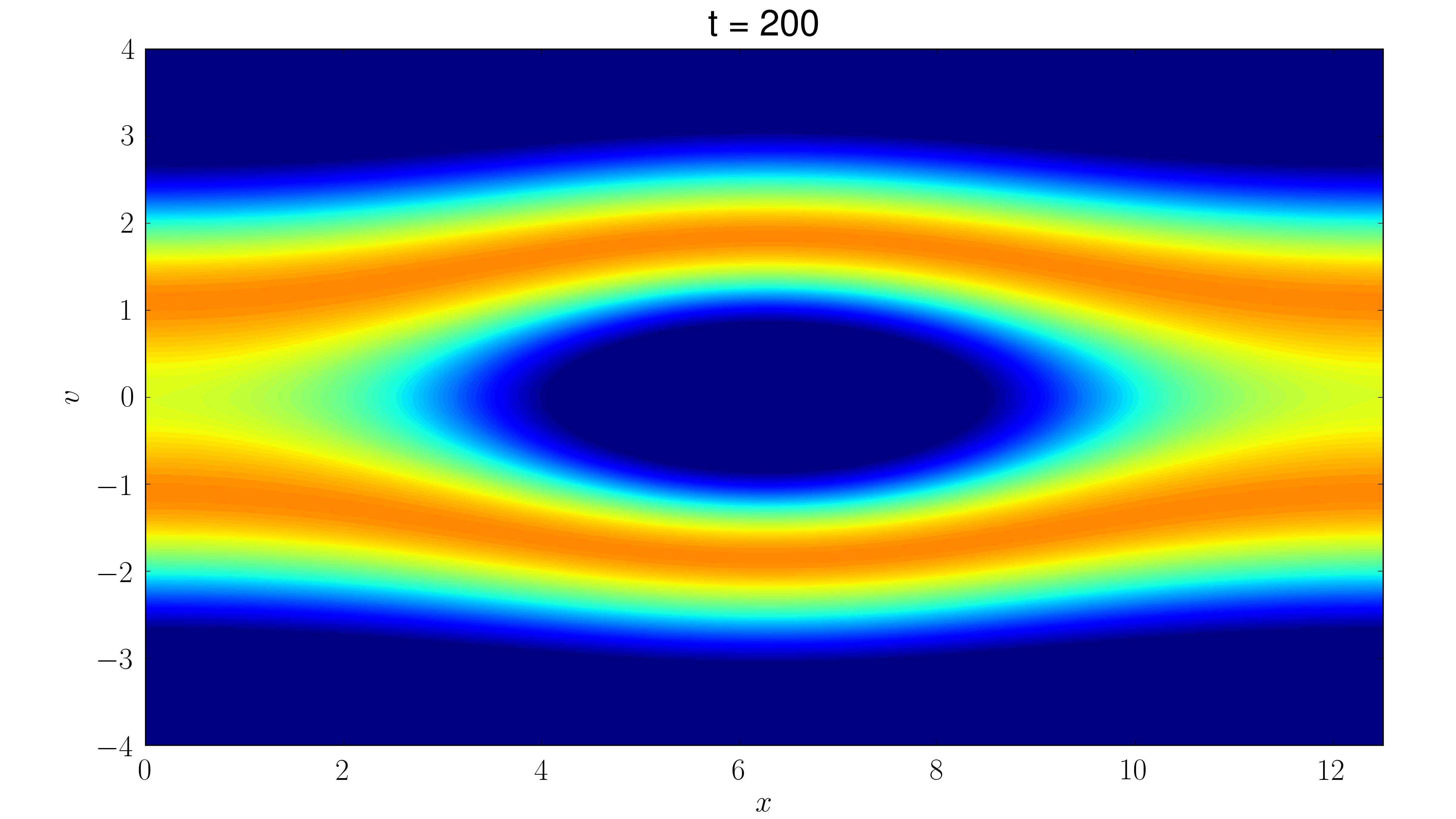}
}
\caption{Twostream instability with collision frequency $\nu = 4 \times 10^{-4}$. Contours of the distribution function in phasespace. Contours are linear and constant.}
\label{fig:vlasov_twostream_nu4E-4_F}
\end{figure}

The distribution function, figure \ref{fig:vlasov_twostream_nu4E-4_F}, exhibits the correct qualitative behaviour.
After an initial growth of the instability, particles become trapped and a hole in phasespace forms.
Between $t = 100$ and $t = 200$, the distribution function takes an almost steady state.

\subsection{Jeans Instability}

Finally, we are considering a test case from gravitational dynamics, the Jeans instability \cite{BinneyTremaine:2008, ChengGamba:2012}.
The only difference compared with plasma dynamics is that the gravitational field is always attractive. This results in a change of sign in the Poisson equation.
The distribution function is initialised as
\begin{align}
f &= f_{M} \, \big( 1 + A \, \cos(k x) \big) ,
\end{align}

with $A = 0.01$ and $k = 0.8$. The spatial domain is $[0, 2.5 \, \pi]$, $v_{\text{max}} = 10$, $n_{x} = 201$, $n_{v} = 401$, $h_{t} = 0.1$.
For $k < 1$, the distribution function is unstable and collapses towards the centre of the simulation domain. For $k > 1$ initial perturbations are damped.

\begin{figure}
\centering
\includegraphics[width=.95\textwidth]{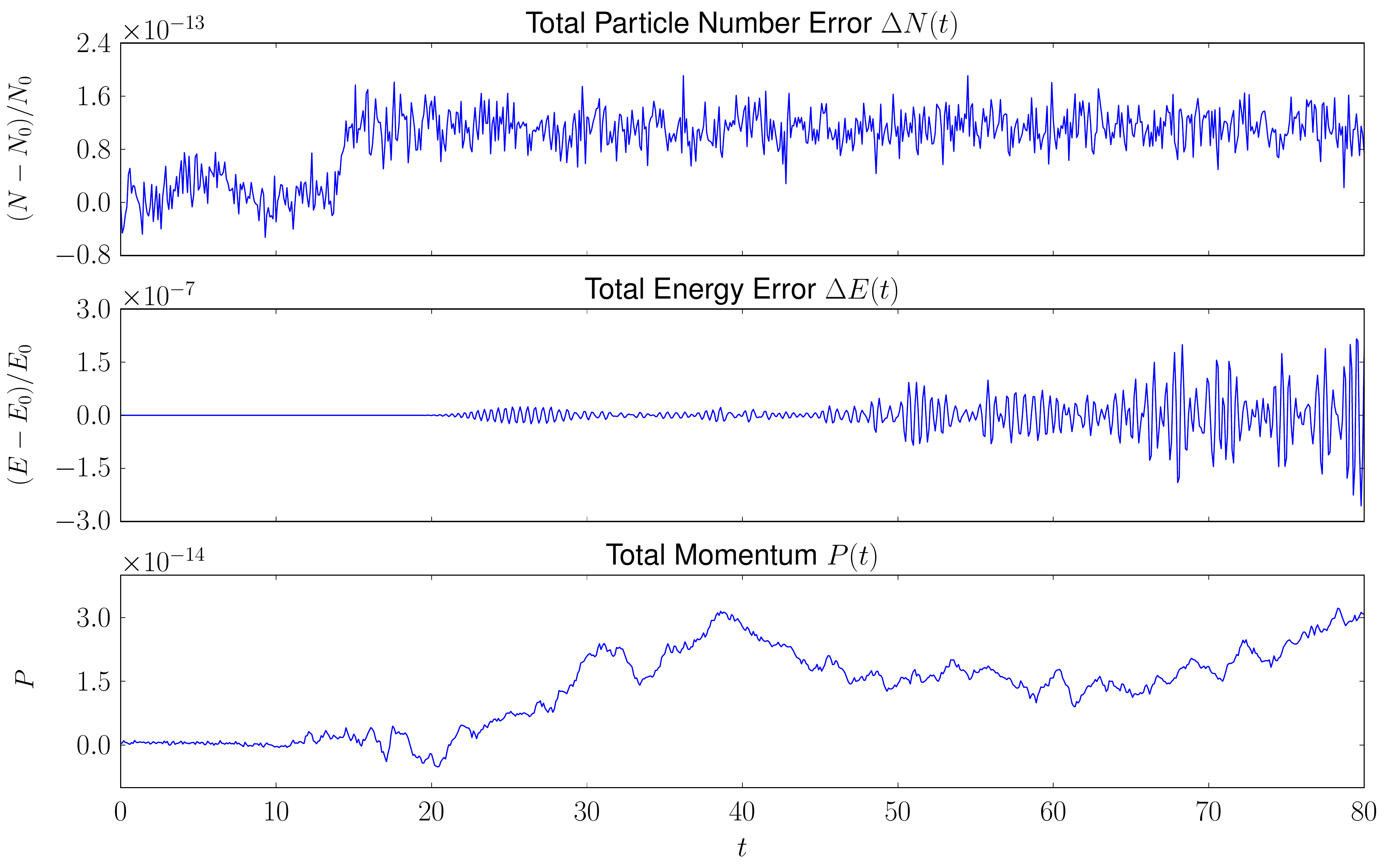}
\caption{Jeans instability without collisions. Total particle number and linear momentum are well preserved, but energy conservation is violated due to subgrid modes.}
\label{fig:vlasov_jeans_weak_nu0_NEP}
\end{figure}

\begin{figure}
\centering
\includegraphics[width=.95\textwidth]{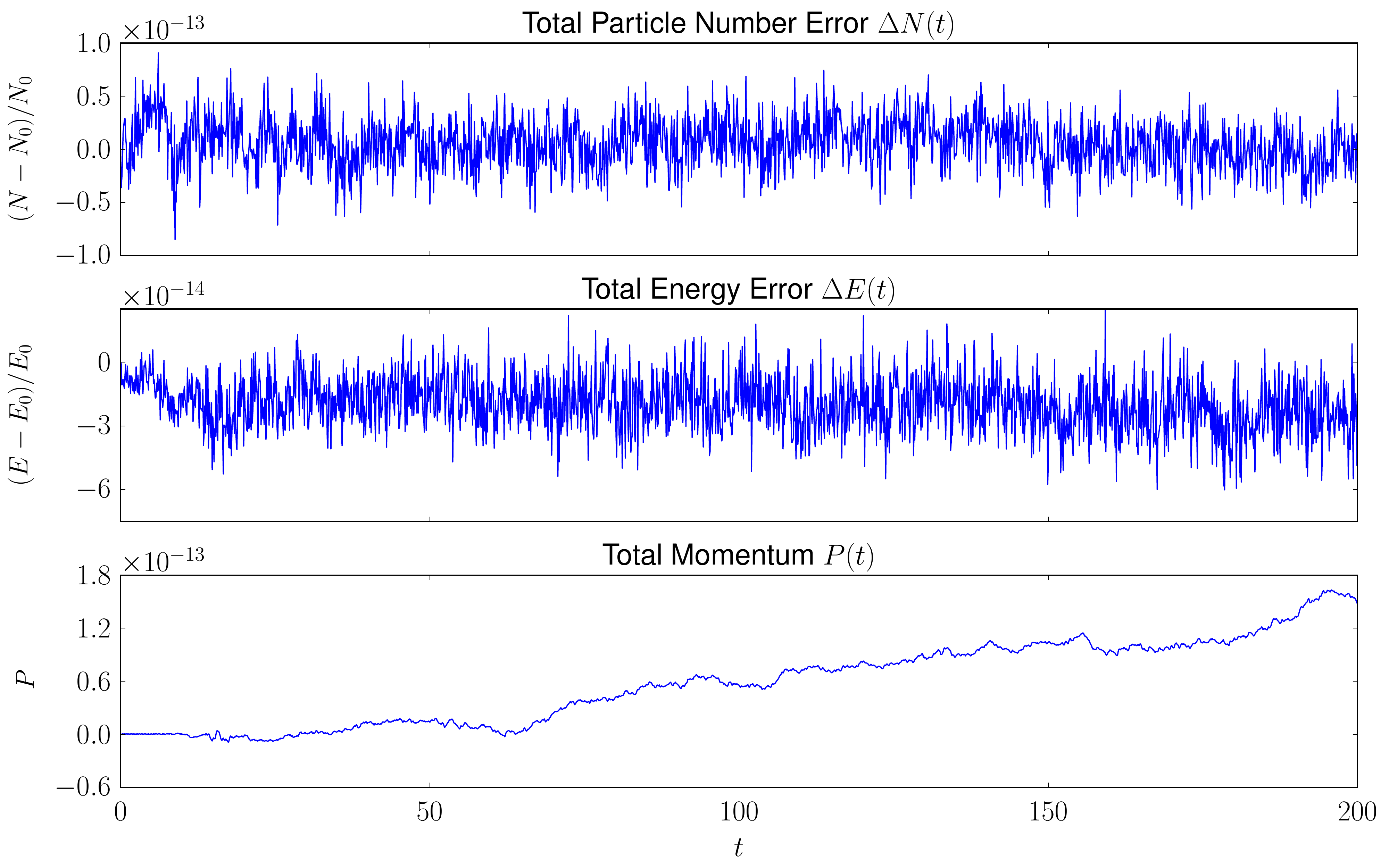}
\caption{Jeans instability with collision frequency $\nu = 4 \times 10^{-4}$. Collisions retain exact energy conservation in addition to exact preservation of the total particle number and linear momentum.}
\label{fig:vlasov_jeans_weak_nu4E-4_NEP}
\end{figure}

The conservation properties are very similar as in the case of the twostream instability.
Without collisions, the total particle number and the total linear momentum are well conserved but not the total energy (figure \ref{fig:vlasov_jeans_weak_nu0_NEP}).
Adding collisions, with the same frequency of $\nu = 4 \times 10^{-4}$ as before ($h_{v}$ is still the same), energy conservation is retained (figure \ref{fig:vlasov_jeans_weak_nu4E-4_NEP}).

\begin{figure}
\centering
\subfloat{
\includegraphics[width=.48\textwidth]{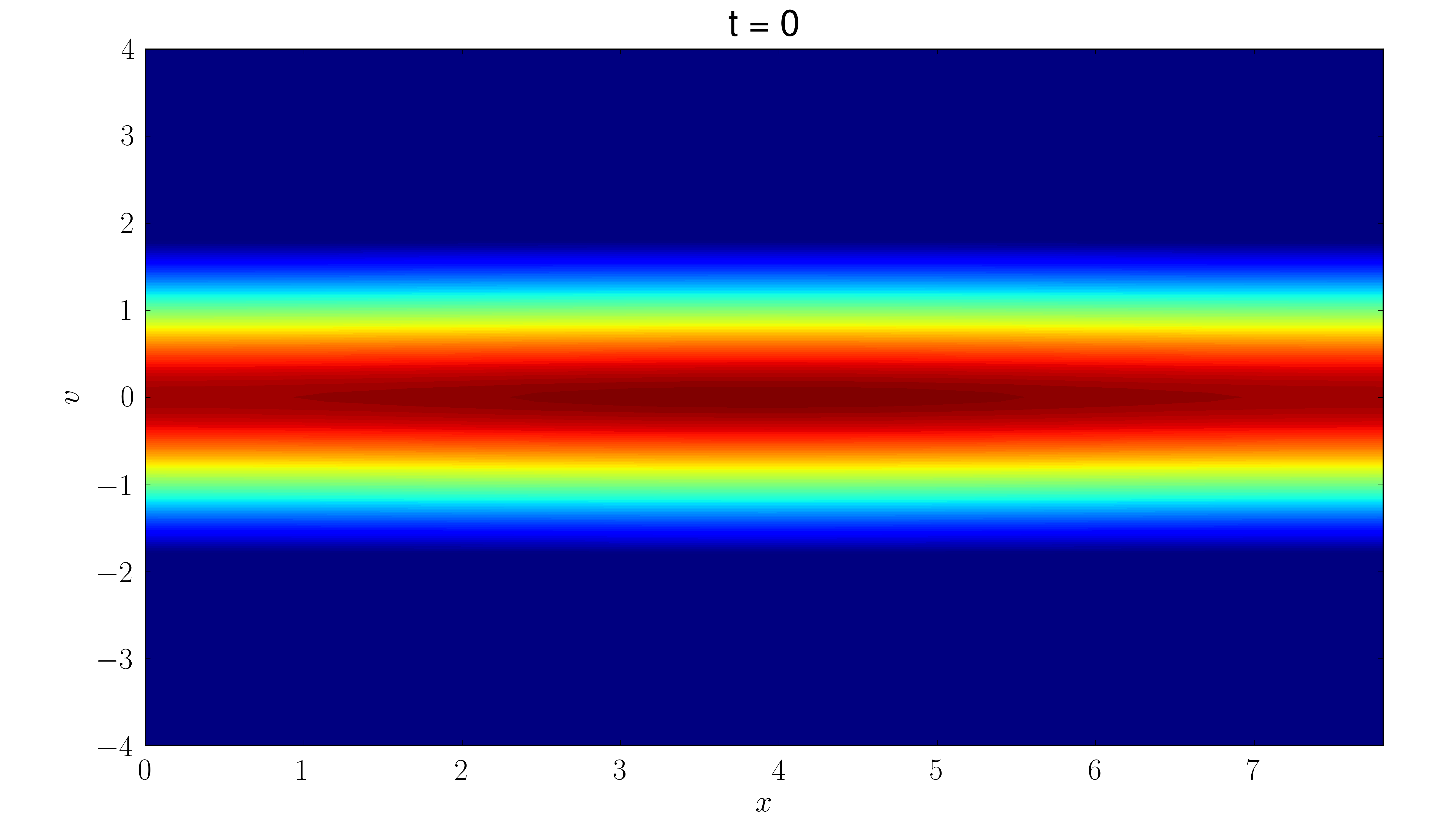}
}
\subfloat{
\includegraphics[width=.48\textwidth]{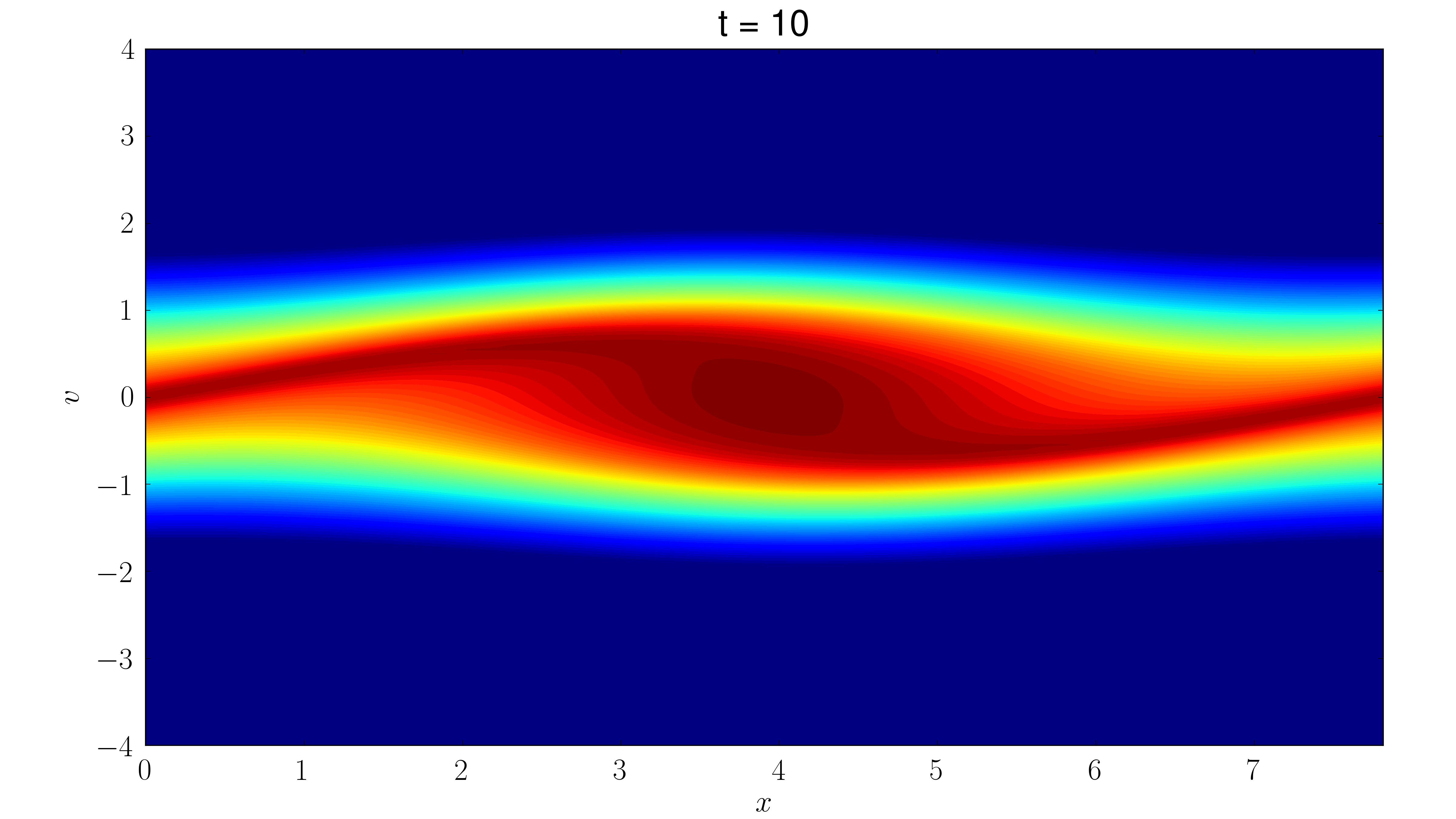}
}

\subfloat{
\includegraphics[width=.48\textwidth]{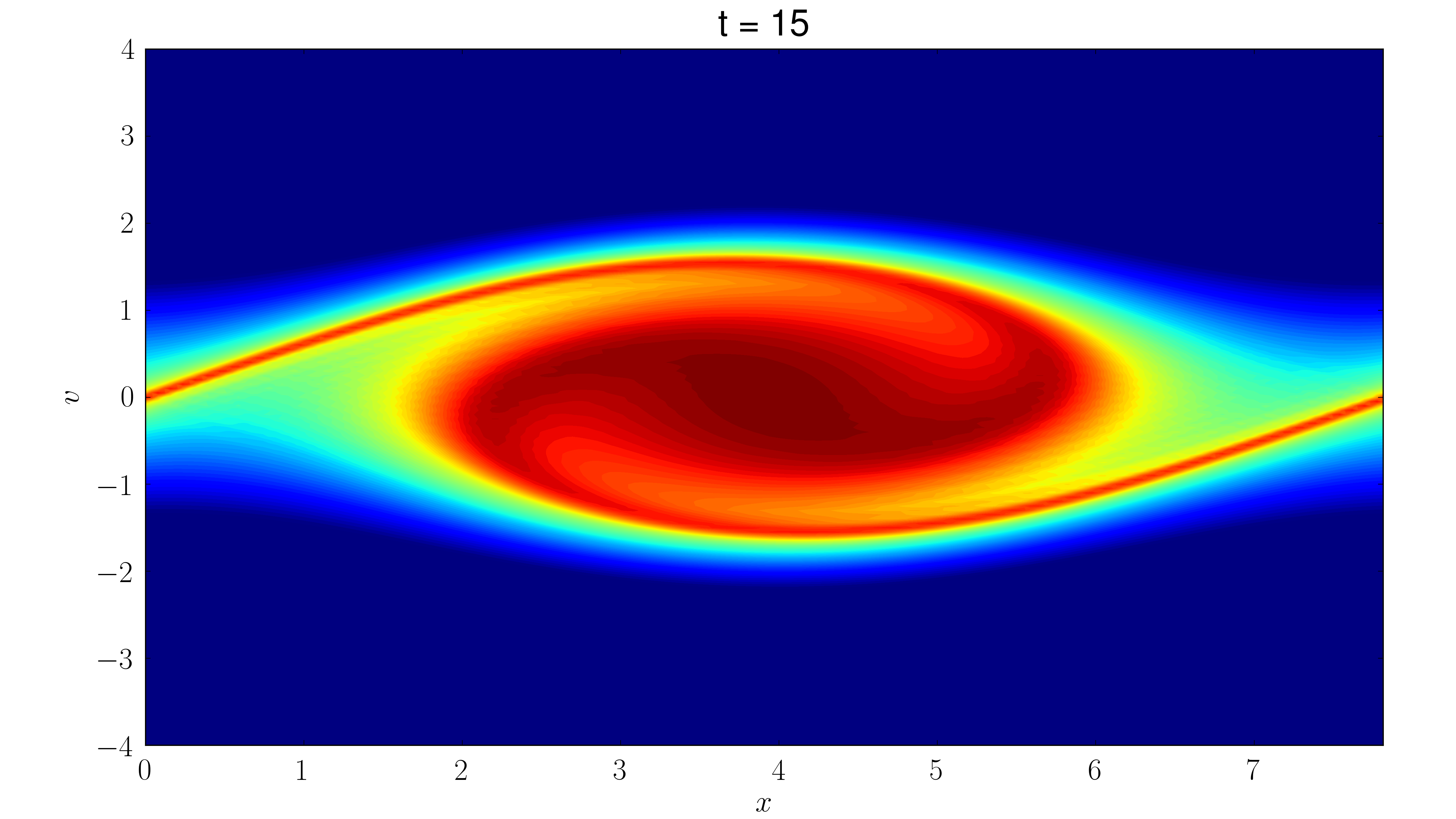}
}
\subfloat{
\includegraphics[width=.48\textwidth]{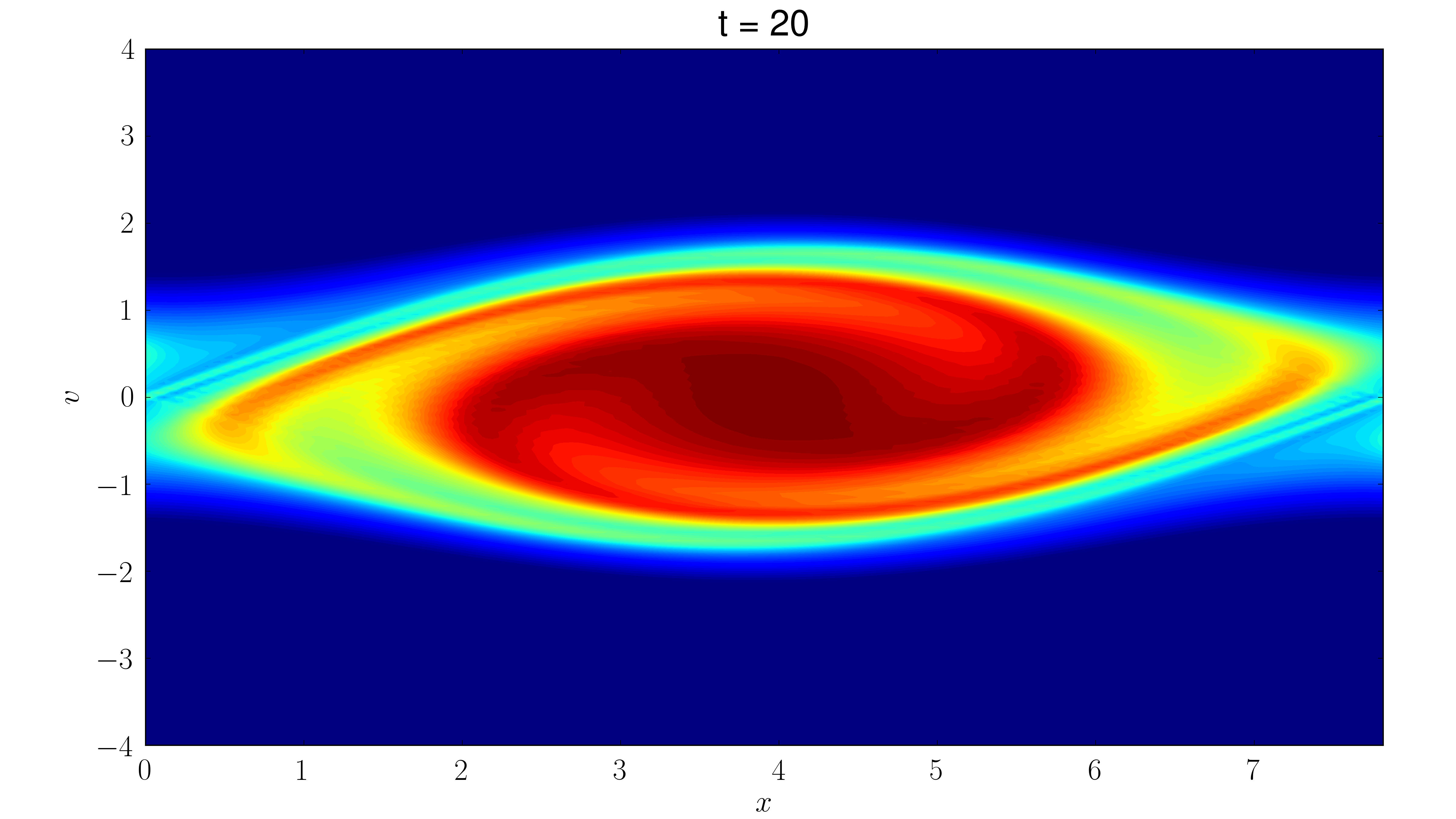}
}

\subfloat{
\includegraphics[width=.48\textwidth]{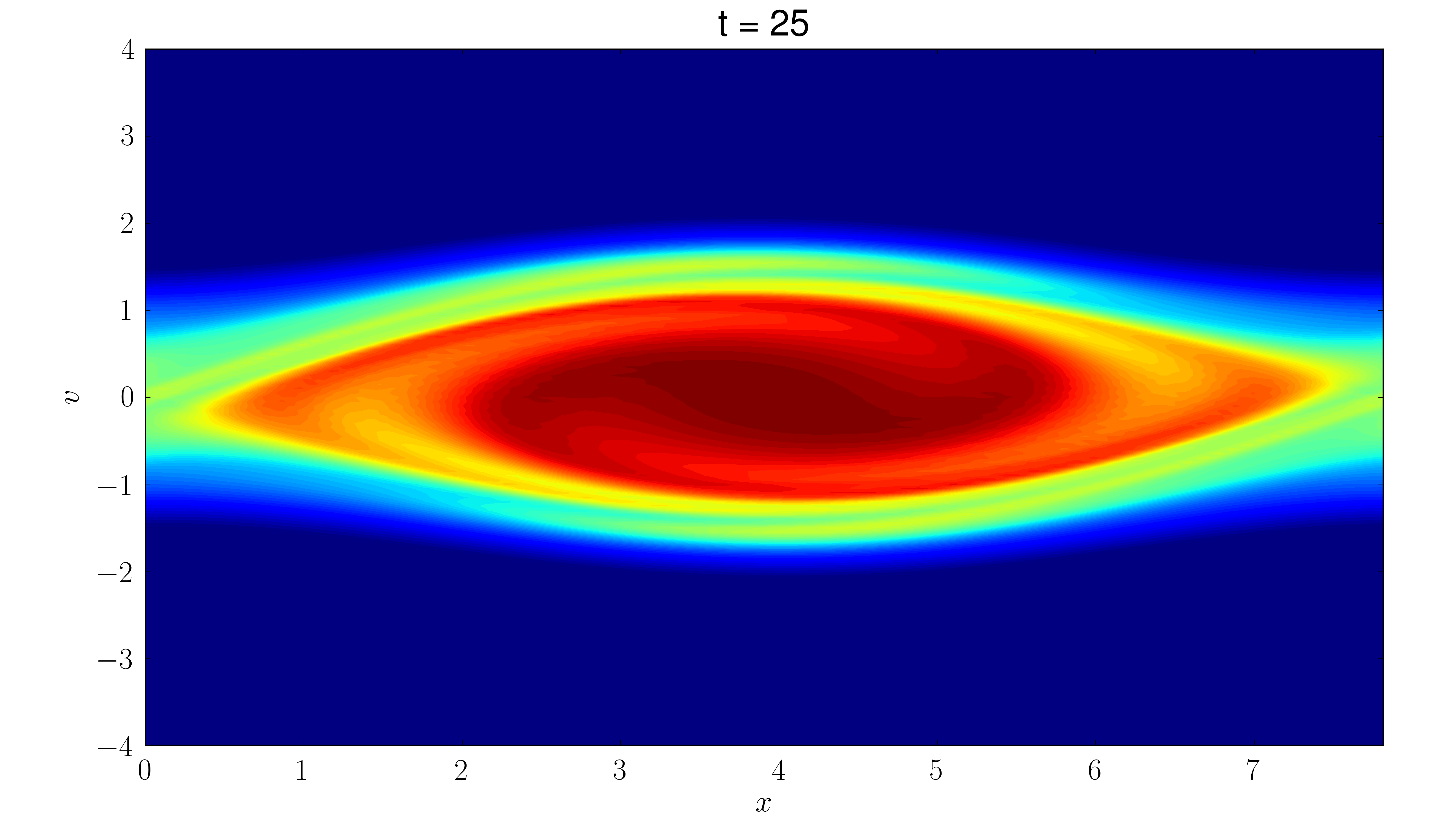}
}
\subfloat{
\includegraphics[width=.48\textwidth]{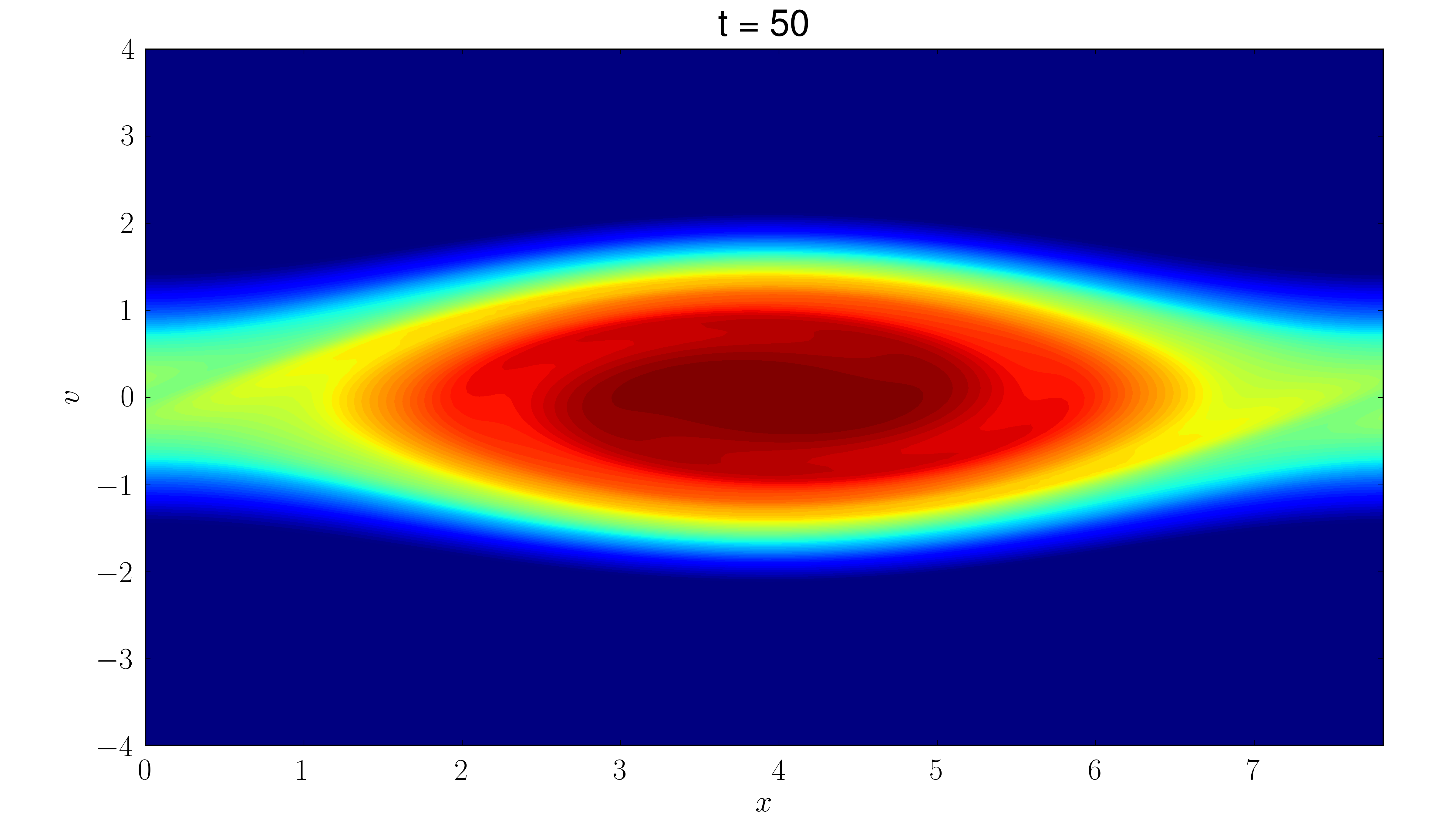}
}

\subfloat{
\includegraphics[width=.48\textwidth]{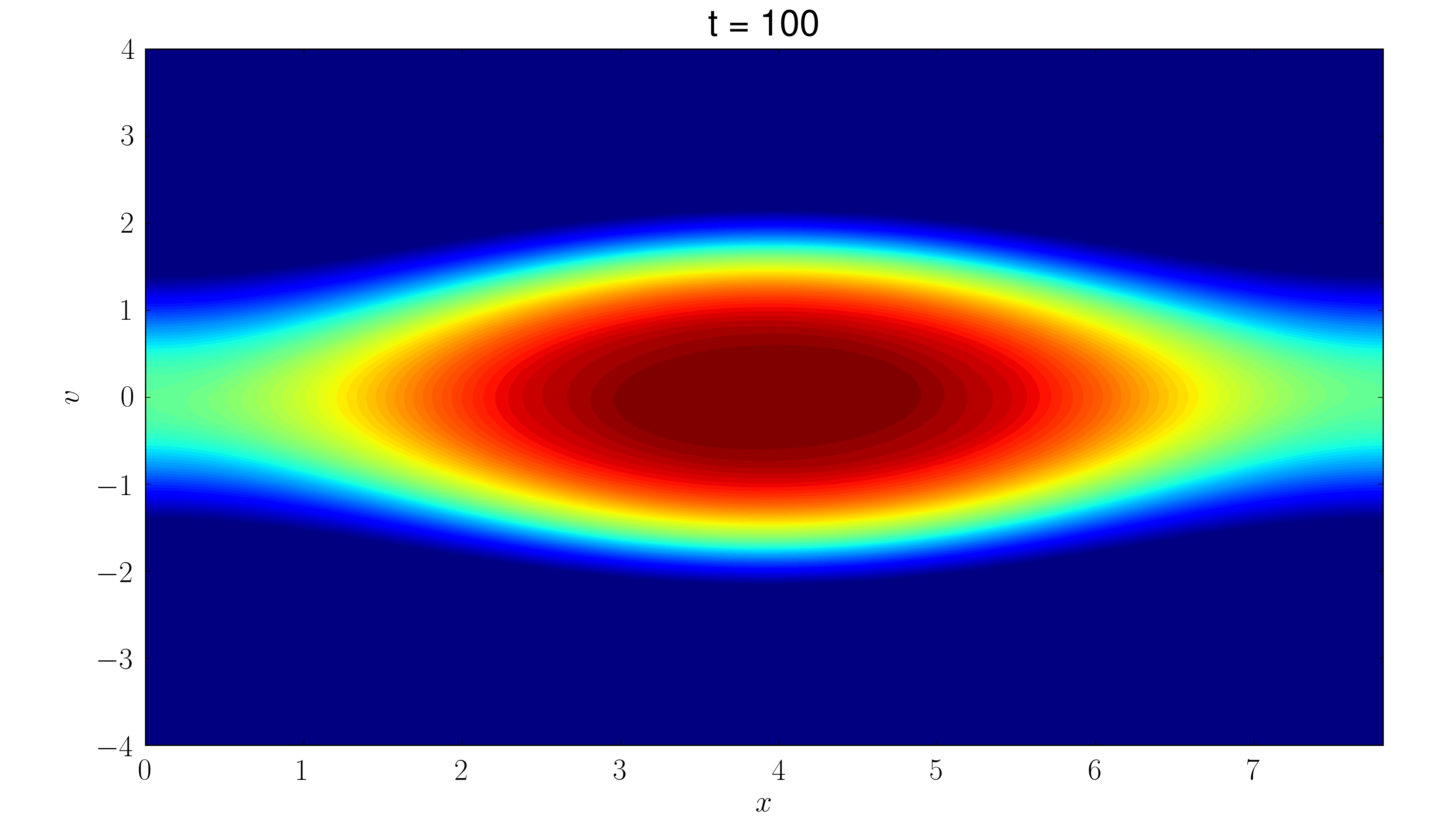}
}
\subfloat{
\includegraphics[width=.48\textwidth]{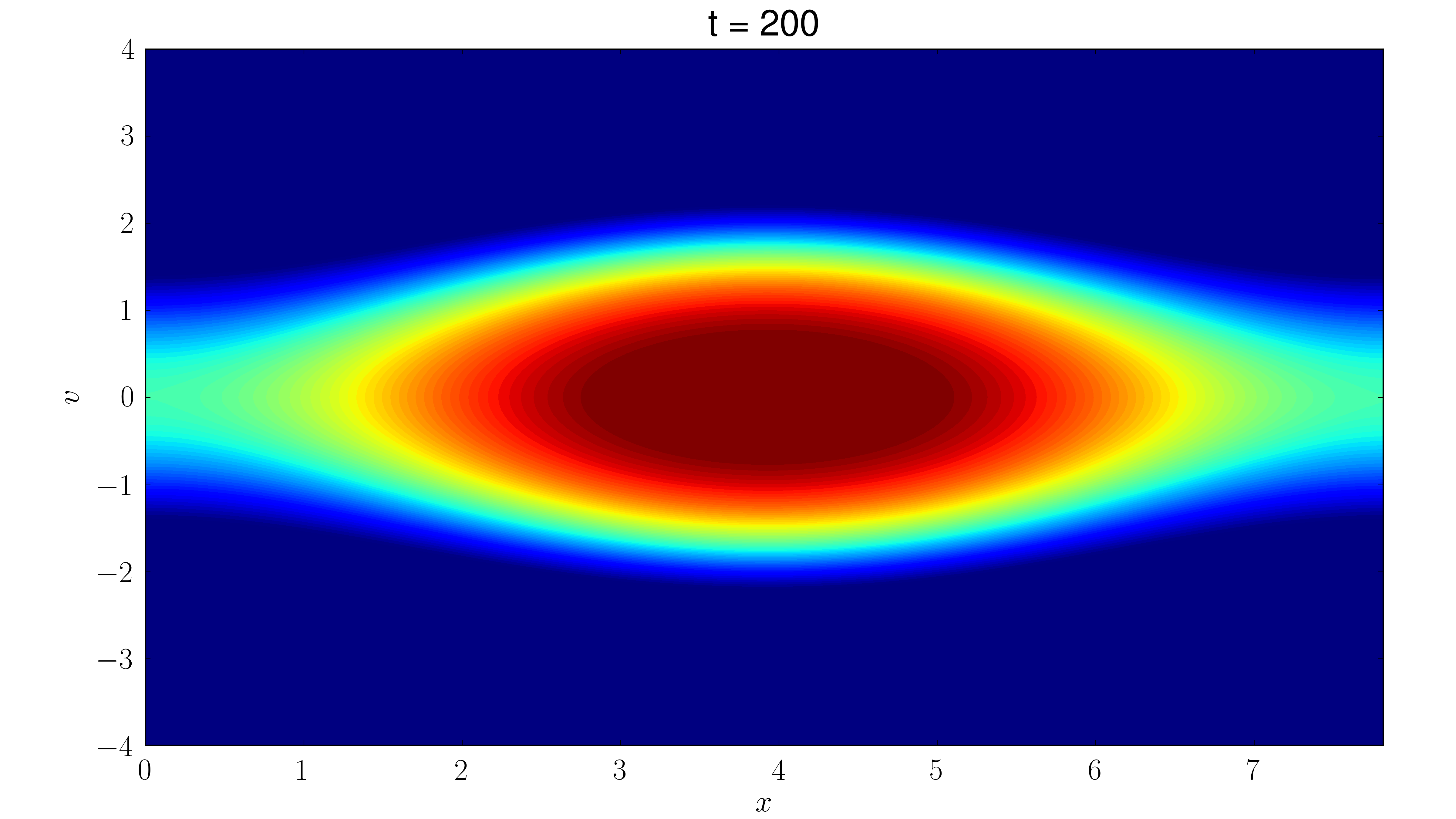}
}
\caption{Jeans instability with collision frequency $\nu = 4 \times 10^{-4}$. Contours of the distribution function in phasespace. Contours are linear and constant.}
\label{fig:vlasov_jeans_weak_nu4E-4_F}
\end{figure}

The qualitative behaviour of the distribution function (figure \ref{fig:vlasov_jeans_weak_nu4E-4_F}) meets the expectations. As we choose $k < 1$, the initial perturbation is unstable and develops a swirl about the centre of the simulation domain, corresponding to a gravitational collapse. Between $t = 50$ and $t = 100$ the distribution functions reaches a steady state and barely changes until $t = 200$.

\newpage

\clearpage

\clearpage

\begin{figure}
\centering
\includegraphics[width=\textwidth]{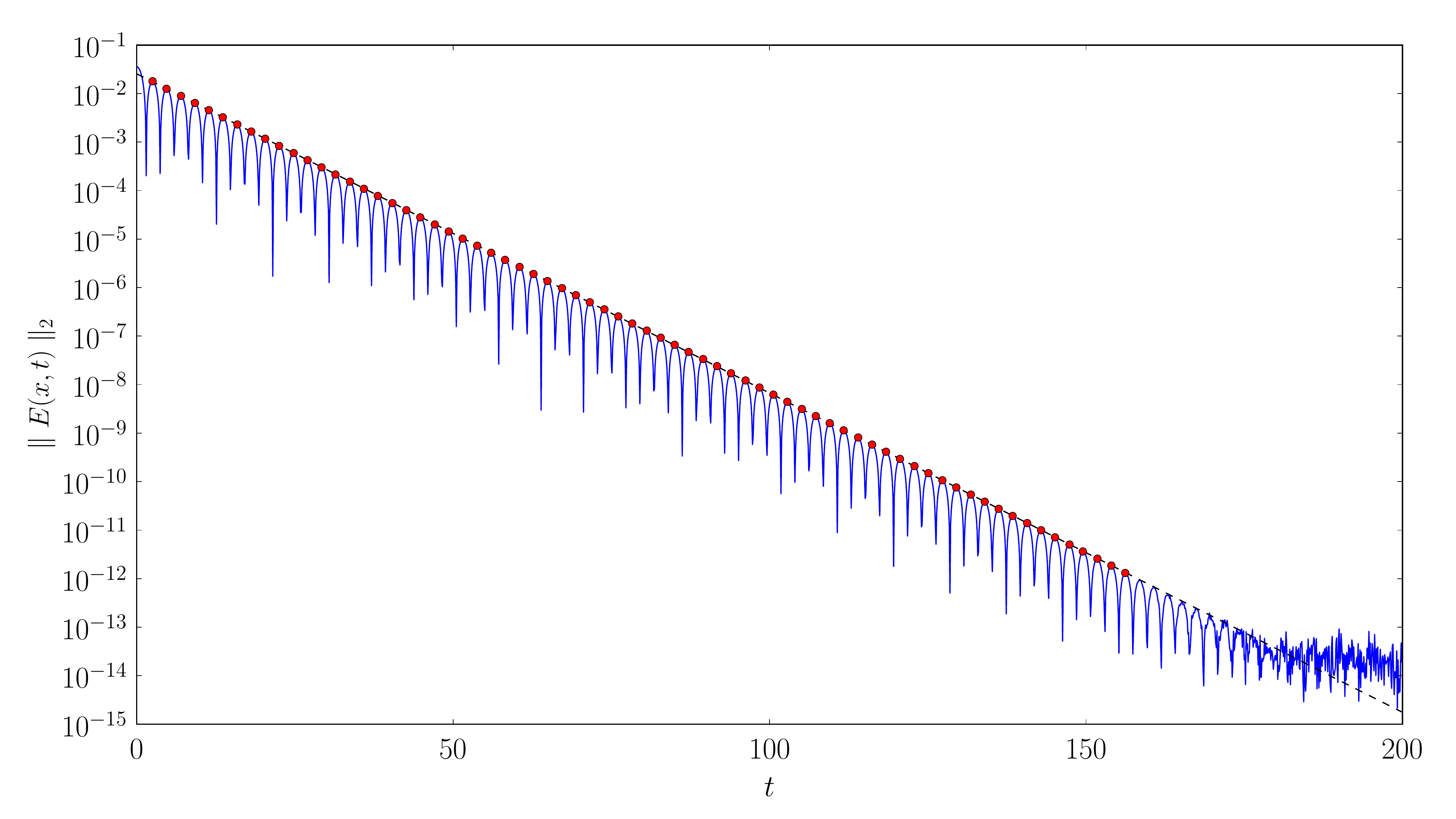}
\caption{Linear Landau damping with collision frequency $\nu = 4 \times 10^{-4}$. Evolution of the electrostatic energy. Both, the correct qualitative behaviour and the correct damping rate, $\gamma = - 0.152$, are obtained for $\nu = 4 \times 10^{-4}$.}
\label{fig:vlasov_landau_linear_nu4E-4_potential}
\end{figure}

\clearpage

\clearpage

\clearpage

\clearpage

\clearpage

\clearpage

\clearpage

\clearpage

\clearpage

\clearpage

\chapter{Magnetohydrodynamics}\label{ch:mhd}

Magnetohydrodynamics (MHD) describes the dynamics of electrically conducting fluids like plasmas or liquid metals.
It is one of the most widely applied theories in laboratory as well as astrophysical plasmas physics \cite{Schnack:2009,GoedbloedPoedts:2004,Biskamp:2003,Davidson:2001,Freidberg:1987}, used to describe macroscopic phenomena 	like equilibrium states in tokamaks or stellarators, large scale turbulence, or dynamos that generate magnetic fields of stars and planets.
The structure of the equations is very similar to hydrodynamics, albeit in MHD the fluid equations are coupled with Maxwell's equations, thereby allowing for an even richer variety of phenomena.

\section{Incompressible ideal MHD}

The equations of magnetohydrodynamics result from the combination of the Navier-Stokes equation for an incompressible fluid,
\begin{align}\label{eq:mhd_eqs_navier_stokes}
\dfrac{\partial V}{\partial t} + ( V \cdot \nabla ) V = - \nabla p + \mu \, \nabla^{2} V + F ,
\end{align}

with Maxwell's equations of electrodynamics.
As usual for an incompressible flow, the mass densities are taken constant. $V$ is the fluid velocity, satisfying $\nabla \cdot V = 0$, $p$ is the pressure, $F$ is a force term, and $\mu$ is the viscosity of the fluid. The right-hand side is normalised to the density $\rho$.
The fluid carries an electric current density $J$ and is immersed into a magnetic field $B$, which produces the force $F = J \times B$.
Throughout this chapter, we use natural units for electromagnetic quantities, i.e., $\mu_{0} = \eps_{0} = c = 1$.

To obtain an evolution equation for the magnetic field, we combine Faraday's law
\begin{align}\label{eq:mhd_eqs_faraday}
\dfrac{\partial B}{\partial t} = - \nabla \times E
\end{align}

with Ohm's law for a resistive plasma
\begin{align}\label{eq:mhd_eqs_ohm}
E + V \times B = \eta J ,
\end{align}

$\eta$ being the plasma resistivity and $E$ the electric field,
and Ampere's law
\begin{align}\label{eq:mhd_eqs_ampere}
J = \nabla \times B
\end{align}

with the displacement current neglected, 
to get
\begin{align}\label{eq:mhd_eqs_induction1}
\dfrac{\partial B}{\partial t}
&= \nabla \times ( V \times B ) - \eta \, \nabla \times \nabla \times B .
\end{align}

With the identities
\begin{align}
\nabla \times \nabla \times B &= \nabla \, ( \nabla \cdot B ) - \nabla^{2} B , \\
\nabla \times ( A \times B ) &= A \, ( \nabla \cdot B ) - B \, ( \nabla \cdot A ) + ( B \cdot \nabla ) \, A - ( A \cdot \nabla ) \, B ,
\end{align}

the induction equation (\ref{eq:mhd_eqs_induction1}) becomes
\begin{align}\label{eq:mhd_eqs_induction2}
\dfrac{\partial B}{\partial t}
&= V \, ( \nabla \cdot B) - B \, ( \nabla \cdot V ) + ( B \cdot \nabla ) V - ( V \cdot \nabla ) B
- \eta \, \nabla ( \nabla \cdot B ) + \eta \, \nabla^{2} B .
\end{align}

Both $V$ and $B$ are divergence-free, so the induction equation simplifies to
\begin{align}\label{eq:mhd_eqs_induction3}
\dfrac{\partial B}{\partial t} + ( V \cdot \nabla ) B &= ( B \cdot \nabla ) V + \eta \, \nabla^{2} B .
\end{align}

The force term $F$ in the Navier-Stokes equation is rewritten as
\begin{align}\label{eq:mhd_eqs_lorentz}
F = J \times B = ( \nabla \times B ) \times B = ( B \cdot \nabla ) B - \dfrac{1}{2}
\nabla ( B \cdot B ) .
\end{align}

Again, we used Ampere's law and the identity
\begin{align}
B \times ( \nabla \times B) &= \tfrac{1}{2} \, \nabla ( B \cdot B ) - ( B \cdot \nabla ) B .
\end{align}

Summing up, we obtain the following system of equations

\rimpeq[Incompressible MHD Equations]{
\label{eq:mhd_eqs_V}
\dfrac{\partial V}{\partial t} + ( V \cdot \nabla ) V &= ( B \cdot \nabla ) B + \mu \, \nabla^{2} V - \nabla P , &
\nabla \cdot V &= 0 , \\
\label{eq:mhd_eqs_B}
\dfrac{\partial B}{\partial t} + ( V \cdot \nabla ) B &= ( B \cdot \nabla ) V  + \eta \, \nabla^{2} B , &
\nabla \cdot B &= 0 ,
}

with the generalised pressure $P$ being the sum of the kinetic gas pressure and the magnetic pressure 
\begin{align}\label{eq:mhd_eqs_total_pressure}
P \equiv p + \dfrac{1}{2} \, B^{2} .
\end{align}

The first equation (\ref{eq:mhd_eqs_V}) is called the \emph{momentum equation}, the second equation (\ref{eq:mhd_eqs_B}) the \emph{induction equation}.
Both $V$ and $B$ are divergence-free, $V$ as we are considering an incompressible fluid, and $B$ as there are no magnetic monopoles.
But while $\nabla \cdot B = 0$ is implied by the induction equation (provided that the initial magnetic field $B(t=0)$ is divergence-free), $\nabla \cdot V = 0$ is a dynamical constraint, for which the pressure $P$ is a Lagrange multiplier.

The left-hand sides of (\ref{eq:mhd_eqs_V},\ref{eq:mhd_eqs_B}) represent the advective derivatives of the velocity field $V$ and the magnetic field $B$.
The force term $( B \cdot \nabla ) B$ is the directional derivative of $B$ in direction of $B$. It describes the magnetic tension force, a restoring force that straightens magnetic field lines. This force is perpendicular to $B$ and inversely proportional to the radius of the field line curvature so that the fluid is accelerated towards the local centre of the curvature. As the magnetic field is advected by the fluid, the field lines are dragged with the fluid and thus straightened.

The term $\mu \, \nabla^{2} V$ describes diffusion of the fluid due to viscosity and results from the divergence of the anisotropic part of the stress tensor. As we are dealing with an incompressible fluid, this term only amounts to shear stress.

The pressure gradient $\nabla p$ arises from the isotropic part of the stress tensor which describes normal forces. The effect of this term is that fluid flows from regions of high pressure to regions of low pressure. The magnetic pressure has the same effect, fluid flows from regions of high magnetic pressure to regions of low magnetic pressure. Due to the advection of the magnetic field with the fluid velocity this then leads to the field lines being pushed apart (imagine a bundle of field lines driven apart).

When $\eta = 0$, (\ref{eq:mhd_eqs_B}) states that the magnetic field is advected with the fluid flow, which implies the conservation of the magnetic flux through a surface moving with the fluid \cite{ArnoldKhesin:1998}.
In a resistive plasma, $\eta \, \nabla^{2} B$ describes diffusive effects, for which the magnetic field are not just dragged along with the field, but are free to change their topology.

In ideal MHD, viscosity and resistivity are neglected, thus $\mu = \eta = 0$.
As an effect, the topology of the magnetic field lines is conserved. They are not allowed to open up and reconnect. A property that we would like to maintain on the discrete level.
Two important conserved quantities of ideal MHD in two dimensions \cite{ArnoldKhesin:1998} are the total energy
\begin{align}
E = \dfrac{1}{2} \int \Big[ \norm{V}^{2} + \norm{B}^{2} \Big] \, dx \, dy ,
\end{align}

and cross helicity
\begin{align}
H = \int V \cdot B \, dx \, dy .
\end{align}

Conservation of both quantities are desirable in numerical simulations.

\subsection{Lie Derivative Formulation}\label{sec:mhd_lie_derivative}

To elucidate the link with work related to ours \cite{Gawlik:2011}, we rewrite the MHD equations with Lie derivatives, thereby also emphasising the advective character of the equations a bit further.
Write the ideal MHD equations in component form, use covariant components in the momentum equation, and add and subtract $V^{j} \partial_{i} V_{j}$ and $B^{j} \partial_{i} B_{j}$
\begin{subequations}\label{eq:mhd_eqs2}
\begin{align}
& \partial_{t} V_{i} + V^{j} \partial_{j} V_{i} + V^{j} \partial_{i} V_{j} = B^{j} \partial_{j} B_{i} + B^{j} \partial_{i} B_{j} + V^{j} \partial_{i} V_{j} - B^{j} \partial_{i} B_{j} - \partial_{i} P , \\
& \partial_{t} B^{i} + V^{j} \partial_{j} B^{i} - B^{j} \partial_{j} V^{i} = 0 ,
\end{align}
\end{subequations}

where $V_{i}$ and $V^{i}$ are co- and contravariant components, respectively, and analogously for $B_{i}$ and $B^{i}$.
The second and third term on the left-hand side of the momentum equation are the Lie derivative of a 1-form, $V^{\flat} = V_{i} \, dx^{i}$, along its corresponding vector field $V$. The first two terms on the right-hand side are the Lie derivative of another 1-form, $B^{\flat} = B_{i} \, dx^{i}$, also along its corresponding vector field $B$.
In the induction equation, the second and third term are the Lie derivative of a vector field, $B$, along $V$ (recall section \ref{sec:geometry_lie_derivative}).
The reformulated ideal MHD equations read

\rimpeq[Lie Derivative Formulation of Ideal Magnetohydrodynamics]{
\dfrac{\partial V^{\flat}}{\partial t} &+ \lie_{V} V^{\flat} = \lie_{B} B^{\flat} - \ext \Big[ P + \dfrac{1}{2} \, \norm{B}^{2} - \dfrac{1}{2} \, \norm{V}^{2} \Big] , \\
\dfrac{\partial B}{\partial t} &+ \lie_{V} B = 0 .
}

Now the interpretation of the time evolution of the fields becomes even more apparent. The evolution of the fluid velocity is determined by three different mechanisms.
The Lie derivative of the velocity along itself describes the change of velocity along the fluid flow. The Lie derivative of the magnetic field along itself describes how the magnetic field changes along field lines.
The resulting force, as already discussed, pushes the velocity field towards the local centre of the magnetic field-line curvature, and thus, as the magnetic field is advected by the fluid flow, balances variations in the magnetic field.
The gradient of the pressure states that fluid flows from regions of high (kinetic and magnetic) pressure to regions of low (kinetic and magnetic) pressure.
The Lie derivative in the induction equation just states that the magnetic field is advected along the fluid flow.

But what is the meaning of the additional terms? Essentially they just remove physics that we added with the Lie derivative but that was not present in the original equations.
The Lie derivative describes all actions that happen along the given vector field, namely a (rigid) translation, a (rigid) rotation, and a deformation, but the original term, e.g. $(u \cdot \nabla) u$, describes only a translation.
Let us try to better understand this by having a look at the external derivative of the velocity term (in component form)
\begin{align}
\dfrac{1}{2} \, \partial_{i} \norm{V}^{2}
= V^{j} \, \partial_{i} V_{j}
= \dfrac{1}{2} \, \big( V^{j} \, \partial_{i} V_{j} + V^{j} \, \partial_{j} V_{i} \big)
+ \dfrac{1}{2} \, \big( V^{j} \, \partial_{i} V_{j} - V^{j} \, \partial_{j} V_{i} \big)
= \mrm{T}_{ij} V^{j} + \mrm{S}_{ij} V^{j} .
\end{align}

Here we split $\partial_{i} V_{j}$ in a symmetric part $\mrm{T}_{ij}$ and an antisymmetric part $\mrm{S}_{ij}$.
The symmetric tensor $\mrm{T}$ describes the rate of stretching of the vector $V$ along the direction of the eigenvectors or $\mrm{T}$. As we are discussing incompressible fluids only, the divergence of $V$ vanishes and the trace of $\mrm{T}$ is zero. The fluid element gets deformed but its volume stays constant.
Taking the dot-product of a vector with the antisymmetric tensor $\mrm{S}$ describes the rate of rotation of the vector $V$ with angular velocity vector $\tfrac{1}{2} \xi$, with $\xi = \omega = \nabla \times V$ being the vorticity.
Analogously, $\tfrac{1}{2} \, \partial_{i} \norm{B}^{2}$ splits into two contributions, describing the rate of stretching and rotation of $B$, where now $\xi = J = \nabla \times B$.

\subsection{Potential Formulation in Two Dimensions}\label{sec:mhd_potential}

Another formulation of magnetohydrodynamics, especially popular in reconnection studies, is the so called potential formulation.
Here, the dynamics is not described in terms of the velocity field and the magnetic field, but in terms of their potentials, the streaming function $\psi$ and the magnetic vector potential $A$.

\rimpeq[Potential Formulation of Ideal Magnetohydrodynamics]{
\dfrac{\partial (\Delta \psi)}{\partial t} + [ \psi , \Delta \psi ] &= [ A , \Delta A ] , \\
\dfrac{\partial A}{\partial t} + [ A , \psi ] &= 0 .
}

Here, $[ \cdot , \cdot ]$ are Poisson brackets with respect to the spatial variables $(x,y)$.
The velocity and the magnetic field are computed as
\begin{align}
V &= \nabla \times \psi , &
B &= \nabla \times A , &
& \text{where} &
\psi &= (0,0,\psi) &
& \text{and} &
A &= (0,0,A) , &
\end{align}

such that in this formulation, the constraints $\nabla \cdot V$ and $\nabla \cdot B$ are automatically fulfilled.
With the help of vorticity and current density
\begin{align}
\omega &= \nabla \times V , &
J &= \nabla \times B , &
& \text{where} &
\omega &= (0,0,\omega) &
& \text{and} &
J &= (0,0,J) , &
\end{align}

we can rewrite the above equations as
\begin{align}
\dfrac{\partial \omega}{\partial t} + [ \omega , \psi ] &= [ J , A ] , &
- \Delta \psi &= \omega , \\
\dfrac{\partial J}{\partial t} + [ \omega , A ] &= [ J , \psi ] + 2 \, [ \nabla \psi , \nabla A ] , &
- \Delta A &= J , &
\end{align}

thereby reducing the highest order of derivatives that appear in the equations from three to two. An important point, since this simplifies the derivation of variational integrators.
Looking at this formulation, one might get the expression that we can directly apply the discretisation from the last chapter to this formulation. After all, we already discretised the time derivative, Poisson brackets, and the Laplace operator.
Unfortunately, it is not that simple, as in the above potential formulation, additional derivatives appear within the Poisson brackets. To account for those properly, one has to do so at the level of the discrete action. Consequently, one has to repeat the whole derivation, work that is left for future research.

\section{Variational Discretisation}

We will base the derivation of the variational integrator for ideal MHD on the equations (\ref{eq:mhd_eqs_V}) and (\ref{eq:mhd_eqs_B}). Their respective components are
\begin{subequations}\label{eq:mhd_vi_equations_1}
\begin{align}
& \dfrac{\partial V_{x}}{\partial t} + V_{x} \, \partial_{x} V_{x} + V_{y} \, \partial_{y} V_{x} - B_{x} \, \partial_{x} B_{x} - B_{y} \, \partial_{y} B_{x} + \partial_{x} P = 0 , \\
& \dfrac{\partial V_{y}}{\partial t} + V_{x} \, \partial_{x} V_{y} + V_{y} \, \partial_{y} V_{y} - B_{x} \, \partial_{x} B_{y} - B_{y} \, \partial_{y} B_{y} + \partial_{y} P = 0 , \\
& \dfrac{\partial B_{x}}{\partial t} + V_{x} \, \partial_{x} B_{x} + V_{y} \, \partial_{y} B_{x} - B_{x} \, \partial_{x} V_{x} - B_{y} \, \partial_{y} V_{x} = 0 , \\
& \dfrac{\partial B_{y}}{\partial t} + V_{x} \, \partial_{x} B_{y} + V_{y} \, \partial_{y} B_{y} - B_{x} \, \partial_{x} V_{y} - B_{y} \, \partial_{y} V_{y} = 0 .
\end{align}
\end{subequations}

As we will see in the next section, we have to use a staggered grid approach for the discretisation of the MHD equations.
To be able to discretise (\ref{eq:mhd_vi_equations_1}) on a single grid cell as depicted in figure \ref{fig:mhd_staggered_grid}, we have to transform these equations as follows.
At first, consider the momentum equations.
Like in section (\ref{sec:mhd_lie_derivative}), add and subtract $V_{y} \, \partial_{x} V_{y}$ in the $x$ component, add and subtract $V_{x} \, \partial_{y} V_{x}$ in the $y$ component, and do the same for the corresponding $B$ terms, such that
\begin{subequations}\label{eq:mhd_vi_equations_2}
\begin{align}
\label{eq:mhd_vi_equations_2a}
& \partial_{t} V_{x} + V_{y} \, \big( \partial_{y} V_{x} - \partial_{x} V_{y} \big) - B_{y} \, \big( \partial_{y} B_{x} - \partial_{x} B_{y} \big) + \partial_{x} \Big[ P + \dfrac{1}{2} \norm{V}^{2} - \dfrac{1}{2} \norm{B}^{2} \Big] = 0 , \\
\label{eq:mhd_vi_equations_2b}
& \partial_{t} V_{y} + V_{x} \, \big( \partial_{x} V_{y} - \partial_{y} V_{x} \big) - B_{x} \, \big( \partial_{x} B_{y} - \partial_{y} B_{x} \big) + \partial_{y} \Big[ P + \dfrac{1}{2} \norm{V}^{2} - \dfrac{1}{2} \norm{B}^{2} \Big] = 0 .
\end{align}
\end{subequations}

In the round brackets, we find the rotation operator, which is the $z$ component of the curl, i.e.,
\begin{align}\label{eq:mhd_vi_equations_3}
(\nabla \times V)_{z} \equiv \partial_{x} V_{y} - \partial_{y} V_{x} ,
\end{align}

and in the square brackets, we find a modified pressure term
\begin{align}\label{eq:mhd_vi_equations_4}
\otilde{P} = P + \dfrac{1}{2} \norm{V}^{2} - \dfrac{1}{2} \norm{B}^{2} = p + \dfrac{1}{2} \norm{V}^{2} ,
\end{align}

that contains the kinetic instead of the magnetic energy in addition to the gas pressure $p$.
Next, consider the induction equation. By the divergence free constraint of $V$ and $B$ make the replacements $\partial_{x} V_{x} = - \partial_{y} V_{y}$ in equation (\ref{eq:mhd_vi_equations_2a}), $\partial_{y} V_{y} = - \partial_{x} V_{x}$ in equation (\ref{eq:mhd_vi_equations_2b}), and corresponding terms for $B$, such that
\begin{subequations}\label{eq:mhd_vi_equations_5}
\begin{align}
& \partial_{t} B_{x} + V_{y} \, \partial_{y} B_{x} - V_{x} \, \partial_{y} B_{y} + B_{x} \, \partial_{y} V_{y} - B_{y} \, \partial_{y} V_{x} = 0 , \\
& \partial_{t} B_{y} + V_{x} \, \partial_{x} B_{y} - V_{y} \, \partial_{x} B_{x} + B_{y} \, \partial_{x} V_{x} - B_{x} \, \partial_{x} V_{y} = 0 ,
\end{align}
\end{subequations}

which can be condensed into
\begin{subequations}\label{eq:mhd_vi_equations_6}
\begin{align}
& \partial_{t} B_{x} + \partial_{y} \big( V_{y} B_{x} - V_{x} B_{y} \big) = 0 , \\
& \partial_{t} B_{y} + \partial_{x} \big( V_{x} B_{y} - V_{y} B_{x} \big) = 0 ,
\end{align}
\end{subequations}

where the term in brackets is of course the $z$ component of $V \times B$.
To make the following derivations more tractable, we define two operators with components
\begin{subequations}\label{eq:mhd_vi_equations_7}
\begin{align}
\psi_{x} (V)   &\equiv - V_{y} \, \big( \partial_{x} V_{y} - \partial_{y} V_{x} \big) , &
\phi_{x} (V,B) &\equiv - \partial_{y} \big( V_{x} B_{y} - V_{y} B_{x} \big) , \\
\psi_{y} (V)   &\equiv + V_{x} \, \big( \partial_{x} V_{y} - \partial_{y} V_{x} \big) , &
\phi_{y} (V,B) &\equiv + \partial_{x} \big( V_{x} B_{y} - V_{y} B_{x} \big) ,
\end{align}
\end{subequations}

which is the same notation used by \citeauthor{Gawlik:2011} \cite{Gawlik:2011}.
With that, the incompressible, ideal MHD equations read
\begin{subequations}\label{eq:mhd_vi_equations_8}
\begin{align}
& \partial_{t} V + \psi(V) - \psi(B) + \nabla \otilde{P} = 0, \\
& \partial_{t} B + \phi(V,B) = 0 .
\end{align}
\end{subequations}

The corresponding extended Lagrangian, is readily written upon introducing three auxiliary variables $\alpha$, $\beta$, $\gamma$, where $\alpha$ and $\beta$ are vector fields and $\gamma$ is a scalar field.
The extended Lagrangian for (\ref{eq:mhd_vi_equations_8}) reads
\begin{align}\label{eq:mhd_vi_equations_9}
\mcal{L} ( V, B, \otilde{P}, \alpha, \beta, \gamma)
= \alpha \cdot \Big[ \partial_{t} V + \psi(V) - \psi(B) + \nabla \otilde{P} \Big]
+ \beta  \cdot \Big[ \partial_{t} B + \phi(V,B) \Big]
+ \gamma \, \Big[ \nabla \cdot V \Big] .
\end{align}

This will be the basis for the derivation of the variational integrator.
It is worth mentioning that $\nabla \cdot V = 0$ is a dynamical equation determining the pressure, c.f. comments after equation (\ref{eq:mhd_eqs_total_pressure}).
The Ibragimov multipliers $\alpha$, $\beta$ and $\gamma$ correspond to the physical variables $V$, $B$ and $\otilde{P}$, respectively.

\subsection{Staggered Grid}

Straight forward centred finite difference discretisations of the Navier-Stokes equation, where the components of the velocity vector and the pressure are located at the same grid points, are known to be prone to instabilities (see e.g. \citeauthor{Langtangen:2002} \cite{Langtangen:2002} or \citeauthor{McDonough:2007} \cite{McDonough:2007}).
The pressure often becomes highly oscillatory as a symmetric difference operator, e.g., with stencil $[ -1 \; \hphantom{-}0 \; +1 ]$, annihilates pressures which oscillate between $+1$ and $-1$ between neighbouring grid points. This is often referred to as \emph{checkerboarding}.

\begin{figure}[tbh]
\centering
\subfloat[Divergence Constraint]{\label{fig:mhd_staggered_grid_div}
\includegraphics[width=.43\textwidth]{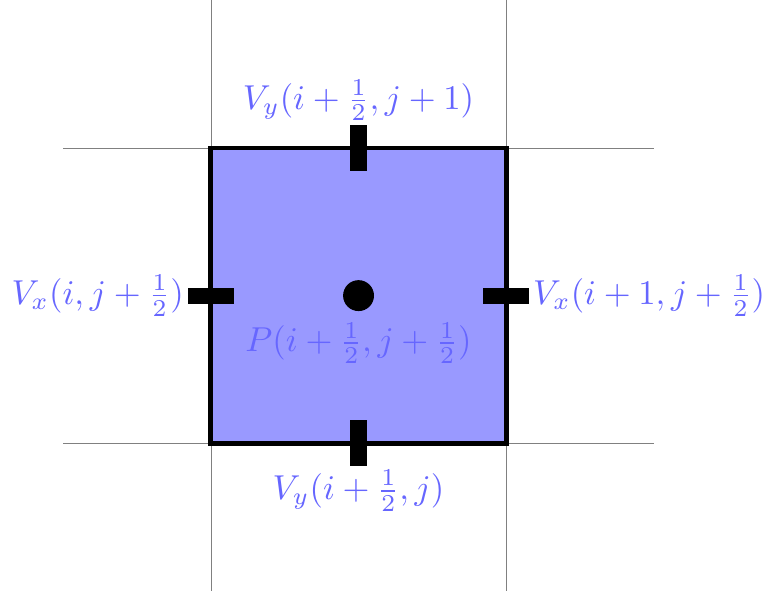}
}
\subfloat[Momentum Equation]{\label{fig:mhd_staggered_grid_momentum}
\includegraphics[width=.44\textwidth]{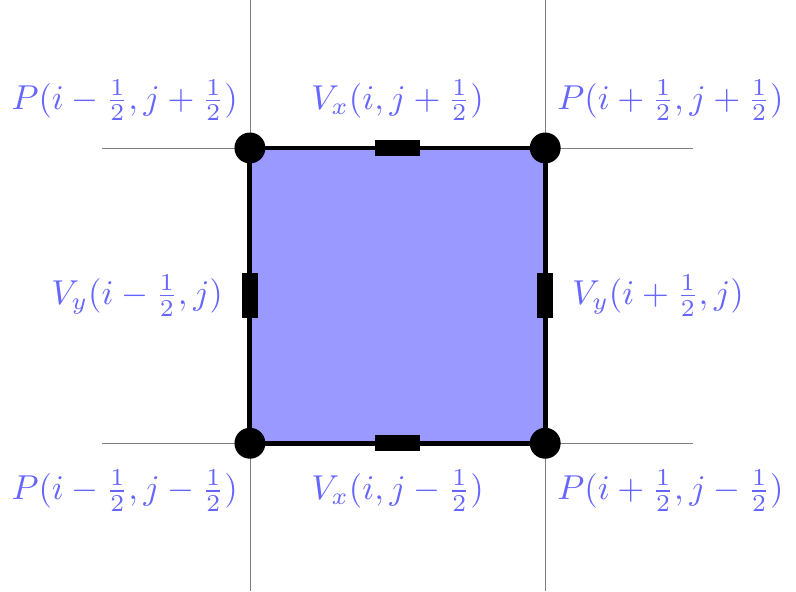}
}
\caption{Staggered Grid. Left: Natural positions for the pressure and the velocity components. Right: Dual grid.}
\label{fig:mhd_staggered_grid}
\end{figure}

An efficacious remedy for this problem is the introduction of a staggered grid, where the pressure is located at the centre of a grid cell and the velocity components at the vertices, like it is depicted in figure \ref{fig:mhd_staggered_grid_div}.
The location of the physical quantities comes natural when viewed as differential forms. The pressure is a zero-form and is therefore collocated at the centre of a cell. The velocity (and in two dimensions also the magnetic field) is a one-form and is therefore collocated at the edges of a cell.

This can also be seen by considering the discrete divergence-free constraint of the velocity field
\begin{align}\label{eq:mhd_staggered_grid_divergence}
\dfrac{V_{x} (i, j+\tfrac{1}{2}, k) - V_{x} (i-1, j+\tfrac{1}{2}, k)}{h_{x}} + \dfrac{V_{y} (i+\tfrac{1}{2}, j, k) - V_{y} (i+\tfrac{1}{2}, j-1, k)}{h_{y}} = 0 ,
\end{align}

which is defined such that the logical location of the divergence coincides with the location of the pressure.
The function of the pressure in incompressible fluid dynamics can be described as taking care of the divergence of the velocity field.
By this discretisation, only one pressure point takes care of the divergence of the neighbouring velocity points. And as the divergence is computed by a simple forward finite difference, without symmetric stencil, checkerboarding will not be an issue.

On this grid, the time derivatives are defined pointwise (without spatial averaging)
\begin{subequations}
\begin{alignat}{5}
& \dfrac{\partial V_{x}}{\partial t} &
& \quad \rightarrow \quad &
& \bigg( \dfrac{\partial V_{x}}{\partial t} \bigg)_{i, j+\tfrac{1}{2}, k+\tfrac{1}{2}} &
& \equiv \; &
& \dfrac{V_{x} (i , j+\tfrac{1}{2} , k+1) - V_{x} (i , j+\tfrac{1}{2}, k)}{h_{t}} ,
\\
& \dfrac{\partial V_{y}}{\partial t} &
& \quad \rightarrow \quad &
& \bigg( \dfrac{\partial V_{y}}{\partial t} \bigg)_{i+\tfrac{1}{2}, j, k+\tfrac{1}{2}} &
& \equiv \; &
& \dfrac{V_{y} (i+\tfrac{1}{2} , j , k+1) - V_{y} (i+\tfrac{1}{2} , j , k)}{h_{t}} .
\end{alignat}
\end{subequations}

For the spatial derivatives of the vectors, we use a midpoint averaging with respect to time, i.e.,
\begin{subequations}
\begin{alignat}{5}
\nonumber
& \dfrac{\partial V_{x}}{\partial x} &
& \quad \rightarrow \quad &
& \bigg( \dfrac{\partial V_{x}}{\partial x} \bigg)_{i+\tfrac{1}{2}, j+\tfrac{1}{2}, k+\tfrac{1}{2}} &
& \equiv \; &
& \dfrac{1}{2} \, \bigg[ \dfrac{V_{x} (i+1, j+\tfrac{1}{2}, k  ) - V_{x} (i, j+\tfrac{1}{2}, k  )}{h_{x}}
\\
&&&&&&
& + \; &
& \dfrac{V_{x} (i+1, j+\tfrac{1}{2}, k+1) - V_{x} (i, j+\tfrac{1}{2}, k+1)}{h_{x}} \bigg] ,
\hspace{-3em}
\\
\nonumber
& \dfrac{\partial V_{x}}{\partial y} &
& \quad \rightarrow \quad &
& \bigg( \dfrac{\partial V_{x}}{\partial y} \bigg)_{i, j, k+\tfrac{1}{2}} &
& \equiv \; &
& \dfrac{1}{2} \, \bigg[ \dfrac{V_{x} (i, j+\tfrac{1}{2}, k  ) - V_{x} (i, j-\tfrac{1}{2}, k  )}{h_{y}}
\\
&&&&&&
& + \, &
& \dfrac{V_{x} (i, j+\tfrac{1}{2}, k+1) - V_{x} (i, j-\tfrac{1}{2}, k+1)}{h_{y}} \bigg] ,
\hspace{-3em}
\end{alignat}
\end{subequations}
\begin{subequations}
\begin{alignat}{5}
\nonumber
& \dfrac{\partial V_{y}}{\partial x} &
& \quad \rightarrow \quad &
& \bigg( \dfrac{\partial V_{y}}{\partial x} \bigg)_{i, j, k+\tfrac{1}{2}} &
& \equiv \; &
& \dfrac{1}{2} \, \bigg[ \dfrac{V_{y} (i+\tfrac{1}{2}, j, k  ) - V_{y} (i-\tfrac{1}{2}, j, k  )}{h_{x}} \\
&&&&&&&&& \hspace{1em}
+ \dfrac{V_{y} (i+\tfrac{1}{2}, j, k+1) - V_{y} (i-\tfrac{1}{2}, j, k+1)}{h_{x}} \bigg] ,
\hspace{-3em}
\\
\nonumber
& \dfrac{\partial V_{y}}{\partial y} &
& \quad \rightarrow \quad &
& \bigg( \dfrac{\partial V_{y}}{\partial y} \bigg)_{i+\tfrac{1}{2}, j+\tfrac{1}{2}, k+\tfrac{1}{2}} &
& \equiv \; &
& \dfrac{1}{2} \, \bigg[ \dfrac{V_{y} (i+\tfrac{1}{2}, j+1, k  ) - V_{y} (i+\tfrac{1}{2}, j, k  )}{h_{y}}
\\
&&&&&&&&& \hspace{1em}
+ \dfrac{V_{y} (i+\tfrac{1}{2}, j+1, k+1) - V_{y} (i+\tfrac{1}{2}, j, k+1)}{h_{y}} \bigg] ,
\hspace{-3em}
\end{alignat}
\end{subequations}

Note, that the $x$ derivative of $V_{x}$ and the $y$ derivative of $V_{y}$ are defined on the grid in figure \ref{fig:mhd_staggered_grid_div}, while the $y$ derivative of $V_{x}$ and the $x$ derivative of $V_{y}$ are defined on the dual grid in figure \ref{fig:mhd_staggered_grid_momentum}.
The indices of the derivatives denote the logical collocation of the derivative, which is always the cell centre.

Derivatives of the pressure can only be defined on the dual grid, figure \ref{fig:mhd_staggered_grid_momentum}.
They are naturally defined on the edges of the cells.
The staggering approach is applied to $P$ also with respect to time, i.e., the pressure nodes are $(i+\tfrac{1}{2}, j+\tfrac{1}{2}, k+\tfrac{1}{2})$.
Taking all of this into account, we define
\begin{subequations}
\begin{alignat}{9}
& \dfrac{\partial P}{\partial x} && \quad \rightarrow \quad &
& \bigg( \dfrac{\partial P}{\partial x} \bigg)_{i, j+\tfrac{1}{2}, k+\tfrac{1}{2}} &
& \equiv \; &
& \dfrac{P (i+\tfrac{1}{2}, j+\tfrac{1}{2}, k+\tfrac{1}{2}) - P (i-\tfrac{1}{2}, j+\tfrac{1}{2}, k+\tfrac{1}{2})}{h_{x}} ,
\\
& \dfrac{\partial P}{\partial y} && \quad \rightarrow \quad &
& \bigg( \dfrac{\partial P}{\partial y} \bigg)_{i+\tfrac{1}{2}, j, k+\tfrac{1}{2}} &
& \equiv \; &
& \dfrac{P (i+\tfrac{1}{2}, j+\tfrac{1}{2}, k+\tfrac{1}{2}) - P (i+\tfrac{1}{2}, j-\tfrac{1}{2}, k+\tfrac{1}{2})}{h_{y}} .
\end{alignat}
\end{subequations}

Averages are needed only on the dual grid, so we are defining them only there and only for the vector fields. For $V$, the averaging is applied with respect to both, space and time,
\begin{subequations}
\begin{align}
\bracket{ V_{x} }_{i, j, k+\tfrac{1}{2}}
&\equiv
\dfrac{V_{x} ( i , j-\tfrac{1}{2} , k  ) + V_{x} ( i , j+\tfrac{1}{2} , k  )}{4}
+ \dfrac{V_{x} ( i , j-\tfrac{1}{2} , k+1) + V_{x} ( i , j+\tfrac{1}{2} , k+1)}{4} ,
\\
\bracket{ V_{y} }_{i, j, k+\tfrac{1}{2}}
&\equiv
\dfrac{V_{y} ( i-\tfrac{1}{2} , j , k  ) + V_{y} ( i+\tfrac{1}{2} , j , k  )}{4}
+ \dfrac{V_{y} ( i-\tfrac{1}{2} , j , k+1) + V_{y} ( i+\tfrac{1}{2} , j , k+1)}{4} ,
\end{align}
\end{subequations}

but as $\alpha$ will be collocated at $k+\tfrac{1}{2}$ (see comment in the next section), its averages do not feature a time average,
\begin{subequations}
\begin{align}
\bracket{ \alpha_{x} }_{i, j, k+\tfrac{1}{2}}
&\equiv
\dfrac{ \alpha_{x} ( i , j-\tfrac{1}{2} , k+\tfrac{1}{2}) + \alpha_{x} ( i , j+\tfrac{1}{2} , k+\tfrac{1}{2}) }{2} ,
\\
\bracket{ \alpha_{y} }_{i, j, k+\tfrac{1}{2}}
&\equiv
\dfrac{ \alpha_{y} ( i-\tfrac{1}{2} , j , k+\tfrac{1}{2}) + \alpha_{y} ( i+\tfrac{1}{2} , j , k+\tfrac{1}{2}) }{2} .
\end{align}
\end{subequations}

\subsection{Navier-Stokes Equation}

We start the derivation of the variational integrator by restricting our attention to the incompressible Navier-Stokes equation
\begin{align}
\partial_{t} V + \psi (V) + \nabla P &= 0 , &
\nabla \cdot V &= 0 , &
\end{align}

neglecting the force term through the magnetic field. The generalisation to magnetohydrodynamics is straight forward, as the magnetic field appears with the same advection term $\psi(B)$ as the velocity field, and the analytical expression of the pressure $P$ in terms of $v$ and $B^{2}$ does not play any role.
The action integral of the extended Lagrangian (\ref{eq:mhd_vi_equations_9}), reduced to this subsystem, is
\begin{align}\label{eq:mhd_navier_stokes_lagrangian}
\mcal{A} &= \int \Big[ ... + \alpha \cdot \big[ \partial_{t} V + \psi(V) + \nabla P \big] + \gamma \, \big[ \nabla \cdot V \big] + ... \Big] \, dt \, dx \, dy .
\end{align}

To be able to discretise all of the derivatives in the first term of the Lagrangian, we have to switch to the dual grid, as depicted in figure \ref{fig:mhd_staggered_grid_momentum}.
The time derivatives are approximated using the trapezoidal rule,
\begin{subequations}
\begin{align}
\alpha_{x} \, \big( \partial_{t} V_{x} \big) \quad \rightarrow \quad
& \alpha_{x} ( i , j-\tfrac{1}{2} , k+\tfrac{1}{2}) \, \bigg( \dfrac{\partial V_{x}}{\partial t} \bigg)_{i, j-\tfrac{1}{2}, k+\tfrac{1}{2}}
+ \alpha_{x} ( i , j+\tfrac{1}{2} , k+\tfrac{1}{2}) \, \bigg( \dfrac{\partial V_{x}}{\partial t} \bigg)_{i, j+\tfrac{1}{2}, k+\tfrac{1}{2}} , \\
\alpha_{y} \, \big( \partial_{t} V_{y} \big) \quad \rightarrow \quad
& \alpha_{y} ( i-\tfrac{1}{2} , j , k+\tfrac{1}{2}) \, \bigg( \dfrac{\partial V_{y}}{\partial t} \bigg)_{i-\tfrac{1}{2}, j, k+\tfrac{1}{2}}
+ \alpha_{y} ( i+\tfrac{1}{2} , j , k+\tfrac{1}{2}) \, \bigg( \dfrac{\partial V_{y}}{\partial t} \bigg)_{i+\tfrac{1}{2}, j, k+\tfrac{1}{2}} .
\end{align}
\end{subequations}

Here, we are exploiting the same ideas, explained in section \ref{sec:vlasov_variational}, but here that is implemented directly into the Lagrangian rather than in the final integrator. We omit the time averaging of $\alpha$, which implies that $\alpha$ is collocated at $k+\tfrac{1}{2}$, just as the time derivative.
We use a trapezoidal approximation to avoid spatial averaging of the time derivatives in the resulting scheme, as that might again lead to grid oscillations (checkerboarding), this time in the velocity field.
We apply the same approximation to the pressure gradient term, for the same reason, i.e., to avoid oscillations, and as the structure of the terms is identical (e.g., $\partial_{t} V_{x}$ and the $x$ component of the pressure gradient $\partial_{x} P$ are both objects collocated at the same logical position),
\begin{align}
\alpha_{x} \, \big( \partial_{x} P \big) \quad \rightarrow \quad
& \alpha_{x} ( i , j-\tfrac{1}{2} , k+\tfrac{1}{2}) \, \bigg( \dfrac{\partial P}{\partial x} \bigg)_{i, j-\tfrac{1}{2}, k+\tfrac{1}{2}}
+ \alpha_{x} ( i , j+\tfrac{1}{2} , k+\tfrac{1}{2}) \, \bigg( \dfrac{\partial P}{\partial x} \bigg)_{i, j+\tfrac{1}{2}, k+\tfrac{1}{2}} ,
\\
\alpha_{y} \, \big( \partial_{y} P \big) \quad \rightarrow \quad
& \alpha_{y} ( i-\tfrac{1}{2} , j , k+\tfrac{1}{2}) \, \bigg( \dfrac{\partial P}{\partial y} \bigg)_{i-\tfrac{1}{2}, j, k+\tfrac{1}{2}}
+ \alpha_{y} ( i+\tfrac{1}{2} , j , k+\tfrac{1}{2}) \, \bigg( \dfrac{\partial P}{\partial y} \bigg)_{i+\tfrac{1}{2}, j, k+\tfrac{1}{2}} .
\end{align}

As previously mentioned, the pressure is collocated at $k+\tfrac{1}{2}$, such that no time average of $P$ is needed.
The $\psi$ operator (\ref{eq:mhd_vi_equations_7}) is discretised by a midpoint approximation, both with respect to space and time, i.e.,
\begin{align}
\alpha_{x} \psi_{x} (V) \quad \rightarrow \quad
&- \bracket{ \alpha_{x} }_{i, j, k+\tfrac{1}{2}}
\, \bracket{ V_{x} }_{i, j, k+\tfrac{1}{2}}
\, \bigg[
\bigg( \dfrac{\partial V_{y}}{\partial x} \bigg)_{i, j, k+\tfrac{1}{2}}
- \bigg( \dfrac{\partial V_{x}}{\partial y} \bigg)_{i, j, k+\tfrac{1}{2}}
\bigg] ,
\\
\alpha_{y} \psi_{y} (V) \quad \rightarrow \quad
&+ \bracket{ \alpha_{y} }_{i, j, k+\tfrac{1}{2}}
\, \bracket{ V_{y} }_{i, j, k+\tfrac{1}{2}}
\, \bigg[
\bigg( \dfrac{\partial V_{y}}{\partial x} \bigg)_{i, j, k+\tfrac{1}{2}}
- \bigg( \dfrac{\partial V_{x}}{\partial y} \bigg)_{i, j, k+\tfrac{1}{2}}
\bigg] .
\end{align}

As, e.g., $\alpha_{x}$ and $V_{y}$ are objects collocated at different logical positions, they cannot be multiplied directly, but products can only be computed of their averages, which are collocated at the cell centres. Therefore, we have to use a midpoint approximation in this term. The rotation is logically located at the centre of the cell, therefore posing no problems.

The discretisation of the second term in (\ref{eq:mhd_navier_stokes_lagrangian}) is implemented on the grid in figure \ref{fig:mhd_staggered_grid_div}.
Recognising that $\gamma$ is a scalar field and thus collocated at the same position as the pressure, the discretisation follows directly from (\ref{eq:mhd_staggered_grid_divergence}), i.e.,
\begin{align}
\gamma \, \big( \nabla \cdot V \big) \quad \rightarrow \quad
\gamma ( i+\tfrac{1}{2} , j+\tfrac{1}{2} , k ) \, \bigg[
\bigg( \dfrac{\partial V_{x}}{\partial y} \bigg)_{i, j-\tfrac{1}{2}, \obar{k}}
+ \bigg( \dfrac{\partial V_{x}}{\partial y} \bigg)_{i, \obar{j}, \obar{k}}
\bigg] .
\end{align}

Summing up all contributions gives the discrete Lagrangian.

\subsection{Induction Equation}

Now we consider those terms of the action that will yield the induction equation.
\begin{align}
\mcal{A} &= \int \Big[ ... + \beta \cdot \big[ \partial_{t} B + \phi(V,B) \big] + ... \Big] \, dt \, dx \, dy .
\end{align}

To find a discretisation of the operator $\phi(V,B)$ on a single grid cell, we have to first do a partial integration, such that
\begin{align}
\mcal{A}
\nonumber
&= \int \Big[ ... + \beta_{x} \, \big[ \partial_{t} B_{x} - \partial_{y} ( V_{x} B_{y} - V_{y} B_{x} ) \big]
+ \beta_{y} \, \big[ \partial_{t} B_{y} + \partial_{x} ( V_{x} B_{y} - V_{y} B_{x} ) \big] + ... \Big] \, dt \, dx \, dy \\
\nonumber
&= \int \Big[ ... + \beta_{x} \, \partial_{t} B_{x} + ( \partial_{y} \beta_{x} ) ( V_{x} B_{y} - V_{y} B_{x} ) \\
& \hspace{8em}
+ \beta_{y} \, \partial_{t} B_{y} - ( \partial_{x} \beta_{y} ) ( V_{x} B_{y} - V_{y} B_{x} ) + ... \Big] \, dt \, dx \, dy
.
\end{align}

The discretisation of the time derivative is the same as in the case of the momentum equation, i.e., using a trapezoidal rule,
\begin{align}
\beta_{x} \, \big( \partial_{t} B_{x} \big) \;\; \rightarrow \;\;\;
& \beta_{x} ( i , j-\tfrac{1}{2} , k+\tfrac{1}{2}) \, \bigg( \dfrac{\partial B_{x}}{\partial t} \bigg)_{i, j-\tfrac{1}{2}, k+\tfrac{1}{2}}
+ \beta_{x} ( i , j+\tfrac{1}{2} , k+\tfrac{1}{2}) \, \bigg( \dfrac{\partial B_{x}}{\partial t} \bigg)_{i, j+\tfrac{1}{2}, k+\tfrac{1}{2}} , \\
\beta_{y} \, \big( \partial_{t} B_{y} \big) \;\; \rightarrow \;\;\;
& \beta_{y} ( i-\tfrac{1}{2} , j , k+\tfrac{1}{2}) \, \bigg( \dfrac{\partial B_{y}}{\partial t} \bigg)_{i-\tfrac{1}{2}, j, k+\tfrac{1}{2}}
+ \beta_{y} ( i+\tfrac{1}{2} , j , k+\tfrac{1}{2}) \, \bigg( \dfrac{\partial B_{y}}{\partial t} \bigg)_{i+\tfrac{1}{2}, j, k+\tfrac{1}{2}} .
\end{align}

The factors of the operator  $\phi$ are collocated at different positions of the grid, e.g., in the first equation, $\partial_{y} \beta_{x}$ is collocated at $(i, j)$, $V_{x}$ and $B_{x}$ are collocated at $(i, j+\tfrac{1}{2})$, and $V_{y}$ and $B_{y}$ are collocated at $(i+\tfrac{1}{2}, j)$.
Therefore, we have to use a midpoint rule, i.e.,
\begin{align}
( \partial_{y} \beta_{x} ) ( V_{x} B_{y} - V_{y} B_{x} ) \; \rightarrow \;\;
&  \bigg( \dfrac{\partial \beta_{x}}{\partial y} \bigg)_{i, j, k+\tfrac{1}{2}}
\, \bigg[
\bracket{ V_{x} }_{i, j, k+\tfrac{1}{2}} \, \bracket{ B_{y} }_{i, j, k+\tfrac{1}{2}}
- \bracket{ V_{y} }_{i, j, k+\tfrac{1}{2}} \, \bracket{ B_{x} }_{i, j, k+\tfrac{1}{2}}
\bigg] ,
\\
( \partial_{x} \beta_{y} ) ( V_{x} B_{y} - V_{y} B_{x} ) \; \rightarrow \;\;
&  \bigg( \dfrac{\partial \beta_{y}}{\partial x} \bigg)_{i, j, k+\tfrac{1}{2}}
\, \bigg[
\bracket{ V_{x} }_{i, j, k+\tfrac{1}{2}} \, \bracket{ B_{y} }_{i, j, k+\tfrac{1}{2}}
- \bracket{ V_{y} }_{i, j, k+\tfrac{1}{2}} \, \bracket{ B_{x} }_{i, j, k+\tfrac{1}{2}}
\bigg] .
\end{align}

With that, we have all the ingredients for a complete discretisation of the action integral corresponding to (\ref{eq:mhd_vi_equations_9}).

\subsection{Variational Integrator}

Computing the variation of the discrete action with respect to $\alpha$, $\beta$ and $\gamma$, we obtain the discrete ideal MHD equations
\begin{align}
\nonumber
& \dfrac{V_{x}^{i,j+1/2,k+1} - V_{x}^{i,j+1/2,k}}{h_{t}} + \psi_{x}^{i,j+1/2} \bigg( \dfrac{V^{k} + V^{k+1}}{2} \bigg) - \psi_{x}^{i,j+1/2} \bigg( \dfrac{B^{k} + B^{k+1}}{2} \bigg) \\
& \hspace{18em}
+ \dfrac{P^{i+1/2,j+1/2,k+1/2} - P^{i-1/2,j+1/2,k+1/2}}{h_{x}} = 0 ,
\\
\nonumber
& \dfrac{V_{y}^{i+1/2,j,k+1} - V_{y}^{i+1/2,j,k}}{h_{t}} + \psi_{y}^{i+1/2,j} \bigg( \dfrac{V^{k} + V^{k+1}}{2} \bigg) - \psi_{y}^{i+1/2,j} \bigg( \dfrac{B^{k} + B^{k+1}}{2} \bigg) \\
& \hspace{18em}
+ \dfrac{P^{i+1/2,j+1/2,k+1/2} - P^{i+1/2,j-1/2,k+1/2}}{h_{y}} = 0 ,
\\
& \dfrac{B_{x}^{i,j+1/2,k+1} - B_{x}^{i,j+1/2,k}}{h_{t}} + \phi_{x}^{i,j+1/2} \bigg( \dfrac{V^{k} + V^{k+1}}{2} , \dfrac{B^{k} + B^{k+1}}{2} \bigg) = 0 ,
\\
& \dfrac{B_{y}^{i+1/2,j,k+1} - B_{y}^{i+1/2,j,k}}{h_{t}} + \phi_{y}^{i+1/2,j} \bigg( \dfrac{V^{k} + V^{k+1}}{2} , \dfrac{B^{k} + B^{k+1}}{2} \bigg) = 0 ,
\\
\nonumber
& \bigg[ \dfrac{V_{x}^{i+1,j+1/2,k} - V_{x}^{i,j+1/2,k}}{h_{x}} + \dfrac{V_{x}^{i+1,j+1/2,k+1} - V_{x}^{i,j+1/2,k+1}}{h_{x}} \\
& \hspace{10em}
+ \dfrac{V_{y}^{i+1/2,j+1,k} - V_{y}^{i+1/2,j,k}}{h_{y}} + \dfrac{V_{y}^{i+1/2,j+1,k+1} - V_{y}^{i+1/2,j,k+1}}{h_{y}} \bigg] = 0
.
\end{align}

with the discrete operators defined by
\begin{align}
\psi_{x}^{i,j+1/2} (V^{k+1/2})
\nonumber
&=
- \bracket{ V_{x} }_{i, j,   k+\tfrac{1}{2}} \,
\bigg[
\bigg( \dfrac{\partial V_{y}}{\partial x} \bigg)_{i, j,   k+\tfrac{1}{2}}
- \bigg( \dfrac{\partial V_{x}}{\partial y} \bigg)_{i, j,   k+\tfrac{1}{2}}
\bigg] \\
& \hspace{8em}
- \bracket{ V_{x} }_{i, j+1, k+\tfrac{1}{2}} \,
\bigg[
\bigg( \dfrac{\partial V_{y}}{\partial x} \bigg)_{i, j+1, k+\tfrac{1}{2}}
- \bigg( \dfrac{\partial V_{x}}{\partial y} \bigg)_{i, j+1, k+\tfrac{1}{2}}
\bigg] ,
\\
\psi_{y}^{i+1/2,j} (V^{k+1/2})
\nonumber
&=
+ \bracket{ V_{y} }_{i,   j, k+\tfrac{1}{2}} \,
\bigg[
\bigg( \dfrac{\partial V_{y}}{\partial x} \bigg)_{i,   j, k+\tfrac{1}{2}}
- \bigg( \dfrac{\partial V_{x}}{\partial y} \bigg)_{i,   j, k+\tfrac{1}{2}}
\bigg]  \\
& \hspace{8em}
+ \bracket{ V_{y} }_{i+1, j, k+\tfrac{1}{2}} \,
\bigg[
\bigg( \dfrac{\partial V_{y}}{\partial x} \bigg)_{i+1, j, k+\tfrac{1}{2}}
- \bigg( \dfrac{\partial V_{x}}{\partial y} \bigg)_{i+1, j, k+\tfrac{1}{2}}
\bigg] ,
\end{align}

and
\begin{align}
\phi_{x}^{i,j+1/2} (V^{k+1/2}, B^{k+1/2})
\nonumber
&=
+ \bigg[
\bracket{ V_{x} }_{i, j+1, k+\tfrac{1}{2}} \, \bracket{ B_{y} }_{i, j+1, k+\tfrac{1}{2}}
- \bracket{ V_{y} }_{i, j+1, k+\tfrac{1}{2}} \, \bracket{ B_{x} }_{i, j+1, k+\tfrac{1}{2}}
\bigg] \\
& \hspace{5em}
- \bigg[
\bracket{ V_{x} }_{i, j,   k+\tfrac{1}{2}} \, \bracket{ B_{y} }_{i, j,   k+\tfrac{1}{2}}
- \bracket{ V_{y} }_{i, j,   k+\tfrac{1}{2}} \, \bracket{ B_{x} }_{i, j,   k+\tfrac{1}{2}}
\bigg] ,
\\
\phi_{y}^{i+1/2,j} (V^{k+1/2}, B^{k+1/2})
\nonumber
&=
- \bigg[
\bracket{ V_{x} }_{i+1, j, k+\tfrac{1}{2}} \, \bracket{ B_{y} }_{i+1, j, k+\tfrac{1}{2}}
- \bracket{ V_{y} }_{i+1, j, k+\tfrac{1}{2}} \, \bracket{ B_{x} }_{i+1, j, k+\tfrac{1}{2}}
\bigg] \\
& \hspace{5em}
+ \bigg[
\bracket{ V_{x} }_{i,   j, k+\tfrac{1}{2}} \, \bracket{ B_{y} }_{i,   j, k+\tfrac{1}{2}}
- \bracket{ V_{y} }_{i,   j, k+\tfrac{1}{2}} \, \bracket{ B_{x} }_{i,   j, k+\tfrac{1}{2}}
\bigg] .
\end{align}

This discretisation of the operators $\psi$ and $\phi$ is the very same as those found by \citeauthor{Gawlik:2011} \cite{Gawlik:2011} and \citeauthor{LiuWang:2001} \cite{LiuWang:2001}.
However, \citeauthor{Gawlik:2011} resolve the nonlinearity in a different way, such that in their scheme, cross helicity is preserved exactly, but the energy error shows a monotonic growth. In our scheme, energy is preserved exactly as is cross helicity (up to machine accuracy).
They follow a different but related path in their derivation, based on Euler-Poincaré reduction. The crucial difference is, that in their method, only the velocity field is treated variational, and the magnetic field is treated as a quantity advected with the velocity field, whereas in our method, the velocity field and the magnetic field are treated on equal footing, fully variationally.
The scheme of \citeauthor{LiuWang:2001} uses an explicit Runge-Kutta method for time integration, so that conservation laws are broken in long time simulations.

Note the absence of any spatial averaging of the time derivatives and the pressure gradient. This is on purpose, as we wanted to prevent the emergence of grid-scale oscillations in the fields by the introduction of the staggered grid.

\section{Numerical Examples}

In this section, we consider four quite different examples of ideal magnetohydrodynamics problems taken from the previous literature \cite{CordobaMarliani:2000, GardinerStone:2005, Gawlik:2011}:
Alfvénic waves, the passive advection of a magnetic loop, the development of current sheaths in an Orszag-Tang vortex, and the development of magnetic islands along a current sheath.

\subsection{Diagnostics}

\subsubsection{Energy}

The total energy of the system is the sum of kinetic energy and magnetic energy, which are computed by
\begin{align}
E_{\text{kin},d}^{k} = h_{x} \, h_{y} \, \dfrac{1}{2} \sum \limits_{i,j} V_{x}^{2} (i, j+1/2, k) + h_{x} \, h_{y} \, \dfrac{1}{2} \sum \limits_{i,j} V_{y}^{2} (i+1/2, j, k) \\
E_{\text{mag},d}^{k} = h_{x} \, h_{y} \, \dfrac{1}{2} \sum \limits_{i,j} B_{x}^{2} (i, j+1/2, k) + h_{x} \, h_{y} \, \dfrac{1}{2} \sum \limits_{i,j} B_{y}^{2} (i+1/2, j, k) .
\end{align}

As there is no dissipation term in the ideal MHD equations, the total energy is always preserved.

\subsubsection{Cross Helicity}

The cross helicity is the integral of the scalar product of the velocity and magnetic fields,
\begin{align}
H_{d}^{k}
\nonumber
&= h_{x} \, h_{y} \, \dfrac{1}{2} \sum \limits_{i,j} V_{x} (i, j+1/2, k) \, B_{x} (i, j+1/2, k) \\
&+ h_{x} \, h_{y} \, \dfrac{1}{2} \sum \limits_{i,j} V_{y} (i+1/2, j, k) \, B_{y} (i+1/2, j, k) .
\end{align}

In ideal MHD, the parallel components of the velocity and magnetic fields do not interact, such that the integral of their product over the spatial domain stays constant.

\subsubsection{Vector Potential}

In two dimensions, the magnetic field is given by
\begin{align}
B_{x} &= \partial_{y} A &
& \text{and} &
B_{y} &= - \partial_{x} A . & &&
\end{align}

The magnetic potential is collocated at the cell centres of figure \ref{fig:mhd_staggered_grid_momentum}.
Therefore these equations discretise as
\begin{subequations}
\begin{align}
B_{x} (i, j+\tfrac{1}{2}) &=   \dfrac{A (i, j+1) - A (i,j)}{h_{y}}   &
& \text{and} &
B_{y} (i, j+\tfrac{1}{2}) &= - \dfrac{A (i+1, j) - A (i,j)}{h_{x}} . &
\end{align}
\end{subequations}

The potential field can be reconstructed by fixing $A(1,1)$ and looping over the whole grid, computing
\begin{align}
A (i, j+1) &= A (i,j) + h_{y} \, B_{x} (i,j+\tfrac{1}{2})   &
& \text{and} &
A (i+1, j) &= A (i,j) - h_{x} \, B_{x} (i,j+\tfrac{1}{2}) , &
\end{align}

using the first equation to compute columns and the second to jump between rows, or the other way around.
To which value $A(1,1)$ is fixed is not important as $A$ is determined only up to a constant.
The contour lines of the magnetic potential $A$ correspond to field lines of the magnetic field $B$. Hence, $A$ is an important diagnostic.

\subsubsection{Current Density}

The current is given by the curl of the magnetic field, or in two dimensions by the $z$ component of the curl. The discrete version of that is
\begin{align}\label{eq:mhd_diagnostics_current_density}
J (i,j,k)
= \dfrac{B_{y} (i+1/2, j, k) - B_{y} (i-1/2, j, k)}{h_{x}}
- \dfrac{B_{x} (i, j+1/2, k) - B_{x} (i, j-1/2, k)}{h_{y}} .
\end{align}

As the vector potential, the current is collocated at cell centres.

\subsection{Alfvén Waves}

In the first example, we consider a travelling Alfvén wave, initialised by
\begin{align*}
V_{x} &= 0 , &
V_{y} &= V_{0} \, \sin (\pi x) , &
B_{x} &= B_{0} , &
B_{y} &= B_{0} \, \sin (\pi x) , &
P &= 0.1 , &
\end{align*}

with $V_{0} = 1$ and $B_{0} = 1$.
The simulation domain is $[0,2] \times [0,2]$ with periodic boundaries and a resolution of $n_{x} \times n_{y} = 30 \times 30$. The timestep is $h_{t} = 0.1$ in units of the Alfvén time (i.e., the Alfvén velocity is one).

\begin{figure}[tb]
\centering
\includegraphics[width=\textwidth]{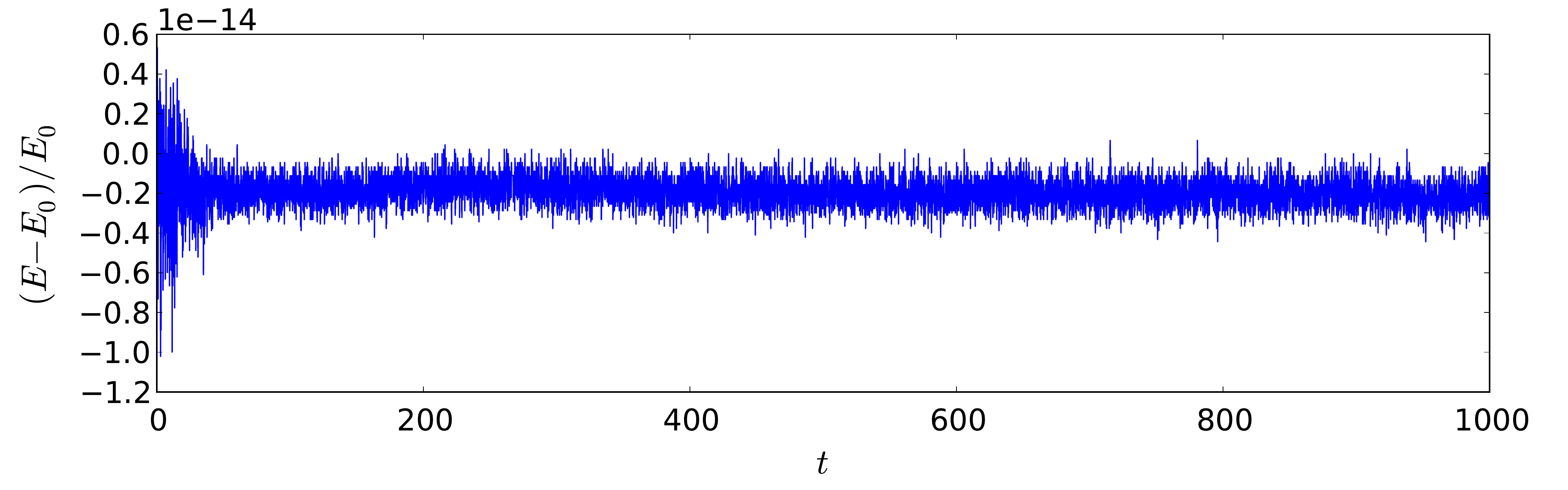}
\includegraphics[width=\textwidth]{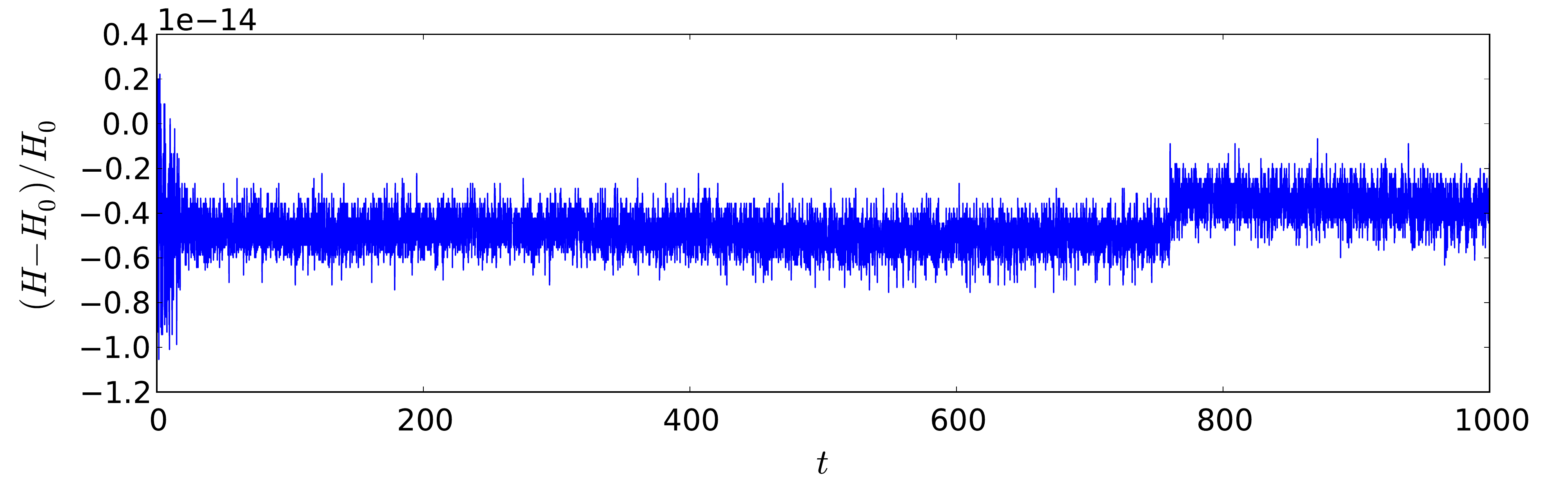}
\caption{Conservation of energy and cross helicity for a travelling Alfvén wave.}
\label{fig:mhd_alfven_wave_travelling}
\end{figure}

Although this example is rather simple, the results of our variational integrator are already remarkable.
Figure \ref{fig:mhd_alfven_wave_travelling} shows the time traces of the errors in the total energy and the cross helicity. For most of the simulation, the amplitude of the oscillations is $\mcal{O}(10^{-15})$, i.e., machine precision.
The sudden jump in the error of the cross helicity at about $t = 800$ might be explained by the residual of the Newton iteration being slightly larger for some timesteps than it is during the rest of the simulation.
We want to stress that during the runtime of $1000$ characteristic times there is no change in the energy within the machine accuracy.
It appears as if the Alfvén wave would continue travelling practically forever.
It also worth mentioning, that this is a fully nonlinear wave, i.e., the amplitudes of the perturbations of the magnetic field as well as the velocity field are $\mcal{O}(1)$.

\subsection{Loop Advection}

We now consider a case with very small magnetic field, such that the momentum and induction equations are almost decoupled and the magnetic field is passively advected by the fluid.
The initial conditions are
\begin{align*}
V_{x} &= V_{0} \, \cos (\theta) , &
V_{y} &= V_{0} \, \sin (\theta) , &
B_{x} &= \hphantom{-}\partial_{y} A, &
B_{y} &= - \partial_{x} A, &
P &= 1.0 , &
\end{align*}

with
\begin{align*}
A &= A_{0} \, ( R - \sqrt{x^{2} + y^{2}} ) , &
V_{0} &= \sqrt{5}, &
\theta &= \tan^{-1} (0.5) , &
A_{0} &= 10^{-3}, &
R &= 0.3 . &
\end{align*}

The spatial domain is $[-1, +1] \times [-0.5, +0,5]$ with periodic boundaries. We consider two resolutions, $128 \times 64$ and $256 \times 128$, and a timestep $h_{t} = 0.01$.

\begin{figure}
\centering
\includegraphics[width=\textwidth]{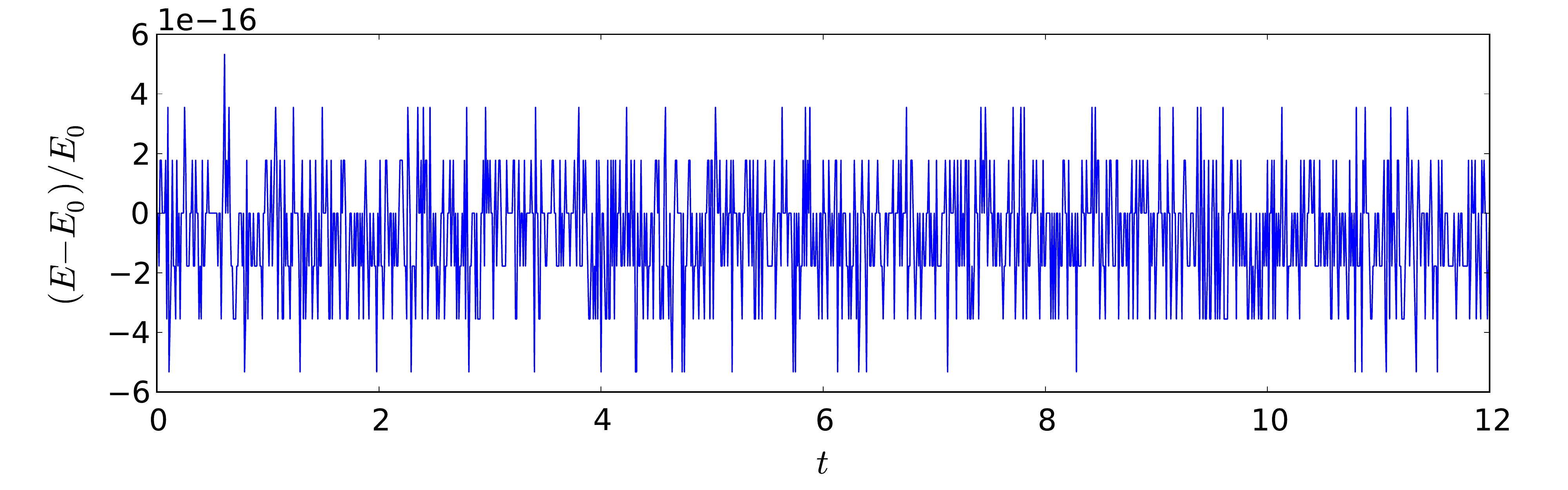}
\includegraphics[width=\textwidth]{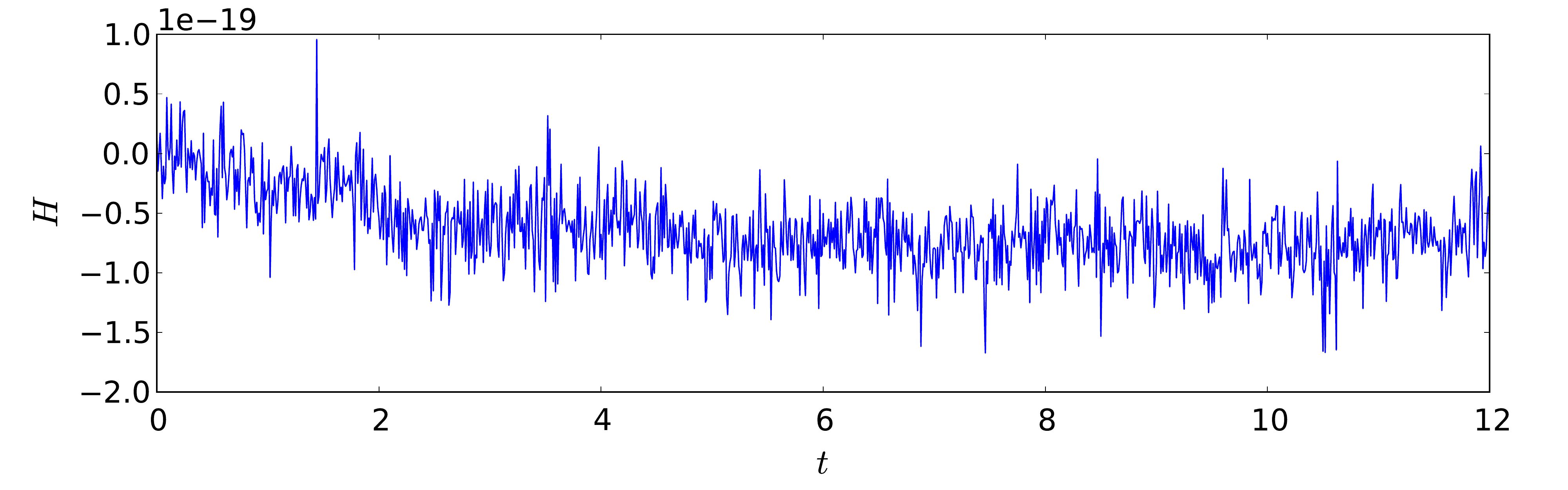}
\caption{Loop advection, $128 \times 64$ grid points. Conservation of energy and cross helicity.}
\label{fig:loop_advection_64x128_errors}
\end{figure}

\begin{figure}
\centering
\includegraphics[width=\textwidth]{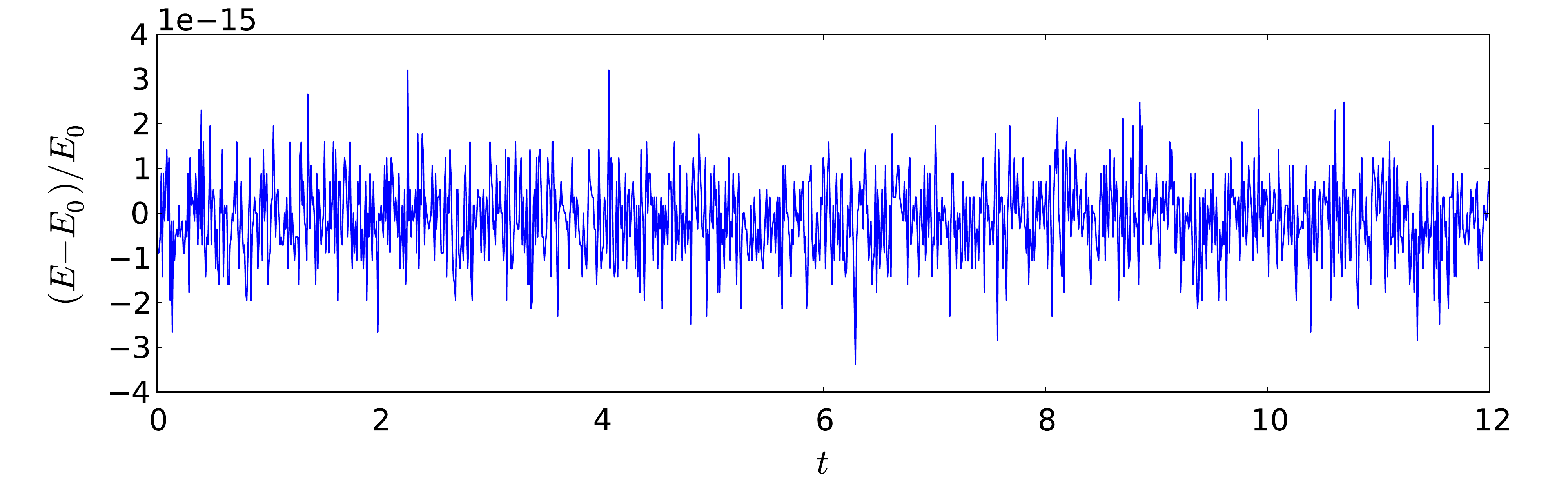}
\includegraphics[width=\textwidth]{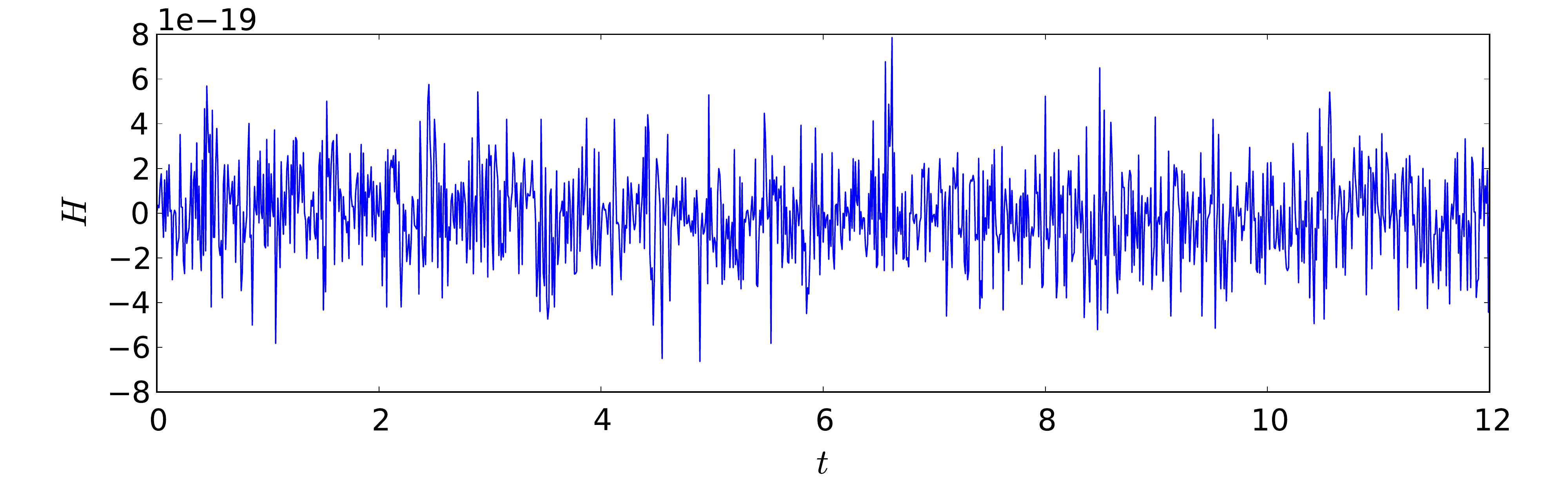}
\caption{Loop advection, $256 \times 128$ grid points. Conservation of energy and cross helicity.}
\label{fig:loop_advection_128x256_errors}
\end{figure}

\begin{figure}
\centering
\includegraphics[width=.48\textwidth]{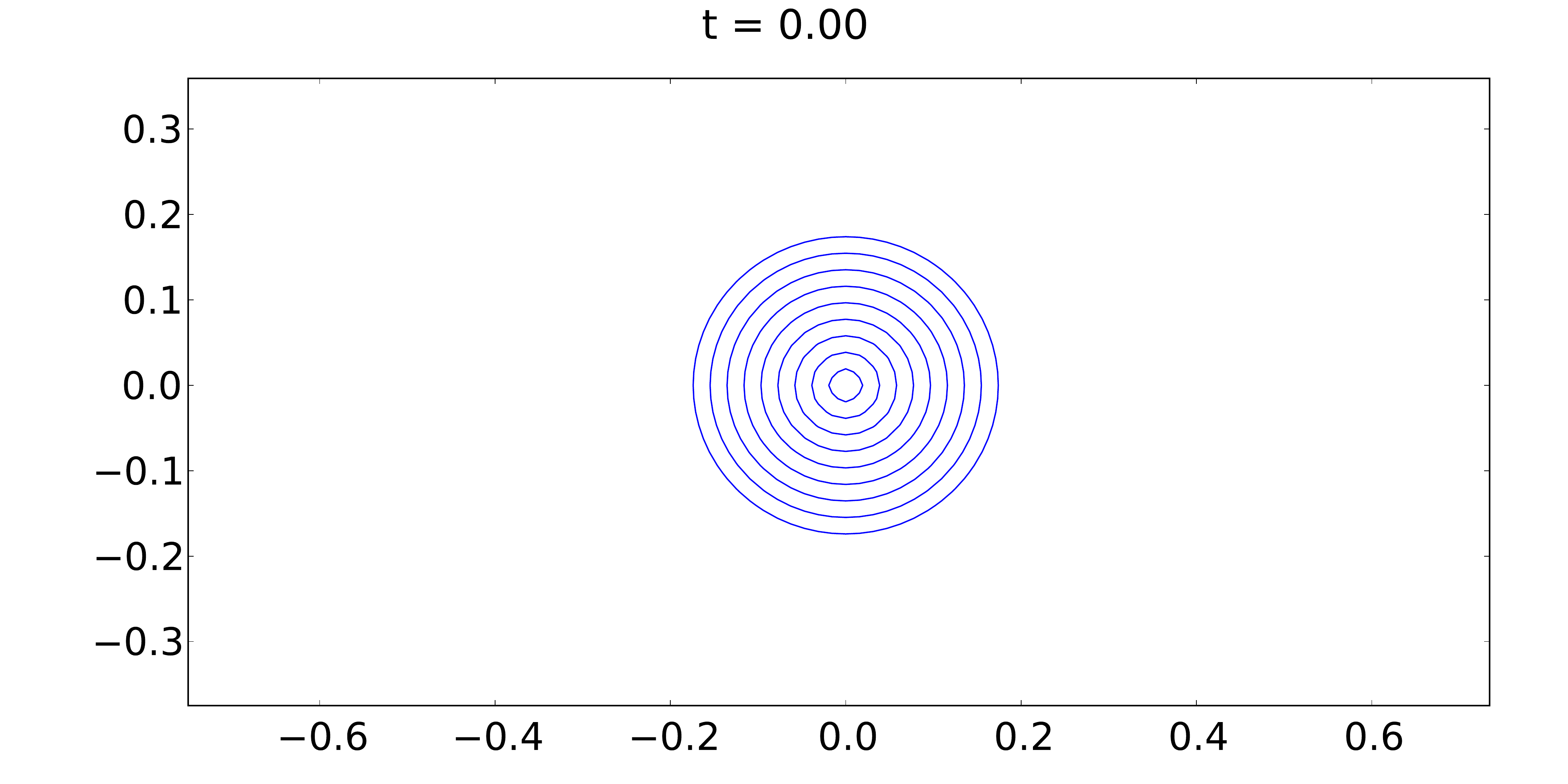}
\includegraphics[width=.48\textwidth]{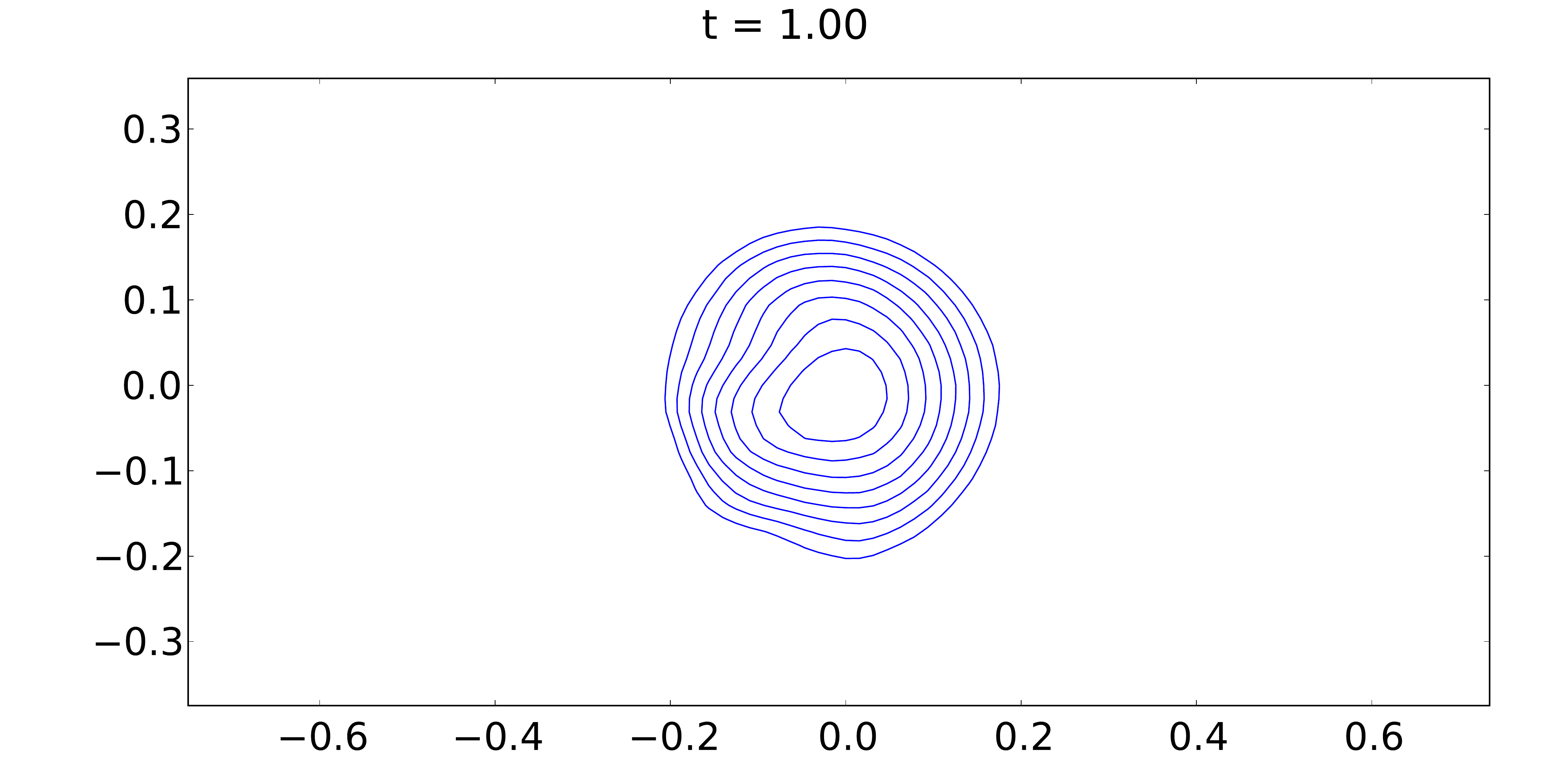}

\includegraphics[width=.48\textwidth]{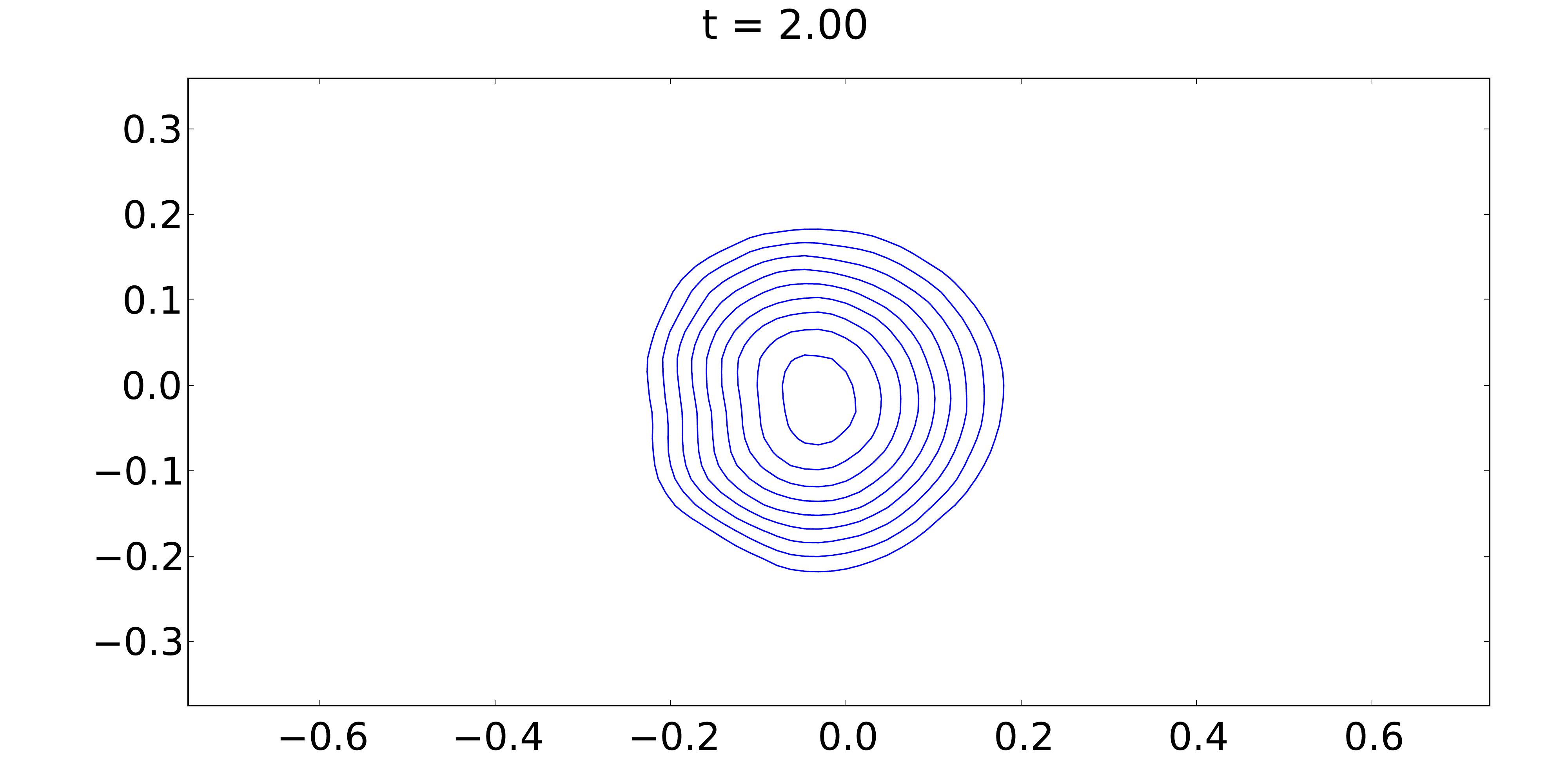}
\includegraphics[width=.48\textwidth]{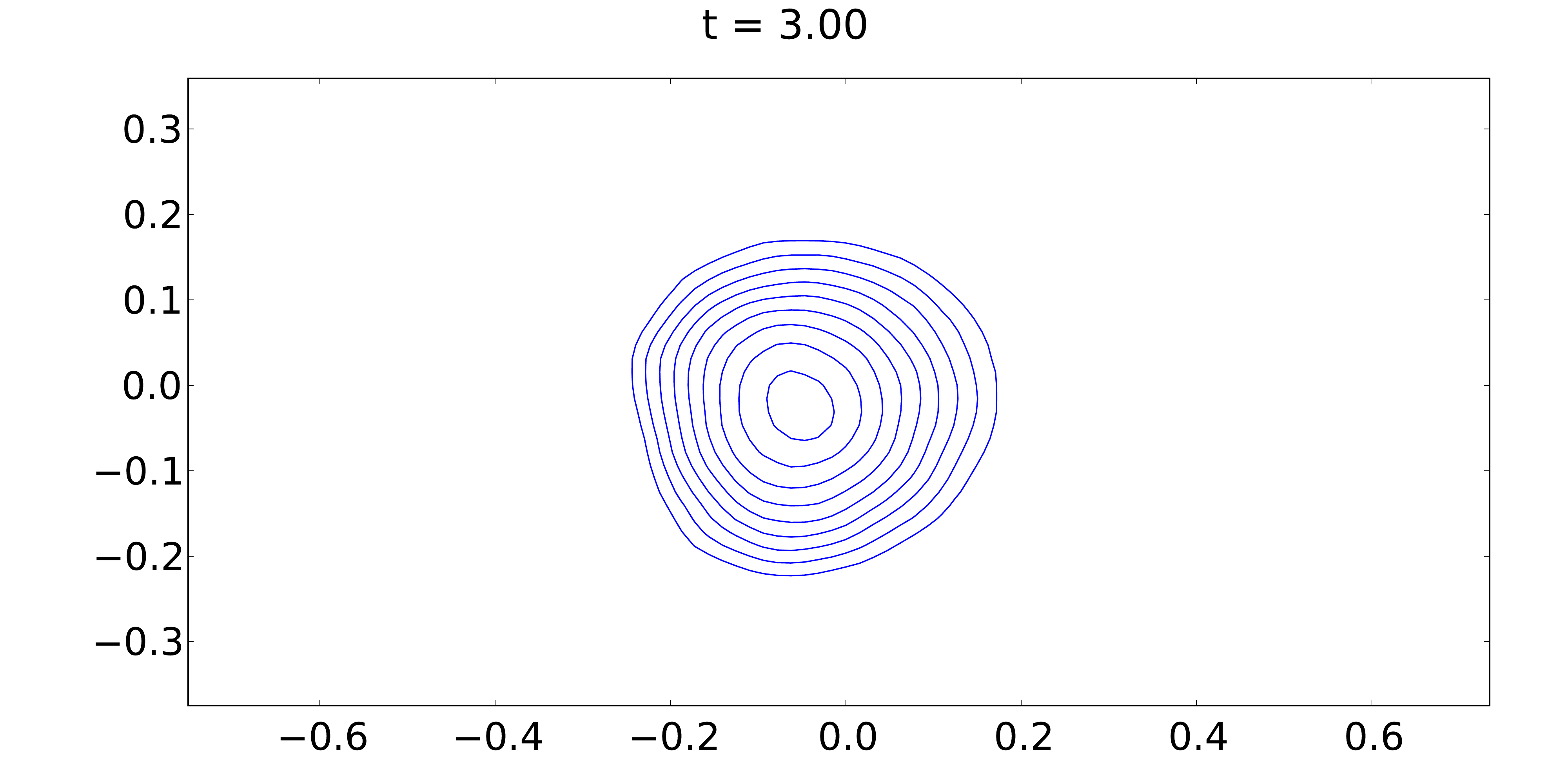}

\includegraphics[width=.48\textwidth]{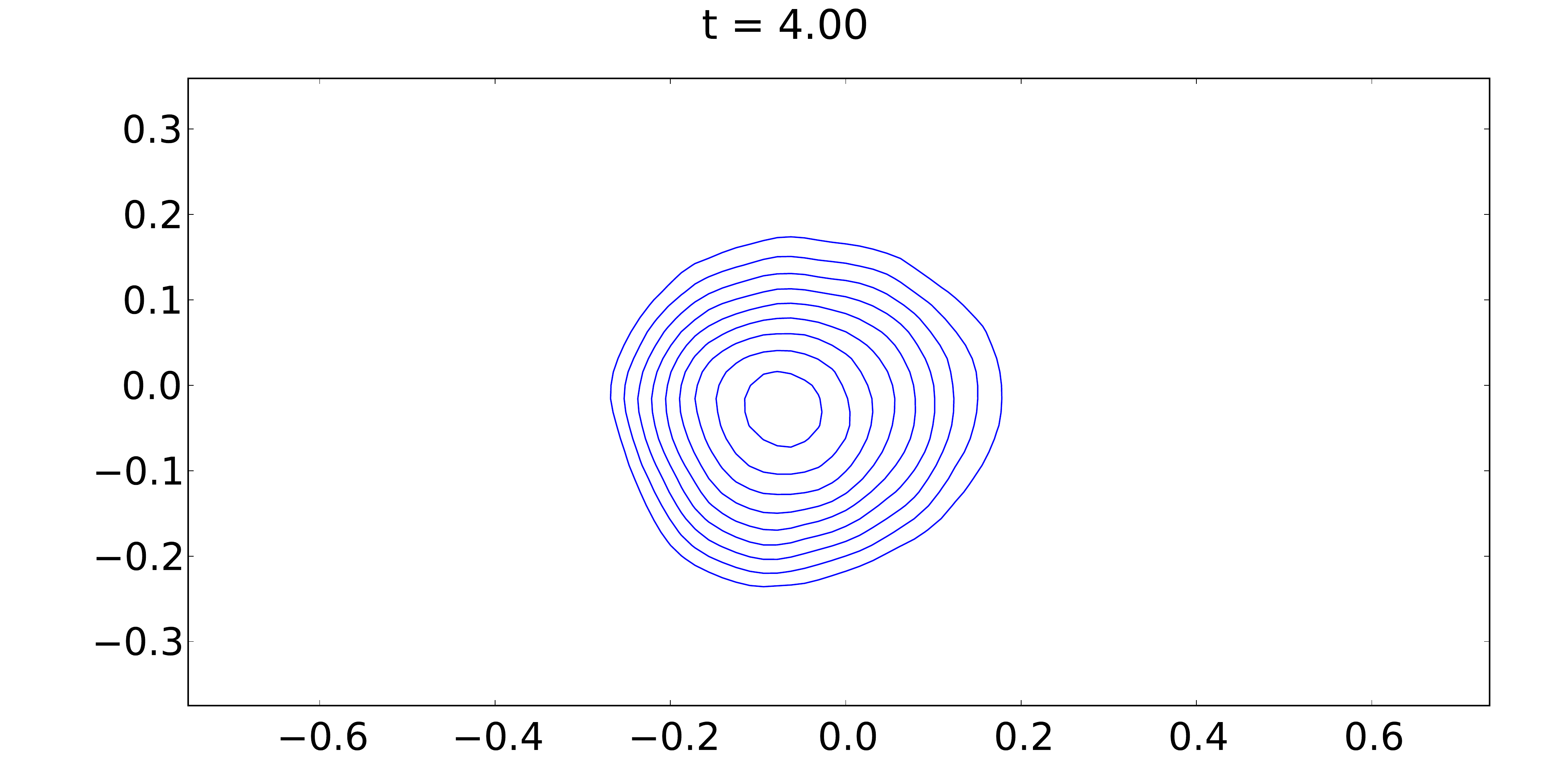}
\includegraphics[width=.48\textwidth]{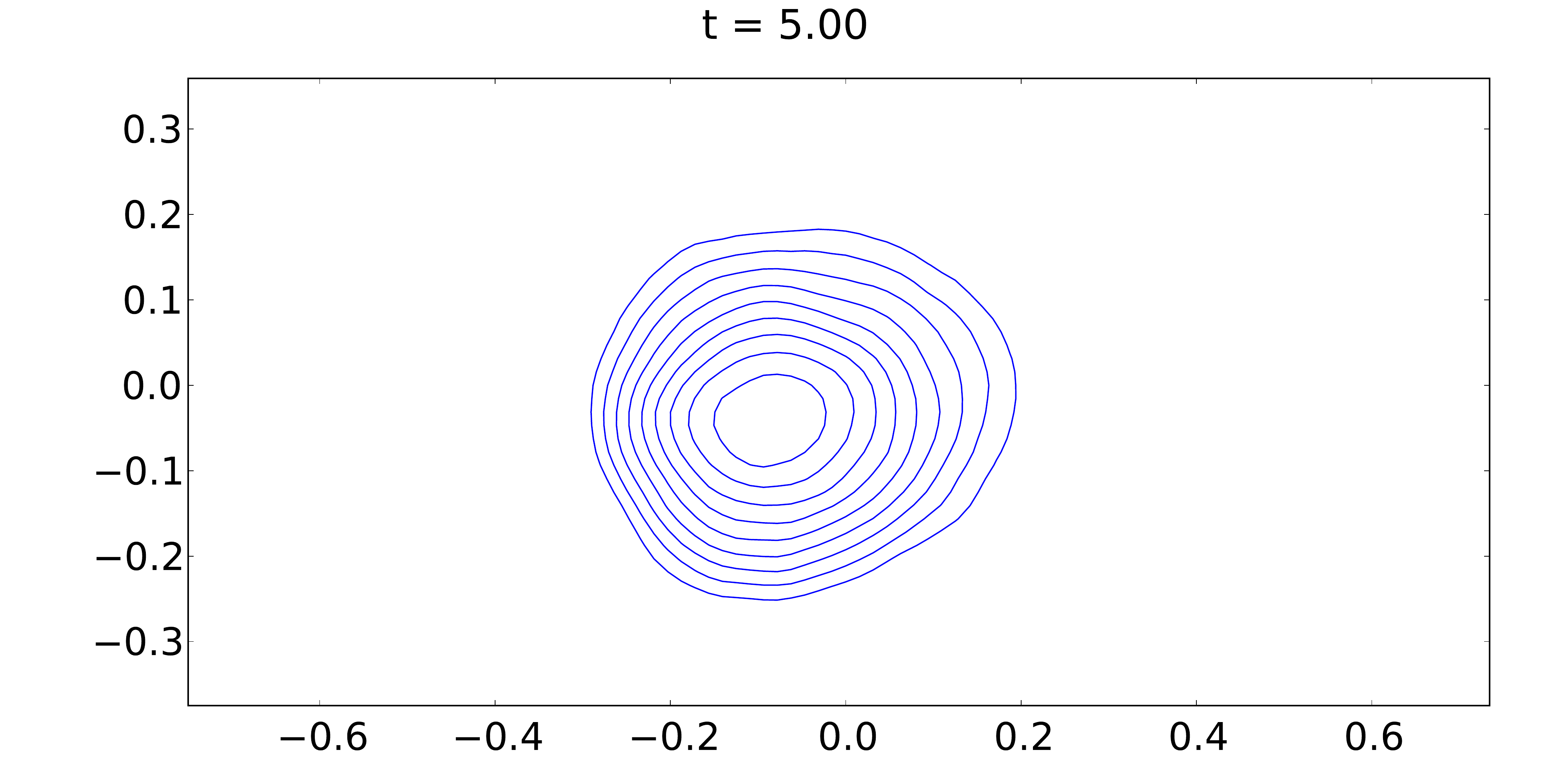}

\includegraphics[width=.48\textwidth]{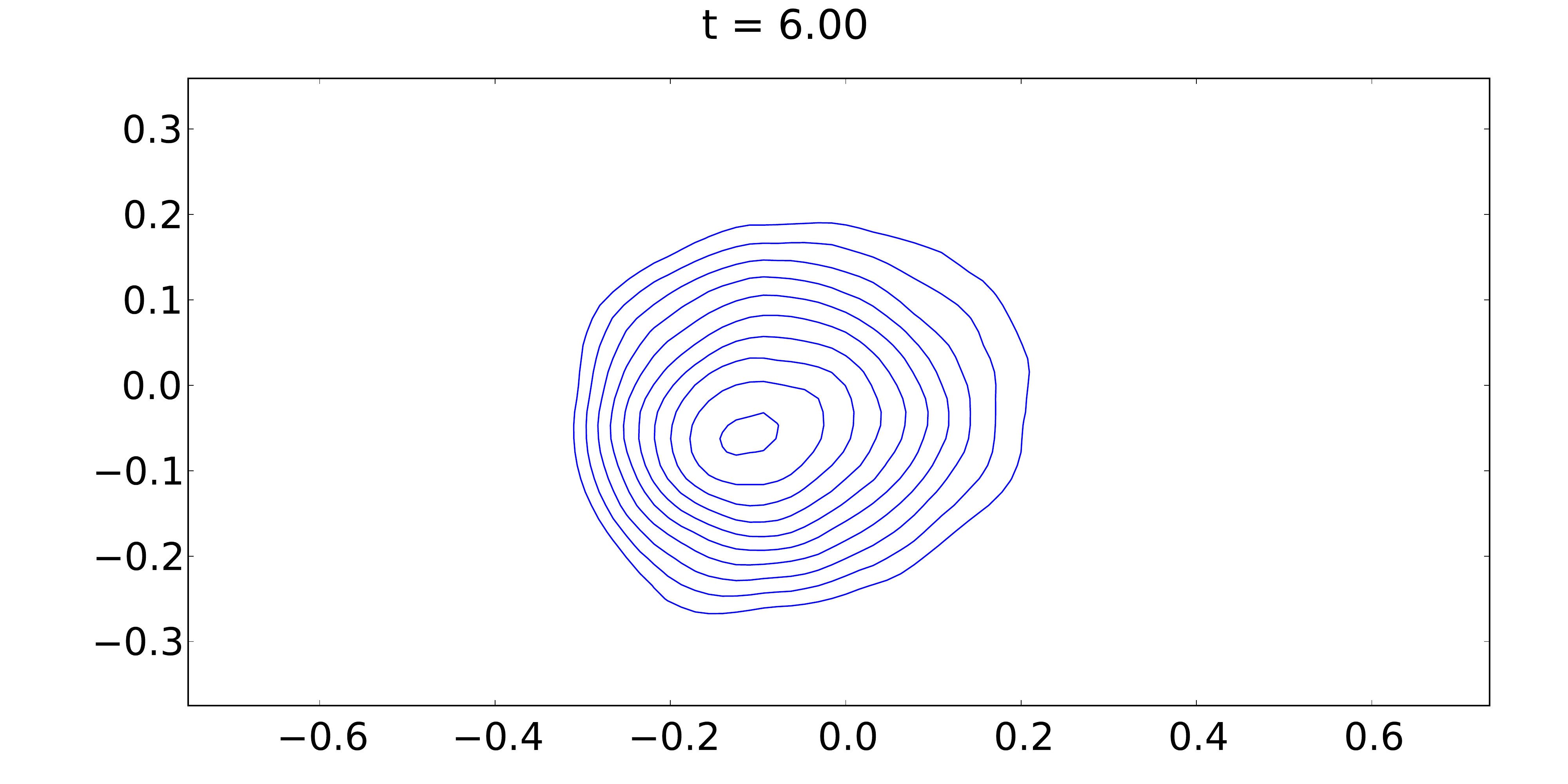}
\includegraphics[width=.48\textwidth]{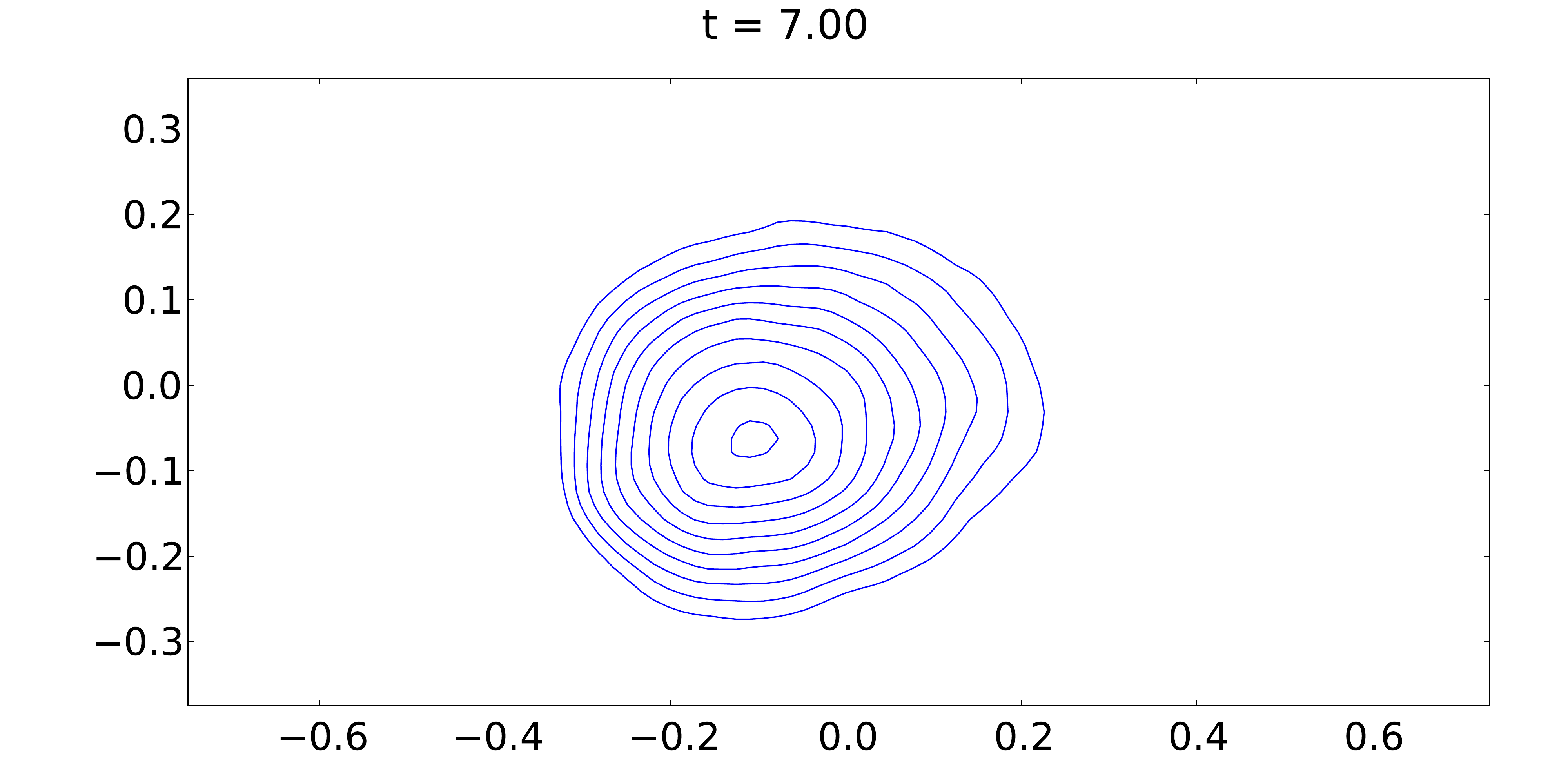}

\includegraphics[width=.48\textwidth]{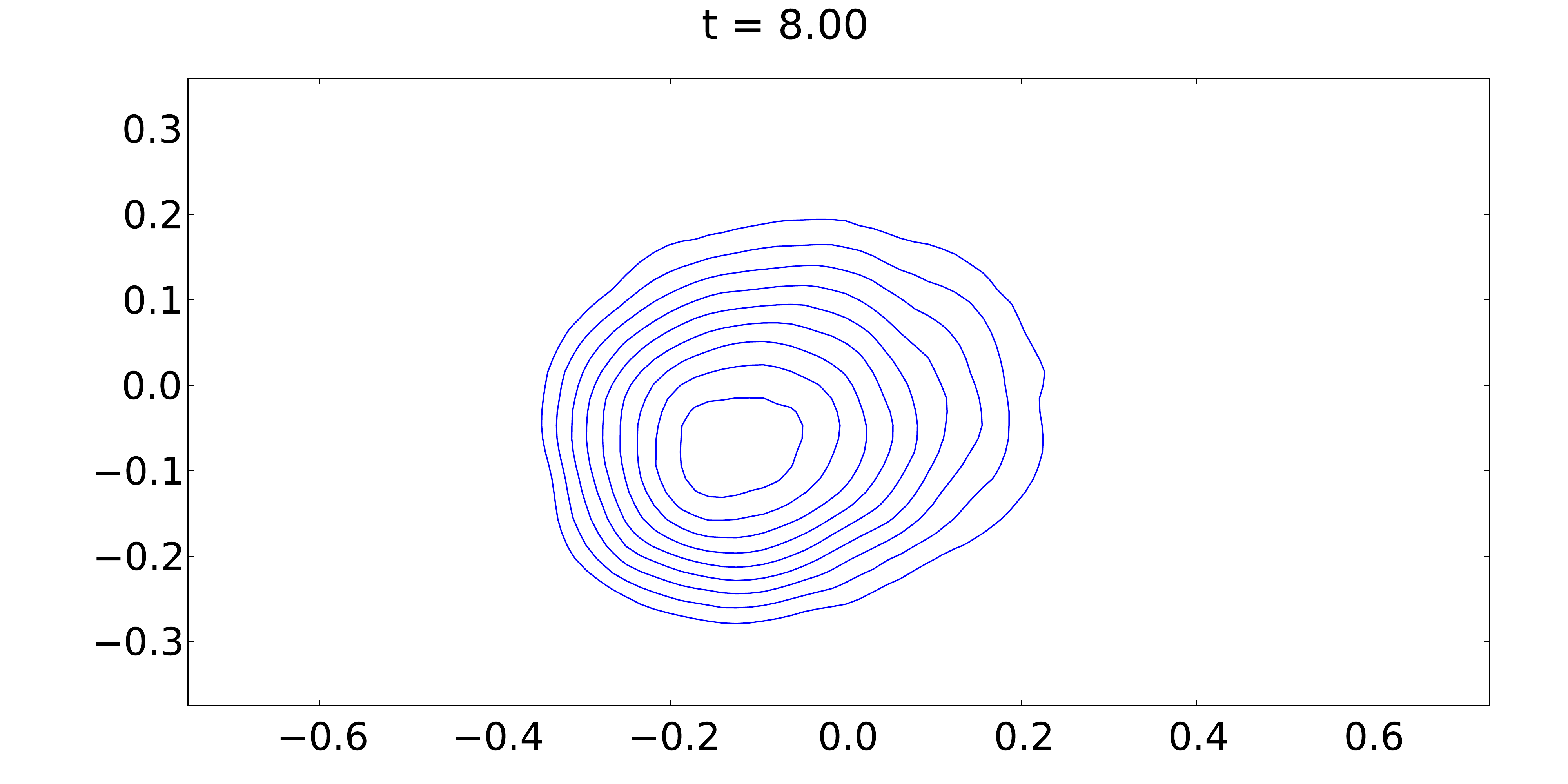}
\includegraphics[width=.48\textwidth]{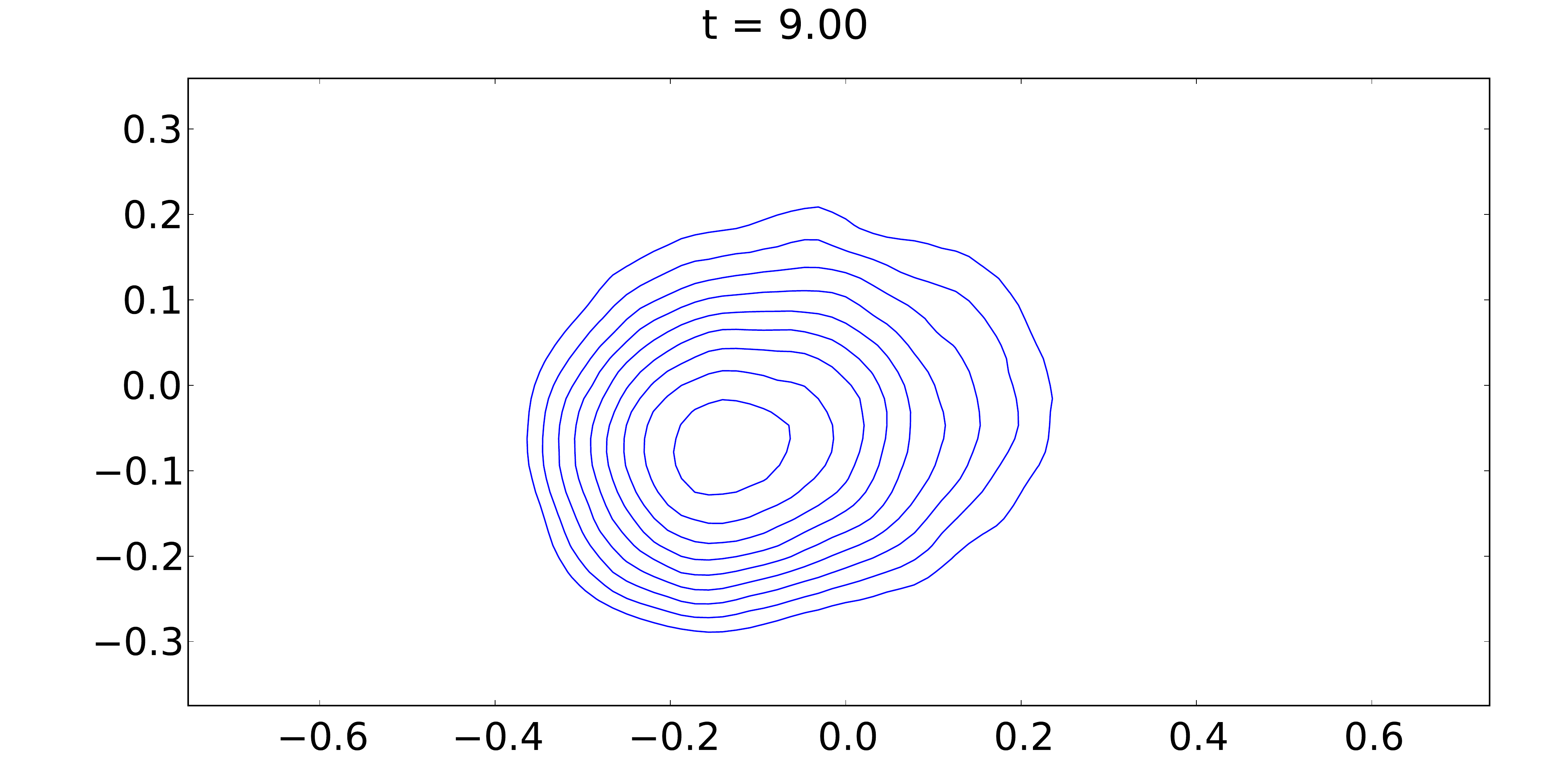}

\includegraphics[width=.48\textwidth]{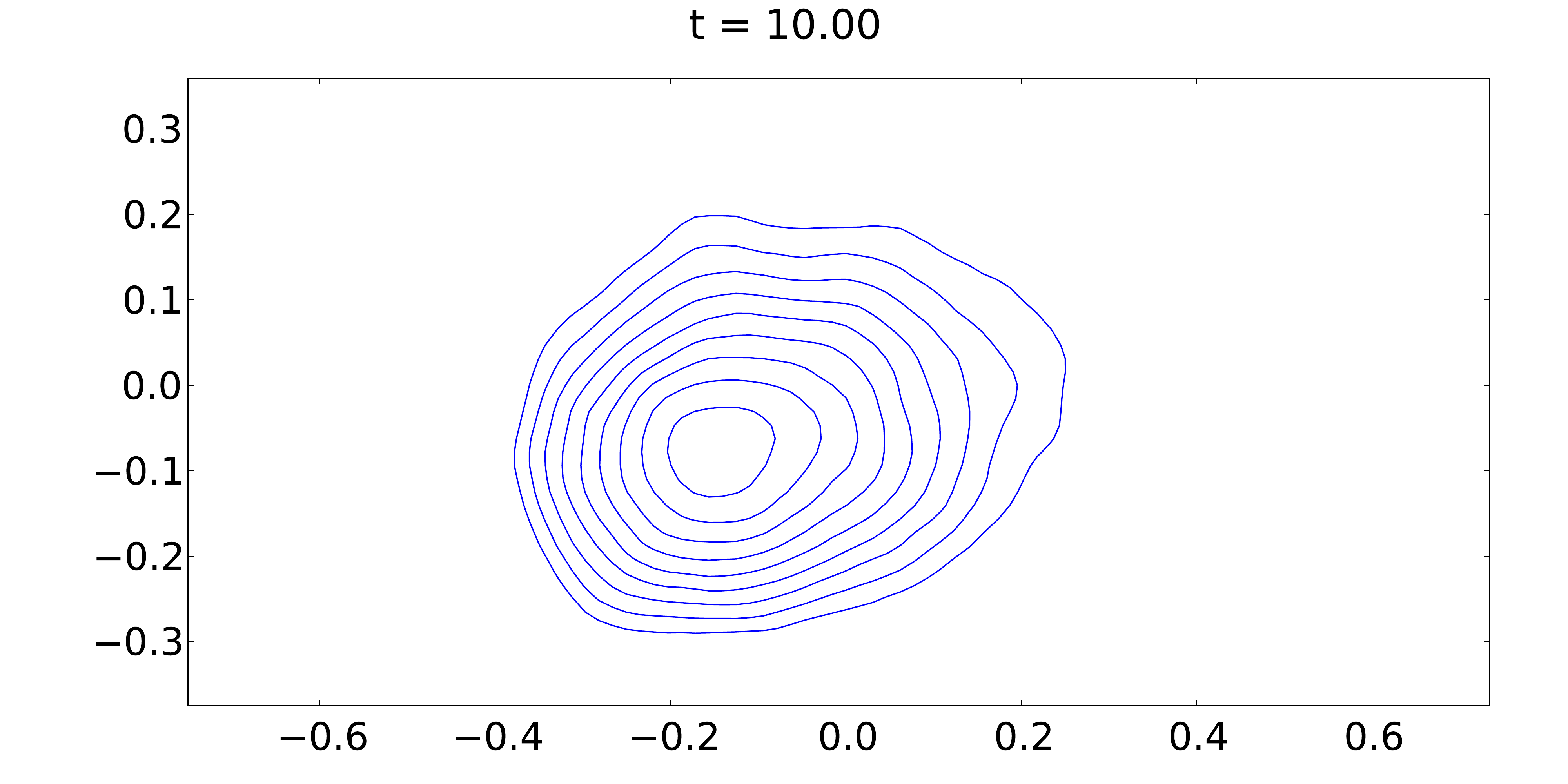}
\includegraphics[width=.48\textwidth]{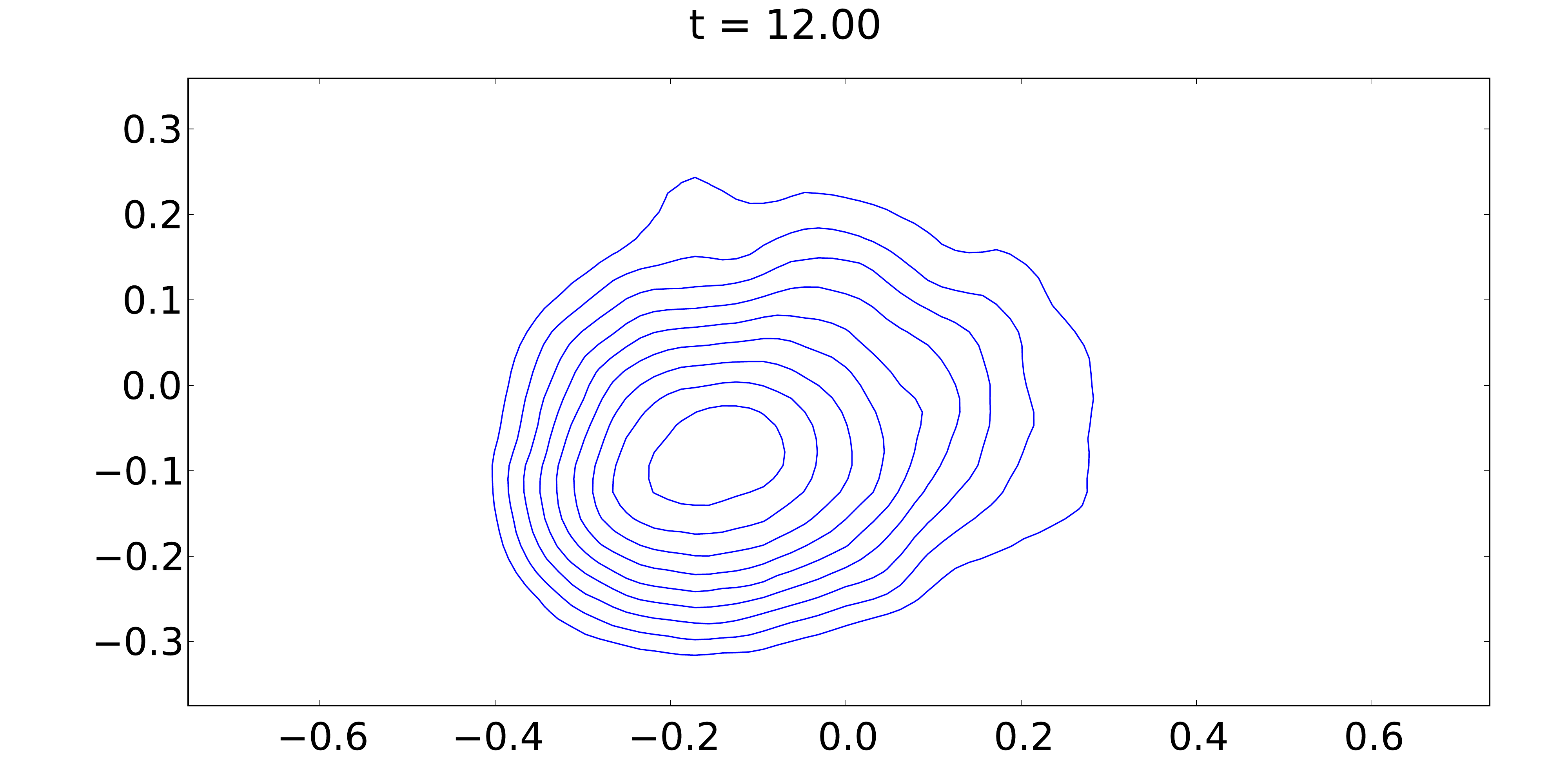}

\caption{Loop advection, $128 \times 64$ grid points. Evolution of the magnetic loop.}
\label{fig:loop_advection_64x128_field_lines}
\end{figure}

The problem is setup such that the loop should return to its initial position after integer times $t = 1, 2, 3, ...$.
In figure \ref{fig:loop_advection_64x128_field_lines} it can be seen that this is initially the case, but after some time, the loop gets deformed, such that its centre is slightly displaced from its initial position at integer times.
Note, that in figure \ref{fig:loop_advection_64x128_field_lines} only contours in the interval $[1.2 \times 10^{-4} , 3.0 \times 10^{-4}]$ are plotted. Otherwise, the plot would be too polluted since the velocity and magnetic field are of course not entirely decoupled but back react onto each other.
A possible reason for this behaviour is that our scheme is only of second order and therefore not extremely accurate. Higher order variational integrators might achieve better accuracy, but this is a topic left for future research.
Nevertheless, the energy and cross helicity are preserved optimally throughout the whole simulation (figures \ref{fig:loop_advection_64x128_errors} and \ref{fig:loop_advection_128x256_errors}). Interestingly, the errors are somewhat larger, when the resolution is increased. The amplitude of the energy error, for example, is about $5 \times 10^{-16}$ on a grid of $128 \times 64$ points, but $3 \times 10^{-15}$ on a grid of $256 \times 128$ points. First of all, one must not forget how close to machine precision these values are, so slight deviations for similar but not identical simulations are no surprise. However, this issue can be explained more specifically by the larger number of degrees of freedom in the case with higher resolution that results in a stronger error accumulation.

\subsection{Orszag-Tang Vortex}

Next we consider the evolution of current sheets in an Orszag-Tang vortex, where we use the same initial conditions as \citeauthor{CordobaMarliani:2000} \cite{CordobaMarliani:2000}, namely
\begin{align*}
V_{x} &=   \partial_{y} \psi, &
V_{y} &= - \partial_{x} \psi, &
B_{x} &=   \partial_{y} A, &
B_{y} &= - \partial_{x} A , &
P &= 0.1 . &
\end{align*}

with
\begin{align*}
\psi &= 2 \sin (y) - 2 \cos (x) , &
A &= \cos (2y) - 2 \cos (x) . & &&
\end{align*}

The spatial domain is $[0,2 \pi] \times [0,2 \pi]$ with periodic boundaries. We consider two resolutions, $64 \times 64$ and $128 \times 128$, and a timestep $h_{t} = 0.01$.

\begin{figure}
\centering
\includegraphics[width=\textwidth]{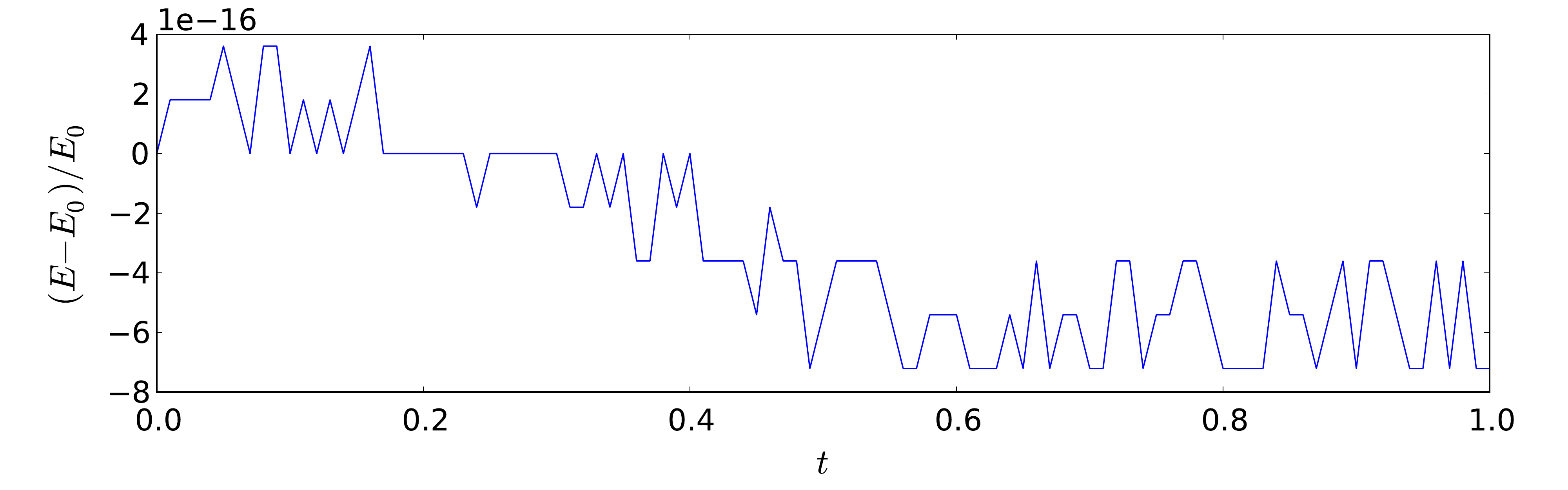}
\includegraphics[width=\textwidth]{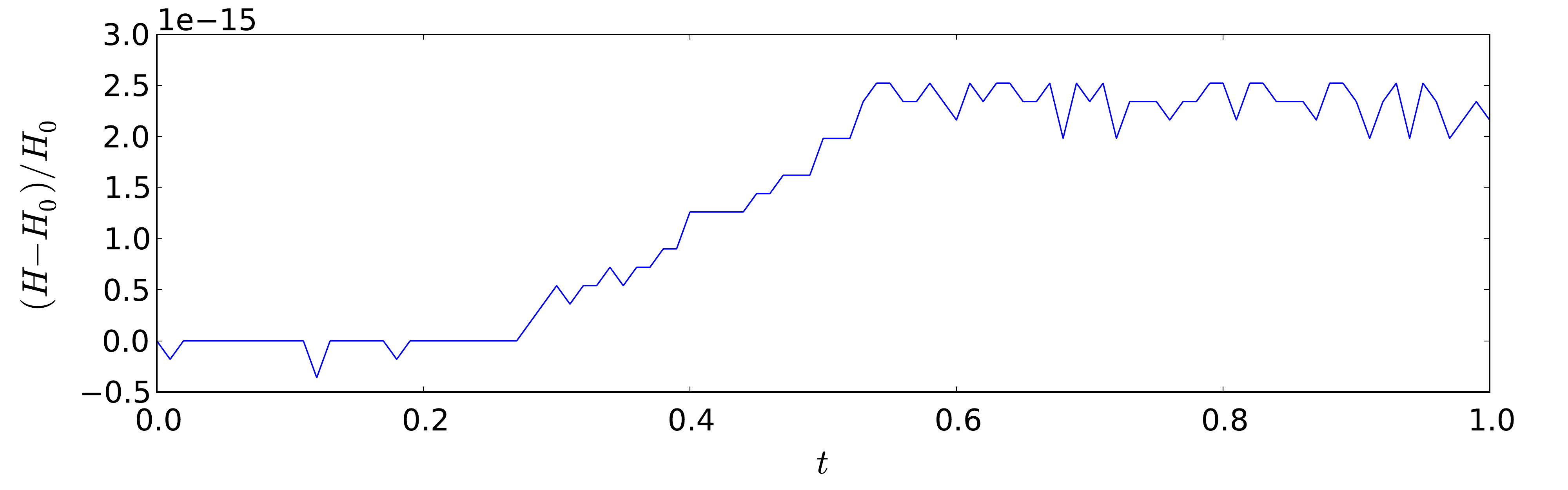}
\caption{Orszag Tang Vortex, $64 \times 64$ grid points. Conservation of energy and cross helicity.}
\label{fig:orszag_tang_64x64_errors}
\end{figure}

\begin{figure}
\centering
\includegraphics[width=\textwidth]{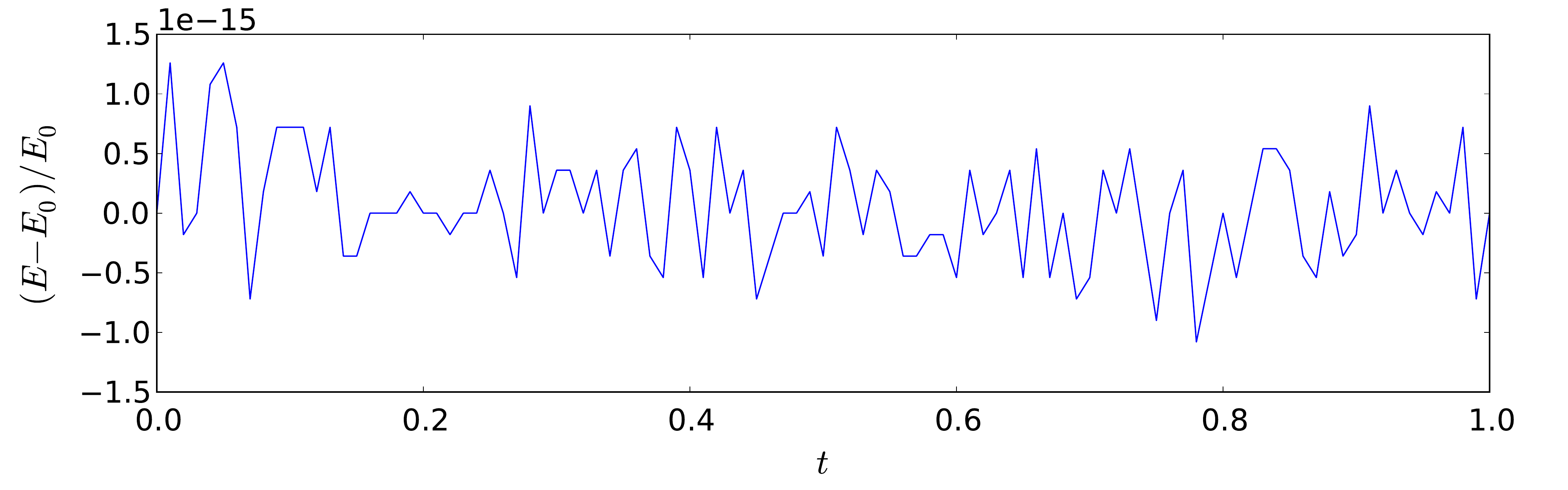}
\includegraphics[width=\textwidth]{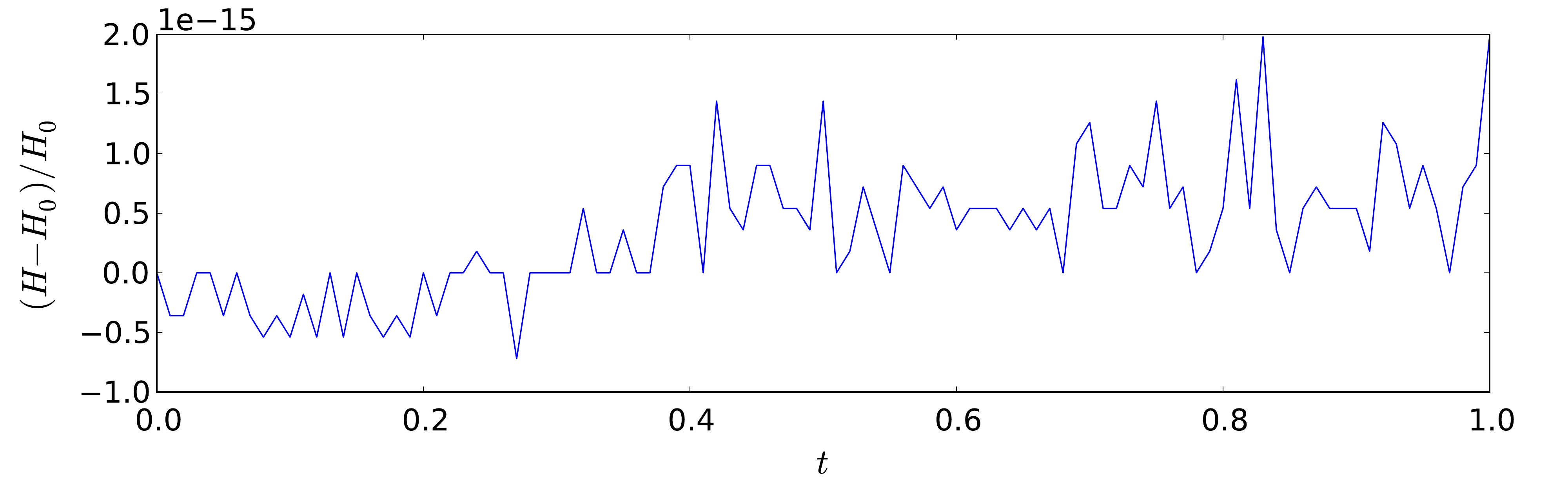}
\caption{Orszag Tang Vortex, $128 \times 128$ grid points. Conservation of energy and cross helicity.}
\label{fig:orszag_tang_128x128_errors}
\end{figure}

\begin{figure}
\centering
\subfloat{
\includegraphics[width=.32\textwidth]{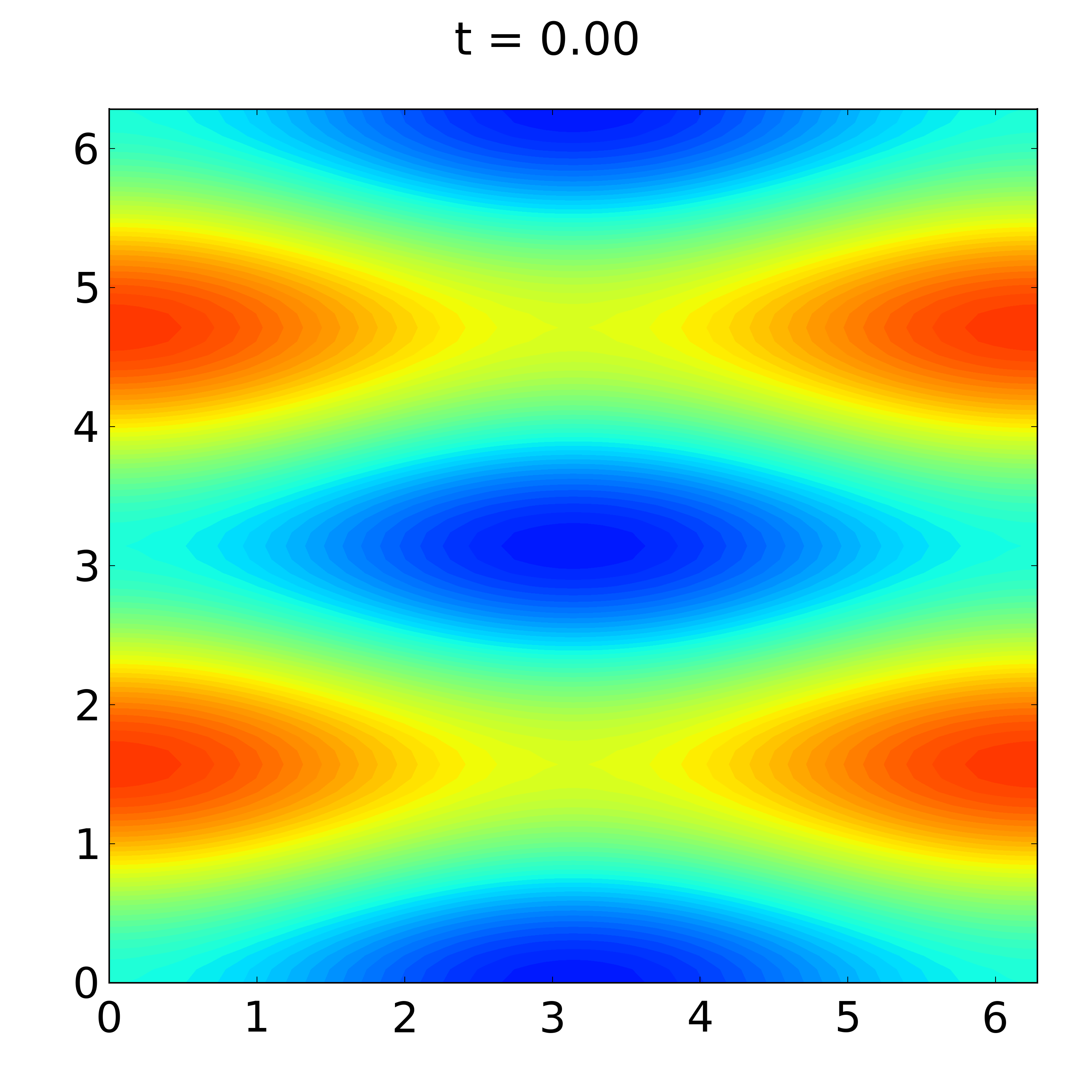}
}
\subfloat{
\includegraphics[width=.32\textwidth]{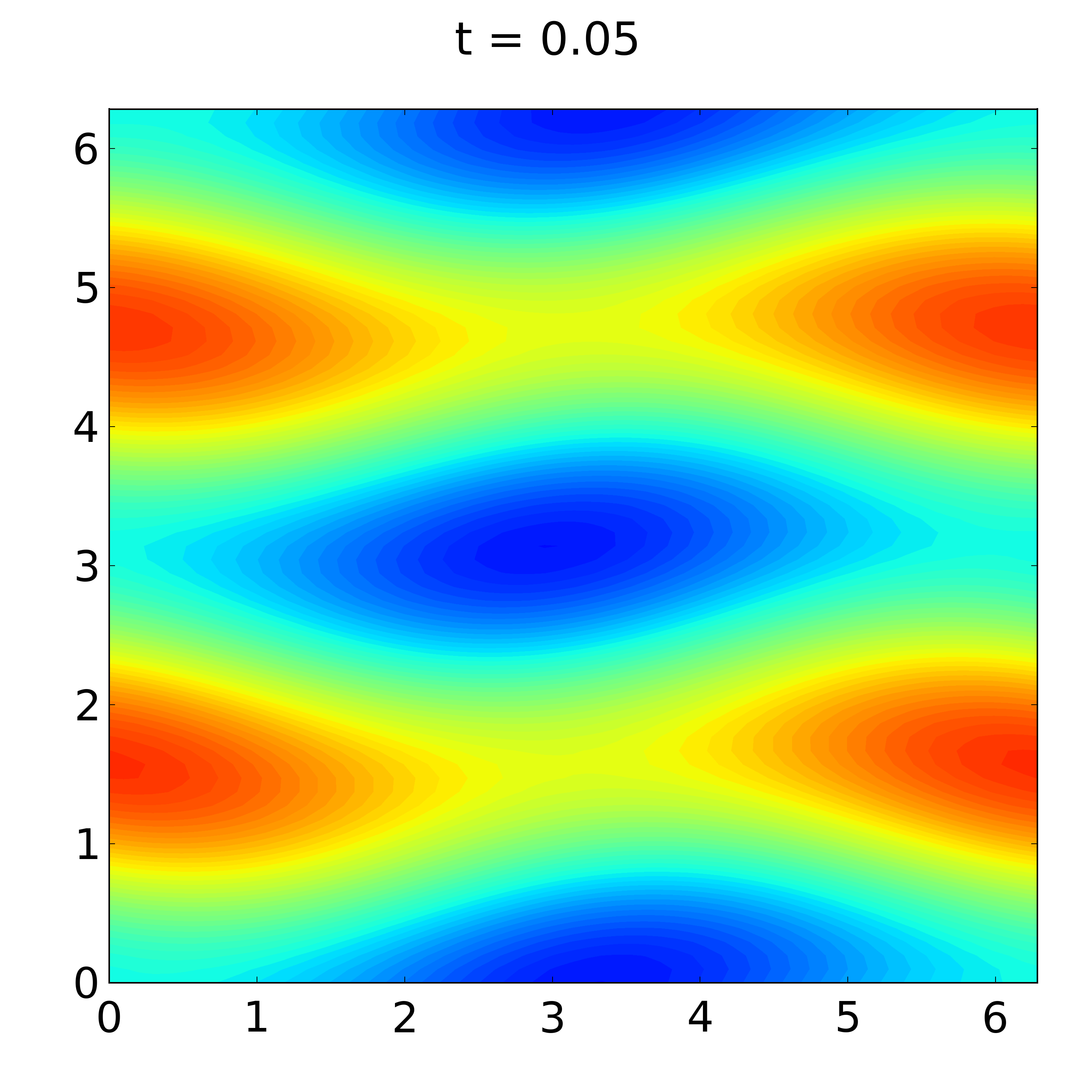}
}
\subfloat{
\includegraphics[width=.32\textwidth]{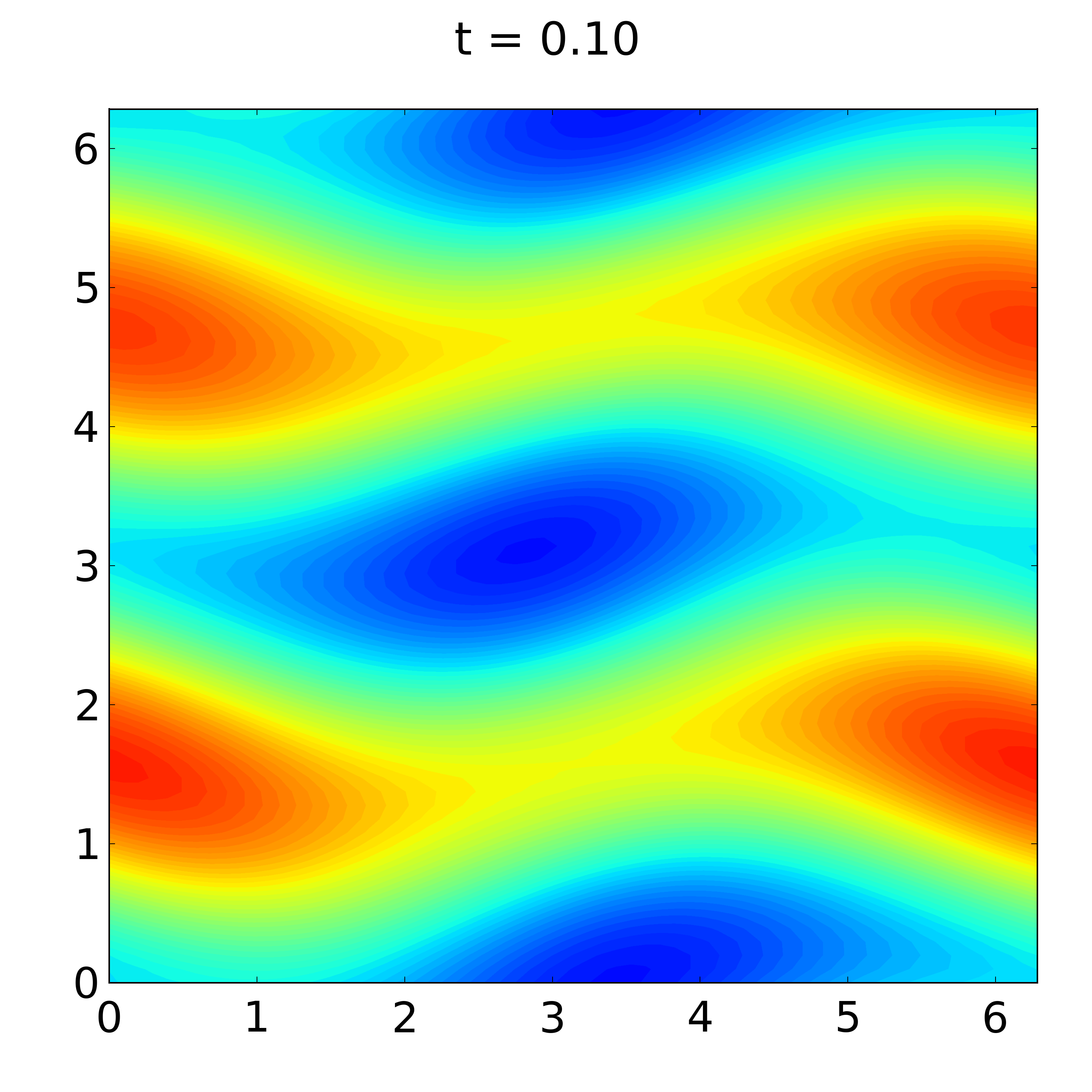}
}

\subfloat{
\includegraphics[width=.32\textwidth]{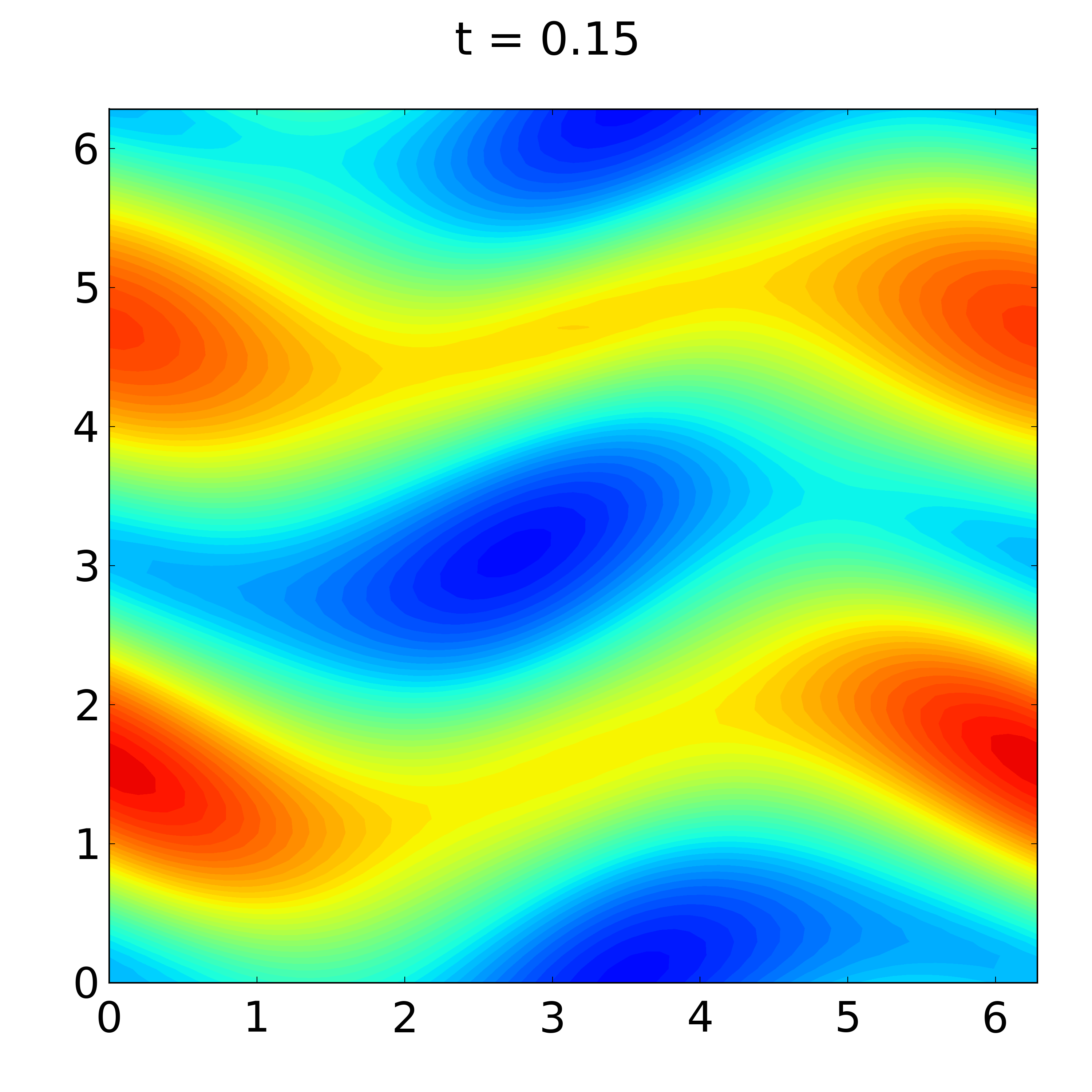}
}
\subfloat{
\includegraphics[width=.32\textwidth]{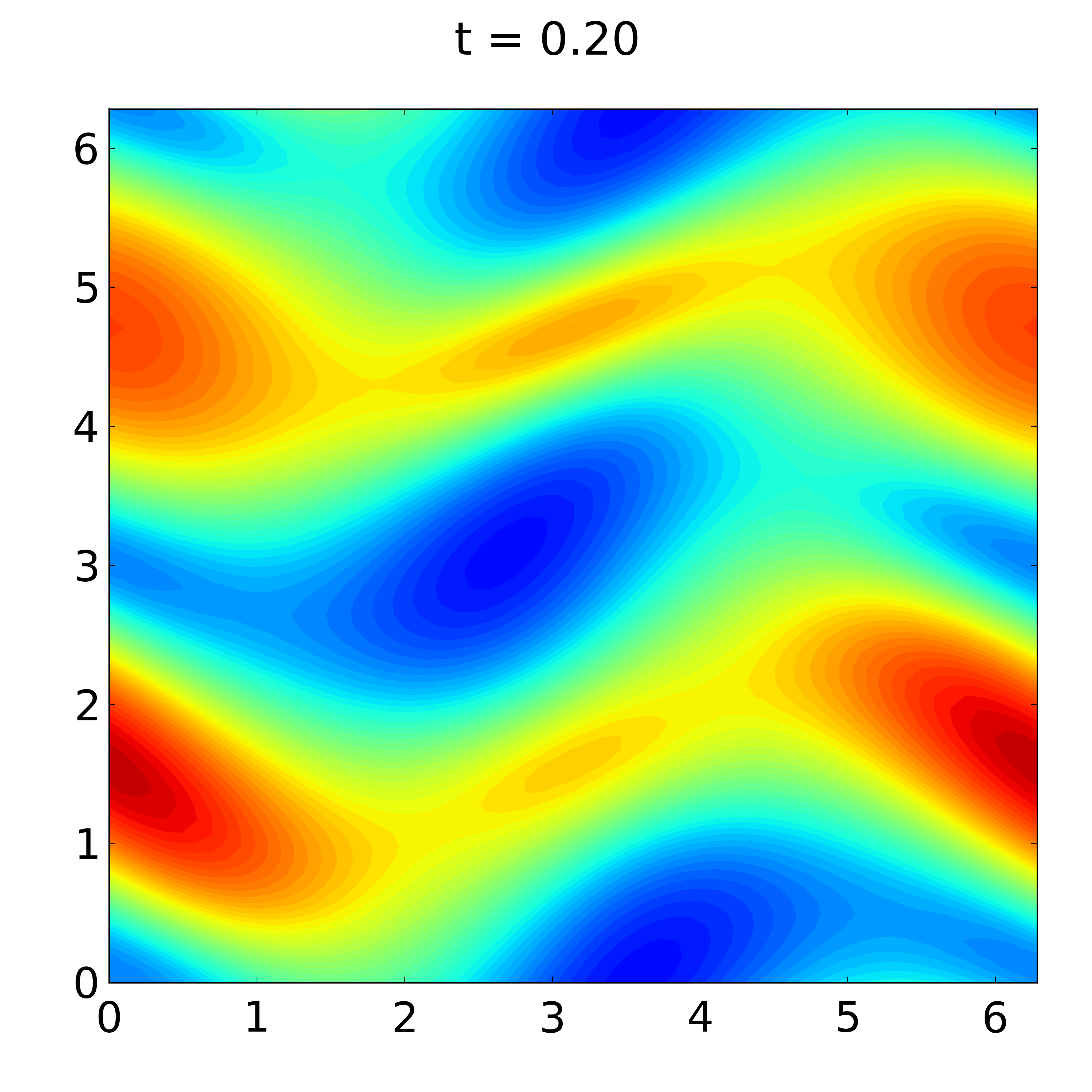}
}
\subfloat{
\includegraphics[width=.32\textwidth]{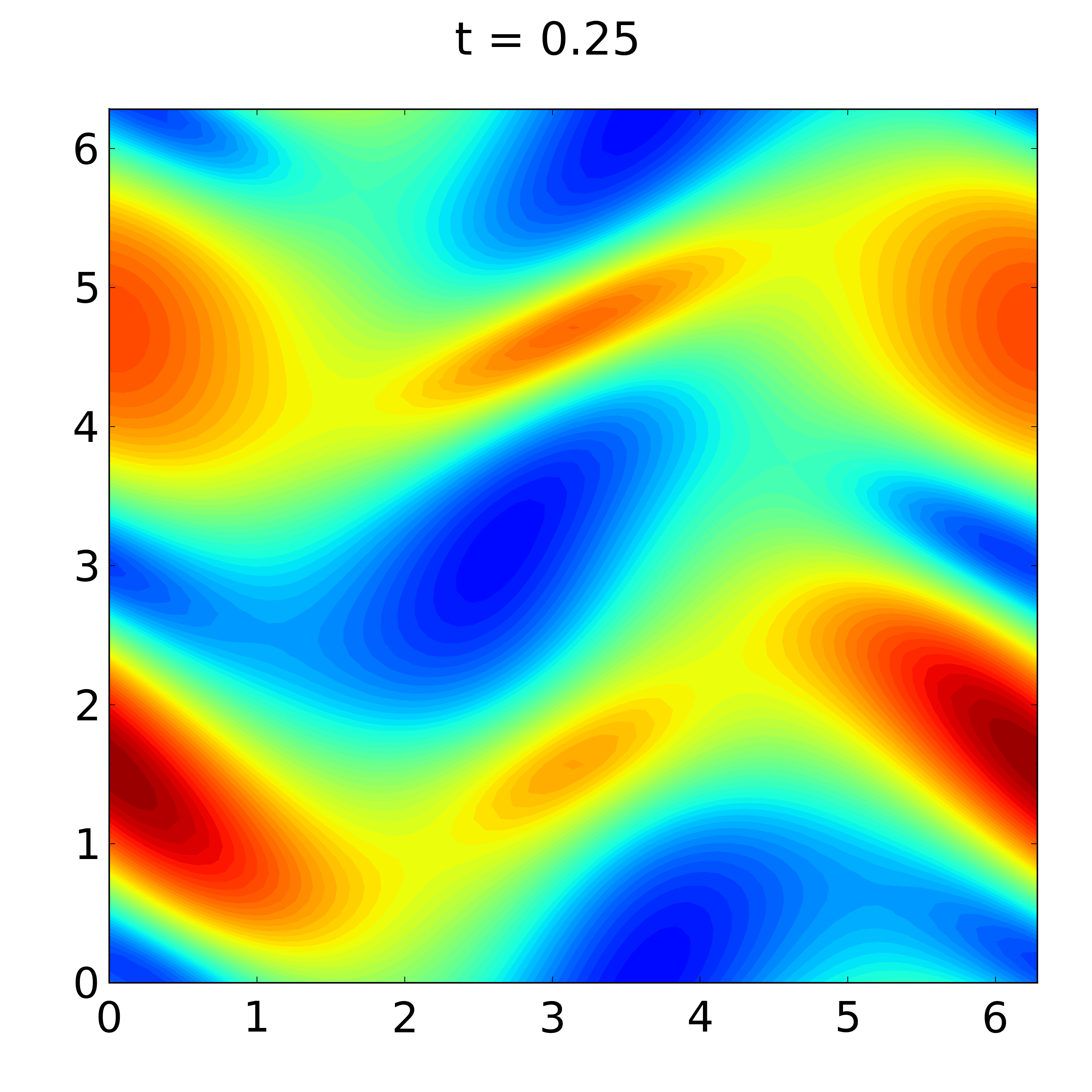}
}

\subfloat{
\includegraphics[width=.32\textwidth]{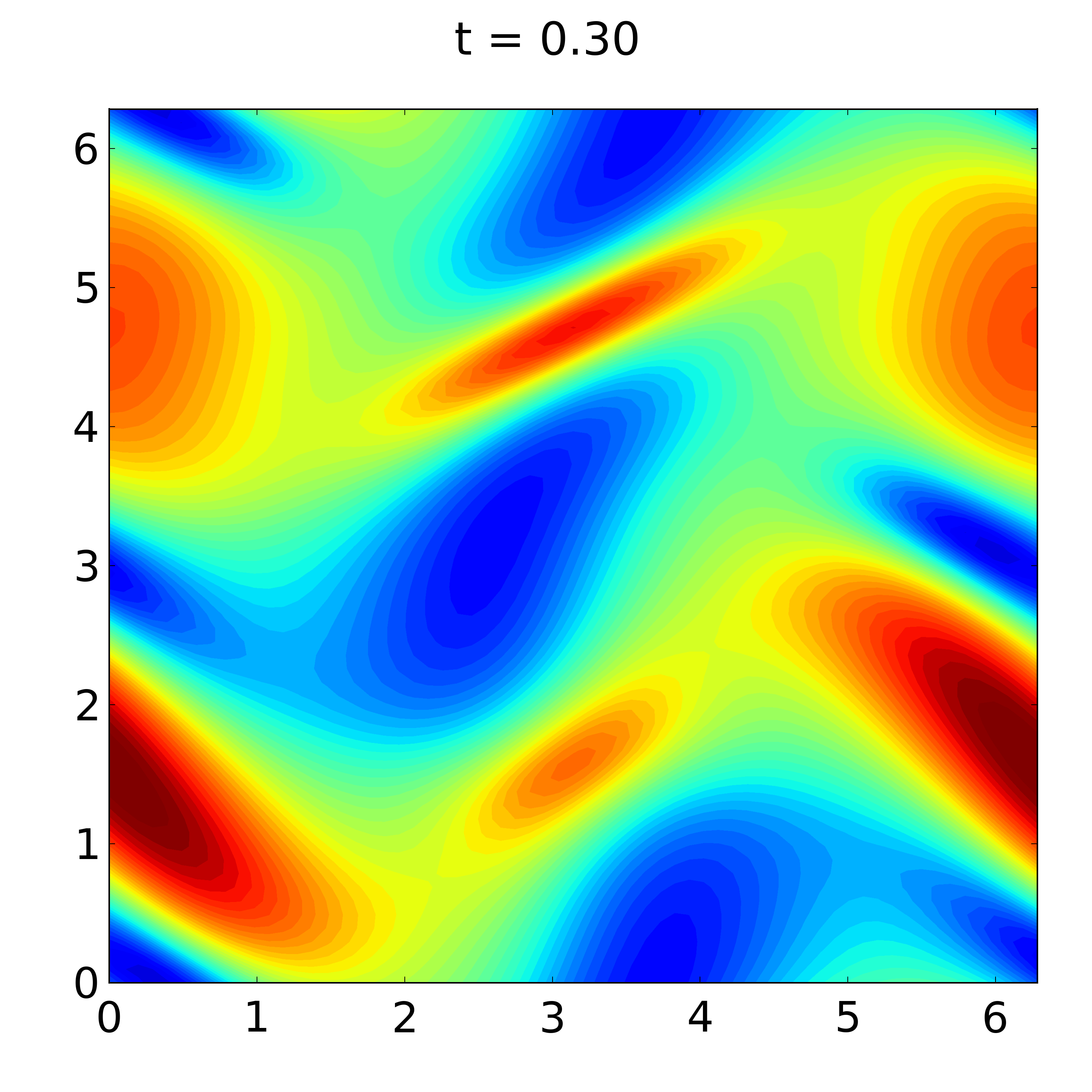}
}
\subfloat{
\includegraphics[width=.32\textwidth]{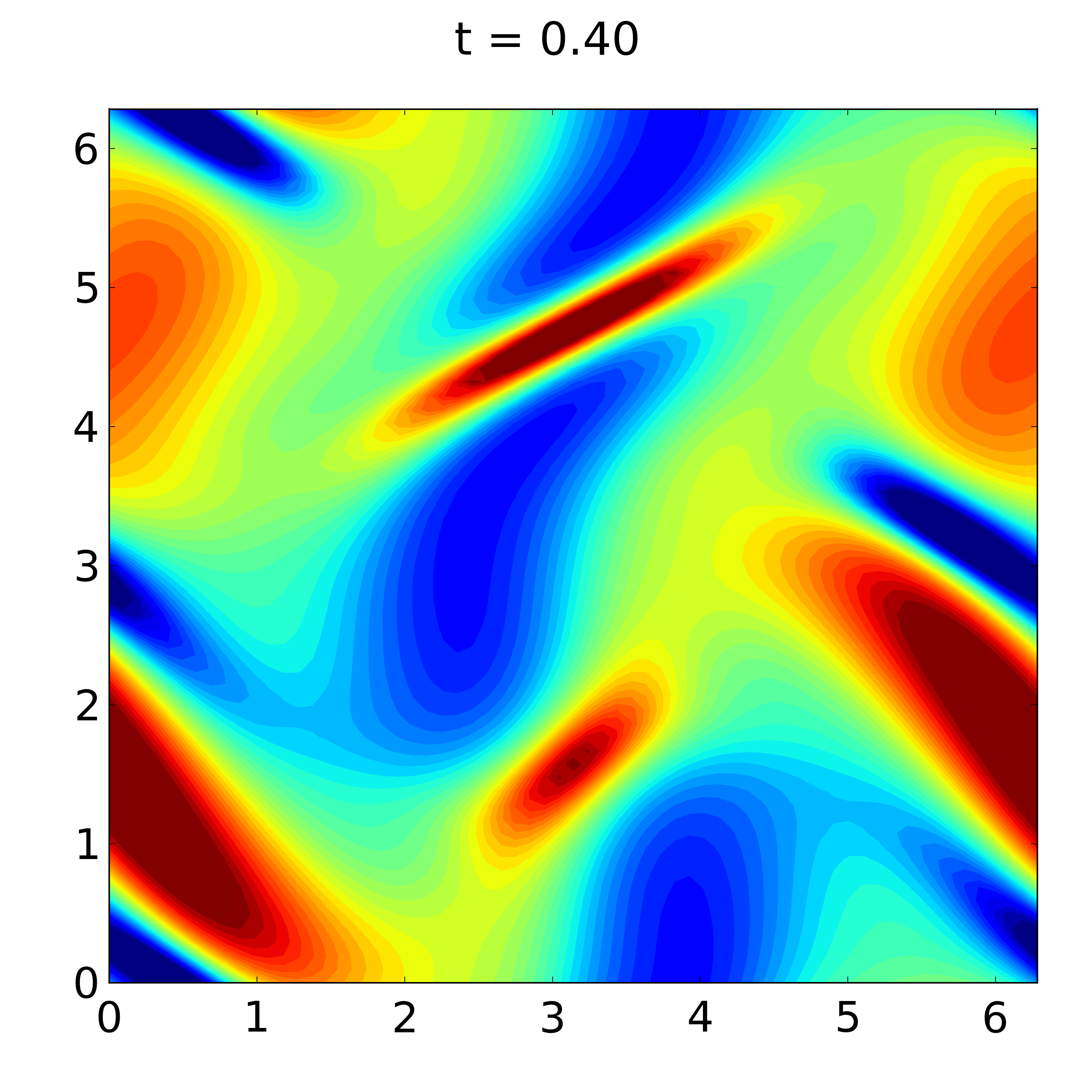}
}
\subfloat{
\includegraphics[width=.32\textwidth]{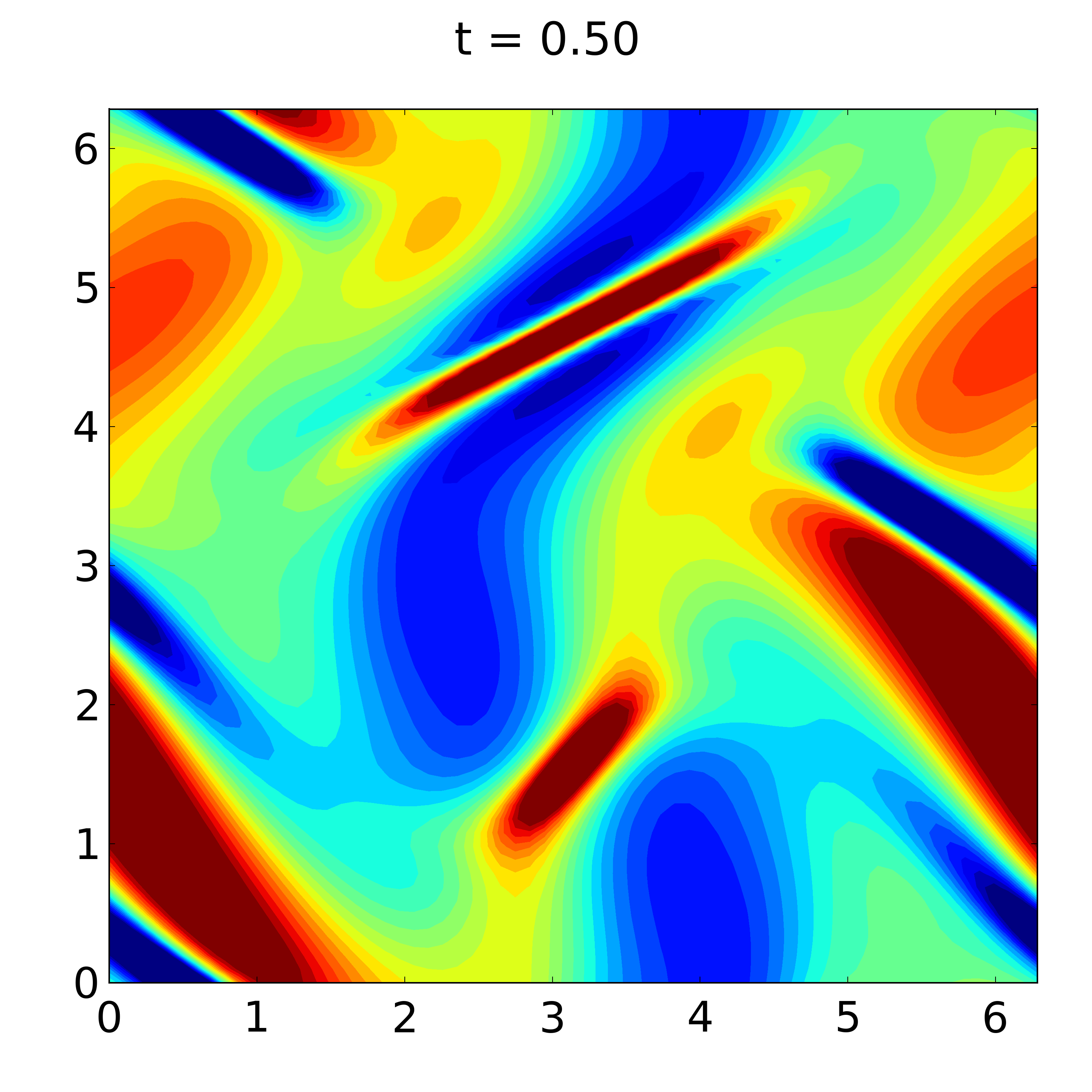}
}

\subfloat{
\includegraphics[width=.32\textwidth]{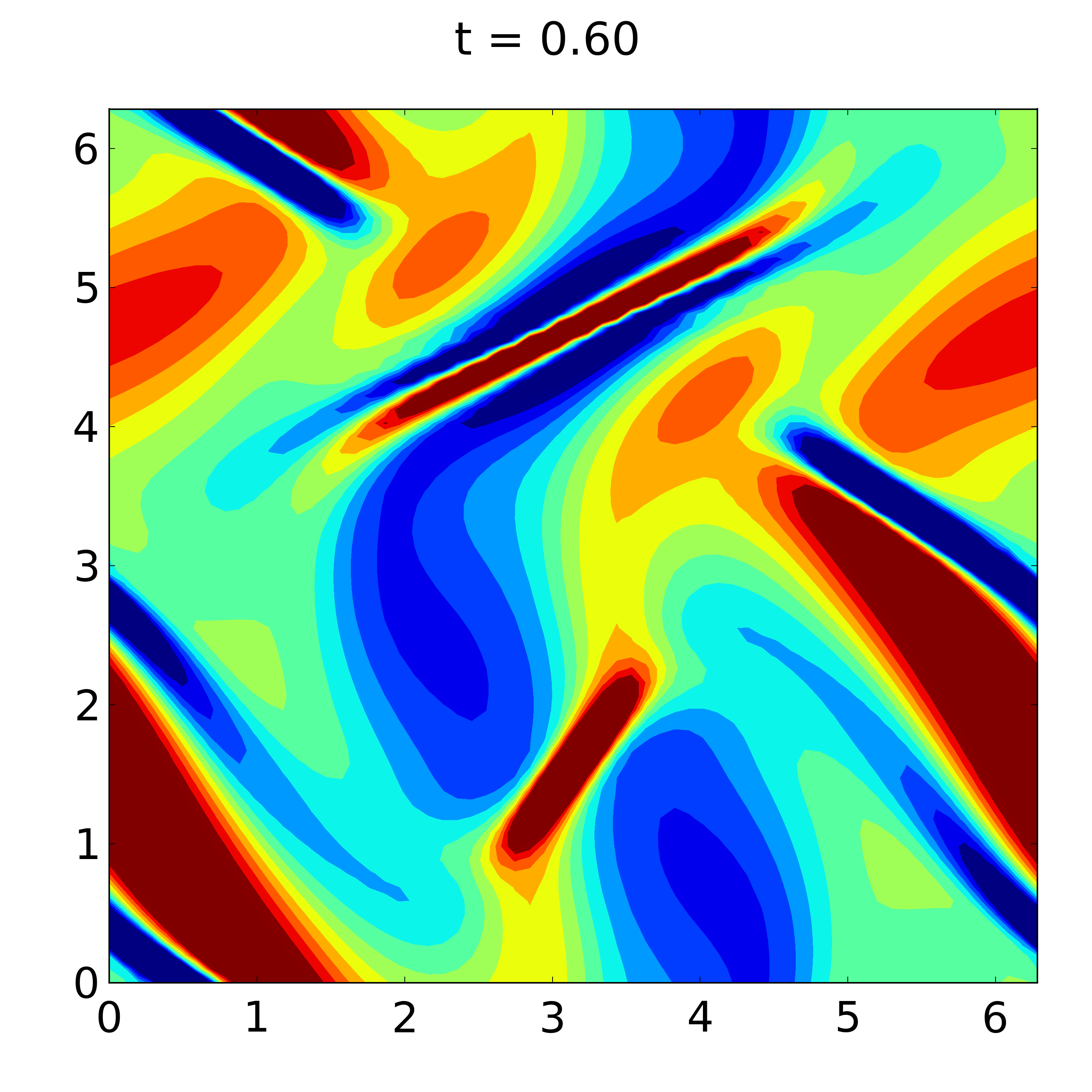}
}
\subfloat{
\includegraphics[width=.32\textwidth]{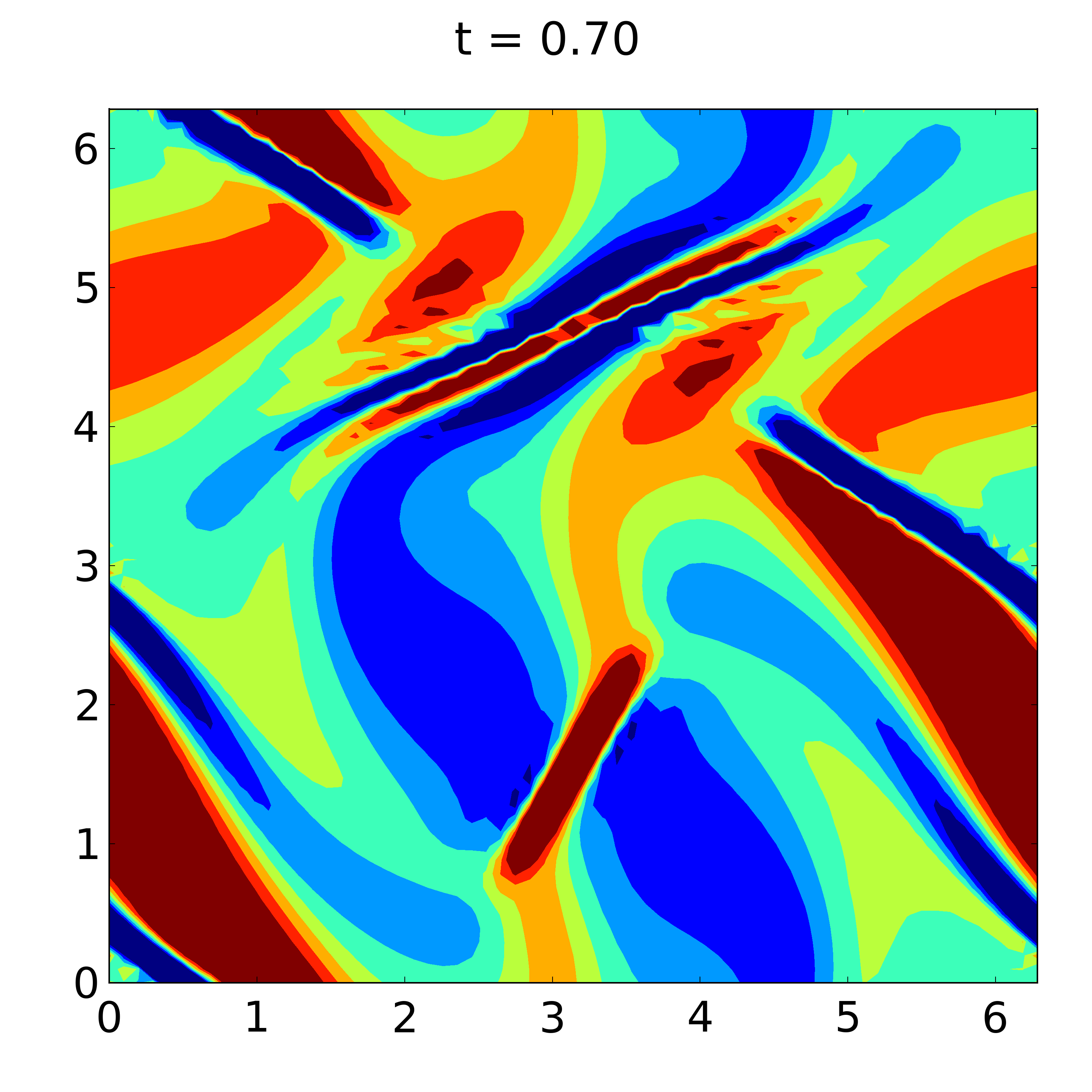}
}
\subfloat{
\includegraphics[width=.32\textwidth]{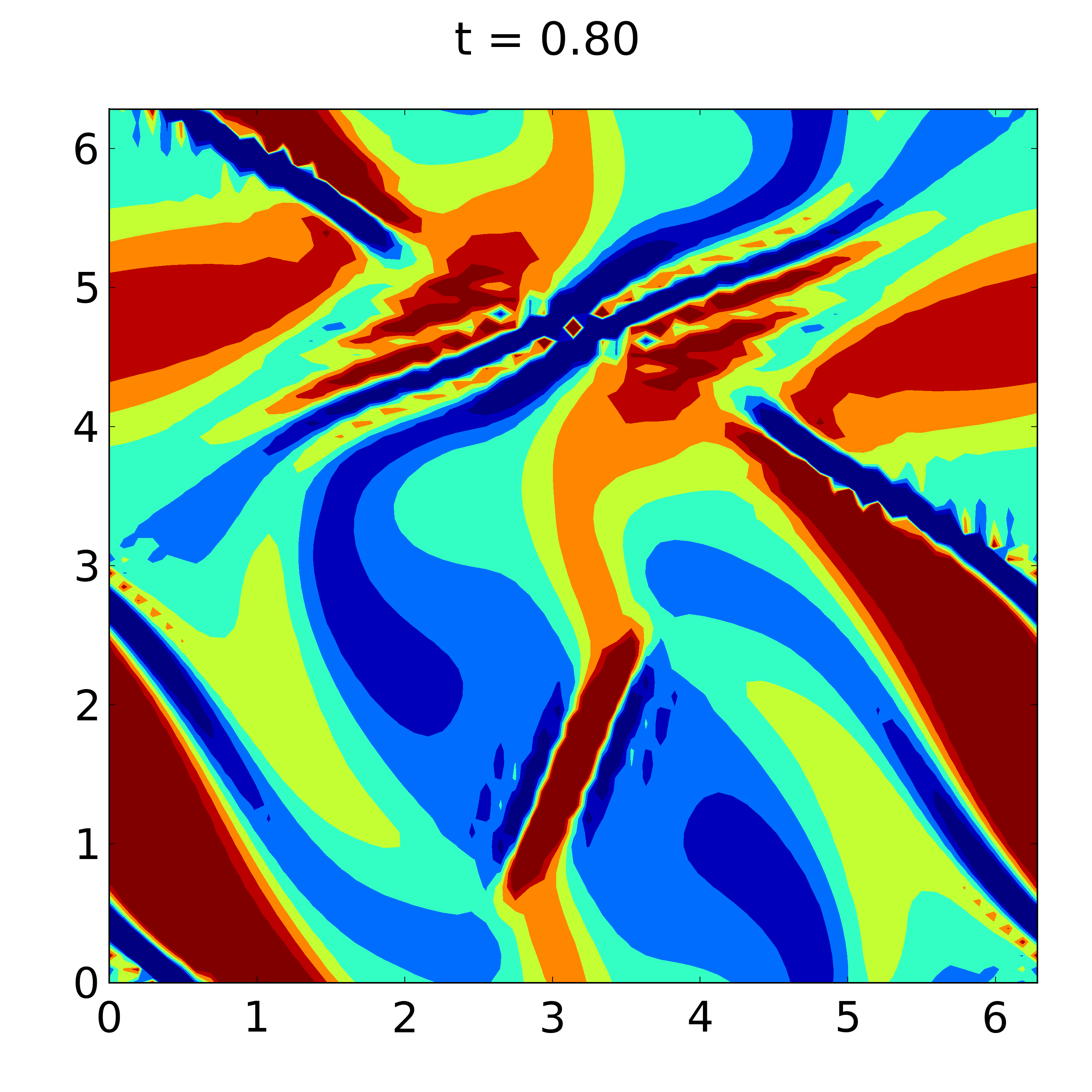}
}

\caption{Orszag Tang Vortex, $64 \times 64$ grid points. Current density. Fixed colour scale.}
\label{fig:orszag_tang_vortex_64x64_current_density}
\end{figure}

The Orszag Tang vortex constitutes a turbulent setting that leads to the growths of current sheaths. These are areas of large current density due to a change of sign in the magnetic field.
In figure \ref{fig:orszag_tang_vortex_64x64_current_density}, the current density computed by (\ref{eq:mhd_diagnostics_current_density}) is plotted.
The current sheaths are located in those parts of the plot where the colour changes from blue to red within a small region.
Starting from about $t = 60$, the simulation is under resolved and subgrid modes start to appear. The situation is similar using $128 \times 128$ grid points.
Note that in the original work, \citeauthor{CordobaMarliani:2000} \cite{CordobaMarliani:2000} use an adaptive mesh refinement approach with an initial resolution of $1024 \times 1024$ points.
The important observation is, that even with low resolution energy and cross helicity are preserved optimally (see figures \ref{fig:orszag_tang_64x64_errors} and \ref{fig:orszag_tang_128x128_errors}). During the simulation, only a slight growth in the errors is observed, probably due to the subgrid modes, but the errors stay $\mcal{O} (10^{-15})$ for the whole time.
As before, we observe that the energy error is slightly larger in the case of higher resolution. As before, this is probably due to the higher number of degrees of freedom, and therefore increased error accumulation.

\subsection{Current Sheath}

In the following, we consider as initial conditions for the magnetic field three different current sheath models that appear in reconnection studies: a sharp jump of the magnetic field
\begin{align*}
B_{y}^{\text{sharp}} &=
\begin{cases}
+ B_{0} & x < x_{1} \\
- B_{0} & x_{1} \leq x \leq x_{2} \\
+ B_{0} & x > x_{2}
\end{cases} &
& \text{with} &
x_{1} &= 0.5 , &
x_{2} &= 1.5 , &
\end{align*}

a $\tanh$ profile
\begin{align*}
B_{y}^{\tanh} = B_{0} \, \tanh (\pi x) ,
\end{align*}

and a $\cosh$ profile
\begin{align*}
B_{y}^{\cosh} = \dfrac{B_{0}}{\cosh (\pi x)} ,
\end{align*}

with $B_{x} = 0$ and $B_{0} = 1$ in all three cases.
The initial conditions for the fluid quantities are
\begin{align*}
V_{x} &= V_{0} \, \sin (\pi y), &
V_{y} &= 0, &
P &= 0.1  &
\end{align*}

for all three cases.
The spatial domains are $[0,2] \times [0,2]$ for the sharp jump, $[-4 , +4] \times [0 , 4]$ for the $\tanh$ case and $[-1 , +1] \times [-1 , +1]$ for the $\cosh$ case. For the sharp case we consider different resolutions, namely $32 \times 32$ and $64 \times 64$ grid points, to analyse the qualitative and conservative behaviour. To compare the field line evolution for the different models, we use a common resolution of $n_{x} \times n_{y} = 30 \times 30$. In all three cases we use periodic boundaries and a timestep of $h_{t} = 0.1$.

\begin{figure}
\centering
\includegraphics[width=\textwidth]{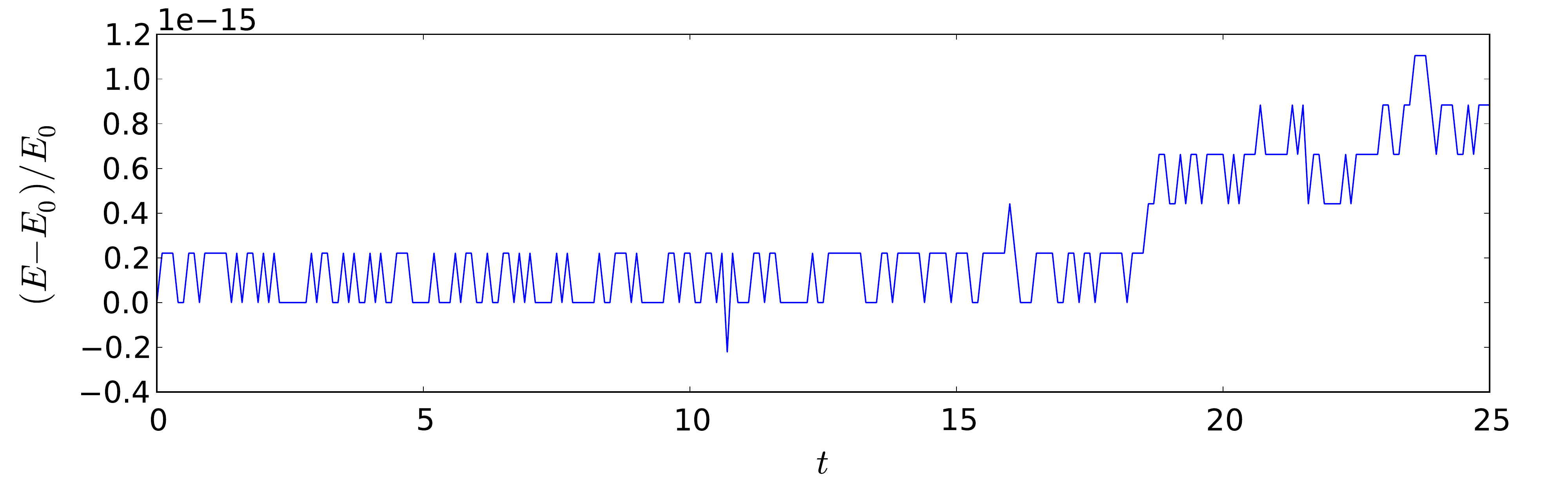}
\includegraphics[width=\textwidth]{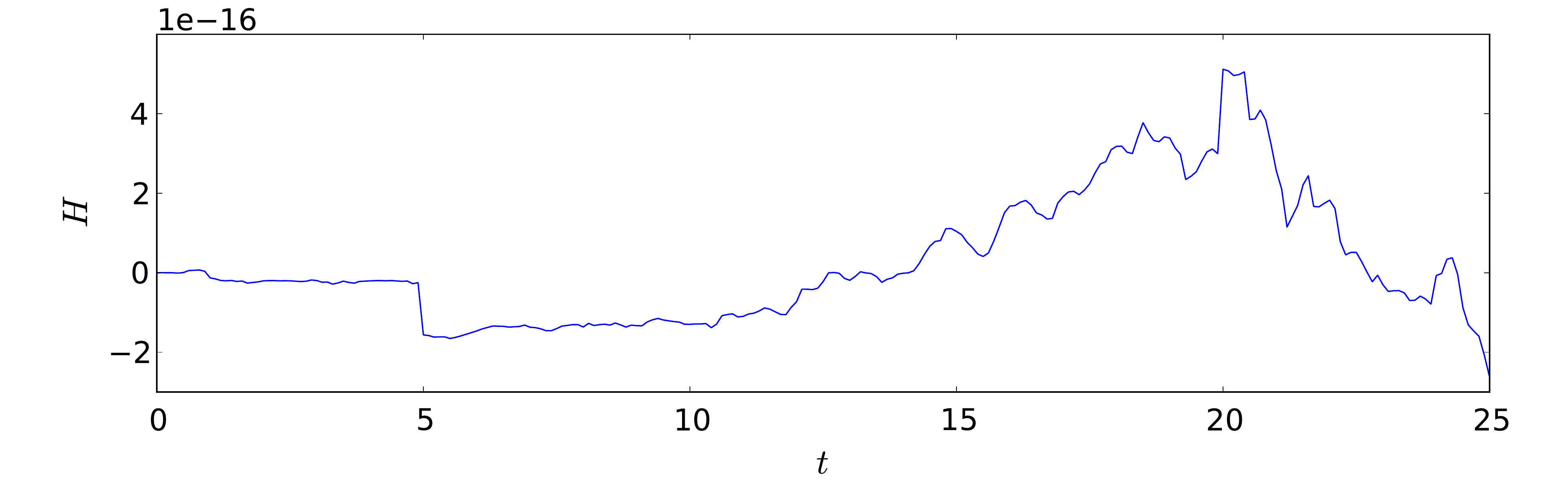}
\caption{Current sheath, $32 \times 32$ grid points. Conservation of energy and cross helicity.}
\label{fig:current_sheath_32x32_errors}
\end{figure}

\begin{figure}
\centering
\includegraphics[width=\textwidth]{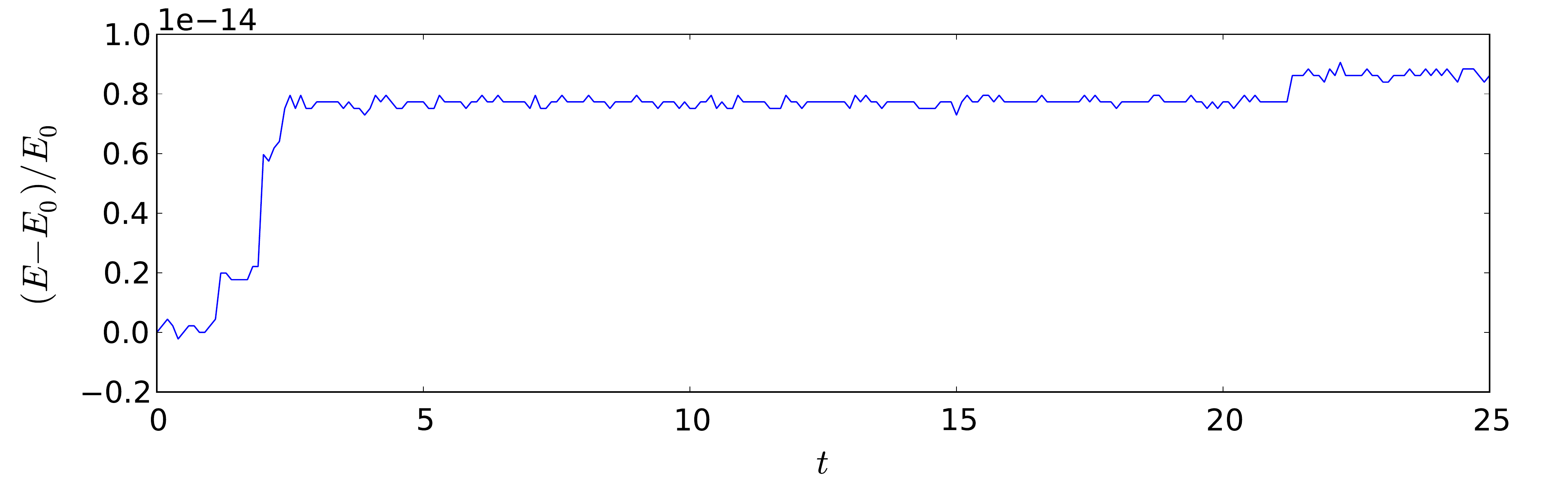}
\includegraphics[width=\textwidth]{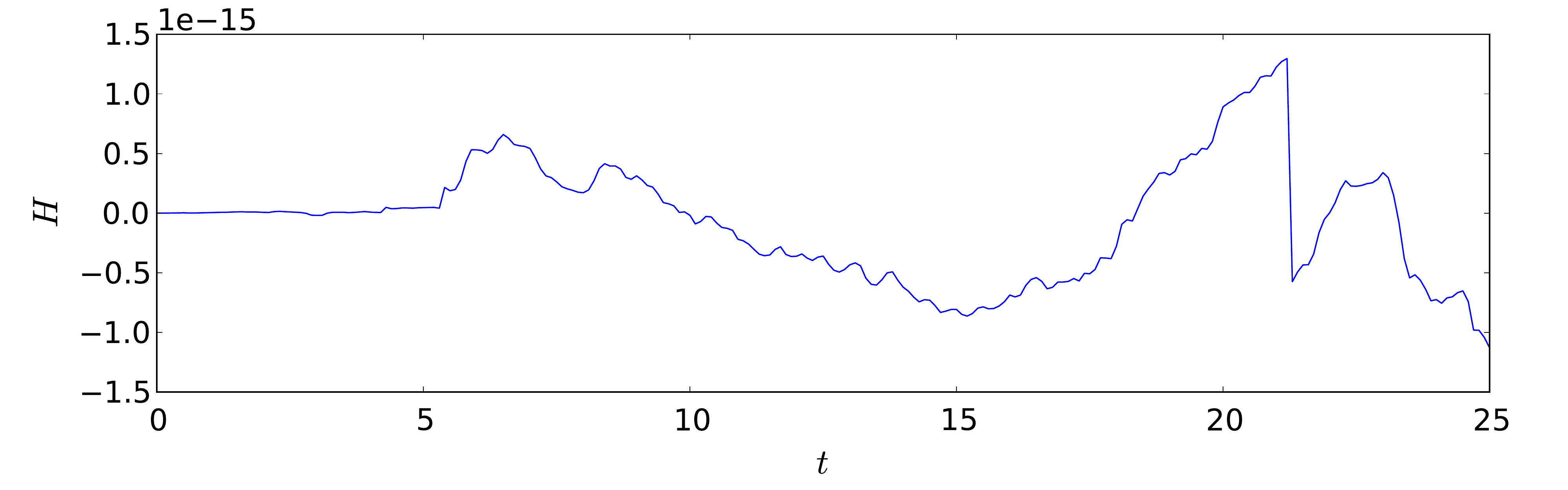}
\caption{Current sheath, $64 \times 64$ grid points. Conservation of energy and cross helicity.}
\label{fig:current_sheath_64x64_errors}
\end{figure}

\begin{figure}
\centering
\includegraphics[width=\textwidth]{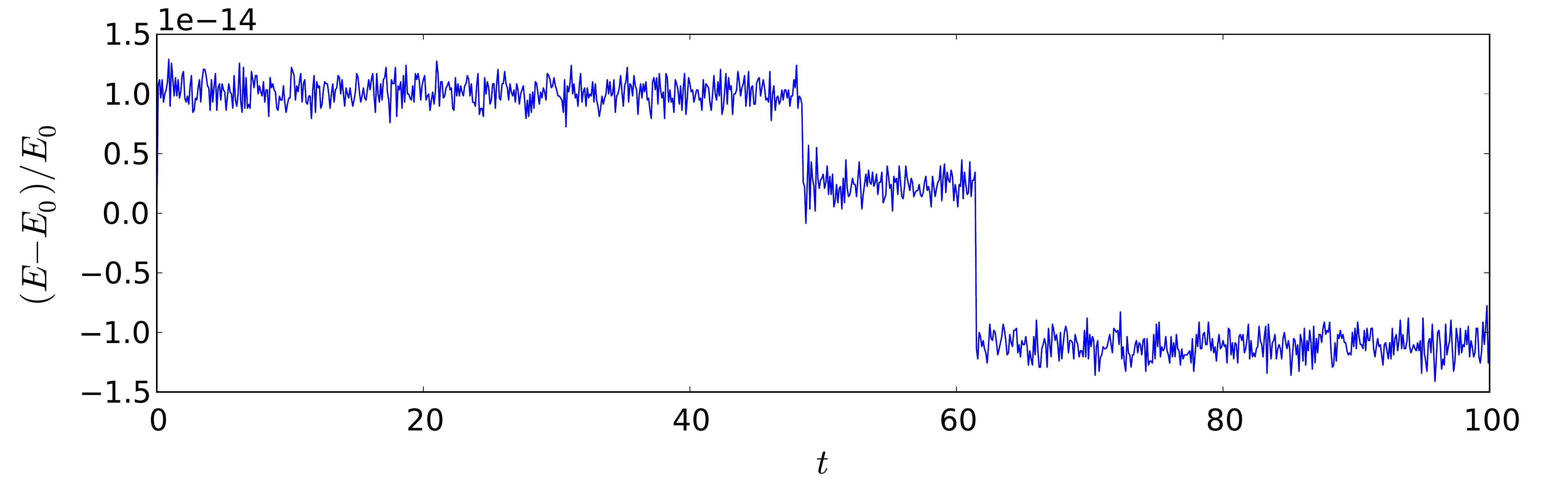}
\includegraphics[width=\textwidth]{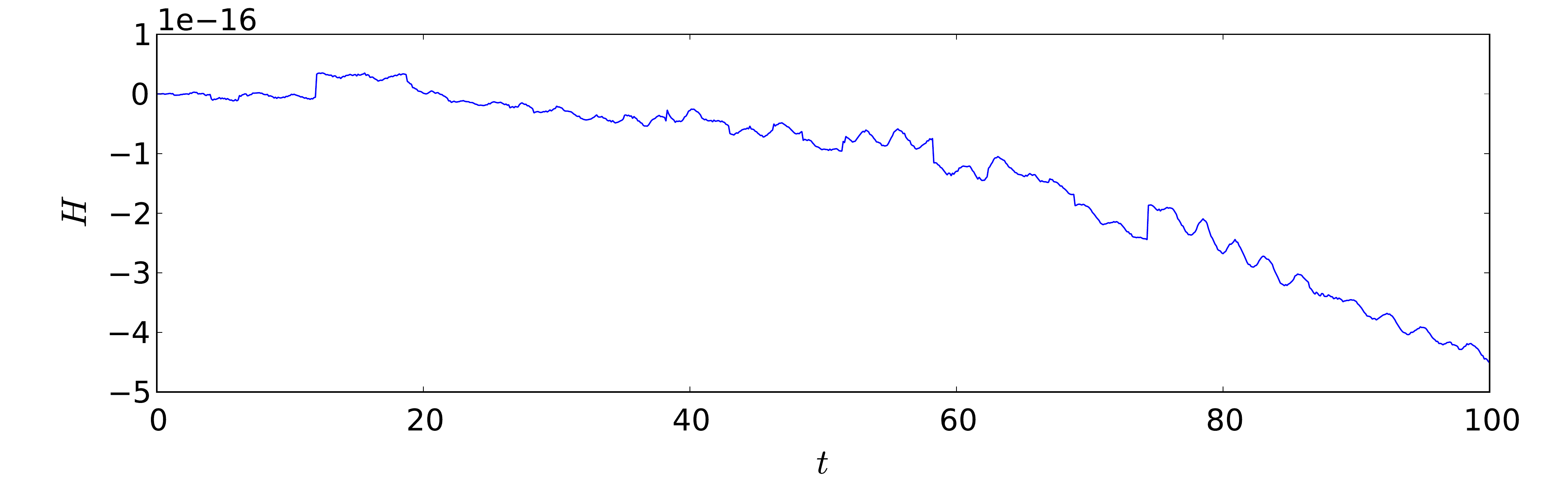}
\caption{Cosh current sheath, $30 \times 30$ grid points. Conservation of energy and cross helicity.}
\label{fig:current_sheath_cosh_errors}
\end{figure}

\begin{figure}
\centering
\includegraphics[width=\textwidth]{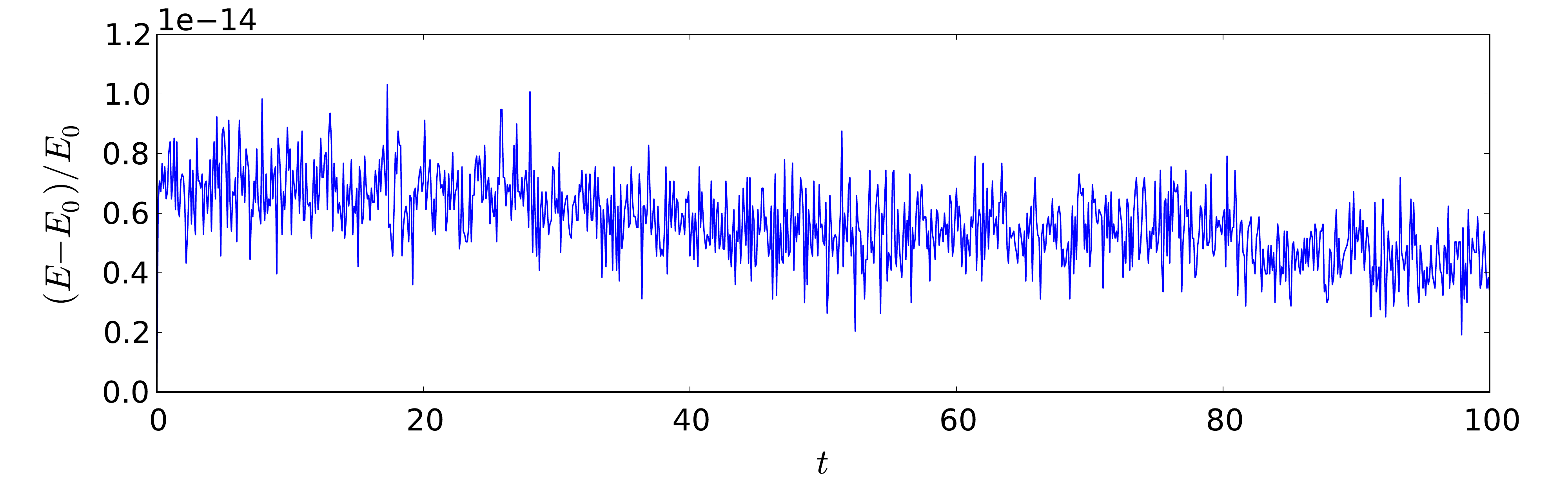}
\includegraphics[width=\textwidth]{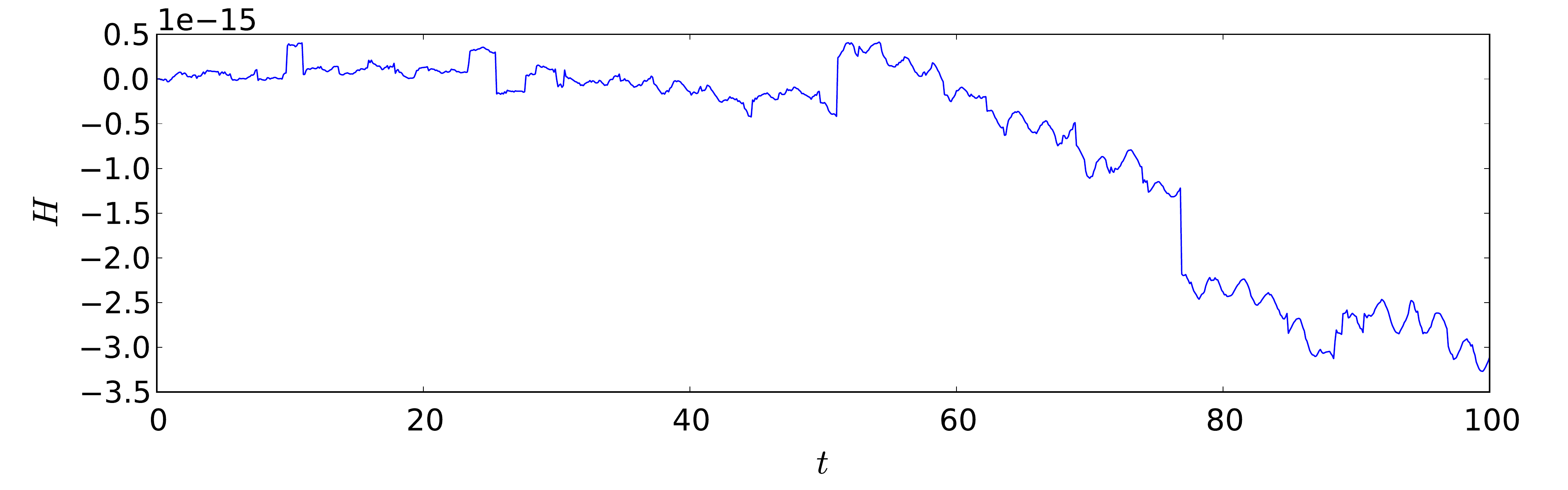}
\caption{Tanh current sheath, $30 \times 30$ grid points. Conservation of energy and cross helicity.}
\label{fig:current_sheath_tanh_errors}
\end{figure}

In all of the considered cases, the conservation of energy and cross helicity is optimal (see figures \ref{fig:current_sheath_32x32_errors} - \ref{fig:current_sheath_tanh_errors}).
The jumps in the energy error of the $\cosh$ case (figure \ref{fig:current_sheath_cosh_errors}) are probably due to the Newton solver converging to a slightly larger residual than during the rest of the simulation. Nevertheless, the energy error is extraordinary small.
Thus here we want to focus on the conservation of field line topology.

\begin{figure}
\centering
\subfloat{
\includegraphics[width=.32\textwidth]{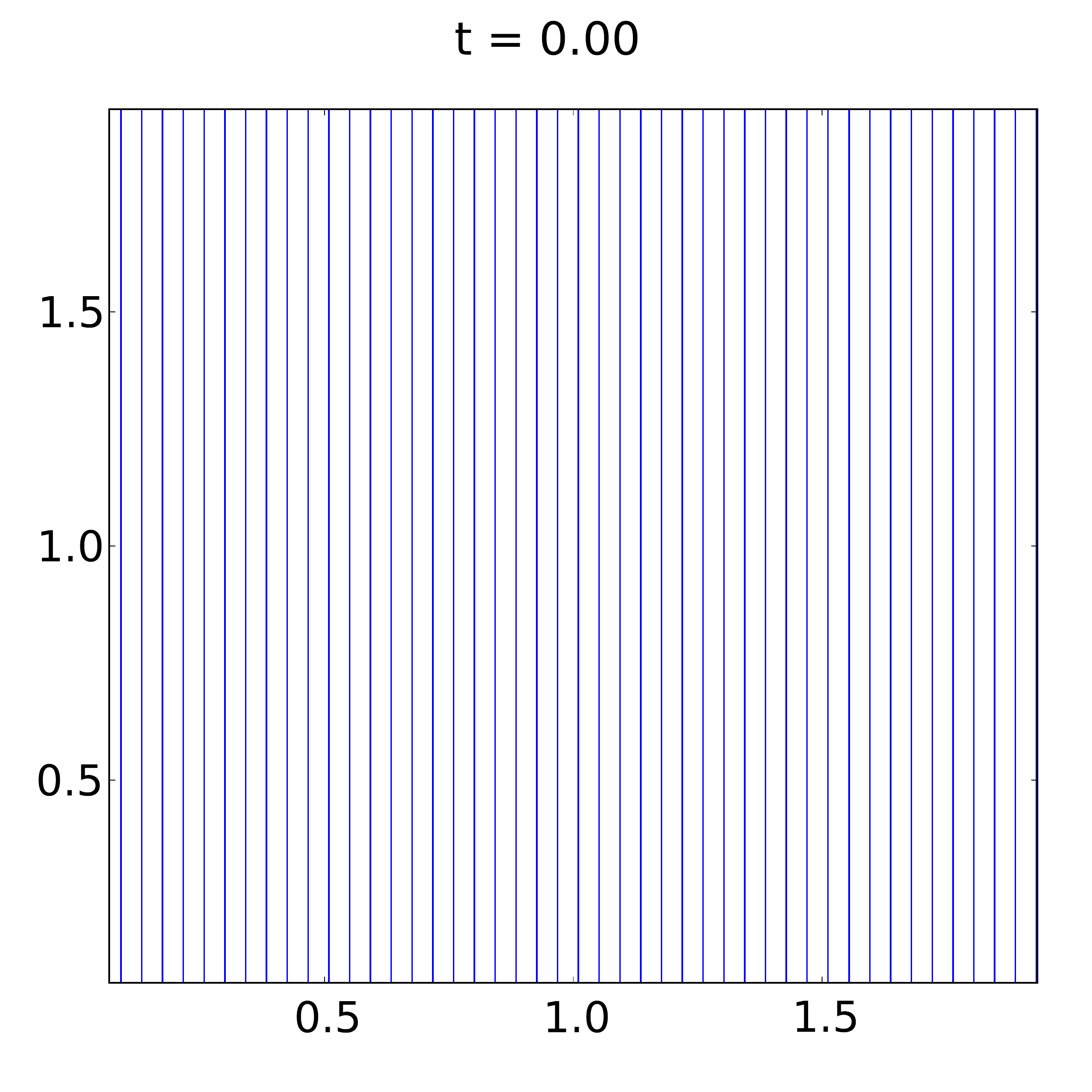}
}
\subfloat{
\includegraphics[width=.32\textwidth]{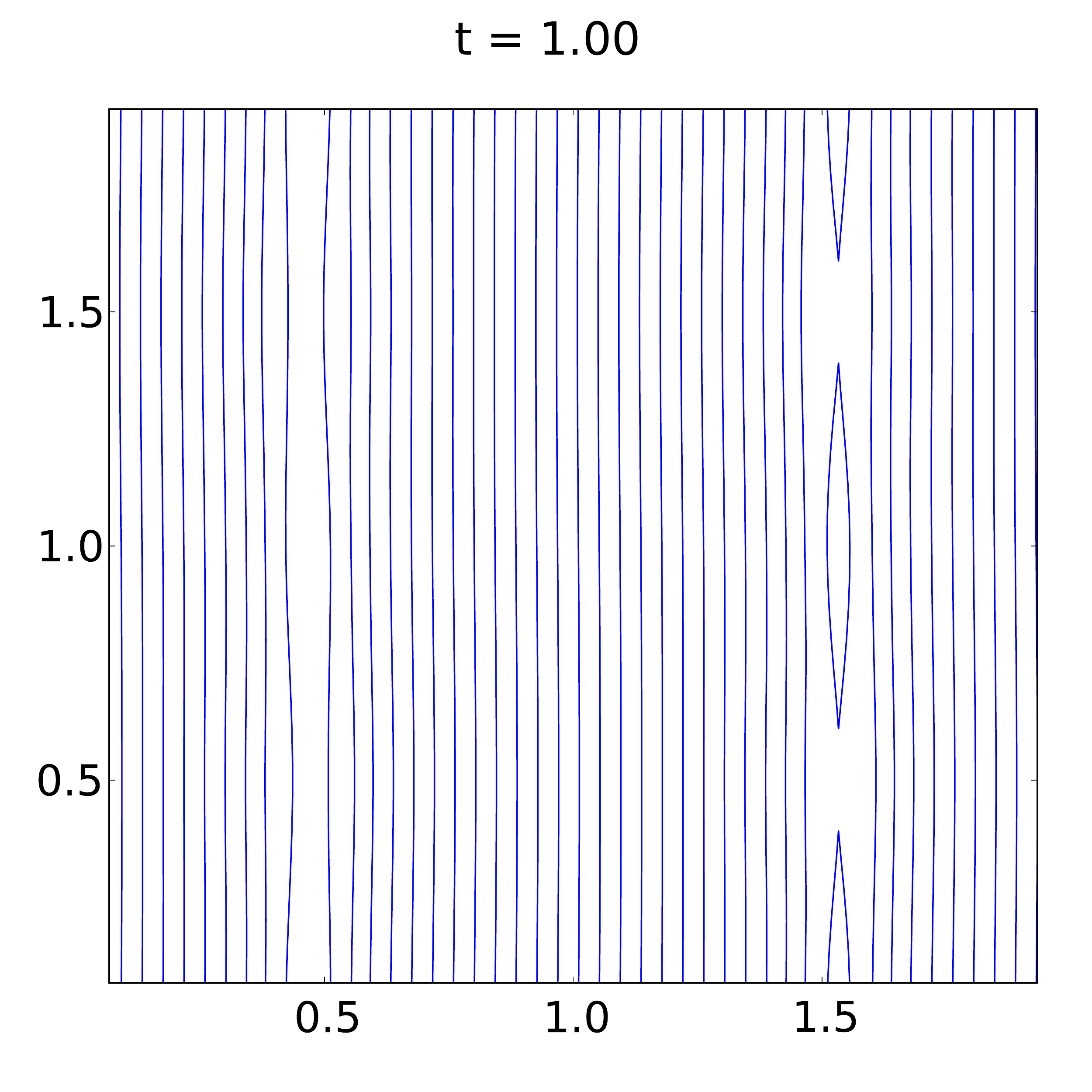}
}
\subfloat{
\includegraphics[width=.32\textwidth]{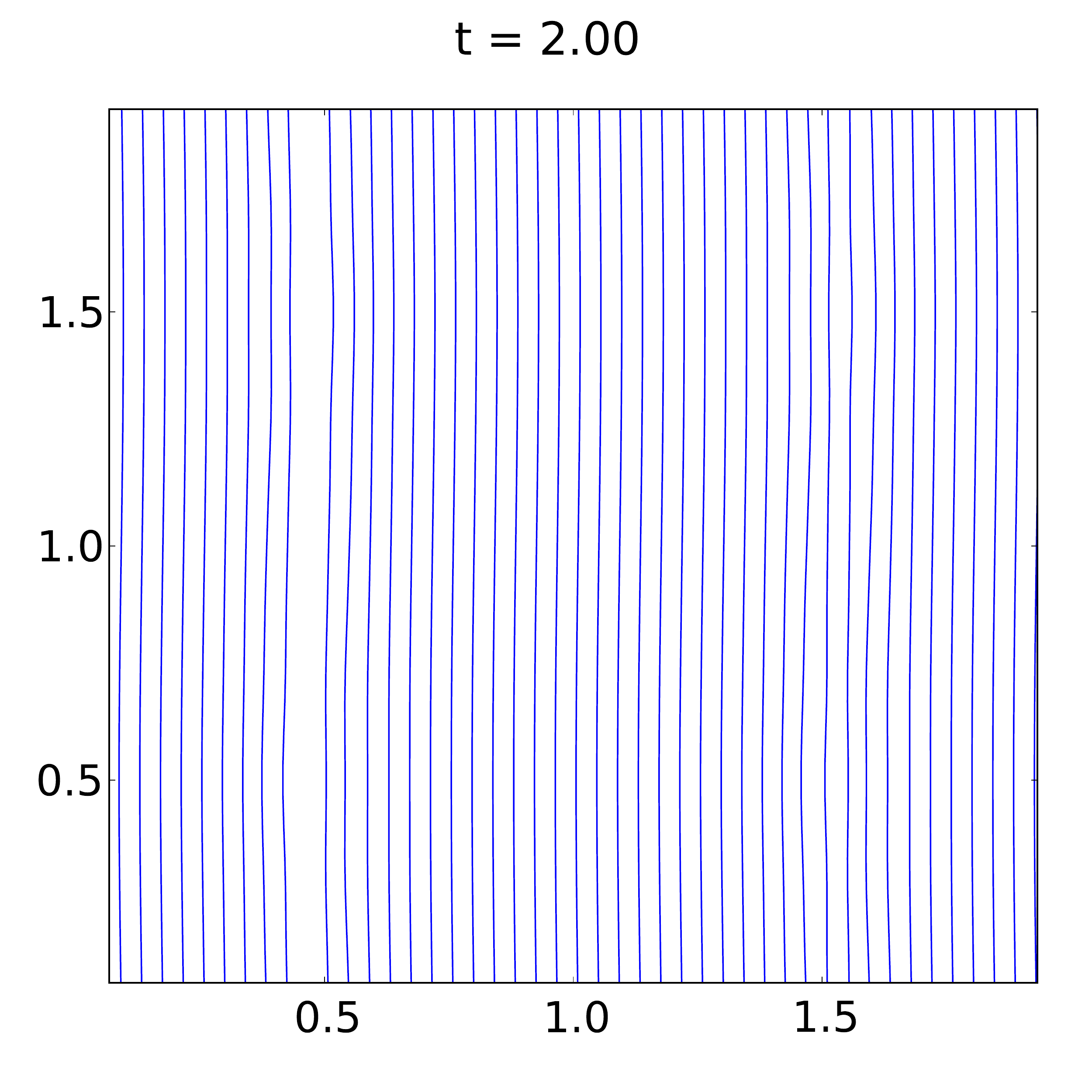}
}

\subfloat{
\includegraphics[width=.32\textwidth]{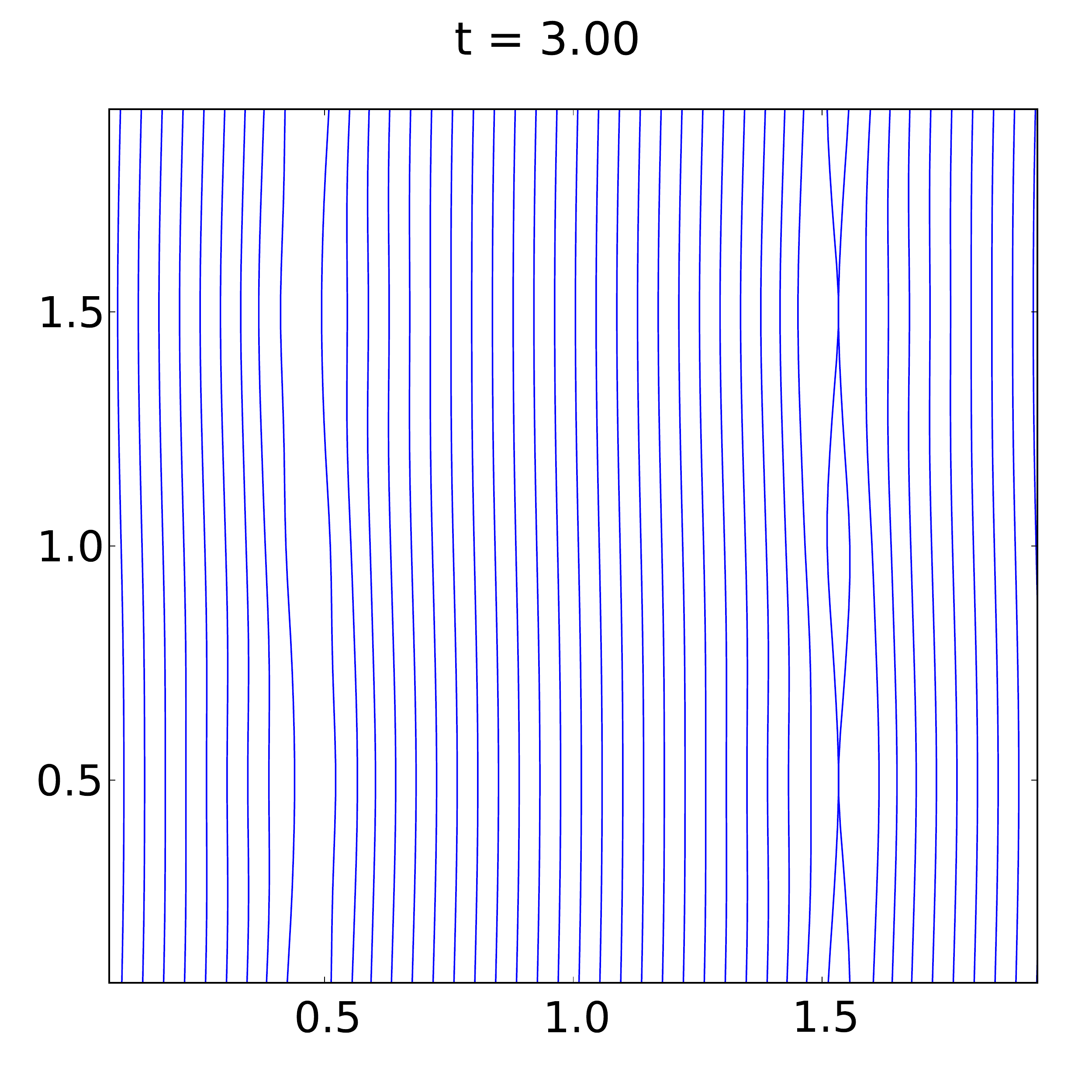}
}
\subfloat{
\includegraphics[width=.32\textwidth]{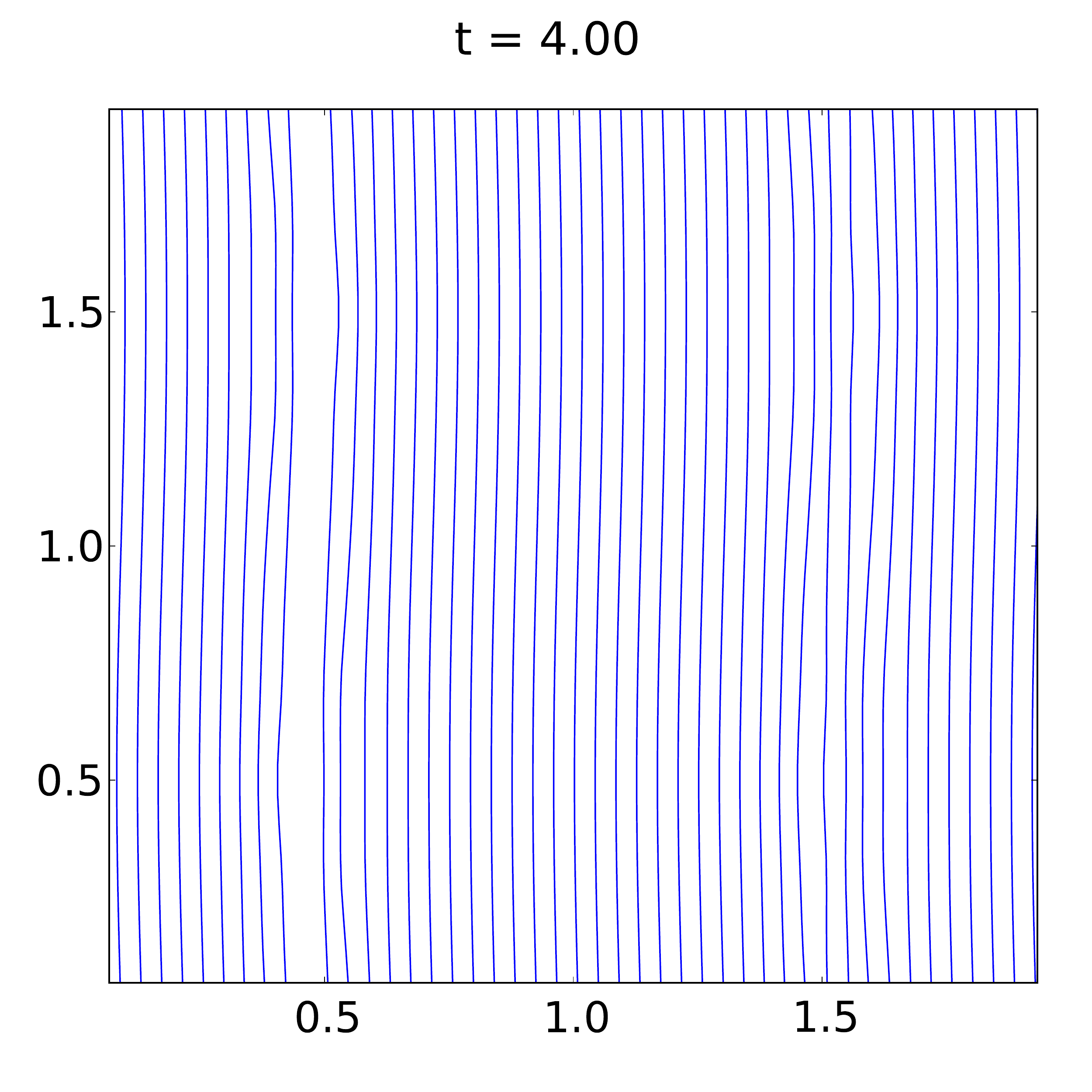}
}
\subfloat{
\includegraphics[width=.32\textwidth]{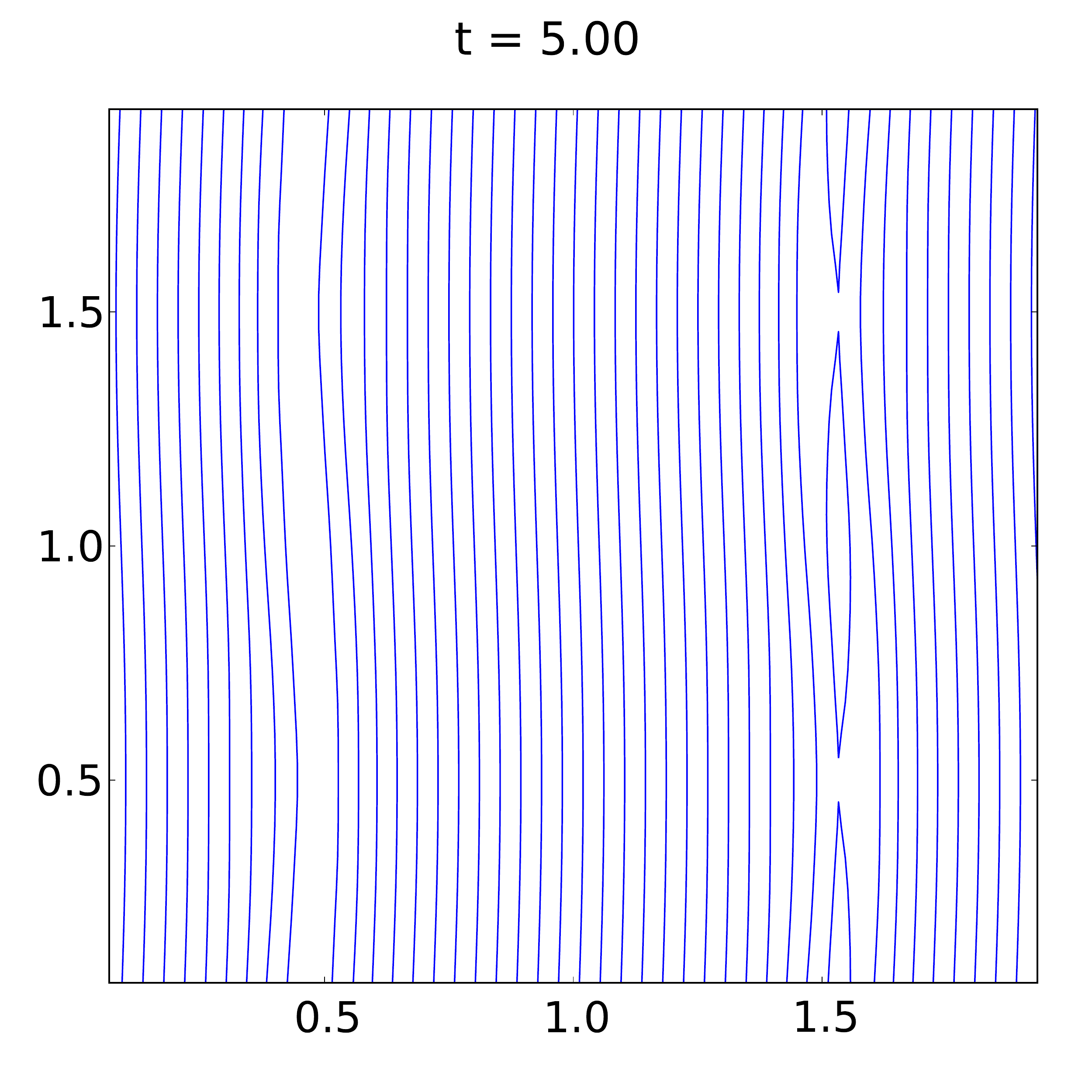}
}

\subfloat{
\includegraphics[width=.32\textwidth]{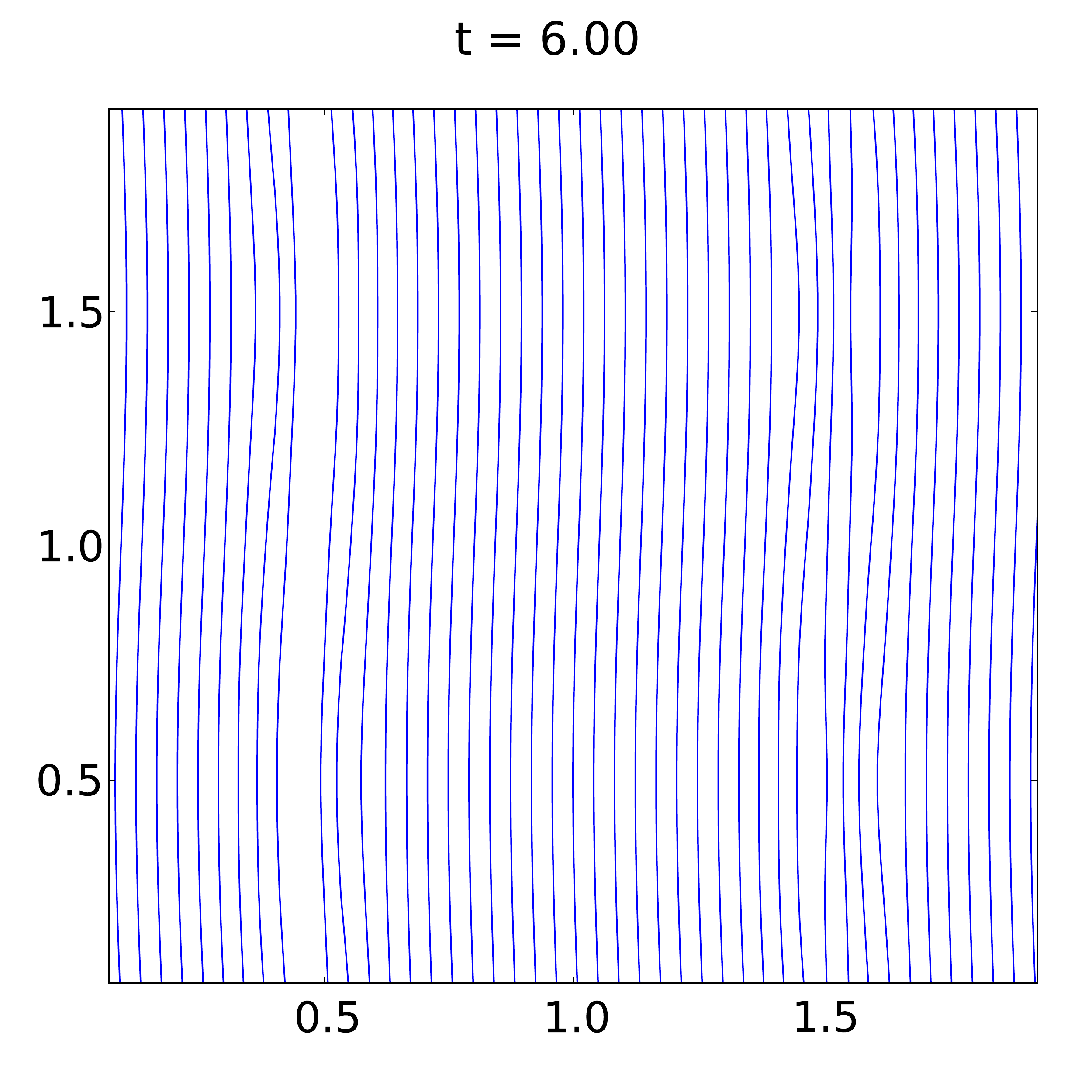}
}
\subfloat{
\includegraphics[width=.32\textwidth]{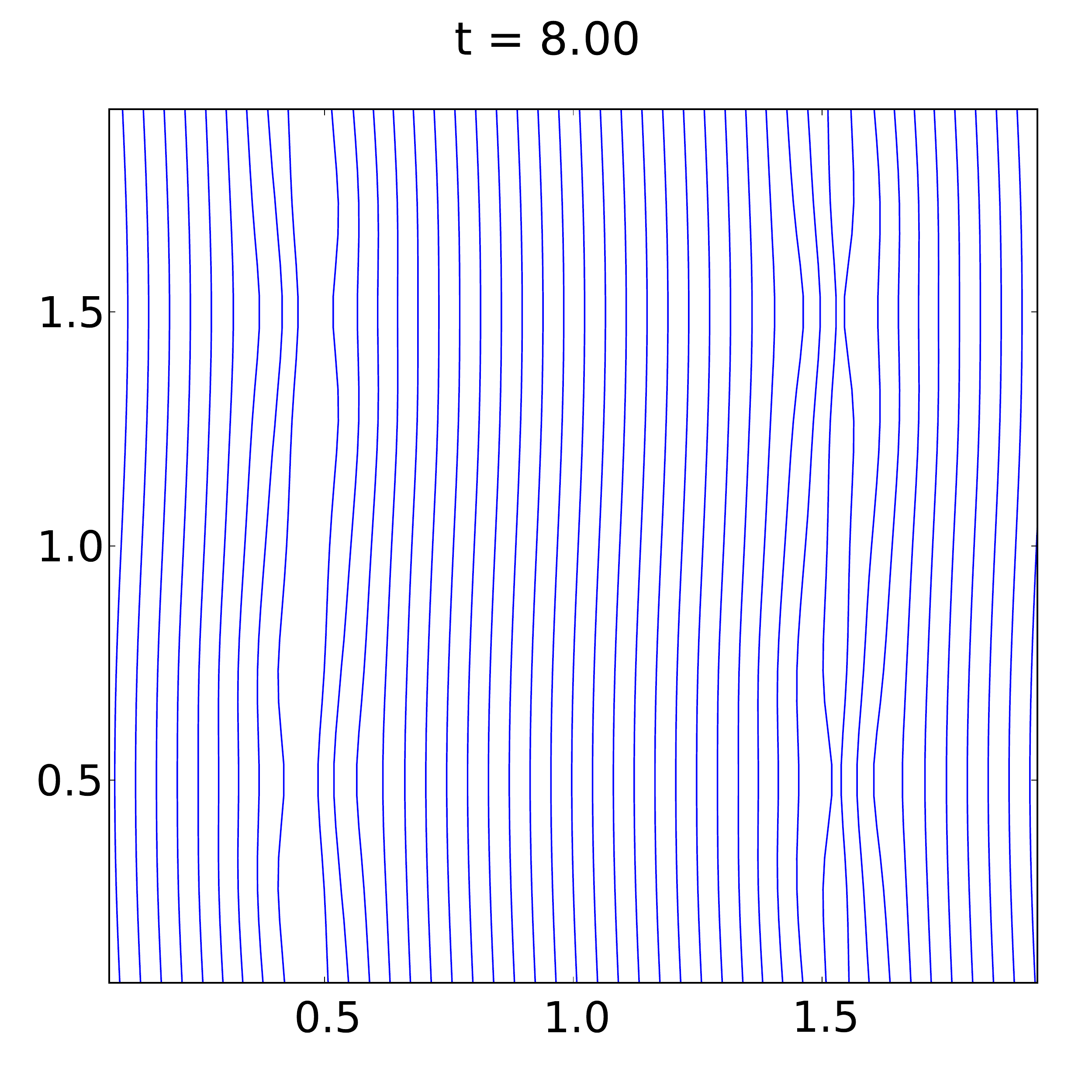}
}
\subfloat{
\includegraphics[width=.32\textwidth]{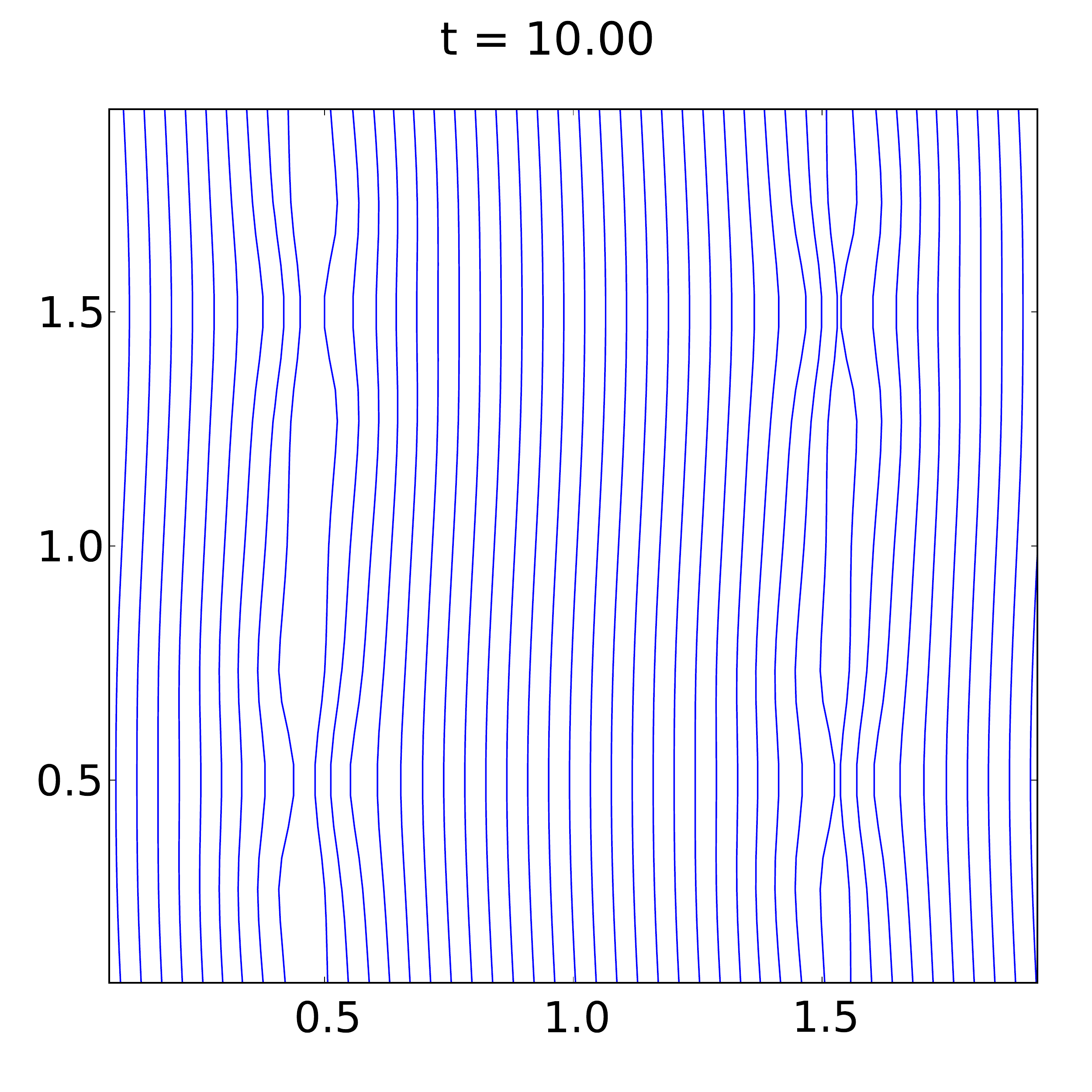}
}

\subfloat{
\includegraphics[width=.32\textwidth]{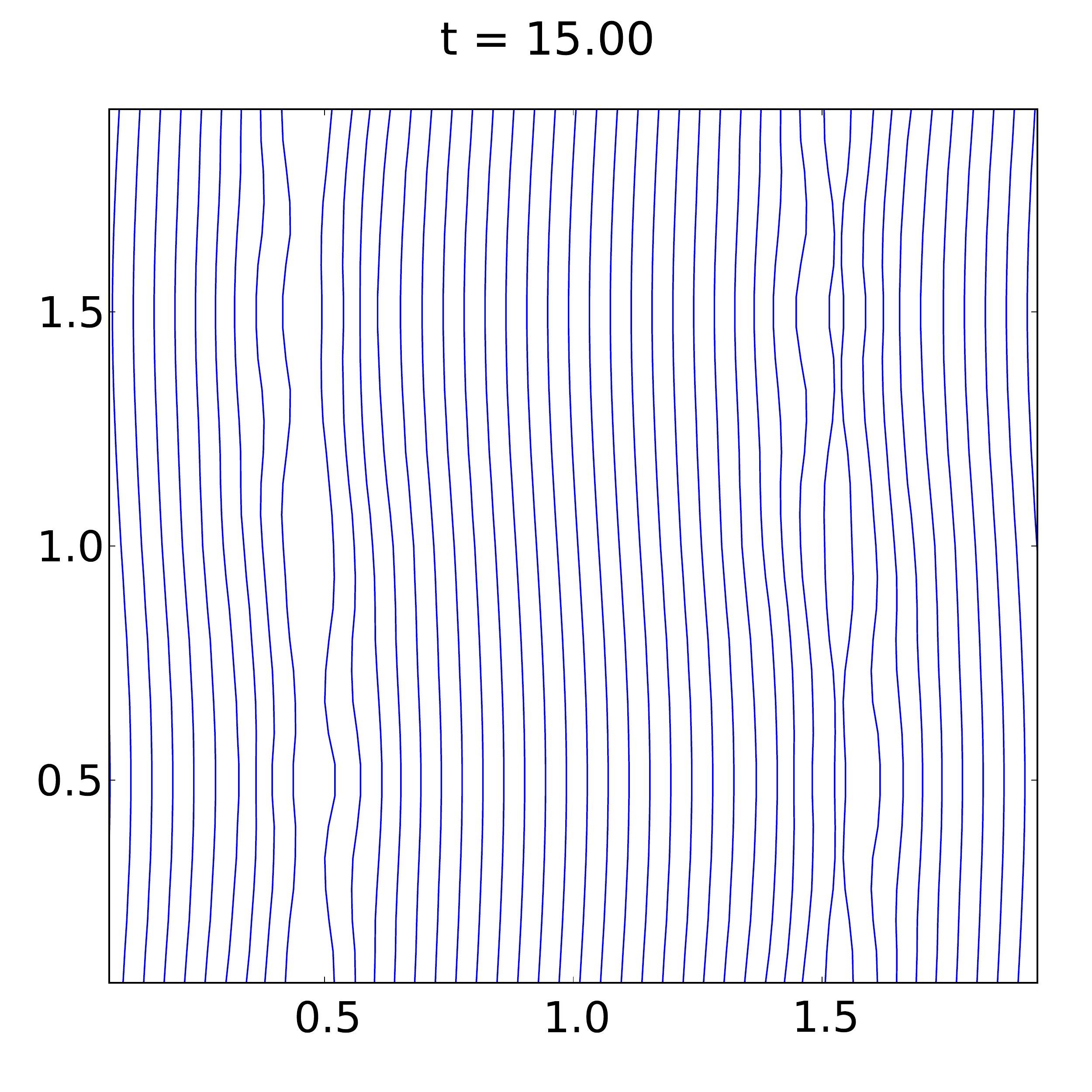}
}
\subfloat{
\includegraphics[width=.32\textwidth]{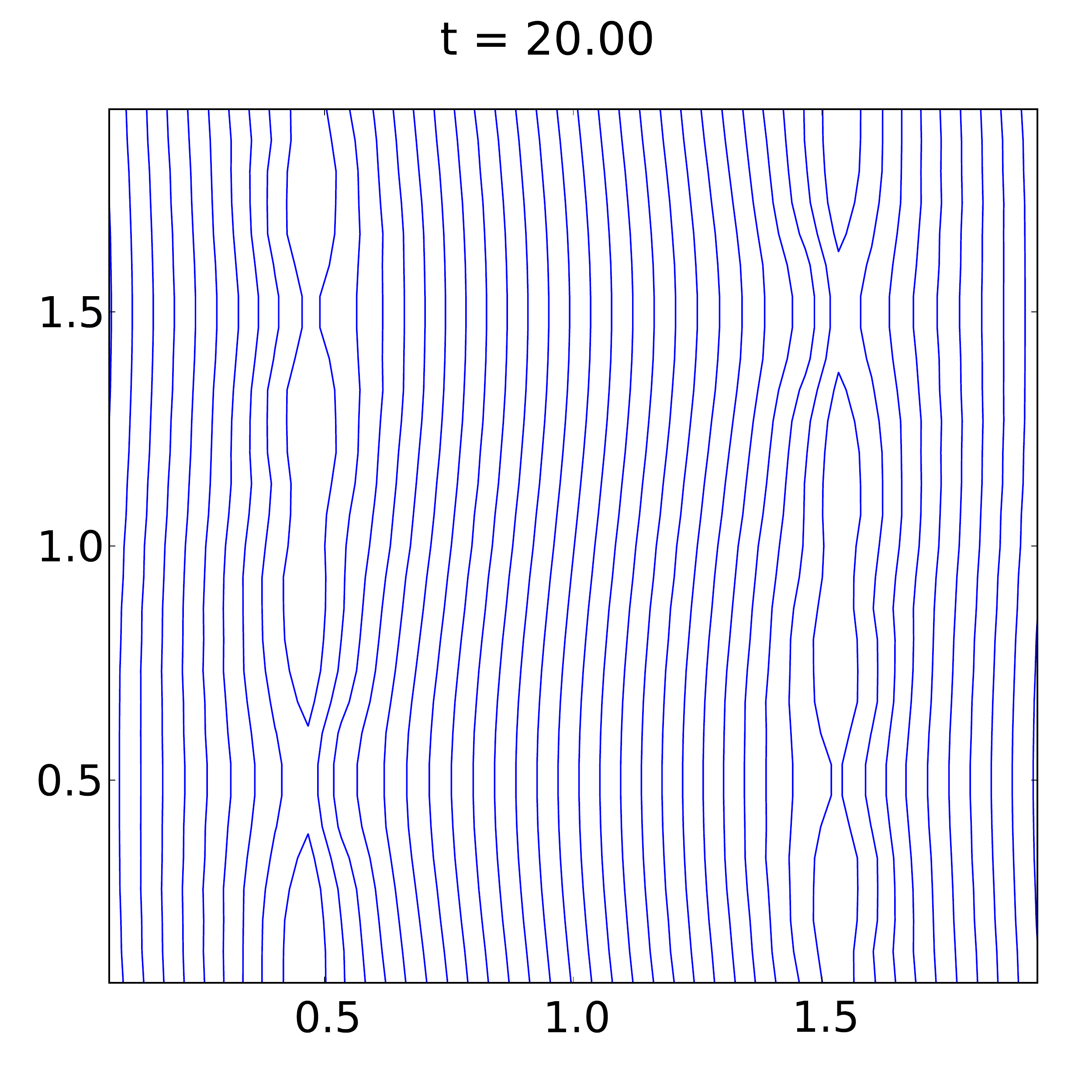}
}
\subfloat{
\includegraphics[width=.32\textwidth]{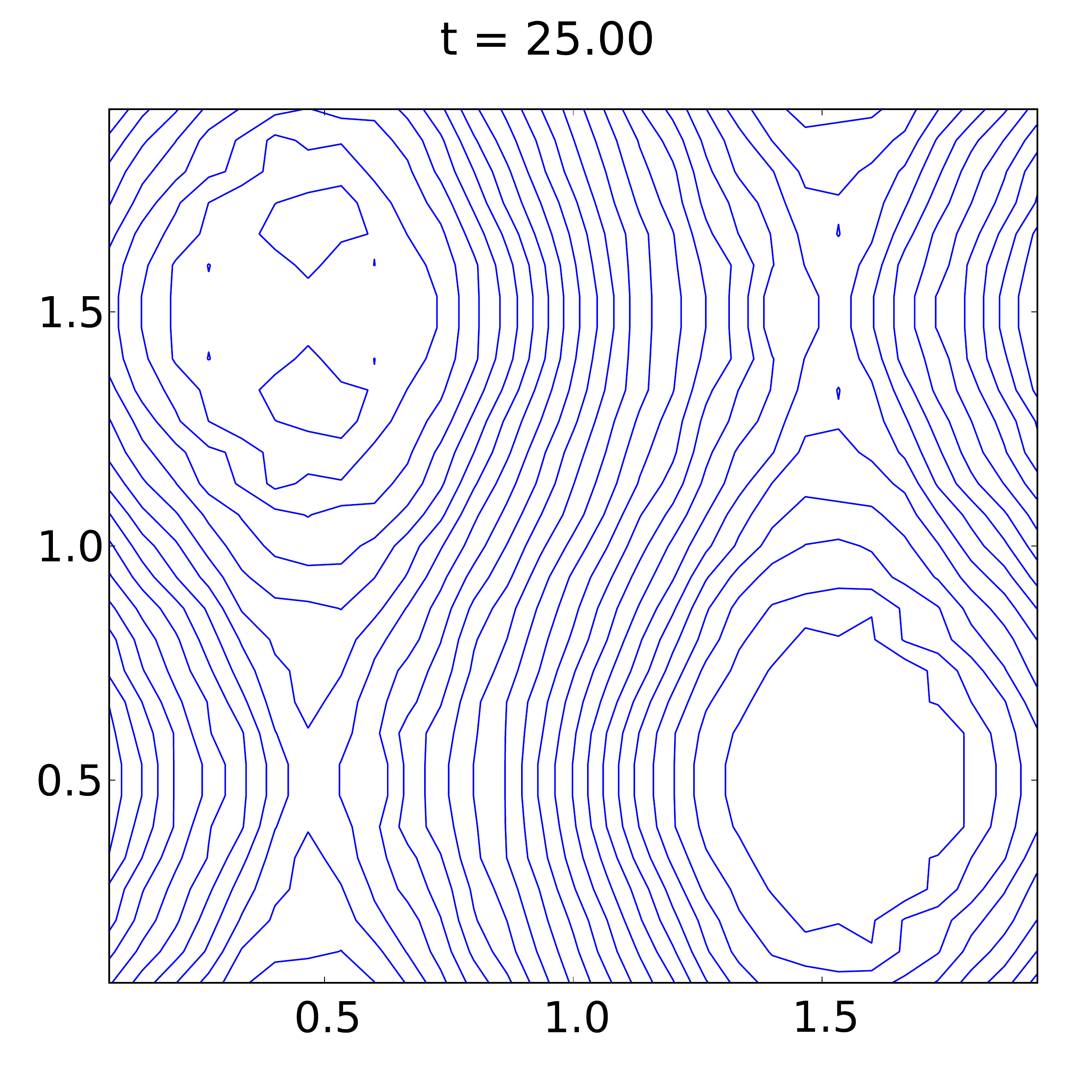}
}

\caption{Sharp current sheath, $30 \times 30$ grid points. Magnetic field lines.}
\label{fig:current_sheath_30_field_lines}
\end{figure}

In figure \ref{fig:current_sheath_30_field_lines}, the field line evolution for the sharp jump is plotted. Initially all field lines are parallel. Due to the perturbation in the velocity field, the magnetic field lines get bend, but for most of the time, they do not break up and reconnect. At $t = 1$ and $t = 5$ we see what appears to be reconnection events, but as these islands disappear very soon after they form, this is likely to be attributed to artefacts of the plotting function (like an inaccurate interpolation or integration of the contour lines).
At about $t = 20$, however, real islands form an consecutively grow as can be seen at $t = 25$. At this point, the solution can not be regarded physical anymore.

\begin{figure}
\centering
\subfloat{
\includegraphics[width=.32\textwidth]{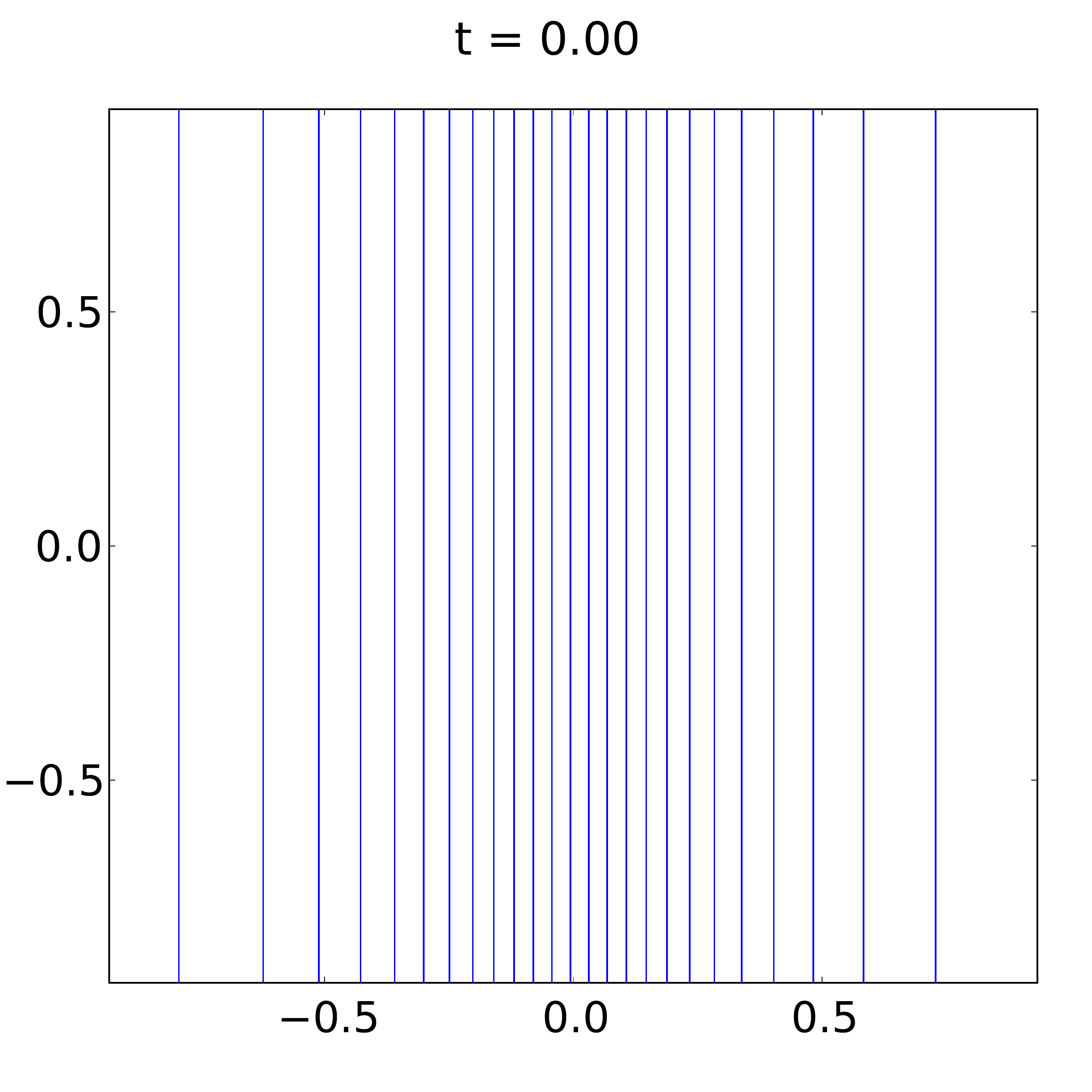}
}
\subfloat{
\includegraphics[width=.32\textwidth]{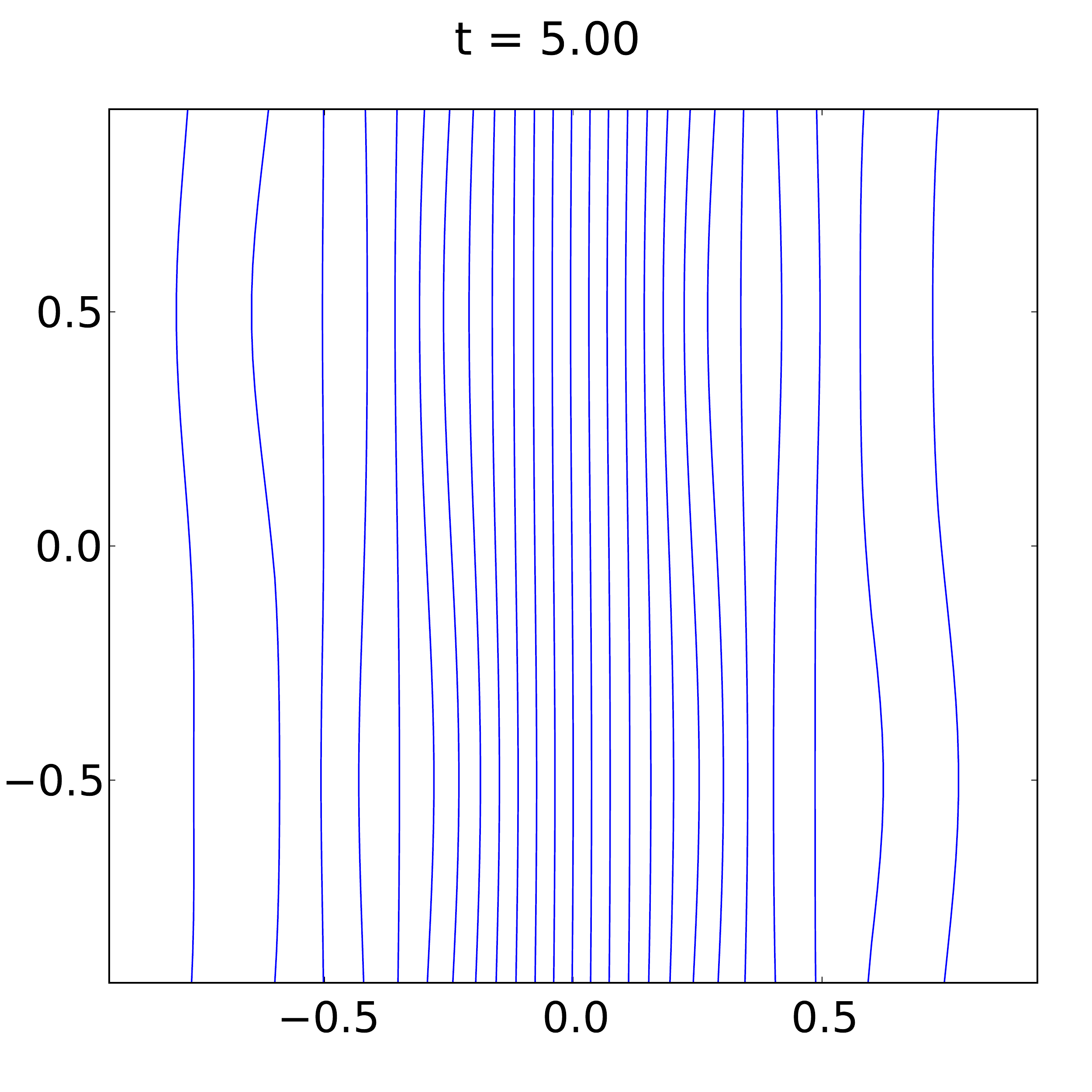}
}
\subfloat{
\includegraphics[width=.32\textwidth]{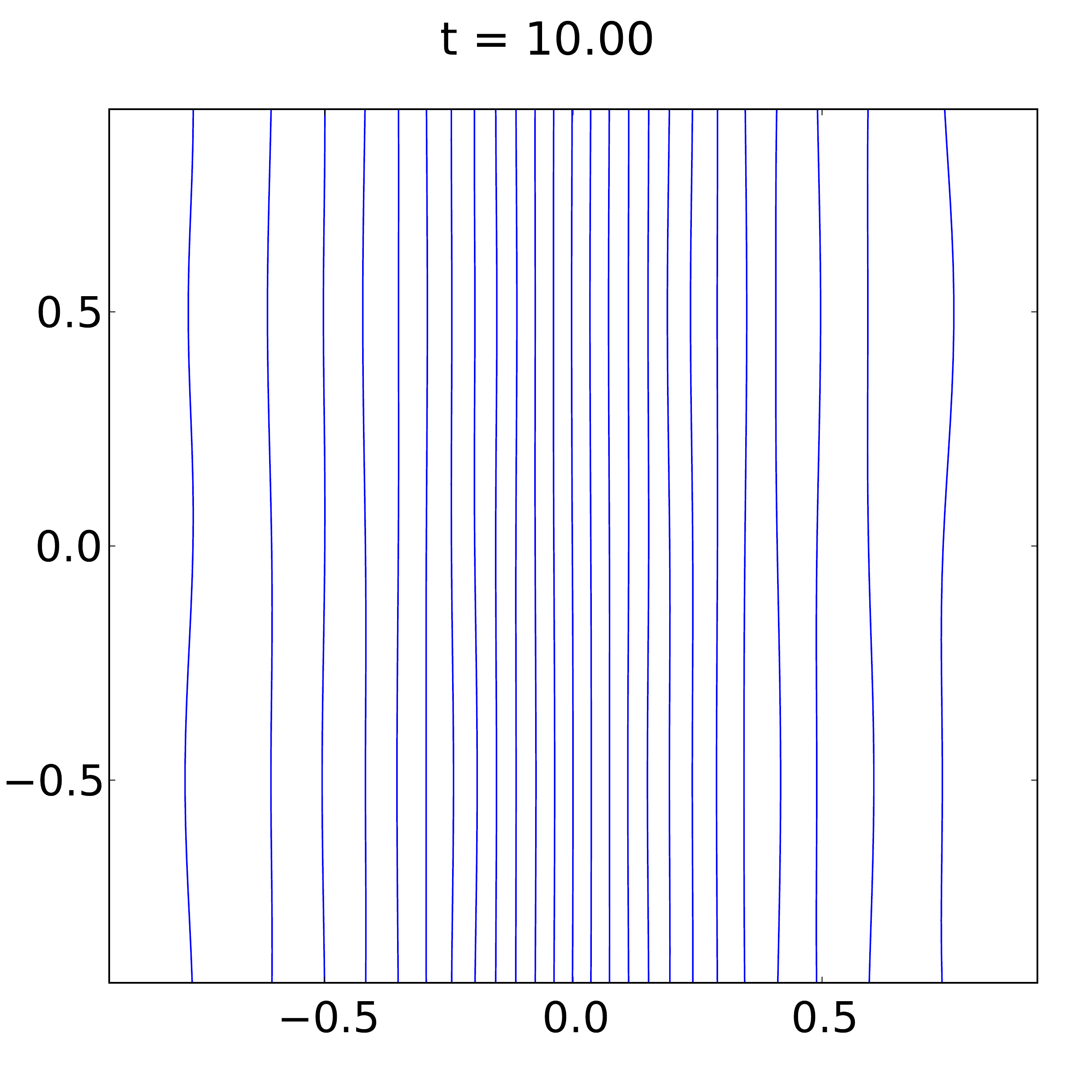}
}

\subfloat{
\includegraphics[width=.32\textwidth]{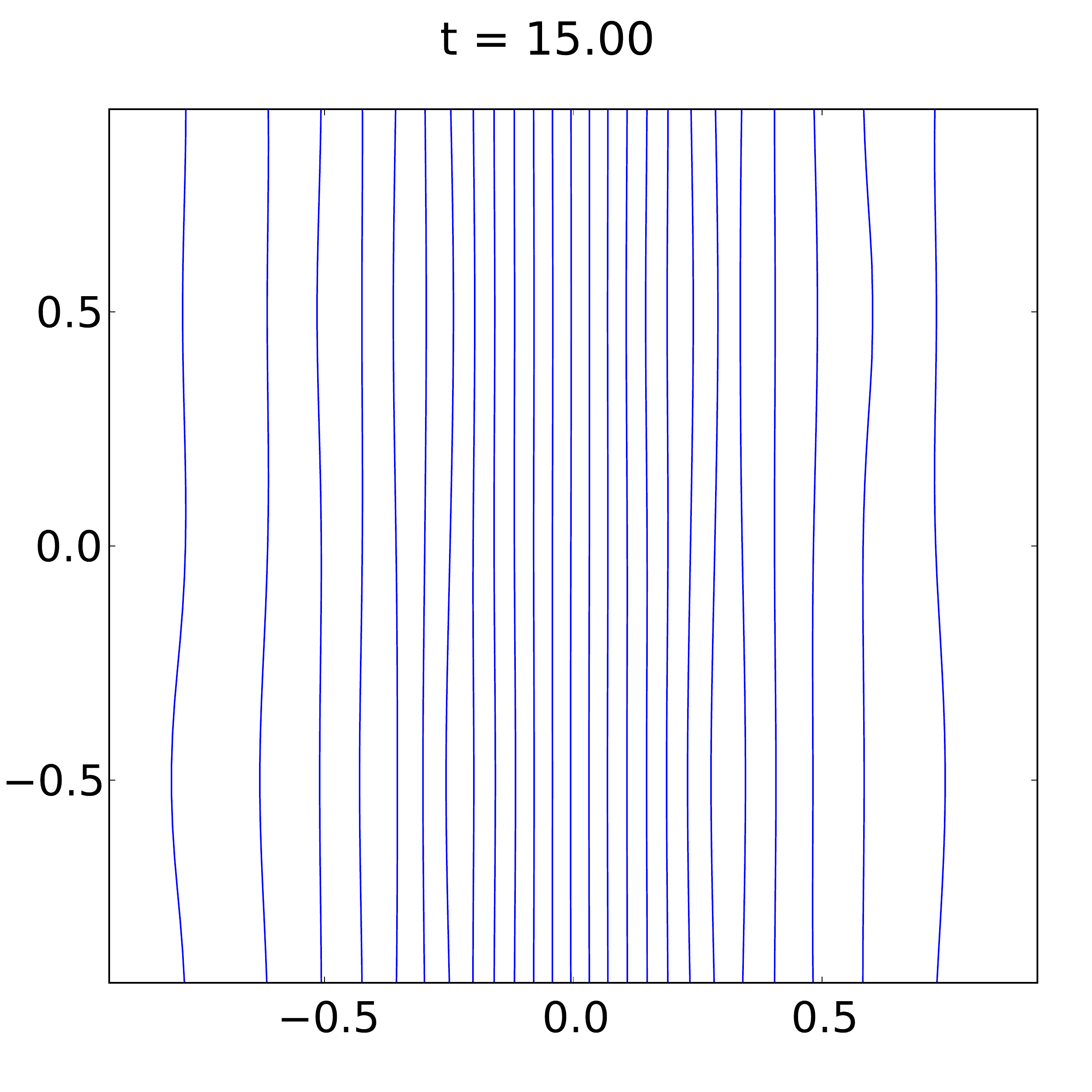}
}
\subfloat{
\includegraphics[width=.32\textwidth]{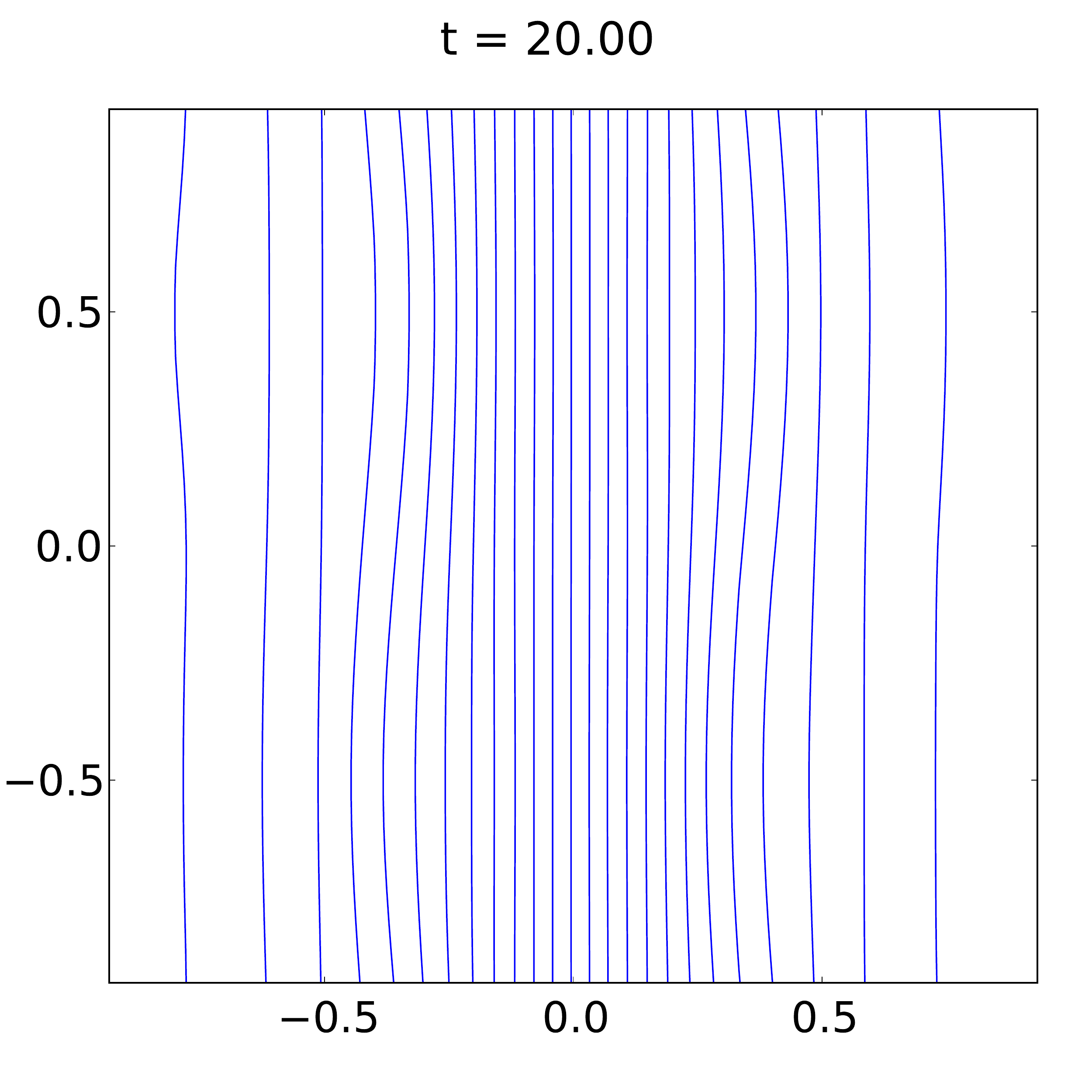}
}
\subfloat{
\includegraphics[width=.32\textwidth]{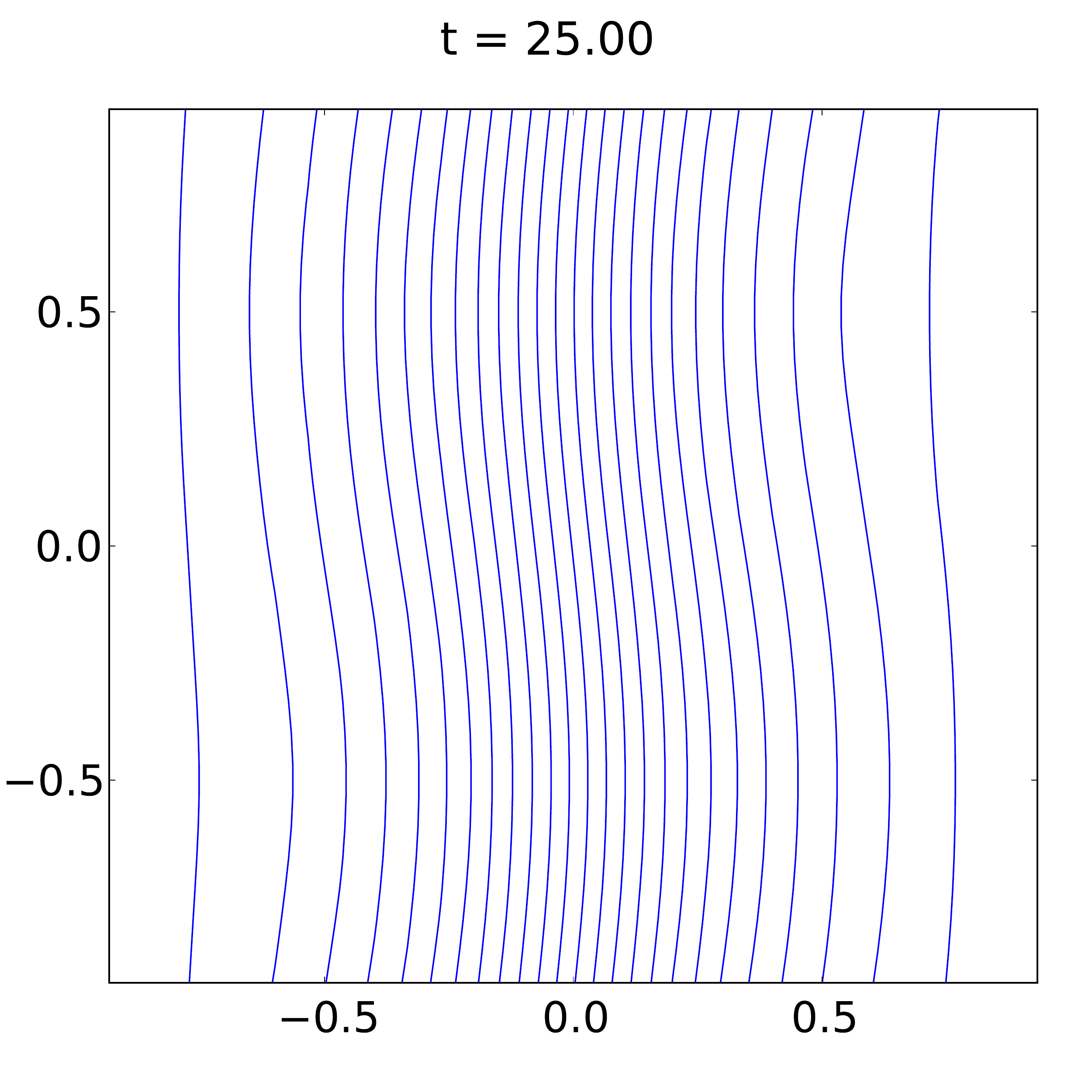}
}

\subfloat{
\includegraphics[width=.32\textwidth]{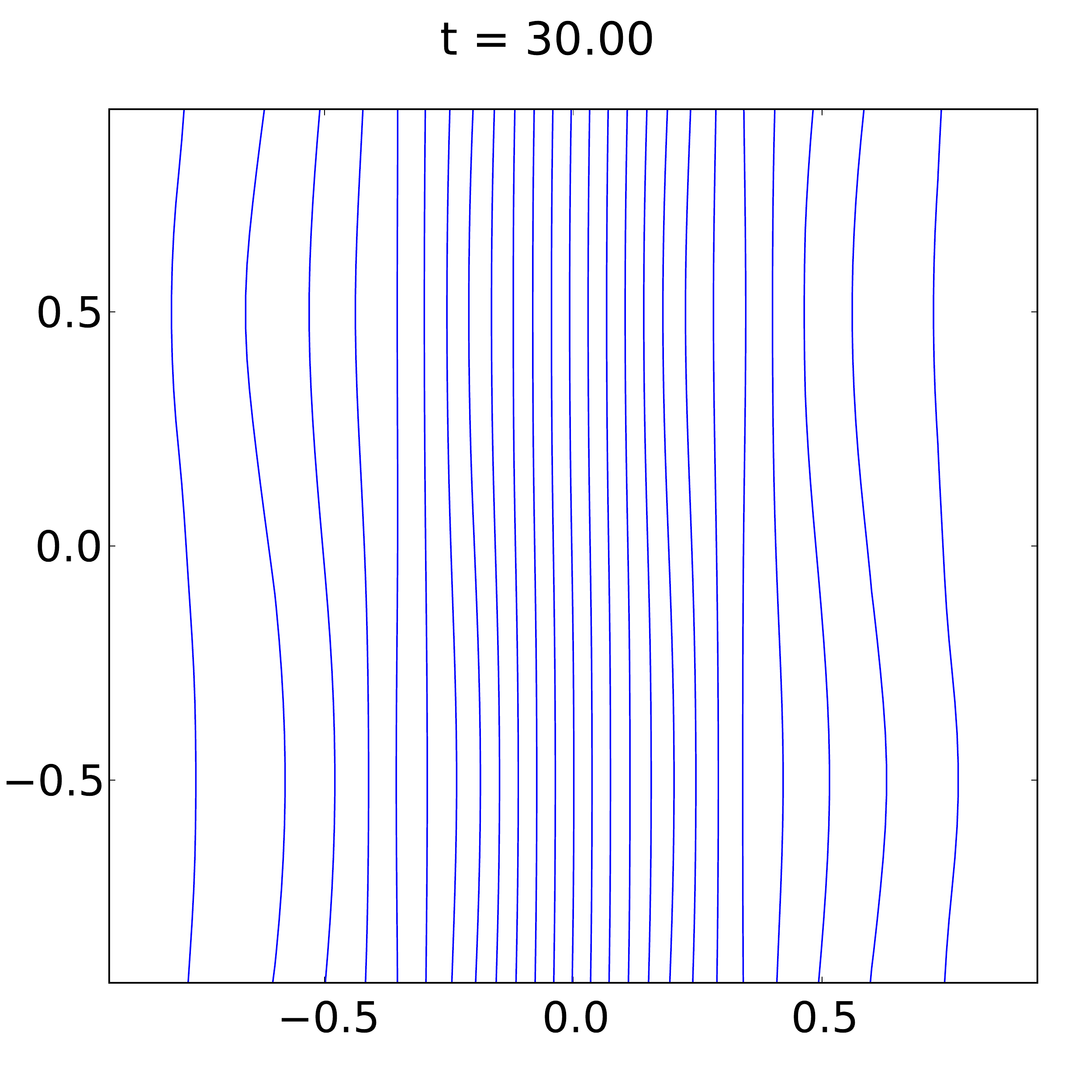}
}
\subfloat{
\includegraphics[width=.32\textwidth]{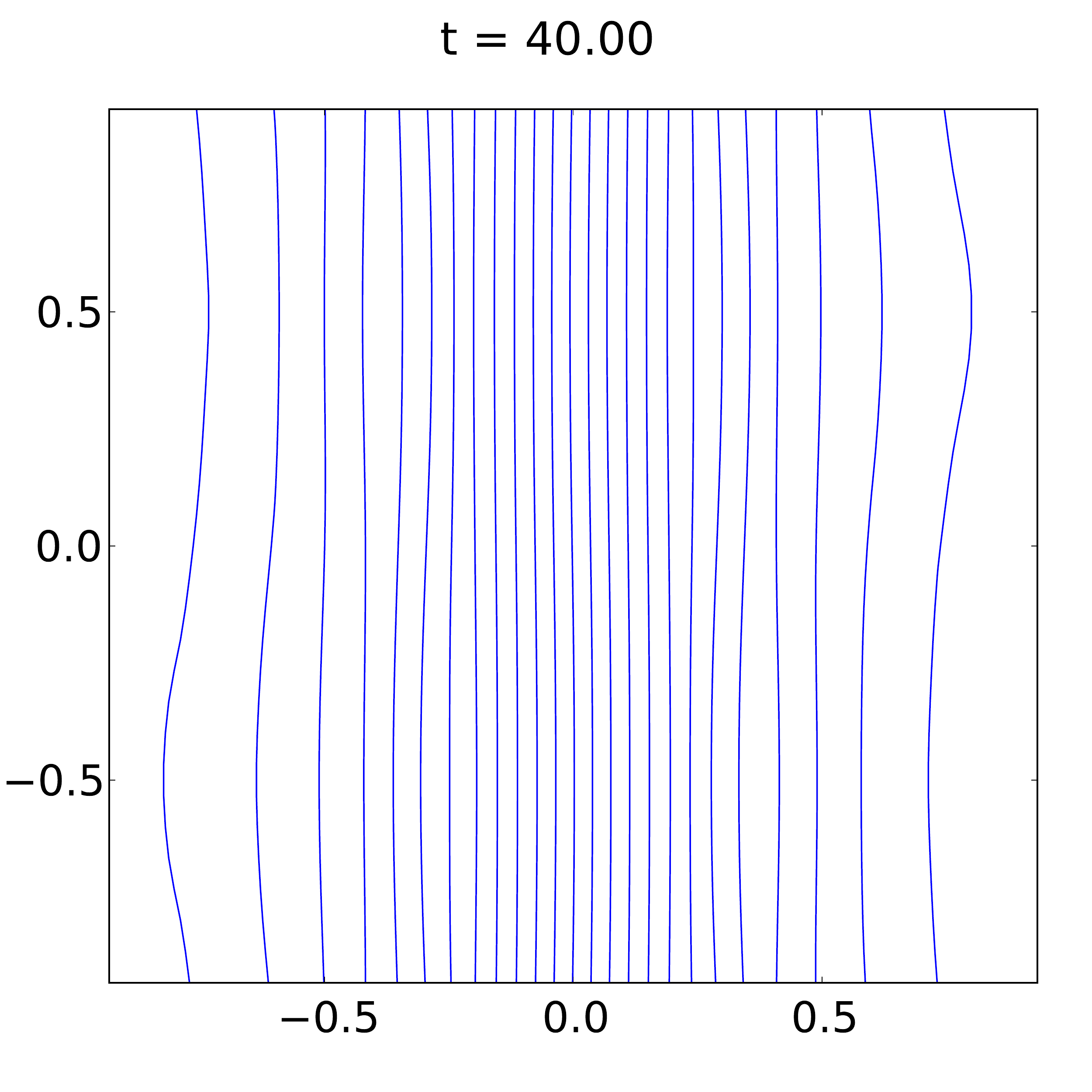}
}
\subfloat{
\includegraphics[width=.32\textwidth]{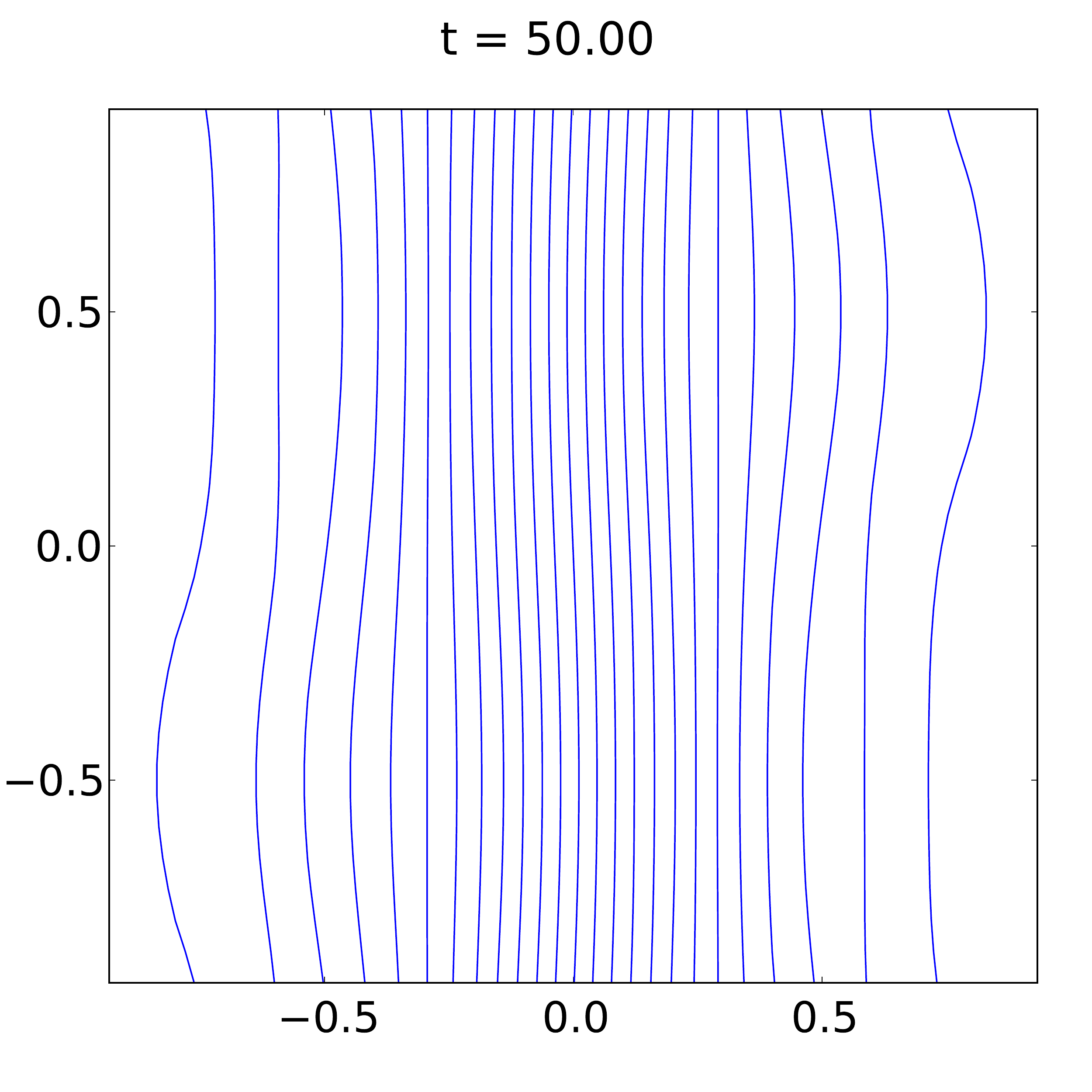}
}

\subfloat{
\includegraphics[width=.32\textwidth]{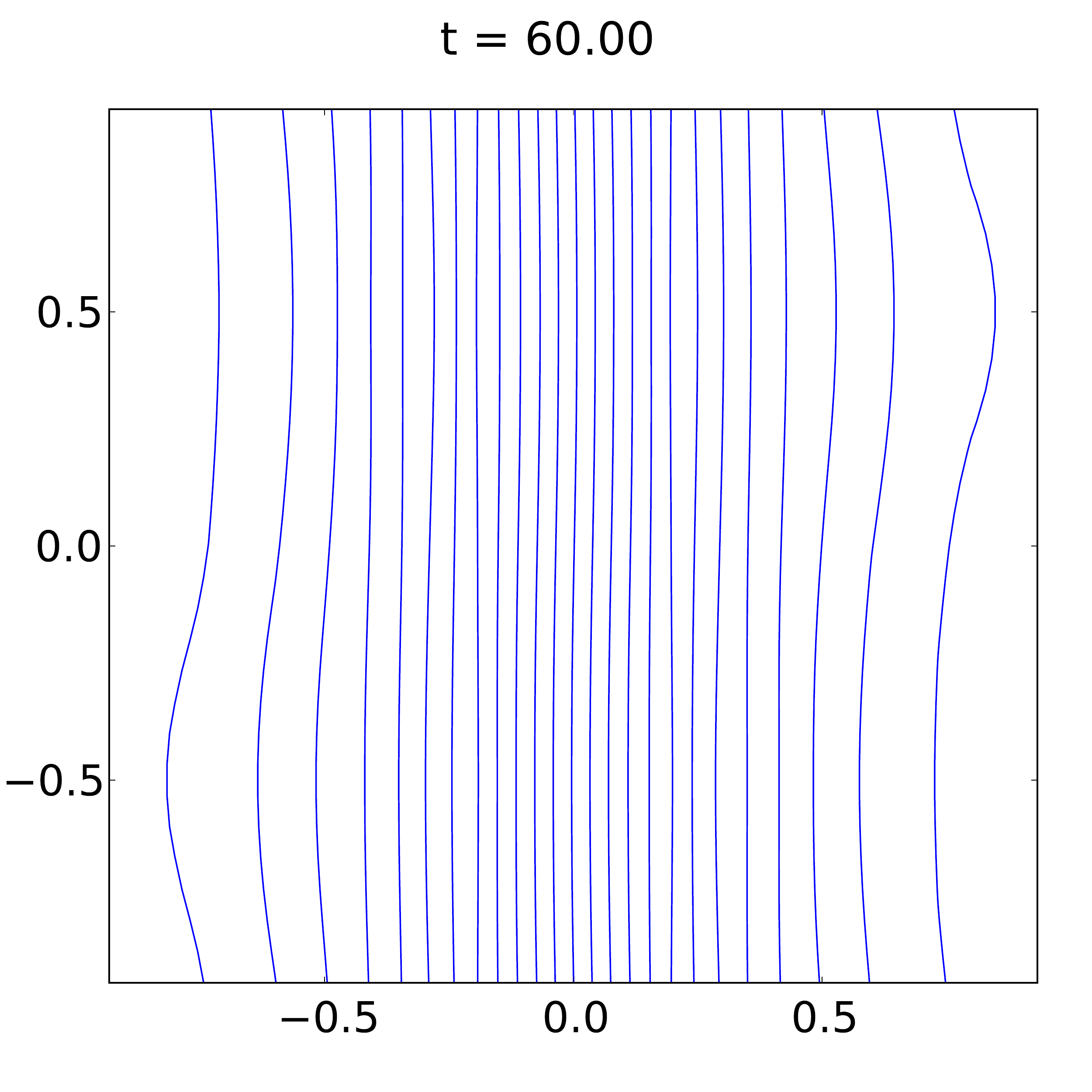}
}
\subfloat{
\includegraphics[width=.32\textwidth]{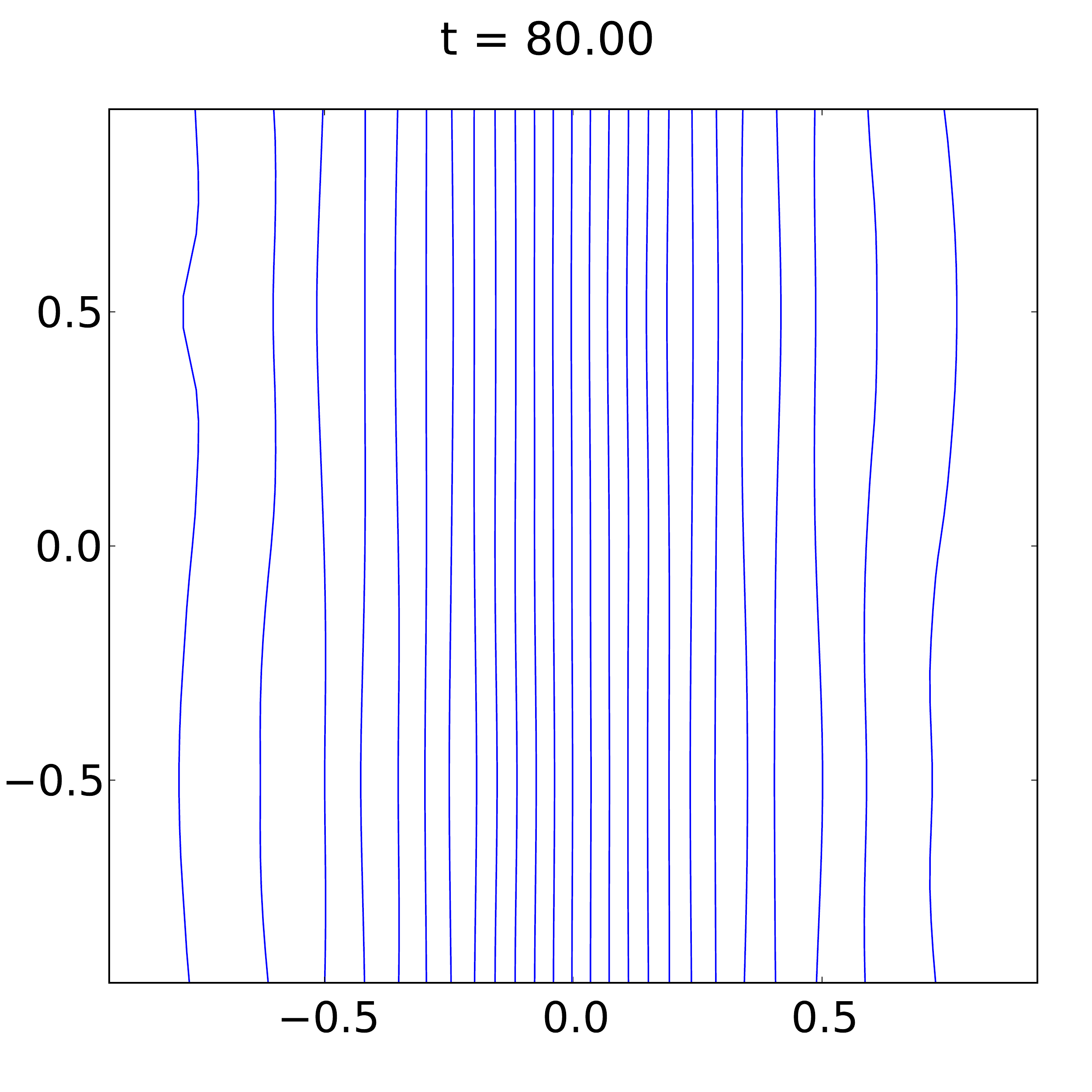}
}
\subfloat{
\includegraphics[width=.32\textwidth]{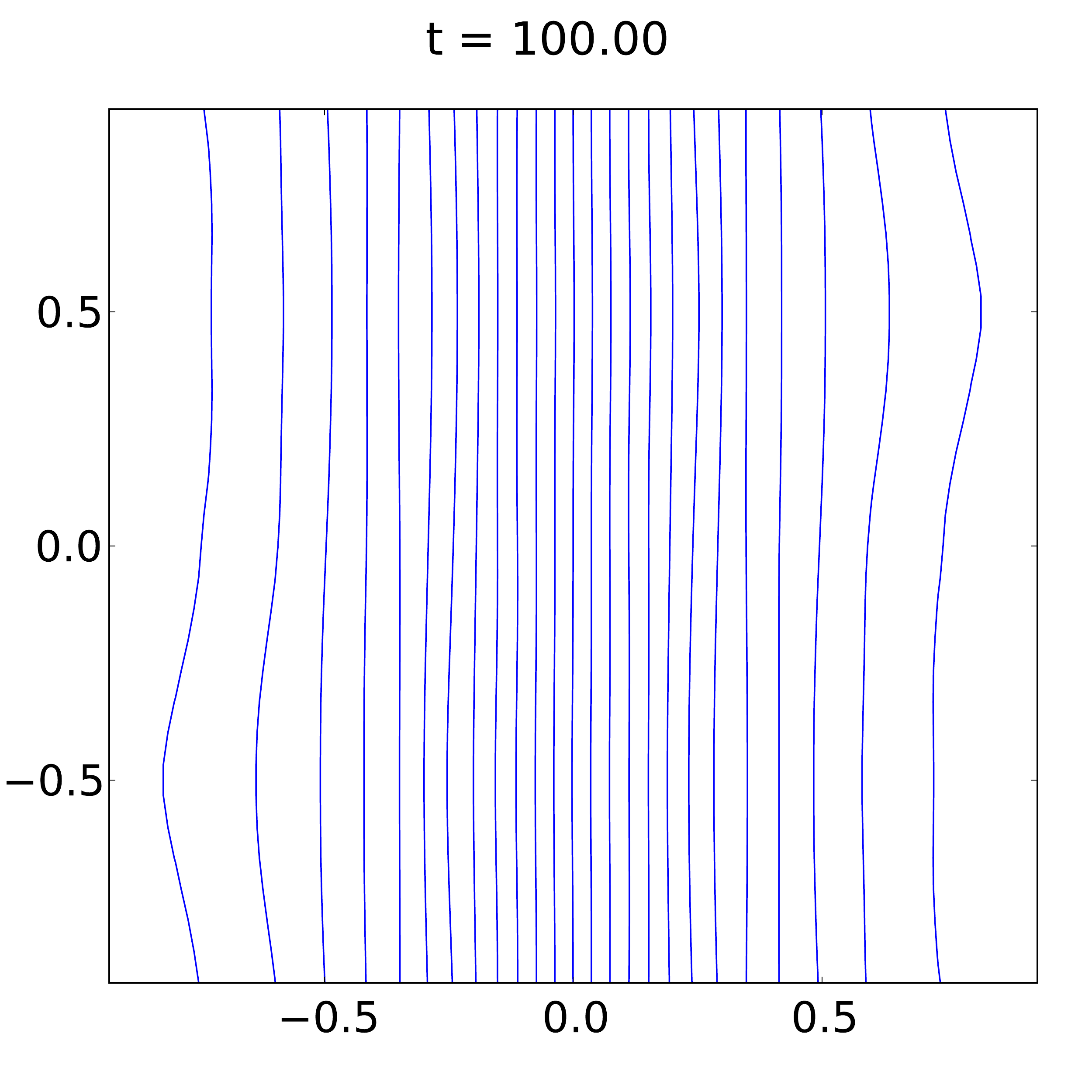}
}

\caption{Cosh current sheath, $30 \times 30$ grid points. Magnetic field lines.}
\label{fig:current_sheath_cosh_field_lines}
\end{figure}

\begin{figure}
\centering
\subfloat{
\includegraphics[width=.32\textwidth]{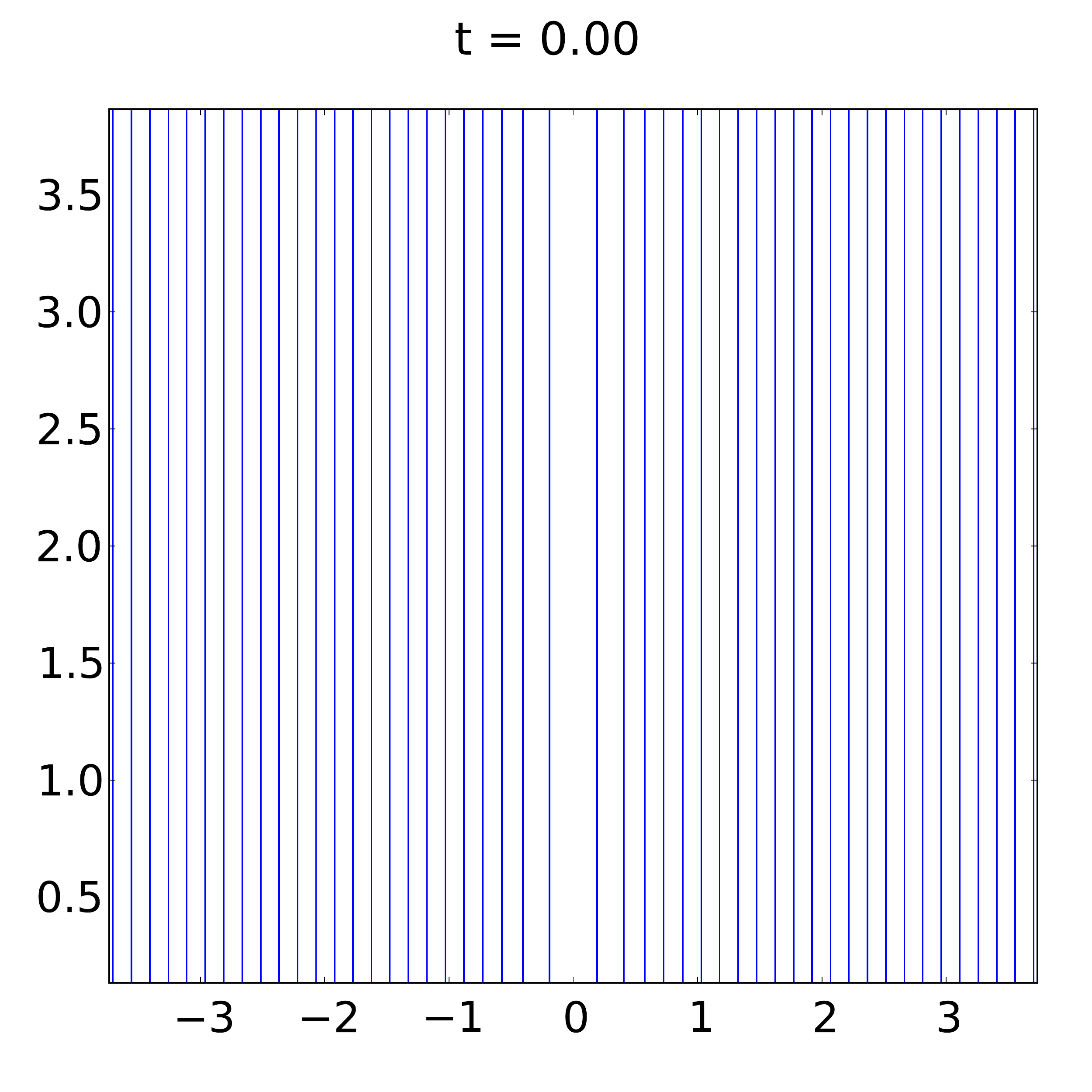}
}
\subfloat{
\includegraphics[width=.32\textwidth]{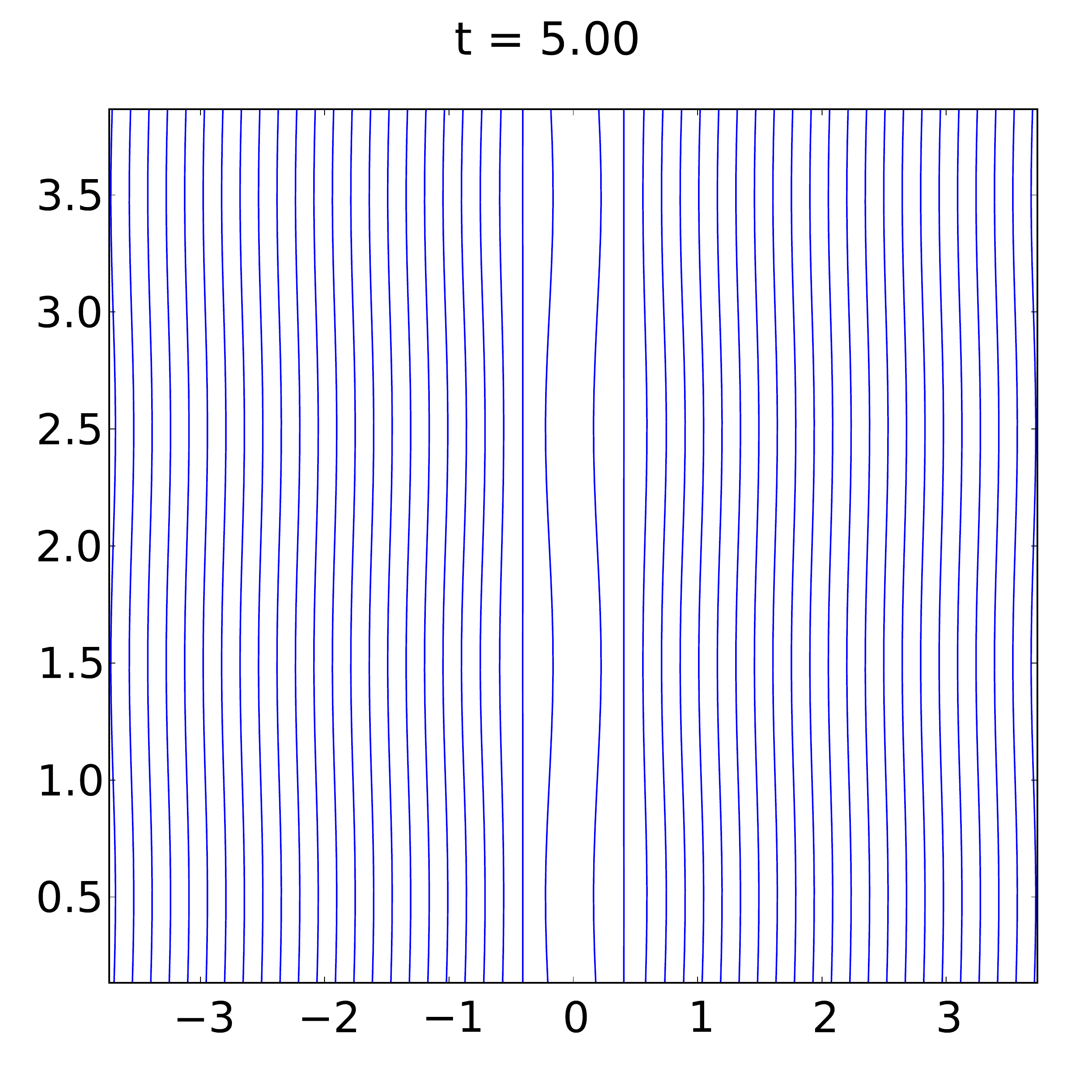}
}
\subfloat{
\includegraphics[width=.32\textwidth]{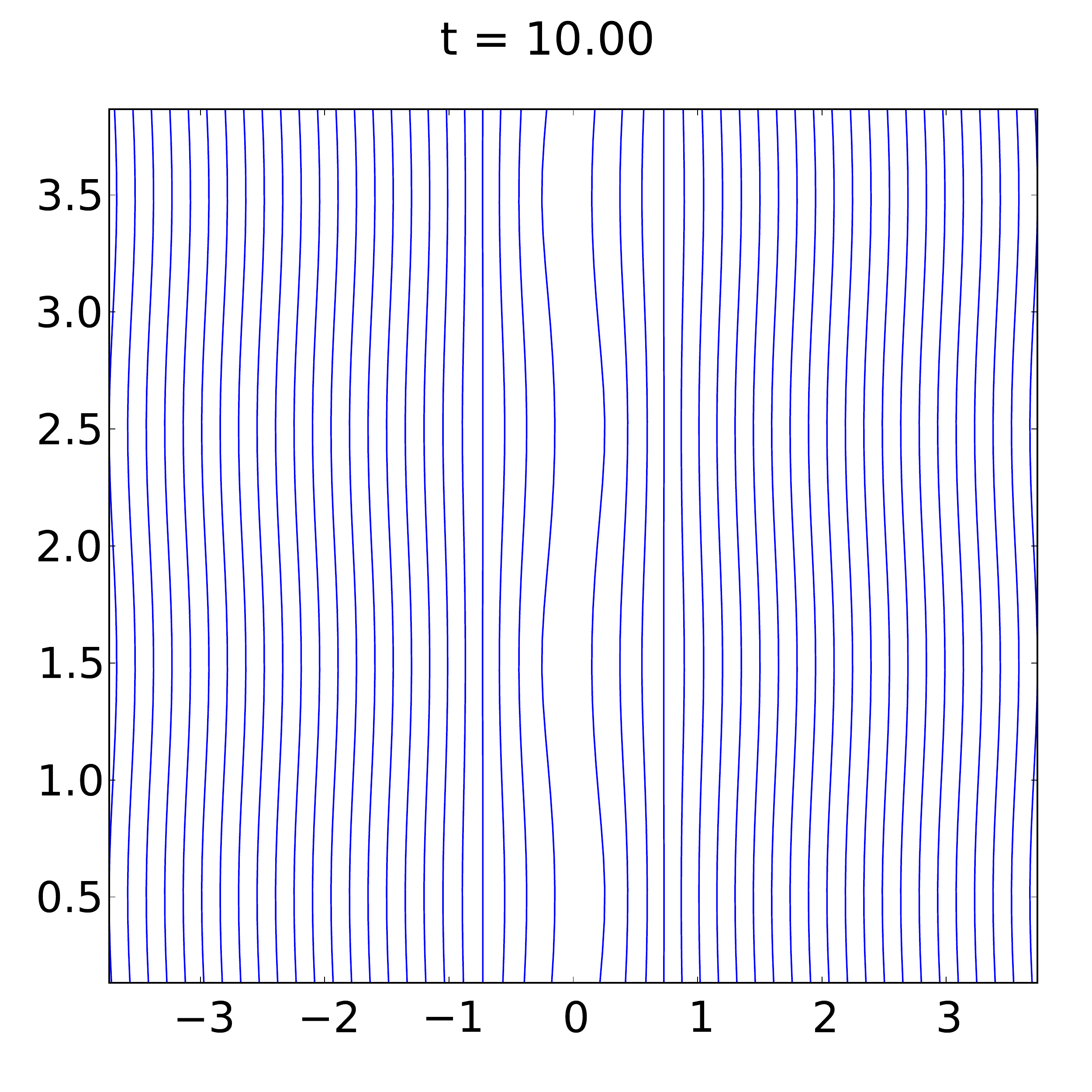}
}

\subfloat{
\includegraphics[width=.32\textwidth]{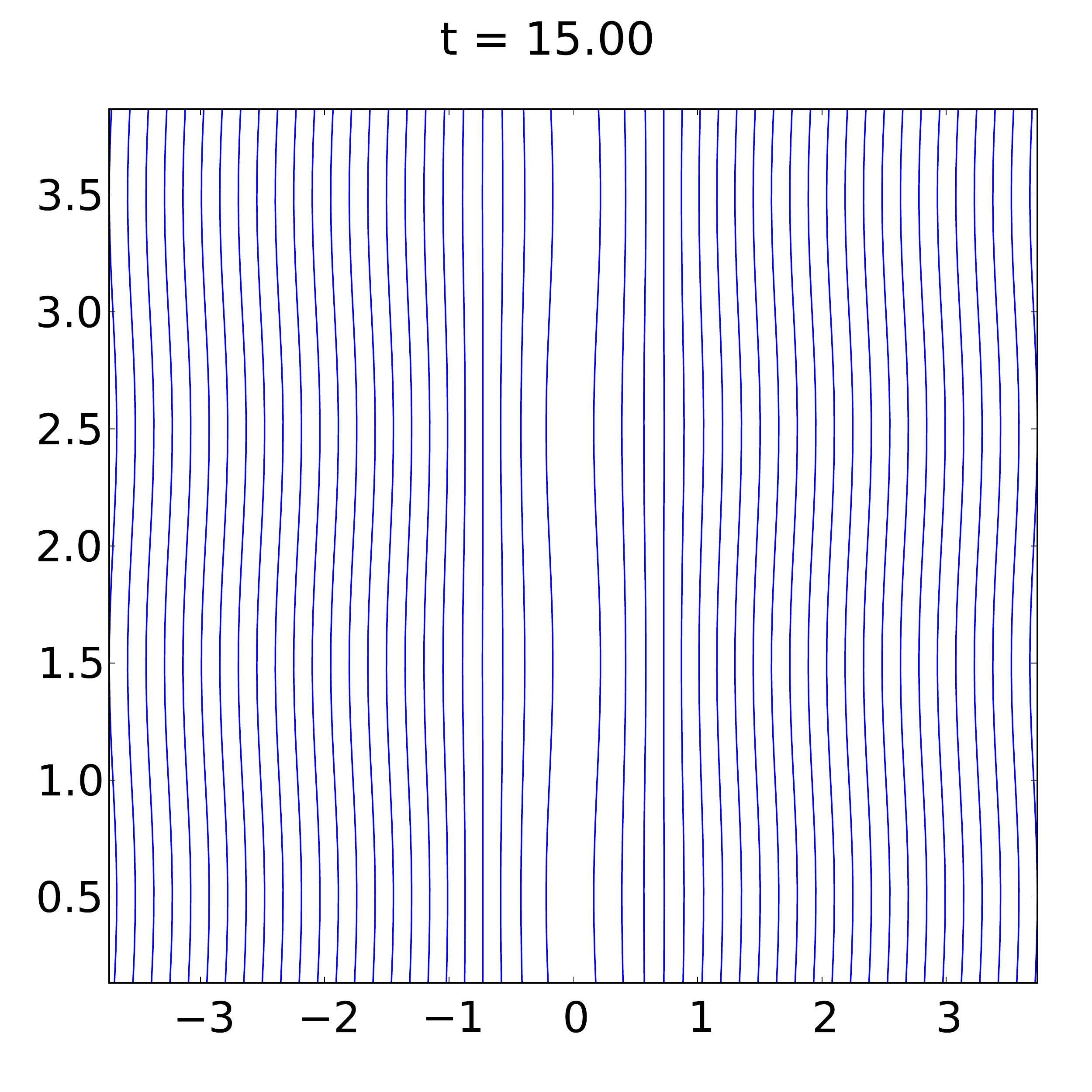}
}
\subfloat{
\includegraphics[width=.32\textwidth]{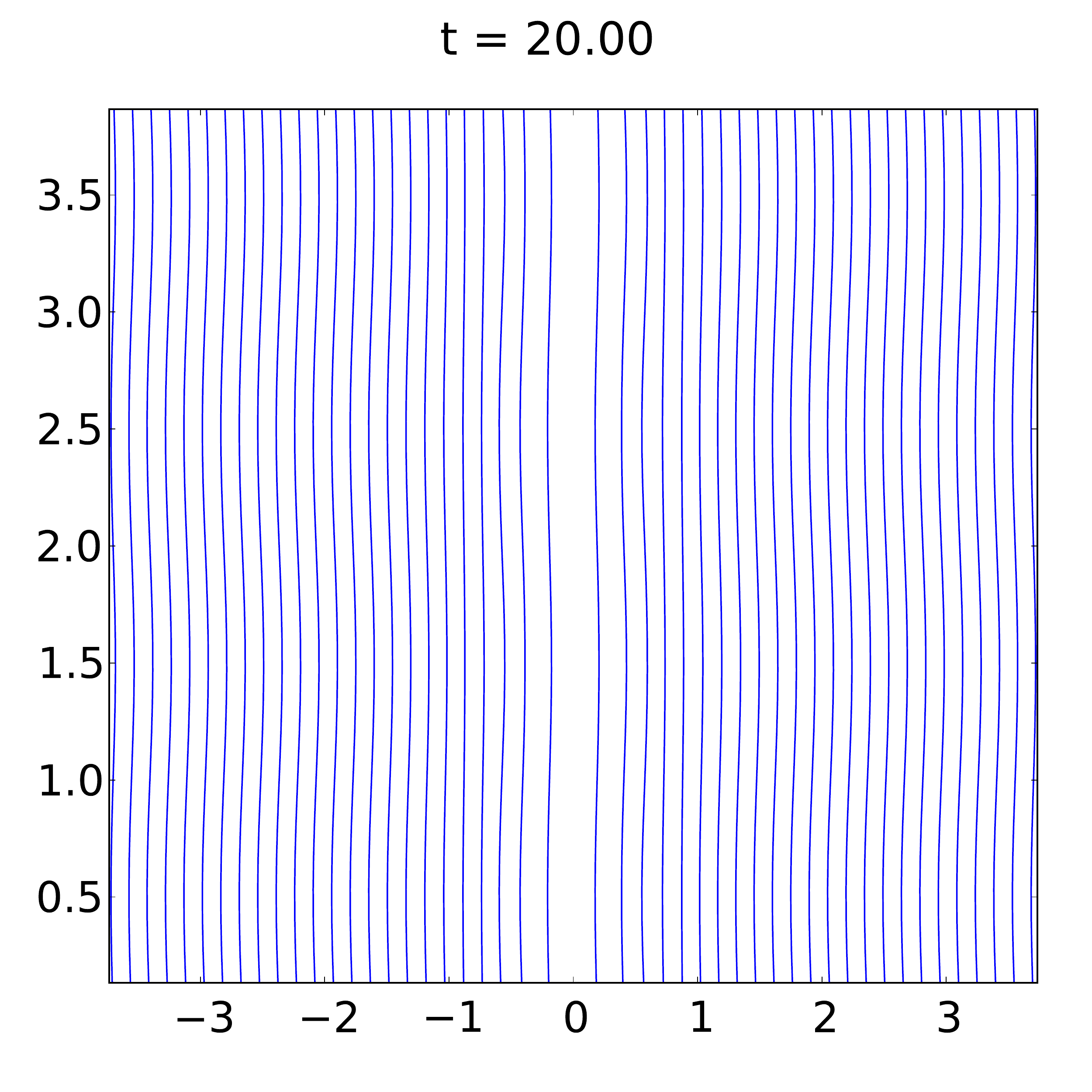}
}
\subfloat{
\includegraphics[width=.32\textwidth]{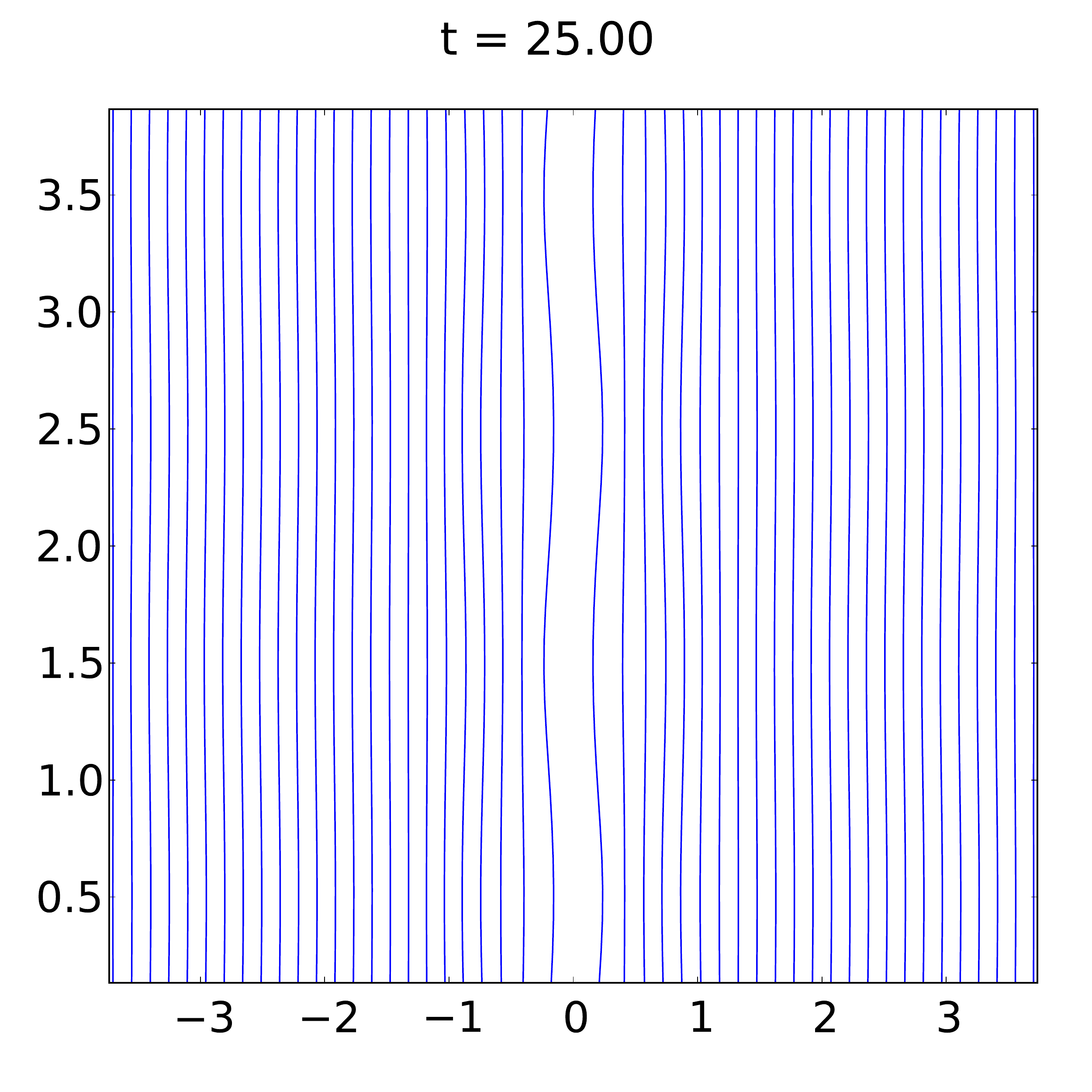}
}

\subfloat{
\includegraphics[width=.32\textwidth]{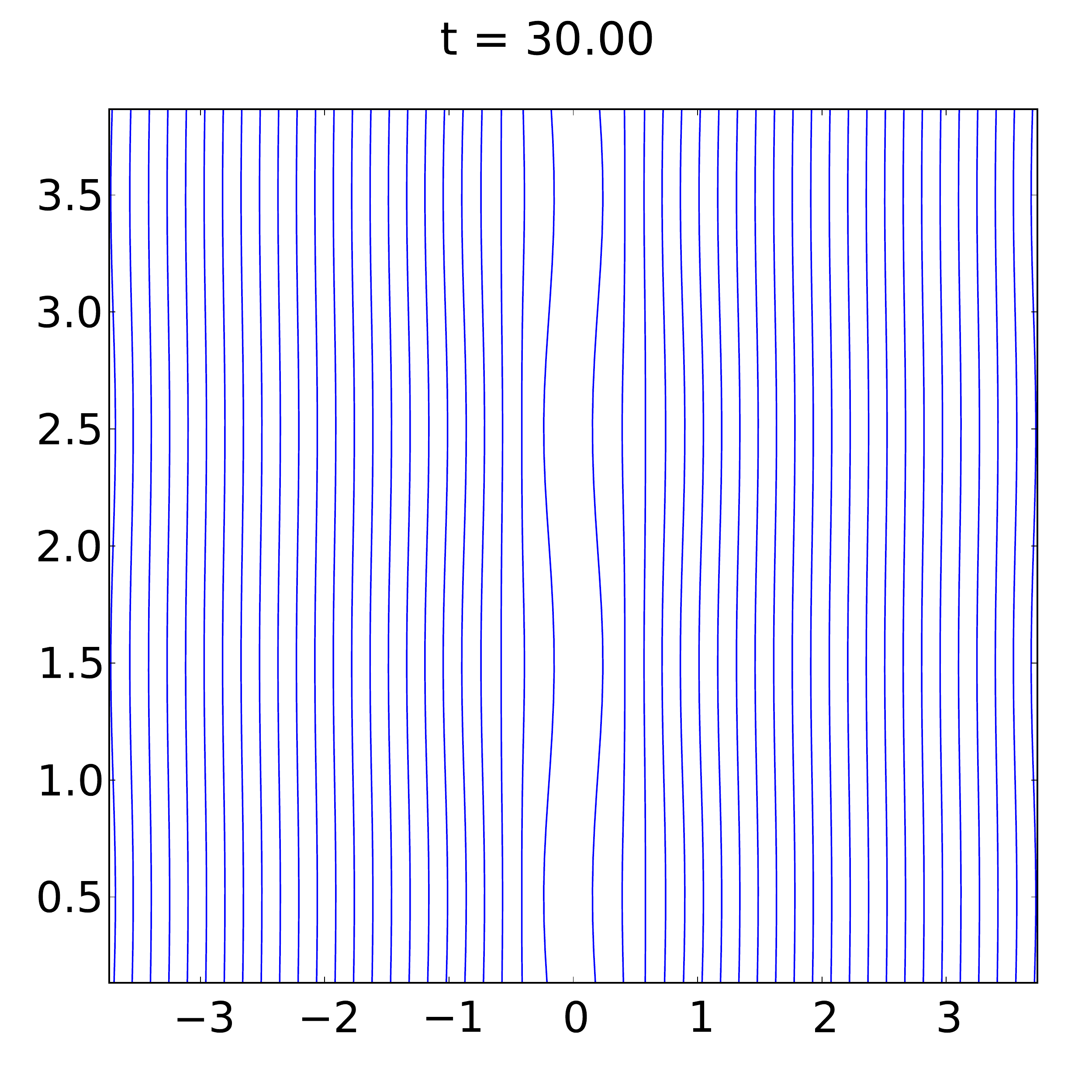}
}
\subfloat{
\includegraphics[width=.32\textwidth]{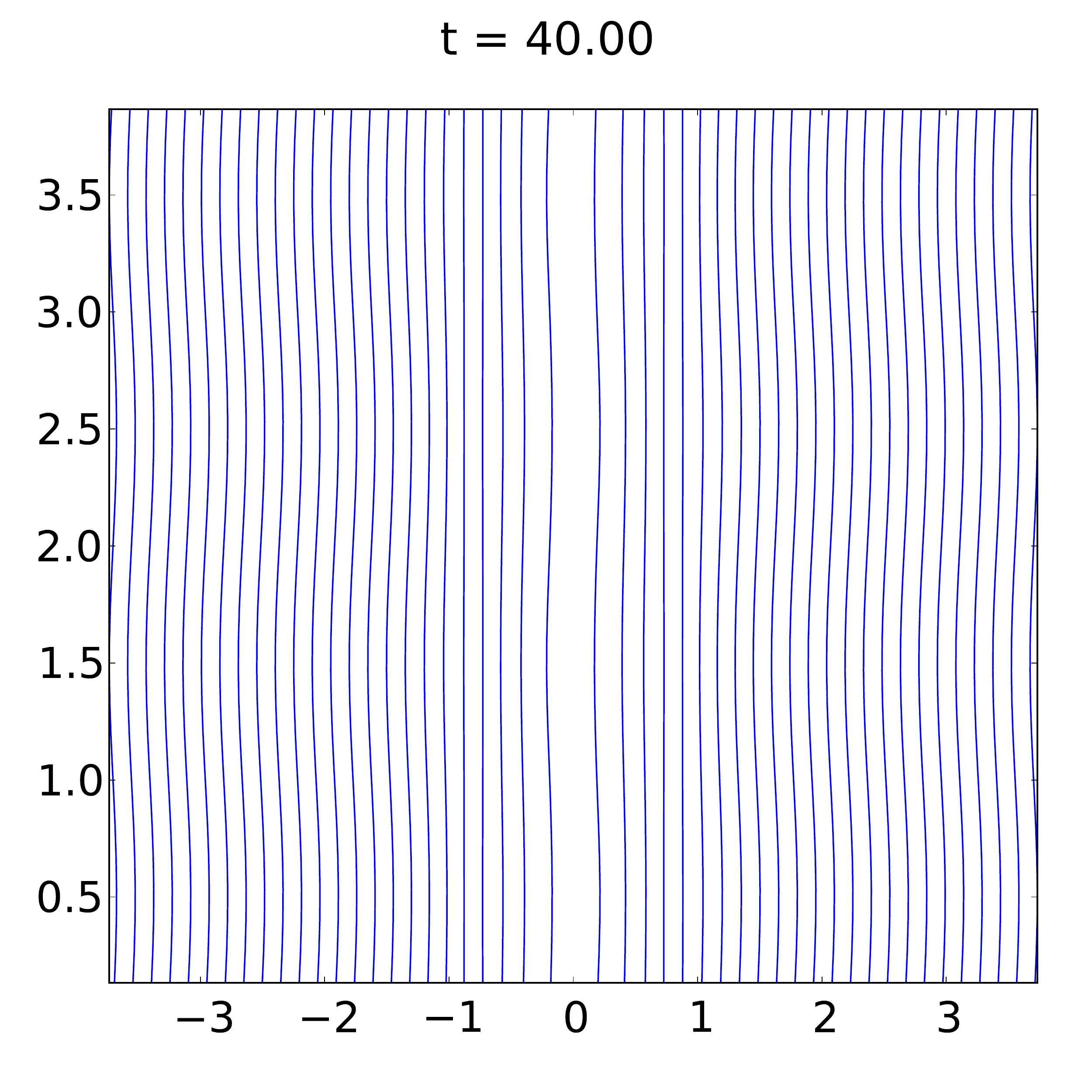}
}
\subfloat{
\includegraphics[width=.32\textwidth]{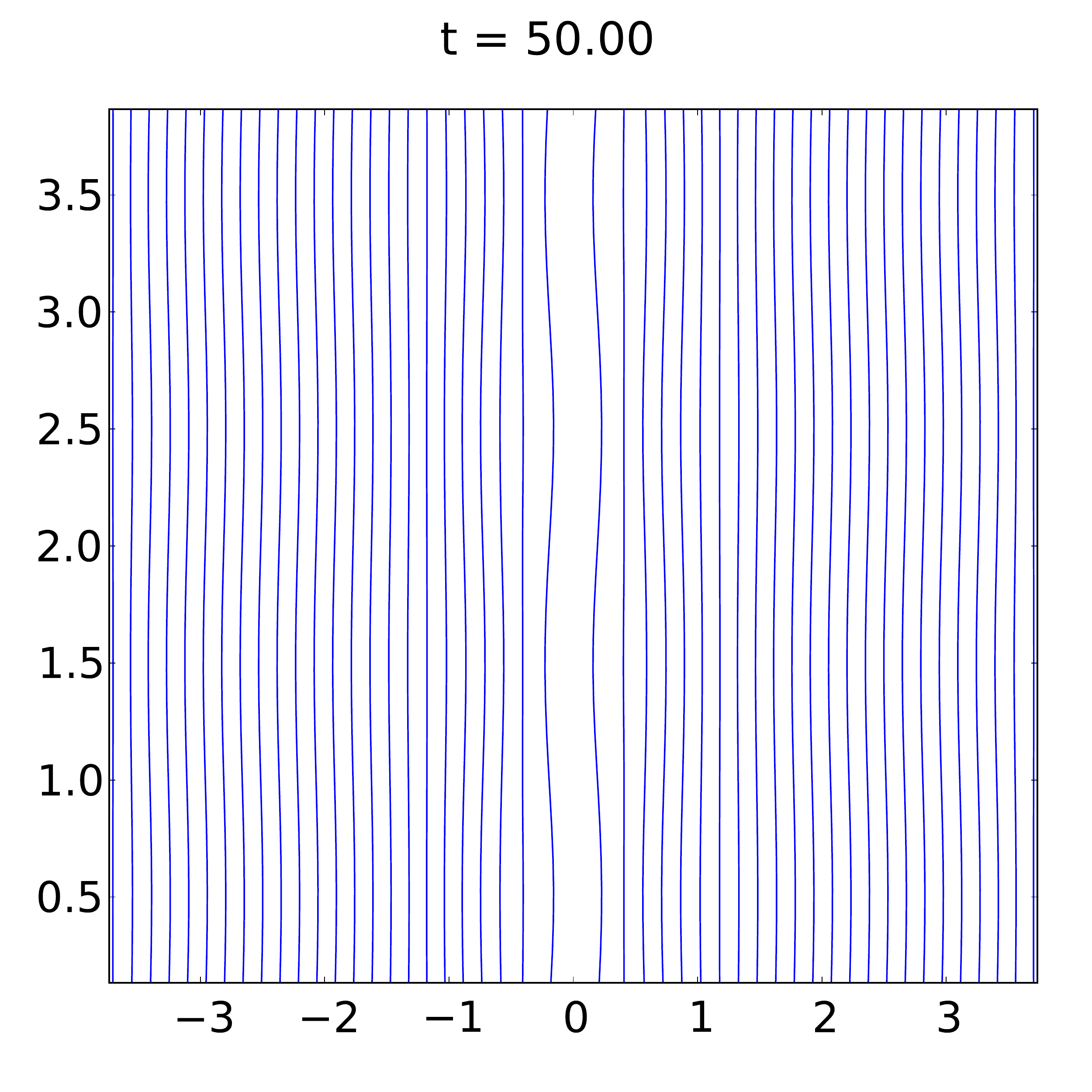}
}

\subfloat{
\includegraphics[width=.32\textwidth]{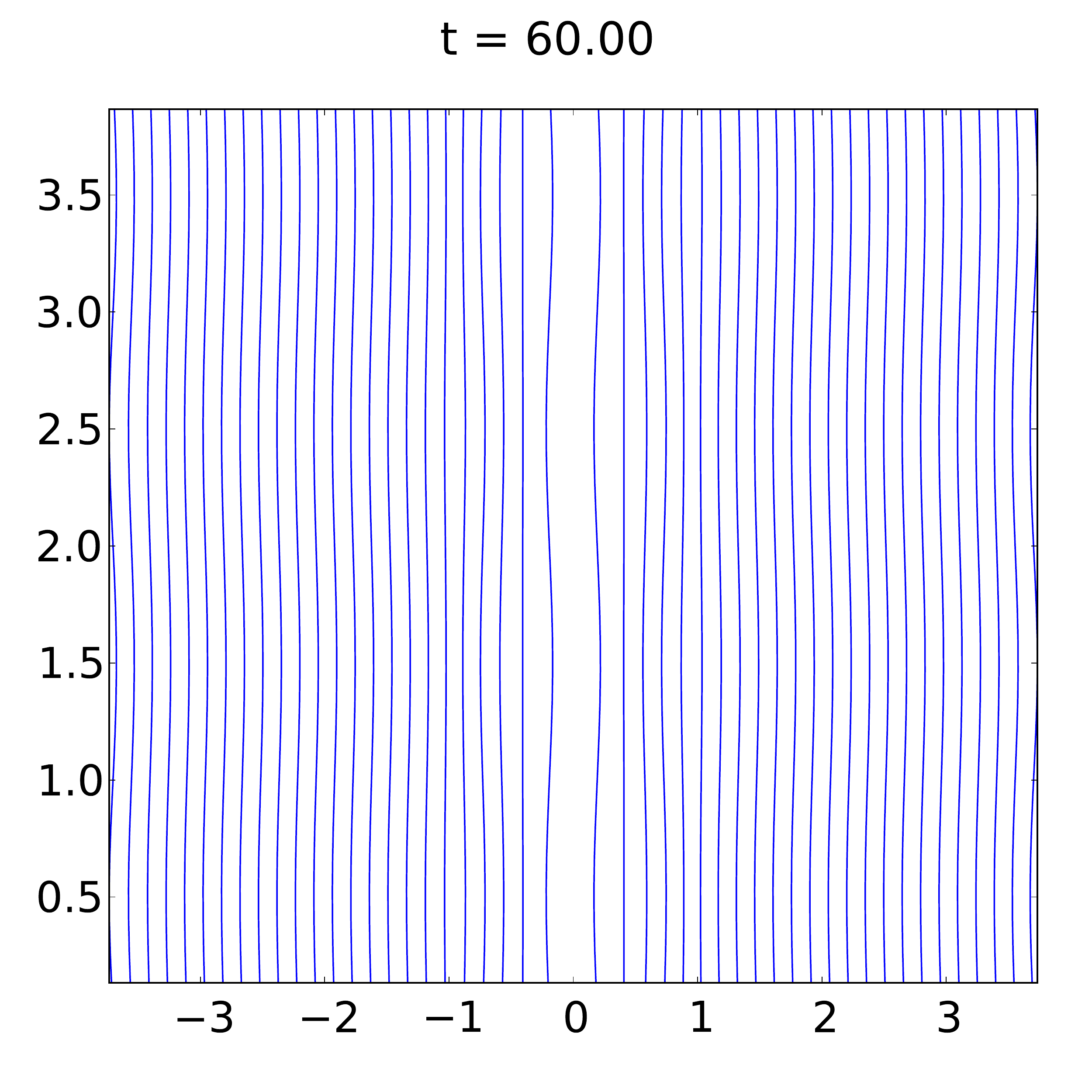}
}
\subfloat{
\includegraphics[width=.32\textwidth]{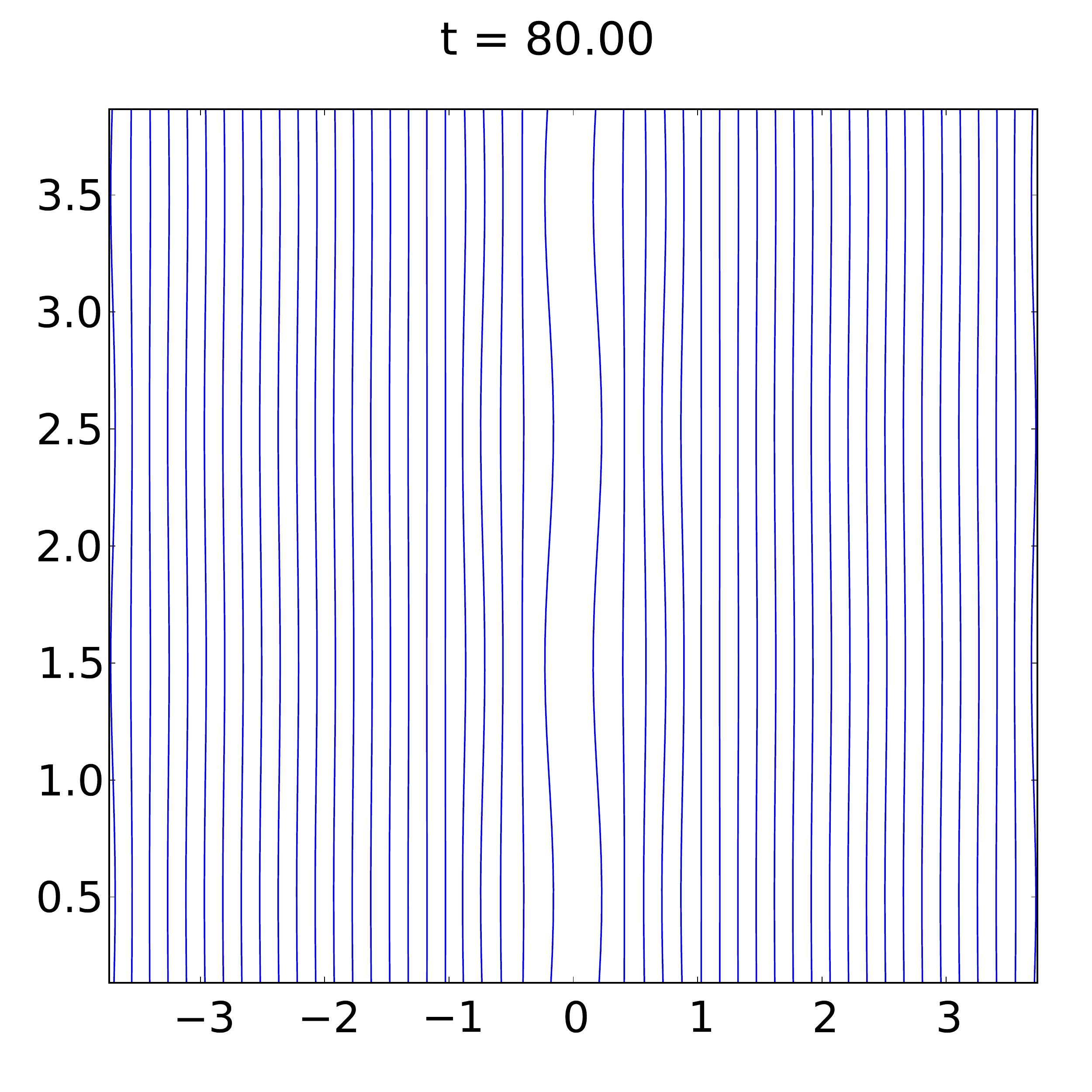}
}
\subfloat{
\includegraphics[width=.32\textwidth]{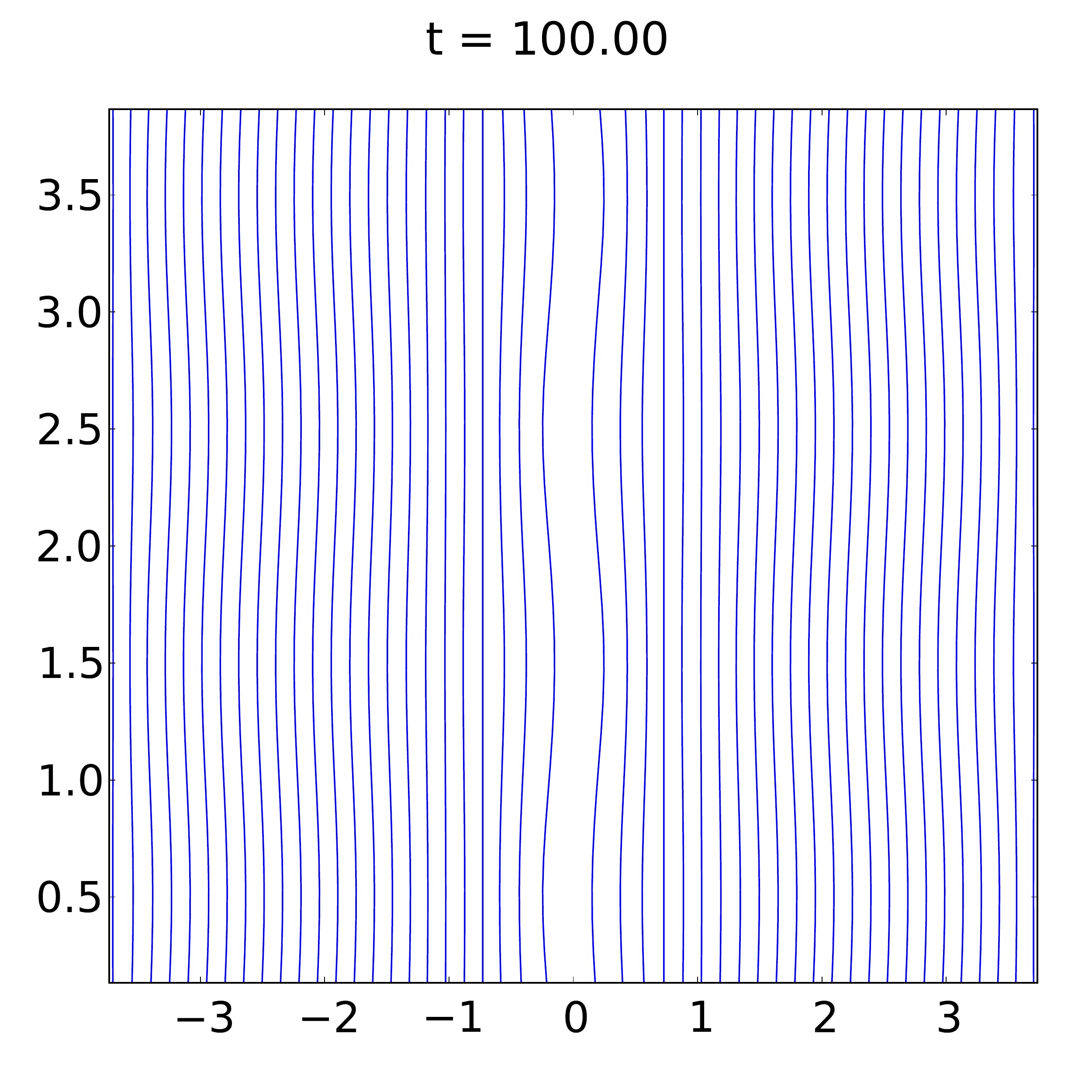}
}

\caption{Tanh current sheath, $30 \times 30$ grid points. Magnetic field lines.}
\label{fig:current_sheath_tanh_field_lines}
\end{figure}

We have to stress here, that this set of initial conditions is quite challenging for most numerical schemes due to the discontinuity and other methods break down much earlier.
To investigate the preservation of the magnetic field line topology on longer time scales, we consider therefore the two less severe examples of current sheaths, defined by $\cosh$ and $\tanh$ profiles, that are frequently used in reconnection studies \cite{GrassoCalifano:2001, GrassoTassi:2010, TassiMorrison:2010}. In both, the $\cosh$ and the $\tanh$ case, the magnetic field changes sign not suddenly but more smoothly.
Under this condition, we can run the simulation much longer. Both, figures \ref{fig:current_sheath_cosh_field_lines} and \ref{fig:current_sheath_tanh_field_lines}, show the field line evolution for the case of the $\cosh$ sheet and the $\tanh$ sheet, respectively, up to $t = 100$. We observe, that for a smooth magnetic field, magnetic field lines are only bend but do not reconnect, as is expected from the theory but rarely observed in numerical simulations, especially on the time scales we are considering here.
In our description of magnetohydrodynamics, reconnection can only occur if the resistive term, proportional to $\nabla^{2} B$, in (\ref{eq:mhd_eqs_B}) is present.  But in ideal MHD, $\eta = 0$, such that the topology of the magnetic fields lines is preserved. Most numerical schemes, however, do feature a certain amount of numerical resistivity, leading to unphysical reconnection. In the variational integrator, such spurious resistivity appears to be completely absent, at least in the case of a continuous magnetic field.

\newpage

\clearpage

\clearpage

\clearpage

\clearpage

\clearpage

\clearpage

\clearpage

\clearpage

\chapter{Summary and Outlook}

\section{Results}

In this work, it was successfully demonstrated that the application of variational integrators to models from plasma physics has great potential in the development of novel numerical algorithms that map the real world to the discrete more realistically than most standard methods.
In the following our results will be summarised.

\subsection{Theory}

Most systems from plasma physics, especially infinite dimensional ones (field theories), do not admit a natural Lagrangian formulation to which the variational integrator method is directly applicable.
To be able to apply variational integrators to plasma physics problems it was therefore necessary to extend the theory.
A combination of Ibragimov's theory of extended Lagrangians with the theory of variational integrators made it possible to use this method for systems not envisaged in the original theory. This includes several systems, like advection type problems, that are typical for plasma physics.
The impact of this result reaches far beyond the realm of the examples considered in this work, as the class of problems that is admissible to the variational integrator theory is now greatly enlarged.

An important feature of Ibragimov's theory is the recovery of Noether's theorem, which allows us to analyse the conservation properties of a system on the continuous as well as on the discrete level, thereby allowing us to find the exact discrete expressions of the quantities that are conserved by the variational integrators.

\subsection{Particle Dynamics}

In chapter four, previous work on the development of variational integrators for guiding centre dynamics in two dimensions was extended. A new integrator based on a midpoint discretisation was developed and compared with the previous one based on a trapezoidal discretisation.
The properties of both integrators are then evaluated against the standard explicit fourth order Runge-Kutta method. The midpoint integrator appears to be more accurate and more stable at smaller timesteps than the trapezoidal integrator. But both are by far superior to the Runge-Kutta method, for which the particle is found to severely deviate from its expected orbit.
The variational integrators exhibit a much better long term stability as the global error of the energy is limited and does not grow in the course of a simulation, i.e., there is no numerical dissipation. The particle orbit is preserved for millions of characteristic times and hundreds of millions of timesteps. Indeed, the energy is not constant but found to oscillate about a fixed value. The important point is that the amplitude of the oscillation stays constant during the whole simulation.

These integrators were then extended to four dimensions where an additional conserved quantity, the toroidal momentum, is present. The performance of the integrators with respect to energy conservation and accuracy of the particle orbit were similar to those of the two dimensional case. The error in the toroidal momentum is determined by the residual of the nonlinear solver.

All of the derived integrators are nonlinearly implicit such that iterative methods have to be used to solve the corresponding system of equations.
In the two dimensional case, an analytic solution of the linear system comprising each Newton step can be used with a fixed number of iterations, namely three plus an explicit initial guess, such that the computational effort is about the same as with the fourth order Runge-Kutta method.
In the four dimensional case, a matrix solver is used to solve the system of equations, which is computationally more demanding than the corresponding Runge-Kutta integrator, but most probably more efficient ways of solving this system could be found. That, however, is an issue that was not part of this work and is left for future research.

\subsection{Kinetic Theory}

In chapter five, several variational integrators for the Vlasov-Poisson system in one dimension (actually, one spatial plus one velocity dimension) were derived. This system can be regarded as the stepping stone to more complicated models like gyrokinetics, which is the backbone of nowadays' small scale turbulence simulations in plasma physics.

A fully nonlinear integrator was constructed, that has outstanding conservational properties, preserving the total particle number, the total energy, linear momentum and the $L^{2}$ norm up to machine precision. The only limit to the performance of this integrator is the grid resolution. If not sufficient, subgrid modes may develop, leading to large gradients in the distribution function, therefore altering the solution unphysically and spoiling the conservation properties. This is, however, a well known phenomenon with low-order finite difference methods like ours.
To counteract these effects, a collision operator was constructed, that dissipates the $L^{2}$ norm and damps the subgrid modes, thereby retaining the conservation of the total particle number, energy and momentum.
It was shown that with the collision operator the variational integrator allows for long time simulations of different standard test cases from kinetic theory, linear and nonlinear ones, for example Landau damping and the twostream instability.

Furthermore, a linearly implicit integrator was constructed, which is computationally less demanding but retains most of the conservation properties, even though it does not reach the accuracy in the errors of the conserved quantities as the nonlinear integrator.
It is important to note that the linearisation is implemented at the level of the Lagrangian. It is not the discrete Euler-Lagrange equations that are linearised. This is a crucial point, as it guarantees that the properties of variational integrators, conservation of discrete momenta and the multisymplectic form, are kept intact.
What is lost is exact energy preservation. Instead, an oscillating energy error is found, as it is typical for symplectic and multisymplectic methods. The amplitude of the oscillation depends largely on the strength of the nonlinearity as well as on the timestep.
In most cases, however, the linearised system still yields very good results. The energy error is still bounded and there is still no numerical dissipation.
This integrator can be used to compute an initial guess for an iterative solver of the nonlinear integrator but also standalone, if the solution of a nonlinearly implicit scheme is computationally too expensive.

Albeit implementing efficient solution techniques for the derived schemes was not an objective of this work, first results suggest that efficient iterative methods (GMRES) can be used effectively if a sufficiently accurate initial guess is available.
If the linear integrator could be solved efficiently with the GMRES method, the computational effort could be in the vicinity of explicit methods.

\subsection{Fluid Dynamics}

In chapter six, a variational discretisation of a fluid plasma model, namely ideal incompressible magnetohydrodynamics, was developed and analysed. The ideal model shares many features with more sophisticated and physically comprehensive descriptions used in simulations of large scale turbulence or magnetic reconnection, to name just two examples.

The scheme is similar to the one for the Vlasov-Poisson system in that it is fully nonlinear and exhibits excellent conservation properties. Both, the total energy and the cross helicity are preserved up to machine precision over long simulation times. Nonlinear Alfvénic waves, a core phenomenon of magnetohydrodynamics, were found to keep travelling for thousand characteristic times without loosing energy. In turbulent settings, where small scale structures can develop, problems with subgrid modes were observed if the resolution was not sufficient, similar to the Vlasov-Poisson integrator. As in realistic simulations usually some form of dissipation is present in the model, this is not expected to be a problem impairing the applicability of the derived scheme.
In these cases, the use of a conservative scheme can ensure that the dissipation mechanism is consistent with the physical dissipation and not numerically generated.

The analysis of several reconnection models has shown the preservation of magnetic field line topology by the variational integrator for very long times, provided that the magnetic field is continuous. For a discontinuous magnetic field, the topology is still preserved for times longer than with most standard methods, but at some point the solution was found to become spurious, probably due to error accumulation. This is, however, not discouraging as the variational integrator was not constructed with discontinuities in mind for which usually specially designed methods are employed.

\subsection{Semi-Discretisations}

In the appendix we describe two different approaches of using variational integrators or closely related methods to obtain semi-discretisations for the Vlasov equation and similar systems.

In appendix A, the spatial dimensions are transformed into Fourier space and only time or time and the velocity dimensions are discretised variationally.
In appendix B, discretisation methods for different kinds of brackets are presented, where only the phasespace dimensions are discretised but not time. A new discretisation for Nambu three brackets is derived, that has potential use in simulations of the gyrokinetic Vlasov equation in axisymmetric tokamaks.
Furthermore, it is shown that discretisation strategies for Lie-Poisson brackets are applicable also to problems from plasma physics.
Interestingly, these bracket discretisations are quite similar to variational integrators, e.g., for the Vlasov-Poisson system.

\section{Future Work}

This work has shown that the application of variational integrators in plasma physics is both viable and attractive. Thereby it has lead the way to a number of possible directions for future research, some of which are sketched below.

\subsection{Theory}

The theory of variational integrators is still young and therefore not fully developed. So far, there are no clear guidelines on which discretisations of the Lagrangian lead to ``good'' (e.g., stable) numerical schemes.
A possible solution to this might be the combination of the present methods with Arnold's theory of finite element exterior calculus (more on this below).

In general, the variational integrator method does not lead to schemes that are (multi)symplectic, momentum and energy preserving at the same time, without implementing one or another form of timestep adaption, either global timestep adaption or asynchronous variational integrators \cite{Kane:1999, Lew:2003}. However, the methods we derived for the Vlasov-Poisson system as well as for ideal magnetohydrodynamics turned out to be exactly energy preserving (to machine precision), momentum preserving (also to machine precision), and they are multisymplectic by design. It would be interesting and important to understand such properties and find general criteria for which discretisations of the Lagrangian lead to such optimal integrators.

\subsubsection{Finite Element Discrete Exterior Calculus}

Recently, a discrete theory of exterior calculus based on finite elements has been developed by \citeauthor{Arnold:2006} \cite{Arnold:2006, Arnold:2010}.
In contrast to the theory of variational integrators, this theory of finite element exterior calculus is embedded in an abstract Hilbert space framework that makes the analysis of stability and convergence of the derived discretisations much more systematic.
A connection between these two theories might allow for a rigorous numerical analysis of variational integrators and therefore help to find answers to some of the questions raised in the previous paragraph.

\subsubsection{Analysis of Discrete Conservation Laws and Preservation of the Multisymplectic Form}

To judge the performance of the variational integrators derived in chapters five and six, only heuristic diagnostics were used.
Instead, a detailed analysis of the discrete conservation laws using the discrete Noether theorem of section \ref{sec:vi_infinite_noether_theorem} should be carried out to find exact expressions of the conserved quantities.

It was shown how the solution space of a system described by an extended Lagrangian (section \ref{ch:classical_extended_lagrangians}) can be restricted to the solution space of the original system to recover the conserved quantities.
An open problem is to find out if the multisymplectic structure of the extended system endows the physical system with a compatible multisymplectic structure and if so how the structure of the extended system can be restricted to obtain the one of the physical system.
The variational integrators preserve a discrete counterpart of the multisymplectic structure by construction, but so far no statement about the conservation of the multisymplectic structure of the original system is possible.

\subsection{Vlasov-Poisson and Vlasov-Maxwell}

Some obvious extensions of the variational integrator for the Vlasov-Poisson system from chapter five include higher dimensional domains, the Vlasov-Maxwell system, and higher order schemes. Not many surprises are expected in deriving variational integrators for the Vlasov-Poisson system in higher dimensions, especially extending the integrator of the Vlasov equation should prove to be straight forward.
Treating the electrodynamic equations of the Vlasov-Maxwell system properly might however necessitate the use of a staggered grid approach \cite{Yee:1966} and therefore suggests a treatment based of discrete exterior calculus \cite{Stern:2009a, Stern:2009b}.

The derivation of higher order schemes as well appears to be a straight forward generalisation of the results presented in this work. The most interesting question in this respect is if Arakawa's fourth order discretisation of the Poisson brackets can be derived by a variational method similarly to his second order discretisation.

\subsubsection{Nambu Bracket Discretisation}

In appendix B, it was shown that the Lie-Poisson bracket formulation of the Vlasov equation can be used to obtain semi-discretisations of the phasespace part of the Vlasov equation by first transforming them to Nambu brackets.
A Lie-Poisson bracket also exists for the Vlasov-Maxwell system \cite{MarsdenRatiu:2002}
\begin{align}
\{ F, G \}
\nonumber
&= \int f \, \left[ \dfrac{\delta F}{\delta f} , \dfrac{\delta G}{\delta f} \right] \, dx \, dv
+ \int \left( \dfrac{\delta F}{\delta E} \cdot \left( \nabla \times \dfrac{\delta G}{\delta B} \right) - \dfrac{\delta G}{\delta E} \cdot \left( \nabla \times \dfrac{\delta F}{\delta B} \right) \right) dx \, dv \\
&+ \int \left( \dfrac{\delta F}{\delta E} \cdot \dfrac{\delta f}{\delta v} \dfrac{\delta G}{\delta f} - \dfrac{\delta G}{\delta E} \cdot \dfrac{\delta f}{\delta v} \dfrac{\delta F}{\delta f} \right) dx \, dv
+ \int f \, B \cdot \left( \dfrac{\partial}{\partial v} \dfrac{\delta F}{\delta f} \times \dfrac{\partial}{\partial v} \dfrac{\delta G}{\delta f} \right) dx \, dv ,
\end{align}

where $f$ is the distribution function, $E$ and $B$ the electric and magnetic fields, respectively, and $H$ is the Hamiltonian energy functional, given by
\begin{align}
H (f, E, B)
= \dfrac{1}{2} \int \norm{v}^{2} \, f (t, x, v) \, dx \, dv
+ \dfrac{1}{2} \int \Big( \norm{E(t,x)}^{2} + \norm{B(t,x)}^{2} \Big) \, dx ,
\end{align}

such that the evolution of all functionals $F(f, E, B)$ is determined by
\begin{align}
\dot{F} = \{ F, H \} .
\end{align}

Preliminary results suggest that this Lie-Poisson bracket has a Nambu bracket formulation as well which could be used to derive semi-discretisations for the Vlasov-Maxwell system, where in contrast to the Vlasov-Poisson system not only the distribution function but also the electric and magnetic fields obey dynamical equations.

\subsubsection{Dirac Bracket Discretisation}

Recently, there has been some effort to derive Dirac brackets for several models of plasma physics \cite{Chandre:2012a, Chandre:2012b}, including the Vlasov-Poisson and Vlasov-Maxwell systems.
These brackets are constructed to respect constraints of the dynamics automatically, e.g. $\Delta \phi - \int f \, dv = 0$ for the Vlasov-Poisson system or $\nabla \cdot V = 0$ for incompressible fluids.
The Dirac brackets for the Vlasov-Poisson system,
\begin{align}
\{ F, G \}_{*} (f, E) &= \int f \, \bigg[ \dfrac{\delta F}{\delta f} - \Delta^{-1} \nabla \cdot \dfrac{\delta F}{\delta E} , \dfrac{\delta G}{\delta f} - \Delta^{-1} \nabla \cdot \dfrac{\delta G}{\delta E} \bigg] \, dx \, dv ,
\end{align}

can be transformed into Lie-Poisson brackets,
\begin{align}
\{ F, G \} (f, \Phi) &= \int f \, \bigg[ \dfrac{\delta F}{\delta f} - \dfrac{\delta F}{\delta \Phi} , \dfrac{\delta G}{\delta f} - \dfrac{\delta G}{\delta \Phi} \bigg] \, dx \, dv,
\end{align}

where
\begin{align}
\Phi = \Delta \phi = - \nabla \cdot E .
\end{align}

Again it might be possible to find a relation with Nambu brackets and use that relation to obtain semi-discretisations.
An interesting peculiarity of this formulation is that due to the Poisson equation being included as a constraint two dynamical equations, the usual one for the distribution function $f$ and another one for the potential vorticity $\Phi$, have to be solved. This might turn out to be an advantage for parallel implementations of the resulting methods.

\subsubsection{Euler-Poincaré Action Principle}\label{sec:outlook_euler_poincare}

The Euler-Poincaré action principle reviewed in section \ref{sec:kinetic_theory_euler_poincare} might pose the starting point for the derivation of variational integrators in a similar way as was done by \citeauthor{Pavlov:2011} \cite{Pavlov:2011, Pavlov:2009} and \citeauthor{Gawlik:2011} \cite{Gawlik:2011} for incompressible fluids.
As already noted by \citeauthor{Squire:2013} \cite{Squire:2013} this might prove nontrivial as, e.g., a discretisation of the group of symplectomorphisms, which describes the dynamics of the Vlasov-Poisson and Vlasov-Maxwell systems, has to be found (for a more thorough discussion see \cite{Squire:2013}).

Despite the possible difficulties it seems worthwhile to pursue this path as the Euler-Poincaré reduced system appears to be the most natural, most geometric formulation of the family of Vlasov systems known to date.

\subsubsection{Gyrokinetics}

From the point of view of applying variational integrators in large scale plasma physics codes, the kinetic model appears to be less attractive than the gyrokinetic model. Gyrokinetics is a version of kinetic theory reduced to five phasespace dimensions instead of six, thereby lowering the computational burden.

An extension of the variational integrators derived in this work to gyrokinetics is thus an important topic of future work, overlapping with the extension to higher dimensions and the Vlasov-Maxwell system.

\subsection{Ideal and Reduced MHD}

The most obvious extensions of the variational integrator for magnetohydrodynamics from chapter six are similar to those of the Vlasov-Poisson case: moving to three dimensional domains, higher order discretisations, and more comprehensive models. Some others are shortly explained below.

\subsubsection{Potential Formulation}

In studies of magnetic reconnection, the potential formulation presented in section (\ref{sec:mhd_potential}) is very popular.
To derive variational integrators for this formulation, an extended version of the theory for second or third order field theories has to be applied \cite{KouranbaevaShkoller:2000, Kouranbaeva:1999}. Apart from this there appear to be no obvious obstacles.

It would be quite interesting to compare the performance of variational integrators for the two different approaches, i.e., the description in terms of the potentials $A$ and $\psi$ and the description in terms of the fields $B$ and $V$.

\subsubsection{Nambu and Dirac Brackets}

Several flavours of magnetohydrodynamics can also be described by Lie-Poisson brackets, e.g., for the ideal case in potential formulation that is
\begin{align}
\{ F, G \} (A, \omega) = \int A \, \bigg( \bigg[ \dfrac{\delta F}{\delta A} , \dfrac{\delta G}{\delta \omega} \bigg] - \bigg[ \dfrac{\delta F}{\delta \omega} , \dfrac{\delta G}{\delta A} \bigg] \bigg) \, dx \, dy + \int \omega \, \bigg[ \dfrac{\delta F}{\delta \omega} , \dfrac{\delta G}{\delta \omega} \bigg] \, dx \, dy ,
\end{align}

with the Hamiltonian energy functional
\begin{align}
H (A, \omega) = \dfrac{1}{2} \int \Big( \psi \omega - A j \Big) \, dx \, dy .
\end{align}

Here, $A$ is the magnetic vector potential, $j = - \Delta A$ the current density, $\omega = - \Delta \psi$ is the vorticity, and $\psi$ the streaming function.
With the Casimir invariant
\begin{align}
C = \int A \omega \, dx \, dy ,
\end{align}

this can be transformed into a Nambu bracket
\begin{align}
\{ F, G, C \} (A, \omega) = \int \dfrac{\delta C}{\delta \omega} \, \bigg( \bigg[ \dfrac{\delta F}{\delta A} , \dfrac{\delta G}{\delta \omega} \bigg] - \bigg[ \dfrac{\delta F}{\delta \omega} , \dfrac{\delta G}{\delta A} \bigg] \bigg) \, dx \, dy + \int \dfrac{\delta C}{\delta A} \, \bigg[ \dfrac{\delta F}{\delta \omega} , \dfrac{\delta G}{\delta \omega} \bigg] \, dx \, dy .
\end{align}

A discretisation approach as described in appendix \ref{sec:brackets_discrete_nambu} should be straight forwardly applicable to this expression.

As for the family of Vlasov systems, there has been recent research on Dirac bracket formulations for magnetohydrodynamics as well \cite{Chandre:2012a, Chandre:2012b}, which might open new possibilities for deriving discretisations if a relation between Dirac brackets and Nambu brackets could be drawn as outlined above.

\appendix

\chapter{Mixed Spectral-Variational Schemes}

In this appendix, we sketch how to derive mixed spectral-variational methods for the Vlasov-Poisson system.
We do the derivation for a system very similar to Vlasov-Poisson, but somewhat simpler, namely the vorticity equation in two spatial dimensions.
The generalisation to the Vlasov-Poisson equation, e.g., in two spatial and two velocity dimensions is mostly straight forward.

\section{The Vorticity Equation in 2D}

The vorticity equation describes the evolution of the vorticity of a fluid element in an incompressible ideal fluid
\begin{align}\label{eq:spectral_vorticity_1}
\dfrac{\partial \omega}{\partial t} + \dfrac{\partial \omega}{\partial x} \dfrac{\partial \psi}{\partial y} - \dfrac{\partial \omega}{\partial y} \dfrac{\partial \psi}{\partial x} = 0 ,
\end{align}
where $\omega$ is the vorticity of the fluid and $\psi$ is the streaming function, determined by
\begin{align}\label{eq:spectral_vorticity_2}
- \Delta \psi = \omega .
\end{align}

In two dimensions, the vorticity equations takes a particularly interesting form which has a structure similar to the one of the Vlasov-Poisson system.
The analogy is not exact, as the Poisson equation (\ref{eq:spectral_vorticity_2}) is two-dimensional. \\

In Fourier representation $\omega$ and $\psi$ take the form
\begin{subequations}\label{eq:spectral_vorticity_3}
\begin{align}
\omega (r, t) &= \sum \limits_{k} \exp \left\{ - i \, k \cdot r \right\} \, \ohat{\omega} (k, t) , \\
\psi (r, t) &= \sum \limits_{k} \exp \left\{ - i \, k \cdot r \right\} \, \ohat{\psi} (k, t) ,
\end{align}
\end{subequations}

where $k \in \msp{Z}^{2}$ and $r = (x,y) \in \msp{I}^{2}$.
With the shorthand notation $\omega_{k}$ for $\ohat{\omega} (k, t)$ and $\psi_{k}$ for $\ohat{\psi} (k, t)$, the vorticity and Poisson equations become
\begin{align}
\label{eq:spectral_vorticity_4}
\dfrac{\partial \omega_{k}}{\partial t} &= \ohat{z} \cdot \sum \limits_{p + q = k} ( p \times q ) \; \omega_{q} \, \psi_{p} , \\
\label{eq:spectral_vorticity_5}
- k^{2} \, \psi_{k} &= \omega_{k} .
\end{align}

The reality condition
\begin{align}\label{eq:spectral_vorticity_6}
\omega (- k, t) = \omega^{*} (k, t)
\end{align}

determines the Fourier modes in the lower half plane and will play an important role in the analysis of conserved quantities.

\section{Conservation Laws}

We will concentrate on a particular invariant of the spectral vorticity equation, the $L^{2}$ norm of the vorticity
\begin{align}
G = \sum \limits_{k} \abs{ \omega_{k} }^{2} = \int \limits_{\msp{I}^{2}} \omega^{2} \, dx \, dy .
\end{align}

In physical space, $dG/dt = 0$ is a consequence of the anti-symmetry of the Poisson brackets in (\ref{eq:spectral_vorticity_1}).
In Fourier space, the conservation of $G$ follows by multiplying the spectral vorticity equation with $\omega_{- k}$
\begin{align}
\omega_{- k} \, \dfrac{\partial \omega_{k}}{\partial t} = \ohat{z} \cdot \sum \limits_{p + q = k} p \times q \; \omega_{- k} \, \omega_{q} \, \psi_{p}
\end{align}
and adding of the result to the same equation with the sign of $k$ flipped
\begin{align}
\omega_{k} \, \dfrac{\partial \omega_{- k}}{\partial t} = \ohat{z} \cdot \sum \limits_{p + q = - k} p \times q \; \omega_{k} \, \omega_{q} \, \psi_{p}
\end{align}
to get
\begin{align}
\omega_{- k} \, \dfrac{\partial \omega_{k}}{\partial t}
+ \omega_{k} \, \dfrac{\partial \omega_{- k}}{\partial t}
= \ohat{z} \cdot \bigg[
\sum \limits_{p + q = k} p \times q \; \omega_{- k} \, \omega_{q} \, \psi_{p}
+ \sum \limits_{p + q = - k} p \times q \; \omega_{k} \, \omega_{q} \, \psi_{p}
\bigg] .
\end{align}

With the reality condition, the left-hand side becomes
\begin{align}
\omega_{- k} \, \dfrac{\partial \omega_{k}}{\partial t}
+ \omega_{k} \, \dfrac{\partial \omega_{- k}}{\partial t}
= \dfrac{\partial}{\partial t} \big( \omega_{k}^{*} \omega_{k} \big)
= \dfrac{\partial}{\partial t} \abs{ \omega_{k} }^{2}
.
\end{align}

We sum the full equation over $k$ and rewrite the sums on the right-hand side
\begin{align}
\dfrac{\partial}{\partial t} \sum \limits_{k} \abs{ \omega_{k} }^{2}
\nonumber
&= \ohat{z} \cdot \sum \limits_{k} \sum \limits_{q} \bigg[
(k - q) \times q \; \omega_{- k} \, \omega_{q} \, \psi_{k - q}
+ (- k - q) \times q \; \omega_{k} \, \omega_{q} \, \psi_{- k - q}
\bigg] \\
&= \ohat{z} \cdot \sum \limits_{k} \sum \limits_{q}
k \times q \, \bigg[
\omega_{- k} \, \omega_{q} \, \psi_{k - q} - \omega_{k} \, \omega_{q} \, \psi_{- k - q}
\bigg] .
\end{align}

On the right-hand side, exchange $k$ and $q$ (we can do this as both, the sum of $k$ and the sum of $q$, run over the whole wave number space) and add the result to the original equation
\begin{align}
\dfrac{\partial}{\partial t} \sum \limits_{k} \abs{ \omega_{k} }^{2}
\nonumber
&= \dfrac{1}{2} \, \ohat{z} \cdot \sum \limits_{k} \sum \limits_{q} \Bigg\lgroup
k \times q \, \bigg[
\omega_{- k} \, \omega_{q} \, \psi_{k - q}
- \omega_{k} \, \omega_{q} \, \psi_{- k - q}
\bigg]
+ q \times k \, \bigg[
\omega_{- q} \, \omega_{k} \, \psi_{q - k}
- \omega_{q} \, \omega_{k} \, \psi_{- q - k}
\bigg] \Bigg\rgroup \\
\nonumber
&= \dfrac{1}{2} \, \ohat{z} \cdot \sum \limits_{k} \sum \limits_{q}
k \times q \, \bigg[
\omega_{- k} \, \omega_{q} \, \psi_{k - q}
- \omega_{k} \, \omega_{q} \, \psi_{- k - q}
+ \omega_{q} \, \omega_{k} \, \psi_{- q - k}
- \omega_{- q} \, \omega_{k} \, \psi_{q - k}
\bigg] \\
&= \dfrac{1}{2} \, \ohat{z} \cdot \sum \limits_{k} \sum \limits_{q}
k \times q \, \bigg[
\omega_{- k} \, \omega_{q} \, \psi_{k - q}
- \omega_{- q} \, \omega_{k} \, \psi_{q - k}
\bigg] .
\end{align}

Now change the sign of both $k$ and $q$, and add the result to the original equation
\begin{align}
\dfrac{\partial}{\partial t} \sum \limits_{k} \abs{ \omega_{k} }^{2}
\nonumber
&= \dfrac{1}{4} \, \ohat{z} \cdot \sum \limits_{k} \sum \limits_{q}
k \times q \, \bigg[
\omega_{- k} \, \omega_{q} \, \psi_{k - q}
- \omega_{q} \, \omega_{- k} \, \psi_{- q + k}
+ \omega_{k} \, \omega_{- q} \, \psi_{- k + q}
- \omega_{- q} \, \omega_{k} \, \psi_{q - k}
\bigg] \\
&= 0 .
\end{align}

Thus we have proved that $G$ is conserved.

\section{Extended Lagrangian}

The extended Lagrangian one-form for the spectral vorticity equation is
\begin{align}
\mcal{L} \big( \omega_{k}, \xi_{k} \big) = \xi_{k} \, \bigg[ \dfrac{\partial \omega_{k}}{\partial t} - \ohat{z} \cdot \sum \limits_{p + q = k} p \times q \; \omega_{q} \, \psi_{p} \bigg] \, dt .
\end{align}

The Poisson equation (\ref{eq:spectral_vorticity_5}) is not explicitly time dependent and therefore not included.
The variational principle applied to the action
\begin{align}
\mcal{A} = \int \mcal{L} \big( \omega_{k}, \xi_{k} \big)
\end{align}

yields two equations
\begin{subequations}
\begin{align}
\dfrac{\partial \omega_{k}}{\partial t} &= \ohat{z} \cdot \sum \limits_{p + q = k} p \times q \; \omega_{q} \, \psi_{p} , \\
\dfrac{\partial \xi_{k}}{\partial t} &= \ohat{z} \cdot \sum \limits_{p + q = k} p \times q \; \xi_{q} \, \psi_{p} ,
\end{align}
\end{subequations}

that are identical. So we can assume that $\xi$ has the same solution as $\omega$, which is important for the analysis of conserved quantities with Noether's theorem.

\section{Variational Integrator}

In the discrete system we apply the Fourier-Galerkin truncation where the sum is taken only over a finite subset of all possible wave numbers $k$. We will only discretise time with the help of a discrete variational principle. \\

The discrete extended Lagrangian density is
\begin{align}
\mcal{L}_{d} \big( \omega_{k}^{1}, \omega_{k}^{2}; \xi_{k}^{1}, \xi_{k}^{2} \big)
= \dfrac{1}{2} \, \Big( \xi_{k}^{1} + \xi_{k}^{2} \Big) \bigg[ \dfrac{\omega_{k}^{2} - \omega_{k}^{1}}{h_{t}} - \dfrac{1}{4} \, \ohat{z} \cdot \sum \limits_{p + q = k} p \times q \; \Big( \omega_{q}^{1} + \omega_{q}^{2} \Big) \Big( \psi_{p}^{1} + \psi_{p}^{2} \Big) \bigg] \, h_{t} .
\end{align}

With that one finds the discrete Euler-Lagrange equations
\begin{align}
\dfrac{\partial \mcal{L}_{d}}{\partial \xi^{1}} (\omega^{j}, \omega^{j+1}) + \dfrac{\partial \mcal{L}_{d}}{\partial \xi^{2}} (\omega^{j-1}, \omega^{j}) = 0
\end{align}
to be
\begin{align}
\dfrac{\omega_{k}^{j+1} - \omega_{k}^{j}}{2 h_{t}}
+\dfrac{\omega_{k}^{j} - \omega_{k}^{j-1}}{2 h_{t}}
= \dfrac{1}{8} \, \ohat{z} \cdot \sum \limits_{p + q = k} p \times q \; \bigg[
\Big( \omega_{q}^{j-1} + \omega_{q}^{j} \Big) \Big( \psi_{p}^{j-1} + \psi_{p}^{j} \Big)
+ \Big( \omega_{q}^{j} + \omega_{q}^{j+1} \Big) \Big( \psi_{p}^{j} + \psi_{p}^{j+1} \Big)
\bigg] .
\end{align}

This corresponds to the sum of two equations (c.f. the discussion in section \ref{sec:kinetic_theory_vi_simplifications}),
\begin{subequations}
\begin{align}
\dfrac{\omega_{k}^{j+1} - \omega_{k}^{j}}{2 h_{t}}
&= \dfrac{1}{8} \, \ohat{z} \cdot \sum \limits_{p + q = k} p \times q \;
\Big( \omega_{q}^{j} + \omega_{q}^{j+1} \Big) \Big( \psi_{p}^{j} + \psi_{p}^{j+1} \Big) ,
\\
\dfrac{\omega_{k}^{j} - \omega_{k}^{j-1}}{2 h_{t}}
&= \dfrac{1}{8} \, \ohat{z} \cdot \sum \limits_{p + q = k} p \times q \;
\Big( \omega_{q}^{j-1} + \omega_{q}^{j} \Big) \Big( \psi_{p}^{j-1} + \psi_{p}^{j} \Big) .
\end{align}
\end{subequations}

If $\omega_{k}$ is a solution of the first equation, it is also a solution of the second equation, as well as of the original equation (for a detailed discussion see section \ref{sec:vlasov_variational}).
So the discrete equation we use to advance the spectral vorticity in time  is
\begin{align}
\dfrac{\omega_{k}^{j+1} - \omega_{k}^{j}}{h_{t}}
&= \dfrac{1}{4} \, \ohat{z} \cdot \sum \limits_{p + q = k} p \times q \;
\Big( \omega_{q}^{j} + \omega_{q}^{j+1} \Big) \Big( \psi_{p}^{j} + \psi_{p}^{j+1} \Big) .
\end{align}

The proof of conservation of the discrete $L^{2}$ norm follows exactly the same path as in the continuous case and is therefore omitted.

\section{The Vlasov-Poisson System}

The Vlasov-Poisson system (see also chapter \ref{ch:kinetic_theory})
\begin{align}
& \dfrac{\partial f}{\partial t} + v \cdot \dfrac{\partial f}{\partial x} - \dfrac{\partial \phi}{\partial x} \cdot \dfrac{\partial f}{\partial v} = 0 \\
& \Delta \phi = - \int f \, dv .
\end{align}

describes the dynamics of a charged particle system, characterised by the particle distribution function $f$, in an electrostatic potential $\phi$.
In Fourier representation, $f$ and $\phi$ read
\begin{align}
f (t, x, v) &= \sum \limits_{k} \exp \left\{ - i \, k \cdot x \right\} \, \ohat{f} (t, k, v) , \\
\phi (t, x) &= \sum \limits_{k} \exp \left\{ - i \, k \cdot x \right\} \, \ohat{\phi} (t, k) ,
\end{align}

and the Vlasov and  Poisson equations become
\begin{align}
& \dfrac{\partial \ohat{f}}{\partial t} - i \, v \cdot k \, \ohat{f} + i \, \ohat{\phi} \, k \cdot \dfrac{\partial \ohat{f}}{\partial v} = 0 , \\
& k^{2} \, \ohat{\phi} = - \int \ohat{f} \, dv .
\end{align}

The variational integrator is derived in the same way as for the vorticity equation, only that now, in addition to time, also the velocity dimensions are considered in the derivation.
However, the analysis of conservation laws is supposedly more complicated than in the previous case.

\chapter{Discretisation of Brackets}\label{ch:brackets}

In this appendix we want to discuss the discretisation of Poisson brackets, Lie-Poisson brackets and both finite and infinite dimensional Nambu brackets \cite{Nambu:1973}.
This treatment is based on ideas of \citeauthor{SalmonTalley:1989}. In \cite{SalmonTalley:1989}, they describe a general way of discretising Poisson brackets $[ \cdot , \cdot ]$ by a method based on a discrete functional derivative that is very similar to the variational integrator formalism. In \cite{Salmon:2005}, \citeauthor{Salmon:2005} generalises these ideas to infinite dimensional Nambu brackets (Nambu field brackets), which are related to Lie-Poisson brackets as they often appear in the Hamiltonian description of plasma physics models.

\section{Canonical Poisson Brackets}

The starting point for \citeauthor{SalmonTalley:1989} \cite{SalmonTalley:1989} is the rephrasing of the equation at hand, e.g.,
\begin{align}\label{eq:brackets_poisson_1}
\mcal{D} f(t,x,p) = 0 ,
\end{align}

where $\mcal{D}$ is any operator, in a weak formulation, that is
\begin{align}\label{eq:brackets_poisson_2}
\int g (x,p) \, \mcal{D} f(t,x,p) \, dx \, dp &= 0 &
& \text{for any test function $g(x,p)$} . &
\end{align}

If the solutions of (\ref{eq:brackets_poisson_2}) are regular and (\ref{eq:brackets_poisson_2}) vanishes for any $g(x,y)$ it is equivalent to the original equation (\ref{eq:brackets_poisson_1}).
Consider as an example the Vlasov equation from chapter \ref{ch:kinetic_theory},
\begin{align}\label{eq:brackets_poisson_3}
\partial_{t} f + [ f, h ] = 0 ,
\end{align}

for which the corresponding weak formulation reads
\begin{align}\label{eq:brackets_poisson_4}
\int g (x,p) \, \big( \partial_{t} f + [ f, h ] \big) \, dx \, dp = 0 .
\end{align}

Here, $f(t,x,p)$ is the distribution function and $h$ is the particle Hamiltonian.
This is almost identical to the extended Lagrangian (\ref{eq:vlasov_vi_action}) from chapter \ref{ch:kinetic_theory}.
What is missing is the integral over time. So from this formulation, a semi-discretisation of the phasespace part of the equation is obtained.
Such a semi-discretisation has also been considered by \citeauthor{Leon:2008} \cite{Leon:2008} in the framework of variational integrators.

In the following we want to concentrate on the brackets, ignoring the time derivative. Therefore, we consider the energy functional
\begin{align}\label{eq:brackets_poisson_5}
\mathcal{E} = \int g \, [f, h] \, dx \, dp ,
\end{align}

such that the Poisson bracket $[f, h]$ is retained by a functional derivative with respect to $g$, that is
\begin{align}\label{eq:brackets_poisson_6}
\dfrac{\delta \mathcal{E}}{\delta g} = [f, h] .
\end{align}

It is now important to realise, integrating by parts with appropriate boundary conditions, that the following expressions are identical
\begin{align}\label{eq:brackets_poisson_7}
\int g \, [ f, h ] \, dx \, dp
= \int f \, [ h, g ] \, dx \, dp
= \int h \, [ g, f ] \, dx \, dp .
\end{align}

The energy functional (\ref{eq:brackets_poisson_5}) can thus be written as a convex combination of these three expressions
\begin{align}\label{eq:brackets_poisson_8}
\mathcal{E} = \int \Big[ \alpha \, g \, [ f, h ] + \beta \, f \, [ h, g ] + \gamma \, h \, [ g, f ] \Big] \, dx \, dp ,
\end{align}

where $\alpha + \beta + \gamma = 1$.
This observation is important for the discretisation of the brackets. To retain the antisymmetry properties of the brackets on the discrete level \citeauthor{SalmonTalley:1989} observed that equal factors, $\alpha = \beta = \gamma = 1/3$, have to be used, such that
\begin{align}\label{eq:brackets_poisson_9}
\mathcal{E} = \dfrac{1}{3} \int \Big[ g \, [ f, h ] +  f \, [ h, g ] + h \, [ g, f ] \Big] \, dx \, dp .
\end{align}

\subsection{Discrete Poisson Brackets on a Rectangular Mesh}

The next step in the derivation is the discretisation of this integral.
\citeauthor{SalmonTalley:1989} discretise the derivatives along the diagonals of a grid cell as depicted below.

\begin{center}
\includegraphics[width=.3\textwidth]{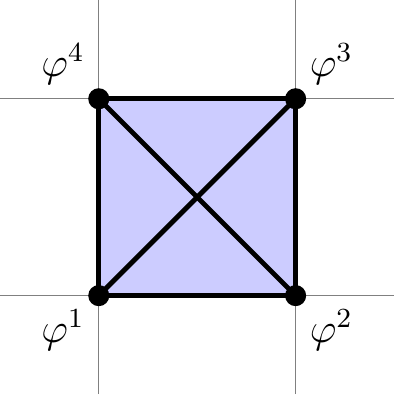}
\end{center}

Where the discrete derivatives are defined as
\begin{align}\label{eq:brackets_poisson_square_1}
\phy_{x}^{\boxtimes} (x,p) &\approx \dfrac{\phy^{2} - \phy^{4}}{h_{x}} , &
\phy_{p}^{\boxtimes} (x,p) &\approx \dfrac{\phy^{3} - \phy^{1}}{h_{p}} . &
\end{align}

This definition of the derivatives appears quite unnatural and ambiguous as both the $x$ and the $p$ derivative could be defined along each of the two diagonals.
The actual choice seems mostly motivated by the desired result, to obtain Arakawa's discretisation of the brackets (details follow below).
However, we will see that this specific discretisation is also obtained by using the discrete derivatives from chapter \ref{ch:variational} instead. I.e., we define the derivatives along the edges of the grid cell, that is
\begin{align}\label{eq:brackets_poisson_square_2}
\phy_{x}^{\square} (x,p) &\approx \dfrac{1}{2} \bigg( \dfrac{\phy^{2} - \phy^{1}}{h_{x}} + \dfrac{\phy^{3} - \phy^{4}}{h_{x}} \bigg) , &
\phy_{p}^{\square} (x,p) &\approx \dfrac{1}{2} \bigg( \dfrac{\phy^{4} - \phy^{1}}{h_{p}} + \dfrac{\phy^{3} - \phy^{2}}{h_{p}} \bigg) . &
\end{align}

As in chapter \ref{ch:variational}, the fields are approximated by averages over all four vertices of the grid cell,
\begin{align}\label{eq:brackets_poisson_square_3}
\phy^{\square} (x,p) &\approx \dfrac{1}{4} \Big( \phy^{1} + \phy^{2} + \phy^{3} + \phy^{4} \Big) .
\end{align}

Upon defining
\begin{multline}\label{eq:brackets_poisson_square_4}
S^{\square} (g,f,h)
\equiv \dfrac{1}{4} \Big( g^{1} + g^{2} + g^{3} + g^{4} \Big)
\bigg[
\dfrac{1}{2} \bigg( \dfrac{f^{2} - f^{1}}{h_{x}} + \dfrac{f^{3} - f^{4}}{h_{x}} \bigg)
\dfrac{1}{2} \bigg( \dfrac{h^{4} - h^{1}}{h_{p}} + \dfrac{h^{3} - h^{2}}{h_{p}} \bigg) \\
- \dfrac{1}{2} \bigg( \dfrac{f^{4} - f^{1}}{h_{p}} + \dfrac{f^{3} - f^{2}}{h_{p}} \bigg)
\dfrac{1}{2} \bigg( \dfrac{h^{2} - h^{1}}{h_{x}} + \dfrac{h^{3} - h^{4}}{h_{x}} \bigg)
\bigg]
,
\end{multline}

the discrete equivalent of the energy functional (\ref{eq:brackets_poisson_9}) becomes
\begin{align}\label{eq:brackets_poisson_square_5}
\mcal{E}_{d} = \dfrac{1}{3} \sum \limits_{\square} \Big( S^{\square} (g,f,h) + S^{\square} (f,h,g) + S^{\square} (h,g,f) \Big) .
\end{align}

The discrete Poisson brackets at a grid point $(i,j)$ are obtained by computing a discrete functional derivative, the same way as we computed the discrete variation, that is
\begin{align}\label{eq:brackets_poisson_square_6}
[f,h]_{i,j} = \dfrac{\partial \mcal{E}_{d}}{\partial g_{i,j}} .
\end{align}

Only four terms of the sum (\ref{eq:brackets_poisson_square_6}) include $g_{i,j}$, such that those four terms define the discrete Poisson bracket
\begin{align}\label{eq:brackets_poisson_square_7}
[f,h]_{i,j}
\nonumber
&= \dfrac{\partial \mcal{E}_d}{\partial g^1} \Big( \phy_{i,  j  }, \phy_{i+1,j  }, \phy_{i+1,j+1}, \phy_{i,  j+1} \Big)
+ \dfrac{\partial \mcal{E}_d}{\partial g^2} \Big( \phy_{i-1,j  }, \phy_{i,  j  }, \phy_{i,  j+1}, \phy_{i-1,j+1} \Big) \\
&+ \dfrac{\partial \mcal{E}_d}{\partial g^3} \Big( \phy_{i-1,j-1}, \phy_{i,  j-1}, \phy_{i,  j  }, \phy_{i-1,j  } \Big)
+ \dfrac{\partial \mcal{E}_d}{\partial g^4} \Big( \phy_{i,  j-1}, \phy_{i+1,j-1}, \phy_{i+1,j  }, \phy_{i,  j  } \Big)
.
\end{align}

The result of this computation is the well-known Arakawa scheme \cite{Arakawa:1966}.

\subsection{Arakawa's Discretisation}

In his original work \cite{Arakawa:1966}, \citeauthor{Arakawa:1966} considers different discretisations of the Poisson brackets with the aim of preserving the total number of particles
\begin{align}\label{eq:brackets_poisson_arakawa_1}
\int [f,h] \, dx \, dp &= 0 &
& \rightarrow &
\int f (t) \, dx \, dp &= \int f (0) \, dx \, dp , &
\end{align}

the total energy
\begin{align}\label{eq:brackets_poisson_arakawa_2}
\int [f,h] \, h \, dx \, dp &= 0 &
& \rightarrow &
\int f (t) \, h(t) \, dx \, dp &= \int f (0) \, h (0) \, dx \, dp , &
\end{align}

and the $L^{2}$ norm of the distribution function
\begin{align}\label{eq:brackets_poisson_arakawa_3}
\int [f,h] \, f \, dx \, dp &= 0 &
& \rightarrow &
\int f^{2} (t) \, dx \, dp &= \int f^{2} (0) \, dx \, dp . &
\end{align}

He defines four different discretisations of $J(f,h) = [f,h]$, that is
\begin{subequations}\label{eq:brackets_poisson_arakawa_4}
\begin{align}
J^{++}
&= \dfrac{1}{4 \, h_x \, h_p} \Big\lgroup \big( f_{+0} - f_{-0} \big) \big( h_{0+} - h_{0-} \big) - \big( f_{0+} - f_{0-} \big) \big( h_{+0} - h_{-0} \big) \Big\rgroup , \\
J^{+ \times}
\nonumber
&= \dfrac{1}{4 \, h_x \, h_p} \Big\lgroup f_{+0} \big( h_{+-} - h_{++} \big) - f_{-0} \big( h_{--} - h_{-+} \big) - f_{0+} \big( h_{-+} - h_{++} \big) + f_{0-} \big( h_{--} - h_{+-} \big) \Big\rgroup \\
&= \dfrac{1}{4 \, h_x \, h_p} \Big\lgroup h_{++} \big( f_{0+} - f_{+0} \big) - h_{--} \big( f_{-0} - f_{0-} \big) - h_{-+} \big( f_{0+} - f_{-0} \big) + h_{+-} \big( f_{+0} - f_{0-} \big) \Big\rgroup , \\
J^{\times +}
\nonumber
&= \dfrac{1}{4 \, h_x \, h_p} \Big\lgroup f_{++} \big( h_{+0} - h_{0+} \big) - f_{--} \big( h_{0-} - h_{-0} \big) - f_{-+} \big( h_{-0} - h_{0+} \big) + f_{+-} \big( h_{0-} - h_{+0} \big) \Big\rgroup \\
&= \dfrac{1}{4 \, h_x \, h_p} \Big\lgroup h_{+0} \big( f_{++} - f_{+-} \big) - h_{-0} \big( f_{-+} - f_{--} \big) - h_{0+} \big( f_{++} - f_{-+} \big) + h_{0-} \big( f_{+-} - f_{--} \big) \Big\rgroup , \\
J^{\times \times}
&= \dfrac{1}{8 \, h_x \, h_p} \Big\lgroup \big( f_{++} - f_{--} \big) \big( h_{-+} - h_{+-} \big) - \big( f_{-+} - f_{+-} \big) \big( h_{++} - h_{--} \big) \Big\rgroup ,
\end{align}
\end{subequations}

where the subscript $00$ refers to the grid point where the brackets are defined, $+0$ the grid point on the right, $-0$ the grid point on the left, and so on, as depicted below.

\begin{center}
\includegraphics[width=.3\textwidth]{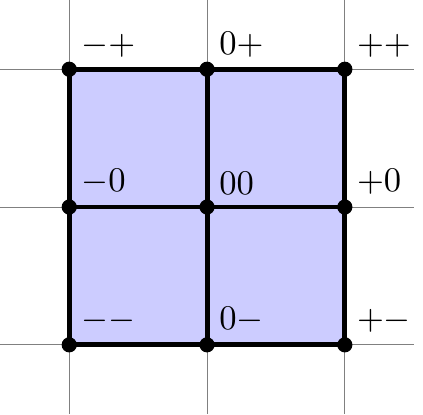}
\end{center}

Arakawa then explores different combinations of those expressions to find that all three of the above conservation properties (\ref{eq:brackets_poisson_arakawa_1}) - (\ref{eq:brackets_poisson_arakawa_3}) are only respected for
\begin{align}\label{eq:brackets_poisson_arakawa_5}
J = \dfrac{1}{3} \Big( J^{++} + J^{+ \times} + J^{\times +} \Big) ,
\end{align}

which is exactly the expression we obtain from (\ref{eq:brackets_poisson_square_7}).
The above expressions correspond to a second order discretisation of the Poisson bracket. Arakawa also provides fourth order expressions. It would be very interesting to see, if these can be derived by a similar approach. This, however, is a problem left for future research.

\subsection{Discrete Poisson Brackets on a Triangular Mesh}

A similar discretisation of the brackets can be performed on a mesh of triangles, leading to the scheme of \citeauthor{Sadourny:1968} \cite{Sadourny:1968}.
The only complication comes with the fact that we have two kinds of triangles, namely those pointing upward and those pointing downward, and we have, of course, to consider the contribution of both.

\begin{center}
\includegraphics[width=.18\textwidth]{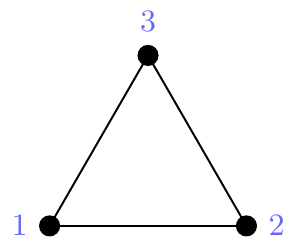}
\hspace{2cm}
\includegraphics[width=.18\textwidth]{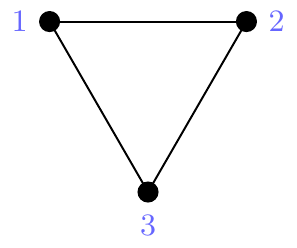}
\end{center}

Therefore we need to define cell averages and discrete derivatives for each type of triangle separately.
The grid points of the triangles are given by
\begin{align}\label{eq:brackets_poisson_triangle_1}
\triangle &= \Big( (i,j), (i+1,j), (i,j+1) \Big) &
\triangledown &= \Big( (i,j), (i+1,j), (i+1,j-1) \Big) , &
\end{align}

such that the vertices are denoted by
\begin{subequations}\label{eq:brackets_poisson_triangle_2}
\begin{align}
\triangle^{1} &= (i,j) , &
\triangle^{2} &= (i+1,j) , &
\triangle^{3} &= (i,j+1) , & \\
\triangledown^{1} &= (i,j) , &
\triangledown^{2} &= (i+1,j) , &
\triangledown^{3} &= (i+1,j-1) . &
\end{align}
\end{subequations}

Field averages are the same on both kind of triangles but the derivatives are different, that is
\begin{subequations}\label{eq:brackets_poisson_triangle_3}
\begin{align}
\phy^{\triangle}         &= \dfrac{1}{3} \left( \phy^{\triangle^{1}} + \phy^{\triangle^{2}} + \phy^{\triangle^{3}} \right) , &
\phy^{\triangledown}     &= \dfrac{1}{3} \left( \phy^{\triangledown^{1}} + \phy^{\triangledown^{2}} + \phy^{\triangledown^{3}} \right) ,
\\
\phy_{x}^{\triangle}     &= \dfrac{\phy^{\triangle^{2}} - \phy^{\triangle^{1}}}{h_x} , &
\phy_{x}^{\triangledown} &= \dfrac{\phy^{\triangledown^{2}} - \phy^{\triangledown^{1}}}{h_x} ,
\\
\phy_{p}^{\triangle}     &= \dfrac{1}{2} \bigg( \dfrac{\phy^{\triangle^{3}} - \phy^{\triangle^{1}}}{h_p} + \dfrac{\phy^{\triangle^{3}} - \phy^{\triangle^{2}}}{h_p} \bigg) , &
\phy_{p}^{\triangledown} &= \dfrac{1}{2} \bigg( \dfrac{\phy^{\triangledown^{1}} - \phy^{\triangledown^{2}}}{h_p} + \dfrac{\phy^{\triangledown^{1}} - \phy^{\triangledown^{3}}}{h_p} \bigg) .
\end{align}
\end{subequations}

Upon defining
\begin{subequations}\label{eq:brackets_poisson_triangle_4}
\begin{multline}\label{eq:brackets_poisson_triangle_4a}
S^{\triangle} (g,f,h)
\equiv \dfrac{1}{3} \Big( \phy^{\triangle^{1}} + \phy^{\triangle^{2}} + \phy^{\triangle^{3}} \Big)
\bigg[
\dfrac{1}{2} \bigg( \dfrac{f^{\triangle^{2}} - f^{\triangle^{1}}}{h_x} \bigg)
\bigg( \dfrac{h^{\triangle^{3}} - h^{\triangle^{1}}}{h_p} + \dfrac{h^{\triangle^{3}} - h^{\triangle^{2}}}{h_p} \bigg) \\
- \dfrac{1}{2} \bigg( \dfrac{f^{\triangle^{3}} - f^{\triangle^{1}}}{h_p} + \dfrac{f^{\triangle^{3}} - f^{\triangle^{2}}}{h_p} \bigg)
\bigg( \dfrac{h^{\triangle^{2}} - h^{\triangle^{1}}}{h_x} \bigg)
\bigg]
,
\end{multline}
\begin{multline}\label{eq:brackets_poisson_triangle_4b}
S^{\triangledown} (g,f,h)
\equiv \dfrac{1}{3} \Big( \phy^{\triangledown^{1}} + \phy^{\triangledown^{2}} + \phy^{\triangledown^{3}} \Big)
\bigg[
\dfrac{1}{2} \bigg( \dfrac{f^{\triangledown^{2}} - f^{\triangledown^{1}}}{h_x} \bigg)
\bigg( \dfrac{h^{\triangledown^{1}} - h^{\triangledown^{2}}}{h_p} + \dfrac{h^{\triangledown^{1}} - h^{\triangledown^{3}}}{h_p} \bigg) \\
- \dfrac{1}{2} \bigg( \dfrac{f^{\triangledown^{1}} - f^{\triangledown^{2}}}{h_p} + \dfrac{f^{\triangledown^{1}} - f^{\triangledown^{3}}}{h_p} \bigg)
\bigg( \dfrac{h^{\triangledown^{2}} - h^{\triangledown^{1}}}{h_x} \bigg)
\bigg]
,
\end{multline}
\end{subequations}

the discrete energy functional (\ref{eq:brackets_poisson_9}) becomes
\begin{align}\label{eq:brackets_poisson_triangle_5}
\mcal{E}_{d}
\nonumber
&= \sum \limits_{\triangle    } \dfrac{1}{3} \Big( S^{\triangle    } (g,f,h) + S^{\triangle    } (f,h,g) + S^{\triangle    } (h,g,f) \Big) \\
&+ \sum \limits_{\triangledown} \dfrac{1}{3} \Big( S^{\triangledown} (g,f,h) + S^{\triangledown} (f,h,g) + S^{\triangledown} (h,g,f) \Big) .
\end{align}

The discrete Poisson brackets at a grid point $(i,j)$ are obtained by computing the discrete functional derivative as in (\ref{eq:brackets_poisson_square_6}), that is
\begin{align}\label{eq:brackets_poisson_triangle_6}
[f,h]_{i,j} = \dfrac{\partial \mcal{E}_{d}}{\partial g_{i,j}} .
\end{align}

Six terms of the sum (\ref{eq:brackets_poisson_triangle_5}) include $g_{i,j}$, such that those terms define the discrete Poisson bracket
\begin{align}\label{eq:brackets_poisson_triangle_7}
[f,h]_{i,j}
\nonumber
&= \dfrac{\partial \mcal{E}_d}{\partial g^{\triangle^1    }} \Big( \phy_{i  , j  }, \phy_{i+1, j  }, \phy_{i  , j+1} \Big)
+ \dfrac{\partial \mcal{E}_d}{\partial g^{\triangle^2    }} \Big( \phy_{i-1, j  }, \phy_{i  , j  }, \phy_{i-1, j+1} \Big) \\
\nonumber
& \hspace{4em}
+ \dfrac{\partial \mcal{E}_d}{\partial g^{\triangle^3    }} \Big( \phy_{i  , j-1}, \phy_{i+1, j-1}, \phy_{i  , j  } \Big)
+ \dfrac{\partial \mcal{E}_d}{\partial g^{\triangledown^1}} \Big( \phy_{i  , j  }, \phy_{i  , j+1}, \phy_{i-1, j+1} \Big) \\
& \hspace{8em}
+ \dfrac{\partial \mcal{E}_d}{\partial g^{\triangledown^2}} \Big( \phy_{i  , j-1}, \phy_{i  , j  }, \phy_{i-1, j  } \Big)
+ \dfrac{\partial \mcal{E}_d}{\partial g^{\triangledown^3}} \Big( \phy_{i+1, j-1}, \phy_{i+1, j  }, \phy_{i  , j  } \Big)
.
\end{align}

The result of this computation is the scheme of \citeauthor{Sadourny:1968} \cite{Sadourny:1968},
\begin{align}\label{eq:brackets_poisson_triangle_9}
[f, h]_{d} &= \dfrac{1}{A_{\hexagon}} \sum \limits_{a=1}^{6} \dfrac{1}{2} \, f_{a} \, \Big( h_{a+1} - h_{a-1} \Big)
\end{align}

where the area of the hexagon, $A_{\hexagon}$, is
\begin{align}\label{eq:brackets_poisson_triangle_10}
A_{\hexagon} = 6 \, A_{\triangle} = 3 \, h_{x} \, h_{p} \, \cos 30 \degree ,
\end{align}

$a$ denotes the vertices of the hexagon as depicted below, $a = 6+1$ is replaced with $a = 1$, and $a = 1-1$ is replaced with $a=6$.

\begin{center}
\includegraphics[width=.4\textwidth]{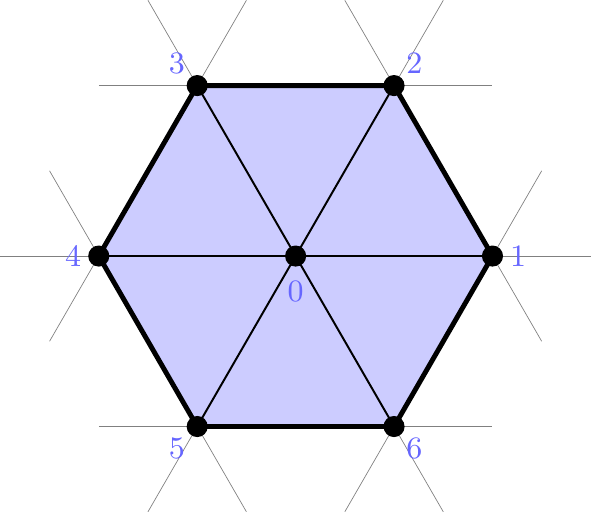}
\end{center}

\section{Nambu Three Brackets}

The ideas of the previous section can be extended to the three brackets defined by \citeauthor{Nambu:1973} \cite{Nambu:1973}
\begin{align}
[ f, g, h ]
&\equiv \epsilon^{abc} \, f_{,a} \, g_{,b} \, h_{,c}
= \dfrac{\partial f}{\partial x} \, \big[ g, h \big]_{yz} + \dfrac{\partial f}{\partial y} \, \big[ g, h \big]_{zx} + \dfrac{\partial f}{\partial z} \, \big[ g, h \big]_{xy} , &
& (a,b,c) \in \{ x, y, z \} ,
\end{align}

where $[ \cdot , \cdot ]_{xy}$ denote Poisson brackets with respect to $x$ and $y$, etc..
The energy functional is completely analogous to the previous one, that is
\begin{align}
\mathcal{E}
&= \int k \, [ f, g, h ] \, dx \, dy \, dz \\
&= \dfrac{1}{4} \int \Big\lgroup k \, [ f, g, h ] + f \, [ k, h, g ] + g \, [k, f, h] + h \, [ k, g, f ] \Big\rgroup dx \, dy \, dz ,
\end{align}

where we applied the same symmetrisation ideas as before.

\subsection{Discrete Nambu Brackets}

The discrete derivatives are defined similar as before, only that now we have to consider a three dimensional grid, such that the field averages and discrete derivatives become
\begin{subequations}\label{eq:brackets_nambu_discrete_derivatives}
\begin{align}
\phy^{\square} &= \dfrac{1}{8} \left( \phy^1 + \phy^2 + \phy^3 + \phy^4 + \phy^5 + \phy^6 + \phy^7 + \phy^8 \right) \\
\phy_{x}^{\square} &= \dfrac{1}{4} \left( \dfrac{\phy^2 - \phy^1}{h_x} + \dfrac{\phy^3 - \phy^4}{h_x} + \dfrac{\phy^6 - \phy^5}{h_x} + \dfrac{\phy^7 - \phy^8}{h_x} \right) , \\
\phy_{y}^{\square} &= \dfrac{1}{4} \left( \dfrac{\phy^4 - \phy^1}{h_y} + \dfrac{\phy^3 - \phy^2}{h_y} + \dfrac{\phy^8 - \phy^5}{h_y} + \dfrac{\phy^7 - \phy^6}{h_y} \right) , \\
\phy_{z}^{\square} &= \dfrac{1}{4} \left( \dfrac{\phy^5 - \phy^1}{h_z} + \dfrac{\phy^6 - \phy^2}{h_z} + \dfrac{\phy^7 - \phy^3}{h_z} + \dfrac{\phy^8 - \phy^4}{h_z} \right) .
\end{align}
\end{subequations}

The discrete action follows exactly along the lines of the previous derivations.
The exact form of the discrete brackets is quite complex and therefore not explicitly repeated here.

\subsection{Application to Gyrokinetics}

An application of this formulation and our discretisation is the gyrokinetic Vlasov equation on an extruded triangular mesh as depicted below \cite{Scott:2010},
\begin{align}
\dfrac{\partial f}{\partial t} + \dfrac{1}{\sqrt{g} \, B_{\parallel}^{*}} \, \Big[ h, f, A_{\phy}^{*} \Big]_{xyp_z} = 0 ,
\end{align}

where $(x,y)$ are coordinates in the poloidal plane of an axisymmetric tokamak, $\phy$ is the toroidal coordinate and $p_z$ is the parallel momentum, $g$ is the metric, $f$ the distribution function and $h$ the particle Hamiltonian, $A^{*} = \, A + \tfrac{c}{e} p_{z} b$ is the generalised vector potential and $B^{*}_{\parallel} = \, b \cdot ( \nabla \times A^{*} )$ the parallel magnetic field strength.

\begin{center}
\includegraphics[width=.45\textwidth]{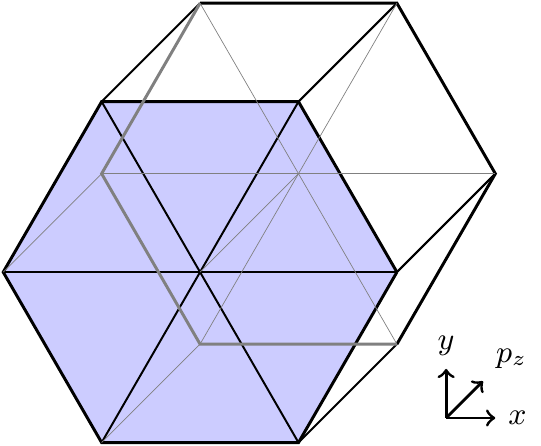}
\end{center}

Here, we have to combine the ideas of the previous sections. The spatial derivatives are defined in the poloidal plane over triangles and averaged in the momentum coordinate. Again we have to consider upward pointing and downward pointing triangles separately. The momentum derivative is defined in the extruded direction and averaged over the vertices of each triangle. Field averages are defined over all vertices of a three-dimensional grid cell. The rest of the derivation follows straight forwardly.

\section{Lie-Poisson and Nambu Field Brackets}

The natural Hamiltonian formulation of the Vlasov equation is via Lie-Poisson brackets. As we discussed in chapter five, the Vlasov-Poisson system does not have a canonical Hamiltonian formulation, nevertheless it is a Hamiltonian system. In this section we want to explore a discretisation approach for such Lie-Poisson brackets.

After sketching the noncanonical Hamiltonian description of systems like Vlasov-Poisson, we review \citeauthor{Salmon:2005}'s approach \cite{Salmon:2005} of discretising Lie-Poisson brackets. He does not discretise the brackets directly, but first finds equivalent infinite dimensional Nambu brackets (Nambu field brackets) and then discretises those by similar ideas as they were presented in the previous section.

\subsection{Noncanonical Hamiltonian Field Theory}

The dynamics of Hamiltonian systems is usually expressed with the help of canonically conjugate variables $(q,p)$ and Hamilton's equations
\begin{align}\label{eq:brackets_noncanonical_1}
\dot{q} &= \dfrac{\partial H}{\partial p} &
& \text{and} &
\dot{p} &= - \dfrac{\partial H}{\partial q} . & &&
\end{align}

The evolution of any functional $F(q,p)$ can be obtained with Poisson brackets, i.e.,
\begin{align}\label{eq:brackets_noncanonical_2}
\dot{F}(q,p) = [F,H] .
\end{align}

A large class of Hamiltonian systems, especially infinite-dimensional ones and especially systems from plasma physics (e.g., Vlasov-Poisson, reduced and ideal MHD, incompressible Fluid dynamics), do not fit into the form of (\ref{eq:brackets_noncanonical_1}).
They can, however, be described by a generalisation of (\ref{eq:brackets_noncanonical_2}).
The dynamics of a functional $F(\xi)$ of state variables $\xi (t, x)$ (e.g., distribution function, vorticity, density, temperature) of a Hamiltonian system is determined by
\begin{align}
\dot{F} (\xi) = \{ F, H \} ,
\end{align}

where $H(\xi)$ is the Hamiltonian functional and $\{ \cdot, \cdot \}$ are generalised Poisson brackets, that means they are antisymmetric and fulfil the Leibniz rule and the Jacobi identity.

\subsection{Lie-Poisson and Nambu Brackets in the Vlasov Equation}

The Vlasov equation can be express in terms of Lie-Poisson brackets \cite{MarsdenRatiu:2002} as
\begin{align}\label{eq:brackets_nambu_vlasov_1}
\dot{F} = \{ F, H \} \equiv \int f \, \bigg[ \dfrac{\delta F}{\delta f}, \dfrac{\delta H}{\delta f} \bigg] \, d x \, d p ,
\end{align}

where $F$ is any functional of $f$ and $H$ is the total energy functional
\begin{align}\label{eq:brackets_nambu_vlasov_2}
H = \int \dfrac{\vert p \vert^2}{2m} \, f(x, p) \, d x \, d p + \dfrac{1}{2} \int \vert \nabla \phi(x) \vert^2 \, d x .
\end{align}

Following \citeauthor{Salmon:2005}'s considerations for the vorticity equation \cite{Salmon:2005}, the Lie-Poisson bracket of the Vlasov equation can be expressed as a Nambu field bracket (see also \cite{BialynickiBirula:1991}). Therefore we just have to replace the single $f$ in (\ref{eq:brackets_nambu_vlasov_1}) with the functional derivative of the $L^{2}$ norm,
\begin{align}\label{eq:brackets_nambu_vlasov_3}
Z = \dfrac{1}{2} \int f^2 \, d x \, d p .
\end{align}

Hence, the Lie-Poisson bracket in (\ref{eq:brackets_nambu_vlasov_1}) becomes a Nambu three bracket,
\begin{align}\label{eq:brackets_nambu_vlasov_4}
\dot{F} = \{ F, H, Z \}
\equiv \int \dfrac{\delta Z}{\delta f} \, \bigg[ \dfrac{\delta F}{\delta f}, \dfrac{\delta H}{\delta f} \bigg] \, d x \, d p .
\end{align}

This bracket is antisymmetric in its three parameters, a property that is important in the discretisation procedure.

\subsection{Discretisation of Nambu Field Brackets}\label{sec:brackets_discrete_nambu}

The functionals $Z$ and $H$ are approximated by a simple quadrature rule as
\begin{align}
Z &= \dfrac{1}{2} \sum \limits_{i,j} f_{i,j}^2, &
H &= \sum \limits_{i,j} f_{i,j} \, h_{i,j} = \sum \limits_{i,j} f_{i,j} \, \Big( p_{j}^2 / m + q \phi_{i,j} \Big) . &
\end{align}

As in the previous sections, the key observation to a successful discretisation is the equivalence of the following expressions (integrating by parts with appropriate boundary conditions)
\begin{align}
\int \dfrac{\delta Z}{\delta f} \, \bigg[ \dfrac{\delta F}{\delta f}, \dfrac{\delta H}{\delta f} \bigg] \, d x \, d p
= \int \dfrac{\delta F}{\delta f} \, \bigg[ \dfrac{\delta H}{\delta f}, \dfrac{\delta Z}{\delta f} \bigg] \, d x \, d p
= \int \dfrac{\delta H}{\delta f} \, \bigg[ \dfrac{\delta Z}{\delta f}, \dfrac{\delta F}{\delta f} \bigg] \, d x \, d p .
\end{align}

The functional derivatives are defined on each vertex of a grid cell as we know it from the previous chapters and we define averages and derivatives just as before, c.f. equation (\ref{eq:brackets_nambu_discrete_derivatives}).

\begin{center}
\includegraphics[width=.3\textwidth]{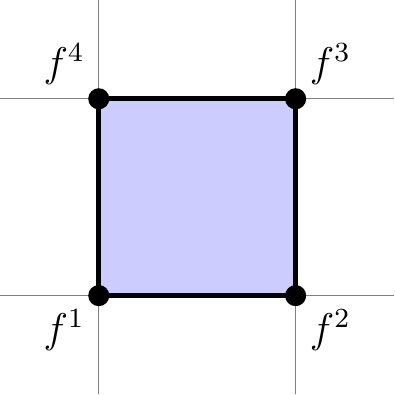}
\end{center}

The functional derivative of the $L^{2}$ norm is therefore discretised as
\begin{align}
\dfrac{\delta Z}{\delta f} \approx \dfrac{1}{4} \bigg( \dfrac{\partial Z}{\partial f^1} + \dfrac{\partial Z}{\partial f^2} + \dfrac{\partial Z}{\partial f^3} + \dfrac{\partial Z}{\partial f^4} \bigg) ,
\end{align}

and the Poisson brackets in \ref{eq:brackets_nambu_vlasov_4} are discretised by
\begin{multline}
\bigg[ \dfrac{\delta F}{\delta f}, \dfrac{\delta H}{\delta f} \bigg]_{d} =
\dfrac{1}{4 h_x h_p} \Bigg\lgroup
\bigg( \dfrac{\partial F}{\partial f^2} - \dfrac{\partial F}{\partial f^1} + \dfrac{\partial F}{\partial f^3} - \dfrac{\partial F}{\partial f^4} \bigg)
\bigg( \dfrac{\partial H}{\partial f^4} - \dfrac{\partial H}{\partial f^1} + \dfrac{\partial H}{\partial f^3} - \dfrac{\partial H}{\partial f^2} \bigg) \\
-
\bigg( \dfrac{\partial F}{\partial f^4} - \dfrac{\partial F}{\partial f^1} + \dfrac{\partial F}{\partial f^3} - \dfrac{\partial F}{\partial f^2} \bigg)
\bigg( \dfrac{\partial H}{\partial f^2} - \dfrac{\partial H}{\partial f^1} + \dfrac{\partial H}{\partial f^3} - \dfrac{\partial H}{\partial f^4} \bigg)
\Bigg\rgroup ,
\end{multline}

such that the discrete Nambu field bracket can be written as
\begin{align}
\{ F, H, Z \}_{\tilde{d}} &= \sum \limits_{\text{grid boxes}} \dfrac{1}{4} \bigg( \dfrac{\partial Z}{\partial f^1} + \dfrac{\partial Z}{\partial f^2} + \dfrac{\partial Z}{\partial f^3} + \dfrac{\partial Z}{\partial f^4} \bigg) \bigg[ \dfrac{\delta F}{\delta f}, \dfrac{\delta H}{\delta f} \bigg]_{d} .
\end{align}

To retain the antisymmetry property of the continuous Nambu bracket on the the discrete level, this expression has to be symmetrised, taking into account all even as well as all odd permutations
\begin{multline}
\{ F, H, Z \}_{d}
= \dfrac{1}{6} \Big( \{ F, H, Z \}_{\tilde{d}} + \{ H, Z, F \}_{\tilde{d}} + \{ Z, F, H \}_{\tilde{d}} \\
- \{ F, Z, H \}_{\tilde{d}} - \{ H, F, Z \}_{\tilde{d}} - \{ Z, H, F \}_{\tilde{d}} \Big) .
\end{multline}

With the discretisations we described above, considering only the even permutations suffices, but in general this is not the case.
The semi-discrete analogue of the Vlasov-Poisson equation is then
\begin{align}
\dfrac{\partial f_{ij}}{\partial t} = \{ f_{ij}, H, Z \}_{d} .
\end{align}

With our discretisation of the derivatives, this will again lead to the Arakawa discretisation.
It is an interesting observation that derivations on the Lagrangian side (variational integrators) and on the Hamiltonian side (discrete Nambu brackets) lead to similar discretisations of the equations of motion.
This is due to the fact that both the extended Lagrangian formulation and the Lie-Poisson brackets (\ref{eq:brackets_nambu_vlasov_1}) are constructed on the basis of the particle brackets $[ \cdot , \cdot ]$. This appears to be a consequence of the use of extended Lagrangians.

It will be interesting to further develop this approach as Lie-Poisson brackets exists for a wide range of system from plasma physics, like the Vlasov-Maxwell system and different flavours of magnetohydrodynamics.

\bibliographystyle{mynat}
\bibliography{thesis}

\end{document}